%% file: body.tex
\documentclass[12pt]{book}
\usepackage{amssymb,fancyheadings,makeidx}
\input{oldgerman}
\input{dissmacros}

\input{diagrams}
\makeindex
\begin{document}
\input{title}

\input{zusammenfassung}

\input{vorwort}

\pagestyle{fancyplain}
\headrulewidth0.2pt
\lhead[\fancyplain{}{\small\bf\sf\thepage}]%
{\fancyplain{}{\scriptsize\sf\rightmark}}
\rhead[\fancyplain{}{\scriptsize\sf\leftmark}]%
{\fancyplain{}{\small\bf\sf\thepage}}
\cfoot{}
\pagestyle{empty}
\input{contents}

\cleardoublepage
\pagestyle{fancyplain}
\headrulewidth0.2pt
\lhead[\fancyplain{}{\small\bf\sf\thepage}]%
{\fancyplain{}{\scriptsize\sf\rightmark}}
\rhead[\fancyplain{}{\scriptsize\sf\leftmark}]%
{\fancyplain{}{\small\bf\sf\thepage}}
\cfoot{}
\setcounter{page}{1}
\input{introduction}

\cleardoublepage
\input{algebra}
\input{rewriting}
\input{buchberger}
\input{presentations}

\cleardoublepage
\input{ideals}

\cleardoublepage
\input{strongreduction}
\input{rightreduction}
\input{prefixreduction}
\input{commutativereduction}

\cleardoublepage
\input{specialgroups}
\input{nilpotent}

\setcounter{lemma}{0}
\cleardoublepage
\input{reductionrings}

\cleardoublepage
\input{conclusions}


\cleardoublepage
\input{references}

\input{index}

\cleardoublepage
\input{lebenslauf}

\end{document}

%% file: oldgerman.tex
\makeatletter
\if \@ptsize 0
   \newfont{\scriptscriptscriptgoth}{ygoth scaled 760}
   \newfont{\scriptscriptgoth}{ygoth scaled 833}
   \newfont{\scriptgoth}{ygoth scaled 912}
   \newfont{\gothnomath}{ygoth}
   \newfont{\Goth}{ygoth scaled \magstephalf}
   \newfont{\GOth}{ygoth scaled \magstep1}
   \newfont{\GOTh}{ygoth scaled \magstep2}
   \newfont{\GOTH}{ygoth scaled \magstep3}

   \newfont{\scriptscriptscriptswab}{yswab scaled 760}
   \newfont{\scriptscriptswab}{yswab scaled 833}
   \newfont{\scriptswab}{yswab scaled 912}
   \newfont{\swab}{yswab}
   \newfont{\Swab}{yswab scaled \magstephalf}
   \newfont{\SWab}{yswab scaled \magstep1}
   \newfont{\SWAb}{yswab scaled \magstep2}
   \newfont{\SWAB}{yswab scaled \magstep3}

   \newfont{\scriptscriptscriptfrak}{yfrak scaled 760}
   \newfont{\scriptscriptfrak}{yfrak scaled 833}
   \newfont{\scriptfrak}{yfrak scaled 912}
   \newfont{\frak}{yfrak}
   \newfont{\Frak}{yfrak scaled \magstephalf}
   \newfont{\FRak}{yfrak scaled \magstep1}
   \newfont{\FRAk}{yfrak scaled \magstep2}
   \newfont{\FRAK}{yfrak scaled \magstep3}

   \newfont{\init}{yinit}
   \newfont{\Init}{yinit scaled \magstephalf}
   \newfont{\INit}{yinit scaled \magstep1}
   \newfont{\INIt}{yinit scaled \magstep2}
   \newfont{\INIT}{yinit scaled \magstep3}
\fi
\if \@ptsize 1
   \newfont{\scriptscriptscriptgoth}{ygoth scaled 833}
   \newfont{\scriptscriptgoth}{ygoth scaled 912}
   \newfont{\scriptgoth}{ygoth}
   \newfont{\gothnomath}{ygoth scaled \magstephalf}
   \newfont{\Goth}{ygoth scaled \magstep1}
   \newfont{\GOth}{ygoth scaled \magstep2}
   \newfont{\GOTh}{ygoth scaled \magstep3}
   \newfont{\GOTH}{ygoth scaled \magstep4}

   \newfont{\scriptscriptscriptswab}{yswab scaled 833}
   \newfont{\scriptscriptswab}{yswab scaled 912}
   \newfont{\scriptswab}{yswab}
   \newfont{\swab}{yswab scaled \magstephalf}
   \newfont{\Swab}{yswab scaled \magstep1}
   \newfont{\SWab}{yswab scaled \magstep2}
   \newfont{\SWAb}{yswab scaled \magstep3}
   \newfont{\SWAB}{yswab scaled \magstep4}

   \newfont{\scriptscriptscriptfrak}{yfrak scaled 833}
   \newfont{\scriptscriptfrak}{yfrak scaled 912}
   \newfont{\scriptfrak}{yfrak}
   \newfont{\fraknomath}{yfrak scaled \magstephalf}
   \newfont{\Frak}{yfrak scaled \magstep1}
   \newfont{\FRak}{yfrak scaled \magstep2}
   \newfont{\FRAk}{yfrak scaled \magstep3}
   \newfont{\FRAK}{yfrak scaled \magstep4}

   \newfont{\init}{yinit scaled \magstephalf}
   \newfont{\Init}{yinit scaled \magstep1}
   \newfont{\INit}{yinit scaled \magstep2}
   \newfont{\INIt}{yinit scaled \magstep3}
   \newfont{\INIT}{yinit scaled \magstep4}
\fi
\if \@ptsize 2
   \newfont{\scriptscriptscriptgoth}{ygoth scaled 912}
   \newfont{\scriptscriptgoth}{ygoth}
   \newfont{\scriptgoth}{ygoth scaled \magstephalf}
   \newfont{\gothnomath}{ygoth scaled \magstep1}
   \newfont{\Goth}{ygoth scaled \magstep2}
   \newfont{\GOth}{ygoth scaled \magstep3}
   \newfont{\GOTh}{ygoth scaled \magstep4}
   \newfont{\GOTH}{ygoth scaled \magstep5}

   \newfont{\scriptscriptscriptswab}{yswab scaled 912}
   \newfont{\scriptscriptswab}{yswab}
   \newfont{\scriptswab}{yswab scaled \magstephalf}
   \newfont{\swab}{yswab scaled \magstep1}
   \newfont{\Swab}{yswab scaled \magstep2}
   \newfont{\SWab}{yswab scaled \magstep3}
   \newfont{\SWAb}{yswab scaled \magstep4}
   \newfont{\SWAB}{yswab scaled \magstep5}

   \newfont{\scriptscriptscriptfrak}{yfrak scaled 912}
   \newfont{\scriptscriptfrak}{yfrak}
   \newfont{\scriptfrak}{yfrak scaled \magstephalf}
   \newfont{\fraknomath}{yfrak scaled \magstep1}
   \newfont{\Frak}{yfrak scaled \magstep2}
   \newfont{\FRak}{yfrak scaled \magstep3}
   \newfont{\FRAk}{yfrak scaled \magstep4}
   \newfont{\FRAK}{yfrak scaled \magstep5}

   \newfont{\init}{yinit scaled \magstep1}
   \newfont{\Init}{yinit scaled \magstep2}
   \newfont{\INit}{yinit scaled \magstep3}
   \newfont{\INIt}{yinit scaled \magstep4}
   \newfont{\INIT}{yinit scaled \magstep5}
\fi

\newcommand{\mswab}[1]{\mbox{\swab#1}}

\makeatother

%% file: dissmacros.tex
\newcommand{\auskommentieren}[1]{}
\newcommand{\betonen}[1]{{\bf #1}}

\newcommand{\m}{{\cal M}}      
\newcommand{\mm}{\circ}        
\newcommand{\cm}{\circ_{\cal T}}  
\newcommand{\factormonoid}{\mswab{M}}
\newcommand{\id}{\equiv} 
\newcommand{\free}{{\cal F}}
\newcommand{\freemonoid}{\Sigma^*}
\newcommand{\freecomm}{{\cal T}}
\newcommand{\class}{\mswab C}
\newcommand{\variety}{\mswab V}
\newcommand{\conc}{{\sf conc}\/}
\newcommand{\irr}{{\rm IRR}\/}
\newcommand{\pred}{\prec}

\newcommand{\g}{{\cal G}}      
\newcommand{\gm}{\circ}        
\newcommand{\inv}[1]{{\sf inv}\/(#1)}
\newcommand{\ord}{{\sf ORD}}

\newcommand{\n}{{\bf N}} 
\newcommand{\z}{{\bf Z}} 

\newcommand{\myk}{{\bf K}}
\newcommand{\q}{{\bf Q}}

\newcommand{\rr}{{\bf R}}

\newcommand{\mrm}{\ast} 
\newcommand{\skm}{\cdot} 
\newcommand{\cd}{\cdot}
\newcommand{\qmult}{\otimes}
\newcommand{\qadd}{\oplus}

\newcommand{\ideal}[2]{{\sf ideal}_{#1}^{#2}}
\newcommand{\spol}[1]{{\sf spol}_{#1}}
\newcommand{\can}{{\sf can}}
\newcommand{\satpoly}{{\sf acan}}
\newcommand{\sgn}{{\sf sgn}}

\newcommand{\hm}{{\sf HM}}
\newcommand{\hc}{{\sf HC}}
\newcommand{\hterm}{{\sf HT}}
\newcommand{\reductum}{{\sf RED}}
\newcommand{\terms}{{\sf T}}
\newcommand{\lcm}{{\sf LCM}}

\newcommand{\tup}{<_{\rm tup}}
\newcommand{\tupgreater}{>_{\rm tup}}
\newcommand{\tupeq}{\geq_{\rm tup}}
\newcommand{\tupleq}{\leq_{\rm tup}}
\newcommand{\syll}{>_{\rm syll}}
\newcommand{\syllless}{<_{\rm syll}}
\newcommand{\sylleq}{\geq_{\rm syll}}

\newcommand{\R}{\Longrightarrow}
\newcommand{\myr}{\longrightarrow}
\newcommand{\la}{\longleftarrow}
\newcommand{\lr}{\longleftrightarrow}

\newcommand{\red}[4]{\mbox{$\,\stackrel{#1}{#2}\!\!\mbox{}^{{\rm #3}}_{#4}\,$}}
\newcommand{\nred}[4]%
{\mbox{$\,\,\,\,{\not\!\!\!\stackrel{#1}{#2}\!\!\mbox{}^{{\rm #3}}_{#4}}\,$}}

\newcommand{\p}{{\cal P}}

\newcommand{\nf}[1]{\!\!\downarrow_{#1}} 
\newcommand{\rnf}[1]{\!\!\Downarrow_{#1}}
\newcommand{\SAT}{{\cal SAT}}

\newcommand{\Ba}[1]{{\bf Proof #1}: 
                      \renewcommand{\baselinestretch}{1.1}\small\normalsize}   
\newcommand{\qed}{\hspace*{\fill} q.e.d. \par\medskip
        \renewcommand{\baselinestretch}{1}\small\normalsize}    

\newtheorem{lemma}{Lemma}[section]
\newtheorem{example}[lemma]{Example}
\newtheorem{remark}[lemma]{Remark}
\newtheorem{definition}[lemma]{Definition}
\newtheorem{theorem}[lemma]{Theorem}
\newtheorem{corollary}[lemma]{Corollary}
\newcommand{\exaend}{\hfill $\diamond$}
\newcommand{\dend}{\hfill $\diamond$}
\newcommand{\remend}{\hfill $\diamond$}
\newcommand{\lemend}{}
\newcommand{\theoend}{}
\newcommand{\ohnebeweis}{\hfill $\square$}

\newcommand{\s}{\mbox{\sc Sat}}
\newcommand{\sm}{\mbox{\footnotesize\sc Sat}} 
\newcommand{\gb}{\mbox{\sc Gb}}

\newcommand{\procedure}[2]%
{\noindent{{\bf Procedure}: {\sc #1}}  \vspace{3mm}
\hrule height 0.5pt \hfill
#2
\hrule height 0.5pt \hfill}

\newcommand{\kommentar}{\small}

\newcommand{\spruch}[4]{
\begin{minipage}[b]{#1cm}
\mbox{   }
\end{minipage}
\begin{minipage}[b]{#2cm}
{\small
#3
}
\end{minipage}

\hfill{\footnotesize\sc #4}

\vspace{1cm}
}

%% file: title.tex
\begin{titlepage}

\begin{center}
\vspace{10cm}
{\huge\bf \mbox{       }    \\[5ex]
       On Gr\"obner Bases   \\[2ex]
              in            \\[3ex]
    Monoid and Group Rings} \\[7ex]

Vom Fachbereich Informatik \\[1ex]
der Universit{\"a}t Kaiserslautern \\[1ex]
zur Verleihung des akademischen Grades \\[1ex]
Doktor der Naturwissenschaften (Dr. rer. nat.)\\[1ex]
genehmigte Dissertation \\[1ex]
von \\[5ex]
 Dipl.-Math. Birgit Reinert \\[10ex]
\begin{tabbing}
XXXXXXXXXXXXX \=\kill\\
Datum der wissenschaftlichen Aussprache: 14. Juni 1995 \\
\\
Dekan:            \> Prof. Dr. Hans Hagen \\
\\
Promotionskommission: \\[1ex]
Vorsitzender: \> Prof. Dr. Theo H\"arder \\
Berichterstatter: \> Prof. Dr. Klaus E. Madlener\\
                  \> Prof. Dr. Volker Weispfenning
\end{tabbing}
\vspace{5mm}
D 386
\end{center}
\end{titlepage}

%% file: zusammenfassung.tex
\pagebreak
\thispagestyle{empty}
{\bf Zusammenfassung}
\\
\\
 Gr\"obnerbasen, entwickelt von Bruno Buchberger f\"ur  kommutative
 Polynomringe, finden  h\"aufig Anwendung bei der L\"osung algorithmischer Probleme.
Beispielsweise  l\"a{\ss}t sich  das Kongruenzproblem f\"ur Ideale
 mit Hilfe der Gr\"obnerbasen l\"osen.
Bis heute wurden diese Ideen auf verschiedene 
 zum Teil nichtkommutative und nichtnoethersche Algebren \"ubertragen.
Die meisten dieser Ans\"atze setzen eine zul\"assige Ordnung
 auf den Termen voraus.

In dieser Dissertation wird das Konzept der Gr\"obnerbasen
 f\"ur endlich erzeugte
 Monoid- und Gruppenringe verallgemeinert.
Dabei werden Reduktionsmethoden sowohl zur Darstellung der Monoid-
 beziehungsweise Gruppenelemente, als auch zur Beschreibung der
 Rechtsidealkongruenz in den entsprechenden Monoid- beziehungsweise
 Gruppenringen benutzt.
Da im allgemeinen Monoide und insbesondere Gruppen keine zul\"assigen
 Ordnungen mehr erlauben, treten bei der Definition einer geeigneten
 Reduktionsrelation wesentliche
 Probleme auf:
Zum einen ist es schwierig, die Terminierung einer Reduktionsrelation
 zu garantieren, zum anderen sind Reduktionsschritte nicht mehr mit
 Multiplikationen vertr\"aglich und daher beschreiben Reduktionen nicht
 mehr unbedingt eine Rechtsidealkongruenz.
In dieser Arbeit werden  verschiedene M\"oglichkeiten Reduktionsrelationen zu
 definieren aufgezeigt und im Hinblick auf die beschriebenen Probleme untersucht.
Dabei wird das Konzept der Saturierung, d.h. eine Polynommenge so zu
erweitern, da{\ss} man die von ihr erzeugte Rechtsidealkongruenz durch
Reduktion erfassen kann, benutzt, um Charakterisierungen
von Gr\"obnerbasen bez\"uglich der verschiedenen Reduktionen durch
s-Polynome zu geben.
Mithilfe dieser Konzepte ist es gelungen f\"ur spezielle Klassen von
Monoiden, wie z.B. endliche, kommutative oder freie,  und verschiedene Klassen von
Gruppen, wie z.B. endliche, freie, plain, kontext-freie oder nilpotente, unter
Ausnutzung struktureller Eigenschaften spezielle Reduktionsrelationen
zu definieren und terminierende Algorithmen zur Berechnung von Gr\"obnerbasen
bez\"uglich dieser Reduktionsrelationen zu 
entwickeln.
\vspace*{4cm}

%% file: vorwort.tex
\pagebreak
\thispagestyle{empty}
\auskommentieren{
\spruch{4.2}{11}
{Im \"ubrigen, mein Sohn, la{\ss} dich warnen!
 Es nimmt kein Ende mit dem vielen B\"ucherschreiben,
 und viel Studieren erm\"udet den Leib.
 Hast du alles geh\"ort, so lautet der Schlu{\ss}:
 F\"urchte Gott, und achte auf seine Gebote!
 Das allein hat jeder Mensch n\"otig.}{Kohelet, 12:12-13}
\vspace{1cm}}
{\bf Vorwort}
\\
\\
Mein Dank gilt allen, die auf vielf\"altige Weise
 zum Gelingen meiner Promotion beigetragen haben.

Professor Madlener hat mir die Gelegenheit zur Erstellung dieser Arbeit gegeben.
Ich danke ihm f\"ur alles, was ich fachlich und menschlich durch ihn lernen durfte.
Professor Weispfenning danke ich f\"ur die \"Ubernahme der Zweitbegutachtung meiner Arbeit.
Professor H\"arder hat freundlicherweise die Leitung der Promotionskommission
 \"ubernommen.

Mein Dank gilt auch den Mitgliedern der
 Arbeitsgruppen von Professor Madlener und Professor Avenhaus,
 sowie Frau Rita Kohl, die das
 Entstehen meiner Arbeit mit Anteilnahme mitverfolgt und unterst\"utzt haben.
Insbesondere Andrea Sattler-Klein, Inger Sonntag,  Thomas Dei{\ss},
 Roland Fettig und Claus-Peter Wirth waren immer f\"ur mich da und haben mir
 stets geduldig zugeh\"ort.

Meinen Freunden und Verwandten danke ich f\"ur die Ermunterungen und Gebete.
Gewidmet ist diese Arbeit meinen Eltern und Joachim.
Ihre Liebe hat mir geholfen, den Blick
 f\"ur das Wesentliche nicht zu verlieren.

\vspace{1.5cm}
\spruch{11}{7}{ Nada te turbe, \\
 nada te espante, \\
 todo se pasa.\\

 Dios no se muda. \\
 La paciencia \\
 todo lo alcanza. \\

 Quien a Dios tiene, \\
 nada le falta: \\
 s\'olo Dios basta.\\}{Teresa de Jes\'us}

\pagebreak
\thispagestyle{empty}
\vspace*{15cm}
\hspace*{\fill}\parbox{6cm}{%
\begin{center}
{\em F\"ur meine Eltern} \\
{\em und Joachim}
\end{center}}

%% file: contents.tex
\tableofcontents
\newpage

%% file: introduction.tex
\chapter{Introduction}
\spruch{11}{7}{Was wir sind, ist nichts,\\
was wir suchen, ist alles.}{Lessing}

One of the amazing features of computers is the ability to discover
new mathematical results due to extensive computations impossible to
be done by hand.
Besides incredible numerical calculations, symbolical mathematical
manipulations are substantial to many fields in mathematics and
physics.
Hence the idea of using a computer to do such manipulations led to
open up whole new areas of mathematics and computer science.
In the wake of these developments has come a new access to abstract
algebra in a computational fashion - computer algebra.
One important contribution is Buchberger's
algorithm for manipulating systems of polynomial equations.
In 1965 Buchberger introduced the theory of Gr\"obner bases for
 polynomial ideals in commutative
 polynomial rings over fields (see \cite{Bu65}).
It established a rewriting approach to the theory of polynomial
 ideals.
Polynomials can be used as rules by giving an admissible\footnote{A term
 ordering $\succeq$ is called admissible if for every term $s,t,u$,
 $s \succeq 1$ holds, and $s \succeq t$ implies $s \mm u \succeq t \mm u$.
 An ordering fulfilling the latter condition is also said to be
 compatible with the respective multiplication $\mm$.} ordering
 on the terms and using the largest monomial according to this
 ordering as a left hand side of a rule.
``Reduction'' as defined by Buchberger then can be compared to
 division of one polynomial by a
 set of finitely many polynomials.
A Gr\"obner basis $G$ is a  set of  polynomials
 such that every polynomial in the polynomial ring has a unique
 normal form with respect to reduction using the polynomials
 in $G$ as rules (especially the
 polynomials in the ideal generated by $G$ reduce to zero using $G$).
Buchberger developed a terminating procedure to transform a finite generating
 set of a polynomial ideal into a finite 
 Gr\"obner basis of the same ideal.

The method of Gr\"obner bases allows to solve many problems related to
 polynomial ideals in a computational fashion.
It was shown by Hilbert (compare Hilbert's basis theorem) that every
 ideal in a polynomial ring has a finite generating set.
However, an arbitrary finite generating set need not provide much
insight into the nature of the ideal.
Let $f_1 = X_1^{2} + X_2$ and $f_2 = X_1^{2} + X_3$ be two polynomials
 in the polynomial ring $\q[X_1,X_2,X_3]$.
Then  $\mswab{i} = \{ f_1 \mrm g_1 + f_2
\mrm g_2 | g_1, g_2 \in \q[X_1,X_2,X_3] \}$ is the ideal they generate 
and it is not hard to see that
the polynomial $X_2 - X_3$ belongs to $\mswab{i}$ since $X_2 - X_3 =
f_1 - f_2$.
But what can be said about the polynomial $f = X_3^{3} + X_1 + X_3$?
Does it belong to $\mswab{i}$ or not?
\\
The problem to decide whether a given polynomial lies in a given
 ideal is called the membership problem for ideals.
In case the generating set is a Gr\"obner basis this problem becomes
immediately solvable, as the membership problem then reduces to checking whether
 the polynomial reduces to zero.
\\
In our example the set $\{ X_1^{2} + X_3, X_2 - X_3 \}$ is a
generating set of $\mswab{i}$ which is in fact a Gr\"obner basis.
Now returning to the polynomial $f = X_3^{3} + X_1 + X_3$ we find that
it cannot
belong to $\mswab{i}$ since neither $X_1^{2}$ nor $X_2$ is a divisor
of a term in  $f$ and hence $f$ cannot be reduced to zero by the
polynomials in the Gr\"obner basis.
\\
Further applications of Gr\"obner bases to algebraic questions can be
 found  e.g. in the work of Buchberger \cite{Bu87}, Becker and
 Weispfenning  \cite{BeWe92} and in
 the book of Cox, Little and O'Shea \cite{CoLiOS92}.

Since the theory of Gr\"obner bases turned out to be of outstanding
 importance for polynomial rings, several generalizations of
 Buchberger's ideas to other structures followed.
We only want to give a brief outline of some of them, mainly of those
 which influenced our work.

A first generalization was given by Buchberger himself and his student 
 Stifter  in characterizing reduction rings by adding additional axioms
 to the ring axioms (\cite{St85,St87}).
Further characterizations
 of such reduction rings were provided by Kapur and Narendran
 (\cite{KaNa85}) and
 Madlener (\cite{Ma86}).

Besides these theoretical studies of reduction rings, the Gr\"obner
 basis theory has been extended to commutative polynomial rings over
 coefficient domains other than fields.
It was shown by authors as  Buchberger, Kandri-Rody, Kapur,
Narendran, Lauer, Stifter and Weispfenning that Buchberger's approach remains valid for polynomial
 rings over the integers, or even Euclidean rings, and over
 regular rings  and reduction rings 
 (see e.g. \cite{Bu83,Bu85,KaKa84,KaKa88,KaNa85,La76,St85,We87}).

Since the development of computer algebra systems for commutative
algebras enabled to perform tedious calculations using computers,
attempts to generalize such systems and especially Buchberger's ideas
to non-commutative algebras followed.
Originating from special problems in physics, Lassner in \cite{La85}
suggested how to extend existing computer algebra systems in order to
handle special classes of non-commutative algebras, e.g. Weyl algebras.
He studied structures where the elements could be represented using
the usual representation of polynomials in commutative variables and
the non-commutative multiplication could be performed by a so-called
``twisted product'' which required only procedures involving
commutative algebra operations and differentiation.
His extensions could be incorporated into computer algebra systems to
solve tasks of interest to physicists.
Later on together with Apel he extended Buchberger's algorithm to
enveloping fields of Lie algebras (see \cite{ApLa88}).
Because these ideas use representations by commutative polynomials,
Dickson's lemma can be carried over.
The existence and construction of Gr\"obner bases for finitely
generated left ideals is ensured.
On the other hand, Mora gave a concept of Gr\"obner bases for a class
of non-commutative algebras by saving an other property of the
polynomial ring while losing the validity of Dickson's lemma.
The usual polynomial ring can be viewed as a monoid ring where the
 monoid is a finitely generated free commutative monoid.
Mora studied the class where the free commutative monoid is substituted
by a free monoid - the class of finitely generated free monoid rings (compare e.g. \cite{Mo85,Mo93}).
The ring operations are mainly performed in the coefficient
 domain while the terms are treated like words, i.e., the variables no
 longer commute with each other.
The definitions of (one- and two-sided) ideals, reduction and Gr\"obner
 bases are carried over from the commutative  case to establish
 a similar theory of Gr\"obner bases in ``free non-commutative
 polynomial rings over fields''.
But these rings are no longer Noetherian if they are generated by more
 than one variable.
Moreover, the word problem for semi-Thue systems can be reduced  to the
membership problem for two-sided ideals.
Mora presented a terminating completion procedure for finitely generated one-sided ideals
and an enumeration procedure for finitely generated two-sided ideals
with respect to some term ordering in free monoid
rings.

Another class of non-commutative rings  where the elements can
be represented by the usual polynomials and which allow the construction of finite
Gr\"obner bases for arbitrary ideals are the so-called solvable rings, a class
intermediate between commutative and general non-commutative
polynomial rings.
They were studied by Kandri-Rody, Weispfenning and Kredel (\cite{KaWe90,Kr93}).
Solvable polynomial rings can be described by ordinary polynomial
rings $\myk[X_1, \ldots, X_n]$ provided with a ``new'' definition of multiplication which coincides with the
ordinary multiplication except for the case that a variable $X_j$ is
multiplied with a variable $X_i$ with lower index, i.e., $i<j$.
In the latter case multiplication can be defined by equations
$$X_j \star X_i = c_{ij} X_iX_j + p_{ij}$$
where $c_{ij} \in \myk^*= \myk \backslash \{ 0 \}$ and $p_{ij}$ is a polynomial
``smaller'' than $X_iX_j$ with respect to a fixed admissible term
ordering on the polynomial ring.

In \cite{We92} Weispfenning showed the existence of finite Gr\"obner
 bases for arbitrary finitely
 generated ideals in non-Noetherian skew polynomial rings over two
 variables $X,Y$ where a ``new'' multiplication $\star$ is introduced
 such that $X \star Y = XY$ and $Y \star X = X^eY$ for some fixed $e
 \in \n^+$.

Most of the approaches mentioned so far fulfill the following requirements:
\begin{itemize}
\item The rings allow admissible well-founded orderings.
\item If a polynomial can be reduced to zero by a set of polynomials,
  so can a multiple of this polynomial by a monomial.
\item The translation lemma holds, i.e., if the difference of two
  polynomials can be reduced to zero using a set of polynomials, then
  the two polynomials are joinable using the same set of polynomials
  for reduction.
\item Two polynomials give rise to finitely many s-polynomials only,
  in general at most one.
\end{itemize}
These statements can be used to characterize Gr\"obner bases in the
respective ring with respect to the corresponding reduction in a
finitary manner and  it is decidable whether a finite set is a Gr\"obner
basis by checking whether the s-polynomials are reducible to
zero\footnote{Note that we always assume that the reduction in the ring is effective.}.

There are rings combined with reduction where these statements cannot
 be accomplished and therefore other concepts to characterize
 Gr\"obner bases have been developed.
For example
 in case the ring contains zero-divisors a well-founded ordering on the ring
 is no
 longer compatible with the ring multiplication\footnote{When studying monoid rings over reduction rings it is possible
 that the ordering on the ring is not compatible with scalar multiplication
 as well as with multiplication with monomials or polynomials.}.
This phenomenon has been studied for the case of zero-divisors in
 the coefficient domain by Kapur and Madlener (\cite{KaMa86}) and by
 Weispfenning for the special case of regular rings (\cite{We87}).
In his PhD thesis \cite{Kr93}, Kredel mentioned that in dropping the 
 axioms guaranteeing the existence
 of admissible orderings in the theory of solvable polynomial rings by
 allowing $c_{ij} = 0$ in the defining equations above, the properties
 mentioned above need no longer hold.
He sketched the idea of using saturation to repair some of the
 problems occurring and for special cases, e.g. for the Grassmann
(exterior) algebras, positive results can be  achieved (compare the paper
 of Stokes \cite{St90}).

Starting point of the current work was the idea to study arbitrary
finitely generated monoid rings similar to Mora's approach to free
monoid rings.
Since we want to treat rings with well-founded but no longer admissible
orderings, we mainly have to deal with the fact that many of the
properties used to give a characterization of a Gr\"obner basis as in the
classical case no longer hold. 
We  show that there are weaker requirements that still enable
characterizations of Gr\"obner bases in terms of s-polynomials in case
we use appropriate reductions combined with appropriate concepts of saturation.
For special reductions we can even give a characterization where
localization of critical situations to one s-polynomial for each pair
of polynomials is possible.
 
Our approach can be characterized by the following tasks: 
\begin{itemize}
\item Combine string rewriting and polynomial rewriting in the field
  of monoid rings.
\item Generalize the concept of Gr\"obner bases to arbitrary monoid
  rings.
\item Find classes of monoids and groups that allow the construction of
  finite Gr\"obner bases for finitely generated one- or even two-sided ideals.
\end{itemize}

The thesis  organizes as follows:

Chapter \ref{chapter.definitions} introduces some of the basic themes of
this work.
We need some definitions and notions from algebra and the theory of
rewriting systems.
Furthermore, as this work is based on Buchberger's ideas, a short
summary of the theory of Gr\"obner bases and Buchberger's algorithm
are given.

Chapter \ref{chapter.ideals} gives a short outline on introducing
Gr\"obner bases to non-commutative structures by sketching Mora's
approach to free monoid rings and Weispfenning's approach to skew
polynomial rings.
We prove that the word problem for semi-Thue systems is equivalent to
a restricted version of the membership problem for free monoid rings and
similarly that the word problem for group presentations is
equivalent to
a restricted version of the membership problem for free group rings.
Hence for free monoids respectively free groups with more than one
generator the ideal membership problem is undecidable.
(It is decidable for one generator.)
Further we show that it is undecidable whether a finite Gr\"obner
basis in a free
monoid ring generated by more than one generator exists.

Chapter \ref{chapter.reduction} gives different approaches to define
reduction in monoid rings.
Since monoid rings in general are not commutative, we are mainly
interested in right ideals.
A well-founded ordering on terms is used to split a polynomial into a
head monomial (the largest monomial) and a reduct.
A natural way to define reduction is to use a multiple of a
polynomial  to replace a monomial in another polynomial in case the result
is smaller.
This is the case if the head term of the multiple equals the term of
the removed monomial.
This reduction is called strong reduction and can be used to express
the congruence of a right ideal.
Although a characterization of Gr\"obner bases with respect to this
reduction in terms of strong s-polynomials is possible, this
characterization is not finitary
and it cannot be used to decide whether a finite set of
polynomials is a strong Gr\"obner basis.
One idea to localize a confluence test, i.e., reduce the number of
s-polynomials to be considered, is to weaken reduction.
The first weakening studied is restricting the multiples of
polynomials used to stable ones, i.e., multiples where the new head
term results from the original head term of the polynomial.
The first problem is that the expressiveness of a right ideal
congruence by reduction using an arbitrary set of generators of the
right ideal is lost.
This can be regained using a concept called saturation (e.g.
mentioned in \cite{Kr93}), which enlarges
the set of polynomials used for reduction. But such saturating sets
need not be finite.
For saturated sets a characterization of Gr\"obner bases with respect
to right reduction in terms of right s-polynomials is provided which
is still not finitary.
Therefore, two weakenings involving syntactical information on the
representatives of the monoid elements (the terms) are given -- prefix
reduction for arbitrary monoid rings and commutative reduction for
Abelian monoid rings.
Saturation concepts with respect to these reductions are provided and
for saturated sets now a finitary characterization of the respective
Gr\"obner bases in terms of special
s-polynomials is possible.
The characterization of prefix Gr\"obner bases can be used to give an
enumerating procedure for such a basis which terminates in case a
finite prefix Gr\"obner basis exists.
For Abelian monoid rings a terminating procedure to complete a finite set of
polynomials is provided.
Furthermore, we introduce interreduction to both approaches and the
existence of unique monic reduced Gr\"obner bases with respect to the respective
ordering on the monoid is shown.

In chapter \ref{chapter.grouprings} we show that the subgroup problem
for a group is equivalent to a restricted version of the right ideal membership
problem in the corresponding group ring.
Hence, only groups with solvable subgroup problem can be expected to
allow the construction of finite Gr\"obner bases.
We apply the concept of prefix reduction to the classes of free, plain
respectively context-free groups and give algorithms to compute finite
reduced prefix Gr\"obner bases for finitely generated right ideals.
Furthermore, a generalization of commutative reduction to nilpotent
groups, namely quasi-commutative reduction, is given and a finitary
characterization of Gr\"obner bases in this setting is provided.
A procedure to compute finite Gr\"obner bases with respect to this
reduction for finitely generated right ideals is given.

Chapter \ref{chapter.reductionrings} gives a sketch how the ideas of
chapter \ref{chapter.reduction} can be carried over to monoid rings
over reduction rings.

We close this introduction by giving a brief overview on the main results of this
thesis:
String rewriting and polynomial rewriting have been combined to
transfer Buchberger's ideas to monoid rings. 
Different definitions of reduction have been studied and Gr\"obner
bases have been defined in the respective settings.
The existence of finite Gr\"obner bases for finitely generated right
ideals has been ensured and procedures for finding them have been
given for the following classes: the class of finite monoids,
 the class of free monoids,
 the class of Abelian\footnote{Here we can 
     successfully treat
     the case of ideals.} monoids,
 the class of finite groups,
  the class of free groups,
 the class of plain groups, 
 the class of context-free groups, and
 the class of nilpotent groups\footnotemark[4].

%% file: algebra.tex
\chapter{Basic Definitions}\label{chapter.definitions}
\spruch{6.3}{11.5}{Content, if hence th' unlearn'd their wants may view\\
The learn'd reflect on what before they knew}{Pope}

The main task of this thesis is to combine string rewriting and polynomial
rewriting in the field of monoid rings.
In this chapter we hence embark on  the basic notions and ideas combined in
our approach to generalize Buchberger's ideas to the more general
setting of monoid rings.

{\bf Section 2.1:} Some of the important algebraic systems --
  monoids, groups, rings and fields -- are introduced and studied.
  Further the main objects of this thesis, ideals  in monoid and
  group rings are specified.

{\bf Section 2.2:} Rewriting is a technique that can be used as
  means of presenting structures as well as for reasoning in
  structures. This will be the foundation for both, representing our
  monoids and the right ideal congruence in our rings later on.
 Hence we introduce an abstract concept of rewriting including 
  normalforms, confluence, termination, completion and the ideas
 related to these terms.

{\bf Section 2.3:} Buchberger introduced the ideas of rewriting to
   commutative polynomial rings. We sketch how a polynomial can be
   used as a rule and how reduction using a set of polynomials
   describes the congruence of the ideal generated by the polynomials.
   Critical pairs of Buchberger's reduction can be localized to a
   special overlap of the head terms of the two polynomials involved,
   namely the least common multiple of the terms. This leads to
   the definition of s-polynomials and Gr\"obner bases can be
   characterized as sets of polynomials where all s-polynomials
   related to pairs of the polynomials in the set can be reduced to
   zero.
   This characterization provides a test whether a finite set is a
   Gr\"obner basis which can be used to give a terminating completion
   procedure called Buchberger's algorithm.

{\bf Section 2.4:} A means of presenting monoids and groups are
  semi-Thue systems and semi-Thue systems modulo commutativity. 
  We give definitions of such systems and distinguish finite complete systems
  for presenting special classes of monoids and groups, namely the
  classes of finite,
  free or Abelian monoids and the classes of finite, free, plain, context-free,
  Abelian or nilpotent groups.
  In using complete presentations, the elements of the respective structures have unique
  representatives and computation in the monoid is possible.
  The equivalence between the irreducible elements of the presenting
  system and the presented structure will be even more important for the
  definition of the ordering of the monoid, which will be assumed to
  be induced by the completion ordering of the presentation.
  This leads to syntactical
   weakenings of reduction in special monoid rings that will
   be studied in chapter \ref{chapter.reduction}.

For more information on algebra and group theory the reader is e.g.
 referred to the books of Herstein \cite{He64} or  Kargapolov and 
 Merzljakov \cite{KaMe79}.
A detailed description on the subject of Gr\"obner bases can be found
for example 
in the books of Becker and Weispfenning \cite{BeWe92} and 
 Geddes,  Czapor and  Labahn \cite{GeCzLa92}.
A good source on rewriting  and monoid presentations is the work of
Book and Otto in \cite{BoOt93}.
%
%
%
\section{Algebra}\label{section.algebra}
Mathematical theories are closely related with the study of two objects, namely
 sets and functions.
Algebra can be regarded as the study of algebraic operations on sets, i.e., 
 functions that take elements from a set to the set itself.
Certain algebraic operations on sets combined with certain axioms are again 
 the objects of independent theories.
This chapter is a short introduction to some of the algebraic systems
  used later on:
 monoids, groups, rings and fields.
\begin{definition}~\\
{\rm
A non-empty set of elements $\m$ together with a binary operation $\mm_{\m}$ is
 said to form a \index{monoid}\betonen{monoid}, if for all $a,b,c$ in
 $\m$      
  \begin{enumerate}
\item $\m$ is closed under $\mm_{\m}$, i.e.,  $a \mm_{\m} b \in \m$,
\item the associative law
  holds for $\mm_{\m}$, i.e., 
  $a \mm_{\m} ( b \mm_{\m} c) =_{\m} (a \mm_{\m} b) \mm_{\m} c$, and
\item there exists  $\lambda_{\m} \in \m$ such
  that $a \mm_{\m} \lambda_{\m} =_{\m} \lambda_{\m} \mm_{\m} a =_{\m} a$.
  The  element $\lambda_{\m}$ is called \index{identity}\betonen{identity}.
\dend
\end{enumerate}
}
\end{definition}
For simplicity of notation we will henceforth drop the index $\m$ and 
 write $\mm$ respectively $=$ if no confusion is likely to arise.
Furthermore, we will often
 talk about a monoid without mentioning its binary operation
 explicitly.
The monoid operation will often be called multiplication or  addition.
Since the algebraic operation is associative we can omit brackets, hence
 the product $a_1 \mm \ldots \mm a_n$ is uniquely defined.
The product of $n \in \n$\footnote{In the following $\n$ denotes the
  set of natural numbers including zero and $\n^+ = \n \backslash \{ 0 \}$.} times the same element $a$ is called the
 \index{n-th power of an element}\betonen{n-th power of $a$} and
 will be denoted by $a^n$,
 where $a^0 = \lambda$.
\begin{definition}~\\
{\rm
An element $a$ of a monoid $\m$ is said to have 
 \index{infinite order}\index{order!infinite}\betonen{infinite
  order}\/
 in case for all $n,m \in \n$, $a^n = a^m$ implies $n=m$.
We say that $a$ has \index{finite
  order}\index{order!finite}\betonen{finite order}\/ 
 in case the set $\{ a^n \mid n \in \n^+ \}$
 is finite and the cardinality of this set is then called the order of
 $a$.
\dend
}
\end{definition}
\begin{definition}~\\
{\rm
For a subset $S$ of a monoid $\m$ we call 
\begin{enumerate}
\item $\ideal{r}{\m}(S) = \{ s \mm m \mid s \in S, m \in \m \}$
       the right ideal,
\item $\ideal{l}{\m}(S) =\{ m \mm s \mid s \in S, m \in \m \}$
       the left ideal, and
\item $\ideal{}{\m}(S) =\{ m \mm s \mm m' \mid s \in S, m, m' \in \m \}$
       the ideal
\end{enumerate}
generated by $S$ in $\m$.
\dend
}
\end{definition}
A monoid $\m$ is called 
 \index{monoid!Abelian}\index{Abelian}\index{Abelian!monoid}\index{monoid!commutative}\index{commutative}\index{commutative!monoid}\betonen{commutative (Abelian)}\/
 if  we have $a \mm b = b \mm a$ for all elements $a,b$ in $\m$.
A natural example for a commutative monoid are the integers together
 with multiplication or addition. \index{addition}\index{multiplication}

A mapping $\phi$ from one monoid $\m_1$ to another monoid $\m_2$
 is called a \index{homomorphism}\betonen{homomorphism},
 if $\phi(\lambda_{\m_1})=
 \lambda_{\m_2}$ and for all $a,b$ in $\m_1$, $\phi(a \mm_{\m_1} b) =
 \phi(a) \mm_{\m_2} \phi(b)$.
In case $\phi$ is surjective we call it an
 \index{epimorphism}\betonen{epimorphism}, in case $\phi$
 is injective a \index{monomorphism}\betonen{monomorphism}\/
 and in case it is both an \index{isomorphism}\betonen{isomorphism}.
The fact that two structures $S_1$, $S_2$ are isomorphic will be denoted by
 $S_1 \cong S_2$.

A monoid is called
 \index{left-cancellative}\betonen{left-cancellative}\/ 
 (respectively \index{right-cancellative}\betonen{right-cancellative}\/) if
for all $a,b,c$ in $\m$, $c \mm a = c \mm b$
 (respectively $a \mm c = b \mm c$) implies $a = b$.
In case a monoid is both, left- and right-cancellative,
 it is called \index{cancellative!left-}\index{cancellative!right-}\index{cancellative}\betonen{cancellative}.
In case $a \mm c = b$ we say that $a$ is a \betonen{left divisor}\/
of $b$ (denoted by $a\;{\sf ldiv}\;b$) and $c$ is called a \betonen{right
  divisor}\/ of $b$ (denoted by $c\;{\sf rdiv}\;b$).
If $c \mm a \mm d =b$ then $a$ is called a \betonen{divisor} of $b$ (denoted
by $a\;{\sf div}\;b$).
A special class of monoids fulfill that for all $a,b$ in $\m$
 there exist $c,d$ in $\m$ such
 that $a \mm c = b$ and $d \mm a = b$, i.e., right and left divisors
 always exist.
These structures are called groups and they can be specified by
 extending the definition of monoids and we do so by
 adding one further axiom.
\begin{definition}~\\
{\rm
A monoid $\m$ together with its binary operation $\mm$ is said to
 form a \index{group}\betonen{group}\/ if additionally 
\begin{enumerate}
\item[4.] for every $a \in \m$ there exists an element $\inv{a} \in \m$ (called
  \index{inverse}\betonen{inverse}\/ of $a$) such
  that $a \mm \inv{a} = \inv{a} \mm
  a = \lambda$. 
\dend
\end{enumerate}
}
\end{definition}
Obviously, the integers form a group with respect to addition, 
but this is no longer
 true for multiplication.

For  n-th powers in groups we can set
 $(\inv{a})^n = \inv{a^n}$.
If for a non-identity element $a \in \g$ there exists an element
 $n \in \n$ such that $a^n=1$
 we call the smallest such number the 
 \index{order!of a group element}\index{order}\betonen{order}\/ of the element
 $a$ in the group $\g$.
Otherwise the order of $a$ is said to be infinite.
A group will be called
\index{group!torsion-free}\index{torsion-free}\betonen{torsion-free}\/
 if every non-identity element of $\g$ has
 infinite order.
\begin{definition}~\\
{\rm
A subset ${\cal H}$ of a group $\g$ is called a  
\index{group!subgroup}\index{subgroup}\betonen{subgroup}\/  of $\g$
 if ${\cal H}$ itself forms a group with respect to the binary operation on $\g$.
The operation on ${\cal H}$ is then said to be 
 \index{induced by}\betonen{induced}\/ by
 the operation on $\g$.
We will write ${\cal H} \leq \g$ and in case ${\cal H}$ is a
 proper subset  ${\cal H} < \g$.\phantom{XX}
\dend
}
\end{definition}
Since intersections of subgroups are again
 subgroups,
 we can define a group in terms of generators as follows:
Given an arbitrary subset $S \subseteq \g$ we define $\langle S \rangle$ to be the
 intersection of all subgroups of $\g$ containing $S$.
Then  $\langle S \rangle$ is a subgroup of $\g$ and $S$ is called a
 \index{group!generating set of a }\index{subgroup!generating
 set of a }\index{generating set}\betonen{generating set}\/
 for $\langle S \rangle$.
In case a group can be generated by some finite subset we call it
 \index{group!finitely generated}\index{subgroup!finitely
  generated}\index{finitely generated}\betonen{finitely generated}.
For example, a cyclic group can be finitely generated by one element.
The following theorem gives a constructive description of this
 generating process.
\begin{theorem}~\\
{\sl
If $S$ is a subset of a group $\g$ then
$$ \langle S \rangle = \{ s_1 \mm \ldots \mm s_n \mid n \in \n, s_i \in 
    S \cup \{ \inv{s} | s \in S \} \}.$$
\ohnebeweis
}
\end{theorem}
Given a subgroup ${\cal H}$ of a group $\g$, for each element $g$ in $\g$ we
 can define special subsets of $\g$, the  
 \index{coset}\index{coset!left}\index{left coset}\betonen{left coset}\/
 $g{\cal H} = \{ g \mm h \mid h \in {\cal H} \}$
 and the \index{coset!right}\index{right!coset}\betonen{right coset}\/
  ${\cal H}g = \{ h \mm g \mid h \in {\cal H} \}$.
Note that for two elements $g_1,g_2$ in $\g$ we have
 $g_1{\cal H} = g_2{\cal H}$ if and only if $\inv{g_1} \mm g_2 \in {\cal H}$
 respectively  ${\cal H}g_1 = {\cal H}g_2$ if and only if
 $g_1 \mm \inv{g_2} \in {\cal H}$.
Thus a subgroup ${\cal H}$ defines a left respectively right 
 congruence by setting $g_1 \sim^{l}_{\cal H} g_2$ if and only if $\inv{g_1}
 \mm g_2 \in {\cal H}$, respectively $g_1 \sim^{r}_{\cal H} g_2$ if
 and only if $g_1 \mm \inv{g_2} \in {\cal H}$.
It can be shown that the sets 
 of all left respectively right cosets are isomorphic.
Their cardinality is called the  \index{index}\betonen{index}\/
 of ${\cal H}$ in $\g$,
 denoted by $|\g : {\cal H}|$.
%

In group theory a particularly important role is played by a special
 kind of subgroup.
\begin{definition}~\\
{\rm
A subgroup ${\cal N}$ of a group $\g$ is called 
 \index{subgroup!normal}\index{normal subgroup}\betonen{normal}\/ if for
 each $g$ in $\g$, we have $g{\cal N} = {\cal N}g$.
We denote this by ${\cal N} \trianglelefteq\g$ or ${\cal N} \triangleleft \g$ 
 in case ${\cal N}$ is a proper subgroup.
\dend
}
\end{definition}
We call two elements $g,h$ of a group
 \index{conjugate}\betonen{conjugate}
 if there exists
 an element $a$ in the same group
 such that $g = \inv{a} \mm h \mm a$ which is sometimes 
 abbreviated by $g = h^{a}$.
Obviously,  ${\cal N} \trianglelefteq \g$ implies that for every
 $g \in \g$ and $h \in {\cal N}$, $\inv{g} \gm h \gm g \in {\cal N}$.
Of course every group contains $\{ \lambda \}$ and itself
 as trivial normal subgroups.
A group containing only these normal subgroups is called
 \index{group!simple}\betonen{simple}.
For example finite cyclic groups of prime order are simple.

The importance of normal subgroups  in group theory stems from the
 fact that they can be used in various ways to built new groups.
\pagebreak
\begin{theorem}\label{theo.factorgroup}~\\
{\sl
The  \index{group!quotient}\index{quotient!of
 groups}\index{group!factor}\index{factor!group}\betonen{quotient}\/ or \betonen{factor group}\/
 of two groups $\g$ and ${\cal N}$, where ${\cal N}$ is a normal subgroup
 of $\g$, is defined by  
 $$\g / {\cal N} = \{ g{\cal N} \mid g \in \g \}$$
 with multiplication
 $$g {\cal N} \mm_{\g / {\cal N}} h {\cal N} = (g \mm_{\g} h) {\cal
   N}.$$
\ohnebeweis
}
\end{theorem}
Moreover, this theorem is closely related to a group construction
method as $\g$ can be seen as an extension of its normal subgroup
${\cal N}$.
\begin{definition}~\\
{\rm
A group $\g$ is said to be an  
 \index{extension}\index{group!extension}\betonen{extension}\/
 of a group  ${\cal N}$ by a group ${\cal H}$ if
\begin{enumerate}
\item ${\cal N} \triangleleft \g$ and
\item $\g / {\cal N} \cong {\cal H}$.
\dend
\end{enumerate}
}
\end{definition}
Thus for two groups ${\cal N} \triangleleft \g$ as described in theorem 
 \ref{theo.factorgroup},
 $\g$ is an extension of ${\cal N}$ by  the quotient $\g / {\cal N}$.
Of course such extensions need not be unique:
 $\z_4$ (which is isomorphic to the structure $\{ 0,1,2,3 \}$ with
 addition modulo 4) as well as the direct product
 $\z_2 \times \z_2$ (which is isomorphic to the structure
 $\{ (0,0),(0,1),(1,0),(1,1) \}$ with component wise addition
 modulo 2) are both extensions of $\z_2$ by $\z_2$ but not isomorphic.
The next definition gives a characterization of groups that are finite
extensions of free groups\footnote{A free group in this
   context is a group where no additional algebraic laws hold except
   the group axioms. A more specific definition will be given later on.}.
\begin{definition}~\\
{\rm
A finitely generated  group $\g$  is called  
 \index{group!context-free}\index{context-free group}\betonen{context-free} if it
 contains a free normal subgroup
 of finite index.
\dend
}
\end{definition}
Another familiar means to construct new groups are products and we
 will give two such constructions that will be of use later on.
\begin{theorem}\label{theo.directproduct}~\\
{\sl
The \betonen{direct}\/ or \index{direct product}\index{product!direct}\index{group!product!direct}\index{Cartesian product}\index{product!Cartesian}\index{group!product!Cartesian}\betonen{Cartesian product}\/
 of two groups $\g$ and ${\cal H}$
 defined by 
 $$\g \times {\cal H} = \{ (g,h) \mid g \in \g, h \in {\cal H} \}$$
 with component wise multiplication
 $$(g,h) \mm_{\g \times {\cal H}} (g',h') = (g \mm_{\g} g',h \mm_{{\cal H}} h' )$$
 is a group.
\ohnebeweis
}
\end{theorem}
It can be shown that the direct product is the unique answer to the
 following question:
Given two groups $\g$ and ${\cal H}$, does there exist a group ${\cal N}$
 and homomorphisms $\phi_1: {\cal N} \myr \g$, $\phi_2: {\cal N} \myr {\cal H}$
 such that for any group ${\cal M}$ and any homomorphisms $\psi_1: {\cal M} \myr \g$,
 $\psi_2: {\cal M} \myr {\cal H}$, there is a unique homomorphism 
 $\theta: {\cal M} \myr {\cal N}$ such that 
 $\theta \circ \phi_1 = \psi_1$ and $\theta \circ \phi_2 = \psi_2$?
This corresponds to the following diagram:

\begin{diagram}[size=2em]
         &                & \g            \\
         &  \ruTo_{\psi_1}  & \uTo_{\phi_1} \\
{\cal M} &  \rTo_{\theta}   & {\cal N}      \\ 
         &  \rdTo_{\psi_2}  & \dTo_{\phi_2} \\
         &                & {\cal H}      \\
\end{diagram}

The answer to the dual question as illustrated in the next diagram is 
 again a unique group called the \betonen{free product}\/ of $\g$ and ${\cal H}$.

\begin{diagram}[size=2em]
\g             & & \\
\dTo_{\phi_1}  & \rdTo_{\psi_1} & \\
{\cal N}       & \rTo_{\theta}  & {\cal M} \\
\uTo_{\phi_2}  & \ruTo_{\psi_2} & \\
{\cal H}       &&\\
\end{diagram}

\auskommentieren{
\begin{theorem}\label{theo.freeproduct}~\\
{\sl 
The \betonen{free product}\/ of two groups $\g_1$ and $\g_2$
 \index{free product}\index{product!free}\index{group!product!free}
 defined by $$\g_1 \ast \g_2 = \langle \g_1,\g_2 \rangle$$
 is a group.
In particular, if $S_1,S_2$ are generating sets of $\g_1,\g_2$ then we get
 $\g_1 \ast \g_2 = \langle S_1,S_2 \rangle$.
\ohnebeweis
}
\end{theorem}}
A very simple type of groups can be built using free products.
\begin{definition}~\\
{\rm 
A \index{plain group}\index{group!plain}\betonen{plain group}\/ is a
 free product of finitely many finite and  free groups.
\dend
}
\end{definition}
One way to learn something about a group is to study subgroups it can be built from.
This can sometimes be done by looking at special chains of subgroups of the group.
Every subgroup ${\cal H}$ of a group $\g$ determines a chain
 $\{ \lambda \} \leq {\cal H} \leq \g$.
A finite chain
 $$ \{ \lambda \} = \g_0 \leq \ldots \leq \g_k = \g$$
 is called a 
 \index{series!normal}\index{normal series}\betonen{normal series of length $k$}\/ for the group $\g$
 if for each $1 \leq i \leq k-1$
 we have $\g_i \trianglelefteq \g$.
In case we only require  $\g_i \trianglelefteq \g_{i+1}$ for all $1 \leq i \leq
 k-1$ such a chain is called a 
 \index{subnormal series}\index{series!subnormal}\betonen{subnormal series}\/ of the group $\g$.
The factor groups $\g_{i+1}/\g_i$ of such a subnormal series
 are called \index{factor}\betonen{factors}.
A subnormal series is called  
\index{polycyclic!series}\index{series!polycyclic}\betonen{polycyclic}\/
 if all factors are cyclic and a group
 possessing a polycyclic subnormal series is called 
 \index{polycyclic!group}\index{group!polycyclic}\betonen{polycyclic}.
A group is called
 \index{polycyclic-by-finite group}\index{group!polycyclic-by-finite}\betonen{polycyclic-by-finite}\/
 if it has a polycyclic normal
 subgroup ${\cal N}$ of finite index.

In order to characterize an interesting subclass of the polycyclic
 groups we need some more definitions.
The following facts are taken from \cite{Ha59} and \cite{KaMe79}.

The \index{group!center}\index{center}\betonen{center}
 of a group $\g$ is defined as 
 ${\cal C}(\g) = \{ g \in \g \mid h \mm g = g \mm h\mbox{ for all } h
 \in \g \}$ and is a characteristic Abelian subgroup of $\g$. 
For two elements $g,h$ of a group their 
 \index{commutator}\betonen{commutator} is defined as
 $[g,h] = \inv{g} \mm \inv{h} \mm g \mm h$.
The \betonen{commutator} of two subsets $S_1,S_2$ of a group $\g$ is denoted by
 $[S_1,S_2] = \{ [g,h] \mid g \in S_1, h \in S_2 \}$.
The special commutator $[\g,\g] = \{ [g,h] \mid g, h \in \g \}$ is a
 subgroup of $\g$ called the 
 \index{derived group}\index{group!derived}\betonen{derived group}\/ of $\g$.
A normal series of a group $\g$ 
 $$ \{ \lambda \} = \g_0 \leq  \ldots \leq \g_k = \g$$
 is called \index{series!central}\index{central
   series}\betonen{central}\/ 
 if all of its factors are \betonen{central}, i.e., if
 $$\g_{i+1} / \g_{i} \leq {\cal C}(\g / \g_{i})
 \mbox{ for all $1 \leq i \leq k-1$}$$
 or equivalently
 $$ [\g_{i+1} , \g] \leq \g_{i}  \mbox{ for all $1 \leq i \leq k-1$}. $$
A group having a (finite) central series is called 
 \index{nilpotent!series}\index{series!nilpotent}\betonen{nilpotent}.
The length of the shortest of all such series for a group is
  called its \index{nilpotency class}\betonen{nilpotency class}.
Note that the Abelian groups are exactly the nilpotent groups of class
 less or equal to 1.
We get the following classification of polycyclic groups by nilpotent groups.
\begin{theorem}~\\
{\sl 
Every polycyclic group
 $\g$ has a normal subgroup ${\cal H}$ of finite index such that ${\cal H}$ is nilpotent-by-Abelian,
 i.e., ${\cal H}$ has a normal subgroup ${\cal N}$ such that ${\cal N}$ is nilpotent and 
 ${\cal H} / {\cal N}$ is Abelian.
\ohnebeweis
}
\end{theorem}
Given a nilpotent group $\g$ we can construct a special central
 series called the
 \index{lower central series}\index{series!lower
   central}\betonen{lower central series of $\g$} by setting:
$$\gamma_1(\g) := \g \mbox{ and } 
\gamma_i(\g) := [\gamma_{i-1}(\g), \g] \mbox{ for } i \geq 2.$$
This gives us a normal series
 $\ldots \leq \gamma_i(\g) \leq \ldots  \leq \gamma_1(\g)=\g$.
If $\gamma_c(\g) \neq \{ \lambda \}$ and 
 $\gamma_{c+1}(\g) = \{ \lambda \}$ then $\g$ is
 nilpotent of class $c$.
A group $\g$ then is
 nilpotent if $\g = \{ \lambda \}$ or
 there is a $c \in \n$ such that $\g$
 is nilpotent of class $c$.

The following theorem shows how nilpotent groups can be regarded
 as polycyclic groups as 
 it establishes the existence of a subnormal series with cyclic factors.
\begin{theorem}~\\
{\sl
A  finitely generated group $\g$ is nilpotent if and only if there exist $\g_0, \ldots, \g_{k} \leq \g$
 and $g_1, \ldots, g_k \in \g$ such that the following statements hold:
\begin{enumerate}
\item $\g_0 =  \{ \lambda \}$, $\g_{k} =\g$ and $\g_{i} \triangleleft \g_{i+1}$
  for $1 \leq i \leq k-1$, 
\item $\g_{i+1} = \langle \{ g_i \} \cup \g_{i} \rangle$ for  $1 \leq
  i \leq k-1$, and
\item $\g_{i+1} / \g_{i} \leq {\cal C}(\g / \g_{i})$  for  $1 \leq i
  \leq k-1$.
\ohnebeweis
\end{enumerate}
}
\end{theorem}
In case we have a finitely generated nilpotent group that is
 additionally torsion-free, even a central series with infinite cyclic
 factors exits.
Moreover, arbitrary finitely generated nilpotent groups are finite
 extensions of torsion-free ones (compare theorem 17.2.2 in \cite{KaMe79}).
\begin{theorem}~\\
{\sl
Every finitely generated nilpotent group has a normal subgroup of
 finite index that is torsion-free.
\ohnebeweis
}
\end{theorem}
This theorem also holds for finitely generated polycyclic groups.

Up to this point we have seen that besides the class of all
 groups there are many 
 interesting subclasses that have become the objects of study
 in various theories:
 the class of finite groups, the class of Abelian groups,
 the class of context-free groups,
 the class of plain groups, the class of nilpotent groups,
 and  the class of polycyclic groups.
We move on now to introduce the concept of free groups\index{free group}.
A group $\free$ belonging to  a class of groups  $\class$ is called
 \index{free group!in the class $\class$}\betonen{free in the class
  $\class$, freely generated by the set
 $X =  \{ x_i \mid i \in I \}$}\/, if for every group $\g \in \class$
 every mapping  $\phi: X \myr \g$ uniquely extends to a homomorphism
 $\Phi: \free \myr \g$.
We will write $\free(X)$ if we want to indicate the generating set.
The cardinality of the index set $|I|$ is called the
 \index{rank}\index{free group!rank}\betonen{rank}\/ of
 $\free$ and the set
 $X$ its \index{basis}\index{free group!basis}\betonen{basis}.
Notice that not all classes of groups contain free groups, e.g. the class of
finite groups.
To give more insight into this concept, we will inspect some classes
 with free groups.
Let us start by showing the existence of free groups in the class of all groups.
Informally, a group $\free$  on a set of generators $X$
 is free in the class of all groups
 if only the trivial relations hold among these generators,
 i.e., only the relations 
 $x \mm \inv{x} = \lambda = \inv{x} \mm x$ hold for $x \in X$.
More precisely we can restate the definition given above for the class
 of all groups as follows.
\begin{definition}~\\
{\rm
A group $\free$ is called a
 \index{free group!in the class of all groups}\betonen{free group in the class of all
  groups generated by a set $X \subseteq \free$}\/ provided the following holds:
If $\phi$ is a mapping from $X$ into an arbitrary group $\g$,
 then there exists a unique
 extension of $\phi$ to a homomorphism from $\free$ into $\g$. 
\dend
}
\end{definition}
We will only sketch how the existence of free groups can be shown. 
For more details the reader is referred to \cite{LySch77,KaMe79}.

Given a set of generators $X$ we add a ``copy'' $X^{-1}$ of $X$
 containing the \index{inverse!formal} \index{formal inverse}\betonen{formal inverses} of
 the generators, i.e., for each $x \in X$ there is an element $x^{-1} \in X^{-1}$.
We can set $\Sigma = X \cup X^{-1}$ and call the elements of this
 set \index{letter}\betonen{letters}.
A \index{word}\betonen{word} $w$ is a finite sequence of letters and
 the  \index{empty word}\index{word!empty}\betonen{empty word} is denoted by
 $\lambda$.
An \index{elementary transformation}\betonen{elementary
  transformation}\/ of a word $w$ consists
 of an application of a relation
 $xx^{-1} = \lambda = x^{-1}x$ for $x \in X$.
Now we can state that two words $w_1,w_2$ are equivalent,
 denoted by $w_1 \sim w_2$, if there
 is a chain of elementary transformations leading from $w_1$ to $w_2$.
Let $\Sigma^{*}$ denote the set of all words generated by
 the letters of  $\Sigma$.
Then one can show that the quotient $\free = \Sigma^{*}/\sim$ is in fact a free group in the
 class of all groups.
This factor group will be denoted by $\free(X)$ and multiplication
 is defined by
 $[w_1]_{\sim} \mm_{\free} [w_2]_{\sim} = [w_1w_2]_{\sim}$.
\begin{theorem}~\\
{\sl
All bases for a given free group have the same cardinality which is equal to the rank of the free
 group.
\ohnebeweis
}
\end{theorem}
\begin{theorem}~\\
{\sl
If a group is generated by a set of $n$ elements ($n$ finite or infinite),
 then it is a quotient of a free group of rank $n$.
\ohnebeweis
}
\end{theorem}
Moreover, if the group $\g$ is generated by a set $S=\{ g_i \mid i \in I \}$ and
 $X= \{ x_i \mid i \in I \}$ is an alphabet, then the mapping $\phi : X \myr S$ defined
 by $\phi(x_i) = g_i$ extends to a unique epimorphism
 $\Phi : \free(X) \myr \g$ and
 the elements of the kernel ${\cal R}$ of this epimorphism are called the \index{relator}\betonen{relators} of $\g$ in terms of 
 alphabet $S$.
If a subset $R$ of ${\cal R}$ is such that the smallest normal subgroup
 containing $R$ is ${\cal R}$ itself, then we call $R$ a
 \betonen{set of defining relators} in the alphabet $S$.
This is a way to determine a group  completely (up to isomorphism).
We call the pair $(S, R)$ a
 \index{presentation!group}\index{group!presentation}\betonen{presentation}\/
 of $\g$. 
A group can
 have several presentations. 
In case there exists a presentation where the sets $S$ and $R$ are 
 finite we call the
 group \index{finitely presented}\betonen{finitely presented}.
For example the finite cyclic group of order $3$ can be presented by
 $(\{ a \} , a^3 = \lambda)$ or $(\{ a,b \} , a^2b^{-1}= ab = ba =
 \lambda )$.

Before we continue with the study of free groups in
 other classes of groups, we briefly
 introduce the concept of varieties, as it provides an elegant way of
 describing special classes of groups together with their free groups.
Let $\{ {\rm x}_i,{\rm x}_i^{-1} \mid i \in \n \}$ be a set of variables.
A word $w$ in those letters is called a \index{law}\betonen{law}\/
 in a class $\class$ 
 of groups, if for every group $\g \in \class$, $w$ becomes trivial for all
 assignments $\{ {\rm x}_i \mid i \in \n \} \myr \g$, i.e.,
 if $w$ contains the variables ${\rm x}_1, \ldots , {\rm x}_n$ we get
 $w(g_1, \ldots, g_n) = \lambda$ for all possible assignments of elements
 $g_i \in \g$ to the variables.
Note that laws are just elements of the free group generated by the
 set $\{ {\rm x}_i \mid i \in \n \}$.
Given a set $V$ of laws and an arbitrary group $\g$ (not
 necessarily in $\class$) we can define
$$V(\g) = \langle v_j(g_1, \ldots , g_{n_j}) \mid v_j \in V, g_i \in \g, n_j
\in \n\rangle\footnote{We assume that the word $v_j$ contains $n_j$ variables.}$$
and $V(\g)$  is a subgroup of $\g$ and in some sense $V(\g)$ ``measures'' the deviation of $\g$
 from the laws in $V$.
For example the law ${\rm x}_1^{-1}{\rm x}_2^{-1}{\rm x}_1{\rm x}_2$ derived from the equation
 ${\rm x}_1{\rm x}_2={\rm x}_2{\rm x}_1$ specifies the class of Abelian groups.
The center and the derived subgroup of a group give a ``measure'' for the
 departure from commutativity, i.e., the bigger the center and the
 smaller the derived subgroup, the nearer the group is to being commutative. 

We proceed now to establish the existence of free groups in classes
 of groups that are \index{variety}\betonen{varieties}, i.e.,
 can be defined by a set
 of laws.
\begin{theorem}\label{theo.freegroups}~\\
{\sl
For a variety $\variety$  of groups defined by a set of laws $V$ and
for any alphabet $X = \{ {\rm x}_i \mid i \in I \}$, define
$$\free_V(X) = \free(X) / V(\free(X)).$$
Then these quotient groups are free in the variety $\variety$.
In particular, every free group in $\variety$ is isomorphic to some $\free_V(X)$.
\ohnebeweis
}
\end{theorem}
Obviously, the class of all groups is a trivial variety defined
 by the empty set of laws and the free groups described by this
 theorem coincide with the  free groups $\free$
 as specified before.

Since the class of Abelian groups is a variety defined by
 $V = \{ {\rm x}_1^{-1}{\rm x}_2^{-1}{\rm x}_1{\rm x}_2 \}$
 we can apply theorem \ref{theo.freegroups} to characterize the free
 groups in this class.
Let $S = \{ a_1, \ldots, a_n \}$ be a finite\footnote{The same
  approach is possible for
 infinite generating sets.} set of generators.
Then the set of all ordered group words on $S$
 is a free group in the class of Abelian groups:
 $\free_V(S) \cong \{ a_1^{i_1} \ldots a_n^{i_n} \mid i_j \in \z \}$. 

The class of nilpotent groups of class $\leq s$ can also be defined as a variety.
Remember that $\g$ is nilpotent of class $\leq s$ if and only if
 $\gamma_{s+1}(\g) = \{ \lambda \}$
 in the lower central series of $\g$.
Hence we get the commutator equation
 $[{\rm x}_1, \ldots, {\rm x}_{s+1}] =[ \dots [{\rm x}_1,{\rm x}_2], \ldots, {\rm x}_{s+1}]= \lambda$
 and thus $V = \{ [{\rm x}_1, \ldots, {\rm x}_{s+1}] \}$.
The free groups in the class of nilpotent groups of class $\leq s$ can
be specified by theorem \ref{theo.freegroups} but
 have no such nice characterization as in the Abelian case as the quotient 
 involves the commutator equation for all generators.
Therefore, we only give an example for a free nilpotent group.
\begin{example}\label{exa.fnp2}~\\
{\rm
Let $\g$ be the free nilpotent group of class 2 with the generators $a$ and $b$.
Further we set $c = [a,b]$.
Then every element $g \in \g$ can be uniquely written as $g \id
a^{i_1}b^{i_2}c^{i_3}$ for some $i_1,i_2,i_3 \in \z$ and
multiplication in $\g$ is given by
$$ a^{i_1}b^{i_2}c^{i_3} \mm a^{j_1}b^{j_2}c^{j_3} =
 a^{i_1 + j_1}b^{i_2 + j_2}c^{i_3 + j_3 + i_2 \cd j_1}.$$
\exaend
}
\end{example}

We end this section by briefly introducing some more algebraic structures that will be used throughout.
\begin{definition}~\\
{\rm
A nonempty set ${\cal R}$ is called an 
 \index{ring}\index{addition!in a ring}\index{multiplication!in a ring}\index{ring!addition}\index{ring!multiplication}\index{ring!unit}\betonen{(associative) ring (with
 unit element)}\/ if there are two binary operations $+$ (addition) and
 $\mrm$ (multiplication) such
 that for all $a,b,c$ in ${\cal R}$
\begin{enumerate}
\item ${\cal R}$ together with $+$ is an Abelian group with zero
  element $0$ and inverse $-a$,
\item ${\cal R}$ is closed under $\mrm$, i.e., $a \mrm b \in {\cal R}$, 
\item $\mrm$ is associative, i.e.,  $a \mrm (b \mrm c) = (a \mrm b ) \mrm c$,
\item the distributive laws hold, i.e., $a \mrm (b + c) = a \mrm b + a \mrm c$
       and $(b + c) \mrm a = b \mrm a + c \mrm a$,
\item there is an element $1 \in {\cal R}$ (called
       \index{unit}\betonen{unit}\/) such that $1 \mrm a = a
       \mrm 1 = a$.
\dend
\end{enumerate}
}
\end{definition}

A ring is called 
 \index{ring!commutative}\index{ring!Abelian}\index{commutative!ring}\index{Abelian!ring}\betonen{commutative (Abelian)}\/ if $a \mrm b = b \mrm a$ for
 all $a,b \in {\cal R}$.
The integers together with addition and multiplication are a well-known  example.
A commutative ring ${\cal R}$ is said to contain 
 \index{zero-divisors}\index{ring!zero-divisors}\betonen{zero-divisors}, if
 there exist not necessarily different elements $a, b$ in ${\cal R}$ such that $a \neq 0$ and $b \neq 0$,
 but $ a \mrm b =0$.
\begin{definition}~\\
{\rm
A commutative ring is called a \index{field}\betonen{field}\/ if
 its non-zero elements form a group under multiplication.
\dend
}
\end{definition}
Similar to our proceeding in group theory we will now look at
 subsets of a ring $R$.
For a subset $U \subseteq R$ to be a subring of $R$ with the operations
 $+$ and $\mrm$ it is necessary and sufficient that
\begin{enumerate}
\item $U$ is a subgroup of $(R,+)$, i.e., for $a,b \in U$ we have $a-b \in U$, and
\item for all $a,b \in U$ we have $a \mrm b \in U$.
\end{enumerate}

We will now take a closer look at special subrings that play a
 role similar to normal subgroups in group theory.
\begin{definition}~\\
{\rm
A nonempty subset $\mswab{i}$ of a ring ${\cal R}$ is called a
  \index{ideal}\index{ideal!right}\index{ideal!left}\index{right!ideal}\index{ideal!two-sided}\index{two-sided ideal}\betonen{right (left) ideal}\/ of ${\cal R}$, if
\begin{enumerate}
\item for all $a, b \in \mswab{i}$ we have $a-b \in \mswab{i}$, and
\item  for every $a \in \mswab{i}$ and $r \in {\cal R}$, the element $a \mrm
  r$ (respectively $r \mrm a$)
  lies in $\mswab{i}$.
\end{enumerate}
A subset that is both, a right and a left ideal, is called a
 \betonen{(two-sided) ideal}\/ of ${\cal R}$.
\dend
}
\end{definition}
For each ring  $\{ 0 \}$ and $R$ are trivial ideals.
Similar to subgroups, ideals can be described in terms of generating sets. 
\index{ideal!trivial}\index{trivial ideal}
\begin{lemma}\label{lem.idealconstruction}~\\
{\sl
Let $A$ be a non-empty subset of ${\cal R}$. Then
\begin{enumerate}
\item  $\ideal{}{}(A) = \{\sum_{i=1}^{n} r_i \mrm a_i \mrm s_i
         \mid a_i \in A, r_i , s_i \in {\cal R}, n \in \n \}$ is an ideal of ${\cal R}$,
\item  $\ideal{r}{}(A) = \{\sum_{i=1}^{n} a_i \mrm r_i
        \mid a_i \in A, r_i \in {\cal R}, n \in \n \}$ is a right ideal
        of ${\cal R}$, and 
\item  $\ideal{l}{}(A) = \{\sum_{i=1}^{n} r_i \mrm a_i
        \mid a_i \in A, r_i \in {\cal R}, n \in \n \}$ is a left ideal
        of ${\cal R}$.
\ohnebeweis
\end{enumerate}
}
\end{lemma}
Notice that the empty sum $\sum_{i=1}^0 a_i$ is zero.

We move on now to combine the algebraic structures introduced so far
to define the rings we are interested in, namely monoid and group rings.
\begin{definition}~\\
{\rm
Let $\myk$ be a field with multiplication denoted by $\skm$ and addition
denoted by $+$ and let $\m$ be a monoid. 
Further let  $\myk[\m]$ denote the set of all mappings \mbox{$f : \m$ $\rightarrow \myk$} where
 the sets ${\sf supp}(f)=\{ m \in \m \mid f(m) \neq 0 \}$ are finite
 and let  $f, g \in \myk[\m]$.
Then the \index{sum!monoid ring}\index{monoid ring!sum}\betonen{sum}
 of $f$ and $g$ is denoted by $f +_{\myk[\m]} 
 g$,  where $(f +_{\myk[\m]} g)(m) = f(m) + g(m)$ and 
 the \index{product!monoid ring}\index{monoid ring!product}\betonen{product}
   is denoted by $f  
 \mrm_{\myk[\m]}  g $, where $(f \mrm_{\myk[\m]} g)(m) = \sum_{x \mm y = m 
 \in \m} f(x) \cd g(y)$.
Notice that $\m$ can be ``embedded'' into $\myk[\m]$ by assigning to every $m \in \m$ 
 a characteristic function $\chi_m: \m \myr \myk$ with
 $\chi_m(m) = 1$ and $\chi_m(m') = 0$ for $m' \in \m\backslash \{ m \}$.
\dend
}
\end{definition}
We will henceforth drop the suffixes of $+$ and $\mrm$ if no confusion
is likely to arise.
Abbreviating $f(m)$ by $ \alpha_{m} \in \myk$ we can express $f$ by the more convenient concept of
 ``polynomials'', i.e., \mbox{$f = \sum_{m \in \m} \alpha_{m} \skm
   m$}.
Notice that since ${\sf supp}(f)$ is finite this sum again is finite. 
This notation now gives us a shorthand for addition and multiplication
in $\myk[\m]$, namely
for  $f = \sum_{m \in \m} \alpha_{m} \skm m$ and $g = \sum_{m \in \m}
 \beta_{m} \skm m$, we get
  $f + g = \sum_{m \in \m} (\alpha_{m} + \beta_{m}) 
 \skm m$ and
 $f \mrm g = \sum_{m \in \m} \gamma_{m} \skm m$ with $\gamma_{m} = \sum_{x 
 \mm y = m  \in \m} \alpha_{x} \skm \beta_{y}$.
$\myk[\m]$ is indeed a ring\footnote{All
  operations mainly involve operations on the
 coefficients in the field $\myk$.} and 
 we call $\myk[\m]$ the \index{monoid ring}\betonen{monoid ring}\/
 of $\m$ over $\myk$, or in case $\m$ is a group the 
 \index{ring!group ring}\index{group!ring}\betonen{group ring}\/ or 
 \index{group!algebra}\betonen{group algebra}\/ of $\m$ over $\myk$.
\begin{remark}\label{rem.zero-divisors}~\\
{\rm
If $\m$ is not cancellative $\myk[\m]$ may contain zero-divisors.
To see this 
let $\Sigma = \{ a,b,c \}$ be the generators of a monoid $\m$ together
with the relations $ab = c$ and $ac = b$.
Then the elements $a^2 -1, c \in {\bf Q}[\m]$ are non-zero, but $(a^2 - 1) \mrm c =
c-c = 0$.
\remend
}
\end{remark}
\begin{example}~\\
{\rm
The polynomial ring over a field $\myk$ in the variables $X_1, \ldots,
 X_n$ is a well-known ring which is also a monoid ring,
 namely the monoid ring
 of the free commutative monoid generated by $X_1, \ldots,X_n$ over $\myk$.
\exaend
}
\end{example}
Since our main interest will be in ideals we give a short
 description of their structure in 
 $\myk[\m]$ in terms of generating sets as described in lemma \ref{lem.idealconstruction}.
For  a non-empty set of polynomials $F$  in $\myk[\m]$ we get 
\index{ideal!in a monoid ring}\index{monoid ring!ideal}
\begin{tabbing}
$\ideal{}{\myk[\m]}(F)\;$ \= $=$ \=
 $\{ \sum_{i=1}^n \alpha_i \skm m_i \mrm p_i \mrm m'_i \mid n \in\n,
 \alpha_i \in \myk, p_i \in F, m_i, m'_i \in \m \}$,\\
$\ideal{r}{\myk[\m]}(F)$ \> $=$ \> 
 $\{ \sum_{i=1}^n \alpha_i \skm p_i \mrm m_i \mid n \in\n,
 \alpha_i \in \myk, p_i \in F, m_i \in \m \}$, and \\
$\ideal{l}{\myk[\m]}(F)$ \> $=$ \>  
 $\{ \sum_{i=1}^n \alpha_i \skm m_i \mrm p_i \mid n \in\n,
 \alpha_i \in \myk, p_i \in F, m_i \in \m \}$.\\
\end{tabbing}
We will simply write $\ideal{}{}(F)$, $\ideal{r}{}(F)$ and
$\ideal{l}{}(F)$ if the context is clear.
Many algebraic problems for rings are related to ideals and we will
close this section by stating two of them\footnote{For more information on
such problems in the special case of commutative polynomial rings see \cite{Bu87}.}.

\begin{tabbing}
XXXXXXXX\=XXXXX \kill
{\bf The Membership Problem\index{membership problem} } \\
\\
{\bf Given:} \> A polynomial $g \in \myk[\m]$ and a set of polynomials 
                $F \subseteq \myk[\m]$. \\
{\bf Question:} \> Is $g$ in the ideal generated by $F$?
\end{tabbing}

\begin{definition}~\\
{\rm
Two elements $f,g \in \myk[\m]$ are said to be
 \index{congruence}\index{ideal!congruence}\betonen{congruent modulo} $\ideal{}{}(F)$,
 denoted by $f \equiv_{\ideal{}{}(F)} g$, if $f = g + h$ for some $h \in
 \ideal{}{}(F)$, i.e., $f-g \in \ideal{}{}(F)$.
\dend
}
\end{definition}
\begin{tabbing}
XXXXXXXX\=XXXXX \kill
{\bf The Congruence Problem\index{congruence problem} } \\
\\
{\bf Given:} \> Two polynomials $f,g \in \myk[\m]$ and a set of polynomials 
                $F \subseteq \myk[\m]$. \\
{\bf Question:} \> Are $f$ and $g$ congruent modulo the ideal generated by $F$?
\end{tabbing}

Note that both problems can similarly be specified for left and right ideals.

%% file: rewriting.tex
\section{The Notion of Reduction}\label{section.rewriting}
%
This section summarizes some important notations and definitions of
reduction relations and basic properties related to them, as can be 
found more explicitly for example in the work of Huet or Book and Otto
(\cite{Hu80,Hu81,BoOt93}).

Let ${\cal E}$ be a set of elements and $\myr$ a binary relation  on
${\cal E}$ called 
 \index{reduction!relation}\index{reduction}\betonen{reduction}.
For $a,b \in {\cal E}$ we will write $a \red{}{\myr}{}{} b$ in case
$(a,b) \in\;\; \myr$.
A pair $({\cal E},\myr)$ will be called a \index{reduction system}\betonen{reduction system}.
Then we can expand the binary relation as follows:
\begin{tabbing}
\hspace{5mm}\=\hspace{1cm}\=\hspace{8mm}\=\hspace{2.5cm}\=\hspace{6cm}\= \kill
\> $\red{0}{\myr}{}{}$\> \>  \> denotes the identity on ${\cal E}$, \\
\> $\red{}{\longleftarrow}{}{}$\phantom{$\red{+}{\lr}{}{}$} \> \> \>denotes the inverse relation for $\myr$, \\
\> $\red{n+1}{\myr}{}{}$\> $:=$ \> $\red{n}{\myr}{}{} \circ \myr$ \> where
                $\circ$ denotes composition of relations and $n \in
                \n$, \\
\> $\red{\leq n}{\myr}{}{}$\> $:=$ \> $\;\!\bigcup_{0 \leq i \leq n}
\red{i}{\myr}{}{}$, \\
\> $\red{+}{\myr}{}{}$ \> $:=$ \> $\;\!\bigcup_{n>0} \red{n}{\myr}{}{}$ \>  denotes the transitive closure
                                               of $\myr$, \\
\> $\red{*}{\myr}{}{}$ \> $:=$ \> $\red{+}{\myr}{}{} \cup \red{0}{\myr}{}{}$ \>  denotes 
   the reflexive transitive closure of $\myr$, \\
\> $\red{}{\lr}{}{}$ \> $:=$ \> $\;\!\longleftarrow \cup \myr$\phantom{$\red{+}{\lr}{}{}$} \>  denotes the
   symmetric closure of $\myr$, \\
\> $\red{+}{\lr}{}{}$ \> \>\>  denotes the symmetric transitive closure
   of $\myr$, \\
\> $\red{*}{\lr}{}{}$\phantom{$\red{+}{\lr}{}{}$} \>\>\>  denotes the reflexive symmetric transitive closure
   of $\myr$.   
\end{tabbing}
%
A well-known decision problem related to a reduction system is the
word problem.
\begin{definition}~\\
{\rm
The \index{word problem}
 \betonen{word problem} for $({\cal E}, \myr)$ is to decide for 
 $a,b$ in ${\cal E}$, whether $a
 \red{*}{\lr}{}{} b$ holds.
\dend
}
\end{definition}
Instances of this problem are well-known in the literature and undecidable in
general.
In the following we will outline sufficient conditions such that $({\cal E}, \myr)$ has
solvable word problem.

An element $a \in {\cal E}$ is said to be
 \index{reducible}\betonen{reducible}\/
 (with respect to $\myr$) if there
 exists an element $b \in {\cal E}$ such that $a \myr b$.
All elements $b \in {\cal E}$ such that $a \red{*}{\myr}{}{} b$ are
 called \index{successor}\betonen{successors}\/
 of $a$ and in case  $a \red{+}{\myr}{}{} b$
 they are called 
 \index{proper!successor}\index{successor!proper}\betonen{proper successors}.
An element which has no proper successors  is called
 \index{irreducible}\betonen{irreducible}.
In case $a \red{*}{\myr}{}{} b$ and $b$ is irreducible, $b$ is called a
 \index{normal form}\betonen{normal form}\/ of $a$.
Notice that for an element $a$ in ${\cal E}$ there can be no, one or many normal
forms.
\begin{definition}~\\
{\rm
A reduction system $({\cal E}, \myr)$ is said to be 
 \index{Noetherian}\index{reduction system!Noetherian}\betonen{Noetherian}\/
 (or \index{reduction system!terminating}\index{terminating}
 \betonen{terminating}\/) in case
 there are no infinitely
 descending reduction chains $a_0 \myr a_1
 \myr \ldots\; $, with  $a_i \in {\cal E}$, $i \in \n$.
\dend
}
\end{definition}
In case $({\cal E}, \myr)$ is Noetherian every element in ${\cal E}$ has
at least one normal form.
\begin{definition}~\\
{\rm
A reduction system $({\cal E}, \myr)$ is called 
 \index{reduction system!confluent}\index{confluence}\betonen{confluent},
 if for all $a, a_1, a_2
 \in {\cal E}$, $a \red{*}{\myr}{}{} a_1$ and $a \red{*}{\myr}{}{} a_2$
implies the existence of $a_3 \in {\cal E}$ such that $a_1
\red{*}{\myr}{}{} a_3$ and $a_2 \red{*}{\myr}{}{} a_3$, which will be
abbreviated by $a_1 \downarrow a_2$ and $a_1$, $a_2$ are called 
\index{joinable}\betonen{joinable}.
\dend
}
\end{definition}
In case $({\cal E}, \myr)$ is confluent every element has at most one
normal form.
We can combine these two properties to give sufficient conditions for
the solvability of the word problem.
\begin{definition}~\\
{\rm
A reduction system $({\cal E}, \myr)$ is said to be 
 \index{reduction system!complete}\index{complete}\betonen{complete}\/
 (or 
  \index{reduction system!convergence}\index{convergent}\betonen{convergent}\/) in case it
 is both, Noetherian and confluent.
\dend
}
\end{definition}
Convergent reduction systems with effective\footnote{By ``effective'' we mean that 
 given an element we can always construct a successor in case one exists.} 
 reduction relations have solvable word problem, as every
 element has a unique normal form and two elements are equal if and
 only if their normal forms are equal.
Of course we cannot always expect $({\cal E}, \myr)$ to be convergent.
Even worse, both properties are undecidable in general.
Nevertheless, there are weaker conditions which guarantee
convergence.
\begin{definition}~\\
{\rm
A reduction system $({\cal E}, \myr)$ is said to be 
 \index{local confluence}\index{confluence!local}\betonen{locally confluent}, 
 if for all $a, a_1, a_2
 \in {\cal E}$, $a \red{}{\myr}{}{} a_1$ and $a \red{}{\myr}{}{} a_2$
 implies the existence of an element $a_3 \in {\cal E}$ such that $a_1
 \red{*}{\myr}{}{} a_3$ and $a_2 \red{*}{\myr}{}{} a_3$.
\dend
} 
\end{definition}
Before stating Newman's lemma which gives a connection between
confluence and local confluence, we introduce the notion of Noetherian
induction that we will use in the proof of the lemma.
\begin{definition}~\\
{\rm
Let $({\cal E}, \myr)$ be a reduction system.
A predicate ${\cal P}$ on ${\cal E}$ is called \betonen{$\myr$-complete}, in case for
 every $a \in {\cal E}$ the following implication holds:  if ${\cal
   P}(b)$ is true for all proper successors of $a$,
 then ${\cal P}(a)$ is true.
\dend
}
\end{definition}
{\sl
{\bf The Principle of Noetherian Induction}:\index{Noetherian!induction} \\
In case  $({\cal E}, \myr)$ is a Noetherian reduction system and ${\cal P}$ is
a predicate that is $\myr$-complete, then  for all $a \in {\cal E}$,
${\cal P}(a)$ is true.
}
\index{Newman's Lemma}
\begin{lemma}[Newman]~\\
{\sl
Let $({\cal E}, \myr)$ be a Noetherian reduction system.
Then $({\cal E}, \myr)$ is confluent if and only if $({\cal E}, \myr)$ is
locally confluent.
\lemend
}
\end{lemma}
\Ba{}~\\
Suppose, first, that  the reduction system $({\cal E}, \myr)$ is confluent.
This immediately implies the local
 confluence of  $({\cal E}, \myr)$ as a special case.
To show the converse, since  $({\cal E}, \myr)$ is Noetherian we can apply the principle of
Noetherian induction to the following predicate:
\begin{center}
${\cal P}(a)$ if and only if for all $a_1,a_2 \in {\cal E}$, $a
\red{*}{\myr}{}{} a_1$ and $a \red{*}{\myr}{}{} a_2$ implies $a_1
\downarrow a_2$.
\end{center}
All we have to do now is to show that ${\cal P}$ is $\myr$-complete.
Let $a \in {\cal E}$ and let ${\cal P}(b)$ be true for all proper successors $b$
of $a$.
We have to prove that ${\cal P}(a)$ is true.
Suppose $a \red{*}{\myr}{}{} a_1$ and $a \red{*}{\myr}{}{} a_2$.
In case $a = a_1$ or $a = a_2$ there is nothing to show.
Therefore, let us assume $a \neq a_1$ and $a \neq a_2$, i.e.,
 $a \myr \tilde{a}_1 \red{*}{\myr}{}{} a_1$ and $a \myr \tilde{a}_2
 \red{*}{\myr}{}{} a_2$. 
Then we can deduce the following  figure
\vspace{-5mm}
\begin{diagram}[size=1.6em]
  & & & & a  & & & &  \\
  & & &\ldTo &    &\rdTo & & &  \\
  & &\tilde{a}_1 & &    & &\tilde{a}_2 & &  \\
  &\ldTo^{*} & &\rdTo^{*} &     &\ldTo^{*}& &\rdTo^{*} &  \\
a_1  & & & & b_0 & & & & a_2  \\
  &\rdTo^{*} & &\ldTo^{*} & & & &\ldTo(4,4)^{*}  &  \\
  & & b_1 & & & && & \\
  & & &\rdTo^{*} & & & & &  \\
  & & & & b  \\
\end{diagram}
\auskommentieren{
\begin{center}
\setlength{\unitlength}{0.0125in}%
\begin{picture}(160,171)(240,540)
\thicklines
\put(315,695){\vector(-1,-1){ 28}}
\put(275,655){\vector(-1,-1){ 28}}
\put(355,655){\vector(-1,-1){ 28}}
\put(315,615){\vector(-1,-1){ 28}}
\put(325,695){\vector( 1,-1){ 28}}
\put(285,655){\vector( 1,-1){ 28}}
\put(245,615){\vector( 1,-1){ 28}}
\put(365,655){\vector( 1,-1){ 28}}
\put(285,575){\vector( 1,-1){ 28}}
\put(395,615){\vector(-1,-1){ 68}}
\put(316,700){\makebox(0,0)[lb]{\raisebox{0pt}[0pt][0pt]{\twlrm $a$}}}
\put(276,660){\makebox(0,0)[lb]{\raisebox{0pt}[0pt][0pt]{\twlrm $\tilde{a}_1$}}}
\put(356,660){\makebox(0,0)[lb]{\raisebox{0pt}[0pt][0pt]{\twlrm $\tilde{a}_2$}}}
\put(236,620){\makebox(0,0)[lb]{\raisebox{0pt}[0pt][0pt]{\twlrm $a_1$}}}
\put(316,620){\makebox(0,0)[lb]{\raisebox{0pt}[0pt][0pt]{\twlrm $b_0$}}}
\put(396,620){\makebox(0,0)[lb]{\raisebox{0pt}[0pt][0pt]{\twlrm $a_2$}}}
\put(276,580){\makebox(0,0)[lb]{\raisebox{0pt}[0pt][0pt]{\twlrm $b_1$}}}
\put(316,540){\makebox(0,0)[lb]{\raisebox{0pt}[0pt][0pt]{\twlrm $b$}}}
\put(255,645){\makebox(0,0)[lb]{\raisebox{0pt}[0pt][0pt]{\twlrm $*$}}}
\put(299,645){\makebox(0,0)[lb]{\raisebox{0pt}[0pt][0pt]{\twlrm $*$}}}
\put(335,645){\makebox(0,0)[lb]{\raisebox{0pt}[0pt][0pt]{\twlrm $*$}}}
\put(379,645){\makebox(0,0)[lb]{\raisebox{0pt}[0pt][0pt]{\twlrm $*$}}}
\put(259,605){\makebox(0,0)[lb]{\raisebox{0pt}[0pt][0pt]{\twlrm $*$}}}
\put(295,605){\makebox(0,0)[lb]{\raisebox{0pt}[0pt][0pt]{\twlrm $*$}}}
\put(299,565){\makebox(0,0)[lb]{\raisebox{0pt}[0pt][0pt]{\twlrm $*$}}}
\put(355,585){\makebox(0,0)[lb]{\raisebox{0pt}[0pt][0pt]{\twlrm $*$}}}
\end{picture}
\end{center}}
where $b_0$ exists, as $({\cal E}, \myr)$ is locally confluent and
$b_1$ and $b$  exist by our induction hypothesis since $a_1$, $b_0$ as well
as $a_2$, $b_1$ are proper successors of $a$.
Hence $a_1 \downarrow a_2$, i.e., the reduction system $({\cal E}, \myr)$ is confluent.
\\
\qed
Therefore, if the reduction system is terminating,
 a check for confluence can be reduced to a check for local
 confluence.
It remains to look for conditions ensuring $({\cal E}, \myr)$ to be
Noetherian.
\begin{definition}~\\
{\rm
A binary relation $\succeq$ on  a set $M$
 is said to be a 
 \index{partial ordering}\index{ordering!partial}\betonen{partial ordering},
 if for all $a,b,c$ in $M$:
\begin{enumerate}
\item $\succeq$ is reflexive, i.e., $a \succeq a$,
\item $\succeq$ is transitive, i.e., $a \succeq b$ and $b \succeq c$
  imply $a \succeq c$, and
\item $\succeq$ is anti-symmetrical, i.e., $a \succeq b$ and $b \succeq
  a$ imply $a=b$.
\dend
\end{enumerate}
}
\end{definition}
A partial ordering is called 
 \index{total ordering}\index{ordering!total}\betonen{total},
 if for all $a, b \in M$
 either $a \succeq b$ or $b \succeq a$ holds.
Further a partial ordering $\succeq$ defines a transitive irreflexive
 ordering $\succ$, where $a \succ b$ if and only if $a \succeq b$ and $a
\neq b$, which is often called a \betonen{proper} or \betonen{strict} ordering.
We call a partial ordering $\succeq$ 
 \index{well-founded
   ordering}\index{ordering!well-founded}\betonen{well-founded},
 if the corresponding strict
 ordering $\succ$ allows  no infinite descending chains $a_0 \succ a_1 \succ
 \ldots\;$, with $a_i \in M$, $i \in \n$.
Now we can give a sufficient condition for a reduction system to be terminating. 
\begin{lemma}~\\
{\sl
Let $({\cal E}, \myr)$ be a reduction system and
suppose there exists a  partial ordering  $\succeq$ on ${\cal E}$ which is
well-founded such that $\myr\;\;\subseteq\;\;\succ$.
Then $({\cal E}, \myr)$ is Noetherian.
}
\end{lemma}
\Ba{}~\\
Suppose the reduction system  $({\cal E}, \myr)$ is not Noetherian.
Then there is an infinite sequence $a_0 \myr a_1 \myr \ldots\;$,  $a_i
\in {\cal E}$, $i \in \n$.
As  $\myr \;\subseteq\; \succ$ this sequence gives us an  infinite sequence
$a_0 \succ a_1 \succ \ldots\;$, with $a_i \in {\cal E}$, $i \in \n$ 
 contradicting our assumption that
$\succeq$ is well-founded on ${\cal E}$.
\\
\qed
We will later on see how a  reduction system with a reduction
 relation fulfilling the requirements of this lemma
 can be
 made convergent by introducing ``completion''.
The partial ordering is then called the completion ordering of the
 reduction system.

The ideas introduced in this section will be specified to special reduction
systems, namely semi-Thue systems to present monoids or groups and
polynomial reduction to present the ideal respectively right ideal congruence in rings.
%


%% file: buchberger.tex
\section{Gr\"obner Bases in Polynomial
  Rings}\label{section.buchberger}
%
The main interest in this section is the study of ideals in polynomial
 rings\index{polynomial ring}\index{ring!polynomial ring} over fields. 
Let $\myk [X_1, \ldots , X_n]$ denote a polynomial ring over the
 (ordered) variables $X_1, \ldots , X_n$.
By ${\cal T} = \{ X_1^{i_1} \ldots X_{n\phantom{1}}^{i_n} \mid i_1, \ldots  i_n \in \n\}$ 
 we define the set of 
 \index{polynomial ring!term}\index{term}\betonen{terms}\/ in this structure.
We recall that a subset
 $F$ of $\myk[X_1, \ldots , X_n]$
 generates an ideal\index{polynomial ring!ideal}\index{ideal!in a polynomial ring}
 $\ideal{}{}(F) = \{ \sum_{i=1}^{k} f_i \mrm g_i \mid k
 \in \n, g_i \in \myk[X_1, \ldots, X_n] \}$
 and $F$ is called a basis\index{basis!of an ideal}\index{ideal!basis} of this ideal. 
It was shown by Hilbert that every ideal in $\myk [X_1, \ldots ,X_n]$ in fact
 has a finite basis, but such a generating set need not allow
 algorithmic solutions for the membership or congruence problem
 related to the ideal.
It was Buchberger who developed a special type of basis, namely the
 Gr\"obner basis, which allows algorithmic solutions for several
 algebraic problems concerning ideals.
He introduced reduction to $\myk [X_1, \ldots , X_n]$ by transforming
 polynomials into ``rules'' and gave a terminating procedure to
 ``complete'' an ideal basis taken as a reduction system.
This procedure is called Buchberger's algorithm. \index{Buchberger's algorithm}
We will give a sketch of his approach below.

Let $\preceq$ be a total well-founded ordering on ${\cal T}$, which is
admissible, i.e., $1 \preceq t$, and $s \prec t$ implies $s \mm u \prec
t \mm u$ for all $s,t,u$ in ${\cal T}$.
In this context $\mm$ denotes the multiplication in ${\cal T}$, i.e.,
 $X_1^{i_1} \ldots X_{n\phantom{1}}^{i_n} \mm X_1^{j_1} \ldots X_{n\phantom{1}}^{j_n} =
  X_1^{i_1+j_1} \ldots X_{n\phantom{1}}^{i_n+j_n}$.
With respect to this multiplication we say that a term $s=X_1^{i_1}
\ldots X_{n\phantom{1}}^{i_n}$ divides a term $t=X_1^{j_1} \ldots X_{n\phantom{1}}^{j_n}$, if
for all $1 \leq l \leq n$ we have $i_l \leq j_l$.
The \index{least common multiple}\index{lcm}\betonen{least common
  multiple}\/ $\lcm(s,t)$ of the terms $s$ and $t$  is the term
 $X_1^{\max \{i_1,j_1\}} \ldots X_{n\phantom{1}}^{\max \{i_n,j_n\}}$.
Note that ${\cal T}$ as mentioned before can be interpreted
 as the free commutative monoid
 generated by $X_1, \ldots, X_n$ with the same multiplication $\cm$ as
 defined above and identity $\lambda = X_1^{0} \ldots
 X_{n\phantom{1}}^{0}$.
We proceed to give an example for a total well-founded admissible
ordering on the set of terms ${\cal T}$.
\begin{example}~\\
{\rm
A 
 \index{total degree ordering}\index{ordering!total degree}\betonen{total degree ordering}\/ $\succ$
 on the terms $X_1^{i_1} \ldots X_{n\phantom{1}}^{i_n} \succ
 X_1^{j_1} \ldots X_{n\phantom{1}}^{j_n} \in {\cal T}$ is specified as follows: 
$X_1^{i_1} \ldots X_{n\phantom{1}}^{i_n} \succ  X_1^{j_1} \ldots X_{n\phantom{1}}^{j_n}$
 if and only if
$ \sum_{s=1}^{n} i_s > \sum_{s=1}^{n} j_s$
 or
$\sum_{s=1}^{n} i_s = \sum_{s=1}^{n} j_s$ and there exists $k$
such that  $i_k > j_k$ and  $i_s=j_s, 1 \leq s < k$.
\dend
}
\end{example}
Henceforth, let $\succeq$ denote a total admissible
ordering on ${\cal T}$ which is of course well-founded.
%
\begin{definition}~\\
{\rm
Let $p = \sum_{i=1}^{k} \alpha_i \skm t_i$ be a non-zero polynomial in $\myk[X_1, \ldots , X_n]$ such
that $\alpha_i \in \myk^* = \myk \backslash \{ 0 \}$, $t_i \in {\cal T}$ and $t_1 \succ \ldots \succ t_n$.
Then we let $\hm(p) = \alpha_{1} \cd t_{1}$ denote the
 \index{head monomial}\betonen{head monomial},
 $\hterm(p) = t_{1}$ the \index{head term}\betonen{head term}\/
 and $\hc(p) = \alpha_{1}$ the 
 \index{head coefficient}\betonen{head coefficient}\/ of $p$.
 $\reductum(p) = p - \hm(p)$ stands for the
 \index{reduct}\betonen{reduct}\/ of $p$.
 We call $p$  \betonen{monic}\index{monic} in case $\hc(p) = 1$.
\mbox{ }\dend
}
\end{definition}
Using the notions of this definition we can recursively extend 
 $\succeq$ on ${\cal T}$ to a partial well-founded admissible ordering
 $\geq$ on $\myk[X_1,\ldots, X_n]$.
\begin{definition}~\\
{\rm
Let $p,q$ be two polynomials in $\myk[X_1,\ldots, X_n]$.
Then we say $p$ is \betonen{greater} than $q$ with respect to a total well-founded admissible
ordering $\succeq$, i.e., $p > q$, if
\begin{enumerate}
\item $\hterm(p) \succ \hterm(q)$ or
\item $\hm(p) = \hm(q)$ and $\reductum(p) > \reductum(q)$.
\dend
\end{enumerate}
}
\end{definition}
We can now split a non-zero polynomial $p$ into a \index{rule}\betonen{rule}
 $\hm(p) \myr - \reductum(p)$
       and we have $\hm(p) > - \reductum(p)$.
Therefore, a set of polynomials gives us a binary relation $\myr$ on
 $\myk[X_1, \ldots , X_n]$ which induces a reduction as follows.
\begin{definition}\label{def.buchberger.red}~\\
{\rm 
Let $p, f$ be two polynomials  in $\myk[X_1, \ldots , X_n]$. 
We say $f$ 
 \index{Buchberger's
   reduction}\index{reduction!Buchberger}\betonen{reduces}
 $p$ to $q$ at a monomial $\alpha \skm t$ of $p$
 in one step, denoted by $p \red{}{\myr}{b}{f} q$, if
\begin{enumerate}
\item[(a)] $\hterm(f) \mm u = t$ for some $u \in {\cal
    T}$\footnote{I.e., $\hterm(f)$ divides $t$.}, and
\item[(b)] $q = p - \alpha \skm \hc(f)^{-1} \skm f \mrm u$.
\end{enumerate}
We write $p \red{}{\myr}{b}{f}$ if there is a polynomial $q$ as defined
above and $p$ is then called reducible by $f$. 
Further, we can define $\red{*}{\myr}{b}{}, \red{+}{\myr}{b}{}$, and
 $\red{n}{\myr}{b}{}$ as usual.
Reduction by a set $F \subseteq \myk[X_1, \ldots , X_n]$ is denoted by
 $p \red{}{\myr}{b}{F} q$ and abbreviates $p \red{}{\myr}{b}{f} q$
 for some $f \in F$,
 which is also written as  $p \red{}{\myr}{b}{f \in F} q$.
\dend
}
\end{definition} 
Note that if $f$ \betonen{reduces}
 $p$ to $q$ at a monomial $\alpha \skm t$ then $t$ is no longer among
 the terms of $q$.
We call a set of polynomials $F \subseteq \myk[X_1, \ldots , X_n]$
 \betonen{interreduced}, if no $f \in F$ is reducible by a polynomial
 in $F \backslash \{ f \}$.
Notice that we have $\myr \;\subseteq\;\; >$ and indeed one can show that
 reduction on $\myk[X_1, \ldots, X_n]$ is Noetherian.
Therefore, we can restrict ourselves to ensuring local confluence when
 describing a completion procedure to compute Gr\"obner bases later on.
\begin{definition}\label{def.buchberger.gb}~\\
{\rm
A set $G \subseteq \myk[X_1,\ldots,X_n]$ is said to be a 
 \index{Gr\"obner basis!Buchberger}\index{Gr\"obner basis}\betonen{Gr\"obner basis}, if
\begin{enumerate}
\item $\red{*}{\lr}{b}{G} = \;\;\equiv_{\ideal{}{}(G)}$, and
\item $\red{}{\myr}{b}{G}$ is confluent.
\dend
\end{enumerate}
}
\end{definition}
The first statement expresses that reduction describes the ideal
congruence and the second one ensures the existence of unique normal forms.
If we additionally require a Gr\"obner basis to be interreduced, such
a basis is unique in case we assume that the polynomials  are monic.
The following lemma gives some properties of reduction, which are
essential in giving a constructive description of a Gr\"obner basis.
\begin{lemma}\label{lem.buchberger.confluence}~\\
{\sl
Let $F$ be a set of polynomials  and $p,q,h$ some polynomials in $\myk[X_1, \ldots, X_n]$.
Then the following statements hold:
\begin{enumerate}
\item
Let $p-q \red{}{\myr}{b}{F} h$.
Then there are polynomials  $p',q' \in \myk[X_1, \ldots, X_n]$ such that 
 $p  \red{*}{\myr}{b}{F} p'$, $q  \red{*}{\myr}{b}{F} q'$ and $h=p'-q'$.
\item
Let $0$ be a  normal form of $p-q$ with respect to $F$.
Then there exists a polynomial  $g \in \myk[X_1, \ldots, X_n]$ such that
 $p  \red{*}{\myr}{b}{F} g$ and $q  \red{*}{\myr}{b}{F} g$.
\item
 $p \red{*}{\lr}{}{F} q \mbox{ if and only if } p - q \in \ideal{}{}(F)$.
\item $p \red{*}{\myr}{b}{F} 0$ implies $\alpha \skm p \mrm u \red{*}{\myr}{b}{F}
  0$ for all $\alpha \in \myk$ and $u \in {\cal T}$.
\item $\alpha \skm p \mrm u \red{}{\myr}{b}{p} 0$ for all $\alpha \in \myk^*$ and $u \in {\cal T}$.
\ohnebeweis
\end{enumerate}
}
\end{lemma}
The second statement of this lemma is often called the translation lemma in the
literature.
Statement 3 shows that Buchberger's reduction always captures the ideal congruence.
Statement 4 is connected to the important fact that
reduction steps are preserved under multiplication with 
coefficients.

Of course we cannot expect an arbitrary ideal basis to be complete\footnote{Note
 that we call a set of polynomials complete (confluent,
 etc.) if the reduction induced by these polynomials used as rules is
 complete (confluent, etc.).}.
But Buchberger was able to show that in order to ``complete'' a given
 basis one only has to add finitely many special polynomials which arise from
 critical situations as described in the context of reduction in the previous
 section.
\begin{definition}\label{def.buchberger.spol}~\\
{\rm
The \index{s-polynomial!Buchberger}\betonen{s-polynomial}\/ for two
 non-zero polynomials
 $p,q \in \myk[X_1, \ldots,X_n]$ is defined as
$$ \spol{}(p,q) = \hc(p)^{-1} \skm p \mrm u - \hc(q)^{-1} \skm q \mrm v,$$
 where $\lcm(\hterm(p),\hterm(q))=\hterm(p) \mm u = \hterm(q)\mm v$
 for some $u,v \in {\cal T}$.
\dend
}
\end{definition}
An s-polynomial will be called non-trivial in case it is not zero and
notice that for non-trivial s-polynomials we always have $\hterm(\spol{}(p,q))
\pred \lcm(\hterm(p),\hterm(q))$.

The following theorem now gives a constructive characterization of
Gr\"obner bases.
\begin{theorem}\label{theo.buchberger.completion}~\\
{\sl
For a set of polynomials $F$ in $\myk[X_1, \ldots,X_n]$,
 the following statements are equivalent:
\begin{enumerate}
\item $F$ is a  Gr\"obner basis.
\item For all polynomials $g \in \ideal{}{}(F)$ we have $g \red{*}{\myr}{b}{F} 0$.
\item For all polynomials $f_{k}, f_{l} \in F$  we have 
  $ \spol{}(f_{k}, f_{l}) \red{*}{\myr}{b}{F} 0$.
\end{enumerate}
}
\end{theorem}
\Ba{}~\\
\mbox{$1 \R 2:$ }
Let $F$ be a Gr\"obner basis and $g \in \ideal{}{}(F)$.
Then $g$ is congruent to $0$ modulo the ideal generated by
 $F$, i.e., $g \red{*}{\lr}{b}{F} 0$.
Thus, as $0$  is irreducible
 we get $g  \red{*}{\myr}{b}{F} 0$.

\mbox{$2 \R 1:$ }
We have to show that reduction with respect to $F$ is confluent.
Since our reduction is terminating it is sufficient to show local
confluence.
Thus, suppose there are three different polynomials $g,h_1,h_2$
 such that $g \red{}{\myr}{b}{F} h_1$ and  $g \red{}{\myr}{b}{F} h_2$.
Then we know $h_1 \equiv_{\ideal{}{}(F)} g \equiv_{\ideal{}{}(F)} h_2$ and hence 
 $h_1 - h_2 \in \ideal{}{}(F)$.
Now by lemma \ref{lem.buchberger.confluence} (the translation lemma), 
 $h_1 - h_2  \red{*}{\myr}{b}{F} 0$ implies the
 existence of  a polynomial  $h \in \myk[X_1, \ldots, X_n]$ such that
 $h_1  \red{*}{\myr}{b}{F} h$ and $h_2  \red{*}{\myr}{b}{F} h$.
Hence, $h_1$ and $h_2$ are joinable.

\mbox{$2 \R 3:$ }
By definition \ref{def.buchberger.spol} the s-polynomial for  two
 non-zero polynomials
 $f_k,f_l \in \myk[X_1, \ldots,X_n]$ is defined as
 $$\spol{}(f_k,f_l) = \hc(f_k)^{-1} \skm f_k \mrm u -
                   \hc(f_l)^{-1} \skm f_l \mrm v,$$
 where $\lcm(\hterm(p),\hterm(q))=\hterm(p) \mm u = \hterm(q)\mm v$  
      and, hence, $\spol{}(f_{k}, f_{l}) \in \ideal{}{}(F)$.
Therefore,  $\spol{}(f_{k}, f_{l})  \red{*}{\myr}{b}{F} 0$ follows immediately.

\mbox{$3 \R 2:$ }
We have to show that every 
 $g \in \ideal{}{}(F) \backslash \{ 0 \}$ is $\red{}{\myr}{b}{F}$-reducible
 to zero.
Remember that for
 $h \in \ideal{}{}(F)$, $ h \red{}{\myr}{b}{F} h'$ implies $h' \in \ideal{}{}(F)$.
As  $\red{}{\myr}{b}{F}$ is Noetherian, thus
 it suffices to show that every  $g \in \ideal{}{}(F) \backslash \{ 0 \}$ 
 is $\red{}{\myr}{b}{F}$-reducible.
Let $g = \sum_{j=1}^m \alpha_{j} \skm f_{j} \mrm w_{j}$ be an arbitrary
 representation of $g$ with $\alpha_{j} \in \myk^*$, $f_j \in F$, and  $w_{j} \in {\cal T}$.
Depending on this representation of $g$ and a total well-founded
admissible  ordering 
 $\succeq$ on ${\cal T}$ we define
 $t = \max \{ \hterm(f_{j}) \mm w_{j} \mid j \in \{ 1, \ldots m \}  \}$ and
 $K$ is the number of polynomials $f_j \mrm w_j$ containing $t$ as a term.
Then $t \succeq \hterm(g)$ and
 in case $\hterm(g) = t$ this immediately implies that $g$ is
 $\red{}{\myr}{b}{F}$-reducible. 
Thus we will prove that $g$ has a representation where every
occurring term is less or equal to $\hterm(g)$, i.e., there exists a
representation such that $t = \hterm(g)$\footnote{Such
  representations are often called standard representations in
  literature (compare \cite{BeWe92}).}.
This will be done by induction
 on $(t,K)$, where
 $(t',K')<(t,K)$ if and only if $t' \prec t$ or $(t'=t$ and
 $K'<K)$\footnote{Note that this ordering is well-founded since $\succ$
                  is well-founded on ${\cal T}$ and $K \in\n$.}.
In case $t \succ \hterm(g)$ there are  two polynomials $f_k,f_l$ in the corresponding 
 representation\footnote{Not necessarily $f_l \neq f_k$.}
 such that  $\hterm(f_k) \mm w_k = \hterm(f_l) \mm w_l=t$.
By definition \ref{def.buchberger.spol} we have an s-polynomial
 \mbox{$\spol{}(f_k,f_l) = \hc(f_k)^{-1} \skm  f_k
 \mrm z_k -\hc(f_l)^{-1}\skm f_l \mrm z_l$} such that
 $\hterm(f_k) \mm z_k = \hterm(f_l) \mm z_l = \lcm(\hterm(f_k),\hterm(f_l))$.
Since $\hterm(f_k) \mm w_k = \hterm(f_l) \mm w_l$ there exists an element 
 $z \in {\cal T}$ such that
 $w_k = z_k \mm z$ and $w_l = z_l \mm z$.
We will now change our representation of $g$ by using the additional
information on this s-polynomial in such a way that for the new
representation of $g$ we either have a smaller maximal term or the 
occurrences of the term $t$
are decreased by at least 1.
Let us assume that $\spol{}(f_k,f_l)$ is not trivial\footnote{In case 
               $\spol{}(f_k,f_l) = 0$, just substitute $0$
               for the sum $\sum_{i=1}^n \delta_i \skm h_i \mrm
               v_i$ in the equations below.}.
Then the reduction sequence $\spol{}(f_k,f_l) \red{*}{\myr}{b}{F} 0 $
results in a representation of the form
 $\spol{}(f_k,f_l) =\sum_{i=1}^n \delta_i \skm h_i \mrm v_i$,
 where $\delta_i  \in \myk^*,h_i \in F,v_i \in {\cal T}$.
As the $h_i$ are due to the reduction of the s-polynomial,
 all terms occurring in the sum are bounded by the term $\hterm(\spol{}(f_k,f_l))$.
Moreover, since $\succeq$ is admissible on ${\cal T}$ this 
 implies that all terms of the sum 
 $\sum_{i=1}^n \delta_i \skm h_i \mrm v_i \mrm z$ are bounded by 
 $\hterm(\spol{}(f_k,f_l)) \mm z \pred t$, i.e., they are strictly bounded
 by $t$\footnote{This can also be concluded by statement four of lemma
   \ref{lem.buchberger.confluence} since $\spol{}(f_k,f_l)
   \red{*}{\myr}{b}{F} 0 $ implies $\spol{}(f_k,f_l) \mrm z
   \red{*}{\myr}{b}{F} 0 $ and $\hterm(\spol{}(f_k,f_l) \mrm z) \prec t$.}.
We can now do the following transformations: 
\begin{eqnarray}
 &  & \alpha_{k} \skm f_{k} \mrm w_{k} + \alpha_{l} \skm f_{l} \mrm w_{l}  \nonumber\\  
 &  &  \nonumber\\                                                             
 & = &  \alpha_{k} \skm f_{k} \mrm w_{k} +
         \underbrace{ \alpha'_{l} \skm \beta_k \skm f_{k} \mrm w_{k}
                     - \alpha'_{l} \skm \beta_k \skm f_{k} \mrm w_{k}}_{=\, 0} 
        + \alpha'_{l}\skm \beta_l  \skm f_{l} \mrm w_{l} \nonumber\\ 
 & = & (\alpha_{k} + \alpha'_{l} \skm \beta_k) \skm f_{k} \mrm w_{k} - \alpha'_{l} \skm 
        \underbrace{(\beta_k \skm f_{k} \mrm w_{k}
        -  \beta_l \skm f_{l} \mrm w_{l})}_{=\, \spol{}(f_k,f_l) \mrm z} \nonumber\\
 & = & (\alpha_{k} + \alpha'_{l} \skm \beta_k) \skm f_{k} \mrm w_{k} -
        \alpha'_{l} \skm (\sum_{i=1}^n \delta_{i} \skm h_{i} \mrm (v_{i} \mm z)
        ) \label{s.buchberger}
\end{eqnarray}
where, $\beta_k=\hc(f_k)^{-1}$, $\beta_l=\hc(f_l)^{-1}$, and  $\alpha'_l \skm \beta_l = \alpha_l$.
By substituting (\ref{s.buchberger}) in our representation of $g$
 either $t$ disappears or $K$ is decreased.
\\
\qed
\begin{remark}~\\
{\rm
A closer inspection of the proof of $3 \R 2$ given above reveals a
 concept that will be initial to the proofs of similar theorems for
 specific monoid rings in the following chapter.
The heart of this proof consists in transforming an arbitrary
 representation of an element $g$ belonging to the ideal generated by
 the set $F$ in such a way that we can deduce a top reduction
 sequence for $g$ to zero, i.e., a reduction sequence where the
 reductions only take part at the respective head term.
Such a representation of $g$ then is  a standard representation
 and hence this technique is closely related to 
the concept of standard bases
as given for example in \cite{BeWe92}.
\remend
}
\end{remark}
As a consequence of theorem \ref{theo.buchberger.completion} it is
decidable whether a finite set of polynomials is a Gr\"obner basis.
Moreover, this theorem gives rise to the following completion procedure
 for sets of polynomials called Buchberger's Algorithm.

\procedure{Buchberger's Algorithm}%
{\vspace{-4mm}\begin{tabbing}
XXXXX\=XXXX \kill
\removelastskip
{\bf Given:} \> A finite set of polynomials $F  \subseteq \myk[X_1, \ldots, X_n]$. \\
{\bf Find:} \> $\gb(F)$, a  Gr\"obner basis of $F$.
\end{tabbing}
\vspace{-7mm}
\begin{tabbing}
XX\=XX\= XXX \= XXX \=\kill
$G$ := $F$; \\
$B$ := $\{ (q_{1}, q_{2}) \mid q_{1}, q_{2} \in G, q_{1} \neq q_{2}
\}$; \\
{\bf while} $B \neq \emptyset$ {\bf do} \\
\>      $(q_{1}, q_{2}) := {\rm remove}(B)$; \\
\>     {\rm\kommentar \% Remove an element from
                     the set $B$}\\
\>      $h := {\rm normalform}(\spol{}(q_{1}, q_{2}),\red{}{\myr}{b}{G})$ \\
\>     {\rm\kommentar \% Compute a normal form of
                     $\spol{}(q_{1}, q_{2})$ with respect to
                     Buchberger's reduction} \\
\>      {\bf if} \>$h \neq 0$\\
\> \>     {\bf then}  \>      $B := B \cup \{ (f,h) \mid f \in G \}$; \\
\> \>                 \>      $G := G \cup \{ h \}$; \\
\>      {\bf endif}\\
{\bf endwhile} \\
$\gb (F):= G$
\end{tabbing}}

Termination can be shown by using a slightly different
 characterization of Gr\"obner bases: 
A subset $G$ of $\ideal{}{\myk[X_1, \ldots, X_n]}(F)$ is a Gr\"obner basis of
 $\ideal{}{\myk[X_1, \ldots, X_n]}(F)$ if and only if
 $\hterm(\ideal{}{\myk[X_1, \ldots, X_n]}(F) \backslash \{ 0 \}) = \ideal{}{\cal T}(\hterm(G))$, i.e., the set of the head
 terms of the polynomials in the ideal generated by $F$ coincides with
 the ideal (in ${\cal T}$) generated by the head terms of the
 polynomials in $G$.
Reviewing the algorithm, we find that every polynomial added 
 in the {\bf while} loop has
 the property that its head term cannot be divided by the head terms of
 the polynomials already in $G$.
By Dickson's lemma or Hilbert's basis theorem,
 the head terms of the polynomials in $G$ will at some step form a
 basis for the set of head terms of the polynomials of the ideal generated by
 $F$ which itself is  the ideal in ${\cal T}$
 generated by the head terms of the polynomials in $G$.
From this time on for every new polynomial $h$ computed by the
 algorithm the head term $\hterm(h)$ must lie in this ideal.
Therefore, its head term must be divisible by at least
 one of the head terms of the polynomials in $G$, i.e., $\hterm(h)$ and hence $h$ cannot be
 in normal form with respect to $G$ unless it is zero.


%% file: presentations.tex
\section{Semi-Thue Systems}\label{section.presentations}
%
In this section we introduce the structures we shall use to present
 our monoids and groups, namely 
 \betonen{semi-Thue systems}\/ (also called 
 \index{string-rewriting system}\betonen{string-rewriting systems}\/).
Let us start with some basic definitions.
\begin{definition}~\\
{\rm
Let $\Sigma$ be a finite alphabet.
\begin{enumerate}
\item By $\Sigma^*$ we will denote the set of all \index{word}\betonen{words}\/ over the
  alphabet $\Sigma$ where $\lambda$ presents the
  \index{empty word}\index{word!empty}\betonen{empty word}, i.e.,
  the word of length zero.
  $\id$ will denote the \index{identity}\betonen{identity}\/ on $\Sigma^*$.
\item Let $u,v$ be words in $\Sigma^*$.
  $u$ is said to be a \index{prefix}\betonen{prefix}\/ of $v$,
  if there exists $w \in
  \Sigma^*$ such that $v \id uw$.
  $w$ is then called a \index{suffix}\betonen{suffix}\/ of $v$.
  In case $w \not\id \lambda$ or $u \not\id \lambda$ we will speak of 
  \index{prefix!proper}\index{suffix!proper}\index{proper!prefix}\index{proper!suffix}\betonen{proper prefixes}\/ respectively \betonen{proper suffixes}.
\item The \index{length}\betonen{length}\/
  of a word is the number of letters it contains,
  i.e., $|\lambda| = 0$ and $|wa|=|w|+1$ for all $w \in \Sigma^*$, $a
  \in \Sigma$.
\item We can define a mapping $\conc : \Sigma^* \times \Sigma^* \myr
  \Sigma^*$  by $\conc(u,v) \id uv$ for
  $u,v \in \Sigma^*$ which will be called
   \index{concatenation}\betonen{concatenation}.
  Then $\conc$ is an associative binary operation on $\Sigma^*$ with
  identity $\lambda$.
  Thus $\Sigma^*$ together with $\conc$ and $\lambda$ is a monoid, namely the 
  \index{free monoid}\index{monoid!free}\betonen{free
  monoid generated by $\Sigma$}.
\item For an element $w \in \Sigma^*$ we define 
\[ \ell (w) =  \left\{ \begin{array}{r@{\quad\quad}l}
              \lambda & \mbox{ if } w \id \lambda \\
              a & \mbox{ if } w \id ua,u \in \Sigma^*,a \in \Sigma,
              \end{array} \right. \]
 i.e., $\ell(w)$ is the last letter of $w$ in case $w$ is not
 the empty word.
\dend
\end{enumerate}
}
\end{definition}
\begin{definition}~\\
{\rm 
Let $\Sigma$ be a finite alphabet.
\begin{enumerate}
\item A \index{semi-Thue system}\betonen{semi-Thue system}
  $T$ over  $\Sigma$ is a subset of
  $\Sigma^* \times \Sigma^*$.
  The elements $(l,r)$ of $T$ are called
  \index{rule}\betonen{rules}\/ and will often be
  written as $l \myr r$.
\item The  \index{single-step reduction relation}\index{reduction!single-step}\betonen{single-step reduction relation}\/
  on $\Sigma^*$ induced
  by a semi-Thue system $T$ is defined as follows:
  For any $u,v$ in $\Sigma^*$, $u \myr_T v$ if and only if there exist
  $x,y$ in $\Sigma^*$ and $(l,r)$ in $T$ such that $u \id xly$ and $v \id
  xry$.
  The \index{reduction relation}\betonen{reduction relation}\/ 
  on $\Sigma^*$ induced by $T$ is the
  reflexive transitive closure of $\myr_T$ and is denoted by
  $\red{*}{\myr}{}{T}$.
  The reflexive transitive symmetric closure   is defined as usual and
  denoted by $\red{*}{\lr}{}{T}$.
\dend  
\end{enumerate}
}
\end{definition}
Recalling section \ref{section.rewriting} we find that the pair
 $(\Sigma^*, \myr_T)$ is a reduction system which is specified by
 $\Sigma$ and $T$.
\begin{definition}~\\
{\rm
A semi-Thue system is called \betonen{normalized} or
\betonen{reduced}\index{normalized}\index{reduced} in case the left
hand sides of the rules can only be reduced by the rule itself and the
right hand sides are irreducible.
\mbox{\phantom{XX}}\dend
}
\end{definition}
Notice that one can even assume that a reduced semi-Thue system does
not contain rules of the form $a \myr b$ or $a \myr \lambda$ where $a, b
\in \Sigma$\footnote{Such rules can be removed using Tietze
  transformations which are a means to change presentations without
  changing the monoid presented.}.
\begin{definition}~\\
{\rm
Let $\Sigma$ be an alphabet. 
A mapping $\imath: \Sigma \myr \Sigma$ is called an
 \index{involution}\betonen{involution}\/
 if $\imath(\imath(a)) = a$
 for all $a \in \Sigma$.
A semi-Thue system is called a 
 \index{group!system}\betonen{group system}\/
 if there exists an
 involution $\imath$ such that for all $a \in \Sigma$
 the rules $(\imath(a)a, \lambda)$ and $(a\imath(a), \lambda)$ are
 included in $T$.
\dend
}
\end{definition}
Note that sometimes we will assume that $\Sigma = \Gamma \cup
\Gamma^{-1}$ where $\Gamma^{-1} = \{ a^{-1} \mid a \in \Gamma \}$ contains the  formal inverses
 of $\Gamma$ and $T$ contains the rules corresponding to the
 trivial relations in a group, namely $\{ (aa^{-1}, \lambda) ,
 (a^{-1}a, \lambda) \mid a \in \Gamma \}$.

An equivalence relation on $\Sigma^*$ is said to be a congruence
relation in case it is admissible, i.e., compatible with
concatenation.\index{congruence}
Since this is obviously true for the reduction relation induced by a
semi-Thue system $T$, the reflexive transitive symmetric closure
$\red{*}{\lr}{}{T}$ is a congruence relation on the set $\Sigma^*$,
the \betonen{Thue
    congruence}\index{Thue congruence}.
Two semi-Thue systems on the same alphabet are called
\betonen{equivalent}\index{equivalent} if they generate the same Thue congruence. 
The congruence classes are denoted by $[w]_T = \{ v \in \Sigma^* \mid v
\red{*}{\lr}{}{T} w \}$ and we can set $\factormonoid_T = \{ [w]_T \mid w \in
\Sigma^* \}$.
In fact $\factormonoid_T$ is the factor monoid of the free monoid
$\Sigma^*$ modulo the congruence induced by $T$ as the following lemma
establishes. 
%
\begin{lemma}~\\
{\sl 
Let $(\Sigma, T)$ be a semi-Thue system.
\begin{enumerate}
\item
The set $\factormonoid_T$ together with the binary operation $\circ :
\factormonoid_T \times \factormonoid_T \myr  \factormonoid_T$ defined by $[u]_T
\circ [v]_T = [uv]_T$ and the identity $[\lambda]_T$ is a monoid, called the
\index{factor monoid}\index{monoid!factor}\betonen{\/factor monoid}\/ of $\Sigma^*$ and $\red{*}{\lr}{}{T}$.
\item
In case $T$ is a group system, the set  $\factormonoid_T$  together with
$\circ$, $[\lambda]_T$ and ${\sf inv}$ is a group, where
$\inv{[w]_T} = [{\rm inv}(w)]_T$, and
${\rm inv}(\lambda) = \lambda$, ${\rm inv}(wa) = \imath(a){\rm inv}(w)$
 for all $w \in \Sigma^*$, $a \in \Sigma$.
\lemend\ohnebeweis
\end{enumerate}
}
\end{lemma}
Hence, semi-Thue systems are means for presenting monoids and groups.
The following definitions are closely related to describing monoids
 and groups in terms of generators and defining relations as given in section
 \ref{section.algebra}.
We call a pair $(\Sigma,T)$ a 
 \index{presentation}\index{monoid!presentation}\index{group!presentation}\betonen{presentation}\/ of a monoid (group) $\m$ if
 $\m \cong \factormonoid_T$.
Note that every monoid can be presented by a (even convergent)
 semi-Thue system.
Just let $\Sigma$  be the set of all elements and $T$  the
 multiplication table.
The problem is that this presentation in general is neither finite nor recursive.
We call a monoid (group) $\m$ 
 \index{finitely generated}\betonen{finitely generated}, if $\m$ has a
 presentation $(\Sigma,T)$ such that $\Sigma$ is finite.
$\m$ is said to be \index{finitely presented}\betonen{finitely
  presented},
 if additionally
 $T$ is finite.
In order to do effective computations in our monoid or group 
 we have to be able
 to compute representatives for the congruence classes of the elements.
A very nice solution occurs in case we are able to give a convergent finite
 semi-Thue system as a presentation, since then every congruence
 class has a unique representative and many problems, e.g.\  the word
 problem, are algorithmically solvable.

After distinguishing  some syntactically restricted semi-Thue
 systems we will conclude this section by giving characterizations
 of certain classes
 of groups and monoids as mentioned in section \ref{section.algebra}
 by special presentations they allow.
A survey on groups allowing convergent presentations can be found
 in \cite{MaOt89}
\begin{definition}~\\
{\rm 
Let $T$ be a semi-Thue system on $\Sigma$.
\begin{enumerate}
\item $T$ is said to be \index{length-reducing}
 \betonen{length-reducing}\/ if for all $(l,r) \in T$,
   $|l|>|r|$.
\item $T$ is said to be \index{monadic}\betonen{monadic}\/ if for all $(l,r) \in T$, $r
  \in \Sigma \cup \{ \lambda \}$ and $|l| \geq |r|$.
\item $T$ is said to be \index{two-monadic, 2-monadic}\index{monadic!2-monadic}\betonen{2-monadic}\/ if for all $(l,r) \in T$,  $r  \in \Sigma \cup \{ \lambda \}$ and $2 \geq |l| \geq |r|$.
\dend
\end{enumerate}
}
\end{definition}
We next introduce 
length-lexicographical orderings which can be used for orienting these 
rule systems.
Then obviously the reduction relation induced by them is Noetherian.
\begin{definition}~\\
{\rm
Let $\Sigma$ be an alphabet and $\succ$ a partial ordering on
 $\Sigma$ called precedence.
Further let $u \id a_1 \ldots a_n$ and $v \id b_1 \ldots b_m$ be two
 words in $\Sigma^*$.
\begin{enumerate}
\item We can define a
       \index{lexicographical ordering}\index{ordering!lexicographical}\betonen{lexicographical ordering} based on $\succ$
       by setting 
      $u >^{\rm lex} v$ if and only if there exists an index $k \in \{
      1, \ldots, \min \{ n,m \} \}$ such that $a_i = b_i$ for all
      $1 \leq i < k$ and $a_k \succ b_k$.
\item We can define a
       \index{length-lexicographical ordering}\index{ordering!length-lexicographical}\betonen{length-lexicographical ordering}
       based on $\succ$ and $>^{\rm lex}$ by setting
      $u >^{\rm llex} v$ if and only if $|u| > |v|$ or $( |u|=|v|$ and
      $u >^{\rm lex} v )$.
\dend
\end{enumerate}
}
\end{definition}
Let us continue with a property of orderings which holds e.g. for  pure 
 lexicographical and length-lexicographical orderings.
\begin{definition}~\\
{\rm 
A partial ordering $\preceq$ on $\Sigma^*$ is called
\betonen{admissible}\index{admissible ordering}\index{ordering!admissible} (with respect to $\conc$)
if for all $u,v,x,y$ in $\Sigma^*$ we have
\begin{enumerate}
\item $1 \preceq u$, and
\item $u \prec v$ implies $xuy \prec xvy$.
\dend
\end{enumerate}
}
\end{definition}
In case $\preceq$ is an admissible well-founded total ordering, then for 
a proper subword $u$ of $v$ (and hence a proper divisor in $\Sigma^*$) we
have $u \prec v$.
This is due to the fact that otherwise $u \succ v \id xuy \succ
x^2uy^2 \succ \ldots$ would give us an infinite descending chain,
contradicting that $\preceq$ is supposed to be well-founded on $\Sigma^*$.
Therefore, in case a monoid $\m$ is presented by a semi-Thue system
$(\Sigma,T)$ which is convergent with respect to an admissible
well-founded total ordering $\succeq$ this yields $uv \succeq u \mm_{\m} v$.

By Newman's lemma we know that under the hypothesis that a
 reduction relation is Noetherian, a reduction system is confluent if
 and only if it is locally confluent.
For semi-Thue systems the global property of being locally
 confluent can be localized to enable a confluence test.
We will now sketch how a finite presentation $(\Sigma, T)$ of a monoid can
 be completed in case we have a total admissible well-founded
 ordering $\succeq$ on $\Sigma^*$ such that for all $(l,r) \in T$
 we have $l \succ r$.
This ordering  then will be called a \betonen{completion ordering}\index{completion ordering}
 and the process of transforming $(\Sigma ,T)$ into a (not necessarily
 finite) convergent presentation of the same monoid is called \betonen{completion}\index{completion}.
Important is that in order to check a finite set $T$ for confluence
 we only have to look at a finite set of critical situations:
for $(l_1,r_1),(l_2,r_2) \in T$ the set of \betonen{critical pairs}\index{critical pairs} is
 $\{ \langle xr_1,r_2y \rangle \mid \mbox{ there are } x,y \in \Sigma^*,
 xl_1 \id l_2y, |x|<|l_2| \} \cup \{ \langle r_1, xr_2y \rangle \mid
 \mbox{ there are } x,y \in \Sigma^*,
 l_1 \id xl_2y, |x|<|l_1| \}$.
Now given a finite semi-Thue system $(\Sigma,T)$ with a completion ordering $\succeq$ we can give a completion process as follows:

\procedure{Knuth Bendix Completion}%
{\vspace{-4mm}\begin{tabbing}
XXXXX\=XXXX \kill
\removelastskip
{\bf Given:} \> $(\Sigma,T), \succeq$ as described above. 
\end{tabbing}
\vspace{-7mm}
\begin{tabbing}
XX\=XX\= XXX \= XXXXX \=\kill
$R_0 := T$; \\
$i:= -1$; \\
{\bf repeat} \\
\> $i:=i+1$ ; \\
\> $R_{i+1} := \emptyset$;\\
\> $B:= {\rm critical.pairs}(R_i)$; \\
\>     {\rm\kommentar \% Compute the critical pairs of all $(l_1,r_1),(l_2,r_2) \in R_i$ as described above} \\
\> {\bf while} $B \neq \emptyset$ {\bf do}\\
\>\> $\langle z_1,z_2 \rangle := {\rm remove}(B)$; \\
\>\>     {\rm\kommentar \% Remove an element using a fair strategy} \\
\>\> $z_1' := {\rm normalform}(z_1, \red{}{\myr}{}{R_i})$; \\
\>\> $z_2' := {\rm normalform}(z_2, \red{}{\myr}{}{R_i})$; \\
\>\> {\bf case }  $z_1' \succ z_2'$:  $R_{i+1}:=R_{i+1} \cup \{(z_1',z_2')\}$; \\
\>\> {\bf case }  $z_2' \succ z_1'$: $R_{i+1}:=R_{i+1} \cup \{(z_2',z_1')\}$; \\
\> {\bf endwhile} \\
\> {\bf case} $R_{i+1} \neq \emptyset$: $R_{i+1} := R_i \cup R_{i+1}$; \\
{\bf until} $R_{i+1} = \emptyset$
\end{tabbing}}
 
Since the word problem for semi-Thue systems is unsolvable, this procedure
 in general will not terminate.
Nevertheless, using a fair startegy to remove elements from the set $B$, it
 always enumerates a convergent semi-Thue system presenting the same
 monoid as the input system.

In the following we will study monoids where finite convergent presentations exist.
Our interest will be in subclasses of the class of finitely presented
 monoids and groups only and the assumption of being finitely presented will often
 be included without being mentioned explicitly.

A very simple subclass of finitely presented monoids is the class of finite monoids.
Obviously every finite monoid $\m$ can be presented by its elements and their 
 multiplication table.
This presentation will be denoted by $(\m ,M_{\m})$ and in fact this gives
 us  presentations for finite monoids by 2-monadic convergent semi-Thue systems.
The same is true for finite groups and we even have finite 2-monadic
convergent group presentations.

In case $(\Sigma', T)$ is a semi-Thue system presenting a free monoid $\m$,
then there is a subset $\Sigma$ of $\Sigma'$ such that $\m$ is freely
generated by $\Sigma$, i.e., $\m \cong \Sigma^*$ and $(\Sigma,
\emptyset)$ is also a presentation of $\m$ and obviously a convergent one.
Of similar simplicity is the characterization of free groups in the
 class of finitely generated groups.
Let $\free$ be a free group generated by $X = \{ x_1, \ldots, x_n \}$.
Then the semi-Thue system 
 $(X \cup X^{-1},I)$ where $I = \{ x_ix_i^{-1} \myr \lambda, x_i^{-1}x_i
 \myr \lambda  \mid 1 \leq i \leq n \}$
 is a presentation of $\free$ which is  2-monadic and convergent.

Presentations for the groups in the class of plain groups can easily be
 constructed from finite and free presentations by using how presentations
 of groups can be combined in order to present free products of these groups.
Let $\g_1,\g_2$ be two groups with  presentations
 $(\Sigma_1,T_1), (\Sigma_2,T_2)$ such that
 $\Sigma_1 \cap \Sigma_2 = \emptyset$\footnote{This can always be
   achieved by renaming the elements of $\Sigma_1$ or $\Sigma_2$.}.
Then $(\Sigma_1 \cup \Sigma_2, T_1 \cup T_2)$ is a presentation of the free
 product $\g_1 \ast \g_2$.
As plain groups are finite free products of finite and free groups, this process
 results in 2-monadic convergent presentations for plain groups.
On the other hand it has been shown (compare \cite{AvMaOt86})
 that the class of plain
 groups is exactly the class of groups allowing finitely generated 2-monadic
 convergent  presentations.
It is obvious that the classes of finite groups and free groups are
 subclasses of this class.

Another combination of finite and free groups occurs in the
 description of context-free groups.
A context-free group $\g$ has a free normal subgroup $\free$ of finite index.
Let the free group $\free$ be generated by $X = \{ x_1, \ldots, x_n \}$ and let
  ${\cal E}$ be a finite group such that  ${\cal E} \cong \g/\free$ and
  $({\cal E} \backslash \{ \lambda \}) \cap (X \cup X^{-1}) = \emptyset$.
Then we can assume that an element  $g \in \g$ can be presented
 as $g \id ew$ for some $e \in {\cal E}$, $w \in \free$.
For all $e \in {\cal E}$ let $\phi_e : X \cup X^{-1} \myr \free$ be a function
 such that $\phi_{\lambda}$ is the inclusion and for all $x \in X \cup X^{-1}$,
 $\phi_e(x) = \inv{e} \mm_{\g} x \mm_{\g} e$, i.e., $\phi_e$ is a
 conjugation homomorphism.
For all $e_1,e_2 \in {\cal E}$ let $z_{e_1,e_2} \in \free$ such that
 $z_{e_1,\lambda} \id z_{\lambda,e_1} \id \lambda$ and for all $e_1,e_2 \in {\cal E}$
 with $e_1 \mm_{\cal E} e_2 =_{\cal E} e_3$, let $e_1 \mm_{\g} e_2 \id e_3z_{e_1,e_2}$. 
Then we can set $\Sigma = ({\cal E} \backslash \{ \lambda \}) \cup X \cup X^{-1}$
 and let $T$ consist of the following rules: 
\begin{tabbing}
XX\=XXXX\=XXX\=XXXXXXX\= XXXXXXXXXXXXXXXXXXXXXXXX\= \kill
\>$xx^{-1}$ \> $\myr$ \> $\lambda$ \> and  \\
\>$x^{-1}x$ \> $\myr$ \> $\lambda$ \> for all $x \in X$,  \\
\>$e_1e_2$ \>  $\myr$ \>  $e_3z_{e_1,e_2}$     \> 
   for all $e_1,e_2 \in {\cal E} \backslash \{ \lambda \}, e_3
   \in {\cal E}$ such that $e_1 \mm_{\cal E} e_2 =_{\cal E} e_3$, \\
\>$xe$ \>  $\myr$ \>  $e \phi_e(x)$ \> and  \\
\>$x^{-1}e$ \>  $\myr$ \>  $e \phi_e(x^{-1})$ \> for all
        $e \in {\cal E} \backslash \{ \lambda \}, x \in X$.
\end{tabbing}
 $(\Sigma , T)$ is called a
 \index{virtually free presentation}\index{presentation!virtually
   free}\betonen{virtually free presentation} and is convergent
 (compare \cite{CrOt94}).
In fact it can be shown that
a group $\g$ has a virtually free presentation if and only if it is
context-free.

Another class of groups allowing special presentations we want to present here are the
polycyclic groups, which include the Abelian and nilpotent
groups (compare  \cite{Wi89}).

Let $\Sigma = \{a_1,a_1^{-1}, \ldots, a_n, a_n^{-1} \}$ be a finite alphabet
 and for $1 \leq k \leq n$ we define the subsets
 $\Sigma_{k} = \{ a_i, a_i^{-1} \mid k \leq i \leq n \}$,
 $\Sigma_{n+1} = \emptyset$.
We first distinguish several particular classes of rules over $\Sigma$.
\begin{definition}~\\
{\rm
Let $i,j \in \{ 1, \ldots, n \}$,  $j>i$ and $\delta, \delta' \in \{ 1, -1 \}$.
\begin{enumerate}
\item A rule $a_j^\delta a_i^{\delta'} \myr  a_i^{\delta'} a_j^\delta$
      is called a \index{CAB-rule}\index{rule!CAB-}\betonen{CAB-rule}.
\item A rule $a_j^\delta a_i^{\delta'} \myr  a_i^{\delta'} a_j^\delta z$,
      $z \in \Sigma_{j+1}^*$ is called a \index{CNI-rule}\index{rule!CNI-}\betonen{CNI-rule}.
\item A rule $a_j^\delta a_i^{\delta'} \myr  a_i^{\delta'}z$,
      $z \in \Sigma_{i+1}^*$ is called a \index{CP-rule}\index{rule!CP-}\betonen{CP-rule}.
\dend
\end{enumerate}
}
\end{definition}
\begin{definition}~\\
{\rm
For ${\rm X} \in \{ \mbox{AB, NI, P} \}$ a subset $C$ of $\Sigma^*
\times \Sigma^*$ is called a \index{commutation-system}\betonen{commutation-system} if
\begin{enumerate}
\item $C$ contains only CX-rules, and
\item for all $1 \leq i < j \leq n$ and for all  $\delta, \delta' \in \{ 1, -1 \}$ there
      is exactly one rule  $a_j^\delta a_i^{\delta'} \myr r$ in $C$.
\dend
\end{enumerate}
}
\end{definition}
\begin{definition}~\\
{\rm
For $1 \leq i \leq n$ a rule $a_i^m \myr r$ where $m \geq 1$, 
 $r \in \Sigma_{i+1}^*$ is called a 
 \index{positive P-rule}\index{rule!positive P-}\betonen{positive P-rule}\/ and  
 a rule $a_i^{-1} \myr uv$ where $u \in \{ a_i \} ^*$ and
 $v \in \Sigma_{i+1}^*$ 
 is called a 
 \index{negative P-rule}\index{rule!negative P-}\betonen{negative P-rule}.
Then a subset  $P$ of $\Sigma^* \times \Sigma^*$ is called a \index{power system}\betonen{power system} if
\begin{enumerate}
\item $P$ contains only positive and negative P-rules.
\item For all $1 \leq i \leq n$ there is a negative P-rule $a_i^{-1}
  \myr uv$ in $P$ if and only if there also
      is a positive P-rule of the form  $a_i^m \myr r$ with $m \geq 1$ in $P$.
\item For all  $1 \leq i \leq n$ there is at most one negative P-rule $a_i^{-1} \myr uv$ and at most
      one positive P-rule  $a_i^m \myr r$ in $P$.
\dend
\end{enumerate}
}
\end{definition}
In combining these rule systems  we can characterize special group presentations.
Let  ${\rm X} \in \{ \mbox{AB, NI, P} \}$.
A  presentation $(\Sigma,T)$ is called a \index{CX-string-rewriting
  system}\index{string-rewriting system!CX-}\betonen{CX-string-rewriting system}\/ (CX-system)
 if $T= C \cup I$ where $C$ is a commutation system and $I$
 contains the trivial rules,
 i.e., $I = \{ a_ia_i^{-1} \myr \lambda, a_i^{-1}a_i \myr \lambda | 1 \leq i \leq n \}$.
It is called a  
 \index{PCX-string-rewriting system}\index{string-rewriting system!PCX-}\betonen{PCX-string-rewriting system}\/ (PCX-system)
 if $T = C \cup P \cup I$ where $T$ additionally includes a power
 system $P$.
The motivation for such presentations stems from the fact that they can
be used to characterize special classes of groups.
\begin{theorem}~\\
{\sl
For  a finitely presented group  $\g$ the following statements hold:
\begin{enumerate}
\item  $\g$ is Abelian if and only if there is a PCAB-system presenting $\g$.
\item  $\g$ is nilpotent if and only if there is a PCNI-system presenting $\g$.
\item  $\g$ is polycyclic if and only if there is a PCP-system presenting $\g$.
\theoend\ohnebeweis
\end{enumerate}
}
\end{theorem}
Using a \betonen{syllable ordering} Wi{\ss}mann has shown that  a PCX-system $(\Sigma,T)$ is
 a Noetherian string-rewriting system and he gave a completion procedure for such systems
 which terminates with an output that is again a PCX-system of the same type.
\begin{definition}\label{def.syllable}~\\
{\rm
Let $\Sigma$ be an alphabet and $\succ$ a partial ordering on
 $\Sigma^*$.
We define an ordering $\succ^{\rm lex}$ on m-tuples over $\Sigma^*$ as
follows:
$$(u_0, \ldots, u_m) \succ^{\rm lex} (v_0, \ldots, v_m)$$
$$\mbox{if and only if}$$
$$\mbox{there exists } 0 \leq k \leq m \mbox{ such that }
  u_i = v_i \mbox{ for all } 0 \leq i < k \mbox{ and } u_k \succ
v_k.$$
Let $a \in \Sigma$.
Then every $w \in \Sigma^*$ can be uniquely decomposed with respect
 to $a$ as $w \id w_0aw_1 \ldots aw_k$, where $|w|_a = k \geq 0$ and
 $w_i \in (\Sigma \backslash \{ a \})^*$.
Given a total precedence\footnote{By a precedence on an alphabet we mean a partial ordering on its letters.} $\triangleright$ on $\Sigma$ we can define a \betonen{syllable ordering} by
$$u >_{{\rm syll}(\Sigma)} v$$
$$\mbox{if and only if}$$
$$|u|_a > |v|_a \mbox{ or} $$
$$|u|_a = |v|_a \mbox{ and } (u_0, \ldots, u_m) >_{{\rm syll}(\Sigma \backslash
  \{ a \})}^{\rm lex}
                                  (v_0, \ldots, v_m)$$
where $a$ is the largest letter in $\Sigma$ according to
$\triangleright$ and $(u_0, \ldots, u_m)$, $(v_0, \ldots, v_m)$
 are the decompositions of $u$ and $v$ with respect to $a$ in case
 $|u|_a = |v|_a = m$.
In case we compare the tuples in reverse order, i.e.,
 $$(u_m, \ldots, u_0) >_{{\rm syll}(\Sigma \backslash \{ a \})}^{\rm lex} (v_m, \ldots, v_0)$$
we say that the syllable ordering has 
 \index{syllable ordering!status right}\betonen{status right}.
\dend
}
\end{definition}
The total precedence used on an alphabet $\Sigma = \{ a_i,
a_i^{-1} \mid 1 \leq i \leq n \}$ in our setting is $a_1^{-1} \succ a_1 \succ \ldots
a_i^{-1} \succ a_i \succ \ldots \succ a_n^{-1} \succ a_n$.
Using the syllable ordering induced by this precedence 
 we can give a characterization of the
 elements of our group as
 a subset of the set of \index{ordered group
   word}\index{group!ordered group word}\index{word!ordered group}
 \betonen{ordered group words} $\ord(\Sigma) = \ord(\Sigma_1)$, where we define
 $\ord(\Sigma_{i})$ recursively by
$\ord(\Sigma_{n+1}) = \{ \lambda \}$, and 
$\ord(\Sigma_i) = \{ w \in \Sigma_i^* \mid w \id uv \mbox{ for some } u \in \{ a_i \}^*
                     \cup \{ a_i^{-1} \}^* ,
                     v \in \ord(\Sigma_{i+1}) \}$.
Further with respect to $T$ we define some constants $\epsilon_T(i)$
for $1 \leq i \leq n$ by setting
$$\epsilon_T(i) = \left\{ \begin{array}{r@{\quad\quad}l}
                  \infty & \mbox{if $T$ contains no P-rules for $a_i$} \\
                  m         & \mbox{if $T$ contains a P-rule  $a_i^m \myr r$ for some unique $m>0$.}
                  \end{array} \right.$$

One can show that using the syllable ordering for orienting $T$ we get 
$$\irr(T) = \{  a_1^{i_1} \ldots a_{n\phantom{1}}^{i_n} | i_1, \ldots, i_n \in \z, \mbox{ and if } \epsilon_R(j) \neq \infty
               \mbox{ then } 0 \leq i_j \leq  \epsilon_T(j) \}.$$
For example the semi-Thue system $(\Sigma,T)$ where $T=C \cup I$ such
that we have 
 $C= \{ a_j^\delta a_i^{\delta'} \myr  a_i^{\delta'} a_j^{\delta} \mid  1 \leq i < j \leq n,
  \delta, \delta' \in \{ 1, -1 \} \}$ and
 $I = \{ a_ia_i^{-1} \myr \lambda, a_i^{-1}a_i \myr \lambda \mid 1 \leq i \leq n \}$
 is a presentation of the free commutative group generated by $\{ a_1,
 \ldots, a_n \}$ and we have $\irr(T) = \ord(\Sigma)$.

In  restricting the syllable ordering introduced in definition
 \ref{def.syllable} to ordered group words this gives us
 $a_1^{i_1} \ldots a_{n\phantom{1}}^{i_n} \syll a_1^{j_1} \ldots
 a_{n\phantom{1}}^{j_n}$ if and only if for some $1 \leq d \leq n$ we
 have $i_l = j_l$ for all $1 \leq l \leq d-1$ and $i_d >_{\z} j_d$
 with
 \[ \alpha <_{\z} \beta \mbox{ iff } \left\{ \begin{array}{l} 
                          \alpha \geq 0 \mbox{ and } \beta <0\\ 
                          \alpha \geq 0 , \beta  >0 \mbox{ and }
                          \alpha < \beta \\
                          \alpha <0, \beta < 0 \mbox{ and } \alpha > \beta
                                   \end{array} \right. \]
where $\leq$ is the usual ordering on $\z$.
We then call $a_d$ the \betonen{distinguishing
  letter}\index{distinguishing letter} between the two ordered group words.

Let us continue by giving some further information on
 nilpotent groups.
The following lemma gives syntactical information on the results of
multiplying a letter by special ordered group words.
\begin{lemma}\label{lem.nilpotentmultiplication}~\\
{\sl
Let $\g$ be a nilpotent group with a convergent PCNI-presentation
$(\Sigma, T)$.
Further for some $1 \leq j < i \leq n$ let $w_1 \in \ord(\Sigma \backslash
\Sigma_{j})$, $w_2 \in \ord(\Sigma_{i+1})$.
Then we have
$a_i \mm w_1 \id w_1a_iz_1$ and
$w_2 \mm a_i \id a_iz_2$
for some $z_1, z_2 \in \ord(\Sigma_{i+1})$.
\lemend\ohnebeweis
}
\end{lemma}
In section \ref{section.algebra} it was stated that an arbitrary
 finitely generated nilpotent group $\g$ has a normal subgroup
 ${\cal N}$ such that ${\cal N}$ is torsion-free nilpotent and 
 $\g/{\cal N}$ is finite.
Furthermore, torsion-free nilpotent groups have a central series
 with infinite cyclic factors and, therefore, we can assume that
 ${\cal N}$ has a CNI-presentation, i.e., a
 presentation containing no power rules (compare \cite{Wi89}).
Now every element $g \in \g$ can be uniquely expressed
 in the form $g \id ew$ where $e \in \g/{\cal N}$ and $w$ is an
 ordered group word in ${\cal N}$.
We can apply the same technique used for context-free groups
 to give a presentation of $\g$ in terms of ${\cal N}$ and 
 $\g/{\cal N}$.
Let $(\Sigma,C \cup I)$ be a CNI-presentation of ${\cal N}$ and 
 ${\cal E} \cong \g/{\cal N}$ such that 
 $({\cal E} \backslash \{ \lambda \}) \cap \Sigma = \emptyset$.
For all $e \in {\cal E}$ let $\phi_e : \Sigma \myr {\cal N}$ be a function
such that $\phi_{\lambda}$ is the inclusion and for all $a \in \Sigma$,
 $\phi_e(a) = \inv{e}\mm_{\g} a \mm_{\g} e$.
For all $e_1,e_2 \in {\cal E}$ let $z_{e_1,e_2} \in {\cal N}$ such that
 $z_{e_1,\lambda} \id z_{\lambda,e_1} \id \lambda$ and for all
 $e_1,e_2, e_3 \in {\cal E}$
 with $e_1 \mm_{\cal E} e_2 =_{\cal E} e_3$, $e_1 \mm_{\g} e_2 \id e_3z_{e_1,e_2}$. 
Let $\Gamma = ({\cal E} \backslash \{ \lambda \}) \cup \Sigma$
 and let $T$ consist of the sets of rules $C$ and $I$, and
 the additional rules: 
\begin{tabbing}
XX\=XXXX\=XXX\=XXXXXXX\= XXXXXXXXXXXXXXXXXXXXXXXX\= \kill
\>$e_1e_2$ \>  $\myr$ \>  $e_3z_{e_1,e_2}$     \> 
   for all $e_1,e_2 \in {\cal E} \backslash \{ \lambda \}, e_3
   \in {\cal E}$ such that $e_1 \mm_{\cal E} e_2 =_{\cal E} e_3$, \\
\>$ae$ \>  $\myr$ \>  $e\phi_e(a)$ \> for all
        $e \in {\cal E} \backslash \{ \lambda \}, a \in \Sigma$.
\end{tabbing}
Then $(\Gamma,T)$ is a convergent presentation of $\g$ as an extension of ${\cal
  N}$ by ${\cal E}$.

Note that finitely generated commutative groups can also be treated as
a special case of nilpotent groups.
But they can also be viewed as special commutative monoids.
Therefore, let us close this section with a short remark on presentations of
finitely presented commutative monoids.
If $\m$ is finitely presented by  a semi-Thue system $(\Sigma,T \cup C)$, where
$C$ is the commutator system for $\Sigma= \{ a_1, \ldots , a_n \}$,
 in general we cannot expect this presentation to be  convergent or to
 allow an equivalent finite convergent
 semi-Thue system\footnote{For example take $\Sigma = \{ a,b,c \}$ and
   $T = \{ (abc, \lambda ) \} \cup \{ (ba, ab), (ca, ac), (cb, bc)
   \}$. Then no equivalent finite convergent semi-Thue system exists.}.
But finitely generated commutative monoids always allow finite
convergent presentations in terms of semi-Thue systems modulo commutativity.
Let us start to specify such presentations by giving some basic definitions.
\begin{definition}\label{def.free.commutative.monoid}~\\
{\rm
For an alphabet $\Sigma = \{ a_1, \ldots, a_n \}$ let
 ${\cal T}_{\Sigma} = \{ a_1^{i_1} \ldots a_{n\phantom{1}}^{i_n} \mid i_j
 \in \n \}$ denote the set of \index{ordered words}\betonen{ordered words}\/ over $\Sigma$.
 We can define a mapping $\mm_{T_{\Sigma}} : {\cal T}_{\Sigma} \times
 {\cal T}_{\Sigma} \myr {\cal T}_{\Sigma}$ by setting 
 $a_1^{i_1} \ldots a_{n\phantom{1}}^{i_n} \mm_{T_{\Sigma}}
  a_1^{j_1} \ldots a_{n\phantom{1}}^{j_n} \id
  a_1^{i_1 + j_1} \ldots a_{n\phantom{1}}^{i_n + j_n}$.
Then $\mm_{T_{\Sigma}}$ is an associative binary operation on ${\cal T}_{\Sigma}$ with
  identity $\lambda = a_1^{0} \ldots a_{n\phantom{1}}^{0}$.
Thus, ${\cal T}_{\Sigma}$ together with $\mm_{T_{\Sigma}}$ and
$\lambda$ is a monoid, namely the 
  \index{free commutative!monoid}\index{monoid!free commutative}\betonen{free
  commutative monoid generated by $\Sigma$}.
\dend
}
\end{definition}
We will write ${\cal T}$ is case no confusion is likely to arise.
\begin{definition}~\\
{\rm
Let $\Sigma$ be a finite alphabet and ${\cal T}$ as in definition \ref{def.free.commutative.monoid}.
\begin{enumerate}
\item A \index{semi-Thue system!modulo commutativity}\betonen{semi-Thue system modulo commutativity}
  $T_c$ is a subset of
  ${\cal T} \times {\cal T}$.
  The elements $(l,r)$ of $T_c$ are called
  \index{rule}\betonen{rules}\/ and will often be
  written as $l \myr r$.
\item The  \index{single-step reduction relation}\index{reduction!single-step}\betonen{single-step reduction relation}\/
  on ${\cal T}$ induced
  by a semi-Thue system  modulo commutativity $T_c$ is defined as follows:
  For any $u,v$ in ${\cal T}$, $u \myr_{T_c} v$ if and only if there exist
  $x$ in ${\cal T}$ and $(l,r)$ in $T_c$ such that $u = l \mm_{{\cal
      T}} x$ and $v = r \mm_{{\cal T}} x$.
  The \index{reduction relation}\betonen{reduction relation}\/ 
  on ${\cal T}$ induced by $T_c$ is the
  reflexive transitive closure of $\myr_{T_c}$ and is denoted by
  $\red{*}{\myr}{}{T_c}$.
  The reflexive transitive symmetric closure is defined as usual and
  denoted by $\red{*}{\lr}{}{T_c}$.
\dend  
\end{enumerate}
}
\end{definition}
Recall that for $u$, $v \in {\cal T}$, $u$ is called a divisor of $v$ if there
exists an element $w \in {\cal T}$ such that $u \mm_{{\cal T}} w = v$.
The least common multiple of two terms $u \id a_1^{i_1} \ldots
a_n^{i_n}$ and $v \id a_1^{j_1} \ldots a_n^{j_n}$ is defined as
$\lcm(u,v) = a_1^{\max\{i_1, j_1 \}} \ldots a_{n\phantom{1}}^{\max \{ i_n, j_n \}}$
and is the ``shortest'' term that has both $u$ and $v$ as divisors.
Hence this term gives rise to a critical situation in case $u$ and $v$
are left hand sides of rules and $\lcm(u,v) \neq u \mm_{{\cal T}} v$.
\begin{definition}~\\
{\rm
A partial ordering $\preceq$ on ${\cal T}$ is called
\betonen{admissible}\index{admissible
  ordering}\index{ordering!admissible} (with respect to $\mm_{T_{\Sigma}}$)
if for all $u,v,x$ in ${\cal T}$
\begin{enumerate}
\item $\lambda \preceq u$, and
\item $u \prec v$ implies $u \mm_{{\cal T}} x \prec v\mm_{{\cal T}} x$.
\dend
\end{enumerate}
}
\end{definition}
Then if $\preceq$ is an admissible total ordering on ${\cal
  T}$ (which is of course well-founded) and $u$ is a proper divisor of $v$ this implies $u \prec v$.
Notice that given an admissible total ordering,
 a finite semi-Thue system modulo commutativity always has
 a finite equivalent convergent  semi-Thue system modulo commutativity with respect to this ordering.
Then in case a commutative monoid $\m$ is presented by a finite semi-Thue system modulo
commutativity $(\Sigma, T_c)$ which is convergent with respect to an admissible
total ordering $\succeq_{T_c}$, 
for $u$, $v\in \m$ we have $u \mm_{{\cal T}} v \succeq_{T_c} u
\mm_{\m} v$.

In the following chapters, if not stated otherwise, our monoids and
groups are presented by finite convergent semi-Thue systems
respectively finite convergent semi-Thue systems modulo commutativity
which are convergent with respect to some admissible well-founded 
total ordering\footnote{This ordering used for completion will
sometimes be called the completion ordering of the system and
we will only consider convergent systems having such an ordering in this thesis. Notice that in general there are convergent systems which allow no such completion ordering (see e.g. \cite{Es86}).}  in the respective setting, i.e., 
 if $\succeq$ is the completion ordering, then for all rules $(l,r)$ in
 the set $T$, $l \succ r$ holds.

%% file: ideals.tex
\chapter{Non-Commutative Gr\"obner Bases}\label{chapter.ideals}
\spruch{10}{8}{Was du tust, bedenke das Ende.}{Sirach 7:40}

In this chapter we will now give a short outline on the results on ideals in
non-commutative rings.
Further we state some undecidability results which give a limit to what can be
achieved in certain non-commutative structures.

{\bf Section 3.1:} Mora generalized Buchberger's ideas to finitely
  generated free monoid rings.
  The notions of reduction, s-polynomials and Gr\"obner bases can be
  carried over and a characterization of Gr\"obner bases in a finitary
  manner is possible.
  Hence it is decidable, whether a finite set of polynomials is a Gr\"obner basis.
  But this characterization can no longer be used to give a terminating
  completion algorithm as none exists.
  Mora developed an enumerating procedure for two-sided ideals and
  in restricting the attention to finitely generated right ideals and
  prefix reduction, a terminating completion algorithm was provided. 

{\bf Section 3.2:} Kandri-Rody and Weispfenning have shown that the
  ideal membership problem for two-sided ideals in a free monoid ring
  with two generators is undecidable in general.
  We here show the equivalence of the word problem for semi-Thue
  systems to a restricted version of the ideal membership problem in a free
  monoid ring.
  Hence the existence of a semi-Thue system over two letters with
  undecidable word problem implies the undecidability of the ideal
  membership problem for a free monoid ring with two generators.
  The same equivalence can be shown for the word problem in groups and
  a restricted version of the ideal membership problem in a free group ring.
  Again the existence of a group system with four letters (two
  generators plus their inverses) with undecidable word problem
  implies the undecidability of the ideal membership problem for free
  group rings with more than one generator.
  Finally we show that it is undecidable whether a finite Gr\"obner
  basis with respect to an admissible ordering in a free monoid ring
  with more than one generator exists.

{\bf Section 3.3:} A class of non-commutative rings where finitely
  generated left, right
  and two-sided ideals have finite Gr\"obner bases are the solvable
  polynomial rings and the skew polynomial rings.
  We sketch Weispfenning's approach to the latter structure as it can
  be viewed as a monoid ring.

For more information on Mora's or Weispfenning's approach
the reader is referred to \cite{Mo85} and \cite{We92}.
\section{The Free Monoid Ring}\label{section.mora}
Let $\freemonoid$ be a free monoid  over a finite alphabet $\Sigma$.
It has been shown by Mora  that Buchberger's ideas generalize
naturally to the free monoid ring over $\myk$ and $\freemonoid$.
In this section we now sketch his approach which can be found more explicitly for
example in
\cite{Mo85}.

Because of non-commutativity it becomes important to distinguish one
 and two-sided ideals and the
results gained differ from the commutative case.
This is e.g. due to the fact that for more than one variable the
corresponding free monoid ring is no longer Noetherian and that for
the set $\freemonoid$ no equivalent to Dickson's lemma in ${\cal T}$
holds, i.e., we cannot prove termination with the methods used in the
commutative case and it can be shown that finite Gr\"obner bases need
not exist, even for finitely generated two-sided ideals.

Let $\succeq$ be a total admissible well-founded ordering
 on $\freemonoid$.
This ordering can be extended to $\myk[\freemonoid]$ and used to
 distinguish the head term, head
 coefficient, head
 monomial and the reduct of a polynomial.
Two-sided reduction then can be defined naturally.
\begin{definition}[Mora]\label{def.red.mora}~\\
{\rm
Let $p,f$ be two non-zero polynomials in $\myk[\freemonoid]$. 
We say $f$ 
 \index{reduction!Mora}\index{reduction!in a free polynomial ring}\index{Mora's reduction}\betonen{reduces} 
 $p$ to $q$ at a monomial $\alpha \skm t$ of $p$ in one step,
 denoted by $p \red{}{\myr}{m}{g} q$, if
\begin{enumerate}
\item[(a)] $u\hterm(f)v \id t$ for some $u,v \in \freemonoid$, and
\item[(b)] $q = p - \alpha \skm \hc(f)^{-1} \skm u \mrm f \mrm v$.
\end{enumerate}
We write $p \red{}{\myr}{m}{f}$ if there is a polynomial $q$ as defined above and $p$ is then called reducible by $f$. 
Further, we can define $\red{*}{\myr}{m}{}, \red{+}{\myr}{m}{}$, and
 $\red{n}{\myr}{m}{}$ as usual.
Reduction by a set $F \subseteq \myk[\freemonoid]$ is denoted by
 $p \red{}{\myr}{m}{F} q$ and abbreviates $p \red{}{\myr}{m}{f} q$
 for some $f \in F$,
 which is also written as  $p \red{}{\myr}{m}{f \in F} q$.
\dend
}
\end{definition}
When studying right ideals this reduction is restricted to prefix
reduction\footnote{A similar approach is possible to study left ideals.}.
\begin{definition}\label{def.redprefix}~\\
{\rm
Let $p, f$ be two non-zero polynomials in $\myk[\freemonoid]$. 
We say $f$ \index{prefix!reduction}\index{reduction!prefix}
 \betonen{prefix reduces}
 $p$ to $q$ at a monomial
 $\alpha \skm t$ of $p$ in one step, denoted by $p \red{}{\myr}{p}{f} q$, if
\begin{enumerate}
\item[(a)] $\hterm(f)w \id t$ for some $w \in \m$,
            i.e., $\hterm(f)$ is a prefix of $t$, and
\item[(b)] $q = p - \alpha \skm \hc(f)^{-1} \skm f \mrm w$.
\end{enumerate}
We write $p \red{}{\myr}{p}{f}$ if there is a polynomial $q$ as defined
above and $p$ is then called  prefix reducible by $f$. 
Further, we can define $\red{*}{\myr}{p}{}, \red{+}{\myr}{p}{}$,
 $\red{n}{\myr}{p}{}$  as usual.
Prefix reduction by a set $F \subseteq \myk[\freemonoid]$ is denoted by
 $p \red{}{\myr}{p}{F} q$ and abbreviates $p \red{}{\myr}{p}{f} q$
 for some $f \in F$,
 which is also written as  $p \red{}{\myr}{p}{f \in F} q$.
\dend
}
\end{definition}
Note that Mora's two-sided reduction and prefix reduction have all nice properties connected to
 reductions in structures with an admissible ordering and
 capture the ideal respectively right ideal congruence in the free monoid ring.
Furthermore, the translation lemma holds and the respective reductions are preserved 
 under two-sided respectively right-sided multiplication with terms.
We can define Gr\"obner bases for two-sided
 ideals and prefix Gr\"obner bases for right ideals and characterize them in a natural way by corresponding s-polynomials.
\begin{definition}~\\
{\rm
A  set $G \subseteq \myk[\freemonoid]$ is called a \betonen{Gr\"obner basis}\/
 with respect to
 the reduction $\red{}{\myr}{m}{}$, if
\begin{enumerate}
\item[(i)] $\red{*}{\lr}{m}{G}  = \;\; \equiv_{\ideal{}{}(G)}$, and
\item[(ii)] $\red{}{\myr}{m}{G}$ is confluent.
\dend
\end{enumerate}
}
\end{definition}
\begin{definition}~\\
{\rm
A  set $G \subseteq \myk[\freemonoid]$ is called a \betonen{Gr\"obner basis}\/
 with respect to
 the reduction $\red{}{\myr}{p}{}$ or a \index{prefix!Gr\"obner basis}
 \index{Gr\"obner basis!prefix}
 \betonen{prefix Gr\"obner basis}, if
\begin{enumerate}
\item[(i)] $\red{*}{\lr}{p}{G}  = \;\; \equiv_{\ideal{r}{}(G)}$, and
\item[(ii)] $\red{}{\myr}{p}{G}$ is confluent.
\dend
\end{enumerate}
}
\end{definition}
\begin{definition}~\\
{\rm
Given two non-zero polynomials $p_{1}, p_{2} \in \myk[\freemonoid]$, 
 if there are $u,v \in \freemonoid$ such that either $\hterm(p_{1})u \id v\hterm(p_{2})$ and $|\hterm(p_1)|<|v|$ or $u\hterm(p_{1})v \id
 \hterm(p_{2})$,
  the corresponding \betonen{s-polynomial}\index{s-polynomial} is defined either as
 $$ \spol{m}(p_{1}, p_{2}, u, v)=\hc(p_1)^{-1} \skm p_1 \mrm u -\hc(p_2)^{-1} \skm
  v \mrm p_2 \mbox{ or as}$$
 $$ \spol{m}(p_{1}, p_{2}, u, v)=\hc(p_1)^{-1} \skm u \mrm p_1 \mrm v
 -\hc(p_2)^{-1} \skm p_2.$$
Let $U_{p_1,p_2}  \subseteq \freemonoid \times \freemonoid$ be the set
containing {\em all} such pairs $u,v \in\freemonoid$. 
\dend
}
\end{definition}
Notice that the sets $U_{p_1,p_2}$ are always finite and that they
 correspond to the critical pairs for semi-Thue systems as mentioned
 in section \ref{section.presentations}.
\begin{definition}~\\
{\rm
Given two non-zero polynomials $p_{1}, p_{2} \in \myk[\freemonoid]$, 
 if there is $w \in \freemonoid$ such that $\hterm(p_{1}) \id \hterm(p_{2})w$
  the \index{prefix!s-polynomial}\index{s-polynomial!prefix}\betonen{prefix
    s-polynomial} is defined as
 $$ \spol{p}(p_{1}, p_{2})=\hc(p_1)^{-1} \skm p_1 -\hc(p_2)^{-1} \skm p_2 \mrm w.$$
\dend
}
\end{definition}
In analogy to theorem \ref{theo.buchberger.completion} we can give the
 following characterizations of Gr\"obner bases respectively prefix
 Gr\"obner bases, which can be used to decide whether a finite set of
 polynomials is a respective basis and to give procedures to enumerate such bases.
\begin{theorem}~\\
{\sl
For a set of polynomials $F \subseteq \myk[\freemonoid]$,
 the following statements are equivalent:
\begin{enumerate}
\item $F$ is a  Gr\"obner basis.
\item For all $f_{k}, f_{l} \in F$, $(u,v) \in
  U_{f_k,f_l}$ we have 
  $ \spol{m}(f_{k}, f_{l}, u, v) \red{*}{\myr}{m}{F} 0$.
\mbox{\phantom{XX}}\ohnebeweis\theoend
\end{enumerate}
}
\end{theorem}
\procedure{Gr\"obner Bases in Free Monoid Rings [Mora]}%
{\vspace{-4mm}\begin{tabbing}
XXXXX\=XXXX \kill
\removelastskip
{\bf Given:} \> A finite set $F \subseteq  \myk[\freemonoid]$. \\
{\bf Find:} \> $G$, a Gr\"obner basis of $\ideal{}{}(F)$.
\end{tabbing}
\vspace{-7mm}
\begin{tabbing}
XX\=XX\=XX\= XXXX\=\kill
$G$ := $F$; \\
$B$ := $\{ (q_{1}, q_{2}) \mid q_{1}, q_{2} \in G \}$; \\
{\bf while} $B \neq \emptyset$ {\bf do} \\
\>      $(q_{1}, q_{2}) := {\rm remove}(B)$; \\
 \> {\kommentar \% Remove an element using a fair strategy}\\
\>      {\bf for all} s-polynomials $h \in S(q_1,q_2)$ {\bf do}  \\
\> \> {\kommentar \% $S(q_1,q_2)= \{\hc(q_1)^{-1} \skm u \mrm q_1 \mrm v -
               \hc(q_2)^{-1} \skm q_2 \mid u\hterm(q_1)v\id\hterm(q_2) \} \cup$} \\
\> \> {\kommentar \%   $\{\hc(q_1)^{-1} \skm q_1 \mrm u -
               \hc(q_2)^{-1} \mrm v \mrm q_2 
             \mid \hterm(q_1)u\id v\hterm(q_2) , |\hterm(q_1)| < |v|\}$} \\
\> \>    $h' := {\rm normalform}(h, \red{}{\myr}{m}{G})$; \\
\> \>    {\kommentar \% Compute a normal form using Mora's two-sided
                         reduction.} \\
\> \>   {\bf if} \>$h' \neq 0$  \\
\> \>            \>{\bf then}  \> $G := G \cup \{ h' \}$; \\
\> \>            \>            \> $B := B \cup \{ (f, h'), (h',f) \mid f \in G \}$; \\
\> \>   {\bf endif}\\
\>  {\bf endfor} \\
{\bf endwhile}
\end{tabbing}}
\begin{theorem}~\\
{\sl
For a set of polynomials $F \subseteq \myk[\freemonoid]$,
 the following statements are equivalent:
\begin{enumerate}
\item $F$ is a  prefix Gr\"obner basis.
\item For all polynomials $f_{k}, f_{l} \in F$  we have 
  $\spol{p}(f_{k}, f_{l}) \red{*}{\myr}{p}{F} 0$.
\mbox{\phantom{XX}}\ohnebeweis\theoend
\end{enumerate}
}
\end{theorem}
\procedure{Prefix Gr\"obner Bases in Free Monoid Rings [Mora]}%
{\vspace{-4mm}\begin{tabbing}
XXXXX\=XXXX \kill
\removelastskip
{\bf Given:} \> A finite set $F \subseteq  \myk[\freemonoid]$.\\
{\bf Find:} \> $G$, a prefix Gr\"obner basis of $\ideal{r}{}(F)$.
\end{tabbing}
\vspace{-7mm}
\begin{tabbing}
XX\=XX\= XXXX \= XXXX\=\kill
$G$ := $F$; \\
{\bf while} there is $g \in G$ such that $\hterm(g)$ is prefix
reducible by $G \backslash \{ g \}$ {\bf do} \\
\> $G$ := $G \backslash \{ g \}$; \\
\> $f$ := ${\rm normalform}(g,\red{}{\myr}{p}{G})$; \\
\> {\kommentar \% Compute a normal form using prefix  reduction.} \\
\> {\bf if} \>$f \neq 0$ \\
\>          \>{\bf then} \>$G$ := $G \cup \{ f \}$; \\
\> {\bf endif} \\
{\bf endwhile}
\end{tabbing}}

While termination for procedure {\sc Prefix Gr\"obner Bases in Free
  Monoid Rings} follows immediately from the fact that polynomials in
 $G$ are replaced by ``smaller'' polynomials and no cycles can occur,
  procedure {\sc Gr\"obner
  Bases in Free Monoid Rings} need not  terminate as the
ideal membership problem for finitely generated two-sided ideals in $\myk[\freemonoid]$ is
unsolvable in general.
The following section will give some insight into undecidability
results related to free monoid rings and free group rings.
\section{Undecidability Results}\label{section.undecidable}
Kandri-Rody and Weispfenning have shown in \cite{KaWe90} that 
 the ideal membership problem for finitely
 generated two-sided ideals is algorithmically unsolvable
 for the free monoid ring $\q[\{ X_1, X_2 \}^*]$ by reducing the 
halting problem for Turing machines to this problem.
Here we state a similar result by showing that  the word problem for
semi-Thue systems is equivalent to a restricted version of the ideal membership
problem in free monoid rings.
\begin{theorem}\label{theo.wp}~\\
{\sl
Let $(\Sigma, T)$ be a finite semi-Thue system and
 $P_T= \{ l-r \mid (l,r) \in T \}$.
Then for $u,v \in \freemonoid$ the following statements are equivalent:
\begin{enumerate}
\item $u \red{*}{\lr}{}{T} v$.
\item $u-v \in \ideal{}{\myk[\freemonoid]}(P_T)$.
\theoend
\end{enumerate}
}
\end{theorem}
\Ba{}~\\
\mbox{$1 \R 2:$ }
 Using induction on $k$ we show that $u \red{k}{\lr}{}{T} v$ implies
        $u - v \in \ideal{}{\myk[\freemonoid]}(P_T)$.
 In the base case $k = 0$ there is nothing to show, since
  $u-u =0 \in \ideal{}{\myk[\freemonoid]}(P_T)$.
 Thus let us assume that  $\tilde{u} \red{k}{\lr}{}{T} \tilde{v}$
  implies $ \tilde{u}-\tilde{v} \in \ideal{}{\myk[\freemonoid]}(P_T)$.
 Then looking at $ u \red{k}{\lr}{}{T} u_{k}  \red{}{\lr}{}{T} v$ we find
 $u_{k}  \red{}{\lr}{}{T} v$ with $(l_j,r_j) \in T$.
 Without loss of generality we can assume $u_k \id xl_jy$ for
  some $x,y \in \freemonoid$  thus
 giving us $v \id xr_jy$, and since multiplication in the
 free monoid is concatenation, $v$ can be expressed in terms of
 polynomials by
 $v = u_k - x\mrm (l_j-r_j) \mrm y$.
 As $u-v = u- u_k + x \mrm (l_j-r_j) \mrm y$ and $u-u_k \in
 \ideal{}{\myk[\freemonoid]}(P_T)$ our induction hypothesis
 yields $u-v \in \ideal{}{\myk[\freemonoid]}(P_T)$.

\mbox{$2 \R 1:$ }
 It remains to show that $u-v \in \ideal{}{\myk[\freemonoid]}(P_T)$
  implies $u \red{*}{\lr}{}{T} v$.
 We know $u-v= \sum_{j=1}^{n} \beta_j \skm x_j \mrm (l_{i_j} - r_{i_j})
 \mrm y_j$, where
  $\beta_j \in \myk^*, x_j,y_j \in \freemonoid$.
 Therefore, by  showing the following stronger result we are done:
 A representation
 $ u - v = \sum_{j=1}^m p_j$ where $p_j = \alpha_j \skm (w_j -w'_j)$,
 $\alpha_j \in \myk^*$ and $w_j \red{+}{\lr}{}{T} w'_j$ implies that
 $u \red{*}{\lr}{}{T} v$.
 Thus let $u - v = \sum_{j=1}^m p_j$ be such a representation.
 Depending on this  representation $\sum_{j=1}^m p_j$ 
 and the ordering $\succeq$ on $\freemonoid$ we can define
  $t = \max \{w_j,w'_j   \mid j = 1, \ldots m  \}$  and
  $K$ is the number of polynomials $p_j$ containing $t$ as a term.
 We will show our claim by induction on $(m,K)$, where 
  $(m',K') < (m,K)$ if and only if $m'<m$ or $(m'=m$ and $K'<K)$.
 In case $m=0$, then $u-v=0$ implies $u \id v$ and hence $u \red{0}{\lr}{}{T} v$.
 Now suppose $m>0$.
\\
 In case $K=1$, let $p_k$ be the polynomial containing
  $t$.
 Since we either have $p_k=\alpha_k \skm (t-w'_k)$
  or $p_k=\alpha_k \skm (w_k - t)$, where $\alpha_k \in \{ 1, -1 \}$,
  without loss of generality we can assume $u \id t$ and $p_k=t-w'_k$.
 Using $p_k$ we can decrease $m$ by subtracting $p_k$ from $u-v$ giving us
  $w'_k-v = \sum_{j=1,j \neq k}^{m} p_j$.
 Since $u \id t \red{*}{\lr}{}{T} w'_k$ and our induction
  hypothesis yields $w'_k \red{*}{\lr}{}{T} v$
  we can conclude $u \red{*}{\lr}{}{T} v$.
\\
 In case $K>1$ there are two polynomials $p_k,p_l$ in
  the corresponding representation containing the term $t$ and without loss
  of generality we can assume $p_k = \alpha_k \skm( t - w'_k)$ and $p_l =
  \alpha_l \skm (t -w'_l)$, as
 the cases where $p_k = \alpha_k \skm( w'_k - t)$ or  $p_l =
 \alpha_l \skm (w'_l - t)$ occur can be treated similarly by modifying
 the respective coefficient.
 If $w'_k \id w'_l$  we can immediately decrease $m$ 
  by substituting the occurrence of $p_k+p_l$ by $(\alpha_k + \alpha_l) \skm p_l$.
 Otherwise we can proceed as follows:
  \begin{eqnarray*}
  p_k + p_l & = & p_k \underbrace{-\alpha_k \skm \alpha_l^{-1} \skm p_l +
                  \alpha_k \skm \alpha_l^{-1} \skm p_l}_{=0}
                  + p_l \\
            & = & (\alpha_k \skm (w'_k - t) - \alpha_k \skm
            \alpha_l^{-1} \skm \alpha_l \skm (w'_l - t)) + (\alpha_k \skm \alpha_l^{-1} +1) \skm p_l\phantom{p_k \underbrace{-\alpha_k \skm \alpha_l^{-1} \skm p_l +
                  \alpha_k \skm \alpha_l^{-1} \skm p_l}_{=0}} \\
            & = & \underbrace{(-\alpha_k \skm w'_k + \alpha_k  \skm
                  w'_l)}_{= p'_k}
                  + (\alpha_k \skm \alpha_l^{-1} +1) \skm p_l \\
  \end{eqnarray*}
  where $p'_k = \alpha_k \skm (w'_l - w'_k)$,
  $w'_k \red{*}{\lr}{}{T} t \red{*}{\lr}{}{T} w'_l$
  and $w'_l \neq w'_k$.
 Therefore, in case $\alpha_k \skm \alpha_l^{-1} + 1=0$, i.e., $\alpha_k=-\alpha_l$, $m$ is decreased.
 On the other hand $p'_k$ does not contain $t$, i.e., $K$ will be
 decreased in any case.
\\
\qed
The existence of a finite semi-Thue system over an alphabet with
two elements having undecidable word problem yields that the ideal
membership problem for free monoid rings with more than one generator is 
undecidable.
In case the free monoid is generated by one element, we have decidable
ideal membership problem.
This is due to the fact that this corresponds to the ordinary polynomial 
 ring in one variable and there e.g. the Euclidean algorithm can be applied to
 solve the ideal membership problem.

Perhaps less obvious is that 
 the word problem for groups is similarly equivalent to a restricted version of the
 membership problem for ideals in a
 free group ring.
Let the group be presented by a semi-Thue system $(\Sigma, T \cup T_I)$ such
 that there exists an involution $\imath : \Sigma \myr \Sigma$ such that
 for all $a \in \Sigma$ we have $\imath(a) \neq a$, $\imath(\imath(a))=a$
 and $(\imath(a)a,\lambda),(a\imath(a),\lambda)\in T_I$.
\begin{theorem}~\\
{\sl
Let $(\Sigma, T \cup T_{I})$ be a finite group system with
  $T_{I}= \{(\imath(a)a, \lambda),(a\imath(a), \lambda) \mid
  a \in\Sigma \}$, i.e.,  $(\Sigma, T_{I})$ is a presentation of a free group $\free$.
Further we can associate a system of polynomials 
 $P_T= \{ l-r \mid (l,r) \in T \}$
 with $T$ and without loss of generality we can assume
 that $l$ and $r$ are in normal form with respect to $T_I$.
Then for $u,v \in \Sigma^*$ the following statements are equivalent:
\begin{enumerate}
\item $u \red{*}{\lr}{}{T \cup T_{I}} v$.
\item $u\nf{T_{I}}  -v\nf{T_{I}} \in \ideal{}{\myk[\free]}(P_T)$.
\theoend
\end{enumerate}
}
\end{theorem}
\Ba{}~\\
\mbox{$1 \R 2:$ }
 Using induction on $k$ we show  $u \red{k}{\lr}{}{T \cup T_{I}} v$ implies
        $u\nf{T_{I}} - v\nf{T_{I}} \in \ideal{}{\myk[\free]}(P_T)$.
 In the base case $k = 0$ we have $u \id v$ and, therefore,
  $u\nf{T_{I}}-u\nf{T_{I}} =0 \in \ideal{}{\myk[\free]}(P_T)$.
 Hence, let us assume that  $\tilde{u} \red{k}{\lr}{}{T \cup T_{I}} \tilde{v}$ 
  implies $ \tilde{u}\nf{T_{I}}-\tilde{v}\nf{T_{I}} \in
  \ideal{}{\myk[\free]}(P_T)$.
 Thus, looking at $ u \red{k}{\lr}{}{T \cup T_{I}} u_{k} 
 \red{}{\lr}{}{T \cup T_{I}} v$ we can distinguish the following cases:
      \begin{enumerate}
      \item $u_{k}  \red{}{\lr}{}{T} v$ with $(l,r) \in T$.
\\
            Without loss of generality we can assume $u_k \id xly$ and
            $v \id xry$ for some words $x,y \in \Sigma^*$.
            Now this gives us $$u\nf{T_{I}}-v\nf{T_{I}} = u\nf{T_{I}}-
            \underbrace{u_k\nf{T_{I}} + xly\nf{T_{I}}}_{=0} -
            xry\nf{T_{I}}$$ and $xly\nf{T_{I}} - xry\nf{T_{I}} = x
            \mrm (l-r) \mrm y$, where $\mrm$ denotes multiplication
            in $\myk[\free]$.
            By our induction hypothesis we know $u\nf{T_{I}}-u_k\nf{T_{I}}
            \in \ideal{}{\myk[\free]}(P_T)$ and, hence,  we get
            $u\nf{T_{I}}-v\nf{T_{I}} \in \ideal{}{\myk[\free]}(P_T)$.
      \item $u_{k}  \red{}{\lr}{}{T_{I}} v$ with $(a\imath(a),
            \lambda) \in T_{I}$\footnote{The case $(\imath(a)a,\lambda)
            \in T$ is similar.}.
\\
             Without loss of generality we can assume
              $u_k \id xa\imath(a)y$ for some $x,y \in \Sigma^*$ and $v
              \id xy$, i.e.,
              $u_k\nf{T_{I}}=v\nf{T_{I}}$ and therefore
              $u\nf{T_{I}}-v\nf{T_{I}} \in \ideal{}{\myk[\free]}(P_T)$.
      \end{enumerate}
\mbox{$2 \R 1:$ }
 It remains  to show that $u\nf{T_{I}}-v\nf{T_{I}} \in \ideal{}{\myk[\free]}(P_T)$
  implies $u \red{*}{\lr}{}{T \cup T_{I}} v$.
 We know $u\nf{T_{I}}-v\nf{T_{I}}= \sum_{j=1}^{n} \beta_j \skm x_j \mrm
  (l_{i_j} - r_{i_j}) \mrm y_j$, where
  $\beta_j \in \myk^*, x_j,y_j \in \free$.
 Therefore, by  showing the following stronger result we are done:
 A representation
 $ u - v = \sum_{j=1}^m p_j$ where $p_j = \alpha_j \skm (w_j -w'_j)$,
 $\alpha_j \in \myk^*$, $u,v,w_j,w'_j \in \free$ and
 $w_j \red{+}{\lr}{}{T} w'_j$ implies that
 $u \red{*}{\lr}{}{T} v$.
 Hence, let $u - v = \sum_{j=1}^m p_j$ be such a representation.
 Depending on this  representation $\sum_{j=1}^m p_j$ 
 and the ordering $\succeq$ on $\freemonoid$ we can define
  $t = \max \{w_j,w'_j   \mid j = 1, \ldots m  \}$  and
  $K$ is the number of polynomials $p_j$ containing $t$ as a term.
 We will show our claim by induction on $(m,K)$, where 
  $(m',K') < (m,K)$ if and only if $m'<m$ or $(m'=m$ and $K'<K)$.
 In case $m=0$, then $u-v=0$ implies $u=v$ and hence $u \red{0}{\lr}{}{T} v$\footnote{Remember that $u,v \in \free$, i.e., they are in normal form with respect to $T_I$.}.
 Now suppose $m>0$.
\\
 In case $K=1$, let $p_k$ be the polynomial containing
  $t$.
 Since we either have $p_k=\alpha_k \skm (t-w'_k)$
  or $p_k=\alpha_k \skm (w_k - t)$, where $\alpha_k \in \{ 1, -1 \}$,
  without loss of generality we can assume $u=t$ and $p_k=t-w'_k$.
 Using $p_k$ we can decrease $m$ by subtracting $p_k$ from $u-v$ giving us
  $w'_k-v = \sum_{j=1,j \neq k}^{m} p_j$.
 Since $u=t \red{*}{\lr}{}{T} w'_k$ and  our induction
  hypothesis yields $w'_k \red{*}{\lr}{}{T} v$
  we get $u \red{*}{\lr}{}{T} v$.
\\
 In case $K>1$ there are two polynomials $p_k,p_l$ in
  the corresponding representation containing the term $t$ and without loss
  of generality we can assume $p_k = \alpha_k \skm( t - w'_k)$ and $p_l =
 \alpha_l \skm (t -w'_l)$, as
 the cases where $p_k = \alpha_k \skm( w'_k - t)$ or $p_l =
 \alpha_l \skm (w'_l - t)$ occur can be treated similarly by modifying
 the respective coefficient.
 If $w'_k = w'_l$  we can immediately decrease $m$ 
  by substituting the occurrence of $p_k+p_l$ by $(\alpha_k + \alpha_l) \skm p_l$.
 Otherwise we can proceed as follows:
  \begin{eqnarray*}
  p_k + p_l & = & p_k \underbrace{-\alpha_k \skm \alpha_l^{-1} \skm p_l +
                  \alpha_k \skm \alpha_l^{-1} \skm p_l}_{=0}
                  + p_l \\
            & = & \underbrace{(-\alpha_k \skm w'_k + \alpha_k  \skm
                  w'_l)}_{= p'_k}
                  + (\alpha_k \skm \alpha_l^{-1} +1) \skm p_l \\
  \end{eqnarray*}
  where $p'_k = \alpha_k \skm (w'_l - w'_k)$,
  $w'_k \red{*}{\lr}{}{T} t \red{*}{\lr}{}{T} w'_l$
  and $w'_l \neq w'_k$.
 Hence, in case $\alpha_k \skm \alpha_l^{-1} + 1=0$, i.e., $\alpha_k=-\alpha_l$, $m$ is decreased.
 On the other hand $p'_k$ does not contain $t$, i.e., $K$ will be
 decreased in any case.
\\
\qed
As before, the existence of a finite group presentation over four
letters (resulting from two generators) with unsolvable word problem
implies that the ideal membership problem for free group rings with
more than one generator is undecidable.
Groups with one generator are known to have decidable word problem 
as this case corresponds to the ring of Laurent polynomials\footnote{Notice that the Laurent polynomials can be treated as a quotient of
 an ordinary polynomial ring.} in the (commutative) free group with one generator.

By theorem \ref{theo.wp} we know that finitely generated ideals in free monoid rings need
 not admit finite Gr\"obner bases.
It even is possible for a finitely generated ideal to admit a finite Gr\"obner basis with respect 
 to one admissible ordering and none with respect to another admissible ordering.
On the other hand, Mora provided a procedure which given an admissible
ordering enumerates a Gr\"obner basis with respect to this ordering.
This procedure terminates in case a finite Gr\"obner basis with
respect to the given ordering exists.
Hence the question might arise, whether it is possible to decide for a 
 given finite set of polynomials whether there exists an admissible ordering 
 such that  procedure {\sc Gr\"obner  Bases in Free Monoid rings} terminates,
 hence whether a finite Gr\"obner basis exists. 
This question turns out to be undecidable.
\begin{theorem}\label{theo.mora.unentscheidbar}~\\
{\sl
It is undecidable, whether a finitely generated ideal
 has a finite Gr\"obner basis in
 the free monoid ring $\myk[\{ s,t \}^*]$ with respect
 to a fixed two-sided reduction as defined in 
 definition \ref{def.red.mora}.
}
\end{theorem}
\Ba{}~\\
Using the technique described by \'{O}'D\'{u}nlaing in \cite{OD83}
Madlener and Otto have shown that the following
 question is undecidable (\cite{Ot93}):

{\sl
Let $\succeq$ be a compatible well-founded partial ordering on $\{s,t \}^*$ such that
 $s \succ \lambda$ and $t  \succ \lambda$ both hold. 
\\
Given a finite Thue system $T$ on $\{ s,t \}$.
Is there a finite and confluent system $T'$ on $\{ s,t \}$
 that is equivalent to $T$ and based on $\succ$?
}

To prove our claim we show that the answer for $T$ is ``yes'' if and only if the
 set of polynomials $P_{T}$
 associated to $T$
 has a finite Gr\"obner basis in   $\myk[\{ s,t \}^*]$ with respect to $\succ$.
If there is an equivalent, finite presentation $(\{s,t\} , T')$  convergent with respect to
 $\succ$, then the set $P_{T'}$ is a finite Gr\"obner basis of $P_{T}$ in 
 $\myk[\{ s,t \}^*]$.
This follows as the Thue reduction $\red{}{\lr}{}{T'}$ can be
simulated by the symmetric closure of the reduction
$\red{}{\myr}{m}{P_T}$ in $\myk[\{ s,t \}^*]$
 (compare definition \ref{def.red.mora}).
Thus it remains 
 to show that in case $P_{T}$ has a finite Gr\"obner basis
 in $\myk[\{ s,t \}^*]$, there exists
 a finite Gr\"obner basis $G$ such that for all $g \in G$
 we have $g = u-v$, where $u,v \in \{s,t\}^*$, and
 $u \red{*}{\lr}{}{T} v$.
Then $(\{s,t\} , T)$ has  an equivalent, convergent,
 finite presentation $(\{s,t\} , T')$, namely
 $T' = \{ (u,v) \mid u-v \in G \}$, since the reduction $\red{}{\myr}{m}{}$
 in $\myk[\{ s,t \}^*]$
 can be compared to a transformation step in a Thue system
 when restricted to polynomials of the form $u-v$.
First we show that in case a finite set $F$ has a finite Gr\"obner basis in  $\myk[\{ s,t \}^*]$ the
 procedure {\sc Gr\"obner Bases in Free Monoid Rings} also computes a
 finite Gr\"obner basis of $F$.
Let $\tilde{G}$ be a finite Gr\"obner basis of $P_{T}$ with
 $\hterm(\tilde{G})=\{\hterm(g) \mid g \in \tilde{G} \} = \{ t_1, \ldots , t_k
 \}$. 
Let $H_{t_i} = \{ xt_iy \mid x,y \in \Sigma^* \}$, then
$\hterm(\ideal{}{}(P_{T}))= \bigcup_{i=1}^k H_{t_i}$, since 
 all polynomials in $\ideal{}{\myk[\{ s,t \}^*]}(P_{T})$
 reduce to zero by $\tilde{G}$.
Further our procedure is correct and, therefore, for
 each $t_i$ there has to be at least 
 one $g_i$ added to $G$ such that $t_i \id x\hterm(g_i)y$ for some $x,y \in
 \Sigma^*$, i.e., $\hterm(g_i)$ ``divides'' $t_i$. 
Note that as soon as all such $g_i$ are added to $G$, we have
 $\hterm(\ideal{}{\myk[\{ s,t \}^*]}(G)) \supseteq \bigcup_{i=1}^k H_{t_i}$
 and all further computed
  s-polynomials must reduce to zero.
Since the procedure is correct, $G$ then is also a Gr\"obner basis of
 $\ideal{}{}(P_T)$.
It remains to show that in case $P_{T}$ has a finite Gr\"obner basis,
the finite output $G$ of our 
 procedure has the desired property that for all $g \in
 G$, $g=u-v$ where $u,v \in \{ s,t \}^*$, and 
 $u \red{*}{\lr}{}{T} v$.
Since all polynomials in $P_T$ have the desired property let us look
at the polynomials added to $G$: 
Let us assume all polynomials in $G$ have the desired structure and a
new polynomial $g$ is added. 
In case $g$ is due to s--polynomial computation of two polynomials
$u_1-v_1$,$u_2-v_2$ we do not 
 lose our structure.
The same is true for computing the normal form of a polynomial
$u-v$ using 
 a set of polynomials having the same structure.
Further $u \red{*}{\lr}{}{T} v$ is inherited within these operations
(compare also the proof of theorem 
 \ref{theo.wp}).
\\
\qed
In this proof we used a result of Madlener and Otto in \cite{Ot93} -
 a strengthening of \'{O}'D\'{u}nlaing's result 
 in \cite{OD83} to alphabets  $\Sigma_2$ containing 2 letters.
Let ${\cal P}$ be a property of  semi-Thue systems
  over $\Sigma_2$ satisfying the following three conditions:
\begin{enumerate}
\item[(P1)] Whenever $T_1$ and $T_2$ are two finite semi-Thue systems on
            the same alphabet $\Sigma_2$ such that $T_1$ and $T_2$ are
            equivalent, then $T_1$ has property ${\cal P}$ if and only
            if $T_2$ has it.
\item[(P2)] Each semi-Thue system
            $T_{\Sigma_2} = \{ a \myr \lambda \mid a \in \Sigma_2 \}$
             has property ${\cal P}$.
\item[(P3)] If a finite semi-Thue system $T$ on $\Sigma_2$ 
            has property ${\cal P}$,
            then $T$ has decidable word problem, i.e., the Thue congruence 
            $\red{*}{\lr}{}{T}$ is decidable.
\end{enumerate}
Then the following problem for ${\cal P}$ is undecidable in general:
\begin{tabbing}
XXXXXXXX\=XXXXX \kill
{\bf Given:} \> A finite semi-Thue system $T$ on $\Sigma_2$. \\
{\bf Question:} \> Does the Thue congruence $\red{*}{\lr}{}{T}$ 
                    have ${\cal P}$?
\end{tabbing}
Using this result we can easily show the following corollary.
\begin{corollary}~\\
{\sl
It is undecidable, whether for a finitely generated ideal in $\myk[\{s,t\}^*]$
 there exists a total, well-founded, admissible ordering on $\{s,t\}^*$
 such that
 the ideal has a finite Gr\"obner basis with respect to reduction as defined
 in \ref{def.red.mora}.
}
\end{corollary}
\Ba{}~\\
This follows using the correspondance between Thue systems and ideal
 bases shown in theorem \ref{theo.mora.unentscheidbar}.
Let us define a property ${\cal P}(T)$ for semi-Thue systems $T$ on
 $\Sigma_2 = \{ s,t \}$ as follows:
${\cal P}(T)$ if and only if there exists a total, well-founded,
 admissible ordering $\succeq$ on $\Sigma_2^*$ such that there exists an
 equivalent finite semi-Thue system $T'$ which is convergent with respect
 to $\succ$.
Then ${\cal P}$ fulfills the conditions (P1), (P2) and (P3) mentioned
 above:
\begin{enumerate}
\item[(P1):] If  ${\cal P}(T_1)$ holds so must ${\cal P}(T_2)$
             as
             the existence of a total, well-founded,
 admissible ordering $\succeq$ on $\Sigma_2^*$ such that there exists an
 equivalent finite semi-Thue system $T'$ which is convergent with respect
 to $\succ$ for $T_1$ at once carries
             over to the equivalent system $T_2$. 
\item[(P2):] The trivial system $\{ s \myr \lambda, t \myr \lambda \}$ has
             property ${\cal P}$.
\item[(P3):] Having property ${\cal P}$ implies decidability of the 
             Thue congruence.
\end{enumerate}
Hence this property is undecidable in general and this result carries
 over to Gr\"obner bases in $\myk[\{s,t\}^*]$ as before.
\\
\qed
Hence, for two-sided ideals the case of free monoids with two generators
is already hard
although free monoids allow simple presentations, namely empty sets of
defining relations.
\section{Skew Polynomial Rings}
Other classes of non-commutative rings and the possibilities of
introducing the theory of Gr\"obner bases to them have been studied
extensively by authors as e.g., Apel and Lassner in \cite{ApLa88},
Gateva-Ivanova in \cite{GaLa88}, Kandri-Rody and Weispfenning in
\cite{KaWe90,We92}, and Kredel in \cite{Kr93}.

A structure where finite Gr\"obner bases for one- and two-sided ideals exist is
the class of solvable polynomial rings which includes the commutative
polynomial rings, enveloping algebras of finite dimensional Lie
algebras and iterated skew polynomial rings.
Solvable polynomial rings can be described by ordinary polynomial
rings provided with a ``new'' multiplication which coincides with the
ordinary multiplication except for the case that a variable $X_j$ is
multiplied with a variable $X_i$ with lower index, i.e., $i<j$.
Then multiplication can be defined by equations
$$X_j \star X_i = c_{ij} X_iX_j + p_{ij}$$
where $c_{ij} \in \myk^*$ and $p_{ij}$ is a polynomial
``smaller'' than $X_iX_j$ with respect to a fixed admissible term
ordering on the polynomial ring.
One-sided reduction relations and one-sided Gr\"obner bases are defined
 naturally but two-sided ideals are not defined by extending the one-sided
 reduction relations. 
Instead two-sided Gr\"obner bases are characterized by one-sided
G\"obner bases, i.e., a set is a two-sided Gr\"obner basis if it is
both, a left and a right Gr\"obner basis.
The proofs for this theory are  more complicated than for
commutative polynomial rings and an extensive study of the Gr\"obner
basis approach to solvable polynomial rings can be found in
Kredel's PhD thesis (\cite{Kr93}).
Note that since we require $c_{ij} \in \myk^*$ and $p_{ij}$ is a polynomial
``smaller'' than $X_iX_j$ with respect to a fixed admissible term
ordering on the polynomial ring, we can use an admissible ordering on
the solvable polynomial ring and reduction is preserved under multiplication. 
This is no longer true if we allow $c_{ij} = 0$ as then we find that the head term
 of a multiple $p \mrm w$ need no longer be the usual commutative  product
 of the two terms $\hterm(p)$ and $w$.
Kredel gives a short discussion on such general solvable polynomial
rings, which include the Grassmann and Clifford algebras.
It is suggested to introduce a concept called saturated reduction to
 remedy the problems arising for general reduction. 

We close this chapter by giving some  details on  Weispfenning's approach to
  skew polynomial rings (compare \cite{We92}), as they
 can be regarded as monoid rings, which will be the
 subject of the next chapter.
Let $\myk[X,Y]$ be the ordinary commutative polynomial ring over $\myk$ in
two variables.
As in the case of solvable polynomial rings, a new multiplication
$\star$ is introduced, such that for some arbitrary but fixed $1 < e \in \n$ the new structure satisfies the
following axioms:
\begin{enumerate}
\item $\myk[X,Y]$ together with $0, 1, +, -$ and $\star$ forms a
  $\myk$-algebra,
\item for all $\alpha \in \myk$, $m,n \in \n$ we have \\
      $X^mY^n \star \alpha = \alpha \star X^mY^n = \alpha \skm X^mY^n$,
      \\
      $X^m \star X^n = X^{m+n}$, \\
      $Y^m \star Y^n = Y^{m+n}$, \\
      $X^m \star Y^n = X^mY^n$, and
\item $Y \star X = X^eY$.
\end{enumerate}
The skew polynomial ring corresponding to $e$ is denote by $R_e$.
Although two-sided ideals in this ring are finitely generated, this is
no longer true for one-sided ideals. 
Nevertheless, Weispfenning shows that for finitely generated one-sided
ideals  finite Gr\"obner bases exist and how they can be characterized
by s-polynomials.
Let us proceed by giving some technical details. 
The next lemma gives some insight into computation in $R_e$ that will
be used later on.
If for some terms $s,t,u$ we have $s \star t = u$, then $s$
left-divides $u$ and $t$ right-divides $u$, denoted by $s \;{\sf ldiv}\;
u$ respectively $t \;{\sf rdiv}\; u$.
Similarly, in case for a further term $s'$ we have $s \star t \star s'
= u$, then $t$ divides $u$, denoted by $t \;{\sf div}\; u$. 
\begin{lemma}~\\
{\sl
For two terms  $s = X^mY^n$ and $t = X^qY^r$ the following statements hold:
\begin{enumerate}
\item $s \;{\sf rdiv}\; t$ if and only if $n \leq r$ and $m \skm e^{r-n}
  \leq q$.  Then $X^{q - m \skm e^{(r-n)}}Y^{(r-n)} \star s = t$. 
\item $s \;{\sf ldiv}\; t$ if and only if $n \leq r$, $m \leq q$ and
  $e^n$ divides $(q-m)$.  Then $s \star X^{\frac{q-m}{e^n}}Y^{(r-n)} = t$.
\item $s \;{\sf div}\; t$ if and only if $m \leq q$ and $n \leq r$ if and
  only if there exist $h,k \in \n$ such that $X^h \star s \star Y^k =
  t$. Note that if furthermore $t = u \star s \star v$ for some terms
  $u,v$, then these terms need not be unique.
\lemend\ohnebeweis
\end{enumerate}
}
\end{lemma}
Weispfenning uses an inverse-lexicographical ordering on the terms and
proves that this ordering is admissible with respect to  the new
multiplication $\star$.
Hence in defining left, right and two-sided reduction, it turns out that
for these reductions essential properties hold, e.g., the translation
lemma, and the reduction steps are preserved under the corresponding
multiplications.
\begin{definition}~\\
{\rm
Let $p, f$ be two non-zero polynomials in $R_e$. 
We say $f$ \betonen{left reduces} $p$ at a monomial
 $\alpha \skm t$ of $p$ in one step, denoted by $p \red{}{\myr}{l}{f} q$, if
\begin{enumerate}
\item[(a)] $w \star \hterm(f) = t$ for some $w \in \m$,
            i.e., $\hterm(f) \;{\sf rdiv}\; t$, and
\item[(b)] $q = p - \alpha \skm \hc(f)^{-1} \skm w \mrm f$.
\end{enumerate}
We say $f$ \betonen{right reduces} $p$ at a monomial
 $\alpha \skm t$ of $p$ in one step, denoted by $p \red{}{\myr}{r}{f} q$, if
\begin{enumerate}
\item[(a)] $\hterm(f) \star w = t$ for some $w \in \m$,
            i.e., $\hterm(f) \;{\sf ldiv}\; t$, and
\item[(b)] $q = p - \alpha \skm \hc(f)^{-1} \skm f \mrm w$.
\end{enumerate}
We say $f$ \betonen{two-sided reduces} $p$ at a monomial
 $\alpha \skm t$ of $p$ in one step, denoted by $p \red{}{\myr}{t}{f} q$, if
\begin{enumerate}
\item[(a)] $X^m \star \hterm(f)\star Y^n = t$ for some $m,n \in \n$, and
\item[(b)] $q = p - \alpha \skm \hc(f)^{-1} \skm X^m \mrm f \mrm Y^m$.
\dend
\end{enumerate}
}
\end{definition}
Notice that the one-sided reductions correspond to the respective one-sided
 ideals.
Since in defining two-sided reduction only special multiples are allowed, this 
 correspondence to ideals no longer holds.
The set connected to two-sided reduction using $F$ is called the
 {\bf restricted ideal} generated by  $F$, which  is the closure of $F$
under addition and multiplication with powers of $X$ from the left and
powers of $Y$ from the right.
Closely related to these reductions are s-polynomials, which  as
usual arise from  special multiples of head terms of polynomials.
For two terms $t = X^mY^n$, $t'=X^{m'}Y^{n'}$ with $n \geq n'$ such
 special multiples are defined as follows:
\begin{enumerate}
\item For $r = \max \{m, m' \skm e^{(n-n')}\}$, the term $X^rY^n$ is
  the least common left multiple of $t$ and $t'$ denoted by ${\sf llcm}(t,t')$.
\item Suppose $t$ and $t'$ have some common right multiple.
      Then put $r = m$ if $m\geq m'$ and else put $r = m' + \lceil
      \frac{m'-m}{e^n} \rceil \skm e^n - m' + m$.
      Then $r \in \n$ and $X^rY^n$ is the least right common multiple
      of $t$ and $t'$ denoted by ${\sf lrcm}(t,t')$.
\end{enumerate}
Now such multiples allow to specify ``overlaps'' between ``rules'' corresponding
 to polynomials.
\begin{definition}~\\
{\rm
Let $p_1,p_2$ be two non-zero polynomials in $R_e$.
Further let $\hterm(p_i) = X^{m_i}Y^{n_i}$, $i \in \{ 1,2 \}$.
For  
 $u = {\sf llcm}(\hterm(p_1), \hterm(p_2))$ with $u_1 \star \hterm(p_1)
 = u_2 \star \hterm(p_2) = u$ we get the \betonen{left s-polynomial}\/
$$\spol{l}(p_1,p_2) = \hc(p_2) \skm u_1 \mrm p_1 - \hc(p_1) \skm u_2 \mrm
p_2.$$
In case   
 $v = {\sf lrcm}(\hterm(p_1), \hterm(p_2))$ exists and  $ \hterm(p_1) \star v_1
 =  \hterm(p_2) \star v_2 = v$ we get the \betonen{right s-polynomial}\/
$$\spol{r}(p_1,p_2) = \hc(p_2) \skm  p_1 \mrm v_1 - \hc(p_1) \skm p_2
\mrm v_2.$$
For $m = \max \{m_1,m_2 \}$ and $n = \max \{ n_1,n_2 \}$ we get the
 \betonen{two-sided s-polynomial}\/
$$\spol{t}(p_1,p_2) = \hc(p_2) \skm X^{m-m_1} \mrm p_1 \mrm Y^{n-n_1}
- \hc(p_1) \skm X^{m-m_2} \mrm p_2 \mrm Y^{n-n_2}.$$
\dend
}
\end{definition}
In analogy to Buchberger, a finite set $G \subseteq R_e$ is called a
\betonen{left, right or restricted Gr\"obner basis} if left, right
respectively two-sided reduction with respect to $G$ is confluent.
These bases can now be characterized as follows:
\begin{theorem}~\\
{\sl
A finite subset $G$ of $R_e$ is a left, right, respectively restricted Gr\"obner basis in
$R_e$ if and only if for all polynomials $f \neq g$ in $G$, $\spol{l}(f,g)
\red{*}{\myr}{l}{G} 0$, $\spol{r}(f,g)
\red{*}{\myr}{r}{G} 0$, $\spol{t}(f,g)
\red{*}{\myr}{t}{G} 0$, respectively. 
\theoend\ohnebeweis
}
\end{theorem}
Weispfenning has further shown the existence of an algorithm to
compute these bases for finite subsets of $R_e$.
Hence, the ideal membership problem for finitely generated left and
right ideals in $R_e$ is solvable.
Weispfenning further showed how restricted Gr\"obner bases can be used
to solve the ideal membership problem for two-sided ideals by constructing
two-sided Gr\"obner bases.
Two-sided Gr\"obner basis are defined as sets $G \subseteq R_e$ such
that $G$ is a restricted Gr\"obner basis and additionally the ideal
and the restricted ideal generated by $G$ coincide.
Such bases now can be characterized by the following lemma which
allows an algorithm to compute them.
\begin{lemma}~\\
{\sl
Let $G \subseteq R_e$ be a restricted Gr\"obner basis with $d = \max
\{ {\sf deg}_Y(g) \mid g \in G \}$\footnote{Here ${\sf deg}_Y$ denotes
  the number of occurrences of the variable $Y$ in the head term of $g$.}.
If for all $g \in G$ and all $0 \leq m \leq d$, $Y \mrm g
\red{*}{\myr}{t}{G} 0$ and $g \mrm X^{e^m} \red{*}{\myr}{t}{G} 0$, then
$G$ is a two-sided Gr\"obner basis.
\lemend\ohnebeweis
}
\end{lemma}
%
Interesting is that two-sided Gr\"obner bases, as in the case of
 solvable polynomial rings, are constructed using a closure of a
 specialized Gr\"obner basis, here a restricted Gr\"obner basis.
Important for termination is that the polynomials considered by the
 completion to compute left, right and restricted Gr\"obner bases have
 a bound on the $Y$-degree of the head terms of the computed
 polynomials. 

The idea of filling up a term that divides another term in order to do
reduction will be used later on in the approaches to commutative
monoids and nilpotent groups (compare chapter \ref{chapter.reduction}
and chapter \ref{chapter.grouprings}). 
In contrary to Weispfenning's approach the structures there will not
allow admissible orderings in general, but as here Dickson's lemma can
be used to prove the existence of finite bases and to ensure
termination of the algorithms. 

%% file: strongreduction.tex
\chapter{Reduction in Monoid and Group Rings}\label{chapter.reduction}
\spruch{12.3}{5}{Docendo discimus.}{Seneca}

In defining reduction for monoid and group rings different approaches
are possible and some will be studied here.
Since we mainly deal with non-commutative structures, we restrict
ourselves to right reductions.

{\bf Section 4.1:} In order to use polynomials as rules a well-founded
  ordering on the rings is necessary and this section gives
  information on how such an ordering can be lifted from the ordering
  induced by the presentation of the respective monoid.
  It is shown why well-founded orderings on the monoid in general
  need not be additionally compatible with the monoid operation.

{\bf Section 4.2:} A natural way to use a set of polynomials as rules
  is to use all term multiples of these polynomials as rules. This is what
  is done e.g. in Buchberger's approach.
  A polynomial $f$ then reduces a polynomial $p$ at a monomial $\alpha
  \skm t$ of $p$ in case there exists an element $w \in \m$ such that
  the head term of the polynomial $f \mrm w$ equals $t$.
 Then $\alpha \skm t$ can be removed from $p$ by subtracting an
 appropriate multiple $\beta \skm f \mrm w$ with  head term
 equal to $\alpha \skm t$.
 This reduction is called strong reduction as it can be used to
 characterize the right ideal congruence generated by the polynomials
 used for reduction.
 In defining s-polynomials related to strong reduction in general it turns out to
 be impossible to give a localization to finitely many candidates for
 the critical situations induced by two polynomials.
 Nevertheless, although some of the properties of Buchberger's reduction no longer
 hold,  we can characterize strong Gr\"obner bases by
 s-polynomials, but this cannot be done in a finitary manner.
 Thus this characterization does not yield a test to decide whether a
 finite set is a strong Gr\"obner basis.

{\bf Section 4.3:} Weakening reduction can provide the means to give
  a finitary confluence test.
  The first idea studied is to restrict the right multiples of a
  polynomial by terms used for reduction to those where the head term is
  preserved by multiplication, i.e., the head term of the multiple
  results from the head term of the original polynomial.
  Defining reduction in this way, the right ideal congruence is no
  longer captured by reduction. 
  This defect can be repaired by a concept called saturation, but 
  saturating sets need not be finite.
  A characterization of right Gr\"obner bases with respect to right
  s-polynomials can be given for saturated sets, but still this
  characterization is not finitary, i.e., it can neither be used to
  decide whether a finite set is a Gr\"obner basis nor does it give
  rise to a  completion procedure.

{\bf Section 4.4 and 4.5:} The next stage is to study weakenings of
  reduction involving syntactical information on the representatives
  of the monoid respectively group elements.
  We introduce the concept of prefix reduction for arbitrary monoid
  rings and the concept of commutative reduction for Abelian monoid
  rings.
  In both cases we have to use a special saturation to regain the
  expressiveness of the right ideal congruence.
  Now for prefix respectively commutatively saturated sets 
  characterizations of the respective Gr\"obner bases by special
  s-polynomials in a finitary manner are given.
  These characterizations can be used to decide whether a finite set
  is an appropriate Gr\"obner basis in case it is saturated in the
  appropriate way which again is decidable.
  A procedure to enumerate a prefix Gr\"obner basis is provided which
  halts in case a finite prefix Gr\"obner basis exists.
  In the commutative case finite Gr\"obner bases
  always can be computed.
  Interreduction is introduced to both settings and the existence of
  unique monic reduced Gr\"obner bases with respect to the respective
  reductions is shown.

\section{Using Polynomials as Rules}
In order to define an effective reduction in monoid and
 group rings we have to ensure
 that certain algebraic operations can be done effectively and certain
 algebraic questions can be solved in our structure.
Therefore, we will have to restrict the presentations allowed for the
 monoids and groups in 
 the approach developed here.
First of all we assume that our monoids are presented by finite convergent
 reduction
 systems (compare section \ref{section.presentations}).
This implies that our monoids have solvable word problem and hence we can decide whether
 two elements are equal which is essential in performing ring
 operations effectively.
When introducing polynomial reduction and later on s-polynomials to our ring, we will find that
 it is important to solve the following two algebraic questions in the
 monoid $\m$:
\begin{enumerate}
\item Given $w_1,w_2 \in \m$, is there an element $m \in \m$ such that $w_1 = w_2 \mm m$?
\item Given $w_1,w_2 \in \m$, are there  elements $m_1,m_2  \in \m$ such that $w_1 \mm m_1 = w_2 \mm m_2$?
\end{enumerate}
Note that these questions in general are undecidable even for monoids
 with convergent presentations.
However, for groups the answer of course is ``yes'' in both cases, as
 the element $m = \inv{w_2} \mm w_1$ is a unique solution to the first
 question and the set $\{ (\inv{w_1} \mm w, \inv{w_2} \mm w) | w \in
 \m \}$ contains all solutions to the second question.

Since the objects of interest will be polynomials and we want to use
 them as rules in a reduction system, we have to introduce an ordering
 on the monoid elements.
This ordering is required to be total and well-founded, but we will
 find that in general we cannot expect it to be admissible as in
 Buchberger's approach and most of the extensions of his ideas to other
 structures.
Because we will assume that the monoids are presented by reduction systems,
 monoid elements can be viewed as syntactical objects over an alphabet.
Let us be more specific now for the case that our presentations are 
 finite convergent semi-Thue systems having a total admissible well-founded completion ordering,
 where the above questions 1. and 2. are solvable and the solutions computable.
If not stated otherwise we will assume our monoids to be presented in
 this manner.
Further we will assume that $\succeq_T$ is a total well-founded
  admissible ordering on $\Sigma^*$ such that the presentation $(\Sigma, T)$
 is convergent with respect to $\succeq_T$.
Then we will always assume that the  well-founded total 
 ordering $\succeq$ 
 on the monoid $\m$ is the restriction of the ordering $\succeq_T$ to
 the irreducible representatives of the monoid elements.
In particular this gives us  $m \succ \lambda$ for all $m \in \m
 \backslash \{\lambda\}$.
Moreover, since we identify the elements of $\m$ with the words in $\irr(T)$
 an essential conclusion used throughout will be that for $u,v \in \m$,
 $uv \succeq u \mm v$ holds, where $uv$ stands for the concatenation
 of $u$ and $v$.
This follows immediately as $\succeq_T$ is admissible on $\Sigma^*$.
Similar properties hold if we use convergent semi-Thue
 systems modulo commutativity to present commutative monoids.

In order to define  reduction in $\myk[\m]$ we want to use polynomials
 as rules.
This can be done by using the well-founded ordering on the monoid to
 give us an ordering on the monomials of a polynomial.
\begin{definition}~\\
{\rm
Let $\succeq$ denote a well-founded total ordering on $\m$.
\begin{enumerate}
\item Let $p \in \myk[\m]\backslash \{ 0 \}$ be denoted by the polynomial 
           $p = \sum_{i = 1}^{n} \alpha_{i} \skm w_{i}$,
           where $\alpha_i \in \myk^*$, $w_i \in \m$ and
             $w_{i} \neq w_{j}$ for $i \neq j$.
           Furthermore, we assume that according to our ordering $\succeq$ we have
            $w_{1} \succ \ldots \succ w_{n}$.
           Then we let $\hm(p) = \alpha_{1} \skm w_{1}$ denote the
            \index{head monomial!in a monoid ring}\betonen{head monomial},
            $\hterm(p) = w_{1}$ the \index{head term!in a monoid
              ring}\betonen{head term}\/ 
            and $\hc(p) = \alpha_{1}$ the \index{head coefficient!in a monoid ring}\betonen{head 
            coefficient}\/ of $p$.
           $\reductum(p) = p - \hm(p)$ stands for the 
            \index{reduct!in a monoid ring}\betonen{reduct}\/ of $p$.
           $\terms(p)= \{ w_1, \ldots, w_n \}$ is the set of
            \index{term!in a monoid ring}\betonen{terms}\/ occurring
            in $p$.
           The polynomial $p$ is called \betonen{monic}\index{monic} in case
            $\hc(p)= 1$.
\item For a set of polynomials $F$ in $\myk[\m]$ we define
      $\hterm(F)  =  \{ \hterm(f) | f \in F \}.$
\dend 
\end{enumerate}
}
\end{definition}
Moreover, we can extend the well-founded total ordering on $\m$
 to an  ordering  on the elements of $\myk[\m]$, which is again well-founded.
\begin{definition}\label{def.order}~\\
{\rm
 Let $p, q$ be two polynomials in $\myk[\m]$. 
 Then we say $p$ is \index{greater!on a monoid ring}\index{ordering!on a monoid ring}\betonen{greater}\/
  than $q$ with respect to an ordering $\succeq$ on $\m$, i.e., $p  > q$, if
  \begin{enumerate}
  \item[(i)]   $\hterm(p) \succ \hterm(q)$ or
  \item[(ii)]  $\hm(p) = \hm(q)$ and $\reductum(p)  > \reductum(q)$.
\dend
  \end{enumerate} 
}
\end{definition}
Note that this ordering is not  total on $\myk[\m]$, e.g.\  
 $5 \skm w$ and $3 \skm w$ are incomparable.
\begin{lemma}\label{lem.well-founded}~\\
{\rm
The ordering $\geq$ on $\myk[\m]$ as given in definition \ref{def.order} is
well-founded.
\lemend
}
\end{lemma}
\Ba{}~\\
The proof of this lemma will use a method that is known as Cantor's
second diagonal argument (compare e.g. \cite{BeWe92} chapter 4).
Let us assume that $\geq$ is not well-founded on $\myk[\m]$.
We will show that this gives us a contradiction to the fact that the
ordering $\succeq$ on $\m$ inducing $\geq$ is well-founded.
Hence, let us suppose $f_0 > f_1 > \ldots > f_k > \ldots\;$, $k \in \n$ is a strictly descending chain in
$\myk[\m]$.
Then we can construct a sequence of sets of pairs $\{ \{ (t_k, g_{kn}) | n
\in \n \} | k \in \n \}$  recursively as follows:
For $k=0$ let $t_0 = \min \{ \hterm(f_i) | i \in \n \}$\footnote{Note
  that this minimum exists as we have a subset
  of $\m$ and $\succeq$ is well-founded on $\m$.}.
Now let $j \in \n$ be the least index such that we have $t_0 =\hterm(f_j)$.
Then $t_0 = \hterm(f_{j + n})$ holds for all $n \in \n$ and we can set
$g_{0n} = f_{j+n} - \hm(f_{j+n})$, i.e., $\hterm(g_{0n}) \pred t_0$
for all $n \in \n$.
For $k+1$ we let $t_{k+1} = \min \{  \hterm(g_{ki}) | i \in \n \}$ and
again let $j \in \n$ be the least index such that $t_{k+1} =
\hterm(g_{kj})$ holds, i.e., $t_{k+1} = \hterm(g_{k(j+n)})$ for all $n \in \n$.
Again we set $g_{(k+1)n} = g_{k(n+j)} - \hm(g_{k(n+j)})$.
\\
Then the following statements hold:
\begin{enumerate}
\item For all $s \in \terms(g_{kn})$ we have $s \prec t_k$.
\item For every $k \in \n$, $g_{k0} > g_{k1} > \ldots\; $ is  a
strictly descending chain in $\myk[\m]$.
\end{enumerate}
Hence we get that $t_0 \succ t_1 \succ \ldots \;$ is a strictly
descending chain in $\m$ contradicting the
fact that $\succeq$ is supposed to be well-founded on $\m$.
\\
\qed
The choice of the well-founded ordering on the reduction ring
 is of great importance for the characteristics
 of reduction in a structure.
As we have seen in the survey on Buchberger's approach,
 in the case of polynomial rings
 (which correspond to free commutative monoid rings)
 a special class of orderings,
 namely term orderings, can be used.
These orderings have the following useful property which is closely
related to the fact that reduction is preserved under multiplication.
\begin{definition}~\\
{\rm
An ordering $\succeq$ on a monoid $\m$ is called 
 \index{compatible ordering}\index{ordering!compatible}\index{ordering!admissible}\index{ordering!monotone} \index{admissible ordering}\index{monotone ordering}\betonen{monotone}\/
 (\betonen{compatible with
  multiplication $\mm$}\/) if
 $u \succeq v$ implies $u \mm w \succeq v \mm w$
  for all $u,v,w \in \m$.
It is called \betonen{admissible} in case we additionally have $w
\succeq \lambda$ for all $w \in \m$.
\dend
}
\end{definition}
This property is essential in Buchberger's approach, since it ensures that 
 $p \red{*}{\myr}{}{F} 0$ implies $\alpha \skm p \mrm w \red{*}{\myr}{}{F}
 0$ for all $p \in \myk[X_1, \ldots, X_n]$, $\alpha \in \myk$ and $w \in
 {\cal T}$.
Unfortunately  in general a monotone total ordering $\succeq$
 on a monoid $\m$ cannot be expected to be well-founded, as the following remark shows.
\begin{remark}~\\
{\rm
  Let $\m \neq \{ \lambda \}$ be a monoid
   with a monotone total ordering $\succeq$.
\begin{enumerate}
\item Then $\m$  
       cannot contain a nontrivial element of finite order.\\
      To see this, suppose $w \in\m \backslash \{ \lambda \}$ is of finite order,
       i.e., there are $n,m \in\n$, $n > m$ such that $w^n = w^m$.
      Without loss of generality let us assume $w \succ \lambda$.
      In case $n = m+1$,
       then (as $\succeq$ is monotone and transitive) we get
       $ w^{m} \succ w^{m-1}$ giving us 
       $w^{m} = w^{n} = w^{m+1} \succ w^{m}$,
       contradicting  our assumption.
      Otherwise we get $w^{n-1} \succ \ldots \succ w^{m+1} \succ
       w^{m}$ likewise giving us $w^{m} = w^{n} \succ w^{m}$ which
       is again a contradiction.
\item The ordering $\succeq$ cannot be well-founded, in case
       there are two elements of infinite order
       $u, v \in \m \backslash \{ \lambda \}$ satisfying
       $u \mm v = \lambda$.
       \\
      Without loss of generality let us assume $u \succ \lambda$.
      Then (as $\succeq$ is monotone) we have $\lambda \succ v$,
       and (as $\succeq$ is transitive)
       $u \succ \lambda \succ v \succ \ldots \succ v^{n}$ for
       all $n \in\n$ gives us an infinite
       descending chain of elements in $\m$.
\remend
\end{enumerate}
}
\end{remark}
Especially non-trivial groups do not allow monotone well-founded total
 orderings.

Now we can move on to discuss reduction in monoid rings.
We will see that
 in defining appropriate reductions we have to be  more cautious than in defining
 reductions in the polynomial ring 
 (compare section \ref{section.buchberger}) where simply the
 head of a polynomial is modified in order to use the polynomial for reduction.
Let us start by  examining a rather natural approach to reduction.
Since we are mainly interested in non-commutative structures, we will
 restrict ourselves to one-sided ideals and investigate
 right ideals and concepts for reduction using right multiplication by monomials only.
\section{The Concept of Strong Reduction}\label{section.reduction}
%
In order to study right ideals generated by a set of polynomials it is
often useful to take a look at special representations of the elements
in this set.
We will start with a first description here and refine this approach
within the following sections.
Such representations of right ideal elements can then be connected to
different definitions of reduction in a monoid ring.
Henceforth, let $\succeq$ denote a well-founded total ordering on $\m$.
\begin{definition}\label{def.sr}~\\
{\rm
Let $F$ be a set of polynomials  and $p$ a non-zero polynomial in $\myk[\m]$.
A representation 
$$ p = \sum_{i=1}^{n} \alpha_i \skm f_{i} \mrm w_i, \mbox{ with }
 \alpha_i \in \myk^*,f_{i} \in F, w_i \in \m $$
is called a \index{standard representation}\betonen{standard representation}\/ of $p$ with respect to $F$,
if for all $1 \leq i \leq n$ we have $\hterm(p) \succeq \hterm(f_{i} \mrm
w_i)$.
\dend
}
\end{definition}
For the reader familiar with the framework of standard and Gr\"obner
 bases, we mention that this is a natural adaption of the term 
 ``standard representation'', as e.g.\  defined in \cite{BeWe92} for
 the polynomial ring, to the case of right ideals in monoid rings.
A standard representation of a polynomial $p$ with respect to a set of
polynomials $F$ is thus a representation where all occurring terms
involved are bounded by the head term of $p$.
Note that for at least one index $1 \leq i \leq n$ we must have $\hterm(p)
= \hterm(f_{i} \mrm w_i)$ and $\hc(p) = \sum_{i,\hterm(p) = \hterm(f_{i} \mrm
w_i)} \alpha_i \skm \hc(f_{i} \mrm w_i)$.
Standard representations can be used to characterize special bases of right ideals.
\begin{definition}~\\
{\rm
A set $F \subseteq \myk[\m]$ is called a 
 \index{right!standard basis}\index{standard basis!right}\index{standard basis}\betonen{(right) standard
   basis}, if every non-zero polynomial in
 $\ideal{r}{}(F)$ has a standard representation with
 respect to $F$.
\dend
}
\end{definition}
One way to characterize such bases arises from introducing
reduction to $\myk[\m]$ and a rather natural approach  is to use a right multiple of a
polynomial as a rule.
\begin{definition}\label{def.reds}~\\
{\rm
Let $p, f$ be two non-zero polynomials in $\myk[\m]$. 
We say  $f$ 
 \index{strong!right reduction}\index{reduction!strong right}\betonen{strongly right reduces} $p$ to $q$ at 
 a monomial $\alpha \skm t$ of $p$ in one step, denoted by 
 $p \red{}{\myr}{s}{f} q$, if
\begin{enumerate}
\item[(a)] $\hterm(f \mrm w) = t$ for some $w \in \m$, and
\item[(b)] $q = p - \alpha \skm \hc(f \mrm w)^{-1} \skm f \mrm w$.
\end{enumerate}
We write $p \red{}{\myr}{s}{f}$ if there is a polynomial $q$ as defined
above and $p$ is then called strongly right reducible by $f$. 
Further, we can define $\red{*}{\myr}{s}{}, \red{+}{\myr}{s}{}$ and
 $\red{n}{\myr}{s}{}$ as usual.
Strong right reduction by a set $F \subseteq \myk[\m]$ is denoted by
 $p \red{}{\myr}{s}{F} q$ and abbreviates $p \red{}{\myr}{s}{f} q$
 for some $f \in F$,
 which is also written as  $p \red{}{\myr}{s}{f \in F} q$.
\dend
}
\end{definition}
Note that in order to strongly right reduce $p$, the polynomial $f$
 need not be smaller than $p$. 
The condition $\hterm(f \mrm w) = t$
 prevents reduction with a polynomial in case $f \mrm w = 0$, i.e., if
 the monomials of $f$ eliminate each other by multiplying $f$ with $w$.
This might happen in case the monoid ring contains zero-divisors 
 (compare remark \ref{rem.zero-divisors}).
Further, in case we have 
 $p \red{}{\myr}{s}{f} q$ at the monomial $\alpha \skm t$, then $t
 \not\in \terms(q)$.
\begin{definition}~\\
{\rm
A set of polynomials $F$ is called
\betonen{interreduced}\index{interreduced} or
\betonen{reduced}\index{reduced} with respect to
$\red{}{\myr}{s}{}$ if for all $f \in F$,$f \red{}{\myr}{s}{F} f'$
 implies $f' = 0$.
\dend
}
\end{definition}
Notice that in literature an interreduced set is often characterized
 by  requiring that none of its elements is reducible by the other
 elements in the set.
In our setting this no longer holds, since strongly right reducing
 a polynomial with {\em itself} need not result in zero as the following
 example shows.
\begin{example}~\\
{\rm
Let $\Sigma = \{ a, b, c \}$ and
       $T = \{ a^{2} \myr \lambda, b^{2} \myr \lambda, ab \myr c, ac \myr b, cb \myr a \}$
       be a presentation of a monoid $\m$ (which is in fact a group)
       with a length-lexicographical ordering induced by $a \succ b
       \succ c$.
\\
Then the polynomial $a+b+c$ is strongly reducible by itself at $a$ as
follows:
$a+b+c \red{}{\myr}{s}{a+b+c} a+b+c - (a+b+c) \mrm b = a+b+c - (c +
\lambda + \underline{a}) = b - \lambda$. \\
Moreover, since $\ideal{r}{}(a+b+c) \neq \ideal{r}{}(b - \lambda)$,
 this shows that  while $f \red{}{\myr}{s}{f} f'$ implies $f' \in
 \ideal{r}{}(f)$ the case $f \not\in \ideal{r}{}(f')$ is possible.
\exaend
}
\end{example}
This example also reveals that interreduced bases of right ideals
 in general need not exist: neither $\{a+b+c\}$ nor $\{a+b+c, b-\lambda\}$
 are interreduced, and while $\{ b - \lambda \}$ is interreduced it
  no longer generates $\ideal{r}{}(a+b+c)$.

In order to decide, whether a polynomial $f$ strongly right reduces a
polynomial $p$ at a monomial $\alpha \skm t$ one has to decide whether
there exist elements  $s \in \terms(p)$ and $w \in \m$ such that $s
\mm w = \hterm(f \mrm w) = t$.
Since this problem reduces to solving equations $s \mm {\rm x} = t$
 in one variable ${\rm x}$ in the monoid $\m$ presented by
 $(\Sigma , T)$, this problem is undecidable in general, even if $\m$
 is presented by a convergent semi-Thue-system.
Note that there can be no, one or even (infinitely) many solutions
depending on $\m$.
\begin{example}~\\
{\rm
Let $\Sigma = \{ a,b \}$ and $T = \{ ab \myr a \}$ be a
 presentation of a monoid $\m$ with a length-lexicographical ordering
 induced by $a \succ b$.
\\
Then the equation $b \mm {\rm x} = a$ has no solution in $\m$,
the equation $b \mm {\rm x} = b$ has one solution in $\m$, namely
${\rm x} = \lambda$, and
the equation $a \mm {\rm x} = a$ has infinitely many solutions in
$\m$, namely the set $\{ b^n | n \in \n \}$.
\exaend
}
\end{example}
The following example  illustrates how different monomials
can become equal when modifying a polynomial in order to use it for
strong right reduction.
\begin{remark}~\\
{\rm
Let $\Sigma = \{ a,b \}$ and $T = \{ ab \myr b \}$ be a
 presentation of a monoid $\m$ with a length-lexicographical ordering
 induced by $a \succ b$.
Furthermore, let $f_1,
f_2, p$ be polynomials in $\q[\m]$ such that $f_1 = a^2 + a$, $f_2 = a^2
- a$ and $p = b + \lambda$.
\\
Then $p$ is strongly right reducible by $f_1$ at $b$, as 
 $\hterm(f_1 \mrm b) = \hterm(2 \skm b) = b$ and 
 $p \red{}{\myr}{s}{f_1} p - \frac{1}{2} \skm f_1 \mrm b = b+\lambda - \frac{1}{2} \skm 2
 \skm b = \lambda$.
On the other hand, although both equations $a^2 \mm
{\rm x} = b$ and $a \mm {\rm x} =b$ have a solution $b$, we get 
 that $p$ is not strongly right reducible by $f_2$, as $f_2 \mrm b = b -  b = 0$.
\remend
}
\end{remark}
In case $\m$ is a right cancellative monoid or a group the phenomenon
 described in this remark can no longer occur, since 
 then $u \mm m = v \mm m$ implies $u = v$ for all $u,w,m \in \m$.
Let us continue to state some  of the properties strong right
 reduction satisfies.
\begin{lemma}\label{lem.red}~\\ 
{\sl 
Let $F$ be a set of polynomials and $p,q, q_1, q_2$ some
 polynomials in $\myk[\m]$.
\begin{enumerate}
\item $p \red{}{\myr}{s}{F} q$ implies $p > q$, in particular $\hterm(p)
  \succeq \hterm(q)$.
\item $\red{}{\myr}{s}{F}$ is Noetherian.
\item $p \red{}{\myr}{s}{q_1} 0$ and $q_1 \red{}{\myr}{s}{q_2} 0$ imply $p \red{}{\myr}{s}{q_2} 0$.
\item\label{lem.red.4} $\alpha \skm p \mrm w \red{{\leq 1}}{\myr}{s}{p} 0$ for
  all $\alpha \in \myk^*$, $w \in \m$.
\lemend
\end{enumerate}
}
\end{lemma}
\Ba{}
\begin{enumerate}
\item This follows from the fact that using a polynomial $f$ together
       with $\alpha \in \myk^*$ and $w \in \m$ for reduction
       we use $\alpha \skm \hm(f \mrm w) \myr -\alpha  \skm \reductum(f \mrm w)$ as a rule
       and we know $\hm(f \mrm w) > -\reductum(f \mrm w)$.
\item This follows from (1), as the ordering $\geq$
       on $\myk[\m]$ is  well-founded.
\item $p \red{}{\myr}{s}{q_1} 0$ implies $p = \alpha_1  \skm q_1 \mrm w_1$
       for some $\alpha_1 \in \myk^*, w_1 \in \m$, and
       $q_1 \red{}{\myr}{s}{q_2} 0$ implies $q_1 = \alpha_2 \skm  q_2 \mrm
       w_2$ for some  $\alpha_2 \in \myk^*, w_2 \in \m$. 
      Combining this information we immediately get  $p \red{}{\myr}{s}{q_2} 0$, as
        $p = \alpha_1 \skm q_1 \mrm w_1
           = \alpha_1 \skm  (\alpha_2 \skm  q_2 \mrm w_2) \mrm w_1
           = (\alpha_1 \skm  \alpha_2) \skm  q_2 \mrm ( w_2 \mm w_1)$
        and thus $\hterm(q_2 \mrm (w_2 \mm w_1))=\hterm(p)$.
\item This follows immediately from definition \ref{def.reds}. 
\\
\qed
\end{enumerate}\renewcommand{\baselinestretch}{1}\small\normalsize
However, a closer inspection of strong right reduction reveals some
dependencies on the monoid and some restrictions in general.
\begin{remark}~\\
{\rm
Let $p,q,q_1$ and $q_2$ be some polynomials in $\myk[\m]$.
\begin{enumerate}
\item $p \red{}{\myr}{s}{q_1} 0$ and $p \red{}{\myr}{s}{q_2} 0$ need
      not imply $q_1 \red{}{\myr}{s}{q_2}$ or $q_2 \red{}{\myr}{s}{q_1}$.
      \\
      Let us recall that $p \red{}{\myr}{s}{q_1} 0$ and
       $p \red{}{\myr}{s}{q_2} 0$ imply
       $p = \alpha_1 \skm q_1 \mrm w_1 =\alpha_2 \skm q_2 \mrm w_2$ for some
       $\alpha_1,\alpha_2 \in \myk^*$ and $w_1,w_2 \in \m$.
      Obviously, if $\m$ is a group we get 
       $q_1 = (\alpha_1^{-1} \skm \alpha_2) \skm q_2 \mrm (w_2 \mm \inv{w_1})$ and
       $q_2 = (\alpha_2^{-1} \skm \alpha_1) \skm q_1 \mrm (w_1 \mm \inv{w_2})$,
       i.e., $q_1 \red{}{\myr}{s}{q_2} 0$ and $q_2 \red{}{\myr}{s}{q_1}
       0$.
      But for arbitrary monoids this does not hold, as it is already
       wrong for Buchberger's reduction in the usual 
       polynomial ring\footnote{Take e.g.\  the polynomials $p =
         X_1X_2$, $q_1 = X_1$ and $q_2 = X_2$.}.
\item $p \red{}{\myr}{s}{q}$ and
       $q \red{}{\myr}{s}{q_1} q_2$ need not imply  $p \red{}{\myr}{s}{\{q_1,q_2\}}$. 
      \\
      Let  $\Sigma = \{ a, b, c \}$ and
       $T = \{ a^2 \myr \lambda, b^2 \myr \lambda, c^2 \myr \lambda \}$  be a
       presentation of a monoid $\m$ (which is in fact a group), 
       with with a length-lexicographical ordering
       induced by $a \succ b \succ c$.
      \\
      Looking at the polynomials $p = ba + b, q = bc + \lambda$
       and $q_1 = ac + b$ we find
       $p \red{}{\myr}{s}{q} p - q \mrm ca = ba + b - (bc + \lambda)
       \mrm ca = ba + b - ba - ca = -ca + b$ and
       $q \red{}{\myr}{s}{q_1} q - q_1 \mrm c = bc + \lambda - (ac + b)
       \mrm c = bc + \lambda - a - bc = -a + \lambda =: q_{2}$,
       but $p \nred{}{\myr}{s}{\{q_1,q_2\}}$, as
       trying to reduce  $ba$ by $q_1$ or $q_{2}$ we get
       $q_1 \mrm a = \underline{aca} + ba, q_1 \mrm caba = ba + \underline{bcaba}$ and $q_{2} \mrm
       aba = -ba + \underline{aba}, q_{2} \mrm ba = -\underline{aba} + ba$
       all violating\footnote{The underlined terms in the polynomial
         multiples are the respective head terms.}
       condition (a) of definition \ref{def.reds}, i.e., there exists
       no $w \in \m$ such that $\hterm(q_1 \mrm w) = ba$ or
       $\hterm(q_2 \mrm w) = ba$.
      Trying to reduce $b$ we get the same problem with
       $q_1 \mrm cab = b + \underline{bcab}, q_2 \mrm ab = -b+\underline{ab}$
       and $q_2 \mrm b = -\underline{ab}+b$.
\remend
\end{enumerate}
}
\end{remark}
Note that the latter property is connected to interreduction and this
example states that in case a polynomial $p$ is strongly right
reducible by a polynomial $q$ and the latter polynomial is reduced to
a new polynomial $q_2$ by a polynomial $q_1$, then $p$ need no longer
be strongly right reducible by the set $\{ q_1, q_2 \}$.
Hence, interreducing a set of polynomials  might affect the set of
polynomials which have been strongly right reducible by the not
interreduced set.
Nevertheless, strong right reduction has the essential properties
which allow us to characterize a right ideal by reduction with respect to a
a set of generators, e.g. the translation lemma holds and the right
ideal congruence can be described by reduction.
\begin{lemma}\label{lem.confluent}~\\
{\sl
Let $F$ be a set of polynomials and $p,q,h$ some
 polynomials in $\myk[\m]$.
\begin{enumerate}
\item
Let $p-q \red{}{\myr}{s}{F} h$.
Then there are  polynomials $p',q' \in \myk[\m]$ such that we have 
 $p  \red{*}{\myr}{s}{F} p', q  \red{*}{\myr}{s}{F} q'$ and $h=p'-q'$.
\item
Let $0$ be a normal form of $p-q$ with respect to $\red{}{\myr}{s}{F}$.
Then there exists a polynomial  $g \in \myk[\m]$ such that
 $p  \red{*}{\myr}{s}{F} g$ and $q  \red{*}{\myr}{s}{F} g$.
\lemend
\end{enumerate}
}
\end{lemma}
\Ba{}
\begin{enumerate}
\item  Let $p-q \red{}{\myr}{s}{F} h = p-q-\alpha \skm f \mrm w$
        with $\alpha \in \myk^*, f \in F, w \in \m$
        and let $\hterm(f \mrm w) = t$, i.e.,
        $\alpha \skm \hc(f \mrm w)$ is the coefficient of $t$ in $p-q$.
       We have to distinguish three cases:
       \begin{enumerate}
         \item $t \in \terms(p)$ and $t \in \terms(q)$:
               Then we can eliminate the term $t$ in the polynomials 
                $p$ respectively $q$ by
                reduction and 
                get $p \red{}{\myr}{s}{f} p - \alpha_1 \skm f \mrm w= p'$,
                $q \red{}{\myr}{s}{f} q - \alpha_2 \skm f \mrm w= q'$,
                with $\alpha_1  -  \alpha_2 =\alpha$,
                where $\alpha_1 \skm \hc(f \mrm w)$ and 
                $\alpha_2 \skm \hc(f \mrm w)$ are
                the coefficients of $t$ in $p$ respectively $q$.
         \item $t \in \terms(p)$ and $t \not\in \terms(q)$:
               Then  we can eliminate the term $t$ in the polynomial 
                $p$  by reduction and  get
                $p \red{}{\myr}{s}{f} p - \alpha \skm f \mrm w= p'$ 
                and $q = q'$.
         \item $t \in \terms(q)$ and $t \not\in \terms(p)$:
               Then  we can eliminate the term $t$ in the polynomial 
                $q$ by reduction and  get
                $q \red{}{\myr}{s}{f} q + \alpha \skm f \mrm w= q'$
                and $p = p'$.
       \end{enumerate}
      In all three cases we have $p' -q' =  p - q - \alpha \skm f \mrm w = h$.
\item We show our claim by induction on $k$, where $p-q \red{k}{\myr}{s}{F} 0$.
      In the base case $k=0$ there is nothing to show.
      Hence, let $p-q \red{}{\myr}{s}{F} h  \red{k}{\myr}{s}{F} 0$.
      Then by (1) there are polynomials $p',q' \in \myk[\m]$ such that 
       $p  \red{*}{\myr}{s}{F} p', q  \red{*}{\myr}{s}{F} q'$ and $h=p'-q'$.
      Now the induction hypothesis for $p'-q' \red{k}{\myr}{s}{F} 0$  yields 
       the existence of a polynomial $g \in \myk[\m]$ such that
       $p  \red{*}{\myr}{s}{F} p' \red{*}{\myr}{s}{F} g$ and
       $q  \red{*}{\myr}{s}{F} q' \red{*}{\myr}{s}{F} g$.
\\
\qed
\end{enumerate}\renewcommand{\baselinestretch}{1}\small\normalsize
\begin{lemma}\label{lem.strong.congruence}~\\
{\sl
Let $F$ be a set of polynomials and $p,q$ some
 polynomials in $\myk[\m]$.
Then
 $$p \red{*}{\lr}{s}{F} q \mbox{ if and only if } p - q \in
 \ideal{r}{}(F).$$
\lemend
}
\end{lemma}
\Ba{}
 \begin{enumerate}
      \item Using induction on $k$ we show that
            $p \red{k}{\lr}{s}{F} q$ implies
            $p - q \in \ideal{r}{}(F)$.
      In the base case $k = 0$ there is nothing to show,
       since $p-p =0 \in \ideal{r}{}(F)$.
      Thus let us assume that  $\tilde{p}\red{k}{\lr}{s}{F}\tilde{q}$
       implies $\tilde{p}-\tilde{q}\in \ideal{r}{}(F)$.
      Then looking at $p \red{k}{\lr}{s}{F} p_{k}  \red{}{\lr}{s}{F} q$
       we can distinguish
       two cases:
        \begin{enumerate}
        \item $p_{k}  \red{}{\myr}{s}{f} q$ using a polynomial $f \in F$.
            \\
            This gives us $q = p_k - \alpha  \skm f \mrm w$,
             where $\alpha \in \myk^*,w \in \m$, and since 
             $p-q = p - p_k + \alpha \skm f \mrm w$ and $p - p_k \in
             \ideal{r}{}(F)$, we get $p-q \in \ideal{r}{}(F)$.
        \item $q  \red{}{\myr}{s}{f} p_k$ using a polynomial
               $f \in F$ can be treated similarly.
        \end{enumerate}
      \item It remains to show that $p-q \in \ideal{r}{}(F)$ implies 
             $p \red{*}{\lr}{s}{F} q$.
            Remember that $p-q \in \ideal{r}{}(F)$ gives us a representation
             $p = q + \sum_{j=1}^{m} \alpha_{j} \skm f_j \mrm w_{j}$
             such that
             $\alpha_{j} \in \myk^*, f_j \in F$, and  $w_{j} \in \m$.
            We will show  $p \red{*}{\lr}{s}{F} q$ by induction on $m$.
            In the base case $m = 0$ there is nothing to show.
            Hence, let  $p = q + \sum_{j=1}^{m} \alpha_{j} \skm f_j \mrm w_{j} +
             \alpha_{m+1} \skm f_{m+1} \mrm w_{m+1}$
             and by our induction hypothesis
             $p \red{*}{\lr}{s}{F} q + \alpha_{m+1} \skm f_{m+1} \mrm w_{m+1}$.
In showing $q + \alpha_{m+1} \skm f_{m+1} \mrm w_{m+1} \red{*}{\lr}{s}{F} q$ we are done.
Notice that $(q + \alpha_{m+1} \skm f_{m+1} \mrm w_{m+1}) - q = \alpha_{m+1} \skm f_{m+1} \mrm w_{m+1} 
 \red{\leq 1}{\myr}{s}{f_{m+1}} 0$ and hence lemma \ref{lem.confluent} implies our claim.
\\
\qed
\end{enumerate}\renewcommand{\baselinestretch}{1}\small\normalsize
In analogy to Buchberger's definition of Gr\"obner bases in
commutative polynomial rings, we can now specify special bases of right
ideals in monoid rings. 
%
\begin{definition}\label{def.gb}~\\
{\rm
A  set $G \subseteq \myk[\m]$ is called a 
 \index{Gr\"obner basis!strong}\index{strong!Gr\"obner basis}\betonen{Gr\"obner basis}\/
 with respect to
 the reduction $\red{}{\myr}{s}{}$ or a \betonen{strong Gr\"obner basis}, if
\begin{enumerate}
\item[(i)] $\red{*}{\lr}{s}{G} = \;\; \equiv_{\ideal{r}{}(G)}$, and
\item[(ii)] $\red{}{\myr}{s}{G}$ is confluent.
\dend
\end{enumerate}
}
\end{definition}
Note that unlike in Buchberger's case a polynomial itself need not be a  Gr\"obner
 basis of the right ideal it generates.
\begin{example}\label{exa.no.strong.gb}~\\
{\rm
Let $\Sigma = \{ a, b, c \}$ and
 $T = \{ a^{2} \myr \lambda, b^{2} \myr \lambda, ab \myr c, ac \myr b, cb \myr a
 \}$ be a presentation of a monoid $\m$ (which is in fact is a
 group) with with a length-lexicographical ordering
 induced by $a \succ b \succ c$.
Further, let us consider the polynomial   $p = a + b + c \in \q[\m]$.
\\
Then 
  $\red{}{\myr}{s}{p}$ is not  confluent on $\ideal{r}{}(p)$, as we can
 reduce
 $a+b+c  \red{}{\myr}{s}{p} b- \lambda$ using $p \mrm b = c+\lambda + \underline{a}$ and
 $a+b+c  \red{}{\myr}{s}{p} 0$,
 but although $b - \lambda \in \ideal{r}{}(p)$, $b- \lambda \nred{*}{\myr}{s}{p} 0$, as for all
 $w \in \m$, $\hterm(p \mrm w) \neq b$.
\exaend
}
\end{example}
The following lemma collects some natural relations between strong right
reduction and standard representations.
\begin{lemma}\label{lem.srprop}~\\
{\sl
Let $F$ be a set of polynomials and
 $p$ a non-zero polynomial in $\myk[\m]$.
\begin{enumerate}
\item Then $p \red{*}{\myr}{s}{F} 0$ implies the existence of a standard
  representation for $p$.
\item In case the polynomial $p$ has a standard representation with respect to $F$,
  then $p$ is strongly right reducible at its head monomial by $F$,
  i.e., $p$ is top-reducible by $F$.
\item\label{lem.srprop.3} In case $F$ is a standard basis, every non-zero polynomial
  $p$ in $\ideal{r}{}(F)$ is top-reducible to zero by $F$.
\lemend
\end{enumerate}
}
\end{lemma}
\Ba{}
\begin{enumerate}
\item This follows directly by adding up the polynomials used in the
  strong right reduction steps  occurring in $p \red{*}{\myr}{s}{F} 0$.
\item This is an immediate consequence of the definition of standard
  representations in \ref{def.sr}
  as the existence of a polynomial $f$ in $F$ and an element
  $w\in\m$ with $\hterm(f \mrm w) = \hterm(p)$ is guaranteed.
\item We show that every non-zero polynomial $p \in \ideal{r}{}(F)
  \backslash \{ 0 \}$ is top-reducible to zero using $F$  by induction
  on $\hterm(p)$. 
  Let $\hterm(p) = \min \{ \hterm(g) | g \in \ideal{r}{}(F) \backslash \{ 0 \} \}$.
  Then, as $p \in \ideal{r}{}(F)$ and $F$ is a  standard basis, we
   have $p = \sum_{i=1}^{k} \alpha_i \skm f_{i} \mrm w_i$, with 
   $\alpha_i \in \myk^*, f_{i} \in F, w_i \in \m$ and 
   $\hterm(p) \succeq \hterm(f_{i} \mrm w_i)$ for all  $1 \leq i \leq k$.
  Without loss of generality, let $\hterm(p) = \hterm( f_1 \mrm w_1)$.
  Hence, $p$ is strongly right reducible by $f_1$.
  Let $p \red{}{\myr}{s}{f_1} q$, i.e.,
   $q = p - \hc(p) \skm \hc(f_1 \mrm w_1)^{-1} \skm f_1 \mrm w_1$, and
   by the definition of strong right reduction the term $\hterm(p)$ is
   eliminated from $p$ implying that $\hterm(q) \pred \hterm(p)$ as $q < p$.
   Hence, as $q \in \ideal{r}{}(F)$ and as $\hterm(p)$ was minimal among
   the head terms of the non-zero elements
   in the right ideal generated by $F$, this implies $q=0$, and, 
   therefore, $p$ is strongly right top-reducible to zero by $f_1$ in 
   one step.
   On the other hand, in case 
    $\hterm(p) \succ \min \{ \hterm(g) | g \in \ideal{r}{}(F)\backslash \{ 0\} \}$, by the
    same arguments used before we can reduce $p$ to a polynomial
    $q$ with $\hterm(q) \pred \hterm(p)$, and, thus, by our induction
    hypothesis we know that $q$ is top-reducible to zero.
   Therefore, as the reduction step $p \red{}{\myr}{s}{f_1} q$ takes
    place at the head term of $p$, $p$ is also top-reducible to
    zero.
\\
\qed
\end{enumerate}\renewcommand{\baselinestretch}{1}\small\normalsize
Now we can prove how standard
 representations can be used to characterize strong Gr\"obner bases, and that in
 fact standard bases are strong Gr\"obner bases.
\begin{theorem}~\\
{\sl
For a set  $F$ of polynomials in $\myk[\m]$,
 the following statements are equivalent:
\begin{enumerate}
\item $F$ is a strong Gr\"obner basis.
\item For all polynomials $g \in \ideal{r}{}(F)$ we have $g \red{*}{\myr}{s}{F} 0$.
\item $F$ is a standard basis.
\theoend
\end{enumerate}
}
\end{theorem}
\Ba{}~\\
\mbox{$1 \R 2:$ }
  By the definition of strong Gr\"obner bases in \ref{def.gb} we know that $g \in
  \ideal{r}{}(F)$ implies $g \red{*}{\lr}{s}{F} 0$, and since
  $\red{}{\myr}{s}{F}$ is confluent and $0$ irreducible, $g \red{*}{\myr}{s}{F} 0$ follows
  immediately.

\mbox{$2 \R 3:$ }
 This follows directly from lemma \ref{lem.srprop}.

\mbox{$3 \R 1:$ }
  In order to show that $F$ is a strong Gr\"obner basis, we have
  to prove two subgoals:
$\red{*}{\lr}{s}{F} = \;\; \equiv_{\ideal{r}{}(F)}$ follows immediately from lemma \ref{lem.strong.congruence}.
It remains to show that $\red{}{\myr}{s}{F}$ is confluent.
             Since $\red{}{\myr}{s}{F}$ is Noetherian, we only have to prove
             local confluence.
             Let us suppose there exist polynomials $g,g_1,g_2 \in \myk[\m]$ such that
             we have $g \red{}{\myr}{s}{F} g_1$, $g \red{}{\myr}{s}{F} g_2$
              and $g_1 \neq g_2$.
             Then $g_1 - g_2 \in \ideal{r}{}(F)$ and, therefore, 
               is top-reducible to
              zero by  $F$ as a result of  lemma \ref{lem.srprop}. 
             Thus lemma \ref{lem.confluent} provides the existence of a
              polynomial 
              $h \in \myk[\m]$ such that
              $g_1 \red{*}{\myr}{s}{F} h$ and $g_2 \red{*}{\myr}{s}{F} h$,
              i.e., $\red{}{\myr}{s}{F}$ is confluent.
\\
\qed
%
In accordance with the terminology used in Buchberger's approach 
 to describe Gr\"obner bases we define critical pairs of polynomials 
 with respect to strong right reduction.
\begin{definition}\label{def.cps}~\\
{\rm
Given two non-zero polynomials $p_{1}, p_{2} \in \myk[\m]$\footnote{Notice
                                        that $p_1=p_2$ is possible.},
 every pair $w_{1}, w_{2} \in \m$ such that
 $\hterm(p_1 \mrm w_{1}) = \hterm(p_2 \mrm w_{2})$, defines a
 \index{strong!s-polynomial}\index{s-polynomial!strong}
 \betonen{strong s-polynomial}\/
$$ \spol{s}(p_{1}, p_{2}, w_{1}, w_{2}) = \hc(p_1 \mrm w_1)^{-1} \skm p_1 \mrm w_1
                                      - \hc(p_2 \mrm w_2)^{-1} \skm p_2 \mrm w_2.$$
Let  $U_{p_1,p_2} \subseteq \m \times \m$ be the set containing {\em
  all}\/ such pairs $w_1,w_2 \in\m$.
\dend
}
\end{definition}
A strong s-polynomial will be called non-trivial in case it is non-zero and
notice that for non-trivial s-polynomials we always have $\hterm(\spol{}(p_{1}, p_{2}, w_{1}, w_{2}))
\pred \hterm(p_1 \mrm w_{1}) = \hterm(p_2 \mrm w_{2})$.
The set $U_{p_1,p_2}$ is  contained in the set of all solutions to the
equations in two variables of the form $u \mm {\rm x} = v \mm {\rm y}$
 where $u \in \terms(p_1)$ and $v \in\terms(p_2)$.
It can be empty, finite or even infinite.
As might be expected, we can give a criterion that implies confluence
for strong right reduction in terms of strong s-polynomials.
\begin{theorem}\label{theo.pcs}~\\
{\sl 
For a set  $F$ of polynomials in $\myk[\m]$,
 the following statements are equivalent:
\begin{enumerate}
\item For all polynomials $g \in \ideal{r}{}(F)$
       we have $g \red{*}{\myr}{s}{F} 0$.
\item For all  not necessarily different polynomials 
       $f_{k}, f_{l} \in F$ and every corresponding pair $(w_{k}, w_{l}) \in U_{f_k,f_l}$  we have 
       $ \spol{s}(f_{k}, f_{l}, w_{k}, w_{l}) \red{*}{\myr}{s}{F} 0$.
\end{enumerate}
}
\end{theorem}
\Ba{}~\\
\mbox{$1 \R 2:$ }
     Let $ (w_{k}, w_{l}) \in U_{f_k,f_l}$ give us a strong
     s-polynomial belonging to the polynomials $f_k,f_l$.
     Then by definition \ref{def.cps} we get
     $$\spol{s}(f_{k}, f_{l}, w_{k}, w_{l}) =
       \hc(f_k \mrm w_k)^{-1} \skm f_{k} \mrm w_{k} -\hc(f_l \mrm w_l)^{-1} \skm f_{l} \mrm w_{l}
       \:\in \ideal{r}{}(F)$$
      and, thus, $\spol{s}(f_{k}, f_{l}, w_{k}, w_{l})
      \red{*}{\myr}{s}{F} 0$.

\mbox{$2 \R 1:$ }
     We have to show that every  non-zero polynomial 
      $g \in \ideal{r}{}(F) \backslash \{ 0 \}$ 
      is $\red{}{\myr}{s}{F}$-reducible to zero.
     Remember that for
      $h \in \ideal{r}{}(F)$, $ h \red{}{\myr}{s}{F} h'$ implies $h' \in \ideal{r}{}(F)$.
     Hence, as  $\red{}{\myr}{s}{F}$ is Noetherian,
      it suffices to show that every  $g \in \ideal{r}{}(F) \backslash
      \{ 0 \}$ is $\red{}{\myr}{s}{F}$-reducible. 
     Now, let $g = \sum_{j=1}^m \alpha_{j} \skm f_{j} \mrm w_{j}$ be a
      representation of a non-zero  polynomial $g$ such that
      $\alpha_{j} \in \myk^*, f_j \in F$, and $w_{j} \in \m$.
     Depending on this  representation of $g$ and the
      well-founded total ordering $\succeq$ on $\m$ we define
      $t = \max \{ \hterm(f_{j} \mrm w_{j}) \mid j \in \{ 1, \ldots m \}  \}$ and
      $K$ is the number of polynomials $f_j \mrm w_j$ containing $t$ as a term.
Then $t \succeq \hterm(g)$ and 
 in case $\hterm(g) = t$ this immediately implies that $g$ is
 $\red{}{\myr}{s}{F}$-reducible. 
So
by lemma \ref{lem.srprop} it is sufficient to  show that
$g$ has a standard representation, as this implies that $g$ is
top-reducible using $F$.
This will be done by induction
 on $(t,K)$, where
 $(t',K')<(t,K)$ if and only if $t' \prec t$ or $(t'=t$ and
 $K'<K)$\footnote{Note that this ordering is well-founded since $\succ$
                  is well-founded on ${\cal T}$ and $K \in\n$.}.
In case $t \succ \hterm(g)$ there are  two polynomials $f_k,f_l$ in the corresponding 
      representation\footnote{Not necessarily $f_l \neq f_k$.}
      such that  $\hterm(f_k \mrm w_k) = \hterm(f_l \mrm w_l)$.
     By definition \ref{def.cps} we have a strong s-polynomial
      $\spol{s}(f_k,f_l,w_k,w_l) = \hc(f_k \mrm w_k)^{-1} \skm  f_k
      \mrm w_k -\hc(f_l \mrm w_l)^{-1}\skm f_l \mrm w_l$ corresponding
      to this overlap.
We will now change our representation of $g$ by using the additional
information on this s-polynomial in such a way that for the new
representation of $g$ we either have a smaller maximal term or the occurrences of the term $t$
are decreased by at least 1.
     Let us assume  $\spol{s}(f_k,f_l,w_k,w_l) \neq 0$\footnote{In case 
               $\spol{s}(f_k,f_l,w_k,w_l) = 0$,
               just substitute $0$ for the sum 
               $\sum_{i=1}^n \delta_i \skm h_i \mrm
               v_i$ in the equations below.}. 
     Hence, the reduction sequence  $\spol{s}(f_k,f_l,w_k,w_l)
     \red{*}{\myr}{s}{F} 0 $ results in a standard representation
      $\spol{s}(f_k,f_l,w_k,w_l) =\sum_{i=1}^n \delta_i \skm h_i \mrm v_i$,
      where $\delta_i \in \myk^*,h_i \in F$, and $v_i \in \m$ and  all terms occurring in the sum are bounded by
       $\hterm(\spol{s}(f_k,f_l,w_k,w_l)) \pred t$.
     This gives us: 
     \begin{eqnarray}
       &  & \alpha_{k} \skm f_{k} \mrm w_{k} + \alpha_{l} \skm f_{l} \mrm w_{l}
             \nonumber\\  
       &  &  \nonumber\\                                                               & = &  \alpha_{k} \skm f_{k} \mrm w_{k} +
              \underbrace{ \alpha'_{l} \skm \beta_k \skm f_{k} \mrm w_{k}
                   - \alpha'_{l} \skm \beta_k \skm f_{k} \mrm w_{k}}_{=\, 0} 
                   + \alpha'_{l}\skm \beta_l  \skm f_{l} \mrm w_{l} \nonumber\\
    & & \nonumber\\ 
       & = & (\alpha_{k} + \alpha'_{l} \skm \beta_k) \skm f_{k} \mrm w_{k} - \alpha'_{l} \skm 
               \underbrace{(\beta_k \skm f_{k} \mrm w_{k}
             -  \beta_l \skm f_{l} \mrm w_{l})}_{=\,
             \spol{s}(f_k,f_l,w_k,w_l)} 
             \nonumber\\
    & & \nonumber\\ 
       & = & (\alpha_{k} + \alpha'_{l} \skm \beta_k) \skm f_{k} \mrm w_{k} - \alpha'_{l} \skm
                   (\sum_{i=1}^n \delta_{i} \skm h_{i} \mrm v_{i}) \label{s1}
     \end{eqnarray}
     where $\beta_k=\hc(f_k \mrm w_k)^{-1}$, $\beta_l=\hc(f_l \mrm w_l)^{-1}$
      and  $\alpha'_l \skm \beta_l = \alpha_l$.
     By substituting (\ref{s1}) in our representation of $g$ 
 either $t$ disappears   or in
 case $t$ remains maximal among the terms occurring in the new
 representation of $g$, $K$ is decreased.
\\
\qed
Note that this theorem, although characterizing a strong Gr\"obner basis by
 strong s-polynomials, does not give a finite test to check
 whether a set is a strong Gr\"obner basis, since in general infinitely
 many strong s-polynomials have to be considered.
The following example shows how already
 two polynomials $p_1,p_2$ can cause infinitely many critical situations.
\begin{example}~\\
{\rm
Let  $\Sigma = \{ a,b,c,d,e,f \}$ and 
 $T = \{ abc \myr ba, fbc \myr bf, bad \myr e \}$ be a presentation of a
 monoid $\m$ with a length-lexicographical ordering
 induced by
 $a \succ b \succ c \succ d \succ e \succ f$.
Further consider two polynomials  $p_1 = a + f, p_2 = bf + a \in \q[\m]$.
\\
Then we get  infinitely many critical situations $\hterm(p_1  \mrm (bc)^idw) = f \mm (bc)^idw = bf \mm (bc)^{i-1} dw = \hterm(p_2 \mrm (bc)^{i-1} dw)$, 
 where $i \in \n^+, w \in \m$, resulting in  infinitely many strong s-polynomials
 $$\spol{s}(p_1,p_2,(bc)^idw,(bc)^{i-1} dw) = (a+f) \mrm (bc)^idw -
 (bf +a) \mrm (bc)^{i-1} dw$$
 and
 $U_{p_1,p_2} = \{ ( (bc)^idw, (bc)^{i-1} dw) | i \in\n^+, w \in\m
 \}$.
\exaend
}
\end{example}
Localization of critical situations might be very hard.
As the previous example shows,
 the set $U_{p_1,p_2}$ need not have a ``suitable'' finite basis,
 i.e., there need not exist a finite set $B \subseteq U_{p_1,p_2}$ such that
 for every  pair $(w_1,w_2) \in U_{p_1,p_2}$ there exists a pair $(u_1,u_2) \in
 B$ and an element $w \in \m$ with $u_1 \mm w = w_1$ and $u_2 \mm w = w_2$. 

One way to reduce the set of critical situations that have to be
 considered is to weaken reduction.
The key idea is that for two reduction relations  $\red{}{\myr}{1}{}$
 and $\red{}{\myr}{2}{}$ on a set ${\cal E}$ such that
 $\red{}{\myr}{1}{} \subseteq \red{}{\myr}{2}{}$ and
 $\red{*}{\lr}{1}{} = \red{*}{\lr}{2}{}$, the confluence of
 $\red{}{\myr}{1}{}$ on ${\cal E}$ implies the
 confluence of $\red{}{\myr}{2}{}$ on ${\cal E}$.
The next section will introduce a way to weaken strong right reduction.

%% file: rightreduction.tex
\section{The Concept of Right Reduction}\label{section.rightreduction}
%
%
%
In the previous section we have introduced standard representations
 to monoid rings and it was shown how they are related to a
 special reduction.
We will now slightly extend this definition in order to reflect
 a possible weakening of strong right reduction and study what can
 be gained by this approach.
\begin{definition}\label{def.rsr}~\\
{\rm
Let $F$ be a set of polynomials  and $p$ a non-zero
polynomial in $\myk[\m]$.
A representation 
$$p = \sum_{i=1}^{n} \alpha_i \skm f_{i} \mrm w_i, \mbox{ with } \alpha_i \in \myk^*,
f_{i} \in F, w_i \in \m $$
is called a 
 \index{stable!standard representation}\index{standard representation!stable}\betonen{stable standard representation}\/ of $p$ with respect to $F$,
if for all $1 \leq i \leq n$ we have $\hterm(p) \succeq \hterm(f_{i}) \mm
w_i \succeq  \hterm(f_{i} \mrm w_i)$.
A set $F \subseteq \myk[\m]$ is called a 
 \index{stable!standard basis}\index{standard basis!stable}\betonen{stable standard basis}, if every
 non-zero polynomial in $\ideal{r}{}(F)$ has a
 stable standard representation with
 respect to $F$.
\dend
}
\end{definition}
Notice that in this definition $\hterm(f_i)
  \mm w_i \succeq \hterm(f_i \mrm w_i)$ need not imply $\hterm(f_i)
  \mm w_i =  \hterm(f_i \mrm w_i)$ as $\hterm(f_i) \mm w_i \not\in
  \terms(f_i \mrm w_i)$ is possible, e.g. if $\m$ is not cancellative.
Nevertheless, in case $\hterm(p) =
 \hterm(f_{i} \mrm w_i)$ we have $\hterm(f_{i} \mrm w_i)=\hterm(f_{i}) \mm w_i$
 and this situation occurs for at least one polynomial in the
 representation. 
In stable standard bases we can refine
the concept of stable standard representations and gain more
information on the head terms of the multiples of the polynomials
involved in the respective representation.
This is due to the fact that if we have $\hterm(f_i) \mm w_i \succeq
\hterm(f_i \mrm w_i)$, but $\hterm(f_i) \mm w_i \neq
\hterm(f_i \mrm w_i)$, i.e., $\hterm(f_i) \mm w_i \not\in \terms(f_i
\mrm w_i)$, then we again have a stable standard representation for
this multiple as it belongs to the right ideal generated by $F$.
This is reflected in the proof of the following lemma.
\begin{lemma}\label{lem.ssr.prop1}~\\
{\sl
Let $F$ be a stable standard basis in $\myk[\m]$.
Then every non-zero polynomial $p$ in $\ideal{r}{}(F)$ has a stable standard
 representation
 $p = \sum_{i=1}^{n} \alpha_i \skm f_{i} \mrm w_i, \mbox{ with } \alpha_i \in \myk^*,
  f_{i} \in F$, and $w_i \in \m$
such that  for all $1 \leq i \leq n$ we even have $$\hterm(p) \succeq 
 \hterm(f_{i}) \mm w_i = \hterm(f_{i} \mrm w_i).$$
\lemend
}
\end{lemma}
\Ba{}~\\
Let $p$ be a non-zero polynomial in $\ideal{r}{}(F)$.
We show our claim by induction on $\hterm(p)$.
In the base case we can assume
 $\hterm(p) = \min \{ \hterm(g) | g \in \ideal{r}{}(F)\backslash \{ 0 \}\}$.
Since $F$ is a stable standard basis, we know that the polynomial $p$ has a stable standard representation
  $p = \sum_{i=1}^{n} \alpha_i \skm f_{i} \mrm w_i, \mbox{ with } \alpha_i \in \myk^*,
 f_{i} \in F, w_i \in \m$
 such that  for all $1 \leq i \leq n$ we have $\hterm(p) \succeq \hterm(f_{i})
 \mm w_i \succeq \hterm(f_{i} \mrm w_i)$.
It remains to show that we can even achieve a representation where $ \hterm(f_{i})
 \mm w_i = \hterm(f_{i} \mrm w_i)$.
Without loss of generality we can assume that $\hterm(p) =
 \hterm(f_1) \mm w_1=\hterm(f_1 \mrm w_1)$.
Then the head term of $p$ can be
 eliminated by subtracting an appropriate right multiple of  $f_1$.
Looking at the polynomial 
 $h = p - \hc(p) \skm \hc(f_1 \mrm w_1)^{-1} \skm f_1 \mrm w_1$ 
 we find that $h$ lies in the right ideal generated by $F$ and, since
 $\hterm(p)$ is minimal and we find $\hterm(h) \prec \hterm(g)$, we can conclude that
 $h=0$.
Thus $p$ has a stable standard
 representation $p = \hc(p) \skm \hc(f_1 \mrm w_1)^{-1} \skm f_1 \mrm
 w_1$ of the desired form.
Now let us suppose $\hterm(p) \succ \min \{\hterm(g) | g \in \ideal{r}{}(F)\backslash \{ 0 \}\}$.
Then again resulting from the existence of a stable
 standard representation for $p$, 
 there exists a polynomial $f_1 \in F$ such that
 $\hterm(p)=\hterm(f_1) \mm w_1= \hterm(f_1 \mrm w_1)$
 for some $w_1 \in \m$.
Hence looking at the polynomial 
 $h = p - \hc(p) \skm \hc(f_1 \mrm w_1)^{-1} \skm f_1 \mrm w_1$  
 we know that $h$ lies in the right ideal generated by $F$ and since
 $\hterm(h) \prec \hterm(p)$ either $h=0$, giving us 
 $p = \hc(p) \skm \hc(f_1 \mrm w_1)^{-1} \skm f_1 \mrm w_1$
 or  our induction hypothesis yields the existence
 of a  stable standard
 representation of the desired form for $h$, say 
 $h = \sum_{j=1}^{m} \beta_j \skm g_j \mrm v_j$, with $\beta_j \in
 \myk^*, g_j \in F$, $v_j \in \m$.
Thus we have a  stable  standard  representation of the polynomial 
 $p$, namely $p =
 \sum_{j=1}^{m} \beta_j \skm g_j \mrm v_j +  \hc(p) \skm \hc(f_1 \mrm
   w_1)^{-1} \skm f_1 \mrm w_1$, which satisfies our requests. 
\\
\qed
Similar as standard representations correspond to strong right
reduction, 
stable standard representations correspond to
a  weakening of strong right reduction.
Instead of using all right multiples of a polynomial  by monomials as rules we
restrict ourselves to those  right
multiples of a polynomial which allow the head
term of the polynomial to keep its head position.
Hence, reduction defined in this way can be called ``stable''.
\begin{definition}\label{def.redr}~\\
{\rm
Let $p, f$ be two non-zero polynomials in $\myk[\m]$. 
We say $f$ 
 \index{right!reduction}\index{reduction!right}\betonen{right reduces} $p$ to $q$ at a monomial
 $
\alpha \skm t$ of $p$ in one step, denoted by $p \red{}{\myr}{r}{f} q$, if
\begin{enumerate}
\item[(a)] $\hterm(f \mrm w) = \hterm(f) \mm w  = t$ for some $w \in \m$, and
\item[(b)] $q = p - \alpha \skm \hc(f \mrm w)^{-1} \skm f \mrm w$.
\end{enumerate}
We write $p \red{}{\myr}{r}{f}$ if there is a polynomial $q$ as defined
above and $p$ is then called  right reducible by $f$. 
Further we can define $\red{*}{\myr}{r}{}, \red{+}{\myr}{r}{}$ and
 $\red{n}{\myr}{r}{}$ as usual.
Right reduction by a set $F \subseteq \myk[\m]$ is denoted by
 $p \red{}{\myr}{r}{F} q$ and abbreviates $p \red{}{\myr}{r}{f} q$
 for some $f \in F$,
 which is also written as  $p \red{}{\myr}{r}{f \in F} q$.
\dend
}
\end{definition}
In case $f$ right reduces $p$ to $q$ at the monomial $\alpha \skm t$, then $t
\not\in \terms(q)$.
Furthermore, as in lemma \ref{lem.red}, we have $p > q$, right
reduction is Noetherian, and $p \red{}{\myr}{r}{q_1} 0$ and $q_1 \red{}{\myr}{r}{q_2} 0$ imply $p \red{}{\myr}{r}{q_2} 0$.

Defining interreduced sets with respect to right reduction we have to
 be careful as in the case of strong right reduction.
\begin{definition}~\\
{\rm
We call a set of polynomials $F \subseteq \myk[\m]$ 
 \index{reduced set}\index{interreduced!(set of polynomials)}
 \betonen{interreduced} or \betonen{reduced} with respect to $\red{}{\myr}{r}{}$, if for all $f \in F$, $f \red{}{\myr}{r}{F} f'$ implies $f'=0$.
\dend
}
\end{definition}
Again right reducing a polynomial with itself need not result in zero
 and interreduced bases for right ideals need not exist.
\begin{example}~\\
{\rm
Let $\Sigma = \{ a,b,c,d \}$ and $T = \{ ab \myr c, b^2 \myr \lambda, cb
\myr d \}$ be a presentation of a cancellative monoid with a
length-lexicographical ordering induced by $a \succ b\succ c \succ d$.
\\
Then the polynomial $a+b+c$ is right reducible by itself as follows:
$a + b + c \red{}{\myr}{r}{a+b+c} a+b+c - (a+b+c) \mrm b = a+b+c -
(\underline{c} + \lambda + d) = a+b - d - \lambda$.
\exaend
}
\end{example}
Note that we only use $\hm(f) \myr -\reductum(f)$ as a rule in case we
 have $\hterm(f \mrm w) = \hterm(f) \mm w$.
We cannot always  use $\hm(f) \myr -\reductum(f)$, since then reduction would no
 longer be Noetherian, i.e., 
 infinite reduction sequences could arise.
This is due to the fact that multiplying $f$ by an element $w \in \m$
can cause $\hm(f) \mm w < \reductum(f) \mrm w$.
The following example  illustrates this phenomenon.
\begin{example}\label{exa.3}~\\
{\rm
Let $\Sigma = \{ a, b \}$ and $T = \{ab \myr \lambda, ba \myr
 \lambda \}$ be a presentation of a group $\g$ with a
 length-lexicographical ordering induced by
  $a \succ b$.
\\
Suppose we use a polynomial $f$ itself as a rule $\hm(f) \myr -\reductum(f)$, i.e.,
 in order to reduce a monomial $\alpha \skm t$ 
 we just require the existence of $w \in\m$ such that $t=\hterm(f) \mm w$
 similar to Buchberger's reduction in definition
 \ref{def.buchberger.red}.
Then  we could  right reduce the
 polynomial $b^2 + 1 \in \q[\m]$ at the monomial $b^2$ by the polynomial
 $a + b$ as $b^2 = a \mm b^3$.
This would give us:
$$ b^{2} + 1 \red{}{\myr}{}{a + b} b^{2} + 1 - (a + b)
 \mrm b^{3} = -b^{4} + 1$$
 and the polynomial $-b^{4} + 1$ likewise would be reducible by $a +
 b$  at the monomial $b^{4}$ causing an 
 infinite reduction sequence.
\exaend
}
\end{example} 
Hence right reduction using a polynomial $p$ corresponds to the following 
 set of
rules $\{ \alpha \skm \hm(p \mrm w) \myr (- \alpha) \skm \reductum(p
\mrm w) \mid \alpha \in \myk^*, w \in \m, \hterm(p \mrm w) = \hterm(p)
\mm w \}$, while strong reduction corresponds to the set
$\{ \alpha \skm \hm(p \mrm w) \myr (- \alpha) \skm \reductum(p
\mrm w) \mid \alpha \in \myk^*, w \in \m, p \mrm w \neq 0 \}$.
Looking at the expressiveness of right reduction we find that
 while  item 3 of lemma \ref{lem.red}  holds for
      our weaker form of reduction, i.e., $p \red{}{\myr}{r}{q_1} 0$
      together with $q_1 \red{}{\myr}{r}{q_2} 0$ yields $p
      \red{}{\myr}{r}{q_2} 0$, we no longer have $\alpha \skm p \mrm w
      \red{{ \leq 1}}{\myr}{r}{p} 0$ for $\alpha \in \myk^*, p \in \myk[\m]$ and $w \in \m$.  
Reviewing example  \ref{exa.3} we find
      that $(a + b) \mrm b = b^2 + 1$ is not right reducible  by $a + b$.
Similar to the case of strong right reduction, $p \red{}{\myr}{r}{q}$ and $q \red{}{\myr}{r}{q_1} q_2$ need not imply $p \red{}{\myr}{r}{\{ q_1,q_2 \}} $ as the following example shows. 
\begin{example}~\\
{\rm
      Let $\Sigma = \{ a, b, c \}$ and
       $T = \{ a^{2} \myr \lambda, b^{2} \myr \lambda, ab \myr c, ac \myr b, cb \myr a \}$
       be a presentation of a monoid $\m$ (which is in fact a group)
       with a
length-lexicographical ordering induced by $a \succ b \succ c$.
      \\
      Looking at the polynomials $p = ba + b, q = a + \lambda$ and
       $q_1 = c^2 + b\in \q[\m]$ we get
       $p \red{}{\myr}{r}{q} p - q \mrm ca = ba + b - (a + \lambda) \mrm
       ca = ba + b - ba - ca = -ca + b$
       and $q \red{}{\myr}{r}{q_1} q - q_1 \mrm bc = a + \lambda - (c^2
       + b) \mrm bc = a + \lambda - a -c = -c + \lambda = q_{2}$,
      but $p \nred{}{\myr}{r}{\{q_1,q_2\}}$.
\\
      Trying to reduce $ba$ by $q_1$ or $q_{2}$ we get
       $q_1 \mrm bc^{2}a = ba + \underline{c^{2}a}$ and
       $q_{2} \mrm bca = -ba + \underline{bca}$ both violating
       condition (a) of definition \ref{def.redr}.
      The same happens trying to reduce $b$, as
       $q_1 \mrm bc^2=b+\underline{c^2}$ and $q_2 \mrm bc=-b +\underline{bc}$.
\exaend
}
\end{example}
Nevertheless, an essential property for reduction  to allow a
 characterization of Gr\"obner bases by
 s-polynomials is still true, as an analogon to lemma
 \ref{lem.confluent} holds for right reduction.
\begin{lemma} \label{lem.confluentr}~\\
{\sl
Let $F$ be a set of polynomials and $p,q,h$ some
 polynomials in $\myk[\m]$.
\begin{enumerate}
\item
Let $p-q \red{}{\myr}{r}{F} h$.
Then there are  polynomials $p',q' \in \myk[\m]$ such that we have
 $p  \red{*}{\myr}{r}{F} p', q  \red{*}{\myr}{r}{F} q'$ and $h=p'-q'$.
\item
Let $0$ be a normal form of $p-q$ with respect to $\red{}{\myr}{r}{F}$.
Then there exists a polynomial  $g \in \myk[\m]$ such that
 $p  \red{*}{\myr}{r}{F} g$ and $q  \red{*}{\myr}{r}{F} g$.
\lemend
\end{enumerate}
}
\end{lemma}
\Ba{}
\begin{enumerate}
\item  Let $p-q \red{}{\myr}{r}{F} h = p-q-\alpha \skm f \mrm w$, where
        $\alpha \in \myk^*, f \in F, w \in \m$
        and $\hterm(f \mrm w) = \hterm(f) \mm w = t$, i.e.,
        $\alpha \skm \hc(f \mrm w)$
        is the coefficient of $t$ in $p-q$.
       We have to distinguish three cases:
       \begin{enumerate}
         \item $t \in \terms(p)$ and $t \in \terms(q)$:
               Then we can eliminate the term $t$ in the polynomials 
                $p$ respectively $q$ by right reduction and 
                get $p \red{}{\myr}{r}{f} p - \alpha_1 \skm f \mrm w= p'$,
                $q \red{}{\myr}{r}{f} q - \alpha_2 \skm f \mrm w= q'$,
                and $\alpha_1  -  \alpha_2 = \alpha$,
                where $\alpha_1 \skm \hc(f \mrm w)$ and 
                $\alpha_2 \skm \hc(f \mrm w)$ are
                the coefficients of $t$ in $p$ respectively $q$.
         \item $t \in \terms(p)$ and $t \not\in \terms(q)$:
               Then  we can eliminate the term $t$ in the polynomial 
                $p$  by right reduction and  get
                $p \red{}{\myr}{r}{f} p - \alpha \skm f \mrm w= p'$ 
                and $q = q'$.
         \item $t \in \terms(q)$ and $t \not\in \terms(p)$:
               Then  we can eliminate the term $t$ in the polynomial 
                $q$ by right reduction and  get
                $q \red{}{\myr}{r}{f} q + \alpha \skm f \mrm w= q'$
                and $p = p'$.
       \end{enumerate}
      In all three cases we have $p' -q' =  p - q - \alpha \skm f \mrm w = h$.
\item We show our claim by induction on $k$, where $p-q \red{k}{\myr}{r}{F} 0$.
      In the base case $k=0$ there is nothing to show.
      Thus let $p-q \red{}{\myr}{r}{F} h  \red{k}{\myr}{r}{F} 0$.
      Then by (1) there are polynomials $p',q' \in \myk[\m]$ such that 
       $p  \red{*}{\myr}{r}{F} p', q  \red{*}{\myr}{r}{F} q'$ and $h=p'-q'$.
      Now the induction hypothesis for $p'-q' \red{k}{\myr}{s}{F} 0$  yields 
       the existence of a polynomial $g \in \myk[\m]$ such that
       $p  \red{*}{\myr}{r}{F} p' \red{*}{\myr}{r}{F} g$ and
       $q  \red{*}{\myr}{r}{F} q' \red{*}{\myr}{r}{F} g$.
\\
\qed
\end{enumerate}\renewcommand{\baselinestretch}{1}\small\normalsize
But, unlike strong right reduction, right reduction no longer captures the right ideal congruence.
\begin{lemma}\label{lem.rdefect}~\\
{\sl
Let  $F$ be a set of polynomials and $p,q$ some polynomials in $\myk[\m]$.
Then  $p \red{*}{\lr}{r}{F} q$ implies $p - q \in \ideal{r}{}(F) $ but not
 vice versa.
\lemend\ohnebeweis
}
\end{lemma}
\begin{example}\label{exa.defect}~\\
{\rm\renewcommand{\baselinestretch}{1.1}\small\normalsize
Let $\Sigma = \{ a, b, c \}$ and
    $T = \{ a^{2} \myr \lambda, b^{2} \myr \lambda, ab \myr c, ac \myr b, cb \myr a \}$
 be a presentation of a monoid $\m$ (which is in fact a group)
  with  a length-lexicographical ordering induced by $a \succ b \succ c$.
\\
Inspecting the polynomials $p = a+b+c, q = b - \lambda\in\q[\m]$ and the set
 $F = \{ a+b+c \}\subseteq \q[\m]$ we get 
 $p-q = a+c+ \lambda = (a+b+c) \mrm b \in \ideal{r}{}(F)$, but
 $a+b+c \nred{*}{\lr}{r}{F}  b - \lambda$.
To prove this claim, let us assume  $a+b+c \red{*}{\lr}{r}{F}  b - \lambda$.
Then, since  $a+b+c  \red{}{\myr}{r}{F} 0$, we get
 $ b - \lambda  \red{*}{\lr}{r}{F} 0$.
Let $n \in\n^+$ be minimal such that $b - \lambda \red{n}{\lr}{r}{F} 0$.
As $b - \lambda \nred{}{\myr}{r}{F} 0$  we know $n>1$.
Thus, let us look at the sequence 
$$b - \lambda =: p_0  \red{}{\lr}{r}{F} p_1 \red{}{\lr}{r}{F}
           \ldots  \red{}{\lr}{r}{F} p_{n-1} \red{}{\lr}{r}{F} 0,$$
 where for all  $1 \leq i \leq n-1$,
 $p_i = p_{i-1} + \alpha_i \skm (a+b+c) \mrm w_i$, $\alpha_i \in \myk^*, w_i \in\m$ and
 $\hterm((a + b + c) \mrm w_i)=a \mm w_i$.
Further let $t = \max \{ \hterm(p_i) \mid 1 \leq i \leq n-1 \}$.
Then $t \succ b$, as $\hterm((a+b+c) \mrm w) \succ b$ for all $w \in\m$.
Let $p_l$ be the first polynomial with $\hterm(p_l)=t$, i.e.,
 $\hterm(p_j) \prec t$ for all $j < l$,
 and let $p_{l+k}$ be the next polynomial,
 where the occurrence of $t$ is changed.
 Since $\hterm((a+b+c) \mrm w_{l+k})= a \mm w_{l+k} =t=
                                  a \mm w_l =\hterm((a+b+c) \mrm w_l)$
  and $(\Sigma,T)$ presents a group, we can conclude $w_{l+k}=w_l$.
 Further our transformation sequence is supposed to be minimal, i.e.,
  $t$ is not changed by the
  reductions taking place in the sequence $p_l \red{k-1}{\lr}{r}{F} p_{l+k-1}$.
But then, eliminating $p_l$  and substituting 
 $p_{l+j}$ by $p'_{l+j}=p_{l+j} -  \alpha_{l} \skm (a+b+c) \mrm w_l$
 for all $1 \leq j<k$ gives us a shorter sequence 
 $b - \lambda \red{n-1}{\lr}{r}{F} 0$
 contradicting our assumption.
\renewcommand{\baselinestretch}{1}\small\normalsize
\exaend
}
\end{example}
Obviously for a set of polynomials $F$ we have $\red{}{\myr}{r}{F} \subseteq \red{}{\myr}{s}{F}$, but as seen in
 example \ref{exa.defect} in general we cannot expect 
 $\red{*}{\lr}{r}{F} = \red{*}{\lr}{s}{F}$ to hold and right
 reduction does not capture the right ideal congruence in general.
To overcome this problem we will use special sets of polynomials for reduction.
The following lemma shows that for special bases of right ideals right
reduction corresponds to the right ideal congruence.
\begin{lemma}\label{lem.ssb.cong}~\\
{\sl
Let $F$ be a stable standard basis and $p,q,h$ some polynomials in
 $\myk[\m]$.
Then
  $$p \red{*}{\lr}{r}{F} q \mbox{ if and only if } p - q \in
  \ideal{r}{}(F).$$
\lemend
}
\end{lemma}
\Ba{}~\\
In order to prove our claim we have to show two subgoals.
The inclusion $\red{*}{\lr}{r}{F}  \subseteq \;\; \equiv_{\ideal{r}{}(F)}$
 is an immediate consequence of the definition of right reduction and
 can be shown by induction as in lemma \ref{lem.strong.congruence}. 
To prove the converse inclusion 
 $\equiv_{\ideal{r}{}(F)} \: \subseteq \red{*}{\lr}{r}{F} $
 we have to use some additional information.
Remember that $p \equiv_{\ideal{r}{}(F)} q$ implies
 $p = q + \sum_{j=1}^{m} \alpha_{j} \skm f_j \mrm w_{j}$, where
 $\alpha_{j} \in \myk^*, f_j \in F$, and $w_{j} \in \m$.
As  every
 multiple $f_j \mrm w_{j}$ lies in $\ideal{r}{}(F)$ and
 $F$ is a stable standard basis, 
 by lemma \ref{lem.ssr.prop1}, we can assume $\hterm(f_j \mrm w_j)=\hterm(f_j) \mm w_j$ for
 all polynomials occurring in the sum.
Using these assumptions the proof of our claim can now be accomplished
by induction on $m$ as in lemma \ref{lem.strong.congruence}.
\auskommentieren{
In the base case $m = 0$ there is nothing to show.
\\
Let  $p = q + \sum_{j=1}^{m} \alpha_{j} \skm f_j \mrm w_{j} +
 \alpha_{m+1} \skm f_{m+1} \mrm w_{m+1}$
 and by our induction hypothesis
 $p \red{*}{\lr}{r}{F} q + \alpha_{m+1} \skm f_{m+1} \mrm w_{m+1}$.
\\
Let $t=\hterm(f_{m+1} \mrm w_{m+1}) = \hterm(f_{m+1}) \mm w_{m+1}$.
\\
In case $t \not\in \terms(q)$ we get $q +  \alpha_{m+1} \skm f_{m+1}
 \mrm w_{m+1}  \red{}{\myr}{r}{f_{m+1}} q$ and  are done.
\\
In case $t \not\in \terms(p)$ we get $p - \alpha_{m+1} \skm f_{m+1}
  \mrm w_{m+1}  \red{}{\myr}{r}{f_{m+1}} p$.
As $p  - \alpha_{m+1} \skm f_{m+1} \mrm w_{m+1} =
   q + \sum_{j=1}^{m} \alpha_{j} \skm f_j \mrm w_{j}$ the
 induction hypothesis yields 
 $p  - \alpha_{m+1} \skm f_{m+1} \mrm w_{m+1}
 \red{*}{\lr}{r}{F} q$ and hence we are done.
\\ 
Otherwise 
 let $\beta_1 \neq 0$ be the coefficient of $t$ in
 $q +  \alpha_{m+1} \skm f_{m+1} \mrm w_{m+1}$ and
 $\beta_2 \neq 0$ the coefficient of $t$ in $q$.
\\
This gives us a right reduction step

\hspace*{1cm}$q +  \alpha_{m+1} \skm f_{m+1} \mrm w_{m+1} \red{}{\myr}{r}{f_{m+1}}$\\
\hspace*{1cm}$q +  \alpha_{m+1} \skm f_{m+1} \mrm w_{m+1}
                  - \beta_1 \skm \hc(f_{m+1} \mrm w_{m+1})^{-1}
                   \skm  f_{m+1} \mrm w_{m+1} =$ \\
\hspace*{1cm}$q - (\beta_1 \skm \hc(f_{m+1}\mrm w_{m+1})^{-1} -\alpha_{m+1})
             \skm  f_{m+1} \mrm w_{m+1}$

 eliminating the occurrence of $t$ in $q +  \alpha_{m+1} \skm f_{m+1} \mrm w_{m+1}$.
\\
Then obviously 
  $\beta_2 = (\beta_1 \skm \hc(f_{m+1}\mrm w_{m+1})^{-1}
                      -\alpha_{m+1}) \skm \hc(f_{m+1}\mrm w_{m+1})$
  and, therefore, $q  \red{}{\myr}{r}{f_{m+1}}
             q - (\beta_1 \skm \hc(f_{m+1}\mrm w_{m+1})^{-1} -\alpha_{m+1})
             \skm  f_{m+1} \mrm w_{m+1}$, i.e.,
  $q$ and $q +  \alpha_{m+1} \skm f_{m+1} \mrm w_{m+1}$ are joinable.}
\\
\qed
Given an arbitrary set of polynomials, we now intend to enrich the set of
 polynomials used for reduction in order to describe the right
ideal congruence 
by right reduction.
This process will be called \betonen{saturation}\/.\index{saturation}
\begin{definition}\label{def.sat.r}~\\
{\rm
A set of polynomials $F \subseteq \{\alpha \skm p \mrm w \mid \alpha
\in \myk^*, w \in\m \}$ is called
 a \index{saturating set}\betonen{saturating set}\/ for a polynomial $p\in \myk[\m]$,
 if for all $\alpha \in \myk^*$, $w \in \m$,  in case $\alpha \skm p
 \mrm w \neq 0$ then $\alpha \skm p \mrm w
 \red{}{\myr}{r}{F} 0$ holds.
Let $\SAT(p)$ denote the family of all saturating sets for $p$.
\dend
}
\end{definition}
Note that in defining  saturating sets we demand right reducibility to
 zero in {\em one}\/ step in case the multiple is non-zero.
This is done to have some equivalent for the situation $\alpha \skm p \mrm w
\red{{ \leq 1}}{\myr}{s}{p} 0$
  (compare lemma \ref{lem.red}) respectively $\alpha \skm p \mrm w \red{}{\myr}{b}{p} 0$
 in Buchberger's approach.
Furthermore, since $\myk$ is a field it is sufficient to demand $p \mrm w
 \red{{ \leq 1}}{\myr}{r}{F} 0$ for all $w \in \m$.

Let us now proceed by inspecting how
 saturating a polynomial reveals  a  natural connection between strong right
 reduction and right reduction.
\begin{lemma}\label{lem.sconnection}~\\
{\sl
Let $f,g,p$ be some polynomials in $\myk[\m]$ and $S \in\SAT(p)$ a
 saturating set for $p$.
Then
$$f  \red{}{\myr}{s}{p} g \mbox{ if and only if } f  \red{}{\myr}{r}{S}
g.$$
\lemend
}
\end{lemma}
\Ba{}
\begin{enumerate}
\item Suppose we have $f  \red{}{\myr}{s}{p} g$ at a monomial $\alpha \skm t$,
       i.e., $g = f - \alpha \skm \hc(p \mrm w)^{-1} \skm p \mrm w$ for some $w \in \m$
       and $\hterm(p \mrm w)=t$.
      Since $p \mrm w  \red{}{\myr}{r}{S} 0$ there exists a polynomial 
       $p_1 \in S$ such that 
       $p \mrm w =  \beta \skm p_1 \mrm w_1$ for some $\beta \in
       \myk^*$, $w_1 \in \m$.
      Further $t= \hterm(p_1 \mrm w_1) = \hterm(p_1) \mm w_1$ implies
      $f  \red{}{\myr}{r}{p_1 \in S} g$.
\item Suppose $f \red{}{\myr}{r}{p_1 \in S} g$, i.e.,
       $g = f - \alpha \skm p_1 \mrm w_1$
       for some $\alpha \in \myk^*, w_1 \in \m$. 
      Since $p_1 \in S$ we have $\beta \in \myk^*$, $w_2 \in \m$
       such that $p_1= \beta \skm p \mrm w_2$. 
      Further $\hterm(p_1) \mm w_1 = \hterm(p_1 \mrm w_1)= 
       \hterm((p \mrm w_2) \mrm w_1) = \hterm(p \mrm (w_2 \mm w_1))$ implies
       $f  \red{}{\myr}{s}{p} g$.
\\
\qed
\end{enumerate}\renewcommand{\baselinestretch}{1}\small\normalsize
The following lemmata establish that saturation indeed is a
 well-defined concept to repair
 some defects of right reduction, i.e., the right ideal
 generated by $p$ is the same as
 the right ideal generated by a saturating set of $p$ and right
 reduction using a saturating set for a polynomial captures the right
 ideal generated by this polynomial.
\begin{lemma}~\\
{\sl
Let $p$ be a polynomial in $\myk[\m]$ and  $S_1, S_2 \in\SAT(p)$ two
 saturating sets for $p$.
Then
  $$ \red{*}{\lr}{r}{S_1} = \red{*}{\lr}{r}{S_2}.$$
\lemend
}
\end{lemma}
\Ba{}~\\
We restrict ourselves to proving $\red{*}{\lr}{r}{S_1} \subseteq
 \red{*}{\lr}{r}{S_2}$ by induction on $k$ for
  $\red{k}{\lr}{r}{S_1}$ since the case 
              $ \red{*}{\lr}{r}{S_2} \subseteq
              \red{*}{\lr}{r}{S_1}$ is symmetric.
In the base case $k=0$ there is nothing to show.
Thue let us assume $\red{k}{\lr}{r}{S_1}\subseteq\red{*}{\lr}{r}{S_2}$ and 
 $p_0 \red{k}{\lr}{r}{S_1} p_k \red{}{\lr}{r}{S_1} p _{k+1}$.
Without loss of generality we further suppose $p_{k} \red{}{\myr}{r}{q} p_{k+1}$
with $q \in S_1$ (the other case is similar).
      Then $p_{k+1} = p_{k} - \alpha \skm q \mrm w$ for $\alpha \in \myk^*,w \in \m$,
       and since $S_2 $ is a saturating set for $p$ we have 
       a polynomial $q_1 \in
       S_2$ with $q \mrm w \red{}{\myr}{r}{q_1\in S_2} 0$, i.e.,
       $q \mrm w =  \beta \skm q_1 \mrm w_1$, where $\beta \in \myk^*$, $
       w_1 \in \m $ and 
       $p_{k+1} = p_k - (\alpha \skm \beta) \skm q_1 \mrm w_1$.
      Therefore we get $p_k \red{}{\myr}{r}{q_1 \in S_2} p_{k+1}$ and 
       our induction hypothesis yields
       $p_0 \red{*}{\lr}{r}{S_2 }p_k \red{}{\lr}{r}{S_2} p_{k+1}$.
\\
\qed  
\begin{lemma}~\\
{\sl
Let $S \in \SAT(p)$ be a saturating set for a polynomial $p \in
\myk[\m]$.
Then 
 $$\red{*}{\lr}{r}{S} = \;\;\equiv_{\ideal{r}{}(p)} .$$
\lemend
}
\end{lemma}
\Ba{}~\\
This is a consequence of lemma \ref{lem.sconnection}, as we have
 $\red{*}{\lr}{r}{S} = \red{*}{\lr}{s}{p} = \;\;\equiv_{\ideal{r}{}(p)}$.
\qed
The next lemma states that we can in some sense ``simplify''
 saturating sets to ``minimal'' saturating sets such that 
 no polynomial in the set is right
 reducible to zero in one step by the remaining polynomials.
\begin{lemma}~\\
{\sl
Let $p$ be a non-zero polynomial in $\myk[\m]$ and $S \in \SAT(p)$ a saturating
 set for $p$.
\\
If there is a polynomial $q \in S$ such that $q \red{}{\myr}{r}{S
 \backslash \{ q \}} 0$, then $S \backslash \{ q \}$ is again a saturating
 set for $p$.
}
\end{lemma}
\Ba{}~\\
This is an immediate consequence of the fact that $p
\red{}{\myr}{r}{q_1} 0$ and $q_1 \red{}{\myr}{r}{q_2} 0$ imply $p
\red{}{\myr}{r}{q_2} 0$ (compare statement 3 of lemma \ref{lem.red}
which also holds for right reduction).
\auskommentieren{
Let $S$ be a saturating set for the polynomial $p$ and suppose there
 exists a polynomial $q$ in $S$ such that $q \red{}{\myr}{r}{S
 \backslash \{ q \}} 0$.
\\
Then we have to show that for every element $w \in \m$, the polynomial
 $p \mrm w$ is right reducible to zero in one step using $S \backslash \{ q
 \}$ in case $p \mrm w \neq 0$.
\\
Let us assume $p \mrm w \red{}{\myr}{r}{q} 0$, as otherwise we are done
 immediately. 
As $q \red{}{\myr}{r}{S \backslash \{ q \}} 0 $, there exists a
  polynomial $q' \in S \backslash \{ q \}$ such that $q = \alpha' \skm q' \mrm u'$
  for some $\alpha' \in \myk^*$, $u' \in \m$ and $\hterm(q) = \hterm(q') \mm u'$.
Now, $p \mrm w \red{}{\myr}{r}{q \in S} 0$ implies that $p \mrm w =
 \alpha \skm q \mrm u$
 for some $\alpha \in \myk^*$, $u \in \m$ and $\hterm(p \mrm w) =
 \hterm(q \mrm u) = \hterm(q) \mm u$.
Thus 
$$p \mrm w = \alpha \skm q \mrm u = 
             \alpha \skm (\alpha' \skm q' \mrm u') \mrm u = 
             (\alpha \skm \alpha') \skm q' \mrm (u' \mm u)$$
 and
$\hterm(p \mrm w) = \hterm(q \mrm u) = \hterm(q) \mm u =
 \hterm(q' \mrm u') \mm u = (\hterm(q') \mm u') \mm u =
 \hterm(q') \mm (u' \mm u)$
imply $p \mrm w \red{}{\myr}{r}{q'} 0$, and hence $S \backslash \{ q \}$
 is a saturating set for $p$.}
\\
\qed
On the contrary, the property of being a saturating set for a
polynomial is in general destroyed when interreducing it with respect
to right reduction.
\begin{example}~\\
{\rm
Let $\Sigma = \{ a \}$ and $T = \{ a^3 \myr \lambda \}$ be a
presentation of a monoid $\m$ (which is in fact a group).
\\
Looking at the polynomial $p = a^2 + \lambda$ we find that the set $S
= \{ a^2 + \lambda, a+ \lambda \}$ is a saturating set for $p$.
Furthermore, we have $p \red{}{\myr}{r}{a+ \lambda} a^2 + \lambda - (a +
\lambda) \mrm a = a^2 + \lambda - a^2 -a = -a + \lambda
\red{}{\myr}{r}{a+ \lambda}  -a + \lambda  -  ( -1) \skm (a + \lambda) =
2 \skm \lambda$, but the set $S' = \{  a + \lambda, 2 \skm \lambda \}$
is no longer a saturating set for $p$.
\exaend
}
\end{example}
We move on now to examine how saturating sets can be constructed.
Obviously the set $\{p \mrm w \mid w \in\m \}$ itself is a trivial
 saturating set for the polynomial $p$.
Since it is recursively enumerable,
 it can be used to specify
 enumerating procedures for saturating sets.
\auskommentieren{To begin with, let us give a naive procedure for saturating sets.

\procedure{Saturation 1}%
{\vspace{-4mm}\begin{tabbing}
XXXXX\=XXXX \kill
\removelastskip
{\bf Given:} \> A polynomial $p \in \myk[\m]$. \\
{\bf Find:} \>  A (generally infinite) set $S \in\SAT(p)$.
\end{tabbing}
\vspace{-7mm}
\begin{tabbing}
XX\=XX\= XXX \= XXX \=\kill
$S$ := $\{ p \}$; \\
{\bf for all} $q \in \{p \mrm w \mid w \in\m \}$ {\bf do} \\
\>      {\bf if} \> $q \nred{}{\myr}{r}{S} 0$ \\
\>               \>{\bf then}\>      $S$ := $S \cup \{ q \}$; \\
{\bf endfor}
\end{tabbing}}

This procedure enumerates a saturating set for $p$ depending on
 the enumeration of the set
 $\{p \mrm w \mid w \in\m \}$ used in the {\bf for all} loop.}
Because in general such enumerations will not terminate we have to
 try to find criteria to decide
 whether a set is a saturating set for a polynomial without
 having to check all right multiples.
We are interested in a systematic way to distinguish what terms can be
brought to head position.
Since in order to bring another term to head position in multiplying
the head term a cancellation (in computing the product using the
presenting semi-Thue system) has to take part, one  first idea is to 
overlap with rules in the presentation of the
monoid instead of multiplying with all monoid elements.

\procedure{Saturation 1\protect{\label{saturation2}}}%
{\vspace{-4mm}\begin{tabbing}
XXXXX\=XXXX \kill
\removelastskip
{\bf Given:} \> A polynomial $p \in \myk[\m]$ and              
              $(\Sigma, T)$ a convergent  presentation of $\m$. \\
{\bf Find:} \>  A (generally infinite) set $S \in\SAT(p)$.
\end{tabbing}
\vspace{-7mm}
\begin{tabbing}
XX\=XX\= XX \= XXXX \= XX \= XXXX \=\kill
$S$ := $\{ p \}$; \\
$H_0$ := $\{ p \}$; \\
$i$ := $0$; \\
{\bf while} $H_i \neq \emptyset$ {\bf do} \\
\> $q$ := {\rm remove}$(H_i)$;\\
\>{\rm\kommentar \% Remove an element using a fair
  strategy} \\
\> $t$ := $\hterm(q)$; \\
\> {\bf for all} $w \in C(t) = \{ w \in \Sigma^* \mid tw \id t_{1}t_{2}w \id t_{1}l, t_2 \neq \lambda$
                            {\rm for some} $(l, r) \in T \} $ {\bf do}
                            \\
\>\>{\rm\kommentar \% $C(t)$ contains special overlaps between 
    $t$ and  left hand sides of rules in $T$} \\
\> \> $q'$ := $q \mrm w$; \\
\> \> {\bf if} \> $q' \not\in S$; \\
\> \>          \>{\bf then} \> {\bf if} \> $q' \nred{}{\myr}{r}{S} 0$ \\
\> \>          \>           \>          \> {\bf then} \> $S$ := $S \cup \{q' \}$; \\
\> \>          \>           \>          \>            \> $H_i$ := $H_i \cup \{q' \}$; \\
\> \>          \>           \>          \>{\bf else} \> $H_i$ := $H_i \cup \{q' \}$; \\
\> \>          \>           \>{\bf endif}\\
\> \> {\bf endif}\\
\> $H_{i+1}$ := $H_{i}$;\\
\> $i$ := $i+1$;\\
\> {\bf endfor} \\
{\bf endwhile} 
\end{tabbing}}

This procedure  gives us an enumeration of a saturating set for a
polynomial.
This enumeration terminates in case for some index $i$ the set $H_i$ becomes empty.
This can only arise if  for any polynomial $q \in H_i$
 either the set $C(\hterm(q))$ is empty or for all 
 $w \in C(\hterm(q))$, $ q \mrm w \in S$ holds.
Since these requirements are strong, it will be important to find
 additional sufficient conditions when such an
 enumeration can be stopped.
The following example illustrates why in case  $q' \nred{}{\myr}{r}{S}
0$ in procedure {\sc Saturation 1} we have to consider $q'$ for further computations, i.e., add it to
our set $H$.
\begin{example}~\\
{\rm
Let $\Sigma = \{ a,b \}$ and $T = \{ ab \myr \lambda \}$ be a
presentation of a monoid $\m$ with a length-lexicographical ordering
induced by $a \succ b$.
\\
Then in modifying procedure {\sc Saturation 1} to omit adding $q'$ to
$H_i$ in case we already have $q' \red{}{\myr}{r}{S} 0$ we would get the
following computation on input $a^2 + a$:
We start with
$S: = \{ a^2 + a \}$ and $H_0 := \{ a^2 + a \}$.
Then $q := a^2 + a$ is removed from $H_0$ and we find $C(a^2) = \{ b \}$
resulting in $q' = (a^2 + a) \mrm b = a + \lambda$.
Since $a+
\lambda \red{}{\myr}{r}{S} 0$ the procedure then halts with output $S =
\{ a^2 + a \}$.
But the polynomial $(a^2 + a) \mrm b^2 = \lambda + \underline{b}$ is
not right reducible by $S$.
\\
On the other hand, procedure {\sc Saturation 1} terminates and computes
 the set $S = \{ a^2+a,b+\lambda \}$, which is a saturating set for $a^2+a$.
\exaend
}
\end{example}
\begin{lemma}~\\
{\sl
Procedure {\sc Saturation 1} generates a saturating set for a polynomial.
}
\end{lemma}
\Ba{}~\\
Let $p$ be a non-zero polynomial in $\myk[\m]$ and $S$ the set
generated   by procedure {\sc Saturation 1} on input $p$ and
$(\Sigma, T)$.
We have to show that for all $w \in \m$ we have $p \mrm w
\red{}{\myr}{r}{S} 0$ in case $p \mrm w \neq 0$.
Since $p \in H_0$ this follows immediately by showing that for all $q \in 
\bigcup_{i \geq 0} H_i$ and all $w
\in \m$ with $q \mrm w \neq 0$ we have $q \mrm w \red{}{\myr}{r}{S} 0$.
Notice that for all $q \in H_i$ we have $q \red{}{\myr}{r}{S} 0$.
Now let us assume that our claim is not true.
Then we can choose a  non-zero counter-example
 $q \mrm  w$, where $\hterm(q)w$ is minimal
 (according to the ordering $\succeq_T$ on $\Sigma^*$) 
 and $q \mrm w \nred{}{\myr}{r}{S} 0$.
Then $\hterm(q)w$ must be $T$-reducible, as otherwise
 $q \mrm w  \red{}{\myr}{r}{q} 0$ with $\hterm(q \mrm w) \id \hterm(q)w$
 and, as $q \in 
\bigcup_{i \geq 0} H_i$, then  $q \red{}{\myr}{r}{s \in S} 0$ implies $q \mrm w \red{}{\myr}{r}{s \in S} 0$.
Let $\hterm(q)w \id t_1t_2w_1w_2$ where
 $\hterm(q) \id t_1t_2, t_2 \neq \lambda, w \id w_1w_2$ and
 $l \id t_2w_1$ for some $(l,r) \in T$.
Furthermore, $w_1 \in \m$ as it is a prefix of $w \in \m$.
As we use a fair strategy 
 to remove elements from the sets $H_i$, $q$ and $C(\hterm(q))$ must
 be considered and we have $w_1 \in C(\hterm(q))$ by the definition of this set.
Since we have $q \mrm w_1 \in \bigcup_{i \geq 0} H_i$ by construction, 
       $\hterm(q)w \id \hterm(q)w_1w_2 \succ \hterm(q \mrm w_1)w_2$  
       contradicts our assumption that $q \mrm w$ was a minimal
       counter example.
\\
\qed
Another  idea to compute saturating sets might be to gain more insight by taking a more constructive
 look at their definition.
\begin{definition}\label{def.sat}~\\
{\rm 
Given a polynomial  $p$, for each term $t \in \terms(p)$
 let $X_{t} = \{ w \in \m \mid \hterm(p \mrm w) = t \mm w \}$,
 i.e., the set of all elements, which put $t$ into head
 position\footnote{Note that if $\m$ is not right-cancellative
 one $w$ may belong to different sets.}.
Further let $Y_{t} = \{p \mrm w  \mid w \in X_{t} \}$.
By choosing  sets $B_{t} \subseteq Y_{t}$ such that
 for all polynomials $q \in Y_{t}$ we have $q  \red{}{\myr}{r}{B_{t}}
 0$, we get a saturating set $\bigcup_{t \in \terms(p)} B_{t} \in \SAT(p)$.
\dend
}
\end{definition}
This definition does not specify how to choose the sets $B_{t}$, but
 setting  $B_{t} = Y_{t}$ we always get the trivial saturating set 
 $\{p \mrm w \mid w \in\m \} \backslash \{ 0 \}$.
Of course $Y_{\hterm(p)}$ must at least contain $p$,
 but all other $Y_{t}$ can be empty.
In case the multiplication on $\m$ is monotone,  we get
 $Y_{\hterm(p)}= \{p \mrm w \mid w \in \m \}$ and  $Y_{t} = \emptyset$
 for $t \in \terms(p) \backslash \{\hterm(p)\}$.
Then obviously the set
 $B_{\hterm(p)}= \{p \}$ is a finite saturating set for $p$.

Unfortunately, one cannot hope for a terminating saturating procedure because
 finite saturating sets need not exist for certain monoid
 presentations, as the following examples reveal.
\begin{example}\label{exa.kippen}~\\
{\rm
Let $\Sigma = \{ a, b, c, d, e, f \}$ and  
 $T= \{ abc \myr ba, bad \myr e , fbc \myr bf \}$ be a presentation\footnote{Note
 that $( \Sigma , T)$ is a convergent length-reducing presentation  and
 the monoid is cancellative.} of
 a monoid $\m$ 
 with a
 length-lexicographical ordering induced by $ a \succ b \succ c \succ d \succ e \succ f$.
\\
Looking at the polynomial $p = a + f \in \q[\m]$ we find that with respect
 to this presentation we have the set
 $X_f =\{ (bc)^idw \mid i \in \n, w \in\m \}$,
 and it can be seen that the set 
 $Y_f = \{ b^{i+1}fdw + b^iew \mid i \in \n, w \in \m \}$
 has no finite basis.
Since if there were a finite basis $B_f$ for $Y_f$, we
 could choose $k \in \n$ such that
 $ b^{k+1}fd + b^ke \not\in B_f$.
But then we get $ b^{k+1}fd + b^ke \nred{}{\myr}{r}{B_f} 0$ as
 $ b^{i+1}fdw \mm {\rm x } =  b^{k+1}fd $ has no solution in
 $\m$ unless ${\rm x }= \lambda$ and $i=k$.
Note that  every saturating set 
 $S \in\SAT(p)$ must right reduce all the
 elements of the set $Y_f$ to zero in one step.
However, if we change our precedence on $\Sigma$ slightly to
 $f \succ a \succ b \succ c \succ d \succ e$, then the set $\{ f + a \}$
 is a saturating set for $p$.
The phenomenon occurring  is depicted below; the respective head
 terms of the polynomials are underlined.
\exaend
}
\end{example}
\begin{diagram}[size=2em]
\underline{a}+f       &                   &      \\
\dLine_{bc}     &                   &      \\
\underline{ba}+bf     & \rLine_{d}  & e+\underline{bfd}     \\
\dLine_{bc}     &                   &      \\
\underline{b^2a}+b^2f & \rLine_{d}  & be+\underline{b^2fd}    \\
\dLine_{bc}     &        &      \\
\vdots \\
\end{diagram}
%
\begin{example}\label{exa.no.finite.sat.set}~\\
{\rm
Let $\Sigma = \{ a, b, c, d, e, f \}$ and  
 $T= \{ abc \myr ba, bad \myr fe , fbc \myr bf, bfd \myr ae \}$ be a presentation of
 a monoid $\m$ with a
 length-lexicographical ordering induced by an arbitrary precedence
 $\succ$ on 
 $\Sigma$.
\\
Then the polynomial $p = a + f \in \q[\m]$ will allow no finite saturating
 sets independent of the chosen precedence:
We have to distinguish two cases according to the relation of the
 two letters $a$ and $f$:
\begin{enumerate}
\item Case $a \succ f$ (depicted in figure 1 below): 
      Then we get the set $X_f =\{ (bc)^idw \mid i \in \n, w \in\m \}$,
       and the corresponding set $Y_f = \{ b^iaew + b^ifew \mid i \in \n, w \in \m \}$
       has no finite basis.
      Since if there was a finite basis $B_f$,
       we could choose $k \in \n$ such that
       $b^kae + b^kfe \not\in B_f$ for $Y_f$.
      But then we get $ b^kae + b^kfe \nred{}{\myr}{r}{B_f} 0$ as
       $ b^iae \mm {\rm x } =  b^kae $ has no solution in
       $\m$ unless ${\rm x }= \lambda$ and $i=k$.
\item Case $f \succ a$ (depicted in figure 2 below): 
      Then we get the set $X_a =\{ (bc)^idw \mid i \in \n, w \in\m \}$,
       and the corresponding set $Y_a = \{ b^ifew + b^iaew \mid i \in \n, w \in \m \}$
       likewise has no finite basis.
      This can be shown as above.
\exaend
\end{enumerate}
}
\end{example}
\begin{diagram}[size=2em]
\underline{a}+f       &            &              & \phantom{XXXXXXXXXXXXX}&  
\underline{f}+a       &            &  \\
\dLine_{bc}     &                   &        & &  
\dLine_{bc}     &            &  \\
\underline{ba}+bf     & \rLine_{d}  & fe+\underline{ae}     & &  
\underline{bf}+ba     &  \rLine_{d} & ae+\underline{fe}   \\
\dLine_{bc}     &                   &     & &  
\dLine_{bc}     &       &     \\
\underline{b^2a}+b^2f & \rLine_{d}  & bfe+\underline{bae}    & &  
\underline{b^2f}+b^2a & \rLine_{d}  & bae+\underline{bfe}   \\
\dLine_{bc}     &        &       & &
\dLine_{bc}     &        &         \\
\vdots       &   &\mbox{\rm\small Figure 1}              & &
\vdots        & &\mbox{\rm\small Figure 2} \\
\end{diagram}

On the other hand in a group ring finite saturating sets always exist.
\begin{lemma}~\\
{\sl
For any polynomial $p$ in a group ring
 $\myk[\g]$ there exists a saturating set
 $S \in \SAT(p)$ containing at most $|\terms(p)|$ elements.
}
\end{lemma}
\Ba{}~\\
To show our claim let us review the more constructive definition
 of saturating sets as given in definition \ref{def.sat}.
\\
For a polynomial $p$, let us take a closer look at the sets
 $X_{t} = \{ w \in \m \mid \hterm(p \mrm w) = t \mm w \}$ and
 $Y_{t} = \{p \mrm w  \mid w \in X_{t} \}$ for $t \in \terms(p)$.
Then in case $Y_{t} \neq \emptyset$,
 there exists a  polynomial $q = p \mrm w$ in $Y_{t}$.
It remains to show that every polynomial in $Y_t$ then is right
reducible to zero in one step using $q$.
As $w$ is a group element, we know $p = q \mrm \inv{w}$ and for every
other polynomial $p \mrm u \in Y_t$ we find
 $\hterm(q \mrm (\inv{w} \mm u)) = \hterm((q \mrm \inv{w}) \mrm u) 
 = \hterm(p \mrm u) = t \mm u = t \mm (w \mm \inv{w}) \mm u =
 \hterm(q) \mm (\inv{w} \mm u)$ and hence $p \mrm u \red{}{\myr}{r}{q} 0$.
\auskommentieren{
It remains to show that the set $\{ q \}$ is a basis for $Y_{t}$.
Let $p \mrm u$ be an arbitrary element in $Y_{t}$.
Then we know that $\hterm(p \mrm w) = t \mm w$ and 
 $\hterm(p \mrm u) = t \mm u$.
Further, as $w$ is an element of a group, we can conclude
 $p = q \mrm \inv{w}$ and thus
 $p \mrm u = ( q \mrm \inv{w}) \mrm u = q \mrm (\inv{w} \mm u)$.
It remains to show that $\hterm(q \mrm (\inv{w} \mm u)) = \hterm(q)
\mm (\inv{w} \mm u)$.
This follows as
$\hterm(q) \mm (\inv{w} \mm u) = (t \mm w) \mm (\inv{w} \mm u) =
t \mm u = \hterm(p \mrm u) = \hterm(q \mrm (\inv{w} \mm u))$.
}
\\
\qed
The proof of this lemma might give rise to another idea how to compute
 saturating sets using the constructive approach in definition  \ref{def.sat}.

\pagebreak
\procedure{Saturation 2\protect{\label{saturation}}}%
{\vspace{-4mm}\begin{tabbing}
XXXXX\=XXXX \kill
\removelastskip
{\bf Given:} \> A polynomial $p \in \myk[\m]$. \\
{\bf Find:} \> A set $S \in \SAT(p)$.
\end{tabbing}
\vspace{-7mm}
\begin{tabbing}
XX\=XX\= XXXX \= XXX \=\kill
$S$ := $\emptyset$; \\
{\bf for all} $t \in \terms(p)$ {\bf do} \\
\> {\bf if} \> $X_{t} \neq \emptyset$\\
\>          \> {\bf then} \> $S$ := $S \cup \{ \mbox{ a suitable basis of } Y_{t}\; \}$; \\
{\bf endfor} 
\end{tabbing}}

Unfortunately this procedure is not effective.
First of all, although each of the sets $X_{t}$  is recursive, 
 one cannot decide, whether it is empty or not.
On the other hand of course even if one knows that $X_{t} \neq
 \emptyset$, a ``suitable'' finite basis for the set $Y_{t}$ need not
 exist (compare example \ref{exa.no.finite.sat.set}).
We will see later on how this idea is used to compute saturating sets
in special group rings (compare chapter \ref{chapter.grouprings}).

The main reason why this procedure cannot be effective is that the
following uniform problem is not solvable:

\begin{tabbing}
XXXXXXXX\=XXXXX \kill
{\bf Given:} \>  A polynomial $p \in \myk[\m]$, a term $t \in \terms(p)$, and \\
             \> $(\Sigma, T)$ a convergent semi-Thue system presenting $\m$. \\
{\bf Question:} \> $X_t \neq \emptyset$ ?
\end{tabbing}

\begin{lemma}\label{lem.andrea}~\\
{\sl
There exists a monoid $\m$ such that the above problem is undecidable.
}
\end{lemma}
\Ba{}~\\
We give a finite, convergent semi-Thue system
 presenting a monoid $\m$
 such that for
 a polynomial $p=t_1-t_2$, $t_1,t_2 \in \m$ in the corresponding
 monoid ring $\q[\m]$ in general it is undecidable whether $X_{t_2} \neq \emptyset$.
This proof will use the following lemma stated in \cite{Sa90}:

{\sl
For every primitive recursive function ${\rm\bf f}: \n^2 \myr \n$ there
 exists a finite convergent interreduced semi-Thue
 system $T_{\rm\bf f}$ over an alphabet $\Sigma_{\rm\bf f}$
 with $\{a,b,g,c,v,e \} \subseteq \Sigma_{\rm\bf f}$ and a symbol
 $f \in \Sigma_{\rm\bf f}$ such that the following holds:
\begin{enumerate}
\item We have $fb^{n_1}ab^{n_2}age \red{*}{\myr}{}{T_{\rm\bf f}}
       b^{n_1}ab^{n_2}agc^{{\rm\bf f}(n_1,n_2)}ve$ and
      $b^{n_1}ab^{n_2}agc^{{\rm\bf f}(n_1,n_2)}ve \in \irr(T_{\rm\bf f})$ for $n_1,n_2 \in\n$.
\item There exists a precedence $\succ \;\subseteq \Sigma_{\rm\bf f} \times \Sigma_{\rm\bf f}$ satisfying
      $a \succ b \succ g \succ c \succ v \succ e$
      such that for the induced syllable ordering with status right 
      $\syll \;\subseteq \Sigma_{\rm\bf f}^* \times \Sigma_{\rm\bf f}^*$ the
      following holds: $\red{+}{\myr}{}{T_{\rm\bf f}} \subseteq\; \syll$.
\item For all $(l,r) \in T_{\rm\bf f}$ we have  $|l|=2$ and the first letter of
  $l$ is not in $\{ a,b,g,c,v,e \}$.
\end{enumerate}
}
Reviewing the proof of this lemma, the precedence on $\Sigma_{\rm\bf
  f}$ can be chosen such that $f$ is maximal among all letters.
Thus for the primitive recursive functions
\[ {\rm\bf f}_1({\rm x},{\rm y}) = \left\{ \begin{array}{r@{\quad\quad}l}
                      0 & \exists {\rm z} \leq {\rm y} \mbox{ such
                        that the computation of 
                        $\Phi_{\rm x}({\rm x})$ halts within} \\
                      \mbox{ }  & \mbox{ z steps} \\
                      2 & \mbox{otherwise}
                      \end{array} \right. \]
where $\Phi$ is a universal function
and ${\rm\bf f}_2({\rm x},{\rm y}) =1$
there exist finite convergent interreduced semi-Thue systems $(\Sigma_{{\rm f}_1},
 T_{{\rm f}_1})$ and $(\Sigma_{{\rm f}_2}, T_{{\rm f}_2})$ with
 precedences $\succ_1, \succ_2$ as
 described above and without loss of generality we can assume that 
 $\Sigma_{{\rm f}_1} \cap \Sigma_{{\rm f}_2} = \{ a,b,g,c,v,e \}$.
Let $\Sigma = \Sigma_{{\rm f}_1} \cup \Sigma_{{\rm f}_2}$ and
  $T = T_{{\rm f}_1} \cup T_{{\rm f}_2}$ with precedence $\succ$ such
  that $\succ$ extends $\succ_1$ and $\succ_2$ and further $f_1 \succ
  f_2 \succ {\rm x}$ for all ${\rm x} \in \Sigma \backslash \{ f_1,f_2
  \}$.
If we slightly modify $\Sigma$ and $T$ to
\begin{eqnarray}
\Sigma' & = & \Sigma \cup \{ L_1, L_2,Q_1, \ldots,Q_7, A,B,G,E \}
             \nonumber \\
        &   & \mbox{without loss of generality we assume } \Sigma \cap
        \{L_1, L_2,Q_1, \ldots,Q_7, A,B,G,E \} = \emptyset \nonumber \\
T'      & = & T \;\;\;\;\cup  \nonumber \\
        &   & \{ {\sf(1)} \;L_2b \myr BL_2, 
                 {\sf(2)} \;L_2a \myr AQ_1,
                 {\sf(3)} \; Q_1b \myr BQ_1,
                 {\sf(4)} \;  Q_1a \myr AQ_2,\nonumber \\
        &   & \phantom{\{} {\sf(5)} \; Q_2g \myr GQ_3,
                   {\sf(6)} \; Q_3e \myr Q_4e,
                   {\sf(7)} \; GQ_4 \myr Q_5g,
                   {\sf(8)} \; AQ_5 \myr Q_6a, \nonumber \\
        &   & \phantom{\{} {\sf (9)} \; BQ_6 \myr Q_6b, 
                   {\sf(10)} \; AQ_6 \myr Q_7a, 
                   {\sf(11)} \; BQ_7 \myr Q_7b,
                   {\sf(12)} \; L_1Q_7 \myr \lambda \} \nonumber 
\end{eqnarray}
with precedence $\succ'$ any total extension of $\succ$ satisfying
 $L_1 \succ' L_2 \succ' Q_1 \succ' Q_2 \succ' Q_3 \succ' Q_4 \succ' Q_5
 \succ' G \succ' B \succ' A \succ' Q_6 \succ' Q_7$,
then $(\Sigma',T')$ is a finite convergent interreduced semi-Thue
 system such
 that for all $(l,r) \in T'$ we have  $|l|=2$ and the first letter of
  $l$ is not in the set $\{ a,b,g,c,v,e \}$.
\\
Additionally we can state that 
 for $w = b^nab^magew'$, where $w' \in \irr(T')$, we get 
 $f_iL_1L_2w \red{*}{\myr}{}{T'} b^{n}ab^{m}agc^{{\rm\bf f}_i(n,m)}vew'$ and
 this is irreducible.
\\
Furthermore, for the elements $w \in \irr(T')$ which do not have a prefix of the
 form $b^nab^mage$ for some $n,m \in \n$, we get  
 $L_1L_2w \red{*}{\myr}{}{T'} L_1\tilde{w}$ where $L_1\tilde{w}$ is
 irreducible and especially $\tilde{w}$ does not start with the letter $Q_7$.
To prove these reduction sequences, we will use that
 $L_2b^n \red{n}{\myr}{}{(1)} B^nL_2$, $Q_1b^n \red{n}{\myr}{}{(3)} B^nQ_1$,
 $B^nQ_6 \red{n}{\myr}{}{(9)} Q_6b^n$ and $B^nQ_7 \red{n}{\myr}{}{(11)}
 Q_7b^n$ for $n \in \n$.
\\
In case $w = b^nab^magew'$ we have 
\begin{eqnarray}
f_iL_1\underline{L_2b^n}ab^magew' & \red{n}{\myr}{}{(1)} & f_iL_1B^n\underline{L_2a}b^magew' \nonumber \\
   & \red{}{\myr}{}{(2)}  & f_iL_1B^nA\underline{Q_1b^m}agew'   \nonumber \\  
   & \red{m}{\myr}{}{(3)} & f_iL_1B^nAB^m\underline{Q_1a}gew'   \nonumber \\  
   & \red{}{\myr}{}{(4)}  & f_iL_1B^nAB^mA\underline{Q_2g}ew'   \nonumber \\ 
   & \red{}{\myr}{}{(5)}  & f_iL_1B^nAB^mAG\underline{Q_3e}w'   \nonumber \\  
   & \red{}{\myr}{}{(6)}  & f_iL_1B^nAB^mA\underline{GQ_4}ew'   \nonumber \\ 
   & \red{}{\myr}{}{(7)}  & f_iL_1B^nAB^m\underline{AQ_5}gew'   \nonumber \\ 
   & \red{}{\myr}{}{(8)}  & f_iL_1B^nA\underline{B^mQ_6}agew'   \nonumber \\ 
   & \red{m}{\myr}{}{(9)} & f_iL_1B^n\underline{AQ_6}b^magew'   \nonumber \\ 
   & \red{}{\myr}{}{(10)}  & f_iL_1\underline{B^nQ_7}ab^magew'   \nonumber \\ 
   & \red{n}{\myr}{}{(11)} & f_i\underline{L_1Q_7}b^nab^magew'   \nonumber \\ 
   & \red{}{\myr}{}{(12)}  & \underline{f_ib^nab^mage}w'         \nonumber \\ 
   & \red{*}{\myr}{}{T_{{\rm\bf f}_i}} & b^nab^magc^{{\rm\bf f}_i(n,m)}vew' \nonumber 
\end{eqnarray}
and are done.
\\
To show that $L_1L_2w$ reduces to a normal form $L_1\tilde{w}$ where
$\tilde{w}$ does not start with the letter $Q_7$ in case $w$ does not
have a prefix of the form $b^nab^mage$ for some $n,m \in \n$, we use
the following fact:
For  $w \in \irr(T')$ there exists an element $w' \in \irr(T')$
 such that  we have $L_2 w \red{*}{\myr}{}{T'} Q_7w'$ if and only if 
 $w=b^nab^magew'$  for some $n,m \in \n$.
The ``if'' part follows directly as above.
It remains to show that when reducing $L_2w$ for some $w \in
\irr(T')$\footnote{It is sufficient to look at irreducible elements,
  as the semi-Thue system is convergent.}
 we can only reach a word starting with the
 letter $Q_7$  in case
 $w$ has a prefix of the form $b^nab^mage$ for some $n,m \in \n$.
Since $w \in \irr(T')$, the word $L_2w$ is only reducible by $T'$ in
case one of the rules $L_2b \myr BL_2$ or $L_2a \myr AQ_1$ is applicable.
Hence we  get a reduction sequence $L_2w \red{n}{\myr}{}{(1)} B^nL_2v_1$
 in case $w \id b^nv_1$ which terminates as soon as
  no more letters $b$ occur as a prefix of the remaining part $v_1$ of $w$ and
  this reduction sequence is unique.
Since $B^n$ and $v_1$ are irreducible, the word $B^nL_2v_1$ has to be
 reducible by $L_2a \myr AQ_1$,
 as if $v_1$ would start with a letter different from $a$ (and
 different from $b$ by assumption) we would end up with an irreducible
 word $B^nL_2v_1$ that cannot be further reduced to a word starting
 with the letter $Q_7$.
Hence, without loss of generality we can assume,
 $$L_2w \red{n}{\myr}{}{(1)} B^nL_2v_1 \red{}{\myr}{}{(2)} B^nAQ_1v_2$$ and $v_1 \id av_2$.
Since we have not reached a word starting with $Q_7$, again
 we can assume that $v_2$ starts with a prefix $b^ma$ for some $m
 \in \n$, as otherwise we would be stuck with an irreducible word
 starting with $B^nA$.
This gives us the unique reduction sequences
$$B^nAQ_1v_2 \red{m}{\myr}{}{(3)} B^nAB^mQ_1v'_2
   \red{}{\myr}{}{(4)} B^nAB^mAQ_2v_3$$ and $v_2 \id b^mav_3$, $v'_2 \id av_3$. 
Again, to keep the word $B^nAB^mAQ_2v_3$ reducible by rules in $T'$,
 we find that $ge$ must be a prefix of $v_3$, i.e., $v_3 \id gv'_3
 \id gev_4$, giving us the only
 possible reductions
 $$B^nAB^mAQ_2v_3 \red{}{\myr}{}{(5)} B^nAB^mAGQ_3v'_3
  \red{}{\myr}{}{(6)} B^nAB^mAGQ_4ev_4.$$
Combining these informations we find
 $w \id b^nv_1 \id b^nav_2 \id b^nab^mav_3 \id b^nab^magev_4$
 and $L_2w \red{*}{\myr}{}{T'} B^nAB^mAGQ_4ev_4$.
\\
Now we can proceed with the proof of our initial claim:
\\
For  $n \in \n$ let $p_n =f_1L_1L_2b^na -f_2L_1L_2b^na$ denote a
 polynomial in the corresponding monoid ring over $\q$ where $\m$ is
 presented by
 the semi-Thue system $(\Sigma',T')$.
Then we get
$$X_{f_2L_1L_2b^na} \neq \emptyset $$
\begin{center} if and only if \end{center}
$$\exists w \in \Sigma'^*: (f_2L_1L_2b^naw)\nf{T'}
   \syll (f_1L_1L_2b^naw)\nf{T'}$$ 
By the definition of the syllable ordering, the construction of $T'$
and the  resulting normal forms we know, that only a word
starting with a prefix $b^mage$ can cause the letters $f_1$
respectively $f_2$ to be affected and this is the only possibility to
cause $X_{f_2L_1L_2b^na}$ to be non-empty.
Hence we can conclude
$$\exists w \in \Sigma'^*: (f_2L_1L_2b^naw)\nf{T'}
   \syll (f_1L_1L_2b^naw)\nf{T'}$$ 
\begin{center} if and only if \end{center}
$$\exists m \in\n, w' \in\Sigma'^* \: (f_2L_1L_2b^{n}ab^{m}agew')\nf{T'} \syll
(f_1L_1L_2b^{n}ab^{m}agew')\nf{T'}$$
\begin{center} if and only if \end{center}
$$\exists m \in\n, w' \in\Sigma'^*: b^{n}ab^{m}agc^{{\rm\bf f}_2(n,m)}vew' \syll
     b^{n}ab^{m}agc^{{\rm\bf f}_1(n,m)}vew'$$
\begin{center} if and only if \end{center}
$$\exists m \in\n : {\rm\bf f}_2(n,m) > {\rm\bf f}_1(n,m)$$

Remembering the definitions of ${\rm\bf f}_1$ and ${\rm\bf f}_2$ we
get for all ${\rm x} \in \n$: 
$$\exists {\rm y} \in\n \mbox{ with } {\rm\bf f}_2({\rm x},{\rm y}) > {\rm\bf f}_1({\rm x},{\rm y})
 \mbox{ if and only if } \Phi_{\rm x}({\rm x}) \mbox{ is defined}$$
and this would solve the halting problem for the universal function
 $\Phi$.
\\
\qed
Nevertheless, we will later on see how  saturating sets in special
classes of monoids and groups can be computed by using additional
information on the respective structure.

The following definition carries the advantages of saturating sets for one
polynomial on to sets of polynomials.
\begin{definition}\label{def.sat.set}~\\
{\rm
A set  $F$ of polynomials in $\myk[\m]$ is called \index{saturated set}\betonen{saturated},
 if  $\alpha \skm f \mrm w \red{}{\myr}{r}{F} 0$ holds for
 all $f \in F$ and all  $\alpha \in \myk^*$, $w \in \m$ such that
 $\alpha \skm f \mrm w \neq 0$.
\dend
}
\end{definition}
Note that saturating sets for a polynomial $p$ are saturated.
Moreover, for a set of polynomials $F$ every
 union $S = \bigcup_{f \in F} S_f$ of saturating sets $S_f \in \SAT(f)$ is a saturated
 set. 
The next lemma shows how saturated sets allow special representations
 of the elements belonging to
 the right ideal they generate.
\begin{lemma}\label{lem.prop1}~\\
{\sl
 Let $F$ be a saturated set of polynomials in $\myk[\m]$.
 Then every non-zero polynomial $g \in \ideal{r}{}(F)\backslash \{ 0 \}$
 has a representation\footnote{Note
    that such representations are in some sense ``stable''
    but no stable standard representations as the head term of $g$ need
    not be a bound for the terms involved.
   Especially saturated sets need not be stable standard bases.
   }
 of the form $g =  \sum_{i=1}^{k} \alpha_i \skm f_i \mrm w_i$,
 where $\alpha_i \in \myk^*, f_i \in F, w_i \in \m$, and
 $\hterm(f_i \mrm w_i) = \hterm(f_i) \mm w_i$.
\lemend
}
\end{lemma}
\Ba{}~\\
This follows immediately from definition \ref{def.sat.set}, as
 in case $\hterm(f_i \mrm w_i) \neq \hterm(f_i) \mm w_i$ for some
 $f_i$ in our representation of $g$ assuming that $f_i \mrm w_i \neq 0$,
 we know $f_i \mrm w_i \red{}{\myr}{r}{F} 0$, i.e.,
 $f_i \mrm w_i = \beta_i \skm f'_i \mrm w'_i$ for some
 $\beta_i \in \myk^*, f'_i \in F, w'_i \in\m$ and $\hterm(f'_i \mrm
 w'_i)=\hterm(f'_i)\mm w'_i$.
Thus we can substitute $f_i \mrm w_i$ by
 $\beta_i \skm f'_i \mrm w'_i$ in the given representation of $g$.
\\
\qed
Similar to lemma \ref{lem.sconnection} the following statement holds:
\begin{lemma}\label{lem.sconnection2}~\\
{\sl
For $f,g$ some polynomials in $\myk[\m]$ and $F$ a saturated set with $p \in F$,
 $f \red{}{\myr}{s}{p} g$ implies $f \red{}{\myr}{r}{F} g$.
\ohnebeweis
}
\end{lemma}

Now we are able to show that we can simulate $\red{*}{\lr}{s}{}$ with
$\red{*}{\lr}{r}{}$ and capture the right ideal congruence by using saturated sets.
\begin{theorem}\label{theo.congruence}~\\
{\sl
Let $F$ be a saturated set of polynomials and  $p,q$ some polynomials in $\myk[\m]$.
Then
  $$p \red{*}{\lr}{r}{F} q \mbox{ if and only if } p - q \in \ideal{r}{}(F).$$
\theoend
}
\end{theorem}
\Ba{}~\\
This is an immediate consequence of lemma  \ref{lem.strong.congruence} and
 lemma \ref{lem.sconnection2}.
\\
\qed
Similar to definition \ref{def.gb} we can define  Gr\"obner bases with
respect to $\red{}{\myr}{r}{}$.
\begin{definition}\label{def.rgb}~\\
{\rm
A  set $G \subseteq \myk[\m]$ is called a \betonen{Gr\"obner basis}\/
 with respect to
 the reduction $\red{}{\myr}{r}{}$ or a \index{right!Gr\"obner basis}
 \index{stable!Gr\"obner basis}\index{Gr\"obner basis!right}
 \index{Gr\"obner basis!stable}
 \betonen{right (or stable) Gr\"obner basis}, if
\begin{enumerate}
\item[(i)] $\red{*}{\lr}{r}{G} = \;\; \equiv_{\ideal{r}{}(G)}$, and
\item[(ii)] $\red{}{\myr}{r}{G}$ is confluent.
\dend
\end{enumerate}
}
\end{definition}
Notice that saturating sets for a polynomial $p$ satisfy statement (i) of this definition, 
 but  in general need not be 
 Gr\"obner bases of $\ideal{r}{}(p)$,
 i.e., the Noetherian relation $\red{}{\myr}{r}{}$ induced by
 them need not be confluent, even restricted to the set $\{p \mrm w
 \mid w \in \m \}$ and so the elements in 
 $\ideal{r}{}(p)$ do not necessarily right reduce to zero.
\begin{example}\label{exa.satr}~\\
{\rm
Let $\Sigma = \{ a, b, c \}$ and
 $T = \{ a^{2} \myr \lambda, b^{2} \myr \lambda, ab \myr c, ac \myr b, cb \myr a
 \}$ be a presentation of a monoid $\m$ (which is in fact a group)
  with a length-lexicographical ordering induced by $a \succ b \succ c$.
Further let us consider the polynomial $p = a + b + c \in \q[\m]$.
\\
Then  $S = \{ a + b + c, a + c + \lambda, bc + c^2 + b \} \in\SAT(p)$, but
  $\red{}{\myr}{r}{S}$ is not  confluent on  the set
 $\{ p \mrm w \mid w \in \m \}$.
This follows as $a+b+c  \red{}{\myr}{r}{ a + c + \lambda} b- \lambda$
  and $a+b+c  \red{}{\myr}{r}{ a + b + c} 0$,
  but $ b- \lambda \nred{*}{\myr}{r}{S} 0$.
\exaend
}
\end{example}
Note that this corresponds to the fact that the set $\{ a+b+c \}$ is
no Gr\"obner basis with respect to $\red{}{\myr}{s}{}$ (compare example \ref{exa.no.strong.gb}).
The following example shows that right  Gr\"obner bases are strong Gr\"obner bases
  but not vice versa.
\begin{example}~\\
{\rm
Let $\m$ be presented as in example \ref{exa.satr} above and
 $F = \{ a+c+ \lambda, b- \lambda \} \subseteq \q[\m]$. 
\begin{enumerate}
\item $F$ is a  strong Gr\"obner basis. \\
      We have to show that all strong s-polynomials reduce to zero.
      A quick inspection reveals that
       $U_{a+c+ \lambda, b- \lambda} = \{ (bw,w) \mid w \in \m
       \backslash b\Sigma^*, w \neq \lambda \}$.
      Hence we have
      $\spol{s}(a+c+ \lambda, b- \lambda,bw,w) =
         (a+c+ \lambda) \mrm bw - (b - \lambda) \mrm w
         = cw + a \mm w + \underline{bw} - \underline{bw} + w
         = a \mm w + cw + w
         = (a + c + \lambda) \mrm w 
         \red{}{\myr}{s}{a+c+ \lambda} 0$.
\item $F$ is no right Gr\"obner basis.\\
      We have $(a+c+ \lambda) \mrm ba = ca + \lambda + \underline{ba}
                \in \ideal{r}{}(a+c+ \lambda, b-\lambda)$, and
       hence $ba + ca + \lambda$ is congruent to zero modulo this
       right ideal.
      But $ba + ca + \lambda$ does not right reduce to zero by $F$,
       as the following inspection shows.
      The polynomial $ba + ca + \lambda$ is not right reducible by $a + c + \lambda$:
      Trying to modify $a$ in order to reduce $ba$ or $ca$ we get,
       $(a+c+\lambda)\mrm ca = ba + \underline{c^2a} + ca$ and
       $(a+c+\lambda)\mrm ba = ca + \lambda + \underline{ba}$.
      Furthermore, $b- \lambda$ can only be applied to reduce
       terms  beginning with $b$.
      Hence the only possible right reduction steps to take place are
       $ba + ca + \lambda \red{}{\myr}{r}{b - \lambda} ca + a + \lambda
       \red{}{\myr}{r}{a+c+\lambda} ca - c$
       and this polynomial is $F$-irreducible.
\item The set $G = \{ a+b+c, a+c+ \lambda, bc+c^2+b, b- \lambda \}$ is
      a right Gr\"obner basis (compare example \ref{exa.interreduction} for a
      proof).
\exaend
\end{enumerate}
}
\end{example}
The next example states that there are cases where finite strong Gr\"obner
 bases exist, but no finite
 right Gr\"obner bases.
\begin{example}\label{exa.strong.but.no.right}~\\
{\rm
Let $\Sigma = \{ a, b, c, d, e, f \}$ and  
 $T= \{ abc \myr ba, bad \myr fe , fbc \myr bf, bfd \myr ae \}$ be a presentation of
 a monoid $\m$ with a
 length-lexicographical ordering induced by  an arbitrary precedence
 $\succ$ on 
 $\Sigma$.
\\
We have seen in example \ref{exa.no.finite.sat.set} that the
polynomial $p = a+f \in \q[\m]$ has no finite saturating set and hence
there exist no finite right Gr\"obner bases for the right ideal 
generated by $p$. 
On the other hand the set $\{ p \}$ itself is a strong Gr\"obner basis
(compare example \ref{exa.no.finite.sat.set} to see that no
critical situations exist).
\exaend
}
\end{example}
We move on now to study how right Gr\"obner bases are connected to
stable standard bases and how they can be characterized by right
reduction. 
Let us start by proving an analogon to lemma \ref{lem.srprop}
\begin{lemma}\label{lem.rprop}~\\
{\sl
Let $F$ be a set of polynomials  and $p$  a non-zero polynomial in $\myk[\m]$.
\begin{enumerate}
\item Then $p \red{*}{\myr}{r}{F} 0$ implies the existence of a stable 
  standard representation for $p$.
\item In case $p$ has a stable standard representation with respect to $F$,
  then $p$ is right reducible at its head monomial by $F$,
  i.e., $p$ is top-reducible by $F$.
\item\label{lem.rprop.3} In case $F$ is a stable standard basis,
  every polynomial $p \in \ideal{r}{}(F)\backslash \{ 0 \}$
 is top-reducible to zero by $F$ using right reduction.
\lemend
\end{enumerate}
}
\end{lemma}
\Ba{}
\begin{enumerate}
\item This follows directly by adding up the polynomials used in the
  right reduction steps  occurring in $p \red{*}{\myr}{r}{F} 0$.
\item This is an immediate consequence of definition \ref{def.rsr}
  as the existence of a polynomial $f$ in $F$ and an element
  $w\in\m$ with $\hterm(f \mrm w)= \hterm(f) \mm w= \hterm(p)$ is guaranteed.
\item We show that every non-zero polynomial $p \in
  \ideal{r}{}(F)\backslash \{ 0 \}$ is top-reducible to zero using $F$
  by induction on $\hterm(p)$. 
  Thus let $\hterm(p) = \min \{ \hterm(g) \mid g \in \ideal{r}{}(F)\backslash \{ 0 \}  \}$.
  Then as $p \in \ideal{r}{}(F)$ and $F$ is a stable standard basis, we
   have $p = \sum_{i=1}^{k} \alpha_i \skm f_{i} \mrm w_i$, with 
   $\alpha_i \in \myk^*, f_{i} \in F, w_i \in \m$ and 
   $\hterm(p) \succeq \hterm(f_{i}) \mm w_i \succeq \hterm(f_{i} \mrm w_i)$
   for all  $1 \leq i \leq k$.
  Without loss of generality, let $\hterm(p) = \hterm(f_1) \mm w_1=
  \hterm(f_1 \mrm w_1)$.
  Hence, the polynomial $p$ is right reducible by $f_1$ at its head monomial.
  Let $p \red{}{\myr}{r}{f_1} q$, i.e.,
   $q = p - \hc(p) \skm \hc(f_1 \mrm w_1)^{-1} \skm f_1 \mrm w_1$, and
   by the definition of right reduction the term $\hterm(p)$ is
   eliminated from $p$ implying that $\hterm(q) \pred \hterm(p)$ as $q < p$.
  Thus, as $q \in \ideal{r}{}(F)$ and  $\hterm(p)$ was minimal among
   the head terms of the elements
   in the right ideal generated by $F$, this implies $q=0$, and, 
   therefore, $p$ is top-reducible to zero by $f_1$ in 
   one step.
   On the other hand, in case 
    $\hterm(p) \succ \min \{ \hterm(g) \mid g \in \ideal{r}{}(F)\backslash \{ 0
    \}  \}$, by the
    same arguments used before we can top right reduce $p$ to a polynomial
    $q$ with $\hterm(q) \pred \hterm(p)$, and, thus, by our induction
    hypothesis we know that $q$ and hence $p$ is top-reducible to zero.
\\
\qed
\end{enumerate}\renewcommand{\baselinestretch}{1}\small\normalsize
As before we find that right Gr\"obner bases and stable standard bases
are in fact equivalent.
\begin{theorem}~\\
{\sl
For a set of polynomials $F$ in $\myk[\m]$,
 the following statements are equivalent:
\begin{enumerate}
\item $F$ is a right Gr\"obner basis.
\item For all polynomials $g \in \ideal{r}{}(F)$ we have $g \red{*}{\myr}{r}{F} 0$.
\item $F$ is a stable standard basis.
\theoend
\end{enumerate}
}
\end{theorem}
\Ba{}~\\
\mbox{$1 \R 2:$ }
  By (i) of definition \ref{def.rgb} we know that $g \in
  \ideal{r}{}(F)$ implies $g \red{*}{\lr}{r}{F} 0$ and since
  $\red{}{\myr}{r}{F}$ is confluent and $0$ is irreducible $g \red{*}{\myr}{r}{F} 0$ follows
  immediately.

\mbox{$2 \R 3:$ }
  This follows directly by adding up the polynomials used in the
  right reduction steps  of $g \red{*}{\myr}{r}{f} 0$.

\mbox{$3 \R 1:$ }
  In order to show that $F$ is a Gr\"obner basis we have
  to prove two subgoals:
$\red{*}{\lr}{r}{F}  = \;\; \equiv_{\ideal{r}{}(F)}$ has already been shown in lemma \ref{lem.ssb.cong}.
It remains to show that $\red{}{\myr}{r}{F}$ is confluent.
             Since $\red{}{\myr}{r}{F}$ is Noetherian, we only have to prove
             local confluence.
             Suppose $g \red{}{\myr}{r}{F} g_1, g \red{}{\myr}{r}{F} g_2$
              and $g_1 \neq g_2$.
             Then $g_1 - g_2 \in \ideal{r}{}(F)$ and, therefore,
              is top-reducible to zero by $F$ as a result of 
              lemma \ref{lem.rprop}.
             Thus 
              lemma \ref{lem.confluentr} provides the existence of
              a polynomial $h \in \myk[\m]$ such that
              $g_1 \red{*}{\myr}{r}{F} h$ and $g_2 \red{*}{\myr}{r}{F} h$,
              i.e.,  $\red{}{\myr}{r}{F}$ is confluent.
\\
\qed
%
We continue by examining critical pairs of polynomials with respect to
right reduction and defining corresponding s-polynomials in order to
characterize right Gr\"obner bases.
\begin{definition}\label{def.cpr}~\\ 
{\rm
Given two non-zero polynomials $p_{1}, p_{2} \in \myk[\m]$\footnote{Notice that
  $p_1=p_2$ is possible.},
 every pair of elements $w_{1}, w_{2}$ in $\m$ such that
 $\hterm(p_{1} \mrm w_{1}) = \hterm(p_{1}) \mm w_{1} = \hterm(p_{2}) \mm  w_{2}
  = \hterm(p_{2} \mrm  w_{2})$, defines a 
 \index{right!s-polynomial}\index{s-polynomial!right}\betonen{(right) s-polynomial}\/
 $$ \spol{}(p_{1}, p_{2}, w_{1}, w_{2}) = \hc(p_1 \mrm w_1)^{-1} \skm p_1 \mrm w_1
                                    -  \hc(p_2 \mrm w_2)^{-1} \skm p_2 \mrm w_2.$$
Let $U_{p_1,p_2}  \subseteq \m \times \m$ be the set
containing {\em all} such pairs $w_1,w_2 \in\m$. 
\dend
}
\end{definition}
A right s-polynomial will be called non-trivial in case it is non-zero and
notice that for non-trivial s-polynomials we always have 
$\hterm(\spol{}(p_{1}, p_{2}, w_{1}, w_{2})) \prec \hterm(p_1) \mm w_1
= \hterm(p_2) \mm w_2$.
The set $U_{p_1,p_2}$ can be empty, finite or even infinite
depending on $\m$ as the following example reveals.
\begin{example}~\\
{\rm
Let $\Sigma = \{ a,b,c \}$ and $T = \{ a^2 \myr \lambda, b^2 \myr \lambda,
 ab \myr c, ac \myr b, cb \myr a,  \}$ be a presentation of a monoid $\m$
 (which is in fact a group) with a length-lexicographical ordering
 induced by
 $a \succ b \succ c$.
\\
Then for the polynomials $c+ \lambda,b + \lambda, a+c$ and
 $ c^2 +b \in \q[\m]$ we find that
 the set $U_{c+ \lambda, b + \lambda}$ is empty,
 the set
 $U_{a+c, c + \lambda} = \{ (\lambda,b)\}$ is finite
 and the set $U_{c+\lambda , c^2 + b}$ contains all pairs
 $(c^{n+1},c^n)$, $n \in \n$, i.e., is infinite.
\exaend
}
\end{example}

Unlike in Buchberger's approach and the previous section,
 s-polynomials as defined above are no longer
 strong enough to characterize Gr\"obner bases as they can only be
 used to give a confluence test but not to ensure that the right ideal
 congruence is expressible by reduction.
\begin{example}~\\
{\rm 
Let $\Sigma = \{ a,b \}$ and $T = \{ ab \myr \lambda, ba \myr \lambda\}$
 be a presentation of a group $\g$  with a length-lexicographical ordering
 induced by $a \succ b$.
\\
Examining  the set $F= \{ b + \lambda\}$ we find that $F$ has no non-trivial
 s-polynomials.
On the other hand we get
 $(b + \lambda) - (b + \lambda) \mrm a = a - b$
 and hence $a - b \in \ideal{r}{}(F)$.
 But the polynomial $a-b$ does not right reduce to zero by $F$.
\exaend
}
\end{example}
This phenomenon is due to the fact that there are critical situations
between the polynomial $b+\lambda$ and the rules in $T$ viewed as
polynomials in the free monoid ring, i.e., $ab - \lambda, ba -
\lambda$.
These situations are considered in the saturating process in procedure
{\sc Saturation 1} on page \pageref{saturation2} (compare the definition
of the sets $C(t)$ there).
Nevertheless, if we require our set of polynomials to be saturated, we can
characterize Gr\"obner bases in a familiar way.
\begin{theorem}\label{theo.rcp}~\\
{\sl
For  a saturated set  $F$ of polynomials in $\myk[\m]$,
 the following statements are equivalent:
\begin{enumerate}
\item For all polynomials $g \in \ideal{r}{}(F)$ we have $g \red{*}{\myr}{r}{F} 0$.
\item For all not necessarily different polynomials 
      $f_{k}, f_{l} \in F$ and every corresponding pair 
      $(w_{k}, w_{l}) \in U_{f_k,f_l}$  we have 
      $ \spol{}(f_{k}, f_{l}, w_{k}, w_{l}) \red{*}{\myr}{r}{F} 0$.
\theoend
\end{enumerate}
}
\end{theorem}
\Ba{}~\\
\mbox{$1 \R 2:$ }
 Let $ (w_{k}, w_{l}) \in U_{f_k,f_l}$ give us an
 s-polynomial belonging to the polynomials $f_k,f_l$.
 Then by definition \ref{def.cpr} we get
  $$\spol{}(f_{k}, f_{l}, w_{k}, w_{l}) =  \hc(p_1 \mrm w_1)^{-1} \skm p_1 \mrm w_1
      -  \hc(p_2 \mrm w_2)^{-1} \skm p_2 \mrm w_2 \:\in \ideal{r}{}(F)$$
 and hence $\spol{}(f_{k}, f_{l}, w_{k}, w_{l}) \red{*}{\myr}{r}{F} 0$.

\mbox{$2 \R 1:$ }
     We have to show that every non-zero polynomial 
      $g \in \ideal{r}{}(F)\backslash \{ 0 \}$ is $\red{}{\myr}{r}{F}$-reducible to zero.
     Remember that for
      $h \in \ideal{r}{}(F)$, $ h \red{}{\myr}{r}{F} h'$ implies $h' \in \ideal{r}{}(F)$.
     Thus as  $\red{}{\myr}{r}{F}$ is Noetherian
      it suffices to show that every 
      $g \in \ideal{r}{}(F) \backslash \{ 0 \}$ is $\red{}{\myr}{r}{F}$-reducible.
     Now, let $g = \sum_{j=1}^m \alpha_{j} \skm f_{j} \mrm w_{j}$  be a
      representation of a non-zero polynomial $g$ such that  $\alpha_{j} \in \myk, f_j \in F,
      w_{j} \in \m$.
     By lemma \ref{lem.prop1} we can assume  
      $\hterm(f_{i} \mrm w_{i}) = \hterm(f_{i}) \mm w_{i}$ since
      $F$ is saturated.
     Depending on this  representation of $g$ and a well-founded total ordering $\succeq$ on $\m$ we define
      $t = \max \{ \hterm(f_{j}) \mm w_{j} \mid j \in \{ 1, \ldots m \}  \}$ and 
      $K$ is the number of polynomials $f_j \mrm w_j$ containing $t$ as a term.
Then $t \succeq \hterm(g)$ and
in case $\hterm(g) = t$ this immediately implies that $g$ is
$\red{}{\myr}{r}{F}$-reducible.
So by lemma \ref{lem.rprop} it is sufficient to  show that
$g$ has a stable standard representation, as this implies that $g$ is
top-reducible using $F$.
This will be done by induction on $(t,K)$, where
      $(t',K')<(t,K)$ if and only if $t' \prec t$ or $(t'=t$ and
      $K'<K)$\footnote{Note that this ordering is well-founded since $\succ$ is and $K \in\n$.}.
In case $t \succ \hterm(g)$
      there are two polynomials $f_k,f_l$ in the corresponding 
      representation\footnote{Not necessarily $f_l \neq f_k$.}
      such that  $t = \hterm(f_k \mrm w_k) = \hterm(f_k) \mm w_k =
      \hterm(f_l) \mm w_l= \hterm(f_l \mrm w_l)$.
     Then by definition \ref{def.cps} we have a corresponding s-polynomial
      $\spol{}(f_k,f_l,w_k,w_l) = \hc(f_k \mrm w_k)^{-1} \skm  f_k
      \mrm w_k - \hc(f_l \mrm w_l)^{-1} \skm f_l \mrm w_l$.
We will now change our representation of $g$ by using the additional
information on this s-polynomial in such a way that for the new
representation of $g$ we either have a smaller maximal term or the occurrences of the term $t$
are decreased by at least 1.
     Let us assume 
      $\spol{}(f_k,f_l,w_k,w_l) \neq 0$\footnote{In case 
        $\spol{}(f_k,f_l,w_k,w_l) = 0$,
               just substitute $0$ for $\sum_{i=1}^n \delta_i \skm h_i \mrm
               v_i$ in the equations below.}.
     Hence, the reduction sequence $\spol{}(f_k,f_l,w_k,w_l)
     \red{*}{\myr}{s}{F} 0 $ results in a  stable standard representation
     $\spol{}(f_k,f_l,w_k,w_l) =\sum_{i=1}^n \delta_i \skm h_i \mrm v_i$,
     where $\delta_i \in \myk^*,h_i \in F,v_i \in \m$ and
      all terms occurring in the sum are bounded by
      $\hterm(\spol{}(f_k,f_l,w_k,w_l)) \prec t$.
     This gives us: 
     \begin{eqnarray}
       &   & \alpha_{k} \skm f_{k} \mrm w_{k} + \alpha_{l} \skm f_{l} \mrm w_{l}  \nonumber\\ 
       &   &  \nonumber\\
       & = &  \alpha_{k} \skm f_{k} \mrm w_{k} + 
              \underbrace{ \alpha'_{l} \skm \beta_k \skm f_{k} \mrm w_{k}
                   - \alpha'_{l} \skm \beta_k \skm f_{k} \mrm w_{k}}_{=\, 0} 
                   + \alpha'_{l}\skm \beta_l  \skm f_{l} \mrm w_{l} \nonumber\\
       &   &  \nonumber\\
       & = & (\alpha_{k} + \alpha'_{l} \skm \beta_k) \skm f_{k} \mrm w_{k} - \alpha'_{l} \skm
             \underbrace{(\beta_k \skm f_{k} \mrm w_{k}
             -  \beta_l \skm f_{l} \mrm w_{l})}_{=\,
             \spol{}(f_k,f_l,w_k,w_l)}                               \nonumber\\
       & = & (\alpha_{k} + \alpha'_{l} \skm \beta_k) \skm f_{k} \mrm w_{k} - \alpha'_{l} \skm
                   (\sum_{i=1}^n \delta_{i} \skm h_{i} \mrm v_{i}) \label{s2}
     \end{eqnarray}
     where  $\beta_k = \hc(f_k \mrm w_k)^{-1}$, $\beta_l = \hc(f_l \mrm
     w_l)^{-1}$ and  $\alpha'_l \skm \beta_l = \alpha_l$.
     By substituting (\ref{s2}) in our representation of $g$ 
 either $t$ disappears or $K$ is decreased.  
\\
\qed
Note that reducing a Gr\"obner basis need not preserve the properties
 of Gr\"obner bases, since for polynomials $p,q,q_1,q_2 \in \myk[\m]$,
 $p \red{}{\myr}{r}{q}$ and $q \red{}{\myr}{r}{q_1} q_2$  need
 not imply  $p \red{}{\myr}{r}{\{ q_1, q_2 \}}$. 
\begin{example}\label{exa.interreduction}~\\
{\rm
Let $\Sigma = \{ a,b,c \}$ and $T = \{ a^2
 \myr \lambda, b^2 \myr \lambda, ab \myr c, ac \myr b, cb \myr a  \}$ be a
 presentation of a monoid $\m$ (which is in fact a group) with a
 length-lexicographical ordering induced by
 $a \succ b \succ c$.
\\
Then the set $F = \{ a+b+c, a+c+ \lambda, bc+c^2+b, b- \lambda \}$ is
a Gr\"obner basis, as $F$ is saturated and all possible s-polynomials
right reduce to zero.
To see this we have to examine the possible s-polynomials.
\begin{enumerate}
\item $U_{a+b+c,a+c+\lambda} = \{ (\lambda,\lambda) \}$ 
      as for all other solutions of the
      equation $a \mm {\rm x} = a \mm {\rm y}$ at least one of the polynomials
      is not stable and
      $\spol{}(a+b+c,a+c+\lambda,\lambda, \lambda) = a+b+c -
      a-c-\lambda = b  - \lambda
      \red{}{\myr}{r}{b-\lambda} 0$.
\item $U_{a+b+c, bc+c^2+b} = \emptyset$ and
      $U_{a+c + \lambda , bc+c^2+b} = \emptyset$
      as for all solutions of the
      equation $a \mm {\rm x} = bc \mm {\rm y}$ one of the polynomials
      is not stable\footnote{The solutions to this equation can be 
      written as $\{ (a \mm w,bc \mm w), (c \mm w, ba
      \mm w), (w, bca \mm w) \mid w \in \m \}$.}. 
\item $U_{a+b+c, b - \lambda} = \emptyset$ and
      $U_{a+c + \lambda , b - \lambda} = \emptyset$
      as for all solutions of the
      equation $a \mm {\rm x} = b \mm {\rm y}$ one of the polynomials
      is not stable\footnote{The solutions to this equation can be 
      written as $\{  (a \mm w. b \mm w) \mid  w \in \m \}$.}.
\item $U_{bc+c^2+b, b - \lambda} = \{ (b,a), (w,c \mm w) \mid  w \in \m, c \mm w
  = cw \}$ and
      $\spol{}(bc+c^2+b, b - \lambda, b, a) = ba + ca + \lambda - ba +
      a = ca + a + \lambda \red{}{\myr}{r}{a+b+c} 0$,
      $\spol{}(bc+c^2+b, b - \lambda,w, cw) = bcw + c^2w + bw - bcw +
      w = c^2w + bw + c \mm w \red{}{\myr}{r}{a+b+c} 0$. 
\end{enumerate}
Now let us look at two possible ways of interreducing $F$:
\\
If we first remove
$bc+c^2+b$ as
$$bc+c^2+b \red{}{\myr}{r}{b - \lambda} c^2+b+c \red{}{\myr}{r}{a+b+c}
0$$
and then $a+b+c$ since
$$a+b+c \red{}{\myr}{r}{a+c+\lambda} b - \lambda \red{}{\myr}{r}{b -
  \lambda} 0$$
 this gives us a set $F' = \{ a+c+\lambda, b- \lambda \}$ which is
no longer a right Gr\"obner basis.
This is due to the fact that we have $(a+c+\lambda) \mrm c - (b -
\lambda) = b+c^2+c - b - \lambda = c^2+c+\lambda 
 \in \ideal{r}{}(F) = \ideal{r}{}(F')$ and this
 polynomial is not right reducible by the polynomials in the set $F'$.
\\
On the other hand , if we first remove the polynomial $a+b+c$, then 
$bc+c^2+b$ no longer right reduces to zero with $\{ a+c+\lambda, b-
\lambda \}$.
Instead we get $$bc+c^2+b \red{}{\myr}{r}{b-\lambda} c^2 + b+c
\red{}{\myr}{r}{b-\lambda} c^2 + c + \lambda$$
and the set $F'' = \{ a+c+\lambda, b-\lambda, c^2+c+\lambda \}$,
although not saturated, is a right Gr\"obner basis.

Notice that in strongly interreducing $F$ we also get a set
 which is no longer a strong Gr\"obner basis:
$$a+b+c \red{}{\myr}{s}{a+c+\lambda} b-\lambda \red{}{\myr}{s}{b-\lambda} 0,$$
$$a+c+\lambda \red{}{\myr}{s}{ba + c^2 + b} b-\lambda \red{}{\myr}{s}{b-\lambda} 0,$$
$$bc+c^2+b \red{}{\myr}{s}{b-\lambda}c^2+b+c \red{}{\myr}{s}{b-\lambda} c^2+c+\lambda$$
leaves us with $F' = \{c^2+c+\lambda,b-\lambda\}$ which is no strong Gr\"obner
 basis, as we are no longer able to reduce $a+b+c$.
\exaend
}
\end{example}
%
Unfortunately, theorem \ref{theo.rcp} is only of theoretical interest
 since in general the following uniform problem is undecidable, 
 even in monoids where the solvability of
 equations of the form $u \mm {\rm x} = v \mm {\rm y}$ is decidable.

\begin{tabbing}
XXXXXXXX\=XXXXX \kill
{\bf Given:} \>  Two polynomials $p,q \in \myk[\m]$, and \\
             \> $(\Sigma, T)$ a convergent semi-Thue system presenting $\m$. \\
{\bf Question:} \> Does there exist an s-polynomial for $p$ and $q$?
\end{tabbing}

To see this we need a construction introduced in lemma
 \ref{lem.andrea}.
\begin{example}~\\
{\rm
Let  $\Phi$ be an universal function and 
 ${\rm\bf f}_1,{\rm\bf f}_2, {\rm\bf f}_3 :  \n^2 \myr \n$  three
 primitive recursive functions, with
\\
$ {\rm\bf f}_1({\rm x},{\rm y}) = \left\{ \begin{array}{r@{\quad\quad}l}
                      0 & \exists {\rm z} \leq {\rm y} \mbox{ such
                        that the computation of 
                        $\Phi_{\rm x}({\rm x})$ halts within  z steps} \\
                      2 & \mbox{otherwise,}
                      \end{array} \right. $

$ {\rm\bf f}_2({\rm x},{\rm y}) = 1$ and $ {\rm\bf f}_3({\rm x},{\rm y}) = 3$.
\\
As described in the proof of lemma
 \ref{lem.andrea} we can construct a finite convergent interreduced semi-Thue
 system $T'$ over an alphabet $\Sigma'$.
Let $\m$ be the monoid presented by $(\Sigma', T')$.
Then for $n \in \n$ and the corresponding polynomials $p_n = f_1L_1L_2b^na - f_2L_1L_2b^na$ and
 $q_n = f_1L_1L_2b^na - f_3L_1L_2b^na$ it is undecidable whether an
 s-polynomial exists in $\q[\m]$.
This is due to the fact that although the equation
 $f_1L_1L_2b^na \mm {\rm x} = f_1L_1L_2b^na \mm {\rm y}$ has all monoid
 elements $(w,w)$ as trivial solutions, 
 $p_n$ and $q_n$ have an s-polynomial if and only if $f_1L_1L_2b^na \mm w
 \syll f_2L_1L_2b^na \mm w$ for such a trivial solution.
But as we have seen in the proof of lemma \ref{lem.andrea}, it is
 not uniformly decidable given an arbitrary $n \in \n$ whether $f_2L_1L_2b^na \mm w
 \syll f_1L_1L_2b^na \mm w$.  
\exaend
}
\end{example}
This example reveals how closely related the problems of saturation
 and s-polynomials are. 
All these problems stem from the fact that the ordering on the monoid
 need not be compatible and hence reduction need not be preserved under multiplication.
\auskommentieren{
Trying to localize this test severe problems arise due to the
 lacks of our reduction relation,
 especially as $p \red{*}{\myr}{r}{F} 0 $ does not imply $p
 \mrm w \red{*}{\myr}{r}{F} 0$.
\begin{example}\label{exa.px}~\\
{\rm
Let $\Sigma = \{ a, b, c, d, e, f,g \}$ and
 $T = \{ ac \myr d, bc \myr e, dg \myr b, eg \myr f \}$ be a presentation
 of a monoid $\m$  with  a  length-lexicographical ordering induced by $ a \succ b \succ c \succ d \succ e \succ f
 \succ g$.
Further consider the set
 $F = \{ a +b, d+e, b+f, fc+e, d+ \lambda, b+g, gc+e, e+g, g^2+f,
 g+ \lambda \}$.
\\
Then we get $e - \lambda \red{}{\myr}{r}{e+g} -g-\lambda \red{}{\myr}{r}{g + \lambda} 0$,
 but $(e- \lambda) \mrm g = f - g$ is not right reducible to zero using $F$.
\exaend
}
\end{example}
Note that the set $F$ in this example is saturated and hence
 saturation does not enable localization.

In ordinary polynomial rings as $\myk[X_{1}, \ldots X_{n}]$ one can select a
 ``smallest'' critical pair by taking the least common multiply of 
 the head terms and it is sufficient to examine this case \cite{Bu85}.
In $\myk[\m]$ the situation is more complicated.
Reviewing  definition \ref{def.cpr} we see that it is important to solve the
 equation $\hterm(p_{1}) \mm {\rm x} = \hterm(p_{2}) \mm {\rm y}$.
Therefore, we are looking for a suitable ``basis'' of a set
 $$ U_{\hterm(p_{1}),\hterm(p_{2})} = \{ (w_1,w_2) \in \m^2 \mid \hterm(p_{1}) \mm w_1 = \hterm(p_{2}) \mm w_2 \}.$$
One idea might be  to look at a basis $B \subseteq U_{t_1,t_2}$
 such that for all $(w_1,w_2) \in  U_{\hterm(p_{1}),\hterm(p_{2})}$ we have $(b_1,b_2) \in B, w \in \m$ fulfilling
 $w_1 = b_1 \mm w, w_2 = b_2 \mm w$.\\
Another approach might be to look for a suitable ``basis'' of a set
 $$ U_{p_1,p_2} = \{ (w_1,w_2)\in\m^2 \mid
    \hterm(p_1 \mrm w_1)=\hterm(p_1) \mm w_1 = \hterm(p_2) \mm w_2=\hterm(p_2 \mrm w_2) \}.$$
$U_{p_1,p_2}$ describes real critical situations in the sense that
 $\hterm(p_1) \mm w_1 = \hterm(p_2) \mm w_2$ is an overlap,
 where both $p_1$ and $p_2$ can be applied for reduction.
Reviewing example \ref{exa.px} we see that  none of these approaches to define ``bases''  is  sufficient.
\begin{example} \label{exa.px2}~\\
{\rm
Let $\Sigma = \{ a, b, c, d, e, f,g \}$ and
 $T = \{ ac \myr d, bc \myr e, dg \myr b, eg \myr f \}$ be a presentation of a
 monoid $\m$ with  a  length-lexicographical ordering induced by $ a \succ b \succ c \succ d \succ e \succ f  \succ g$.
Further take $F = \{ a +b, d+e, b+f, fc+e, d+ \lambda, b+g, gc+e, e+g, g^2+f, g+ \lambda \}$.
\\
Looking at the polynomials  $a+b$ and $d+ \lambda$ we get a ``real'' critical situation in $d$, which leads to
 the s-polynomial
 $e - \lambda \red{}{\myr}{r}{e+g} -g-\lambda \red{}{\myr}{r}{g + \lambda} 0$,
 but $(e- \lambda) \mrm g = f - g$ is not reducible to zero using $F$.
\exaend
}
\end{example}
If this equation is solvable, then a pair of solutions $w_1, w_2$
 such that $\hterm(p_1 \mrm w_1) = \hterm(p_1) \mm w_1 = \hterm(p_2) \mm w_2 =
 \hterm(p_2 \mrm w_2)$ defines an s-polynomial.

One idea to localize the confluence test given in theorem
\ref{theo.rcp}
might be to look for a basis $B \subseteq U_{p_1,p_2}$ such that
for all pairs of solutions $w_1,w_2) \in U_{p_1,p_2}$ there exists a
pair $(b_1,b_2) \in B$ and an element $w \in \m$ such that $w_1 = b_1
\mm w$ and $w_2 = b_2 \mm w$.
Unfortunately, this approach does not work in general, as the
following example shows.
The crucial point is that $p \red{*}{\myr}{r}{F} 0 $ does not imply $p
 \mrm w \red{*}{\myr}{r}{F} 0$.
Localization in Buchberger's approach is possible since we have an
 admissible term ordering which enables us to prove that $p \red{*}{\myr}{}{F} 0$
 implies $p \mrm w \red{*}{\myr}{}{F} 0$ for arbitrary terms $w$.
}

We will end this section by a remark on the algebraic characterization of
 Gr\"obner bases in terms of ideals in the set of terms,
 e.g., the free commutative
 monoid generated by the indeterminants in the usual polynomial ring.
In the polynomial ring $\myk[X_1, \ldots,X_n]$ we know that for a set
of polynomials $F$ we have that the set
$\hterm(\ideal{}{\myk[X_1, \ldots,X_n]}(F) \backslash \{ 0 \})$ itself is an ideal in the set of
 terms ${\cal T}$,
 in fact if $F$ is a Gr\"obner basis then for the ideal generated by
 $\hterm(F)$ in ${\cal T}$, we have  $\hterm(\ideal{}{\myk[X_1,
   \ldots,X_n]}(F) \backslash \{ 0 \})=\ideal{}{{\cal T}}(\hterm(F))$.
This is crucial when proving termination of Buchberger's algorithm.
Unfortunately, this no longer holds for arbitrary monoid rings.
\begin{example}\label{exa.aa.bb.ba=ab}~\\
{\rm
Let $\Sigma = \{ a,b \}$ and $T = \{ a^2 \myr \lambda, b^2 \myr \lambda ,
 ba \myr ab \}$ be a presentation of  a monoid $\m$
 (which is in fact a commutative group) with $a \pred b$  
 inducing a  length-lexicographical ordering on $\m$.
Note that $\m$ is a finite group consisting of the elements $\{
\lambda, a,b, ab \}$.
\\
For $p=ab+ \lambda \in \q[\m]$, we get
 $\ideal{r}{\q[\m]}(p) = \{ \alpha \skm (ab + \lambda) + \beta \skm (
 b + a) \mid  \alpha, \beta \in \q \}$.
Then the set $\{ p \}$ itself is a right Gr\"obner basis,
  but we have $\hterm(\ideal{r}{\q[\m]}(p)\backslash \{ 0 \})=
 \{ b, ab \} \neq \ideal{r}{\m}(ab) = \{ \lambda, a, b, ab \}$.
\exaend
}
\end{example}
In the next section we will introduce weakenings of right
 reduction which provide enough
 information to localize critical situations and characterize
 Gr\"obner bases in some way
 by the right ideals generated by their head terms.
But, nothing really comes for free and we will have to do
 saturation with respect to these
 weaker reductions in order to establish the right ideal congruence.

%% file: prefixreduction.tex
\section{The Concept of Prefix Reduction}\label{section.prefixreduction}
%
In the previous section we have investigated stable 
 standard representations
 of polynomials and we have seen how they are connected to right reduction.
Hence, we start this section by  refining our view on representations of
 polynomials which will lead to a refinement of right reduction.
We will see later on that for certain classes of groups this enables us to compute
 finite Gr\"obner bases for finitely generated right ideals (compare
 chapter \ref{chapter.grouprings}).
\begin{definition}\label{def.psr}~\\
{\rm
Let $F$ be a set of polynomials and $p$ a non-zero
polynomial in $\myk[\m]$.
A representation 
$$ p = \sum_{i=1}^{n} \alpha_i \skm f_{i} \mrm w_i, \mbox{ with } \alpha_i \in \myk^*,
f_{i} \in F, w_i \in \m $$
is called a  
 \index{prefix!standard representation}\index{standard representation!prefix}\betonen{prefix standard
 representation}\/ with respect to the set of polynomials $F$,
 in case for all $1 \leq i \leq n$ we have $\hterm(p) \succeq \hterm(f_{i})w_i$.
A set $F \subseteq \myk[\m]$ is called a 
  \index{prefix!standard basis}\index{standard basis!prefix}\betonen{prefix standard
  basis}\/ if every non-zero polynomial
$g \in \ideal{r}{}(F)$ has a prefix  standard representation with
respect to $F$.
\dend
}
\end{definition}
Notice that $\hterm(p) \succeq \hterm(f_{i})w_i$ immediately implies
 $\hterm(p) \succeq \hterm(f_{i})w_i \succeq \hterm(f_{i} \mrm w_i)$
 and in case $\hterm(f_i) \mm w_i \id \hterm(f_i)w_i$ even 
 $\hterm(p) \succeq \hterm(f_{i})w_i \id \hterm(f_{i} \mrm w_i)$.
On the other hand,  in case $\hterm(p) =
 \hterm(f_{i} \mrm w_i)$ we must have $\hterm(f_{i} \mrm
 w_i) \id \hterm(f_{i})w_i$, i.e., $\hterm(f_i)$ is a prefix of $\hterm(p)$,
 and this situation occurs for at least one polynomial in the
 representation. 

Before we move on, we give two technical but nevertheless
 useful lemmata comparable to lemma \ref{lem.ssr.prop1} which state
 when even more restricted prefix standard
 representations for polynomials exist.
\begin{lemma}\label{lem.psr.prop2}~\\
{\sl
Let $F$ be a set of polynomials in $\myk[\m]$ such that for all $f \in
F$ and all $w \in \m$ the polynomial $f \mrm w$ has a prefix standard
representation with respect to $F$ in case it is non-zero.
Then there even exists a prefix standard representation 
 $f \mrm w = \sum_{i=1}^{n} \alpha_i \skm f_{i} \mrm w_i$, with $\alpha_i \in \myk^*,
 f_{i} \in F, w_i \in \m$ such that $\hterm(f \mrm w)
 \succeq \hterm(f_{i} \mrm w_i) \id \hterm(f_{i})w_i$.
\lemend
}
\end{lemma}
\Ba{}~\\
We will prove this lemma by contradiction.
Let us assume this is not true.
Then there exists a counter-example  $f \mrm w$ such that 
 $\hterm(f \mrm w)$ is minimal among all counter-examples.
By our assumption $f \mrm w$ has a prefix standard representation,
e.g.\  $f \mrm w = \sum_{i=1}^{m} \alpha_i \skm g_{i} \mrm w_i$,
 with $\alpha_i \in \myk^*, g_{i} \in F, w_i \in \m$.
Without loss of generality we can assume that for some  $k \leq m$, 
 $g_1, \ldots, g_k$, are the polynomials involved in the head
 term of $f \mrm  w$, i.e., $\hterm(f \mrm w) \id \hterm(g_i)w_i$ for
 all $1 \leq i \leq k$.
Hence, we know $k < m$, as otherwise we would get a contradiction to
$f \mrm w$ being a counter-example.
Furthermore, for all $k+1 \leq j \leq m$ we know $\hterm(g_j \mrm w_j)
\prec \hterm(f \mrm w)$ and hence every such polynomial has a prefix
standard representation of the desired form, say
 $g_j \mrm w_j = \sum_{l=1}^{n_l} \alpha'_{j_l} \skm g_{j_l} \mrm
 w'_{j_l}$,
 with $\alpha'_{j_l} \in \myk^*$, $g_{j_l} \in F$ and $w'_{j_l} \in \m$.
Thus the representation
 $f \mrm w = \sum_{i=1}^k \alpha_i \skm g_i \mrm w_i + 
 \sum_{i=k+1}^{n} \alpha_i \skm
 (\sum_{l=1}^{n_l} \alpha'_{j_l} \skm g_{j_l} \mrm w'_{j_l})$
is a prefix standard representation having the desired property, contradicting
our assumption.
\\
\qed
\begin{lemma}\label{lem.psr.prop1}~\\
{\sl
Let $F$ be a prefix standard basis in $\myk[\m]$.
Then every non-zero polynomial $p \in \ideal{r}{}(F)$ has a prefix standard
 representation
 $p = \sum_{i=1}^{n} \alpha_i \skm f_{i} \mrm w_i, \mbox{ with } \alpha_i \in \myk^*,
  f_{i} \in F$, and $w_i \in \m$
such that  for all $1 \leq i \leq n$ we even have $$\hterm(p) \succeq \hterm(f_{i})w_i
 \id \hterm(f_{i} \mrm w_i).$$
\lemend
}
\end{lemma}
\Ba{}~\\
Since $p \in \ideal{r}{}(F) \backslash \{ 0 \}$, the polynomial $p$
has a prefix standard representation with respect to $F$, say 
  $p = \sum_{i=1}^{n} \alpha_i \skm f_{i} \mrm w_i$ with $\alpha_i \in \myk^*$,
 $f_{i} \in F$, and $w_i \in \m$.
Moreover, by lemma \ref{lem.psr.prop2} every multiple $f_i \mrm w_i$
has a prefix standard representation with respect to $F$, say
 $f_i \mrm w_i = 
 \sum_{j=1}^{m} \beta_j \skm g_j \mrm v_j$ with $\beta_j \in \myk^*$,
 $g_{j} \in F$, and $v_j \in \m$, such that
 $$\hterm(p) \succeq\hterm(f_i \mrm w_i) \succeq \hterm(g_j \mrm v_j) \id \hterm(g_j)v_j.$$
\qed
Notice that a prefix standard basis is a stable standard
basis, as by lemma \ref{lem.psr.prop1} the equation $\hterm(p) \succeq  \hterm(f_i)w_i =
\hterm(f_i) \mm w_i = \hterm(f_i \mrm w_i)$ holds.
But the following example shows that the converse is not true.
\begin{example}~\\
{\rm
Let $\Sigma = \{ a,b \}$ and $T = \{a^2 \myr \lambda, b^2 \myr \lambda\}$
 be a presentation of a monoid $\m$ (which is in fact a group), with $a \succ b$ inducing a
 length-lexicographical ordering on $\m$.
Further take the set $F = \{ ab \} \subseteq \q[\m]$.
\\
Then  all non-zero polynomials $g \in \ideal{r}{}(F)$ obviously have a stable
 standard representation in $F$, but
 e.g.\  the polynomial $a \in \ideal{r}{}(F)$ has no prefix standard
 representation.
\phantom{XX}
\exaend
}
\end{example}
Prefix standard representations provide
 us with enough information
 to characterize prefix standard bases
 (which are closely related to special Gr\"obner
 bases as we will see later on) by their head terms in a way similar to
 the case of usual polynomial rings.
\begin{theorem}\label{theo.equiv1}~\\
{\sl
Let $F$ be a set of polynomials in $\myk[\m]$ and $G \subseteq
\ideal{r}{\myk[\m]}(F)\backslash \{ 0 \}$.
Then  the following statements are equivalent:
\begin{enumerate}
\item $G$ is a prefix standard basis for
  $\ideal{r}{\myk[\m]}(F)$\footnote{I.e., $\ideal{r}{}(F) =
    \ideal{r}{}(G)$ and $G$ is a prefix standard basis.}.
\item $\ideal{r}{\freemonoid}(\hterm(G)) \cap \m = \hterm(\ideal{r}{\myk[\m]}(F)
  \backslash \{ 0 \})$.
\end{enumerate}
Notice that the set $\hterm(\ideal{r}{\myk[\m]}(F)\backslash \{ 0 \})$
 in general is no right ideal in $\m$.
\theoend
}
\end{theorem}
\Ba{}~\\
\mbox{$1 \R 2:$ }
  The inclusion $\ideal{r}{\freemonoid}(\hterm(G)) \cap \m
   \subseteq \hterm(\ideal{r}{\myk[\m]}(F)\backslash \{ 0 \})$ follows
   at once as $\hterm(g)u \in \m$ for some $g \in G$, $u \in \Sigma^*$ implies
   $u \in \m$ and $\hterm(g \mrm u) \id \hterm(g)u$, and as the
   multiple $g \mrm u$ belongs to $\ideal{r}{\myk[\m]}(F)$.
  It remains to show that $\ideal{r}{\freemonoid}(\hterm(G)) \cap \m \supseteq
  \hterm(\ideal{r}{\myk[\m]}(F)\backslash \{ 0 \})$ holds.
  To see this, let $g \in \ideal{r}{\myk[\m]}(F)\backslash \{ 0 \}$.
  Then as $G$ is  a prefix standard basis  for $\ideal{r}{\myk[\m]}(F)$,
  there exists a prefix standard representation
   $g=\sum_{i=1}^{n} \alpha_i \skm g_{i} \mrm w_i$ with $\alpha_i \in \myk^*, g_i
   \in G$ and $w_i \in \m$ such that 
   $\hterm(g) \succeq \hterm(g_i)w_i$.
  Furthermore there exists $1 \leq k \leq n$ such that 
   $\hterm(g) \id \hterm(g_{k})w_k$, i.e., $\hterm(g) \in
   \ideal{r}{\freemonoid}(\hterm(G)) \cap \m$.

\mbox{$2 \R 1:$ }
 We have to show that every $g \in \ideal{r}{\myk[\m]}(F)\backslash \{ 0
 \}$ has a prefix standard representation with respect to $G$.
 This will be done by induction on the term $\hterm(g)$.
 In the base case we can assume
   $\hterm(g) = \min \{ w | w \in \hterm(\ideal{r}{\myk[\m]}(F)\backslash \{ 0 \})\}$.
  Then since  $\ideal{r}{\freemonoid}(\hterm(G)) \cap \m =
   \hterm(\ideal{r}{\myk[\m]}(F)\backslash \{ 0 \})$ there exists a polynomial 
   $f \in G$ such that $\hterm(g)\id \hterm(f)w$
   for some $w \in \m$.
  Eliminating the head term of $g$ by subtracting an appropriate right
   multiple of $f$ we get 
   $h = g - \hc(g) \skm \hc(f)^{-1} \skm f \mrm w$. 
  As $g \in \ideal{r}{\myk[\m]}(F)$,  $h$ lies in the right ideal
   generated by $F$.
  Moreover, since
   $\hterm(g)$ is minimal and $\hterm(h) \prec \hterm(g)$, we can conclude
   $h=0$ and $g$ has a prefix standard
   representation $g = \hc(g) \skm \hc(f)^{-1} \skm f \mrm w$.
  Now let us suppose $\hterm(g) \succ \min \{w | w \in
  \hterm(\ideal{r}{\myk[\m]}(F)\backslash \{ 0 \})\}$.
  Then again there exists a polynomial $f \in G$ such that $\hterm(g)\id\hterm(f)w$
   for some $w \in \m$.
  Hence looking at the polynomial 
   $h = g - \hc(g) \skm \hc(f)^{-1} \skm f \mrm w$ 
   we know that $h$ lies in the right ideal generated by $F$ and since
   $\hterm(h) \prec \hterm(g)$ either $h=0$, giving us that
   $g = \hc(g) \skm \hc(f)^{-1} \skm f \mrm w$,
   or  our induction hypothesis yields the existence
   of a  prefix standard representation for $h$ with respect to $G$,
   say $h = \sum_{j=1}^{m} \beta_j \skm g_{j} \mrm v_j$ where $\beta_j
   \in \myk^*$, $g_j \in G$ and $v_j \in \m$.
  Thus we have a  prefix  standard  representation of the polynomial 
   $g$, namely $g =
   \sum_{j=1}^{m} \beta_j \skm g_{j} \mrm v_j  +  \hc(g) \skm \hc(f)^{-1} \skm f \mrm w$. 
\\
\qed
We continue by giving a weakening of right reduction that will
correspond to the concepts of prefix standard representations and
prefix standard bases.
\begin{definition}\label{def.redp}~\\
{\rm
Let $p, f$ be two non-zero polynomials in $\myk[\m]$. 
We say $f$ \index{prefix!reduction}\index{reduction!prefix}\betonen{prefix reduces}
 $p$ to $q$ at a monomial
 $\alpha \skm t$ of $p$ in one step, denoted by $p \red{}{\myr}{p}{f} q$, if
\begin{enumerate}
\item[(a)] $\hterm(f)w \id t$ for some $w \in \m$,
            i.e., $\hterm(f)$ is a prefix of $t$, and
\item[(b)] $q = p - \alpha \skm \hc(f)^{-1} \skm f \mrm w$.
\end{enumerate}
We write $p \red{}{\myr}{p}{f}$ if there is a polynomial $q$ as defined
above and $p$ is then called  prefix reducible by $f$. 
Further we can define $\red{*}{\myr}{p}{}, \red{+}{\myr}{p}{}$,
 $\red{n}{\myr}{p}{}$  as usual.
Prefix reduction by a set $F \subseteq \myk[\m]$ is denoted by
 $p \red{}{\myr}{p}{F} q$ and abbreviates $p \red{}{\myr}{p}{f} q$
 for some $f \in F$,
 which is also written as  $p \red{}{\myr}{p}{f \in F} q$.
\dend
}
\end{definition}
Notice that  in the above definition the equation in 
 (a) has at most one solution and
 we then always have $\hc(f \mrm w) = \hc(f)$.
This is due to the fact that $t \id \hterm(f)w$ implies $\hterm(f)w =
\hterm(f \mrm w)$ and $\hterm(f)w \succ s \mm w$ for all $s \in
\terms(\reductum(f))$.
Further, in case $f$ prefix reduces $p$ to $q$ at the monomial $\alpha \skm
t$, we have $t \not\in \terms(q)$ and $p > q$.
Moreover, $\red{}{\myr}{p}{} \subseteq \red{}{\myr}{r}{}$ and the statements 1 to 3 of lemma \ref{lem.red} can
 be carried over to prefix reduction.
\begin{lemma}\label{lem.Noetherian}~\\
{\sl
Prefix reduction with respect to an arbitrary (possibly infinite) set
$U \subseteq \myk[\m]$ is Noetherian.
\lemend
}
\end{lemma}
\Ba{}~\\
This is an immediate consequence of the fact that $\geq$ is well-founded
 on $\myk[\m]$ and $f \red{}{\myr}{p}{g} f'$ implies $f > f'$.
Hence an infinite reduction sequence $f \red{}{\myr}{p}{g_1} f_1
 \red{}{\myr}{p}{g_2} f_2 \ldots \;$, $g_j \in U$, would imply the existence of an
 infinite strictly descending sequence of polynomials $f > f_1> \ldots \;$ in $\myk[\m]$
 contradicting the fact that $\geq$ is well-founded on $\myk[\m]$.
\\
\qed
Unlike in the case of strong right and right reduction now prefix
reducing a polynomial using itself must result in zero.
Therefore, we can define interreduced sets as follows.
\begin{definition}~\\
{\rm
We call a set of polynomials $F \subseteq \myk[\m]$  
 \index{reduced set}\index{interreduced!(set of polynomials)}
 \betonen{interreduced} or \betonen{reduced} with respect to
 $\red{}{\myr}{p}{}$, if no polynomial $f$ in
 $F$ is prefix reducible usinf the set $F \backslash \{
 f \}$.
\dend
}
\end{definition}
As for right reduction, $p \red{}{\myr}{p}{q_1} 0$ and $q_1
\red{}{\myr}{r}{q_2} 0$ imply $p \red{}{\myr}{p}{q_2} 0$.
Furthermore, prefix reduction  gives us additional
 information on the reduction
 step essential to the concept interreduction.
\begin{remark}\label{rem.transitiv}~\\
{\rm
Let $p \red{}{\myr}{p}{q} $ and $q \red{}{\myr}{p}{q_1} q_2 $.
In case $\hterm(q)=\hterm(q_2)$ we immediately get $p \red{}{\myr}{p}{q_2}$.
Otherwise $\hterm(q) \id \hterm(q_1)w$ implies $p \red{}{\myr}{p}{q_1}$.
Hence we have $p \red{}{\myr}{p}{\{ q_1, q_2 \}}$.
\remend
}
\end{remark}
Note that this property of prefix reduction corresponds to the fact
that the existence of prefix standard representations with respect to
a set of polynomials remains true for an interreduced version of the
set.
\begin{lemma}\label{lem.psr.prop3}~\\
{\sl
Let $F$ and $G$ be two sets of polynomials in $\myk[\m]$ such that
 every polynomial in $F$ has a prefix standard representation with
 respect to $G$.
Then if a polynomial $p$ has a prefix standard representation with
 respect to $F$ it also has one with respect to $G$. 
\lemend
}
\end{lemma}
\Ba{}~\\
Let $p = \sum_{i=1}^{n} \alpha_i \skm f_{i} \mrm w_i$,  with  $\alpha_i \in \myk^*,
 f_{i} \in F, w_i \in \m $ be a prefix standard representation of a
 polynomial $p$ with respect to the set of polynomials $F$, i.e., for
 all $1 \leq i \leq n$ we have $\hterm(p)\succeq \hterm(f_{i})w_i$.
Furthermore, every polynomial $f_i$ occurring in this sum has a prefix standard
 representation with respect to the set of polynomials $G$, say
 $f_i = \sum_{j=1}^{n_i} \beta_{i_j} \skm g_{i_j} \mrm v_{i_j}$,
 with  $\beta_{i_j} \in \myk^*,
 g_{i_j} \in G, v_{i_j} \in \m$ such that for all $1 \leq j \leq n_i$ we
 have $\hterm(f_i) \succeq \hterm(g_{i_j})v_{i_j}$. 
\\
These representations can be combined in the sum
 $$p = \sum_{i=1}^{n} \alpha_i \skm (\sum_{j=1}^{n_i} \beta_{i_j} \skm g_{i_j} \mrm
 v_{i_j}) \mrm w_i.$$
It remains to show that this in fact is a prefix standard
 representation with respect to $G$, i.e., to prove that for all 
 $1 \leq i \leq n$ and all $1 \leq j \leq n_i$,
 we get $\hterm(p) \succeq  \hterm(g_{i_j})v_{i_j}w_i$. 
\\
This now follows immediately as  for all 
 $1 \leq i \leq n$ and all $1 \leq j \leq n_i$ we have
 $$\hterm(p) \succeq \hterm(f_i)w_i \succeq \hterm(g_{i_j})v_{i_j}w_i.$$
\qed
One can even show that unique monic reduced standard bases exist.
\begin{lemma}\label{lem.prefixbasis}~\\
{\sl
Let $M$ be a subset of $\freemonoid$.
Then there exists a unique subset $M' \subseteq M$ such that
\begin{enumerate}
\item for all $m \in M$ there exists an element $m' \in M'$ and an
  element $w \in \freemonoid$ such that $m \id m'w$, and
\item for all $m \in M'$ no element $m' \in M' \backslash \{ m \}$ is
  a prefix of $m$.
\lemend\ohnebeweis
\end{enumerate}
}
\end{lemma}
%
Note that the subset $M'$ need not be finite , e.g., the set $M = \{ ab^ic | i \in \n \}$ in $\{a,b,c\}^*$
 contains no  finite subset satisfying the properties above.
\begin{theorem}\label{theo.monic.reduced.psb}~\\
{\sl
Every right ideal  in $\myk[\m]$ contains
 a unique monic reduced prefix standard basis.
\theoend
}
\end{theorem}
\Ba{}~\\
Let $\mswab{i}_r$ be a right ideal in $\myk[\m]$ and $G$ a subset of
$\mswab{i}_r$ such that we have
 $$\ideal{r}{\freemonoid}(\hterm(G)) \cap \m = \hterm(\mswab{i}_r \backslash \{
 0 \}).$$
Then by  theorem \ref{theo.pgb.psr} we
know that $G$ is a prefix standard basis of $\mswab{i}_r$.
By lemma \ref{lem.prefixbasis}, as the set $\hterm(G)$ is a  subset
of $\m$ which can be regarded as a subset of $\freemonoid$, there
exists a  subset $H \subseteq \hterm(G)$ such that
\begin{enumerate}
\item for all $m \in \hterm(G)$ there exists an element $m' \in H$ and an
  element $w \in \freemonoid$ such that $m \id m'w$, 
\item for all $m \in H$ no element $m' \in H \backslash \{ m \}$ is
  a prefix of $m$, and
\item $\ideal{r}{\freemonoid}(H) \cap \m = 
       \ideal{r}{\freemonoid}(\hterm(G)) \cap \m =
       \hterm(\mswab{i}_r \backslash \{  0 \}).$
\end{enumerate}
Further for each term $t \in H$ there exists at least one polynomial in
 $G$ with head term $t$.
Thus we can choose one of them, say $g_t$, for every
 $t \in H$.
If we then  set $G' = \{ g_t | t \in H \}$, by theorem
 \ref{theo.pgb.psr} this is a prefix
 standard basis.
Moreover,  all polynomials in $G'$ have different head terms and no head
term is prefix reducible by the other polynomials in $G'$.
Furthermore, if we prefix interreduce the set $G'$ giving us another set of polynomials
$G'' = \{ {\rm normalform}(g, \red{}{\myr}{p}{G' \backslash \{ g \}}) \mid g \in G' \}$, we know $\hterm(G') = \hterm(G'')$ and this set is a prefix
standard basis.
To see the latter, we show that for $f \in G$ with $f \red{}{\myr}{p}{g'
  \in G \backslash \{ f \}} f'$, the set $G' = (G \backslash \{ f \})
\cup \{f' \}$ is a prefix standard basis.
Since $f$ has a prefix standard representation with respect to $(G
 \backslash \{ f \}) \cup \{ f' \}$,  by lemma \ref{lem.psr.prop3} 
 we can conclude immediately that every polynomial in the
 right ideal generated by $G$ also has a prefix standard representation
 with respect to $G'$.
\\
It remains to show the uniqueness of the reduced prefix standard
basis in case it is monic.
Let us assume $S$ is another monic reduced prefix standard  basis of
 $\mswab{i}_r$.
Further let $f \in S \bigtriangleup G'' = (S \backslash G'') \cup
 (G'' \backslash S)$ be a polynomial such that $\hterm(f)$ is minimal in
 the set of terms $\hterm(S \bigtriangleup G'')$.
Without loss of generality we can assume that $f \in S \backslash
G''$.
As $G''$ is a reduced prefix standard  basis and $f \in \mswab{i}_r$ there
exists a polynomial $g \in G''$ such that $\hterm(f) \id \hterm(g)w$ for some $w
\in \m$.
We can even state that $g \in G'' \backslash S$ as otherwise $S$ would
not be prefix interreduced.
Since $f$ was chosen such that $\hterm(f)$ was minimal in
 $\hterm(S \bigtriangleup G'')$, we get $\hterm(f) = \hterm(g)$\footnote{Otherwise
   $\hterm(f) \succ \hterm(g)$ would contradict our assumption.}.
This gives us $\hterm(f-g) \prec \hterm(f) = \hterm(g)$ and $\hterm(f-g) \in \terms(f) \cup
 \terms(g)$ and without loss of generality let us assume $\hterm(f-g) \in \terms(f)$.
But $f-g \in \mswab{i}_r$ and $f-g \neq 0$ implies the existence of a polynomial $h \in
S$ such that $\hterm(f-g) \id \hterm(h)w$ for some $w \in \m$, implying that $f$
is not prefix reduced.
Hence we get that $S$ is not prefix interreduced, contradicting our
assumption.
\\
\qed
\auskommentieren{
A procedure to interreduce a set of polynomials with respect to prefix
reduction can be stated as in Mora's approach to free monoid rings.

\procedure{Prefix Interreduce}%
{\vspace{-4mm}\begin{tabbing}
XXXXX\=XXXX \kill
\removelastskip
{\bf Given:} \> A finite set $F \subseteq  \myk[\freemonoid]$.\\
{\bf Find:} \> $G$, a prefix interreduced set with $\red{*}{\lr}{p}{F}
= \red{*}{\lr}{p}{G}$.
\end{tabbing}
\vspace{-7mm}
\begin{tabbing}
XX\=XXXX\= XX \= XXXX\=\kill
$G$ := $F$; \\
{\bf while} there is $g \in G$ such that $\hterm(g)$ is prefix
reducible by $G \backslash \{ g \}$ {\bf do} \\
\> $G$ := $G \backslash \{ g \}$; \\
\> $f$ := ${\rm normalform}(g,\red{}{\myr}{p}{G})$; \\
\> {\kommentar \% compute a normal form using prefix  reduction.} \\
\> {\bf if} \>$f \neq 0$ \\
\> {\bf then} \>$G$ := $G \cup \{ f \}$; \\
{\bf endwhile}\\
$G$ := $\{ normalform(g, \red{}{\myr}{p}{G \backslash \{ g \}}) | g \in
G \}$
\end{tabbing}}

Notice that this procedure terminates as during the while loop only
polynomials with head terms smaller than the removed ones are added
and no cycles occur since no head term is added more than once.
Lemma \ref{lem.Noetherian} states that no infinite reduction sequence
with respect to the set of all computed polynomials can occur.
 ??????????????????????????????????????????????????????????????}
Before moving on to the study of prefix reduction we give some
 settings where
 bounds on special representations of polynomials are
 preserved under multiplication.
These properties are of importance to establish a weaker form of the
fact that $p \red{*}{\myr}{b}{F} 0$ implies $\alpha \skm p \mrm w
\red{*}{\myr}{b}{F} 0$ used in the proof of Buchberger's
characterization of Gr\"obner bases.
They will be used in the proofs of different characterizations of
prefix Gr\"obner bases later on.
\begin{lemma}\label{lem.redp}~\\
{\sl
Let $F$ be a set of polynomials in $\myk[\m]$ and $p \in \myk[\m]$. 
Further let $p \red{*}{\myr}{p}{F} 0$ and let us assume this reduction sequence results
 in a representation
 $p = \sum_{i=1}^{k} \alpha_i \skm g_i \mrm w_i$, where
 $\alpha_i \in \myk^*$, $g_i \in F$, and  $w_i \in \m$.
Then for every term $t \in \m$ such that $t \succ \hterm(p)$ and
every term $w \in \m$ we get that if
       $s \in \bigcup_{i=1}^k \terms(g_i \mrm w_i \mrm w)$
       then $tw \succ s$ holds.
\lemend
}
\end{lemma}
\Ba{}~\\
As $\sum_{i=1}^{k} \alpha_i \skm g_i \mrm w_i$ belongs to the reduction
 sequence $p \red{*}{\myr}{p}{F} 0$, for all $u \in \bigcup_{i=1}^{k} \terms(g_i \mrm w_i)$
 we have $\hterm(p) \succeq u$ implying $tw \succ \hterm(p)w \succeq uw \succeq u \mm
 w$.
Note that this proof uses the fact that the ordering $\succ$ on $\m$
 is induced by the completion ordering $\succeq_T$ of the presentation
 $(\Sigma, T)$ of $\m$,
 as we need that the ordering is compatible with concatenation, i.e.,
 $uv \succeq_T  (uv) \nf{T} = u \mm v$ for all $u,v \in \m$.
\\
\qed
Similarly, some properties of Buchberger's reduction are regained,
although not his lemma that $p \red{*}{\myr}{b}{F} 0$ implies $\alpha
\skm p \mrm w \red{*}{\myr}{b}{F} 0$.
\begin{lemma}\label{lem.redr}~\\
{\sl
Let $F$ be a set of polynomials in $\myk[\m]$ and $p,q \in \myk[\m]$. 
Further let $p \red{}{\myr}{p}{q} 0$ and $q \red{*}{\myr}{p}{F} 0$.
Let these reduction sequences result in the representations $p = \alpha \skm q \mrm w$ and 
 $q = \sum_{i=1}^{k} \alpha_i \skm g_i \mrm w_i$, where
 $\alpha,\alpha_i \in \myk^*, g_i \in F$, and $w,w_i \in \m$.
Then the following statements hold:
\begin{enumerate}
\item There exists $s \in \{ 1, \ldots , k \}$ such that 
       $\hterm(p)=\hterm(g_s \mrm w_s \mrm w) = \hterm( g_s \mrm w_s)w$.
\item For all remaining terms 
       $t \in \bigcup_{i=1 \atop i \neq s}^k \terms(g_i \mrm w_i \mrm w)$
       we have $\hterm(p) \succ t$.
\end{enumerate}
Notice that this also holds for a representation resulting from
 a right reduction sequence $p \red{*}{\myr}{r}{F} 0$.
\lemend
}
\end{lemma}
\Ba{}~\\
Since $p \red{}{\myr}{p}{q} 0$ and $p = \alpha \skm q \mrm w$ we know $\hterm(p) \id \hterm(q)w$.
As $\sum_{i=1}^{k} \alpha_i \skm g_i \mrm w_i$ belongs to the reduction
 sequence $q \red{*}{\myr}{p}{F} 0$,
 there exists $s \in \{ 1, \ldots , k \}$ such that
 $\hterm(q) \id \hterm(g_s)w_s \id \hterm(g_s \mrm w_s)$,
 $\hc(q)=\alpha_s \skm \hc(g_s)$ and for all terms $t \in \bigcup_{i=1
   \atop i \neq s}^{k} \terms(g_i \mrm w_i)$
 we have $\hterm(q) \succ t$ implying $\hterm(p) \id \hterm(q)w \succ tw \succeq t \mm
 w$.
\\
\qed
A word of caution:
This lemma does {\em not} imply $p \red{*}{\myr}{p}{F} 0$ or $p
\red{*}{\myr}{r}{F} 0$, as there is no information on how multiplying
the polynomials $g_i \mrm w_i$ by $w$ affects them in case $i \neq s$,
especially prefix reduction is not preserved under right multiplication.
\begin{example}~\\
{\rm
Let $\Sigma = \{ a, b \}$ and $T = \{ ab \myr \lambda, ba \myr \lambda \}$
be a presentation of a monoid $\m$ (which is in fact a group) with  a
length-lexicographical ordering induced by $a
\succ b$.
Further let $F = \{ a+ \lambda, b + \frac{1}{2} \} \subseteq \q[\m]$ and
$p=a^2+2 \skm a + 2$, $q = a + 2 \skm b + 2$ be two polynomials in
$\q[\m]$.
\\
Then we have
$$p \red{}{\myr}{p}{q} a^2 + 2 \skm a + 2 - (a + 2 \skm b + 2 ) \mrm a =
0$$
and 
$$q \red{}{\myr}{p}{a+\lambda} a+2 \skm b + 2 - (a + \lambda) = 2 \skm b
+ \lambda \red{}{\myr}{p}{b + \frac{1}{2}} 2 \skm b + \lambda - 2 \skm
(b + \frac{1}{2})
= 0,$$
but $p \red{}{\myr}{p}{a+\lambda} a^2 + 2 \skm a + 2 - (a + \lambda)
\mrm a = a+2 \red{}{\myr}{p}{a+\lambda} \lambda$, i.e., $p \nred{*}{\myr}{p}{F} 0$.
\exaend
}
\end{example}
\auskommentieren{
A similar lemma can be specified in case we use right reduction for
reducing $q$ to zero.
\begin{lemma}\label{lem.redr}~\\
{\sl
Let $F$ be a set of polynomials in $\myk[\m]$ and $p,q \in \myk[\m]$. 
Further let $p \red{}{\myr}{p}{q} 0$ and $q \red{*}{\myr}{{\bf r}}{F} 0$.
Let these reduction sequences result in the representations $p = \alpha \skm q \mrm w$ and 
 $q = \sum_{i=1}^{k} \alpha_i \skm g_i \mrm w_i$, where
 $\alpha,\alpha_i \in \myk^*, g_i \in F, w,w_i \in \m$.
Then the following statements hold:
\begin{enumerate}
\item There exists $s \in \{ 1, \ldots , k \}$ such that 
       $\hterm(p)=\hterm(g_s \mrm w_s \mrm w) = \hterm( g_s \mrm w_s)w$.
\item For all remaining terms 
       $t \in \bigcup_{i=1 \atop i \neq s}^k \terms(g_i \mrm w_i \mrm w)$
       we have $\hterm(p) \succ t$.
\lemend
\end{enumerate}
}
\end{lemma}
\Ba{}~\\
The proof is similar to the one of the previous lemma using that $\hterm(q) = \hterm(g_s) \mm
 w_s = \hterm(g_s \mrm w_s)$.
\\
\qed}
As before, we show that the translation lemma holds for prefix reduction.
\begin{lemma} \label{lem.confluentp}~\\
{\sl
Let $F$ be a set of polynomials and $p,q,h$ some
 polynomials in $\myk[\m]$.
\begin{enumerate}
\item
Let $p-q \red{}{\myr}{p}{F} h$.
Then there are  $p',q' \in \myk[\m]$ such that 
 $p  \red{*}{\myr}{p}{F} p', q  \red{*}{\myr}{p}{F} q'$ and $h=p'-q'$.
\item
Let $0$ be a normal form of $p-q$ with respect to $\red{}{\myr}{p}{F}$.
Then there exists a polynomial  $g \in \myk[\m]$ such that
 $p  \red{*}{\myr}{p}{F} g$ and $q  \red{*}{\myr}{p}{F} g$.
\lemend
\end{enumerate}
}
\end{lemma}
\Ba{}
\begin{enumerate}
\item  Let $p-q \red{}{\myr}{p}{F} h = p-q-\alpha \skm f \mrm w$, where
        $\alpha \in \myk^*, f \in F, w \in \m$
        and $\hterm(f)w = t$, i.e., $\alpha \skm \hc(f)$ is
        the coefficient of $t$ in $p-q$.
       We have to distinguish three cases:
       \begin{enumerate}
         \item $t \in \terms(p)$ and $t \in \terms(q)$:
               Then we can eliminate the term $t$ in the polynomials 
                $p$ respectively $q$ by prefix
                reduction and get
                $p \red{}{\myr}{p}{f} p - \alpha_1 \skm f \mrm w= p'$,
                $q \red{}{\myr}{p}{f} q - \alpha_2 \skm f \mrm w= q'$,
                with $\alpha_1  -  \alpha_2 =\alpha$,
                where $\alpha_1 \skm \hc(f)$ and 
                $\alpha_2 \skm \hc(f)$ are
                the coefficients of $t$ in $p$ respectively $q$.
         \item $t \in \terms(p)$ and $t \not\in \terms(q)$:
               Then  we can eliminate the term $t$ in the polynomial 
                $p$  by prefix reduction and  get
                $p \red{}{\myr}{p}{f} p - \alpha \skm f \mrm w= p'$ 
                and $q = q'$.
         \item $t \in \terms(q)$ and $t \not\in \terms(p)$:
               Then  we can eliminate the term $t$ in the polynomial 
                $q$ by prefix reduction and  get
                $q \red{}{\myr}{p}{f} q + \alpha \skm f \mrm w= q'$
                and $p = p'$.
       \end{enumerate}
      In all three cases we have $p' -q' =  p - q - \alpha \skm f \mrm w = h$.
\item We show our claim by induction on $k$, where $p-q \red{k}{\myr}{p}{F} 0$.
      In the base case $k=0$ there is nothing to show.
      Hence, let $p-q \red{}{\myr}{p}{F} h  \red{k}{\myr}{p}{F} 0$.
      Then by (1) there are polynomials $p',q' \in \myk[\m]$ such that 
       $p  \red{*}{\myr}{p}{F} p', q  \red{*}{\myr}{p}{F} q'$ and $h=p'-q'$.
      Now the induction hypothesis for $p'-q' \red{k}{\myr}{p}{F} 0$  yields 
       the existence of a polynomial $g \in \myk[\m]$ such that
       $p  \red{*}{\myr}{p}{F} p' \red{*}{\myr}{p}{F} g$ and
       $q  \red{*}{\myr}{p}{F} q' \red{*}{\myr}{p}{F} g$.
\\
\qed
\end{enumerate}\renewcommand{\baselinestretch}{1}\small\normalsize
Notice that prefix reduction like right reduction in general does not capture the right ideal
congruence, but for special bases
 of right ideals this can be regained.
\begin{lemma}\label{lem.psb.cong}~\\
{\sl
Let $F$ be a prefix standard basis and $p,q,h$ some polynomials in
 $\myk[\m]$.
Then
  $$p \red{*}{\lr}{p}{F} q \mbox{ if and only if } p - q \in
  \ideal{r}{}(F).$$
\lemend
}
\end{lemma}
\Ba{}~\\
In order to prove our claim we have to show two subgoals.
The inclusion $\red{*}{\lr}{p}{F}  \subseteq \;\; \equiv_{\ideal{r}{}(F)}$
 is an immediate consequence of the definition of prefix reduction and
 can be shown by induction as in lemma \ref{lem.strong.congruence}. 
To prove the converse inclusion 
 $\equiv_{\ideal{r}{}(F)} \: \subseteq \red{*}{\lr}{p}{F} $
 we can  modify the proof given 
 in lemma \ref{lem.strong.congruence}.
Remember that $p \equiv_{\ideal{r}{}(F)} q$ implies
 $p = q + \sum_{j=1}^{m} \alpha_{j} \skm f_j \mrm w_{j}$, where
 $\alpha_{j} \in \myk^*, f_j \in F, w_{j} \in \m$.
Since  every
 multiple $f_j \mrm w_{j}$ lies in $\ideal{r}{}(F)$ and 
 $F$ is a prefix standard basis, 
 by lemma \ref{lem.psr.prop1}, we can assume $\hterm(f \mrm w) \id \hterm(f)w$ for
 all polynomials occurring in the sum.
Now we can prove our claim straightforward as in lemma
\ref{lem.strong.congruence} by induction on $m$.
\auskommentieren{\\
In the base case $m = 0$ there is nothing to show.
\\
Let  $p = q + \sum_{j=1}^{m} \alpha_{j} \skm f_j \mrm w_{j} +
 \alpha_{m+1} \skm f_{m+1} \mrm w_{m+1}$
 and by our induction hypothesis
 $p \red{*}{\lr}{p}{F} q + \alpha_{m+1} \skm f_{m+1} \mrm w_{m+1}$.
\\
Let $t \id \hterm(f_{m+1})w_{m+1}$.
\\
In case $t \not\in \terms(q)$ we get $q +  \alpha_{m+1} \skm f_{m+1}
 \mrm w_{m+1}  \red{}{\myr}{p}{f_{m+1}} q$ and  are done.
\\
In case $t \not\in \terms(p)$ we get $p - \alpha_{m+1} \skm f_{m+1}
  \mrm w_{m+1}  \red{}{\myr}{p}{f_{m+1}} p$.
As $p  - \alpha_{m+1} \skm f_{m+1} \mrm w_{m+1} =
   q + \sum_{j=1}^{m} \alpha_{j} \skm f_j \mrm w_{j}$ the
 induction hypothesis yields 
 $p  - \alpha_{m+1} \skm f_{m+1} \mrm w_{m+1}
 \red{*}{\lr}{p}{F} q$ and hence we are done.
\\ 
Otherwise 
 let $\beta_1 \neq 0$ be the coefficient of $t$ in
 $q +  \alpha_{m+1} \skm f_{m+1} \mrm w_{m+1}$ and
 $\beta_2 \neq 0$ the coefficient of $t$ in $q$.
\\
This gives us a prefix reduction step

\hspace*{1cm}$q +  \alpha_{m+1} \skm f_{m+1} \mrm w_{m+1} \red{}{\myr}{p}{f_{m+1}}$\\
\hspace*{1cm}$q +  \alpha_{m+1} \skm f_{m+1} \mrm w_{m+1}
                  - \beta_1 \skm \hc(f_{m+1})^{-1}
                   \skm  f_{m+1} \mrm w_{m+1} =$  \\
\hspace*{1cm}$q - (\beta_1 \skm \hc(f_{m+1})^{-1} -\alpha_{m+1})
             \skm  f_{m+1} \mrm w_{m+1}$

 eliminating the occurrence of $t$ in $q +  \alpha_{m+1} \skm f_{m+1} \mrm w_{m+1}$.
\\
Then obviously 
  $\beta_2 = (\beta_1 \skm \hc(f_{m+1})^{-1}
                      -\alpha_{m+1}) \skm \hc(f_{m+1})$
  and, therefore, we have $q  \red{}{\myr}{p}{f_{m+1}}
             q - (\beta_1 \skm \hc(f_{m+1})^{-1} -\alpha_{m+1})
             \skm  f_{m+1} \mrm w_{m+1}$, i.e.,
  $q$ and $q +  \alpha_{m+1} \skm f_{m+1} \mrm w_{m+1}$ are joinable.}
\\
\qed
We can define Gr\"obner bases with respect to prefix reduction
 by slightly changing our previous definitions.
\begin{definition}\label{def.pgb}~\\
{\rm
A  set $G \subseteq \myk[\m]$ is called a \betonen{Gr\"obner basis}\/
 with respect to
 the reduction $\red{}{\myr}{p}{}$ or a \index{prefix!Gr\"obner basis}
 \index{Gr\"obner basis!prefix}
 \betonen{prefix Gr\"obner basis}, if
\begin{enumerate}
\item[(i)] $\red{*}{\lr}{p}{G}  = \;\; \equiv_{\ideal{r}{}(G)}$, and
\item[(ii)] $\red{}{\myr}{p}{G}$ is confluent.
\dend
\end{enumerate}
}
\end{definition}
As in the previous section there is a natural connection between prefix
standard bases and prefix reduction.
\begin{lemma}\label{lem.pprop}~\\
{\sl
Let $F$ be a set of polynomials and $p$ a non-zero polynomial in $\myk[\m]$.
\begin{enumerate}
\item Then $p \red{*}{\myr}{p}{F} 0$ implies the existence of a prefix 
  standard representation for $p$.
\item In case $p$ has a prefix standard representation with respect to $F$,
  then $p$ is prefix reducible at its head monomial by $F$,
  i.e., $p$ is prefix top-reducible by $F$.
\item\label{lem.pprop.3} In case $F$ is a prefix standard basis,
  every non-zero polynomial $p$ in $\ideal{r}{}(F)\backslash \{ 0 \}$
  is prefix top-reducible to zero by $F$.
\lemend
\end{enumerate}
}
\end{lemma}
\Ba{}
\begin{enumerate}
\item This follows directly by adding up the polynomials used in the
  prefix reduction steps  occurring in $p \red{*}{\myr}{p}{F} 0$.
\item This is an immediate consequence of definition \ref{def.psr}
  as the existence of a polynomial $f$ in $F$ and an element
  $w\in\m$ with $\hterm(p) \id \hterm(f)w$ is guaranteed.
\item We show that every non-zero polynomial $p \in
  \ideal{r}{}(F)\backslash \{ 0 \}$ is top-reducible to zero using $F$ by induction on the term $\hterm(p)$.
  Let $\hterm(p) = \min \{ \hterm(g) | g \in \ideal{r}{}(F)\backslash \{ 0 \} \}$.
  Then, as $p \in \ideal{r}{}(F)$ and $F$ is a prefix  standard basis, we
   have $p = \sum_{i=1}^{k} \alpha_i \skm f_{i} \mrm w_i$, with 
   $\alpha_i \in \myk^*, f_{i} \in F, w_i \in \m$ and 
   $\hterm(p) \succeq \hterm(f_{i})w_i$
   for all  $1 \leq i \leq k$.
  Without loss of generality, let us assume  $\hterm(p) \id \hterm(f_1)w_1$.
  Hence, the polynomial $p$ is  prefix reducible by $f_1$.
  Let $p \red{}{\myr}{r}{f_1} q$, i.e.,
   $q = p - \hc(p) \skm \hc(f_1)^{-1} \skm f_1 \mrm w_1$, and
   by the definition of prefix reduction the term $\hterm(p)$ is
   eliminated from $p$ implying that $\hterm(q) \pred \hterm(p)$ as $q < p$.
  Now, as $\hterm(p)$ was minimal among the head terms of the elements
   in the right ideal generated by $F$, this implies $q=0$, and, 
   therefore, $p$ is prefix top-reducible to zero by $f_1$ in 
   one step.
   On the other hand, in case 
    $\hterm(p) \succ \min \{ \hterm(g) | g \in \ideal{r}{}(F) \backslash \{ 0 \} \}$, by the
    same arguments used before we can prefix reduce the polynomial 
    $p$ to a polynomial
    $q$ with $\hterm(q) \pred \hterm(p)$, and, thus, by our induction
    hypothesis we know that $q$ is prefix top-reducible to zero.
   Therefore, as the reduction step $p \red{}{\myr}{p}{f_1} q$ takes
    place at the head term of $p$, the polynomial 
    $p$ is also prefix top-reducible
    to zero.
\\
\qed
\end{enumerate}\renewcommand{\baselinestretch}{1}\small\normalsize
Using the results of this lemma we can show that prefix standard bases
 in fact are  prefix Gr\"obner bases.
\begin{theorem}\label{theo.pgb.psr}~\\
{\sl
For a set $F$ of polynomials in $\myk[\m]$,
 the following statements are equivalent:
\begin{enumerate}
\item $F$ is a prefix Gr\"obner basis.
\item For all polynomials $g \in \ideal{r}{}(F)$ we have $g \red{*}{\myr}{p}{F} 0$.
\item $F$ is a prefix standard basis.
\theoend
\end{enumerate}
}
\end{theorem}
\Ba{}~\\
\mbox{$1 \R 2:$ }
  By (i) of definition \ref{def.pgb} we know that $g \in
  \ideal{r}{}(F)$ implies $g \red{*}{\lr}{p}{F} 0$ and since
  $\red{}{\myr}{p}{F}$ is confluent and $0$ is irreducible,
  $g \red{*}{\myr}{p}{F} 0$ follows
  immediately.

\mbox{$2 \R 3:$ }
  This follows directly by adding up the polynomials used in the
  prefix reduction steps  $g \red{*}{\myr}{p}{F} 0$.

\mbox{$3 \R 1:$ }
  In order to show that $F$ is a prefix Gr\"obner basis we have
  to prove two subgoals:
$\red{*}{\lr}{p}{F}   = \;\; \equiv_{\ideal{r}{}(F)}$ has already been shown in lemma \ref{lem.psb.cong}.
It remains to show that $\red{}{\myr}{p}{F}$ is confluent.
             Since $\red{}{\myr}{p}{F}$ is Noetherian, we only have to prove
             local confluence.
             Suppose $g \red{}{\myr}{p}{F} g_1$, 
              $g \red{}{\myr}{p}{F} g_2$ and $g_1 \neq g_2$.
             Then $g_1 - g_2 \in \ideal{r}{}(F)$ and, therefore, is prefix 
              top-reducible to zero as a result of lemma 
              \ref{lem.pprop}. 
             Hence lemma \ref{lem.confluentp} provides the existence
              of a polynomial $h \in \myk[\m]$ such that
              $g_1 \red{*}{\myr}{p}{F} h$ and $g_2 \red{*}{\myr}{p}{F} h$,
              i.e.,  $\red{}{\myr}{p}{F}$ is confluent.
\\
\qed
%
Since in general for a set of polynomials $F$ we get
 $\red{*}{\lr}{r}{F} \neq \red{*}{\lr}{p}{F} \neq \red{*}{\lr}{s}{F}$,
 we again enrich our set of polynomials
 used for reduction in order to regain the expressiveness of strong 
 right reduction respectively right reduction combined with saturation.
\begin{definition}\label{def.sat.p}~\\
{\rm
A set of polynomials $F\subseteq \{\alpha \skm p \mrm w \mid \alpha
\in \myk^*, w \in\m \}$
 is called a 
 \index{prefix!saturating set}\index{saturating set!prefix}\index{prefix!saturation}\index{saturation!prefix}\betonen{prefix saturating set}\/ for a non-zero polynomial $p\in
 \myk[\m]$, if for all $\alpha \in \myk^*$, $w \in \m$, in case $\alpha
 \skm p \mrm w \neq 0$ then $\alpha \skm p
 \mrm w  \red{}{\myr}{p}{F} 0$ holds\footnote{Since $\myk$ is a
   field it is sufficient to demand $p \mrm w \red{\leq 1}{\myr}{p}{F} 0$ for
   all $w \in \m$.}.
$\SAT_p(p)$ denotes the family of all prefix saturating sets for $p$.
We call a set $F \subseteq \myk[\m]$ \index{prefix!saturated}\index{saturated set!prefix}\betonen{prefix saturated}, if for all
 $f \in F$ and all $\alpha \in \myk^*$,  $w \in \m$, $\alpha \skm f \mrm
 w \red{}{\myr}{p}{F} 0$ holds in case $\alpha \skm f \mrm w \neq 0$.
\dend
}
\end{definition}
As in the case of right reduction, for a set of polynomials $F$ every
 union $S = \bigcup_{f \in F} S_f$ of prefix saturating sets $S_f \in \SAT_p(f)$
 is a prefix saturated set.
But in general this union contains too many polynomials, as the
 following example shows.
\begin{example}~\\
{\rm
Let $\Sigma = \{ a,b \}$ and $T = \{  ba \myr ab \}$ be a presentation
 of a monoid $\m$   with  a  length-lexicographical ordering induced
 by $b \succ a$.
\\
For the set $F = \{ a,b \}$ we find that $b$ has no finite prefix
 saturating set, 
 since there exists no finite set
 $S \subseteq \{ a^nb \mid n \in \n \}$ such that 
 all polynomials $a^nb$, $n \in \n$ are
 prefix reducible to zero in one step using $S$.
But the set $F$ itself is prefix saturated, since obviously
 $a^nb \red{}{\myr}{p}{a} 0$ holds.
\exaend
}
\end{example}
At this point, before we continue to give a characterization of
 Gr\"obner bases in this context, let us take a look
 at the relations between the reductions studied so far and the
 concepts of saturation induced by them.
Remember that saturation enabled us to simulate strong right
 reduction by right reduction (compare lemma \ref{lem.sconnection}).
The same is true for prefix saturation.
\begin{lemma}\label{lem.connection}~\\
{\sl
Let $f,g,p$ be some polynomials in $\myk[\m]$, $S \in\SAT(p)$,
 and $S_p \in\SAT_p(p)$.
Then the following statements hold: 
\begin{enumerate}
\item $f  \red{}{\myr}{r}{S} g$ if and only if $f  \red{}{\myr}{r}{S_p} g$.
\item $f  \red{}{\myr}{r}{S} g$ if and only if $f  \red{}{\myr}{p}{S_p}
  g$.
\lemend
\end{enumerate}
}
\end{lemma}
\Ba{}~\\
We will only prove statement 1, as  2 can be shown analogously.
First, suppose $f  \red{}{\myr}{r}{p_1 \in S} g$, i.e.,
             $g = f - \alpha_1 \skm p_1 \mrm w_1$
             for some $\alpha_1 \in \myk^*, w_1 \in \m$.
            Since $S_p$ is a prefix saturating set for $p$ and $p_1
            \in S$ implies $p_1 = \alpha \skm p \mrm u$ for some
            $\alpha \in \myk^*$, $u \in \m$ this
            gives us $p_1 \mrm w_1 = (\alpha \skm p \mrm u) \mrm w_1 =
            \alpha \skm p \mrm (u  \mm w_1)$, i.e., $p_1 \mrm w_1$ is
            a multiple of $p$.
            Thus
             $p_1 \mrm w_1  \red{}{\myr}{p}{p_2 \in S_p} 0$ and 
             $p_1 \mrm w_1 = \beta \skm p_2 \mrm w_2$ for some $\beta
             \in \myk^*$, $w_2 \in \m$.
             Furthermore, as 
             $\hterm(p_1 \mrm w_1) = \hterm(p_1) \mm w_1 = \hterm(p_2)w_2$,
             we get
             $f  \red{}{\myr}{r}{p_2 \in S_p} g$ and even $f \red{}{\myr}{p}{p_2} g$.
On the other hand, suppose  $f  \red{}{\myr}{r}{p_1 \in S_p} g$,
             i.e., $g =  f - \alpha_1 \skm p_1 \mrm w_1$
             for some $\alpha_1 \in \myk^*, w_1 \in \m$ and
             $\hterm(p_1 \mrm w_1) = \hterm(p_1) \mm w_1$.
            As $p_1 = \alpha \skm p \mrm u$ for some $\alpha \in
            \myk^*$, $u \in \m$  and $p_1 \mrm w_1 = (\alpha \skm p \mrm u) \mrm w_1 =
            \alpha \skm p \mrm (u  \mm w_1)$ i.e., $p_1 \mrm w_1$ is
            a multiple of $p$.
            Hence we get
            $p_1 \mrm w_1  \red{}{\myr}{r}{p_2 \in S} 0$ and $p_1 \mrm w_1
            = \beta \skm p_2 \mrm w_2$ for some $\beta
             \in \myk^*$, $w_2 \in \m$ with
             $\hterm(p_1 \mrm w_1)=\hterm(p_2) \mm w_2$ implying
             $f  \red{}{\myr}{r}{p_2 \in S} g$.
\\
\qed
%
Note that prefix saturated sets are also saturated sets.
Furthermore, they give us additional information as they
 allow special representations of elements
 in the right ideals they generate which are weaker than
 prefix standard representations,
 but sufficient to give a localized confluence criteria.
The following lemma is an analogon to lemma \ref{lem.prop1}.
\begin{lemma}\label{lem.prop1.1}~\\
{\sl
 Let $F \subseteq \myk[\m]$ be a prefix saturated set.
 Then every non-zero polynomial $g$ in $\ideal{r}{}(F)$ has a
 representation\footnote{Note that such a representation need not be a
   prefix standard representation as we cannot conclude that
   $\hterm(g) \succeq \hterm(f_i)w_i$ holds.} of the form
   $g =  \sum_{i=1}^{k} \alpha_i \skm f_i \mrm w_i$
 with $\alpha_i \in \myk^*, f_i \in F, w_i \in \m$,
 and
 $\hterm(f_i \mrm w_i) \id \hterm(f_i)w_i$.
\lemend\ohnebeweis
}
\end{lemma}
\vspace{-1mm}
Prefix reduction combined with prefix saturation is strong
enough to capture the right ideal congruence.
\vspace{-1mm}
\begin{lemma}\label{lem.congruencep}~\\
{\sl
Let $F$ be a prefix saturated set of polynomials in
 $\myk[\m]$ and $p,q,h \in\myk[\m]$.
Then
  $$p \red{*}{\lr}{p}{F} q \mbox{ if and only if } p - q \in
  \ideal{r}{}(F).$$
\lemend
}
\end{lemma}
\Ba{}~\\
This lemma follows directly from theorem \ref{theo.congruence}  and lemma
\ref{lem.connection}.
\\
\vspace{-2mm}
\qed
In the following we will give a procedure, which similar to
{\sc Saturation 1}  on page
\pageref{saturation2} enumerates a 
 prefix saturating set for a polynomial in $\myk[\m]$ depending on a
 convergent presentation $(\Sigma,T)$ of $\m$.
We compute critical situations between the head terms of
multiples of the  polynomial being saturated and the left hand sides of the
rules in $T$.
Later on we will see how this can be compared to computing special
s-polynomials between polynomials and the set of ``polynomials'' $\{
l-r \mid (l,r) \in T \}$ in the free monoid ring generated by $\Sigma$.

\procedure{Prefix Saturation\protect{\label{prefix.saturation}}}%
{\vspace{-4mm}\begin{tabbing}
XXXXX\=XXXX \kill
\removelastskip
{\bf Given:} \> A polynomial $p \in \myk[\m]$ and
             $(\Sigma, T)$ a convergent  presentation of $\m$. \\
{\bf Find:} \> $S_p \in\SAT_p(p)$. 
\end{tabbing}
\vspace{-7mm}
\begin{tabbing}
XX\=XX\=XX\=XXXX\=XXX\= \kill
$S_p$ := $\{p \}$; \\
$H$ := $\{p \}$; \\
{\bf while} $H \neq \emptyset$ {\bf do} \\
\> $q$ := ${\rm remove}(H)$; \\
\> {\rm\kommentar \% Remove an element using a fair
  strategy} \\
\> $t$ := $\hterm(q)$; \\
\> {\bf for all} $w \in C(t) = \{ w \in \Sigma^* \mid tw \id t_{1}t_{2}w \id t_{1}l, t_2 \neq \lambda$
                            {\rm for some} $(l, r) \in T \} $ {\bf do}
                            \\
\> \>{\rm\kommentar \% $C(t)$ contains special overlaps between  $t$ and left hand sides of  rules in $T$} \\
\> \> $q'$ := $q \mrm w$; \\
\> \> {\bf if} \> $q' \nred{}{\myr}{p}{S_p} 0$ and $q' \neq 0$\\
\> \>          \> {\bf then} \>  $S_p$ := $S_p \cup \{q' \}$; \\
\> \>           \>           \>  $H$ := $H \cup \{q' \}$; \\
\> \> {\bf endif}\\
\> {\bf endfor} \\
{\bf endwhile}
\end{tabbing}}

Notice that in contrary to procedure {\sc Saturation 1} in case we have $q' \red{}{\myr}{p}{S} 0$ then
$q'$ does not have to be considered for further computations.
\begin{theorem}\label{theo.cor.ps}~\\
{\sl
For a given polynomial $p \in \myk[\m]$ let $S_p$ be the
 set generated by procedure {\sc Prefix Saturation}.
Then for all elements $w \in \m$ the polynomial
 $p \mrm w$ is prefix reducible to zero in one step using $S_p$ in
 case it is non-zero.
\theoend
}
\end{theorem}
\Ba{}~\\
We show that for all $q \in S_p$,$w \in \m$ we have
 $q \mrm w  \red{}{\myr}{p}{S_p} 0$ in case $q \mrm w \neq 0$.
Suppose this is not true.
Then we can choose a  non-zero counter-example
 $q \mrm  w$, where $\hterm(q)w$ is minimal
 (according to the ordering $\succeq_T$ on $\Sigma^*$) 
 and $q \mrm w \nred{}{\myr}{p}{S_p} 0$.
Thus $\hterm(q)w$ must be $T$-reducible, as otherwise
 $q \mrm w  \red{}{\myr}{p}{q \in S_p} 0$.
Let $\hterm(q)w \id t_1t_2w_1w_2$ such that
 $\hterm(q) \id t_1t_2, t_2 \neq \lambda, w \id w_1w_2$ and
 $l \id t_2w_1$ for some $(l,r) \in T$.
Furthermore, $w_1 \in \m$ as it is a prefix of $w \in \m$.
Since $q \in S_p$ the polynomial $q$ must have been
 added to the set $H$ at some step and as we use a fair strategy 
 to remove elements from $H$, $q$ and $C(\hterm(q))$ are considered.
Thus, we have $w_1 \in C(\hterm(q))$ by the definition of this set and
we can distinguish two cases.
If we have $q \mrm w_1 \in S_p$ then 
       $q \mrm w = (q \mrm w_1) \mrm w_2  \red{}{\myr}{p}{S_p} 0$, since
       $w_1 \in \m$ and 
       $\hterm(q)w \id \hterm(q)w_1w_2 \succ \hterm(q \mrm w_1)w_2$,  contradicting our assumption.
On the other hand, $q \mrm w_1 \not\in S_p$ implies
       $q \mrm w_1   \red{}{\myr}{p}{q' \in S_p} 0$ and we know
       $\hterm(q)w_1 \succ \hterm(q \mrm w_1) \id \hterm(q')z$ for some $z \in\m$.
      Further  
       $q \mrm w = (q \mrm w_1) \mrm w_2 = (\alpha \skm q' \mrm z) \mrm w_2$,
       and $\hterm(q)w \succ \hterm(q')zw_2 \succeq \hterm(q')(z \mm
       w_2)$.
      Therefore, we have 
       $q \mrm w = (\alpha \skm q' \mrm z) \mrm w_2 =
        \alpha \skm q' \mrm (z \mm w_2)  \red{}{\myr}{p}{S_p} 0$,
       contradicting our assumption.
\\
\qed
%
Hence, procedure {\sc Prefix Saturation} enumerates a prefix
saturating set for a polynomial. 
The next lemma states that this process will terminate in case a
finite prefix saturating set exists.
\begin{lemma}\label{lem.sattermination}~\\
{\sl
In case a polynomial has a finite prefix saturating set, then 
 procedure {\sc Prefix Saturation} terminates.
\lemend
}
\end{lemma}
\Ba{}~\\
Let $p \in\myk[\m]$ be the  polynomial which is being saturated and
 $S \in \SAT_p(p)$ finite.
Further let $S_p$ be the set generated by the procedure.
Since we have a correct enumeration of a prefix saturating set for $p$,
 each polynomial $q \in S$ has to be prefix reducible to zero by a 
 polynomial in $S_p$\footnote{Especially there is a polynomial
  $q' \in S_p$ such that $\hterm(q) \id \hterm(q')z$ for some $z \in\m$.}.
Therefore, there exists a finite set $S' \subseteq S_p$ such that for
 every polynomials $q \in S$ there exists a polynomial $q' \in S'$
 such that $q \red{}{\myr}{p}{q'} 0$.
Thus as soon as all polynomials in $S'$ have been enumerated, 
 we have the situation that for every remaining polynomial $h \in H$ 
 on one hand 
 $h \red{}{\myr}{p}{s \in S} 0$ and on the other hand 
 $s \red{}{\myr}{p}{s' \in S'} 0$ hold, implying $h \red{}{\myr}{p}{S'} 0$.
Hence the {\bf while} loop terminates, as no more elements are added to the
 set $H$.
\\
\qed
The following lemma gives some more information on the
 structure of a prefix saturating set for
 a polynomial in case our monoid is presented by a
 convergent monadic semi-Thue system with a 
 length-lexicographical completion ordering ensuring that in this case
 finite prefix saturating sets exist.
\begin{lemma} \label{lem.prop2}~\\
{\sl
Let $(\Sigma,T)$ be a  convergent monadic presentation of a monoid
 $\m$.
For a non-zero polynomial $p$ in $\myk[\m ]$,
 let $S \in \SAT_p(p)$.
Then for each non-zero right multiple $\tilde{q} = p \mrm w$, $w \in \m$
 there is a $q \in S$ such that 
\begin{enumerate}
\item $|\hterm(q)| \leq |\hterm(p)| + \max  \{ |l| \mid (l,a) \in T \} -1 =:  K$.
\item $\tilde{q} \red{}{\myr}{p}{q} 0$.
\lemend
\end{enumerate}
}
\end{lemma}
\Ba{}~\\
Let us assume that our ordering on $\m$ is length-lexicographical.
As $T$ is monadic, for a polynomial $p = \sum_{i = 1}^{n} \alpha_{i}
\skm t_{i}$ the right multiplication $p \mrm w$ results in the terms
 of the form $t_1 \mm w \id t_1'a_1w_1 , \ldots , t_n \mm w \id t_n'a_nw_n $,
 where $t_{i}'$ is a prefix of $t_{i}$, 
 $a_i \in \Sigma \cup \{ \lambda \}$ and 
 $w_1, \ldots , w_n$ are (possibly empty) suffixes of $w$.
Now let us assume there exists a polynomial $\tilde{q}=p \mrm w$
 with $|\hterm(\tilde{q})| > K$.
Then we can decompose $w \id w_1w_2w_3$ in such a way that
 $\hterm(\tilde{q}) = t_j \mm w \id t_j'a_jw_2w_3$, i.e., $t_j \mm w_1 =
 t'_ja_j$ and $|t_j'a_jw_2| = K$.
Let us consider the polynomial $p \mrm w_1w_2$.
We claim that $\hterm(p \mrm w_1w_2) \id t_j'a_jw_2$.
Suppose this is not true.
Then there exists a term $t_i \in \terms(p)$ such that
 $\hterm(p \mrm w_1w_2)= t_i \mm w_1w_2 \succ t_j'a_jw_2$.
Hence  $|t_i \mm w_1w_2| \geq K$ gives us
 $t_i \mm w_1w_2 \id t_i'a_iw'$, where $|t_i'a_i| \leq |\hterm(p)|$
 and thus $|w'| \geq \max \{ |l| \mid (l,a) \in T \} -1$.
Since furthermore $w'$ is a suffix of $w_1w_2$ and $w_1w_2w_3 \in IRR(T)$
 we get that $t_i'a_iw'w_3$ is $T$-irreducible
 giving us $t_i'a_iw'w_3 \succ t_j \mm w_1w_2w_3 = \hterm( \tilde{q})$
 contradicting our assumption.
Therefore, we can conclude $\hterm(p \mrm w_1w_2) \id t_j'a_jw_2$.
In case $p \mrm w_1w_2 \in S$ we can set $q= p \mrm w_1w_2$
 and get $|\hterm(q)| = |\hterm(p \mrm w_1w_2)| = |t'_ja_jw_2| = K$ and we are
 done.
On the other hand, since $S$ is a prefix saturating set for $p$,
 there exists an element $q\in S$ such that $p \mrm
 w_1w_2 \red{}{\myr}{p}{q} 0$ and $|\hterm(q)| \leq |\hterm(p \mrm w_1w_2)| =
 K$.
As $\hterm(p \mrm w) \id \hterm(p \mrm w_1w_2)w_3$ we know $p \mrm w
 \red{}{\myr}{p}{q} 0$ and are done.
\\
\qed
\begin{corollary}~\\
{\sl
Procedure {\sc Prefix Saturation} terminates for  monoids $\m$ with
a convergent monadic presentation.
}
\end{corollary}
\Ba{}~\\
This follows immediately from lemma \ref{lem.sattermination} since
lemma \ref{lem.prop2} provides the existence of finite prefix
saturating sets for polynomials.
\\
\qed

%
%
\begin{corollary}~\\
{\sl
Procedure {\sc Prefix Saturation} terminates for finite monoids $\m$.
}
\end{corollary}
\Ba{}~\\
Let $\m$ be a finite monoid and $p$ a polynomial in $\myk[\m]$.
Then obviously the set $S= \{ p \mrm w \mid w \in \m \}$ is finite and hence the
procedure must terminate.
\\
\qed
The next lemma states the existence of minimal prefix saturating sets.
\begin{lemma}~\\
{\sl
Let $p$ be a polynomial in $\myk[\m]$ and $S \in \SAT_p(p)$ a prefix saturating
 set for $p$.
Then if there is a polynomial $q \in S$ such that $q \red{}{\myr}{p}{S
 \backslash \{ q \}} 0$, the set $S \backslash \{ q \}$ is a prefix saturating
 set for $p$.
\lemend
}
\end{lemma}
\Ba{}~\\
This follows immediately as $p \red{}{\myr}{p}{q_1} 0$ and $q_1
\red{}{\myr}{p}{q_2} 0$ imply $p \red{}{\myr}{p}{q_2} 0$ (compare item 3
of lemma \ref{lem.red} which also holds for prefix reduction). 
\\
\qed
It is now possible to introduce simplification to procedure
 {\sc Prefix Saturation} by removing polynomials
 which are prefix reducible to zero in one step by later
 generated polynomial multiples.

\procedure{Prefix Saturation using Simplification\protect{\label{prefix.sauration.using.simplification}}}%
{\vspace{-4mm}\begin{tabbing}
XXXXX\=XXXX \kill
\removelastskip
{\bf Given:} \> A polynomial $p \in \myk[\m]$, and \\
             \> $(\Sigma, T)$ a convergent semi-Thue system presenting
                $\m$. \\
{\bf Find:} \> $S_p \in\SAT_p(p)$. 
\end{tabbing}
\vspace{-7mm}
\begin{tabbing}
XX\=XX\=XX\=XXXX\=XXX\= \kill
$S_0 := \{p \}$; \\
$H := \{p \}$; \\
$i$ := $0$;\\
{\bf while} $H \neq \emptyset$ {\bf do} \\
\> $i$ := $i+1$; \\
\> $S_i$ := $S_{i-1}$; \\
\> $q$ := ${\rm remove}(H)$;\\
\> {\rm\kommentar \% Remove an element using a fair
  strategy} \\
\> $t := \hterm(q)$; \\
\> {\bf for all} $w \in C(t) = \{ w \in \Sigma^* \mid tw \id t_{1}t_{2}w \id t_{1}l, t_2 \neq \lambda$
                            {\rm for some} $(l, r) \in T \} $ {\bf do}
                            \\
\> \>{\rm\kommentar \% $C(t)$ contains special overlaps between
  $t$ and left hand sides of rules in $T$} \\
\> \> $q' := q \mrm w$; \\
\> \> {\bf if} \> $q' \nred{}{\myr}{p}{S_i} 0$ and $q' \neq 0$\\
\> \> \>{\bf then} \>  $S_i := {\rm simplify}(S_i, q') \cup \{q' \}$; \\
\> \> \>           \>  {\rm\kommentar \% 
                     Simplify removes elements $s$ from $S_i$ in
                     case $s \red{}{\myr}{p}{q'} 0$}\\
\> \> \>           \>  $H := H \cup \{q' \}$; \\
\>\> {\bf endif} \\
\> {\bf endfor} \\
{\bf endwhile} \\
$S_p$ := $S_i$
\end{tabbing}}

\begin{theorem}\label{theo.cor.psws}~\\
{\sl
Let $S_p$ be the set generated by procedure {\sc Prefix
  Saturation using Simplification} for a given polynomial $p \in
 \myk[\m]$.
Then for all elements $w \in \m$ the polynomial
 $p \mrm w$ is prefix reducible to zero using $S_p$.
\theoend
}
\end{theorem}
\Ba{}~\\
First we specify the output of procedure {\sc Prefix Saturation using
  Simplification}.
In case $H$ becomes empty in some iteration $k$ we have $S_p = S_k$.
Otherwise, as a fair strategy is used to remove elements from the set
$H$, this
 guarantees that all polynomials added to $H$ are also considered and
 we can characterize the output as
 $S_p = \bigcup_{i \geq 0} \bigcap_{j \geq i} S_i$.
\\
The fact that no polynomial is entered twice into $H$ is due to the
 following observation:
Since for polynomials $q, q_1, q_2$, $q \red{}{\myr}{p}{q_1} 0$ and $q_1
 \red{}{\myr}{p}{q_2} 0$ yield $q \red{}{\myr}{p}{q_2} 0$, we get that for
 all polynomials 
 $f \red{}{\myr}{p}{S_i} 0$ implies $f \red{}{\myr}{p}{S_{i+n}} 0$, $n \in \n$.
In particular, as $p \in S_0$, we get $p \red{}{\myr}{p}{S_i} 0$ for all
$i \in \n$.
\auskommentieren{
To see this let us look at $f \red{}{\myr}{p}{S_i} 0$ where
 $f = \alpha \skm s \mrm u$ for some
 $\alpha \in \myk^*, u \in \m$ and $\hterm(f) \id \hterm(s)u$.
In case $s \in S_{i+1}$ we are done.
Else $s$ is removed from the set $S_i$ by simplification in the {\bf for
 all} loop, i.e., there exists a polynomial $s' \in S_i \backslash \{ s
 \}$ such that $s \red{}{\myr}{p}{s'} 0$ with $s = \alpha' \skm s' \mrm
  u'$ where $\alpha' \in \myk^*$, $u' \in \m$ and $\hterm(s) \id \hterm(s')u'$.
Then either $s' \in S_{i+1}$ or another similar simplification step
 takes place.
Since the for {\bf all loop} terminates, without loss of generality we can
 proceed assuming $s' \in S_{i+1}$.
\\
We find
$$f = \alpha \skm s \mrm u = \alpha \skm ( \alpha' \skm s' \mrm u')
\mrm u = (\alpha \skm \alpha') \skm s' \mrm (u' \mm u)$$
and
$$\hterm(f) \id \hterm(s)u \id (\hterm(s')u')u \id \hterm(s')u'u,$$
and thus $f \red{}{\myr}{p}{s' \in S_{i+1}} 0$.
\\
In particular, $p \in S_0$ implies $p \red{}{\myr}{p}{S_i} 0$ for all $i
\in \n$.
\\}
\\
We continue by proving that for all $q \in \bigcup_{i \geq 0} S_i$,
 $w \in \m$ we have
 $q \mrm w  \red{}{\myr}{p}{S_p} 0$ as this implies $p \mrm v
 \red{}{\myr}{p}{S_p} 0$ for all $v \in \m$.
Suppose this is not true.
Then we can choose a counter-example $q \mrm  w$
 such that $\hterm(q)w$ is minimal
 (according to the ordering $\succeq_T$ on $\Sigma^*$) 
 among all counter-examples and $q \mrm w \nred{}{\myr}{p}{S_p} 0$.
Then $\hterm(q)w$ must be $T$-reducible, as otherwise $q \in S_j$ for
some $j \in \n$ and $q \not\in S_p$ implies $q \red{}{\myr}{p}{s \in S_{j+k}}
0$ for some $k \in \n$ and either $s \in S_p$ contradicting our
assumption  or $s$ is again removed
by simplification.
Now the latter cannot occur infinitely often, as in case a polynomial
$s$ is removed due to simplification with a polynomial $q'$ we
know $\hterm(q')$ is a proper prefix of $\hterm(s)$ as otherwise we
would have $q' \red{}{\myr}{p}{s} 0$ contradicting the fact that $q'$ is
used for simplification of the set containing $s$.
Therefore, the existence of a polynomial $s \in S_p$ such that $q
\mrm w \red{}{\myr}{p}{s} 0$ is guaranteed contradicting our assumption. 
\auskommentieren{
We will first show that then
 $\hterm(q)w$ must be $T$-reducible.
\\
Suppose that this is not true.
Now $q \red{}{\myr}{p}{S_l} 0$ for some $l \in \n$ as $q \in \bigcup_{i
  \geq 0} S_i$ and we can conclude $q \red{}{\myr}{p}{S_{l + j}} 0$ for
  all $j \in \n$.
In case $S_p = S_k$ for some $k \in \n$ this would give us $q
 \red{}{\myr}{p}{S_p} 0$ implying $q \mrm w \red{}{\myr}{p}{S_p} 0$
 contradicting our assumption.
It remains to look at the case $S_p = \bigcup_{i \geq 0} \bigcap_{j
  \geq i} S_i$.
Let us assume that there is no polynomial $q' \in S_p$ such that $q
\red{}{\myr}{p}{q'} 0$.
This implies that there are infinitely many different polynomials
$q_s$, $s \in \n$ (as no polynomial is entered twice to $H$ or a new $S_i$)
 such that $q \red{}{\myr}{p}{q_s} 0$ and $q_s \in S_{l+j_s}$, $q_s
 \not\in S_{l+j_s+i}$ for some $j_s \in \n$, and all $i \in \n^+$.
Since all $q_s$ are right multiples of $p$, say $q_s = p \mrm u_s$ for
some $u_s \in \m$, and only finitely many different prefixes of
$\hterm(q)$ can occur among the head terms of the polynomials $q_s$,
one head term must occur infinitely many times resulting from the same
term $t \in \terms(p)$.
Now this implies $t \mm u_{s_1} = t \mm u_{s_2}$ for at least two of
the polynomials and, as $\m$ is left cancellative, $u_{s_1} =u_{s_2}$
follows contradicting the fact that the corresponding polynomials are
supposed to be different.
Hence we get $q \mrm w \red{}{\myr}{p}{S_p} 0$ and likewise a
 contradiction.
}
Therefore, we can assume that for our counter-example $\hterm(q)w$
 is $T$-reducible.
Hence, let $\hterm(q)w \id t_1t_2w_1w_2$ such that
 $\hterm(q)\id t_1t_2, t_2 \neq \lambda, w \id w_1w_2$ and
 $l \id t_2w_1$ for some rule $(l,r) \in T$.
Since $q \in \bigcup_{i \geq 0} S_i$ the polynomial $q$ must have been
 added to the set $H$ at some step
 and as we use a fair strategy 
 to remove elements from $H$, $q$ and $C(\hterm(q))$ are considered.
Thus, we can conclude $w_1 \in C(\hterm(q))$ by the definition of this
set.
Now we have to take a closer look at what happens to $q \mrm w_1$.
In case  $q \mrm w_1$ is added to the respective set $S_j$, then 
  $\hterm(q)w \id \hterm(q)w_1w_2 \succ \hterm(q \mrm w_1)w_2$ and $q
  \mrm w_1 \in \bigcup_{i \geq 0} S_i$ imply  
  $q \mrm w = (q \mrm w_1) \mrm w_2  \red{}{\myr}{p}{S_p} 0$
  contradicting our assumption.
Otherwise we have  
 $q \mrm w_1 \red{}{\myr}{p}{S_j} 0$ for the set $S_j$ actual when
 considering $q \mrm w_1$ in the {\bf for all} loop.
But then $q \mrm w_1 = \alpha \skm q' \mrm z$ for some $\alpha \in
 \myk^*$, $q' \in S_j$, $z \in \m$ and $\hterm(q)w_1  \succ \hterm(q \mrm
 w_1) \id \hterm(q')z$.
Moreover, $\hterm(q)w \id \hterm(q)w_1w_2 \succ \hterm(q')zw_2 \succeq
\hterm(q')(z \mm w_2)$ and  $q \mrm w = (q \mrm w_1) \mrm w_2 =
(\alpha \skm q' \mrm z) \mrm w_2 = \alpha \skm q' \mrm (z \mm w_2)$.
Thus $q' \mrm (z \mm w_2) \red{}{\myr}{p}{S_p} 0$ implies $q \mrm w
\red{}{\myr}{p}{S_p} 0$ contradicting our assumption.
\\
\qed
Before we move on to show how the property of being prefix saturated
can be used to characterize prefix Gr\"obner bases, we prove that this
property is decidable for finite sets of polynomials.
\begin{lemma}~\\
{\sl 
It is decidable, whether a finite subset $F$ of $\myk[\m]$ 
 is prefix saturated.
\lemend
}
\end{lemma}
\Ba{}~\\
We can slightly modify the procedure {\sc Prefix Saturation} to give
us a decision procedure, whether a finite set of polynomials is prefix
saturated.


\renewcommand{\baselinestretch}{1}\small\normalsize
\procedure{Prefix Saturated Check\protect{\label{is.prefix.saturated}}}%
{\vspace{-4mm}\begin{tabbing}
XXXXXXX\=XXXX \kill
\removelastskip
{\bf Given:} \>  A finite set $F \subseteq \myk[\m]$ and 
                $(\Sigma,T)$ a convergent presentation of $\m$. \\
{\bf Answer:} \> {\em yes}, \= if $F$ is prefix saturated, \\
        \> {\em no}, \>otherwise. 
\end{tabbing}
\vspace{-7mm}
\begin{tabbing}
XX\=XX\=XX\=XXXX\=XXX\= \kill
answer $:=$ yes \\
{\bf for all}  $q \in F$ {\bf do} \\
\> $t := \hterm(q)$; \\
\> {\bf for all} $w \in C(t) = \{ w \in \Sigma^* \mid tw \id t_{1}t_{2}w \id t_{1}l, t_2 \neq \lambda$
                            for some $(l, r) \in T \} $ {\bf do} \\
\> {\rm\kommentar \% $C(t)$ contains words that will lead to cancellation
  when right multiplied to $t$}\\
\> \> $q'$ := $q \mrm w$ \\
\> \> {\bf if} \>$q' \neq 0$ and  $q' \nred{}{\myr}{p}{F} 0$ \\
\> \>          \>{\bf then} \> answer := no \\
\> \> {\bf endif}\\
\> {\bf endfor} \\
{\bf endfor} 
\end{tabbing}\vspace{-2mm}}
\renewcommand{\baselinestretch}{1.1}\small\normalsize\\
It remains to show that the answer of our procedure is ``no'' if and only if
 $F$ is not prefix saturated.
Obviously, the answer ``no'' implies the existence of an element $w \in \m$ such that for some 
 $f \in F$, $f \mrm w \nred{}{\myr}{p}{F} 0$.
On the other hand, let us assume that our procedure gives us ``yes'', but $F$ is not prefix saturated.
Then there exist $w \in \m$ and  $f \in F$ such that $\hterm(f)w$
 is minimal according to the ordering $\succeq_T$ on $\freemonoid$, $f
 \mrm w \neq 0$
 and $f \mrm w \nred{}{\myr}{p}{F} 0$.
In case $w \in C(\hterm(f))$ the procedure would give us ``no''
contradicting our assumption.
Thus suppose $w \not\in C(\hterm(f))$.
$\hterm(f)w$ must be $T$-reducible as otherwise $f \mrm w \red{}{\myr}{p}{F} 0$.
Let $\hterm(f)w\id\hterm(f)w_1w_2$ such that $w_1 \in C(\hterm(f))$.
Now $f \mrm w_1$ is considered by our procedure and since the answer
 given is ``yes'', we either get $f \mrm w_1 = 0$ contradicting that
 $f \mrm w \neq 0$, or $f \mrm w_1 \red{}{\myr}{p}{F} 0$.
This gives us the existence of a polynomial $f' \in F$ such that
$\hterm(f)w_1 \succ \hterm(f \mrm w_1) \id  \hterm(f')z$
 for some $z \in \m$.
Now $\hterm(f)w \id \hterm(f)w_1w_2 \succ \hterm(f')zw_2 \succeq
\hterm(f')(z \mm w_2)$  gives us $f \mrm w = f' \mrm (z \mm w_2) \red{}{\myr}{p}{F} 0$
contradicting our assumption.
\\
Further this procedure terminates, as the sets $C(t)$ are
always finite.
\\
The complexity  in the number of
monoid multiplications can be described as follows:
Let $n = \max \{ |l| \mid (l,r) \in T \}$, $m = |T|$, and $k = \max \{
|\terms(f)| \mid f \in F \}$.
Then the first {\bf for all} loop is executed $|F|$-times.
A computation of a set $C(t)$ can be bounded by $n \cdot m$ 
and the second {\bf for all} loop is then executes
$|C(t)|$-times.
Within this loop there is one multiplication of a polynomial by a term
involving at most $k$ monoid multiplications and a test whether the
result is zero or prefix reducible to zero in one step.
The latter involves at most $|F|$  reduction
tests, i.e., again $|F|$-times multiplying a polynomial with a term
and checking for equality.
Hence a bound in monoid multiplications is
\vspace{-2mm}
$$\underbrace{|F|}_{{\rm 1st\; loop}} \cdot \underbrace{\phantom{|}n \cdot
  m\phantom{|}}_{{\rm 2nd\; loop}}
\cdot\: (\underbrace{\phantom{|}k\phantom{|}}_{q \mrm w} + \underbrace{|F| \cdot
  k}_{{\rm reducibility\;\; check}}).$$
\vspace{-1cm}
\qed
Prefix saturation enriches a polynomial $p$  to a set
 $S \in \SAT_p(p)$ such that we can
 substitute $q \red{}{\myr}{(s,r)}{p}q'$ by
 $q \red{}{\myr}{p}{p' \in S}q'$.
We use this additional information to give a finite confluence
 criterion that will use a refined
 definition of s-polynomials. 
\begin{definition}\label{def.cpp}~\\
{\rm
Given two non-zero polynomials $p_{1}, p_{2} \in \myk[\m]$, 
 such that $\hterm(p_{1}) \id \hterm(p_{2})w$ for some $w \in \m$ 
  the \index{prefix!s-polynomial}\index{s-polynomial!prefix}\betonen{prefix
    s-polynomial} is defined as
 $$ \spol{p}(p_{1}, p_{2})=\hc(p_1)^{-1} \skm p_1 -\hc(p_2)^{-1} \skm p_2 \mrm w.$$
\dend
}
\end{definition}
As before non-zero  prefix s-polynomials are called non-trivial and
for non-trivial s-polynomials  we
have $\hterm(\spol{p}(p_{1}, p_{2})) \pred \hterm(p_{1}) \id
\hterm(p_{2})w$.
Notice that a finite set $F \subseteq \myk[\m]$ defines
 finitely many prefix s-polynomials.
As before, these s-polynomials alone are not sufficient to
 characterize prefix Gr\"obner
 bases, but lemma  \ref{lem.redp} enables us to localize our confluence test 
 in case we demand our set of polynomials to be prefix saturated.
\begin{theorem}\label{theo.pcp}~\\
{\sl
For a prefix saturated set $F$ of polynomials in $\myk[\m]$, the
 following statements are equivalent:
\begin{enumerate}
\item For all polynomials $g \in \ideal{r}{}(F)$ we have $g
  \red{*}{\myr}{p}{F} 0$. 
\item For all polynomials $f_{k}, f_{l} \in F$ we have 
  $\spol{p}(f_{k}, f_{l}) \red{*}{\myr}{p}{F} 0$.
\end{enumerate}
}
\end{theorem}
\Ba{}~\\
\mbox{$1 \R 2:$ }
Let $\hterm(f_{k}) \id \hterm(f_{l})w$ for
 $w \in \m$.
Then by definition \ref{def.cpp} we get
 $$\spol{p}(f_{k}, f_{l}) =
 \hc(f_k)^{-1} \skm f_{k} -\hc(f_l)^{-1} \skm f_{l} \mrm w \:\in \ideal{r}{}(F),$$
 and hence $\spol{p}(f_{k}, f_{l}) \red{*}{\myr}{p}{F} 0$.

\mbox{$2 \R 1:$ }
We have to show that every non-zero element 
 $g \in \ideal{r}{}(F)$ is $\red{}{\myr}{p}{F}$-reducible
 to zero.
Remember that for
 $h \in \ideal{r}{}(F)$, $ h \red{}{\myr}{p}{F} h'$ implies $h' \in \ideal{r}{}(F)$.
Hence as  $\red{}{\myr}{p}{F}$ is Noetherian
 it suffices to show that every  
 $g \in \ideal{r}{}(F)\backslash\{ 0 \}$ is $\red{}{\myr}{p}{F}$-reducible.
Now, let $g = \sum_{j=1}^m \alpha_{j} \skm f_{j} \mrm w_{j}$ be a
      representation of a non-zero polynomial $g$ such that
      $\alpha_{j} \in \myk^*, f_j \in F, w_{j} \in \m$.
By lemma \ref{lem.prop1.1} we can assume 
 $\hterm(f_{i} \mrm w_{i}) \id \hterm(f_{i})w_{i}$.
This will enable a restriction to prefix s-polynomials in order to
modify the representation of $g$.
Depending on the above  representation of $g$ and a 
 well-founded total ordering $\succeq$ on $\m$ we define
 $t = \max \{ \hterm(f_{j}) \mm w_{j} \mid j \in \{ 1, \ldots m \}  \}$ and
 $K$ is the number of polynomials $f_j \mrm w_j$ containing $t$ as a
 term.
Then $t \succeq \hterm(g)$ and
in case $\hterm(g) = t$ this immediately implies that $g$ is
$\red{}{\myr}{p}{F}$-reducible. 
So
by lemma \ref{lem.pprop} it is sufficient to  show that
$g$ has a prefix standard representation, as this implies that $g$ is
top-reducible using $F$.
This will be done by induction on $(t,K)$, where
 $(t',K')<(t,K)$ if and only if $t' \prec t$
 or $(t'=t$ and $K'<K)$\footnote{Note
  that this ordering is well-founded since $\succ$ is and $K \in\n$.}.
In case $t \succ \hterm(g)$
 there are two polynomials $f_k,f_l$ in the corresponding
 representation\footnote{Not necessarily $f_l \neq f_k$.}
 such that  $\hterm(f_k)w_k \id \hterm(f_l)w_l$.
We have either $\hterm(f_k)z\id\hterm(f_l)$ or $\hterm(f_k)\id\hterm(f_l)z$ for some $z \in\m$.
Without loss of generality let us assume
 $\hterm(f_k)\id\hterm(f_l)z$ and hence $w_l \id zw_k$.
Then definition \ref{def.cpp} provides us with a prefix s-polynomial
 $\spol{p}(f_k,f_l) = \hc(f_k)^{-1} \skm  f_k -
  \hc(f_l)^{-1} \skm f_l \mrm z$.
Note that, while in the proofs of theorem \ref{theo.pcs} and
\ref{theo.rcp} the s-polynomials correspond directly to the overlap
$\hterm(f_k \mrm w_k) = \hterm(f_l \mrm w_l)$, i.e., $w_k$ and $w_l$
are involved in the s-polynomial, now we have an s-polynomial
corresponding directly to the two polynomials $f_k$ and $f_l$.
We will see later on that this localization is strong enough because
 this situation has a prefix of the term $t$ as an upper border and
 lemma \ref{lem.redp} can be applied.
We will now change our representation of $g$ by using the additional
information on the above prefix s-polynomial in such a way that for the new
representation of $g$ we either have a smaller maximal term or the occurrences of $t$
are decreased by at least 1.
Let us assume  $\spol{p}(f_k,f_l) \neq 0$\footnote{In case  $\spol{p}(f_k,f_l) = 0$,
 just substitute $0$ for $\sum_{i=1}^n \delta_i \skm h_i \mrm v_i$
 in the equations below.}.
Hence,  the reduction sequence $\spol{p}(f_k,f_l) \red{*}{\myr}{p}{F} 0
$ results in  a prefix standard representation of the form
 $\spol{p}(f_k,f_l) =\sum_{i=1}^n \delta_i \skm h_i \mrm v_i$,
 where $\delta_i \in \myk^*$, $h_i \in F$, $v_i \in \m$ and
 all terms occurring in the sum are bounded by $\hterm(\spol{p}(f_k,f_l))$.
Now as $\hterm(\spol{p}(f_k,f_l)) \prec \hterm(f_k) \preceq t \id \hterm(f_k)w_k$, 
  by lemma \ref{lem.redp} we then can conclude that
 $t$ is a proper bound for all terms occurring
 in the sum $\sum_{i=1}^n \delta_i \skm h_i \mrm v_i \mrm w_k$.
Without loss of generality we can assume that for all polynomials
 occurring in this representation we have 
$\hterm(h_i \mrm v_i \mrm w_k) \id \hterm(h_i)(v_i \mm w_k)$ as $F$ is prefix
 saturated and in case 
     $\hterm(h_i \mrm v_i \mrm w_k) \neq \hterm(h_i)(v_i \mm w_k)$  we can substitute
     the polynomial $h_i \mrm v_i \mrm w_k$ by a product
     $\tilde{\alpha_i} \skm \tilde{h_i}\mrm u_i$ such that 
     $h_i \mrm v_i \mrm w_k = \tilde{\alpha_i} \skm\tilde{h_i}\mrm u_i$
     and $\hterm(h_i \mrm v_i \mrm w_k) \id \hterm(\tilde{h_i})u_i$
     without increasing neither $t$ nor $K$.
This gives us: 
\begin{eqnarray}
&  & \alpha_{k} \skm f_{k} \mrm w_{k} + \alpha_{l} \skm f_{l} \mrm w_{l}  \nonumber\\ &  &  \nonumber\\ 
& = &  \alpha_{k} \skm f_{k} \mrm w_{k} + 
       \underbrace{ \alpha'_{l} \skm \beta_k \skm f_{k} \mrm w_{k}
                   - \alpha'_{l} \skm \beta_k \skm f_{k} \mrm w_{k}}_{=\, 0} 
       + \alpha'_{l}\skm \beta_l  \skm f_{l} \mrm w_{l} \nonumber\\ 
&  &  \nonumber\\ 
& = & (\alpha_{k} + \alpha'_{l} \skm \beta_k) \skm f_{k} \mrm w_{k} - \alpha'_{l} \skm
       \underbrace{(\beta_k \skm f_{k} \mrm w_{k}
       -  \beta_l \skm f_{l} \mrm w_{l})}_{=\, \spol{p}(f_k,f_l) \mrm w_k} \nonumber\\
& = & (\alpha_{k} + \alpha'_{l} \skm \beta_k) \skm f_{k} \mrm w_{k} - \alpha'_{l} \skm
       (\sum_{i=1}^n \delta_{i} \skm h_{i} \mrm v_{i} \mrm w_k) \label{s3}
\end{eqnarray}
where  $\beta_k = \hc(f_k)^{-1}$, $\beta_l = \hc(f_l)^{-1}$ and  $\alpha'_l \skm \beta_l = \alpha_l$.
By substituting (\ref{s3}) in our representation of $g$ 
 either $t$ disappears   or in
 case $t$ remains maximal among the terms occurring in the new
 representation of $g$, $K$ is decreased. 
\\
\qed
Buchberger's characterization of Gr\"obner bases by s-polynomials in
the commutative polynomial ring provided a finite test to decide
whether a finite set of polynomials is a Gr\"obner basis.
Theorem \ref{theo.pcp} only provides such a test for finite prefix
saturated sets.
Since the property of being prefix saturated is also decidable for
finite sets of polynomials, we can hence decide whether a finite set
of polynomials is a prefix saturated Gr\"obner basis.
We will see later on that  prefix Gr\"obner bases need
not be prefix saturated and how they then can be characterized.
Theorem \ref{theo.pcp} gives rise to the following procedure to
compute prefix Gr\"obner bases.

\procedure{Prefix Gr\"obner Bases\protect{\label{prefix.groebner.bases}}}%
{\vspace{-4mm}\begin{tabbing}
XXXXX\=XXXX \kill
\removelastskip
{\bf Given:} \> A finite set of polynomials $F \subseteq \myk[\m]$. \\
{\bf Find:} \> $\gb(F)$, a prefix Gr\"obner basis of $F$. \\
{\bf Using:} \> $\s_p$ a prefix saturating procedure for polynomials.
\end{tabbing}
\vspace{-7mm}
\begin{tabbing}
XX\=XX\=XXXX\=XX\=XXXX\=XXXX\kill 
$G$ := $\bigcup_{f \in F} \s_p(f)$; \\
{\rm\kommentar \% $G$ is prefix saturated}\\
$B$ := $\{ (q_{1}, q_{2}) \mid q_{1}, q_{2} \in G, q_{1} \neq q_{2} \}$; \\
{\bf while} $B \neq \emptyset$ {\bf do} \\
\>{\rm\kommentar \% Test if
  statement 2 of theorem \ref{theo.pcp} is valid}\\
\>  $(q_{1}, q_{2})$ := {\rm remove}$(B)$; \\
\>  {\rm\kommentar \% Remove an element using a fair strategy}\\
\>      {\bf if}  \>  $\spol{p}(q_{1}, q_{2})$ exists \\
\>                \> {\rm\kommentar \% The s-polynomial is not trivial}\\
\>      \>{\bf then}\>   $h$:= {\rm normalform}$(\spol{p}(q_{1},
                              q_{2}),\red{}{\myr}{p}{G})$; \\
\>      \>          \> {\rm\kommentar \% Compute a normal form using
        prefix reduction} \\
\> \>   \>      {\bf if} \> $h \neq 0$ \\
\> \>  \>\>      {\bf then} \> $G$ := $G \cup \s_p(h)$; \\
\> \>  \>\>                 \> {\rm\kommentar \% $G$ is prefix
  saturated} \\
\> \>  \>\>                 \> $B$ := $B \cup \{ (f, {\tilde h}), ( {\tilde h},f) \mid f \in G, {\tilde h} \in  \s_p(h) \}$; \\
\> \>   \>      {\bf endif}\\
\>      {\bf endif}\\
{\bf endwhile}  \\
$\gb(F):= G$
\end{tabbing}}

There are two crucial points, why  procedure {\sc Prefix Gr\"obner Bases} might not terminate:
 prefix saturation of a polynomial need not terminate 
 and the set $B$ need not become empty.

Note that in case prefix saturation does not terminate it is possible
 to  modify this procedure in order to enumerate a (prefix) Gr\"obner
 basis by using fair enumerations of the prefix saturating sets
 needed.
\auskommentieren{
\procedure{Enumerating Prefix Gr\"obner Bases\protect{\label{enumerating.prefix.groebner.bases}}}%
{\begin{tabbing}
XXXXX\=XXXX \kill
\removelastskip
{\bf Given:} \> A finite set of polynomials $F = \{ f_1, \ldots , f_k \} \subseteq \myk[\m]$.\label{prefix.groebner.bases} \\
{\bf Find:} \> $G$, a recursively enumerable prefix Gr\"obner basis of $\ideal{\myr}{}(F)$. \\
{\bf Using:} \> $\s_p$ enumerates a prefix saturating set for a
polynomial without repetitions.
\end{tabbing}
\begin{tabbing}
XX\=XX\=XX\=XX\=XXXX\=XXXX\=XXXX\kill
$G$ := $\{ f_1, \ldots , f_k \}$; \\
$B$ := $\{ (q_{1}, q_{2}) \mid q_{1}, q_{2} \in G, q_{1} \neq q_{2}
\}$; \\
{\bf for} $1 \leq i \leq k$ {\bf do} \\
start a saturating process  $\s_p(f_i)$; \\
\hspace{1cm}{\rm\scriptsize \% We will assume that $h_{f_i,j}$ is the
  $j$-th polynomial of this process in case it exits, i.e.,}\\
\hspace{1cm}{\rm\scriptsize \%  the
  saturation process has not yet terminated.} \\
{\bf endfor} \\
$j$ := $1$; \\
$l$ := $k$; \\
$s$ := $k+1$; \\
{\bf while} not all saturation processes have terminated {\bf do} \\
\> $S_{j}$ := $\{ h_{f_i,j} | 1 \leq i \leq l \}$; \\
\> \hspace{1cm}{\rm\scriptsize \% Collect the $j$-th elements of the
  still running saturation processes.} \\
\> $G$ := $G \cup S_{j}$; \\
\> $B$ := $B \cup \{ (g,h), (h,g) , (h,f)| g \in G , h,f \in S_{j}
\}$; \\
\>{\bf while} $B \neq \emptyset$ {\bf do} \\
\>\>  $(q_{1}, q_{2})$ := {\rm remove}$(B)$;
    \hspace{1cm}{\rm\scriptsize \% using a fair strategy}\\
\>\>      {\bf if}  $\spol{p}(q_{1}, q_{2})$ exists \\
\>\> \>       {\bf then} $f_s$ :=
            {\rm normalform}$(\spol{p}(q_{1},
            q_{2}),\red{}{\myr}{p}{G})$; \\
\>\> \>      \hspace{1cm}{\rm\scriptsize \% compute a normal form using
        prefix reduction} \\
\>\> \> \>        {\bf if} $f_s \neq 0$ \\
\>\> \> \>  \>        {\bf then} \> $s$ := $s+1$; \\
\>\> \> \>  \>   \> start a saturating process $\s_p(f_s)$; \\
\>\> \> \>  \>   \> $B$ := $B \cup \{ (g, f_s), ( f_s, g) \mid g \in G \}$; \\
\> {\bf endwhile} \\
{\bf endwhile}
\end{tabbing}}}

The sets characterized in theorem \ref{theo.pcp} are prefix Gr\"obner
bases and hence right Gr\"obner bases, but they are required to be
prefix saturated.
Reviewing example \ref{exa.aa.bb.ba=ab} we see that there exist right
 Gr\"obner bases which are not prefix saturated.
\begin{example}~\\
{\rm
Let $\Sigma = \{ a,b \}$ and $T = \{ a^2 \myr \lambda, b^2 \myr \lambda ,
 ba \myr ab \}$ be a presentation of a monoid $\m$ (which is in fact a
 group)  with  a  length-lexicographical ordering induced by $b \succ a$.
\\
Then the set $\{ ab + \lambda \} \subseteq \q[\m]$ itself is
  a right Gr\"obner basis,
  but is neither prefix saturated nor a prefix Gr\"obner basis, as 
we have $b + a \in \ideal{r}{}(ab + \lambda)$
  but $b+a \nred{}{\myr}{p}{ab + \lambda} 0$.
\mbox{\phantom{X}}\exaend
}
\end{example}
Note that even a prefix Gr\"obner basis need not be prefix saturated.
\begin{example}\label{exa.ac=d.bc=e}~\\
{\rm
Let $\Sigma = \{ a,b,c,d, e \}$ and $T = \{ ac \myr d, bc \myr e \}$ be a
 presentation of a monoid $\m$ with a length-lexicographical ordering
 induced by $a \succ b \succ c \succ d \succ
 e$.
\\
Then the set $F = \{ a+b,d+ \lambda , e - \lambda \}$ is
 a prefix Gr\"obner basis in ${\bf Q}[\m]$.
This can be seen by studying the right ideal generated by $F$,
 $\ideal{r}{}(F) = \{ \alpha_1 \skm (a + b)\mrm w_1 + \alpha_2 \skm (d
 - \lambda)\mrm w_2
    + \alpha_3 \skm (e - \lambda)\mrm w_3 + \alpha_4 \skm (d + e)\mrm w_4 | \alpha_i \in \q,
    w_i \in \m, w_1 \neq cw' \}$.
But $F$ is not prefix saturated,
 as $(a+b) \mrm c = d+ e \nred{}{\myr}{p}{F} 0$.
We only have $d+e \red{2}{\myr}{p}{F} 0$.
\mbox{\phantom{XX}}\exaend
}
\end{example}
In the previous section we have seen that we can have a finite strong
Gr\"obner basis although no finite right Gr\"obner basis exists
(compare example \ref{exa.strong.but.no.right}).
Similarly  a finite right Gr\"obner basis can exist while there is no
finite prefix Gr\"obner basis.
\begin{example}~\\
{\rm
Let $\Sigma = \{ a, b \}$ and $T = \{ ba \myr ab \}$ be a presentation
of a commutative monoid $\m$.
\\
Then the set $F = \{ b + \lambda \}$ is a right Gr\"obner basis but no
finite prefix Gr\"obner basis exists.
\exaend
}
\end{example}
Next we will give a characterization of prefix Gr\"obner bases
 without demanding that the set of polynomials is prefix saturated.
This is important as interreducing a set of prefix saturated
 polynomials destroys this property, but in case the  set being interreduced
 is a prefix Gr\"obner basis the resulting set will again be a
 prefix Gr\"obner basis (compare theorem \ref{theo.monic.reduced.psb}). 
Remember that this is not true for right Gr\"obner bases in general
(compare example \ref{exa.interreduction}).
\begin{theorem}\label{theo.altpcp}~\\
{\sl
For a set  $F$ of polynomials in $\myk[\m]$,
 equivalent are:
\begin{enumerate}
\item Every polynomials $g \in \ideal{r}{}(F)$ has a prefix standard representation.
\item \begin{enumerate}
      \item For all polynomials $f \in F$ and all elements  $w \in\m$,
             the polynomial $f \mrm w$ has a  prefix standard representation.
      \item For all polynomials $f_{k}, f_{l} \in F$ the non-trivial
        prefix s-polynomials  have prefix standard representations.
      \end{enumerate}
\end{enumerate}
}
\end{theorem}
\Ba{}~\\
\mbox{$1 \R 2:$ }
This follows immediately.

\mbox{$2 \R 1:$ }
We have to show that every non-zero element 
 $g \in \ideal{r}{}(F)$ has a  prefix standard representation.
Let $g = \sum_{j=1}^m \alpha_{j} \skm f_{j} \mrm w_{j}$  be an
 arbitrary   representation of a non-zero polynomial $g$ such that 
  $\alpha_{j} \in \myk^*$, $f_j \in F$,  $w_{j} \in \m$.
By our assumption and lemma \ref{lem.psr.prop2}
 we can assume  that $\hterm(f_{i} \mrm w_{i}) \id
 \hterm(f_{i})w_{i}$ as $f_{i} \in F$ and every $f_{i} \mrm
 w_{i}$ has a standard prefix representation.
Note that these
 prefix standard representations do not yield a prefix standard
 representation for the polynomial $g$,
 as $\hterm(g) \prec \hterm(f_i)w_i$ is possible.
Using statement (a) and (b) we then can proceed straightforward as in
theorem \ref{theo.pcp} to show that such a representation can be
transformed into a prefix standard representation for $g$ with respect
to $F$.
\auskommentieren{
\\
Depending on this  representation of $g$ and a well-founded 
 total ordering $\succeq$ on $\m$ we define
 $t = \max \{ \hterm(f_{j}) \mm w_{j} \mid j \in \{ 1, \ldots m \}  \}$ and
 $K$ is the number of polynomials $f_j \mrm w_j$ containing $t$ as a term.
\\
Then $t \succeq \hterm(g)$ and
in case $\hterm(g) = t$ this immediately implies that our
representation is in fact a prefix standard representation. 
\\
So it is sufficient to  show that
$g$ such a representation with $\hterm(g) = t$.
This will be done by induction on $(t,K)$, where
 $(t',K')<(t,K)$ if and only if $t' \prec t$
 or $(t'=t$ and $K'<K)$\footnote{Note
      that this ordering is well-founded since $\succ$ is and $K \in\n$.}.
\\
In case $t \succ \hterm(g)$
 there are two polynomials $f_k$
 and $f_l$ in the corresponding
 representation\footnote{Not necessarily $f_l \neq f_k$.}
 such that  $\hterm(f_k)w_k \id \hterm(f_l)w_l$.
We have either $\hterm(f_k)z\id\hterm(f_l)$ or $\hterm(f_k)\id\hterm(f_l)z$ for some
 $z \in\m$.
Without loss of generality let us assume $\hterm(f_k)\id\hterm(f_l)z$
 and hence $w_l \id zw_k$.
Then  definition \ref{def.cpp} provides us with a prefix s-polynomial
 $\spol{p}(f_k,f_l) = \hc(f_k)^{-1} \skm  f_k -
  \hc(f_l)^{-1} \skm f_l \mrm z$.
\\
We will now change our representation of $g$ by using the additional
information on this s-polynomial in such a way that for the new
representation of $g$ we either have a smaller maximal term or the occurrences of the term $t$
are decreased by at least 1.
\\
Let us assume  $\spol{p}(f_k,f_l) \neq 0$\footnote{In case 
    $\spol{p}(f_k,f_l) = 0$,
    just substitute $0$ for the sum $\sum_{i=1}^n d_i \skm h_i \mrm v_i$
    in the equations below.}.
\\
Hence, we have a prefix standard representation of the form
 $\spol{p}(f_k,f_l) =\sum_{i=1}^n \delta_i \skm h_i \mrm v_i$,
 where $\delta_i \in \myk^*,h_i \in F,v_i \in \m$ and
 all terms occurring in the sum are bounded by $\hterm(\spol{p}(f_k,f_l))$.
By lemma \ref{lem.redp} we  can conclude that $t$ is a proper bound for all terms occurring
 in the sum $\sum_{i=1}^n \alpha_i \skm h_i \mrm v_i \mrm w_k$.
Further without loss of generality we can assume the
 new representation to have the required form, i.e., assume that
 $\hterm(h_i \mrm v_i \mrm w_k) \id\hterm(h_i)(v_i \mm w_k)$, as again
 the existence of prefix standard representations for
 right multiples of elements in $F$
 allows us to substitute all polynomials $h_i$, where
 $\hterm(h_i \mrm v_i \mrm w_k) \neq \hterm(h_i)(v_i \mm w_k)$  without increasing
 $t$ or $K$.

This gives us: 
\begin{eqnarray}
&  & \alpha_{k} \skm f_{k} \mrm w_{k} + \alpha_{l} \skm f_{l} \mrm w_{l}  \nonumber\\&  &  \nonumber \\
& = &  \alpha_{k} \skm f_{k} \mrm w_{k} +
        \underbrace{ \alpha'_{l} \skm \beta_k \skm f_{k} \mrm w_{k}
                   - \alpha'_{l} \skm \beta_k \skm f_{k} \mrm w_{k}}_{=\, 0} 
       + \alpha'_{l}\skm \beta_l  \skm f_{l} \mrm w_{l} \nonumber\\ 
&  &  \nonumber \\
& = & (\alpha_{k} + \alpha'_{l} \skm \beta_k) \skm f_{k} \mrm w_{k} - \alpha'_{l} \skm
      \underbrace{(\beta_k \skm f_{k} \mrm w_{k}
       -  \beta_l \skm f_{l} \mrm w_{l})}_{=\, \spol{p}(f_k,f_l) \mrm w_k} \nonumber\\
& = & (\alpha_{k} + \alpha'_{l} \skm \beta_k) \skm f_{k} \mrm w_{k} - \alpha'_{l} \skm
       (\sum_{i=1}^n \delta_{i} \skm h_{i} \mrm v_{i} \mrm w_k) \label{s4}
\end{eqnarray}
where  $\beta_k = \hc(f_k)^{-1}$, $\beta_l = \hc(f_l)^{-1}$ and
  $\alpha'_l \skm \beta_l = \alpha_l$.
By substituting (\ref{s4}) in our representation of $g$
 either $t$ disappears  or in
 case $t$ remains maximal among the terms occurring in the new
 representation of $g$, $K$ is decreased.
}
\\
\qed
Note that this theorem gives us a stronger characterization of prefix
 Gr\"obner bases in so far  as it does not require the sets to be prefix saturated.
We can further use it to prove the following lemma.
\begin{lemma}~\\
{\sl 
It is decidable, whether a finite subset $F$ of $\myk[\m]$ 
 is a prefix Gr\"obner basis.
\lemend
}
\end{lemma}
\Ba{}~\\
The following procedure decides, whether a finite set of polynomials
is a prefix Gr\"obner basis.

\renewcommand{\baselinestretch}{1}\small\normalsize
\procedure{Prefix Gr\"obner Basis Check\protect{\label{is.prefix.groebner.basis}}}%
{\vspace{-4mm}\begin{tabbing}
XXXXXXX\=XXXX \kill
\removelastskip
{\bf Given:} \>  A finite set $F \subseteq \myk[\m]$ and 
                $(\Sigma,T)$ a convergent presentation of $\m$. \\
{\bf Answer:} \> {\em yes}, \= if $F$ is a prefix Gr\"obner basis, \\
        \> {\em no}, \>otherwise. 
\end{tabbing}
\vspace{-7mm}
\begin{tabbing}
XX\=XX\=XXXX\=XX\=XXXX\= \kill
answer := yes; \\
$B$ := $\{ (f,g) \mid f,g \in F, f \neq g \}$; \\
{\bf for all}  $(f,g) \in B$ {\bf do} \\
\> {\bf if} \>$\spol{p}(f,g)$ exists \\
\> \> {\bf then} \>$q$ := {\rm normalform}$(\spol{p}(f,g),\red{}{\myr}{p}{F})$; \\
\> \>            \>  {\bf if} \> $q \neq 0$ \\
\> \>            \>           \> {\bf then} \>  answer := no; \\
\> \>            \>  {\bf endif}\\
\> {\bf endif}\\
{\bf endfor}\\
XX\=XX\=XX\=XXXX\=XXX\= \kill
{\bf for all}  $q \in F$ {\bf do} \\
\> $t := \hterm(q)$; \\
\> {\bf for all} $w \in C(t) = \{ w \in \Sigma^* \mid tw \id t_{1}t_{2}w \id t_{1}l, t_2 \neq \lambda$
                            for some $(l, r) \in T \} $ {\bf do} \\
\> {\rm\kommentar \% $C(t)$ contains words that will lead to cancellation
  when right multiplied to $t$}\\
\> \> $q'$ := {\rm normalform}$(q \mrm w,\red{}{\myr}{p}{F})$; \\
\> \> {\bf if} \>$q' \neq 0$ \\
\> \>          \>{\bf then} \> answer := no; \\
\> \> {\bf endif}\\
\> {\bf endfor} \\
{\bf endfor} 
\end{tabbing}}
\renewcommand{\baselinestretch}{1.1}\small\normalsize\\
It remains to show that the answer of our procedure is ``no'' if and only if
 $F$ is no prefix Gr\"obner basis.
Obviously, the answer ``no'' implies the existence of a polynomial in
 the right ideal generated by $F$ which is not prefix reducible to zero
 using $F$, i.e., $F$ is no prefix Gr\"obner basis.
On the other hand, let us assume that our procedure gives us ``yes'',
 although $F$ is no prefix Gr\"obner basis.
We then know that all prefix s-polynomials originating from polynomials in
 $F$ prefix reduce to zero as the answer is not set to ``no''.
Hence, by theorem \ref{theo.altpcp} there must exist $f \in F$ and $w
 \in \m$ such that $f \mrm w$ has no prefix standard representation
 with respect to $F$ as otherwise $F$ would be a prefix Gr\"obner
 basis.
Let us assume that $\hterm(f)w$
 is minimal according to the ordering $\succeq_T$ on $\freemonoid$
 such that the multiple $f \mrm w \neq 0$ has no prefix standard representation
 with respect to $F$.
Notice that  $w \in C(\hterm(f))$ is not possible as it would
 contradict that the answer given by the procedure is supposed to be ``yes''.
Furthermore, $\hterm(f)w$ must be $T$-reducible as otherwise we would
 get a contradiction by $f \mrm w \red{}{\myr}{p}{f} 0$.
Since $w$ is not $T$-reducible there exist $w_1, w_2 \in \m$ such that
 $\hterm(f)w \id \hterm(f)w_1w_2$ and $w_1 \in C(\hterm(f))$.
As the answer is ``yes'' and $w_1$ has been checked, we can conclude
 $f \mrm w_1 \red{*}{\myr}{p}{F} 0$ implying the existence of a prefix
 standard representation for $f \mrm w_1$, say $f \mrm w_1 =
 \sum_{i=1}^k \alpha_i \skm f_i \mrm u_i$ for some $\alpha_i \in \myk$,
 $f_i \in F$ and $u_i \in \m$.
Moreover we know that $\hterm(f)w_1 \succ \hterm(f \mrm w_1) \succeq
 \hterm(f_i \mrm u_i) \id \hterm(f_i)u_i$ for all $1 \leq i \leq k$.
This immediately provides $\hterm(f)w \id \hterm(f)w_1w_2 \succ
 \hterm(f_i)u_iw_2 \succeq \hterm(f_i)(u_i \mm w_2)$.
Therefore, by our assumption that $\hterm(f)w$ was chosen minimal, we
 can assume that  every multiple $f_i \mrm (u_i \mm w_2)$ has a prefix
 standard representation with respect to $F$, say $f_i \mrm (u_i \mm
 w_2) = \sum_{j=1}^{k_i} \beta_{i_j} \skm f_{i_j} \mrm u_{i_j}$ with
 $\beta_{i_j} \in \myk$, $f_{i_j} \in F$ and $u_{i_j} \in \m$.
Hence, we
 find that $f \mrm w = (\sum_{i=1}^k \alpha_i \skm f_i \mrm u_i) \mrm
 w_2 = \sum_{i=1}^k \alpha_i \skm f_i \mrm (u_i \mm w_2) = \sum_{i=1}^k
 \alpha_i \skm ( \sum_{j=1}^{k_i} \beta_{i_j} \skm f_{i_j} \mrm
 u_{i_j})$ is a prefix standard representation of $f \mrm w$ with
 respect to $F$,
 contradicting our assumption that $f \mrm w$ has none.
\\
\qed
The next remark illuminates the differences between the two characterizations
given for prefix Gr\"obner bases so far.
\begin{remark}~\\
{\rm
Let $F$ be a set of polynomials in $\myk[\m]$.
\begin{enumerate}
\item If $F$ is prefix saturated, then for every polynomial
       $f \in F$ and every element $w \in \m$ the polynomial $f \mrm
       w$ has a prefix standard representation.
       \\
      This follows immediately, since $f \mrm w \red{}{\myr}{p}{F} 0$
       implies that there exists a polynomial
       $f' \in F$ such that $f \mrm w \red{}{\myr}{p}{f'} 0$ and 
       $\hterm(f \mrm w) \id \hterm(f')u$ for some $u \in \m$.
      Note that $F$ need not be a prefix standard basis.  
\item On the other hand, if for every polynomial $f \in F$ and every
       element $w \in \m$ the polynomial $f \mrm w$ has a prefix
       standard representation, this need not imply that $F$ is prefix
       saturated. 
      To see this, let us review example \ref{exa.ac=d.bc=e}.
      \\
      Let $\Sigma = \{ a,b,c,d, e \}$ and $T = \{ ac \myr d, bc \myr e \}$
       be a presentation of a monoid $\m$ with a length
       lexicographical ordering induced by
       $a \succ b \succ c \succ d \succ e$.
      \\
      Then  for every polynomials $f$ in the set
       $F = \{ a+b,d+ \lambda , e - \lambda \}$ and every element
       $w \in \m$ we can show that the multiple $f \mrm w$ has a  prefix
       standard representation.
      For the multiples $(a+b) \mrm w = aw + bw$, $(d+ \lambda) \mrm w = dw
       + w$ and $(e - \lambda) \mrm w = ew -w$, these are  prefix
       standard representations.
      It remains to check the case $(a+b) \mrm cw = dw + ew$.
      Since $d+e \red{}{\myr}{p}{d + \lambda}
             e - \lambda \red{}{\myr}{p}{e - \lambda} 0$,
       we have a  prefix standard representation 
       $dw + ew = (d+ \lambda) \mrm w + (e - \lambda ) \mrm w$,
       but as seen before, $dw + ew$ does not prefix reduce to zero using $F$
       in {\em one}\/ step.
\remend
\end{enumerate}
}
\end{remark}
\begin{definition}~\\
{\rm
We call a set of polynomials $F \subseteq \myk[\m]$ \index{weakly prefix
  saturated}\index{prefix saturated!weakly}\betonen{weakly prefix
  saturated}, if for all $f \in F$ and all $\alpha \in \myk^*$, $w \in
 \m$, we have $\alpha \skm f
 \mrm w \red{*}{\myr}{p}{F} 0$.
\dend
}
\end{definition}
If a set of polynomials $F$ is weakly prefix saturated this implies
 that for all $f \in F$ and all $w \in \m$ the polynomial $f \mrm w$,
 in case it is non-zero, has a prefix standard representation.
Thus we can give the following procedure to compute reduced 
 prefix Gr\"obner bases.
Notice that in this procedure no prefix s-polynomials are computed.
This is due to the fact that we are computing a basis of the right
ideal generated by the input set such that no head terms of the
polynomials in the basis are prefix reducible by the other polynomials
in the set and hence no prefix s-polynomials exist.

\procedure{Reduced Prefix Gr\"obner Bases\protect{\label{reduced.prefix.groebner.bases}}}
{\vspace{-4mm}\begin{tabbing}
XXXXX\=XXXX \kill
\removelastskip
{\bf Given:} \> A finite set $F \subseteq \myk[\m]$. \\
{\bf Find:} \> $\gb(F)$, a  (prefix) Gr\"obner basis of $F$. \\
{\bf Using:} \>   $\s_p$ a prefix saturating procedure for polynomials.
\end{tabbing}
\vspace{-7mm}
\begin{tabbing}
XX\=XX\= XXXX \= XXX \=\kill
$G_0$ := $\emptyset$; \\
$S_0$ := $F$; \\
$i$ := $0$; \\
{\bf while} $S_i \neq \emptyset$ do\\
\> $i$ := $i+1$; \\
\> $q_i$ := {\rm remove}$(S_{i-1})$; \\
\>{\rm\kommentar \%   Remove an element using a fair strategy}\\
\> $q_i'$ := ${\rm normalform}(q_i, \red{}{\myr}{p}{G_{i-1}})$; \\
\>        {\rm\kommentar \% Compute a normal form using
           prefix reduction} \\
\>{\bf if} \>$q_i' \neq 0$ \\
\>         \>{\bf then}  \>$H_i$ := $\{ g \in G_{i-1} |
                              \hterm(g) \mbox{ is prefix reducible using }
                              q_i' \}$;\\
\>         \>            \> {\rm \kommentar \% These polynomials would
  have new head terms after prefix reduction} \\
\>         \>            \> {\rm \kommentar \%  using $q_i'$} \\
\>         \>            \>$G_i$ := {\rm reduce}$((G_{i-1} \backslash H_i) \cup \{ q_i'
\})$; \\
\> \>  \>{\rm\kommentar
  \% reduce$(F) = \{ {\rm normalform}(f, \red{}{\myr}{p}{F \backslash \{
    f \}}) | f \in F \}$\footnotemark}\\ 
\> \>  \>{\rm\kommentar \% No head term of a polynomial in $G_i$ is
  prefix reducible by the other }\\
\> \>  \>{\rm\kommentar \% polynomials in $G_i$} \\
\> \>  \>$S_i$ := $S_{i-1}  \cup H_i \cup 
              \bigcup_{g \in (G_i \backslash G_{i-1})}(\s_p(g) \backslash \{ g \})$; \\
\> \> {\bf else}  \> $G_i$ := $G_{i-1}$; \\
\> \>             \> $S_i$ := $S_{i-1}$; \\   
\>{\bf endif}\\
{\bf endwhile} \\
$\gb (F)$:= $G_i$
\end{tabbing}}
\footnotetext{Notice that only the reducts of the polynomials are touched in this procedure.}

Notice that adding a polynomial to a set $G_i$ we also add the
saturating polynomials to $S_i$.
This will ensure that the final set is weakly prefix saturated.
But in doing so, the sets $S_i$ will in general contain  many
unnecessary polynomials, as in removing a polynomial from a set $G_i$
one could also think of changing the set $S_i$.
We will later on realize this idea for the special case of free group
rings.
In order to show that the procedure actually constructs a prefix
reduced Gr\"obner basis we first prove some lemmata.
Let us start by showing that the sets $G_i \cup S_i$ constructed generate the
same right ideal as $F$.
\vspace{-1mm}
\begin{lemma}\label{lem.idealproperty}~\\
{\sl
Let $F$ be a set of polynomials in $\myk[\m]$ and $G_i$, $S_i$, $i \in
\n$ the respective sets in procedure {\sc Reduced Prefix Gr\"obner
  Bases}.
Then  we have 
$$\ideal{r}{}(F) = \ideal{r}{}(G_i \cup S_i).$$
\lemend}
\end{lemma}
\vspace{-1mm}
\Ba{}~\\
For $i = 0$ we have $F = G_0 \cup S_0$ and hence
 $\ideal{r}{}(F) = \ideal{r}{}(G_0 \cup S_0)$.
On the other hand, for $i>0$ let $G_{i-1}$, $S_{i-1}$ be the respective sets before
 entering the {\bf while} loop for its $i$-th iteration.
Further let $q_i$ be the polynomial chosen from $S_{i-1}$ and 
 $q_i'$ a prefix normal from of $q_i$ with respect to $G_{i-1}$.
\\
Then in case $q_i' =0$ we know $q_i \in \ideal{r}{}(G_{i-1})$ and thus
 as $G_i = G_{i-1}$ and $S_i = S_{i-1} \backslash \{ q_i \}$ we can
 conclude
\vspace{-2mm}
 $$\ideal{r}{}(G_i \cup S_i) = \ideal{r}{}(G_{i-1} \cup S_{i-1})=\ideal{r}{}(F).$$
In case $q_i' \neq 0$, then $\ideal{r}{}({\rm reduce}((G_{i-1} \backslash H_i)
\cup \{ q_i' \})) = \ideal{r}{}((G_{i-1} \backslash H_i)
\cup \{ q_i' \})$,
 $\ideal{r}{}(G_{i-1} \cup \{ q_i' \} \cup (S_{i-1} \backslash \{ q_i \}) ) 
 = \ideal{r}{}(G_{i-1} \cup S_{i-1})$,
$\ideal{r}{}(\bigcup_{g \in (G_i \backslash G_{i-1})}(\s_p(g) \backslash \{
g \})) \subseteq \ideal{r}{}(G_i)$ yield
\vspace{-4mm}
\begin{eqnarray}
  &\phantom{=}& \ideal{r}{}(G_i \cup S_i) \nonumber\\ 
  &=& \ideal{r}{}(G_i \cup (S_{i-1} \backslash \{ q_i \}) \cup H_i \cup 
            \bigcup_{g \in (G_i \backslash G_{i-1})}(\s_p(g)
            \backslash \{ g \})) \nonumber\\
  &=& \ideal{r}{}(G_i \cup (S_{i-1} \backslash \{ q_i \}) \cup H_i)  \nonumber\\
  &=& \ideal{r}{}({\rm reduce}((G_{i-1} \backslash H_i)
      \cup \{ q_i' \}) \cup (S_{i-1} \backslash \{ q_i \}) \cup H_i) \nonumber\\   &=& \ideal{r}{}((G_{i-1} \backslash H_i)
      \cup \{ q_i' \} \cup (S_{i-1} \backslash \{ q_i \}) \cup H_i) \nonumber\\   &=& \ideal{r}{}(G_{i-1} 
      \cup \{ q_i' \} \cup (S_{i-1} \backslash \{ q_i \})) \nonumber\\ 
  &=& \ideal{r}{}(G_{i-1} \cup S_{i-1}) \nonumber\\ 
  &=& \ideal{r}{}(F). \nonumber
\end{eqnarray}
\nolinebreak
\qed
\begin{remark}~\\
{\rm
Taking a close look at the construction of the sets $G_i$, $i \in
\n$, we find that the following observations for the head terms of the
polynomials generated during computation hold:
A set $G_i$, $i>0$, is constructed from a set $G_{i-1}$ by removing a
 polynomial $q_i$ from $S_i$, reducing it to $q_i'$  and in case $q_i'$ 
 is non-zero setting $G_i := {\rm
  reduce}((G_{i-1} \backslash H_i) \cup \{ q_i' \})$.
Notice that $H_i$ contains those polynomials in $G_{i-1}$ that have a head term
prefix reducible by $q_i'$, i.e., these polynomials when reduced with
$q_i'$ will lead to polynomials with different head terms.
Thus in removing these polynomials we find that reducing the
polynomials in $(G_{i-1} \backslash H_i) \cup \{ q_i' \}$ does not
touch the head terms and since $q_i'$ was in prefix normal form with
respect to $G_{i-1}$ it is not changed in this step.
Hence we can conclude $q_i' \in G_i$ and $\hterm(G_i) =
\hterm((G_{i-1} \backslash H_i) \cup \{ q_i' \})$.
Moreover, we know that every term in $\hterm(H_i)$ has $\hterm(q_i')$
as a proper prefix, i.e., all terms in $\hterm(G_{i-1})$ have prefixes
in the set $\hterm(G_i)$.
Therefore, if a polynomial is prefix reducible with respect to some
set $G_k$ it will also be prefix reducible with respect to all sets
$G_{k+n}$, $n \in \n$.
In particular, in case a polynomial $g$ is removed from a set $G_k$ no
polynomial with head term $\hterm(g)w$, $w \in \m$, will be added to a
later set $G_{k+n}$, $n \in \n^+$.
This implies that no cycles occur.
\mbox{\phantom{X}}\remend
}
\end{remark}
The next lemma states that these properties of the sets $G_i$
carry over to the set $G$.
This is obvious in case the procedure terminates.
\begin{lemma}\label{lem.prefix.exists}~\\
{\sl
Let $G$ be the set generated by procedure {\sc
  Reduced Prefix Gr\"obner Bases}.
Then if $f \in G_k$ for some $k \in \n$ there exists a polynomial $g
 \in G$ such that $\hterm(g)$ is a prefix of $\hterm(f)$.
\lemend
}
\end{lemma}
\Ba{}~\\
In case our procedure terminates or $f \in G$ we are done at once.
Hence, let us assume there exists a polynomial $f$ such that $f \in G_k$ for
 some $k \in \n$ but no $g \in G$ exists such that $\hterm(g)$ is a
 prefix of $\hterm(f)$.
Further let $f$ be a counter-example with minimal head term.
As $f \not\in G$ there exists an index $j > k$ such that $f \in
 G_{j-1}$ but $f \not\in G_{j}$.
Let $q_j' \neq 0$ be the polynomial computed in this $j$-th execution
 of the {\bf while} loop.
In case $\hterm(f)$ is prefix reducible by $q_j'$, we know that
 $\hterm(q_j')$ is a proper prefix of $\hterm(f)$ and hence
 $q_j' \in G_{j}$ implies the existence of a polynomial $g \in G$
 such that $\hterm(g)$ is a prefix of $\hterm(q_j')$ and hence of
 $\hterm(f)$, contradicting our assumption.
On the other hand the case that $f$ is replaced in $G_{j-1}$ by a
 polynomial $f'$ after reducing the set $(G_{j-1} \backslash H_j) \cup \{ q_j' \}$
 with $\hterm(f) = \hterm(f')$ cannot occur infinitely often since 
this would imply the existence of an infinite prefix reduction
 sequence $f \red{+}{\myr}{p}{} f_0 \red{+}{\myr}{p}{}
 f_1 \ldots\;$ inducing an infinite
 strictly descending chain $f > f_0 > f_1 \ldots\;$ in $\myk[\m]$
 (compare lemma \ref{lem.Noetherian}).
Hence there exists an index $l \geq j$ such that for a  descendant
 $f'$ of $f$ either no more
 changes occur, i.e., $f' \in G$, or $f'$
 is removed from $G_l$ because its head term is prefix reducible by a
 polynomial $q_{l+1}'$ added to $G_{l+1}$ where $\hterm(q_{l+1}')$ is a proper
 prefix of $\hterm(f')$ and then we can now proceed as above and get
 a contradiction.
\\
\qed
\begin{lemma}\label{lem.psr.exists}~\\
{\sl
Let $G$ be the set generated by procedure {\sc
  Reduced Prefix Gr\"obner Bases}.
Then if $f \in S_k \cup G_k$ for some $k \in \n$, $f$ has a prefix
standard representation with respect to $G$.
\lemend
}
\end{lemma}
\Ba{}~\\
In proving our claim we will distinguish two cases.
Suppose first that $f \in G_k$ but $f$ has no prefix
 standard representation with respect to $G$.
Let us further assume that $f$ is a minimal counter-example.
Since $f \not\in G$  there exists an
 index $j >k$ such that $f \in G_{j-1}$, but $f \not\in G_j$.
Let $q_j' \neq 0$ be the polynomial corresponding to this {\bf while} loop.
In case $f$ is removed since $\hterm(q_j')$ is a proper prefix of $\hterm(f)$
 it is put into $S_j$.
As we are using a fair strategy to remove elements from the
 respective sets $S_i$, there exists an index $l > j$ such that $f$
 is chosen to compute $q_{l}'$.
Then $f \red{+}{\myr}{p}{G_{l-1}} q_{l}'$,$\;q_{l}' < f$ and either
 $q_{l}' = 0$ or $q_{l}' \in G_{l}$. 
In both cases $f$ has a prefix standard representation with respect to
 $G_{l-1} \cup \{ q_l' \}$ and by lemma \ref{lem.psr.prop3} also with respect
 to $G_{l} = {\rm reduce}((G_{l-1} \backslash H_{l}) \cup \{ q_{l}' \})$.
Note that all polynomials involved in this prefix standard representation
 are in $G_{l}$ and are smaller than $f$, i.e., they have prefix standard
 representations with respect to $G$ yielding such a representation for
 $f$ contradicting our assumption.
In case $f$ is replaced by a polynomial $f'$ while computing the set
$G_j = {\rm reduce}((G_{j-1} \backslash H_j) \cup \{ q_j' \})$ we
know that all reductions involved take place at $\reductum(f)$ and
hence all polynomials used for prefix reduction are smaller than $f$.
So are again the polynomials used to prefix reduce these polynomials
and so on.
Then lemma \ref{lem.psr.prop3} gives us the existence of a prefix
standard representation with respect to $G_l$ for $f$ and the existence of
prefix standard representations with respect to $G$ for all
polynomials involved in this representation yield a prefix standard
representation with respect to $G$ for $f$.
On the other hand, suppose  that $f \in S_k$.
Then let $l \geq k$ be the iteration where $f$ is chosen to
 compute $q_l'$, i.e., $f \red{*}{\myr}{p}{G_{l-1}} q'$.
Thus $f$ has a prefix standard representation with respect
 to $G_{l} = {\rm reduce}((G_{l-1} \backslash H_l) \cup \{ q_l' \})$ and 
 since every polynomial in this set has a prefix standard 
 representation with respect to $G$, lemma \ref{lem.psr.prop3}
 yields the existence of a prefix standard representation for $f$
\\
\qed
\begin{theorem}\label{theo.corr.rpgb}~\\
{\sl
Let $G$ be the set generated by procedure
 {\sc Reduced Prefix Gr\"obner Bases} on a finite input $F \subseteq \myk[\m]$.
Then $G$ is a reduced prefix Gr\"obner basis.
\theoend
}
\end{theorem}
\Ba{}~\\
In case procedure {\sc Reduced Prefix Gr\"obner Bases} terminates we
 have $G = G_k$ for some $k \in \n$.
Otherwise, as we assume that elements are removed from the respective
 sets $S_i$ in a fair way, we have $\bigcup_{i \geq 0} \bigcap_{j \geq
  i} S_j = \emptyset$ and $G = \bigcup_{i \geq 0} \bigcap_{j \geq i}
 G_j$.
\\
By construction  no
 prefix s-polynomials exist for the polynomials in $G$.
Hence, in order to show that $G$ is a Gr\"obner basis, by theorem
 \ref{theo.altpcp} it remains to show that $g \in G$ implies that for
 all $w \in \m$ the multiple $g \mrm w$ has a prefix standard
 representation with respect to $G$.
In case  $g \in G$, there exists an index $k \in \n$ such that $g$ is
 added to $G_k$ and hence $\s_p(g) \subseteq G_k \cup S_k$.
Hence we know $g \mrm w \red{}{\myr}{p}{g' \in \sm_p(g)} 0$ and $g' \in S_k \cup
G_k$.
Thus  lemma \ref{lem.psr.exists} yields the existence of a prefix
standard representation for $g'$ with respect to $G$ which can be
extended to a prefix standard representation for $g \mrm w$ by lemma
 \ref{lem.psr.prop3}.
\auskommentieren{Hence, $G$ is a prefix Gr\"obner basis and we have to show that it
generates the same right ideal as $F$.
The inclusion $\ideal{r}{}(G) \subseteq \ideal{r}{}(F)$ follows, as
all polynomials generated by the procedures used lie in the right
ideal generated by $F$.
Now let $f \in \ideal{r}{}(F)$, i.e., $f = \sum_{i=1}^m \alpha_i \skm
f_i \skm w_i$ with $\alpha_i \in \myk^*$, $f_i \in F$ and $w_i \in \m$.
Then for all $1 \leq i \leq m$ we have $f_i \in G_0 \cup S_0$,
implying that every $f_i$ has a prefix standard representation with
respect to $G$, say $f_i = \sum_{j=1}^{m_j} \beta_{i_j} \skm g_{i_j}
\mrm v_{i_j}$ with $\beta_{i_j} \in \myk^*$, $g_{i_j} \in G$ and
$v_{i_j} \in \m$.
Thus we get that $f =  \sum_{i=1}^m \alpha_i \skm
f_i \skm w_i = \sum_{i=1}^m \alpha_i \skm (\sum_{j=1}^{m_j} \beta_{i_j} \skm g_{i_j}
\mrm v_{i_j}) \mrm w_i$, and hence $f \in \ideal{r}{}(G)$.}
\\
\qed
\begin{theorem}~\\
{\sl
Let $F$ be a finite set of polynomials in $\myk[\m]$.
In case $\ideal{r}{}(F)$ has a finite reduced prefix Gr\"obner basis,
 procedure {\sc Reduced Prefix Gr\"obner Basis} terminates.
\theoend
}
\end{theorem}
\Ba{}~\\
Theorem \ref{theo.pgb.psr} implies that
 reduced prefix Gr\"obner bases are
 unique up to multiplication with coefficients.
Hence in case a finite reduced prefix Gr\"obner basis exists, all
 reduced prefix Gr\"obner bases must be finite including the
 one computed by procedure {\sc Reduced Prefix Gr\"obner Basis}.
\auskommentieren{
Let $G'$ be a finite reduced prefix Gr\"obner basis of
$\ideal{r}{}(F)$.
Then, as procedure {\sc Reduced Prefix Gr\"obner Basis} generates a
reduced prefix Gr\"obner basis and for every polynomial ever entered
in a set $G_i$ its head term is prefix reducible with respect to every
following set $G_{i+n}$, $n \in \n$, there must exist an index $k \in
\n$ such that every term in $\hterm(G')$ is prefix reducible with
respect to $G_k$.
Then by theorem \ref{theo.equiv1} we have that $G_k$ is a prefix Gr\"obner
basis since $G_k \subseteq \ideal{r}{}(F)$ and $\ideal{r}{\freemonoid}(\hterm(G')) \cap \m = \hterm(\ideal{r}{\myk[\m]}(F)
  \backslash \{ 0 \})$ and by construction $G_k$ is reduced.
}
\\
\qed
\begin{theorem}\label{theo.term.finite}~\\
{\sl
Let $\m$ be a finite monoid presented by the finite convergent
semi-Thue system.
Further let the prefix saturation procedure $\s_p$ be specified e.g. as
 in procedure {\sc Prefix Saturation}.
Then  the procedure {\sc Reduced Prefix Gr\"obner Bases} terminates.
}
\end{theorem}
\Ba{}~\\
This follows immediately as for a reduced prefix Gr\"obner basis all
head terms are different elements of $\m$ and since $\m$ is finite such
a basis must also be finite.
\\
\qed
\begin{theorem}\label{theo.term.free}~\\
{\sl
Let $\m$ be a free  monoid finitely generated by an alphabet $\Sigma$ 
 and presented by the finite convergent
 semi-Thue system $( \Sigma , \emptyset)$.
Further let the prefix saturation procedure $\s_p$ be specified as
$\s_p(p) = \{ p  \}$ for a polynomial $p \in \myk[\m]$.
Then  the procedure  {\sc Reduced Prefix Gr\"obner Bases} terminates.
}
\end{theorem}
\Ba{}~\\
Let $F \subseteq \myk[\m]$ be the finite input set of polynomials.
Reviewing the definition of prefix reduction we find that all
polynomials added must have head terms of length less equal to $\max \{
|\hterm(f)| \mid f \in F \} = K$.
Hence there exists a prefix Gr\"obner basis whose head terms are all
different and bounded in their length by $K$, i.e., there exists a
finite prefix Gr\"obner basis.
\\
\qed
Note that with the modifications described in theorem
\ref{theo.term.free} procedure {\sc Reduced Prefix Gr\"obner Bases}
directly corresponds to Mora's algorithm for computing Gr\"obner bases
for finitely generated right ideals in non-commutative polynomial
rings as specified in \cite{Mo85} (compare also section \ref{section.mora}).

Chapter \ref{chapter.grouprings} will give some more detailed examples where
finite reduced  prefix Gr\"obner bases can be computed, namely
 in group rings of the class of free groups, the class of plain groups and the
class of context-free groups.
\auskommentieren{
\begin{theorem}~\\
{\sl
Finite Gr\"obner bases for finitely generated right ideals in $\myk[\m]$ exist in
 case
\begin{enumerate}
\item $\m$ is finite,
\item $\m$ is a  finitely generated free monoid.
\end{enumerate}
Using an appropriate presentation the procedure 
 {\sc Reduced Gr\"obner Bases} computes
 such bases.
}
\end{theorem}
\Ba{}~\\
We show the existence of finite Gr\"obner bases for finitely generated right ideals in these structures
 by showing the termination of the completion procedure for them.
\begin{enumerate}
\item First of all, prefix saturation terminates as our monoid
       is finite and at most
       $|\m|$ polynomials have to be considered, namely the set
       $\{ p \mrm w | w \in \m \}$ is finite.
      The procedure stops as soon as $S$ becomes empty.
      Now suppose that our procedure does not terminate.
      Then since $\m$ is finite there must be a term $t \in \m$, which
       occurs infinitely many 
       times among the head terms of the  polynomials $h'$ as computed
       in the above procedure and added to $S$ or $G$. 
      But this is a contradiction to the fact that the $h'$ are in
       head normal form with respect to the already computed polynomials
       in $G$.
\item We can set $\s_p(p) = \{p \}$ and, therefore, saturation terminates.
      All polynomials $q$ added have the property
       $|\hterm(q)| \leq \max \{ |\hterm(f)| \mid f \in F \}$.
      Further the head of any added polynomial is in normal form with
       respect to the already computed set $G$, 
       i.e., the head terms of the added polynomials do not occur
       among the head terms of the polynomials 
       already in $G$.
      Hence, the procedure must terminate.
\\
\qed
\end{enumerate}\renewcommand{\baselinestretch}{1}\small\normalsize
\begin{theorem}~\\
{\sl
Finite Gr\"obner bases for finitely generated right ideals in $\myk[\m]$ exist in
 case
\begin{enumerate}
\item $\m$ is finite
\item $\m$ is a free finitely generated monoid
\end{enumerate}
Using the appropriate presentation the procedure
 {\sc Completion With Respect to Prefix Reduction 1} computes
 such bases.
}
\end{theorem}
\Ba{}
We show the existence of finite Gr\"obner bases for finitely generated right ideals in these structures
 by showing the termination
 of the completion procedure for them.
\begin{enumerate}
\item Since we can set $\s_p(p)= \{ p \}$, we have to take a closer look at 
       the s-polynomials in order to prove termination.
      The procedure stops as soon as all s-polynomials reduce to zero.
      Now suppose that our procedure does not terminate.
      Then since $\m$ is finite there must be a term $t \in \m$, which occurs infinitely many
       times among the head terms of the  polynomials $h'$ as computed in the above procedure and added to $G$.
      But this is a contradiction to the fact that the $h'$ are in normal form.
\item Again we can set $\s_p(p) = \{p \}$ and  take a closer look at 
       the s-polynomials in order to prove termination.
      Looking at the critical overlaps we see that all polynomials $q$ added have 
       $|\hterm(q)| \leq max \{ |\hterm(f)| \mid f \in F \}$.
      Further the head of any added polynomial is in normal form with respect to the already computed set $G$,
       i.e. the head terms of the added polynomials do not occur among the head terms of the polynomials
       already in $G$.
      Therefore, the procedure must terminate.
\\
\qed
\end{enumerate}\renewcommand{\baselinestretch}{1}\small\normalsize}

Let us close this section by considering our monoid ring
 as a quotient\index{quotient} structure of a free monoid ring.
As we have  mainly investigated the right ideal congruence in monoid
 rings, we will modify our quotient in order to introduce right
 reduction.
In the following we will assume that
 $(\Sigma , T)$ is a finite, convergent, and even reduced presentation
 of a monoid $\m$.
As before, $\freemonoid$ is a free monoid generated by $\Sigma$ with 
 concatenation as multiplication.
Let $P_T = \{ l - r \mid (l,r) \in T \}$ be the set of polynomials 
 in $\myk[\freemonoid]$ associated to $T$.
Further let ${\tilde T} = \{ (xl,xr) \mid (l,r) \in T, x \in \freemonoid
 \mbox{ and no proper 
 prefix of $xl$ is $T$-reducible}\}$ denote the \index{prefix-rewriting
 system}\betonen{prefix-rewriting
 system} belonging to $T$ as
 described in \cite{KuMa89} and we can associate the set of polynomials 
 $P_{{\tilde T}} = \{ xl-xr \mid (xl,xr) \in {\tilde T} \}$ to $\tilde{T}$.
\begin{lemma}\label{lem.prefix.quotient}~\\
{\sl
Let $\m$ be a monoid with a finite, convergent, reduced presentation
  $(\Sigma , T)$.
Then the following statements hold:
\begin{enumerate}
\item $\red{}{\myr}{p}{P_{{\tilde T}}}$ is confluent on $\myk[\freemonoid]$.
\item $\ideal{}{\myk[\freemonoid]}(P_{T}) =
        \ideal{r}{\myk[\freemonoid]}(P_{{\tilde T}})$ in $\myk[\freemonoid]$.
\item The quotient
      $\myk[\freemonoid]/\ideal{r}{\myk[\freemonoid]}(P_{{\tilde T}})$
      is a ring.  
\item The monoid ring $\myk[\m]$ is isomorphic to the quotient ring 
        $\myk[\freemonoid]/\ideal{r}{\myk[\freemonoid]}(P_{{\tilde T}})$.
\lemend
\end{enumerate}
}
\end{lemma}
\Ba{}~\\
\vspace{-8mm}
\begin{enumerate}
\item Since the prefix reduction induced by  $P_{{\tilde T}}$ is
       Noetherian, in order to show
       confluence, we only have to take a look at critical pairs 
       of different polynomials
       $xl_i-xr_i,x'l_j-x'r_j \in P_{{\tilde T}}$, caused by a
       superposition $xl_i\id x'l_jz$, where $z \in\Sigma^*$.
      We have to distinguish two cases:
      \begin{enumerate}
      \item Let $|x|<|x'|$, i.e., $x'\id xw, w \in\Sigma^* \backslash
             \{ \lambda \}$.
\\
            Looking at $xl_i\id xwl_jz$ we get $l_i \id wl_jz$ contradicting
             the fact that  $T$ is supposed to be reduced.
      \item Let  $|x|\geq |x'|$, i.e., $x \id x'w, w \in\Sigma^*$.
\\
            Looking at $x'wl_i \id x'l_jz$ we have to consider two cases.
            In case $z \neq \lambda$ we find that $x'l_j$ is a proper
             prefix of $xl_i\id x'wl_i$ and as $x'l_j$
             is $T$-reducible, this is a contradiction to the
             definition
             of  ${\tilde T}$.
            The case $z= \lambda$ likewise gives 
             us a contradiction to $T$ being reduced, as $wl_i \id l_j$.
      \end{enumerate}
      Hence the set  $P_{{\tilde T}}$ leads to no superpositions,
       and therefore is confluent.
\item $\ideal{r}{\myk[\freemonoid]}(P_{{\tilde T}}) \subseteq
       \ideal{}{\myk[\freemonoid]}(P_{T})$ follows immediately.
       To show $\ideal{}{\myk[\freemonoid]}(P_{T}) \subseteq 
        \ideal{r}{\myk[\freemonoid]}(P_{{\tilde T}})$ it suffices
        to prove that
        $x \in \freemonoid$ and $l-r \in P_{T}$ implies
        $xl-xr \in \ideal{r}{}(P_{{\tilde T}})$.
       In case $xl-xr \in P_{{\tilde T}}$ there is nothing to show.
       Suppose  $xl-xr \not\in P_{{\tilde T}}$, i.e., we can
        decompose $xl$ into $xl \id x'l'w$
        for some $x',l',w \in \freemonoid$, where $x'l'$ is the first
        proper $T$-reducible prefix of $xl$.
       Then we get the following reduction sequences, where
        $\red{}{\myr}{l}{T}$ stands for leftmost reduction:
        $$xr \red{}{\longleftarrow}{}{T}
          x\underline{l} \id x'\underline{l'}w
          \red{}{\myr}{l}{T} x'r'w.$$
       Since $T$ is convergent this situation is confluent, even when
        we restrict the reduction strategy to left-most reduction.
       As a consequence there exists an element $z \in IRR(T)$
        such that $xr \red{*}{\myr}{l}{T} z$ and
        $xl \red{}{\myr}{l}{T} x'r'w \red{*}{\myr}{l}{T} z$.
       This means that $xl$ and $xr$ have a common normal form
        with respect to left-most reduction.
       Reviewing the definition of left-most reduction,
        we find that a reduction step
        $ul_1v \myr_T ur_1v$ with  $(l_1,r_1) \in T,
        u,v \in \Sigma^*$ is a left-most reduction step if and only if
        and no proper prefix of $ul_1$ is $T$-reducible.
       This can be compared to subtracting the polynomial 
        $ul_1v-ur_1v$ from $ul_1v$ in  $\myk[\freemonoid]$ and by the
        definition of $P_{\tilde{T}}$ we know that
        $ul_1-ur_1 \in  P_{{\tilde T}}$.
       Therefore $ul_1v \red{}{\myr}{l}{T} ur_1v$ can be
        simulated in $\myk[\freemonoid]$ by prefix reduction as follows:
        $$ul_1v \red{}{\myr}{p}{P_{\tilde{T}}} ul_1v - 
          (ul_1 -ur_1) \mrm v = ur_1v.$$
       Hence, $ul_1v = ur_1v + (ul_1 -ur_1) \mrm v$ and any reduction
        sequence $w \red{*}{\myr}{p}{P_{\tilde{T}}} y$ induces
        $w = y +h$ for some $h \in  \ideal{r}{\myk[\freemonoid]}(P_{{\tilde T}})$.
       As a consequence, $xl \red{*}{\myr}{l}{P_{\tilde{T}}} z$ and
        $xr \red{*}{\myr}{l}{P_{\tilde{T}}} z$ imply that 
        $xl = z +g$ and  $xyr = z +g'$ for some
        $g,g' \in  \ideal{r}{\myk[\freemonoid]}(P_{{\tilde T}})$. 
       Therefore,
        $xl -xr= g - g'\in \ideal{r}{\myk[\freemonoid]}(P_{{\tilde T}})$. 
\item This follows immediately as
       $\ideal{}{\myk[\freemonoid]}(P_{T}) =
        \ideal{r}{\myk[\freemonoid]}(P_{{\tilde T}})$ in
        $\myk[\freemonoid]$ implies that
       $\myk[\freemonoid]/\ideal{r}{\myk[\freemonoid]}(P_{{\tilde T}})$
       and $\myk[\freemonoid]/\ideal{}{\myk[\freemonoid]}(P_{T})$
       are equal as sets and even as rings.
\item Let $\varphi : \myk[\freemonoid] \myr \myk[\m]$ be the natural ring homomorphisms
       defined by setting
      $\varphi(\sum_{i=1}^k \alpha_i \skm w_i) = 
       \sum_{i=1}^k \alpha_i \skm [w_i]_{\m}$ with  
      $\alpha_i \in\myk, w_i \in\Sigma^*$.
      We will show that the kernel of $\varphi$ is a two-sided ideal in
      $\myk[\freemonoid]$, namely $\ideal{}{\myk[\freemonoid]}(P_{T}) =
        \ideal{r}{\myk[\freemonoid]}(P_{{\tilde T}})$, and
      hence $\myk[\m]$ is isomorphic to the ring 
       $\myk[\freemonoid]/\ideal{r}{\myk[\freemonoid]}(P_{{\tilde T}})$.
To see our claim we prove  that for any polynomial
 $g \in \myk[\freemonoid]$ we have $\varphi(g) = 0$ if and only if
 $g \in \ideal{r}{\myk[\freemonoid]}(P_{{\tilde T}})$.
First for $g = \sum_{i=1}^k \alpha_i \skm w_i$ with $\alpha_i \neq 0$
 we show that $\varphi(g)=0$ implies $g \in \ideal{r}{\myk[\freemonoid]}(P_{{\tilde T}})$
 by induction on $k$.
In case $k=0$ we have $g =0$ and are done.
The case $k=1$ is not possible, as $g = \alpha \skm w$ and $\alpha \skm [w]_{\m}=0$
 imply $\alpha = 0$ contradicting our assumption.
In the induction step let us assume 
 $g = \sum_{i=1}^k \alpha_i \skm w_i + \alpha_{k+1} \skm w_{k+1}$.
Then $\varphi(g) = \sum_{i=1}^{k+1} \alpha_i \skm [w_i]_{\m}=0$ gives us
 $\sum_{i=1}^k \alpha_i \skm [w_i]_{\m} = - \alpha_{k+1} \skm [w_{k+1}]_{\m}$ and
 since $\alpha_{k+1} \neq 0$ there exists $1 \leq j \leq k$ such that 
 $[w_j]_{\m} = [w_{k+1}]_{\m}$.
Hence $w_j \red{*}{\lr}{p}{P_{{\tilde T}}} w_{k+1}$ and 
 $w_{k+1} - w_j \in \ideal{r}{\myk[\freemonoid]}(P_{{\tilde T}})$.
Now we can set $\sum_{i=1}^{k+1} \alpha_i \skm [w_i]_{\m} =
 \sum_{i=1}^{k} \beta_i \skm [w_i]_{\m}$ with $\beta_i = \alpha_i$ for
 $1 \leq i \leq k$, $i \neq j$ and $\beta_j = \alpha_j + \alpha_{k+1}$.
Thus, as $\varphi(\sum_{i=1}^{k+1} \alpha_i \skm [w_i]_{\m}) =
 \varphi(\sum_{i=1}^{k} \beta_i \skm [w_i]_{\m}) = 0$, the induction
 hypothesis yields $\sum_{i=1}^{k} \beta_i \skm w_i \in 
 \ideal{r}{\myk[\freemonoid]}(P_{{\tilde T}})$.
Furthermore, $\sum_{i=1}^{k} \beta_i \skm w_i = g - \alpha_{k+1} \skm w_{k+1}
 + \alpha_{k+1} \skm w_j$ and hence, as
 $w_{k+1} - w_j \in \ideal{r}{\myk[\freemonoid]}(P_{{\tilde T}})$,
 $g \in  \ideal{r}{\myk[\freemonoid]}(P_{{\tilde T}})$ holds.
To see on the other hand that $g \in \ideal{r}{\myk[\freemonoid]}(P_{{\tilde T}})$ implies
 $\varphi(g)=0$, let 
 $g = \sum_{i=1}^{n} \gamma_i \skm (x_il_i - x_ir_i) \mrm y_i \in
 \ideal{r}{\myk[\freemonoid]}(P_{{\tilde T}})$.
Remember that $[x_il_iy_i]_{\m} = [x_ir_iy_i]_{\m}$ and hence
 $\varphi(x_il_iy_i - x_ir_iy_i)=0$.
This immediately yields $ \varphi(g) = 0$ and we are done.
\\
\qed
\end{enumerate}\renewcommand{\baselinestretch}{1}\small\normalsize
Thus we have set up a one-to-one correspondence between monoid rings
 and certain quotient rings of free monoid rings.
Now let us state how computation in our quotient structure is done.
For a polynomial  $p =
 \sum_{i=1}^m \alpha_i \skm t_i \in \myk[\freemonoid]$ we define 
 $$[p]_{\ideal{r}{}(P_{\tilde{T}})} =  \sum_{i=1}^m \alpha_i \skm [t_i]_{\m} =  \sum_{i=1}^m \alpha_i \skm (t_i)\nf{P_{\tilde{T}}}$$
 and  we will
 write $p$ instead of $[p]_{\ideal{r}{}(P_{\tilde{T}})}$ to denote elements of 
 $\myk[\freemonoid]/\ideal{}{}(P_{\tilde{T}})$ or $p\nf{P_{\tilde{T}}}$ if we want
 to turn a polynomial in $\myk[\freemonoid]$ into an element of the quotient.
\begin{definition}~\\
{\rm 
Let $p,q \in
\myk[\freemonoid]/\ideal{r}{\myk[\freemonoid]}(P_{\tilde{T}})$.
Then we can define \index{addition!in a quotient}\betonen{addition}\/
 and \index{multiplication!in a quotient}\betonen{multiplication}\/
 as follows:
\begin{enumerate}
\item $p \qadd q = [p +_{\myk[\freemonoid]}
  q]_{\ideal{r}{}(P_{\tilde{T}})}=(p +_{\myk[\freemonoid]} q)\nf{P_{\tilde{T}}}$
\item $p \qmult  q = [p \mrm_{\myk[\freemonoid]} q]_{\ideal{r}{}(P_{\tilde{T}})}=(p \mrm_{\myk[\freemonoid]} q)\nf{P_{\tilde{T}}}$
\end{enumerate}
where $+_{\myk[\freemonoid]}$ and $\mrm_{\myk[\freemonoid]}$ are
 the corresponding ring operations in
 $\myk[\freemonoid]$.
\dend
}
\end{definition}
Then $\myk[\freemonoid]/\ideal{r}{\myk[\freemonoid]}(P_{\tilde{T}})$
 together with $\qadd$ and 
 $\qmult$ is a ring with unit.
We will now introduce quotient prefix reduction to the quotient structure by lifting
 prefix reduction.
\begin{definition}~\\
{\rm
Let $p,q,f$ be some polynomials in $\myk[\freemonoid]/\ideal{r}{\myk[\freemonoid]}(P_{\tilde{T}})$.
Then we set\index{prefix reduction!in a quotient}
$p \red{}{\myr}{p}{f} q$ at a monomial $\alpha \skm t$ of $p$ if there exists a polynomial
 $q' \in \myk[\freemonoid]$ such that
 $p \red{}{\myr}{p}{f} q'$ at $\alpha \skm t$ and $q = q'\nf{P_{\tilde{T}}}$,
 in particular
 $p \red{}{\myr}{p}{f} q'\red{*}{\myr}{p}{P_{\tilde{T}}} q$
 where $q$ is in prefix normal form with respect to $P_{\tilde{T}}$.
We can define $\red{*}{\myr}{p}{}, \red{+}{\myr}{p}{}$,
 $\red{n}{\myr}{p}{}$ as usual.
Reduction by a set $F \subseteq \myk[\freemonoid]/\ideal{r}{\myk[\freemonoid]}(P_{\tilde{T}})$
 is denoted by $p \red{}{\myr}{p}{F} q$ and stands for 
 $p \red{}{\myr}{p}{f} q$ for some $f \in F$,
 also written as  $p \red{}{\myr}{p}{f \in F} q$.
\dend
}
\end{definition}
Then for this reduction we can state:
\begin{lemma}~\\
{\sl
Let $F$ be a set of polynomials and
 $p,q$ some polynomials in $\myk[\freemonoid]/\ideal{}{\myk[\freemonoid]}(P_{\tilde{T}})$.
Then the following statements hold:
\begin{enumerate}
\item $\red{}{\myr}{}{F} \subseteq \red{*}{\myr}{p}{F \cup P_{\tilde{T}}}$.
\item $p \red{}{\myr}{}{F} q$ implies $p > q$.
\item $\red{}{\myr}{}{F}$ is Noetherian.
\end{enumerate}
}
\end{lemma}
\Ba{}~\\
\vspace{-8mm}
\begin{enumerate}
\item This follows by the definition of reduction, since
       $p \red{}{\myr}{p}{f \in F}q$ can be simulated by 
       $p \red{}{\myr}{p}{f}q' \red{*}{\myr}{p}{P_{\tilde{T}}} q$.
\item This follows immediately as prefix reduction already has
       this property.
\item This follows from the fact that prefix reduction is
  Noetherian.
\\
\qed
\end{enumerate}\renewcommand{\baselinestretch}{1}\small\normalsize
But many other properties of prefix reduction on $\myk[\freemonoid]$ are lost.
\begin{lemma}~\\
{\sl
Let $p,q,h$ be some polynomials in
 $\myk[\freemonoid]/\ideal{r}{\myk[\freemonoid]}(P_{\tilde{T}})$,
 $h \neq 0$.
\begin{enumerate}
\item $q < p$ no longer implies $q \qmult  h < p \qmult  h$.
\item $p \red{}{\myr}{p}{p} 0$ no longer implies $p \qmult  h
  \red{}{\myr}{p}{p} 0$.
\end{enumerate}
}
\end{lemma}
\begin{example}\label{exa.fqsdefect}~\\
{\rm
Let $\Sigma = \{ a \}$ and $P_T = \{ a^2 - \lambda \}$.
Then $p=a$, $q=\lambda$ and $h = a$ gives us an appropriate
counter-example, since we have $q<p$ but $q \qmult h = a > p \qmult h
= \lambda$, and $p \qmult h = \lambda$ cannot be quotient prefix reduced by $p$.
\exaend
}
\end{example}
We now can give a definition for Gr\"obner bases in this setting.
\begin{definition}~\\
{\rm
A set $G \subseteq \myk[\freemonoid]/\ideal{r}{\myk[\freemonoid]}(P_{\tilde{T}})$ is said to
 be a \index{Gr\"obner basis!in a quotient}
 \betonen{prefix Gr\"obner basis}\/ with respect to $\red{}{\myr}{p}{}$, if
\begin{enumerate}
\item $\red{*}{\Longleftrightarrow}{p}{G} = \;\;\equiv{\ideal{r}{}(G)}$, and
\item $\red{}{\myr}{p}{G}$ is confluent.
\end{enumerate}
}
\end{definition}
Reviewing example \ref{exa.fqsdefect} we find that 
 a polynomial alone need no longer be  a prefix Gr\"obner basis  of
 the right ideal it generates.
The set $\{ a \}$ is no prefix Gr\"obner basis as we do not have
 $\red{*}{\Longleftrightarrow}{p}{a} = \;\;\equiv{\ideal{r}{}(a)}$.

The following lemma gives a sufficient condition for confluence and we
 will see later on that we can also regain the expressiveness of the
 right ideal congruence.
\begin{lemma}\label{lem.fqs.confluent}~\\
{\sl
Let $P_{\tilde{T}}$ be a prefix Gr\"obner basis in $\myk[\freemonoid]$,
 $F \subseteq \myk[\freemonoid]/\ideal{r}{\myk[\freemonoid]}(P_{\tilde{T}})$
 and let $\red{}{\myr}{p}{F \cup P_{\tilde{T}}}$ be confluent on
 $\myk[\freemonoid]$.
Then $\red{}{\myr}{p}{F}$ is confluent on $\myk[\freemonoid]/\ideal{r}{\myk[\freemonoid]}(P_{\tilde{T}})$.
}
\end{lemma}
\Ba{}~\\
Suppose there exist
 $f,h_1,h_2 \in \myk[\freemonoid]/\ideal{r}{\myk[\freemonoid]}(P_{\tilde{T}})$
 such that $f \red{}{\myr}{p}{F} h_1$ and $f \red{}{\myr}{p}{F} h_2$.
Then we can view these polynomials as elements of 
 $\myk[\freemonoid]$ and simulate $\red{}{\myr}{p}{F}$ by 
 $\red{}{\myr}{p}{F\cup P_{\tilde{T}}}$ giving us
 $f \red{*}{\myr}{p}{F \cup P_{\tilde{T}}} h_1$ and 
 $f \red{*}{\myr}{p}{F \cup P_{\tilde{T}}} h_2$.
Hence, as $\red{}{\myr}{p}{F \cup P_{\tilde{T}}}$ is confluent, there exists
 a polynomial $g \in \myk[\freemonoid]$ such that
 $h_1 \red{*}{\myr}{p}{F \cup P_{\tilde{T}}} g$ and
 $h_2 \red{*}{\myr}{p}{F \cup P_{\tilde{T}}} g$.
Since $\red{}{\myr}{p}{F \cup P_{\tilde{T}}}$ is convergent we
 can use the following reduction strategy:
\begin{enumerate}
\item Do as many prefix reduction steps as possible using
  $P_{\tilde{T}}$.
\item If possible apply one prefix reduction step using $F$ and return
  to 1. 
\end{enumerate}
We stop as soon as no more prefix reduction steps are possible.
Note that we can combine these reduction steps and add some more using $P_{\tilde{T}}$
 to get a sequence as required in the definition of $\red{}{\myr}{p}{F}$.
This gives us that $h_1 \red{*}{\myr}{p}{F} \tilde{g}$ and
 $h_2 \red{*}{\myr}{p}{F} \tilde{g}$, where $\tilde{g} = g\nf{P_{\tilde{T}}}$.
\\
\qed 
Reviewing example \ref{exa.fqsdefect} we see that the converse is not true.
\begin{example}~\\
{\rm
Let $\Sigma = \{ a \}$ and $P_T = \{ a^2 - \lambda \}$.
Then $P_{\tilde{T}} = P_T$.
For $F = \{ a \}$ we get that $\red{}{\myr}{p}{a}$ is  confluent
 on $\myk[\freemonoid]/\ideal{r}{\myk[\freemonoid]}(P_{\tilde{T}})$, 
 but for the set $F \cup P_{\tilde{T}} = \{ a, a^2 - \lambda \}$,
 $\red{}{\myr}{p}{\{ a, a^2 - \lambda\}}$ is not confluent on
 $\myk[\freemonoid]$.
\exaend
}
\end{example}
In order to use the previous lemma to sketch how a Gr\"obner basis
 with respect to $\red{}{\myr}{p}{}$ can be computed, we review the
 definition of prefix s-polynomials for the special case of the free
 monoid ring $\myk[\freemonoid]$.
Given two polynomials $p_{1}, p_{2} \in \myk[\freemonoid]$ such that $\hterm(p_{1}) \id \hterm(p_{2})w$
 for some $w \in \freemonoid$, this  gives us the prefix s-polynomial
 $$ \spol{}(p_{1}, p_{2})=
   \hc(p_1)^{-1} \skm p_1 -  \hc(p_2)^{-1} \skm p_2 \mrm w.$$

Then we can compute a Gr\"obner basis as follows:

Compute a Gr\"obner basis $G'$ of $F \cup P_{\tilde{T}}$ in
 $\myk[\freemonoid]$\footnote{Note that as $P_{\tilde{T}}$ in general
 is infinite we have to use additional information on the structure of
 $P_{\tilde T}$ to specify critical overlaps between polynomials in
 $F$ and $P_{\tilde{T}}$.} with respect to prefix
 reduction without changing the polynomials in $P_{\tilde{T}}$.
Then the set $G := G' \backslash P_{\tilde{T}}$ is a Gr\"obner basis
 of $F$ in $\myk[\freemonoid]/\ideal{r}{\myk[\freemonoid]}(P_{\tilde{T}})$.

Computing a prefix Gr\"obner basis of $F \cup P_{\tilde{T}}$
 in $\myk[\freemonoid]$ three kinds of
 prefix s-polynomials can arise:
\begin{enumerate}
\item $f,g \in P_{\tilde{T}}$:
      Then the corresponding prefix s-polynomial can be omitted as
       $P_{\tilde{T}}$ is already a Gr\"obner basis.
\item $f,g \in F$:
      Then the prefix s-polynomial corresponds to the prefix
       s-polynomial as defined in section \ref{section.prefixreduction}.
\item $f \in F$, $g \in P_{\tilde{T}}$:
      Then the s-polynomial corresponds to the process of
       saturating $f$ as described in procedure {\sc Prefix
         Saturation} on page \pageref{prefix.saturation}.
      In particular this can be compared to overlapping
       the head term of $f$ with the rule $l \myr r \in T$
       where $g = xl-xr$.
      Notice that although $P_{\tilde{T}}$ maybe infinite only finitely many 
       such overlaps can arise.
\end{enumerate}
The correctness of this approach follows from the next theorem.
\begin{theorem}~\\
{\sl 
Let  $P_{\tilde{T}} \subseteq \myk[\freemonoid]$ be a prefix Gr\"obner basis
 and $F \subseteq \myk[\freemonoid]/\ideal{r}{\myk[\freemonoid]}(P_{\tilde{T}})$.
Then the following statements  are equivalent:
\begin{enumerate}
\item $F$ is a prefix Gr\"obner basis in $\myk[\freemonoid]/\ideal{r}{\myk[\freemonoid]}(P_{\tilde{T}})$.
\item
\begin{enumerate}
\item For all $f,g \in F$ we have $\spol{}(f,g) \red{*}{\myr}{p}{F} 0$ and
\item for all $f \in F$, $g \in P_{\tilde{T}}$ we have
       $(\spol{}(f,g))\nf{P_{\tilde{T}}} \red{*}{\myr}{p}{F} 0$.
\theoend
\end{enumerate}
\end{enumerate}
}
\end{theorem}
\Ba{}~\\
\mbox{$1 \R 2:$ }
This follows immediately from the fact that 
 the polynomials
 themselves are elements of
 $\ideal{r}{\myk[\freemonoid]/\ideal{r}{\myk[\freemonoid]}(P_{\tilde{T}})}(F)$ 
 and therefore congruent to zero.
Thus the confluence of $F$ implies that they can be quotient prefix reduced to zero
 using $F$.

\mbox{$2 \R 1:$ }
The items (a) and (b) imply that $F \cup P_{\tilde{T}}$ is a prefix Gr\"obner basis
 in $\myk[\freemonoid]$, i.e., $\red{}{\myr}{p}{F \cup P_{\tilde{T}}}$ is
 confluent.
Hence, by lemma \ref{lem.fqs.confluent}, $\red{}{\myr}{p}{F}$ is also
 confluent.
It remains to show that 
 $\red{*}{\Longleftrightarrow}{p}{F} = \;\;\equiv_{\ideal{r}{}(F)}$.
Obviously, $\red{*}{\Longleftrightarrow}{p}{F} \subseteq
 \;\;\equiv_{\ideal{r}{}(F)}$.
On the other hand, let $p$ and $q$ be polynomials in 
 $\myk[\freemonoid]/\ideal{r}{\myk[\freemonoid]}(P_{\tilde{T}})$.
Then $p \equiv_{\ideal{r}{}(F)} q$ in $\myk[\freemonoid]/\ideal{r}{\myk[\freemonoid]}(P_{\tilde{T}})$
 implies $p \equiv_{\ideal{r}{}(F \cup P_{\tilde{T}})} q$ in
 $\myk[\freemonoid]$.
Further, as $F \cup P_{\tilde{T}}$ is a prefix Gr\"obner basis in
 $\myk[\freemonoid]$, we know 
 $\equiv_{\ideal{r}{}(F \cup P_{\tilde{T}})} = \red{*}{\lr}{p}{F \cup P_{\tilde{T}}}$,
 and as $\red{}{\myr}{p}{F \cup P_{\tilde{T}}}$ is confluent,
 $p \red{*}{\lr}{p}{F \cup P_{\tilde{T}}} q$ implies 
 $p \downarrow_{F \cup P_{\tilde{T}}} q$.
Thus, as in lemma \ref{lem.fqs.confluent}, we can conclude
 $p \Downarrow_{F} q$ giving us $p \red{*}{\Longleftrightarrow}{p}{F} q$.
This completes the proof that $F$ is a prefix Gr\"obner basis in 
 $\myk[\freemonoid]/\ideal{r}{\myk[\freemonoid]}(P_{\tilde{T}})$.
\\
\qed
We now move on to compare prefix reduction defined here for the
 quotient structure to prefix and right reduction in the
 corresponding monoid ring.
For a polynomial $p \in \myk[\freemonoid]/\ideal{r}{\myk[\freemonoid]}(P_{\tilde{T}})$ let
$\tilde{p}$ be the corresponding polynomial in the  monoid ring $\myk[\m]$.
Note that  we can identify the elements of $\m$ with their
normalforms with respect to $P_{\tilde{T}}$ in  $\myk[\freemonoid]/\ideal{r}{\myk[\freemonoid]}(P_{\tilde{T}})$.
\begin{lemma}~\\
{\sl
Let $p,q,f$ be some polynomials in
$\myk[\freemonoid]/\ideal{r}{\myk[\freemonoid]}(P_{\tilde{T}})$ and  let $\tilde{p},\tilde{q},\tilde{f}$ be the corresponding polynomials in $\myk[\m]$.
Then $p \red{}{\myr}{p}{f} q$ if and only if $\tilde{p} \red{}{\myr}{p}{\tilde{f}} \tilde{q}$.
\lemend
}
\end{lemma}
\Ba{}~\\
Before entering the proof of our claim let us first take a closer look at
reduction in $\myk[\freemonoid]/\ideal{r}{\myk[\freemonoid]}(P_{\tilde{T}})$.
If $p \red{}{\myr}{p}{f} q$ at a monomial $\alpha \skm t$ with
 $t \id \hterm(f)u$, then we can express this reduction step by
 $p \red{}{\myr}{p}{f} q' \red{*}{\myr}{p}{P_{\tilde{T}}} q$ and we have
 $q' = p - \alpha \skm \hc(f)^{-1} \skm f \mrm_{\myk[\freemonoid]} u$ 
 and
 $q = p \qadd (- \alpha) \skm \hc(f)^{-1} \skm f
  \qmult u$.
In this context it is easy to see that $p \red{}{\myr}{p}{f} q$
 implies $\tilde{p} \red{}{\myr}{p}{\tilde{f}} \tilde{q}$.
On the other hand, prefix reduction  in the monoid ring requires that the head term
 of the polynomial is a prefix of
 the term to be reduced.
Hence, such a prefix reduction step can be split
 into first doing one step using prefix reduction in the free monoid
 ring and then
 normalizing the new monomials using $P_{\tilde{T}}$.
\\
\qed
\begin{corollary}~\\
{\sl
Let $G \subseteq
 \myk[\freemonoid]/\ideal{r}{\myk[\freemonoid]}(P_{\tilde{T}})$
 and let $\tilde{G}$ be
 the corresponding set in $\myk[\m]$.
Then $G$ is a prefix Gr\"obner basis with respect to $\red{}{\myr}{p}{}$ if and only if
 $\tilde{G}$ is a prefix Gr\"obner basis with respect to $\red{}{\myr}{p}{}$.
\ohnebeweis
}
\end{corollary}
Obviously, then $\red{}{\myr}{p}{}$ must be weaker than $\red{}{\myr}{r}{}$.
\begin{corollary}~\\
{\sl
Let $p,q,f$ be some polynomials in
$\myk[\freemonoid]/\ideal{r}{\myk[\freemonoid]}(P_{\tilde{T}})$ and let $\tilde{p},\tilde{q},\tilde{f}$ be the corresponding polynomials in $\myk[\m]$.
Then $p \red{}{\myr}{p}{f} q$ implies $\tilde{p} \red{}{\myr}{r}{\tilde{f}} \tilde{q}$ but not vice versa. \ohnebeweis
}
\end{corollary}
\begin{example}~\\
{\rm
Let $\Sigma = \{ a \}$ and $T = \{ a^3 \myr \lambda \}$
 be a presentation of a group $\g$ with a
 length-lexicographical ordering on $\g$.
Further let $\tilde{p} = a$ and $\tilde{f} = a^2$ be polynomials in
$\myk[\g]$.
\\
Then $\tilde{p} \red{}{\myr}{r}{\tilde{f}} 0$, but for the corresponding polynomials
 $p= a$ and $f=a^2$ in the quotient
 $\myk[\freemonoid]/\ideal{r}{\myk[\freemonoid]}(P_{\tilde{T}})$ we have
 $p \nred{}{\myr}{p}{f}$ as $a^2$ is no prefix of $a$.
Notice that the set $\{ \tilde{f} \}$ itself is a Gr\"obner basis with respect
 to $\red{}{\myr}{r}{}$ but neither with respect to $\R^p$ nor
 $\red{}{\myr}{p}{}$.
\exaend
}
\end{example}
\auskommentieren{
The following example shall illustrate that there is
 a great difference in computing a Gr\"obner basis in
 a monoid ring where the monoid is commutative, depending on
 whether we take the approach using a quotient of a polynomial ring 
 or a quotient of a free monoid ring.
\begin{example}~\\
{\rm
Let the monoid $\m$ be presented by
 $\Sigma = \{ a,b,c \}$ with $c \succ b \succ a$ and
 $T = \{ ba \myr ab, ca \myr ac, cb \myr bc \}$.
We will investigate the set $F= \{ abc + \lambda \}$.
\\
In $\myk[\m]$ the set $F$ itself is a Gr\"obner basis for
 the ideal generated by $F$.
On the other hand, in $\myk[\freemonoid]/\ideal{r}{}(P_{\tilde{T}})$ no finite
prefix Gr\"obner basis with respect to $\R$ exists,
 but the set $G = \{ ab^ic + b^i | i \in \n \}$ is an infinite
 one.
}
\end{example}}

%% file: commutativereduction.tex
\section{The Concept of Commutative Reduction}
The concept of prefix reduction as introduced in the previous
 section is a very strong restriction as only few classes of monoids
 respectively groups allow finite prefix Gr\"obner bases.
For example,  commutative monoids in general cannot be treated by this
 approach.
Let $\Sigma = \{ a,b \}$ be generators for a free commutative monoid
$\freecomm = \{ a^ib^j \mid i,j \in \n \}$ with a length-lexicographical
ordering induced by $b \succ a$.
Then the right ideal generated by the polynomial  $ab + \lambda$ does
not have a finite prefix Gr\"obner basis.
Therefore, we will introduce another way of weakening right reduction
for commutative monoids which  makes use of the fact that they can be
presented by ordered words and semi-Thue systems modulo commutativity.
We will introduce the same ideas used in prefix reduction by
generalizing the term ``prefix'' in order to refine right reduction.
Remember that for an alphabet $\Sigma = \{ a_1, \ldots, a_n \}$,
 $\freecomm = \{ a_1^{i_1} \ldots a_{n\phantom{1}}^{i_n} \mid i_j \in \n
 \}$ is a free commutative monoid with multiplication $\cm$
(compare definition \ref{def.free.commutative.monoid}).
We can define
a tuple-ordering on $\freecomm$ as follows:
\begin{definition}\label{def.tuple}~\\
{\rm
Let $u \id a_1^{i_1} \ldots a_{n\phantom{1}}^{i_n}, v \id a_1^{j_1}  \ldots a_{n\phantom{1}}^{j_n}$
be two elements of $\freecomm$.
 We define $u \tupeq v$ if for each $1 \leq l \leq n$ we
       have  $i_l \geq j_l$.
 Further we define $u \tupgreater v$ if $u \tupeq v$ and $i_l>j_l$
       for some $1 \leq l \leq n$.
\dend
}
\end{definition} 
Notice that for terms $u,v \in \freecomm$, $u$ is a divisor of $v$ if
 and only if $u \tupleq v$.
$u$ then can  be viewed as a ``commutative prefix'' of $v$ and
 we have similar properties as in the case of prefixes in the free monoid $\Sigma^*$.

Let $\m$ be given by a  semi-Thue systems modulo commutativity
$(\Sigma, T_c)$ which is convergent
with respect to an admissible total ordering $\succeq_{\freecomm}$ on
$\freecomm$ (e.g. a
length-lexicographical ordering)\footnote{Remember that such a finite
  convergent  semi-Thue system modulo commutativity always exists for a finitely
  generated commutative monoid.}.
Then the elements of $\m$ are denoted by irreducible ordered words, i.e.,
 $\m \subseteq \freecomm$, and  $\succeq_{\freecomm}$ is an extension
 of $\tupeq$ and its restriction to the irreducible representatives of
 the monoid elements is
 a total, well-founded ordering  $\succeq$ on $\m$.
An important fact is that  we have $u \cm v \succeq_{\freecomm} u
 \mm v$ for all $u,v \in\m$, but
 the ordering $\succeq$ in general need not be admissible on $\m$.

Since for commutative monoids the right ideals and the ideals in the
corresponding monoid ring coincide, we will study ideals.
Let us start with refining our view on representations of polynomials
for commutative monoid rings as follows.
\begin{definition}\label{def.csr}~\\
{\rm
Let $F$ be a set of polynomials  and $p$ a non-zero
polynomial in $\myk[\m]$.
A representation 
$$ p = \sum_{i=1}^{n} \alpha_i \skm f_{i} \mrm w_i, \;\;
  \mbox{ with } \alpha_i \in \myk^*,
  f_{i} \in F, w_i \in \m $$
  is called a \index{free commutative!standard
    representation}\index{standard representation!free
    commutative}\betonen{free commutative standard representation}\/ 
 in case for all $1 \leq i \leq n$ we have $\hterm(p)   
 \succeq \hterm(f_{i}) \cm w_i$.
A set $F \subseteq \myk[\m]$ is called a 
 \index{free commutative!standard basis}\index{standard basis!free
   commutative}\betonen{free commutative standard basis}\/ if every
 non-zero polynomial in
$\ideal{}{}(F)$ has a free commutative standard representation with
respect to $F$.
\dend
}
\end{definition}
Notice that $\hterm(p) \succeq \hterm(f_i) \cm w_i$ immediately
implies $\hterm(p) \succeq \hterm(f_i) \cm w_i \succeq \hterm(f_i \mrm
w_i)$.
Furthermore, $\hterm(f_i) \cm w_i = \hterm(f_i) \mm w_i$ implies 
 $\hterm(f_i) \cm w_i = \hterm(f_i \mrm w_i)$.
On the other hand, in case $\hterm(p) =
 \hterm(f_{i} \mrm w_i)$ this yields
      $\hterm(f_{i} \mrm w_i)=\hterm(f_{i}) \cm w_i$ and  
  $\hterm(p)  \tupeq \hterm(f_i)$.
This situation must occur for at least one polynomial in the
 representation.
The following lemmata give some more information on special free commutative
 standard representations.
\begin{lemma}\label{lem.csr.prop2}~\\
{\sl
Let $F$ be a set of polynomials in $\myk[\m]$ such that for all $f \in
F$ and all $w \in \m$ the polynomial $f \mrm w$ has a free commutative standard
representation with respect to $F$.
Then the polynomial $f \mrm w$ has a free commutative standard representation 
 $f \mrm w = \sum_{i=1}^{n} \alpha_i \skm f_{i} \mrm w_i$, with $\alpha_i \in \myk^*,
 f_{i} \in F, w_i \in \m$ such that $\hterm(f \mrm w)
 \succeq \hterm(f_{i} \mrm w_i)=\hterm(f_{i}) \cm w_i$.
\lemend
}
\end{lemma}
\Ba{}~\\
We will prove this lemma by contradiction.
Let us assume the claim is not true.
Then there exists a counter-example  $f \mrm w$ such that 
 $\hterm(f \mrm w)$ is minimal among all counter-examples.
By our assumption $f \mrm w$ has a free commutative standard representation,
e.g.\  $f \mrm w = \sum_{i=1}^{m} \alpha_i \skm g_{i} \mrm w_i$,
 such that $\alpha_i \in \myk^*, g_{i} \in F, w_i \in \m$.
Without loss of generality we can assume that for some  $k \leq m$, 
 $g_1, \ldots, g_k$, are the polynomials involved in the head
 term of $f \mrm  w$, i.e., $\hterm(f \mrm w) = \hterm(g_i) \cm w_i$ for
 all $1 \leq i \leq k$.
Hence, we know $k < m$, as otherwise we would get a contradiction to
$f \mrm w$ being a counter-example.
Furthermore, for all $k+1 \leq j \leq m$ we know $\hterm(g_j \mrm w_j)
\prec \hterm(f \mrm w)$ and hence every such polynomial has a free commutative
standard representation of the desired form, say
 $g_j \mrm w_j = \sum_{l=1}^{n_l} \alpha'_{j_l} \skm g_{j_l} \mrm
 w'_{j_l}$,
 with $\alpha'_{j_l} \in \myk^*$, $g_{j_l} \in F$ and $w'_{j_l} \in \m$.
Thus the representation
 $f \mrm w = \sum_{i=1}^k \alpha_i \skm g_i \mrm w_i + 
 \sum_{i=k+1}^{n} \alpha_i \skm
 (\sum_{l=1}^{n_l} \alpha'_{j_l} \skm g_{j_l} \mrm w'_{j_l})$
is a free commutative standard representation of the desired form, contradicting
our assumption.
\\
\qed
\begin{lemma}\label{lem.csr.prop1}~\\
{\sl
Let $F$ be a free commutative standard basis in $\myk[\m]$.
Then every non-zero polynomial $p \in \ideal{}{}(F)$ has a free commutative standard
 representation
 $p = \sum_{i=1}^{n} \alpha_i \skm f_{i} \mrm w_i, \mbox{ with } \alpha_i \in \myk^*,
  f_{i} \in F, w_i \in \m$
such that  for all $1 \leq i \leq n$ we even have $$\hterm(p) \succeq
 \hterm(f_{i}) \cm w_i
 = \hterm(f_{i} \mrm w_i).$$
\lemend
}
\end{lemma}
\Ba{}~\\
Since $p \in \ideal{}{}(F) \backslash \{ 0 \}$, $p$ has a free
commutative standard representation with respect to $F$, say 
$p = \sum_{i=1}^{n} \alpha_i \skm f_{i} \mrm w_i$, with $\alpha_i \in \myk^*$,
 $f_{i} \in F$, $w_i \in \m$.
Further, by lemma \ref{lem.csr.prop2} every multiple $f_i \mrm w_i$
has a free commutative standard representation, say
 $f_i \mrm w_i = \sum_{j=1}^{m} \delta_j \skm g_j \mrm v_j$ with $\delta_j \in
 \myk^*$, $g_j \in F$ and $v_j \in \m$ such that
 $$\hterm(p) \succeq \hterm(f_i \mrm w_i) \succeq \hterm(g_j \mrm
 v_j)= \hterm(g_j) \cm v_j.$$
\auskommentieren{
Let $p$ be a non-zero polynomial in $\ideal{}{}(F)\backslash \{ 0 \}$.
\\
We show our claim by induction on $\hterm(p)$.
\\
In the base case we can assume
 $\hterm(p) = \min \{ \hterm(g) | g \in \ideal{}{}(F)\backslash \{ 0 \}\}$.
\\  
Since $F$ is a free commutative standard basis, we know that $p$ has a free commutative standard representation
  $p = \sum_{i=1}^{n} \alpha_i \skm f_{i} \mrm w_i, \mbox{ with } \alpha_i \in \myk^*,
 f_{i} \in F, w_i \in \m$
 such that  for all $1 \leq i \leq n$ we have $\hterm(p) \succeq
 \hterm(f_{i}) \cm w_i$ and furthermore $\hterm(f_i) \cm w_i
 \succeq \hterm(f_{i} \mrm w_i)$.
Without loss of generality we can assume that $\hterm(p) =
 \hterm(f_1)\cm w_1=\hterm(f_1 \mrm w_1)$, i.e., the head term of $p$ can be
 eliminated by subtracting an appropriate  multiple of $f_1$.
Then looking at the polynomial 
 $h = p - \hc(p) \skm \hc(f_1)^{-1} \skm f_1 \mrm w_1$ 
 we find that $h$ lies in the ideal generated by $F$ and, since
 $\hterm(p)$ is minimal and $\hterm(h) \prec \hterm(p)$, we can conclude that
 $h=0$.
Thus $p$ has a free commutative standard
 representation $p = \hc(p) \skm \hc(f_1)^{-1} \skm f_1 \mrm
 w_1$ of the desired form.
\\
Now suppose $\hterm(p) \succ \min \{\hterm(g) | g \in \ideal{}{}(F)\backslash \{ 0 \}\}$.
Then again resulting from the free commutative standard representation for $p$, 
 there exists a polynomial $f_1 \in F$ such that
 $\hterm(p)=\hterm(f_1) \cm w_1= \hterm(f_1 \mrm w_1)$
 for some $w_1 \in \m$.
Hence as before, looking at the polynomial 
 $h = p - \hc(p) \skm \hc(f_1)^{-1} \skm f_1 \mrm w_1$  
 we know that $h$ lies in the ideal generated by $F$ and since
 $\hterm(h) \prec \hterm(p)$ either $h=0$ giving us 
 $p = \hc(p) \skm \hc(f_1)^{-1} \skm f_1 \mrm w_1$
 or  our induction hypothesis yields the existence
 of a  free commutative standard
 representation of the desired form for $h$, say 
 $h = \sum_{j=1}^{m} \delta_j \skm g_j \mrm v_j$ with $\delta_j \in
 \myk^*$, $g_j \in F$ and $v_j \in \m$.
Thus this gives us the desired special free commutative  standard
 representation of $p$, namely $p =
 \sum_{j=1}^{m} \delta_j \skm g_j \mrm v_j +  \hc(p) \skm \hc(f_1)^{-1} \skm f_1 \mrm w_1$. 
\\}
\qed
In particular this lemma implies that a free commutative standard
basis is a stable standard basis but the converse need not hold as the
next example shows.
\begin{example}\label{exa.fcsb}~\\
{\rm
Let $\Sigma = \{ a,b \}$ and
 $T_c = \{a^2 \myr \lambda, b^2 \myr \lambda \}$ 
 be a presentation of a commutative group $\g$
 with a length-lexicographical ordering induced by $a \succ b$.
Further take the set $F = \{ ab \}$.
\\
Then  all polynomials $g \in \ideal{}{}(F)$ have a stable 
 standard representation in $F$, but
 e.g. the polynomial $b$ has no free commutative standard representation.
\exaend
}
\end{example}
Now we can use the lemma of Dickson to show the existence
 of finite commutative standard bases.
\index{Dickson's Lemma}
\begin{lemma}[Dickson]~\\
{\sl
For every infinite sequence of elements $m_s \in \freecomm$, $s \in \n$,
 there exists
 an index $k \in \n$ such that for every index $i>k$ there exists and index
 $j\leq k$ and an element $w \in \freecomm$ such that 
 $m_i = m_j \cm w$.
\lemend\ohnebeweis
}
\end{lemma}
Free commutative standard representations provide
 us with enough information
 to characterize free commutative standard bases
 (which are closely related to special Gr\"obner
 bases as we will see later on) by their head terms in a way similar to
 the case of usual polynomial rings.
\begin{theorem}\label{theo.cequiv}~\\
{\sl
Let $F$ be a set of polynomials in $\myk[\m]$ and $G \subseteq
\ideal{}{\myk[\m]}(F)\backslash \{ 0 \}$.
Then  the following statements are equivalent:
\begin{enumerate}
\item $G$ is a free commutative standard basis for
 $\ideal{}{\myk[\m]}(F)$\footnote{I.e., $\ideal{}{}(F) =
    \ideal{}{}(G)$ and $G$ is a free commutative standard basis.}.
\item $\ideal{}{\freecomm}(\hterm(G)) \cap \m = \hterm(\ideal{}{\myk[\m]}(F)
  \backslash \{ 0 \})$.
\end{enumerate}
Note that the set $\hterm(\ideal{}{\myk[\m]}(F)\backslash \{ 0 \})$
 in general is no  ideal in $\m$.
\theoend
}
\end{theorem}
\Ba{}~\\
\mbox{$1 \R 2:$ }
  The inclusion $\ideal{}{\freecomm}(\hterm(G)) \cap \m
   \subseteq \hterm(\ideal{}{\myk[\m]}(F)\backslash \{ 0 \})$ follows
   as $\hterm(g) \cm u \in \m$ for some $g \in G$, $u \in \freecomm$ implies
   $\hterm(g \mrm u) = \hterm(g) \cm u$, $u \in \m$ and as the multiple
   $g \mrm u$ belongs to $\ideal{}{\myk[\m]}(F)$.
  It remains to show that $\ideal{}{\freecomm}(\hterm(G)) \cap \m \supseteq
  \hterm(\ideal{}{\myk[\m]}(F)\backslash \{ 0 \})$ holds.
  Let $g \in \ideal{}{\myk[\m]}(F)\backslash \{ 0 \}$.
  Then since $G$ is  a free commutative standard basis  for $\ideal{}{\myk[\m]}(F)$,
  there exists a free commutative standard representation
   $g=\sum_{i=1}^{n} \alpha_i \skm g_{i} \mrm w_i$ with $\alpha_i \in \myk^*, g_i
   \in G$ and $w_i \in \m$ such that 
   $\hterm(g) \succeq \hterm(g_i) \cm w_i$.
  Furthermore there exists $1 \leq k \leq n$ such that 
   $\hterm(g)=\hterm(g_{k}) \cm w_k$, i.e., $\hterm(g) \in
   \ideal{}{\freecomm}(\hterm(G)) \cap \m$.

\mbox{$2 \R 1:$ }
 We have to show that every $g \in \ideal{}{\myk[\m]}(F)\backslash \{ 0
 \}$ has a free commutative standard representation with respect to $G$.
 This will be done by induction on the term $\hterm(g)$.
 In the base case we can assume
   $\hterm(g) = \min \{ w \mid w \in \hterm(\ideal{}{\myk[\m]}(F)\backslash \{ 0 \})\}$.
 Then,
  since  $\ideal{}{\freecomm}(\hterm(G)) \cap \m =
   \hterm(\ideal{}{\myk[\m]}(F)\backslash \{ 0 \})$ there exists a polynomial 
   $f \in G$ such that $\hterm(g)=\hterm(f) \cm w$
   for some $w \in \m$.
  Eliminating the head term of $g$ by subtracting an appropriate 
   multiple of $f$ we get the polynomial 
   $h = g - \hc(g) \skm \hc(f)^{-1} \skm f \mrm w$. 
  As $g \in \ideal{}{\myk[\m]}(F)$,  $h$ also lies in the  ideal
   generated by $F$.
  Moreover, since
   $\hterm(g)$ is minimal and $\hterm(h) \prec \hterm(g)$, we can conclude
   $h=0$ and $g$ has a free commutative standard
   representation $g = \hc(g) \skm \hc(f)^{-1} \skm f \mrm w$.
  Now let us suppose $\hterm(g) \succ \min \{w \mid w \in
  \hterm(\ideal{}{\myk[\m]}(F)\backslash \{ 0 \})\}$.
  Then again there exists a polynomial $f \in G$ such that
   $\hterm(g)=\hterm(f) \cm w$
   for some $w \in \m$.
  Hence looking at the polynomial 
   $h = g - \hc(g) \skm \hc(f)^{-1} \skm f \mrm w$ 
   we know that $h$ lies in the ideal generated by $F$ and since
   $\hterm(h) \prec \hterm(g)$ either $h=0$, giving us that
   $g = \hc(g) \skm \hc(f)^{-1} \skm f \mrm w$,
   or  our induction hypothesis yields the existence
   of a free commutative standard representation for $h$ with respect to $G$,
   say $h = \sum_{j=1}^{m} \beta_j \skm g_{j} \mrm v_j$ where $\beta_j
   \in \myk^*$, $g_j \in G$ and $v_j \in \m$.
  Thus we have a  free commutative standard  representation of the polynomial 
   $g$, namely $g =
   \sum_{j=1}^{m} \beta_j \skm g_{j} \mrm v_j  +  \hc(g) \skm \hc(f)^{-1} \skm f \mrm w$. 
\\
\qed
\begin{theorem}\label{theo.finite.fcsb}~\\
{\sl
Every ideal in $\myk[\m]$ has a finite free
 commutative standard basis.
\theoend
}
\end{theorem}
\Ba{}~\\
Let $\mswab{i}$ be an ideal in $\myk[\m]$.
Then we can view $\hterm(\mswab{i}\backslash \{ 0 \})$ as a subset of $\freecomm$.
Further by  Dickson's lemma every (infinite) set in $\freecomm$ is
 finitely generated, i.e., $\ideal{}{\freecomm}(\hterm(\mswab{i}\backslash \{ 0 \}))$ is
 finitely generated, with respect to the multiplication $\cm$.
In case $\hterm(\mswab{i}\backslash \{ 0 \})$ is finite, we can set $S=\hterm(\mswab{i}\backslash \{ 0 \})$.
Furthermore, for each $t \in \hterm(\mswab{i}\backslash \{ 0 \})$ we take a polynomial $g_t \in \mswab{i}$ such
 that $\hterm(g_t)=t$  and let $G= \{ g_t \mid t \in \hterm(\mswab{i}\backslash \{ 0 \}) \}$.
Otherwise if $\hterm(\mswab{i}\backslash \{ 0 \}) = \{ t_i \mid i \in
\n \}$ is infinite, let $S = \{ s_1, \ldots, s_k \}$ be such a subset
of $\hterm(\mswab{i}\backslash \{ 0 \})$ as described in Dickson's lemma.
Then for each $1 \leq j \leq k$ again take a polynomial
 $g_{s_j} \in \mswab{i}$ such
 that $\hterm(g_{s_j})={s_j}$  and let $G= \{ g_{s_j} | s_j \in S\}$.
Now in both cases we get 
$$\ideal{}{\freecomm}(\hterm(G)) =  \ideal{}{\freecomm}(S) =
 \ideal{}{\freecomm}(\hterm(\mswab{i}\backslash \{ 0 \}))$$ and
 $$\hterm(\mswab{i}\backslash \{ 0 \}) = \ideal{}{\freecomm}(\hterm(\mswab{i}\backslash \{ 0 \})) \cap \m =
 \ideal{}{\freecomm}(\hterm(G)) \cap \m,$$
 i.e., by theorem \ref{theo.cequiv}, since $G \subseteq \mswab{i}$,
 $G$ is a free commutative standard basis of $\mswab{i}$.
\\
\qed
As in the previous section we can characterize free commutative
standard bases by
weakening right reduction.
\begin{definition}\label{def.redc}~\\
{\rm
Let $p, f$ be two non-zero polynomials in $\myk[\m]$. 
We say $f$
\index{commutative!reduction}\index{reduction!commutative}\betonen{commutatively
  reduces} $p$ to $q$
 at a monomial $\alpha \skm t$ of $p$ in one step, denoted by
 $p \red{}{\myr}{c}{f} q$, if
\begin{enumerate}
\item[(a)] $\hterm(f) \cm w = t$ for some $w \in \m$, i.e.,  $t \tupeq \hterm(f)$, and
\item[(b)] $q = p - \alpha \skm \hc(f)^{-1} \skm f \mrm w$.
\end{enumerate}
We write $p \red{}{\myr}{c}{f}$ if there is a polynomial $q$ as defined
above and $p$ is then called  commutatively reducible by $f$. 
Further we can define $\red{*}{\myr}{c}{}, \red{+}{\myr}{c}{}$,
 $\red{n}{\myr}{c}{}$ as usual.
Commutative reduction by a set $F \subseteq \myk[\m]$ is denoted by
 $p \red{}{\myr}{c}{F} q$ and abbreviates $p \red{}{\myr}{c}{f} q$
 for some $f \in F$,
 which is also written as  $p \red{}{\myr}{c}{f \in F} q$.
\dend
}
\end{definition}
Notice that if $f$ commutatively reduces $p$ to $q$ at a
 monomial $\alpha \skm t$
 then $t \not\in \terms(q)$ and $p > q$.
Furthermore, commutative reduction is Noetherian and $p \red{}{\myr}{c}{q_1} 0$ and
 $q_1 \red{}{\myr}{c}{q_2} 0$ imply $p \red{}{\myr}{c}{q_2} 0$ (compare lemma
 \ref{lem.red}).
Moreover, as prefix reduction, commutative reduction is terminating
with respect to arbitrary sets of polynomials.
\vspace{-2mm}
\begin{definition}~\\
{\rm
We call a set of polynomials $F \subseteq \myk[\m]$  
 \index{reduced set}\index{interreduced!(set of polynomials)}\betonen{interreduced} or \betonen{reduced} with respect to
 $\red{}{\myr}{c}{}$, if no polynomial $f$ in
 $F$ is commutatively reducible by the other polynomials in $F \backslash \{
 f \}$.
\dend
}
\end{definition}
\vspace{-2mm}
As in the case of prefix reduction, commutatively reducing a
polynomial by itself results in zero and hence so-defined reduced sets
can be compared to the concept of reduced sets in the usual commutative
polynomial ring.
We have $\red{}{\myr}{c}{} \subseteq \red{}{\myr}{r}{}$ and
similar to prefix reduction,
 we have more information on the reduction step.
\vspace{-2mm}
\begin{remark}\label{rem.ctransitiv}~\\
{\rm
     Let $p \red{}{\myr}{c}{q} $ and $q \red{}{\myr}{c}{q_1} q_2 $.
     Then in case $\hterm(q)=\hterm(q_2)$ we immediately get $p \red{}{\myr}{c}{q_2}$.
     Otherwise $\hterm(q)=\hterm(q_1) \cm y$, for some $y \in\m$ 
      implies $p \red{}{\myr}{c}{q_1}$.
     Hence, we have $p \red{}{\myr}{c}{\{ q_1, q_2 \}}$.
\mbox{\phantom{XX}}\remend
}
\end{remark}
\vspace{-1mm}
This property of commutative reduction corresponds to the fact
that the existence of free commutative standard representations with respect to
a set of polynomials remains true for an interreduced version of the
set (compare lemma \ref{lem.csr.prop3}).
Therefore, we will later on be able to compute reduced
Gr\"obner bases in this setting.
%
\begin{lemma}\label{lem.csr.prop3}~\\
{\sl
Let $F$ and $G$ be two sets of polynomials in $\myk[\m]$ such that
 every polynomial in $F$ has a free commutative standard representation with
 respect to $G$.
Then if a polynomial $p$ has a free commutative standard representation with
 respect to $F$ it also has one with respect to $G$. 
\lemend
}
\end{lemma}
\Ba{}~\\
Let $p = \sum_{i=1}^{n} \alpha_i \skm f_{i} \mrm w_i$,  with  $\alpha_i \in \myk^*,
 f_{i} \in F, w_i \in \m $ be a free commutative standard representation of a
 polynomial $p$ with respect to the set of polynomials $F$, i.e., for
 all $1 \leq i \leq n$ we have $\hterm(p)\succeq \hterm(f_{i}) \cm w_i$.
Furthermore, every polynomial $f_i$ occurring in this sum has a free commutative standard
 representation with respect to the set of polynomials $G$, say
 $f_i = \sum_{j=1}^{n_i} \beta_{i_j} \skm g_{i_j} \mrm v_{i_j}$,
 with  $\beta_{i_j} \in \myk^*,
 g_{i_j} \in G$, and  $v_{i_j} \in \m$ such that for all $1 \leq j \leq n_i$ we
 have $\hterm(f_i) \succeq \hterm(g_{i_j}) \cm v_{i_j}$. 
These representations can be combined in the sum
\vspace{-1mm}
 $$p = \sum_{i=1}^{n} \alpha_i \skm (\sum_{j=1}^{n_i} \beta_{i_j} \skm g_{i_j} \mrm
 v_{i_j}) \mrm w_i.$$
It remains to show that this in fact is a free commutative standard
 representation, i.e., to prove that for all 
 $1 \leq i \leq n$ and all $1 \leq j \leq n_i$,
 we get $\hterm(p) \succeq  \hterm(g_{i_j}) \cm (v_{i_j} \mm w_i)$. 
\\
This now follows immediately as  for all 
 $1 \leq i \leq n$ and all $1 \leq j \leq n_i$ we have
\vspace{-2mm}
 $$\hterm(p) \succeq \hterm(f_i) \cm w_i \succeq (\hterm(g_{i_j}) \cm
 v_{i_j}) \cm w_i \succeq \hterm(g_{i_j}) \cm
 (v_{i_j} \mm w_i).$$
\qed
%
Using this lemma we can show that finite monic reduced free commutative standard bases are
 unique with respect to the presentation of the monoid.
\begin{theorem}~\\
{\sl
Every ideal in $\myk[\m]$ contains
 a unique monic finite reduced free commutative standard  basis.
\auskommentieren{
In particular in case $G$ is a  commutative standard  basis,
 then if there is a polynomial $f$ in $G$ such that $f$ is commutatively
 reducible using the set of polynomials $G \backslash \{ f \}$, i.e.,
 $f \red{}{\myr}{c}{G \backslash \{ f \}} f'$, 
 then the set $G' = (G \backslash \{ f \}) \cup \{ f' \}$ is a free commutative
 standard basis of the same  ideal.
In case $f$ is reduced at its head monomial or $f' = 0$ this also holds
for the set $G' = G \backslash \{ f \}$.}
\theoend
}
\end{theorem}
\Ba{}~\\
Let $G$ be a  finite free commutative standard basis of the ideal
$\mswab{i}$ which must exist by theorem \ref{theo.finite.fcsb}.
Then by theorem \ref{theo.cgb.csr} we
know
 $$\ideal{}{\freecomm}(\hterm(G)) \cap \m = \hterm(\mswab{i}\backslash \{
 0 \}).$$
As the set $\hterm(G)$ is finite\footnote{The sets $\hterm(G)$ and $\hterm(\mswab{i}\backslash \{
 0 \})$ of course depend on the presentation of $\m$ chosen, especially on the ordering induced on $\m$.}, there
exists a  subset $H \subseteq \hterm(G)$ such that
\begin{enumerate}
\item for all $m \in \hterm(G)$ there exists an element $m' \in H$ and an
  element $w \in \freecomm$ such that $m = m' \cm w$, 
\item for all $m \in H$ there exists no element $m' \in H \backslash
  \{ m \}$ such that $m' \tup m$, and
\item $\ideal{}{\freecomm}(H) \cap \m = 
       \ideal{}{\freecomm}(\hterm(G)) \cap \m =
       \hterm(\mswab{i}\backslash \{  0 \}).$
\end{enumerate}
Since for each term $t \in H$ there exists at least one polynomial in
 $G$ with head term $t$ we can choose one of them, say $g_t$, for every
 $t \in H$.
Then as in theorem \ref{theo.finite.fcsb} the set $G' = \{ g_t | t \in H \}$
is a free commutative standard basis.
Further all polynomials in $G'$ have different head terms and no head
term is commutatively reducible by the other polynomials in $G'$.
Hence, if we commutatively interreduce $G'$ giving us another set of polynomials
$G''$, we know $\hterm(G') = \hterm(G'')$ and this set is a free
commutative standard basis as well.
\auskommentieren{
The latter follows as every polynomial $g \in \ideal{}{}(G) = \ideal{}{}(G')$ has a free commutative
standard representation with respect to $G$ and $f$ has a free
commutative standard representation with respect to $G' = G \backslash
\{ f \} \cup \{ f' \}$.
Thus lemma \ref{lem.csr.prop3}  yields that $g$ also has a free commutative
standard representation with respect to $G'$.
}
To see the latter we use the fact that for a free commutative standard basis
$G$, if $f \in G$ and $f \red{}{\myr}{c}{G \backslash \{ f \}} f'$, then
$(G\backslash \{ f \}) \cup \{f' \}$ again is a free commutative
standard basis of the same ideal.
This follows immediately by lemma \ref{lem.csr.prop3} as $f$ has a
free commutative standard representation with respect to $(G\backslash \{ f \}) \cup \{f' \}$.
\\
It remains to show the uniqueness of the  reduced free commutative standard
basis if we restrict ourselves to sets of monic polynomials.
Let us assume $S$ is another monic reduced free commutative standard basis of
 $\mswab{i}$.
Further let $f \in S \bigtriangleup G'' = (S \backslash G'') \cup
 (G'' \backslash S)$ be a polynomial such that $\hterm(f)$ is minimal in
 the set of terms $\hterm(S \bigtriangleup G'')$.
Without loss of generality we can assume that $f \in S \backslash
G''$.
As $G''$ is a free commutative standard basis and $f \in \mswab{i}$ there
exists a polynomial $g \in G''$ such that $\hterm(f) = \hterm(g) \cm w$ for some $w
\in \m$.
We can even state that $g \in G'' \backslash S$ as otherwise $S$ would
not be commutatively interreduced.
Since $f$ was chosen such that $\hterm(f)$ was minimal in
 $\hterm(S \bigtriangleup G'')$, we get $\hterm(f) = \hterm(g)$\footnote{Otherwise
   $\hterm(f) \succ \hterm(g)$ would contradict our assumption.}.
As we assume $f \neq g$ this gives us $f-g \neq 0$, $\hterm(f-g) \prec \hterm(f) = \hterm(g)$ and $\hterm(f-g) \in \terms(f) \cup
 \terms(g)$.
But $f-g \in \mswab{i}$ implies the existence of a polynomial $h \in
S$ such that $\hterm(f-g)=\hterm(h) \cm w'$ for some $w' \in \m$, implying that $f$
is not commutatively reduced.
Hence we get that $S$ is not commutatively interreduced, contradicting our
assumption.
\\
\qed
The following example shows that different presentations for the monoid can result in different 
 reduced free commutative standard bases.
\begin{example}~\\
{\rm
Let $\Sigma = \{ a,b,c \}$ and $T_c = \emptyset$ be a presentation of a 
commutative monoid $\m$ with a length-lexicographical ordering induced by
 $a \succ b \succ c$.
Then the set $F = \{ a+c+\lambda, b+c+\lambda \}$ is a reduced free commutative standard basis.
This is no longer true if we assume that the ordering on $\m$ is induced by
 $c \succ a \succ b$.
Then the set $F' = \{ c+b+\lambda, a-b\}$ is a reduced free commutative standard basis.
\exaend
}
\end{example}
Another setting in which a bound of a representation is preserved under
multiplication is specified in the next lemma.
This observations will be a weaker substitute  for the fact that in a
commutative polynomial ring $p \red{*}{\myr}{b}{F} 0$ implies $\alpha
\skm p \mrm w \red{*}{\myr}{b}{F} 0$.
\begin{lemma}\label{lem.redc}~\\
{\sl
Let $F$ be a set of polynomials and $p$ a polynomial in $\myk[\m]$. 
Further let $p \red{*}{\myr}{c}{F} 0$ and this reduction sequence results
 in a representation
 $p = \sum_{i=1}^{k} \alpha_i \skm g_i \mrm w_i$, where
 $\alpha_i \in \myk^*$, $g_i \in F$, and  $w_i \in \m$.
Then for every term $t \in \m$ such that $t \succ \hterm(p)$ and
every term $w \in \m$ we get that 
       $s \in \bigcup_{i=1}^k \terms(g_i \mrm w_i \mrm w)$
       implies $t \cm w \succ s$.
\lemend
}
\end{lemma}
\Ba{}~\\
As $\sum_{i=1}^{k} \alpha_i \skm g_i \mrm w_i$ belongs to the reduction
 sequence $p \red{*}{\myr}{c}{F} 0$, for all $u \in \bigcup_{i=1}^{k} \terms(g_i \mrm w_i)$
 we have $\hterm(p) \succeq u$ implying $t \cm w \succ \hterm(p) \cm w
 \succ u \cm w \succeq u \mm w$.
\\
Note that this proof uses the fact that the ordering $\succ$ on $\m$
 is induced by the completion ordering $\succeq_T$ of the vector
 replacement system
 $(\Sigma, T)$ presenting  $\m$,
 as we need that the ordering is admissible on $\freecomm$, i.e.,
 $u \cm v \succeq_T  (u \cm v) \nf{T} = u \mm v$ for all $u,v \in \m$.
\\
\qed
Let us continue by taking a closer look at commutative reduction.
An essential property in characterizing Gr\"obner bases, the translation lemma holds.
\begin{lemma} \label{lem.confluentc}~\\
{\sl
Let $F$ be a set of polynomials and $p,q,h$ some
 polynomials in $\myk[\m]$.
\begin{enumerate}
\item
Let $p-q \red{}{\myr}{c}{F} h$.
Then there are  $p',q' \in \myk[\m]$ such that 
 $p  \red{*}{\myr}{c}{F} p', q  \red{*}{\myr}{c}{F} q'$ and $h=p'-q'$.
\item
Let $0$ be a normal form of $p-q$ with respect to $\red{}{\myr}{c}{F}$.
Then there exists a polynomial  $g \in \myk[\m]$ such that
 $p  \red{*}{\myr}{c}{F} g$ and $q  \red{*}{\myr}{c}{F} g$.
\lemend
\end{enumerate}
}
\end{lemma}
\Ba{}
\begin{enumerate}
\item  Let $p-q \red{}{\myr}{c}{F} h = p-q-\alpha \skm f \mrm w$,
        where $\alpha \in \myk^*, f \in F, w \in \m$
        and $\hterm(f) \cm w = t$, i.e. $\alpha \skm \hc(f)$ is the
        coefficient of $t$ in $p-q$.
       We have to distinguish three cases:
       \begin{enumerate}
         \item $t \in \terms(p)$ and $t \in \terms(q)$:
               Then we can eliminate the term $t$ in the polynomials 
                $p$ respectively $q$ by commutative
                reduction and 
                get $p \red{}{\myr}{c}{f} p - \alpha_1 \skm f \mrm w= p'$,
                $q \red{}{\myr}{c}{f} q - \alpha_2 \skm f \mrm w= q'$,
                with $\alpha_1  -  \alpha_2 =\alpha$,
                where $\alpha_1 \skm \hc(f)$ and 
                $\alpha_2 \skm \hc(f)$ are
                the coefficients of $t$ in $p$ respectively $q$.
         \item $t \in \terms(p)$ and $t \not\in \terms(q)$:
               Then  we can eliminate the term $t$ in the polynomial 
                $p$  by commutative reduction and  get
                $p \red{}{\myr}{c}{f} p - \alpha \skm f \mrm w= p'$ 
                and $q = q'$.
         \item $t \in \terms(q)$ and $t \not\in \terms(p)$:
               Then  we can eliminate the term $t$ in the polynomial 
                $q$ by commutative reduction and  get
                $q \red{}{\myr}{c}{f} q + \alpha \skm f \mrm w= q'$
                and $p = p'$.
       \end{enumerate}
      In all cases we have $p' -q' =  p - q - \alpha \skm f \mrm w = h$.
\item We show our claim by induction on $k$, where $p-q \red{k}{\myr}{c}{F} 0$.
      In the base case $k=0$ there is nothing to show.
      Hence, let $p-q \red{}{\myr}{c}{F} h  \red{k}{\myr}{c}{F} 0$.
      Then by (1) there are polynomials $p',q' \in \myk[\m]$ such that 
       $p  \red{*}{\myr}{c}{F} p', q  \red{*}{\myr}{c}{F} q'$ and $h=p'-q'$.
      Now the induction hypothesis for $p'-q' \red{k}{\myr}{c}{F} 0$  yields 
       the existence of a polynomial $g \in \myk[\m]$ such that
       $p  \red{*}{\myr}{c}{F} p' \red{*}{\myr}{c}{F} g$ and
       $q  \red{*}{\myr}{c}{F} q' \red{*}{\myr}{c}{F} g$.
\\
\qed
\end{enumerate}\renewcommand{\baselinestretch}{1}\small\normalsize
The following lemma shows that commutative reduction  captures the  ideal
congruence when using free commutative standard bases
 of ideals.
Reviewing example \ref{exa.fcsb} we find that this is not true in general.
\begin{lemma}\label{lem.csb.cong}~\\
{\sl
Let $F$ be a free commutative standard basis and $p,q,h$ some polynomials in
 $\myk[\m]$.
Then
  $$p \red{*}{\lr}{c}{F} q \mbox{ if and only if } p - q \in
  \ideal{}{}(F).$$
\lemend
}
\end{lemma}
\Ba{}~\\
In order to prove our claim we have to show two subgoals.
The inclusion $\red{*}{\lr}{c}{F}  \subseteq \;\; \equiv_{\ideal{}{}(F)}$
 is an immediate consequence of the definition of commutative reduction and
 can be shown by induction as in lemma \ref{lem.strong.congruence}. 
To prove the converse inclusion 
 $\equiv_{\ideal{}{}(F)} \: \subseteq \red{*}{\lr}{c}{F} $
 let us remember that $p \equiv_{\ideal{}{}(F)} q$ implies
 $p = q + \sum_{j=1}^{m} \alpha_{j} \skm f_j \mrm w_{j}$, where
 $\alpha_{j} \in \myk^*, f_j \in F, w_{j} \in \m$ and every
 multiple $f_j \mrm w_{j}$ belongs to $\ideal{}{}(F)$.
As $F$ is a free commutative standard basis, by lemma \ref{lem.csr.prop1}
 we can assume $\hterm(f \mrm w)=\hterm(f) \mm_{\freecomm} w$ for
 all polynomials occurring in the sum.
Under these assumptions we can then show our claim straightforward as
in lemma \ref{lem.strong.congruence} by induction on $m$.
\auskommentieren{
In the base case $m = 0$ there is nothing to show.
\\
Let  $p = q + \sum_{j=1}^{m} \alpha_{j} \skm f_j \mrm w_{j} +
 \alpha_{m+1} \skm f_{m+1} \mrm w_{m+1}$
 and by our induction hypothesis
 $p \red{*}{\lr}{c}{F} q + \alpha_{m+1} \skm f_{m+1} \mrm w_{m+1}$.
\\
Let $t=\hterm(f_{m+1}) \mm_{\freecomm} w_{m+1}$.
\\
In case $t \not\in \terms(q)$ we get $q +  \alpha_{m+1} \skm f_{m+1}
 \mrm w_{m+1}  \red{}{\myr}{c}{f_{m+1}} q$ and  are done.
\\
In case $t \not\in \terms(p)$ we get $p - \alpha_{m+1} \skm f_{m+1}
  \mrm w_{m+1}  \red{}{\myr}{c}{f_{m+1}} p$.
As $p  - \alpha_{m+1} \skm f_{m+1} \mrm w_{m+1} =
   q + \sum_{j=1}^{m} \alpha_{j} \skm f_j \mrm w_{j}$ the
 induction hypothesis yields 
 $p  - \alpha_{m+1} \skm f_{m+1} \mrm w_{m+1}
 \red{*}{\lr}{c}{F} q$ and hence we are done.
\\ 
Otherwise 
 let $\beta_1 \neq 0$ be the coefficient of $t$ in
 $q +  \alpha_{m+1} \skm f_{m+1} \mrm w_{m+1}$ and
 $\beta_2 \neq 0$ the coefficient of $t$ in $q$.
\\
This gives us a commutative reduction step

\hspace*{1cm}$q +  \alpha_{m+1} \skm f_{m+1} \mrm w_{m+1} \red{}{\myr}{c}{f_{m+1}}$\\
\hspace*{1cm}$q +  \alpha_{m+1} \skm f_{m+1} \mrm w_{m+1}
                  - \beta_1 \skm \hc(f_{m+1})^{-1}
                   \skm  f_{m+1} \mrm w_{m+1} = $ \\
\hspace*{1cm}$q - (\beta_1 \skm \hc(f_{m+1})^{-1} -\alpha_{m+1})
             \skm  f_{m+1} \mrm w_{m+1}$

 eliminating the occurrence of $t$ in $q +  \alpha_{m+1} \skm f_{m+1} \mrm w_{m+1}$.
\\
Then obviously 
  $\beta_2 = (\beta_1 \skm \hc(f_{m+1})^{-1}
                      -\alpha_{m+1}) \skm \hc(f_{m+1})$
  and, therefore, we find $q  \red{}{\myr}{c}{f_{m+1}}
             q - (\beta_1 \skm \hc(f_{m+1})^{-1} -\alpha_{m+1})
             \skm  f_{m+1} \mrm w_{m+1}$, i.e.,
  $q$ and $q +  \alpha_{m+1} \skm f_{m+1} \mrm w_{m+1}$ are joinable.
}
\\
\qed
Let us continue by defining Gr\"obner bases with respect to
 commutative reduction.
\begin{definition}\label{def.cgb}~\\
{\rm
A  set $G \subseteq \myk[\m]$ is called a 
 \index{Gr\"obner basis!commutative}\betonen{Gr\"obner basis}\/ with respect to
 the reduction $\red{}{\myr}{c}{}$ or a
 \index{commutative!Gr\"obner basis}\betonen{commutative Gr\"obner basis}, if
\begin{enumerate}
\item[(i)] $\red{*}{\lr}{c}{G}   = \;\; \equiv_{\ideal{}{}(G)}$, and
\item[(ii)] $\red{}{\myr}{c}{G}$ is confluent.
\dend
\end{enumerate}
}
\end{definition}
As in the previous section there is a natural connection between
 free commutative standard bases and commutative reduction.
\begin{lemma}\label{lem.cprop}~\\
{\sl
Let $F$ be a set of polynomials and $p$ a non-zero polynomial in $\myk[\m]$.
\begin{enumerate}
\item Then $p \red{*}{\myr}{c}{F} 0$ implies the existence of a free 
  commutative standard representation for $p$.
\item In case $p$ has a free commutative  standard representation
  with respect to $F$,
  then $p$ is commutatively reducible at its head monomial by $F$,
  i.e., $p$ is commutatively top-reducible by $F$.
\item\label{lem.cprop.3} In case $F$ is a free commutative standard basis,
  every polynomial $p \in \ideal{}{}(F)\backslash \{ 0 \}$
  is commutatively top-reducible to zero by $F$.
\lemend
\end{enumerate}
}
\end{lemma}
\Ba{}
\begin{enumerate}
\item This follows directly by adding up the polynomials used in the
   reduction steps  occurring in $p \red{*}{\myr}{c}{F} 0$.
\item This is an immediate consequence of definition \ref{def.csr}
  as the existence of a polynomial $f$ in $F$ and an element
  $w\in\m$ with $\hterm(f \mrm w)= \hterm(f) \cm w= \hterm(p)$ is guaranteed.
\item We show that every non-zero polynomial $p \in
  \ideal{r}{}(F)\backslash \{ 0 \}$ is top-reducible to zero using $F$ by induction on $\hterm(p)$.
  First let $\hterm(p) = \min \{ \hterm(g)\mid g \in \ideal{}{}(F)\backslash \{ 0 \} \}$.
  Then, as $p \in \ideal{}{}(F)$ and $F$ is a free commutative  standard
  basis, this gives us a representation $p = \sum_{i=1}^{k} \alpha_i
  \skm f_{i} \mrm w_i$, with  
   $\alpha_i \in \myk^*, f_{i} \in F, w_i \in \m$ and 
   $\hterm(p) \succeq \hterm(f_{i}) \cm w_i$
   for all  $1 \leq i \leq k$.
  Without loss of generality, let us assume $\hterm(p) = \hterm(f_1) \cm w_1$.
  Hence, the polynomial $p$ is  commutatively reducible by $f_1$.
  Let $p \red{}{\myr}{c}{f_1} q$, i.e.,
   $q = p - \hc(p) \skm \hc(f_1)^{-1} \skm f_1 \mrm w_1$, and
   by the definition of commutative reduction the term $\hterm(p)$ is
   eliminated from $p$ implying that $\hterm(q) \pred \hterm(p)$ as $q < p$.
  Now, since $\hterm(p)$ was minimal among the head terms of the elements
   in the ideal generated by $F$, this implies $q=0$, and, 
   therefore, $p$ is commutatively top-reducible to zero by $f_1$ in 
   one step.
   On the other hand, in case 
    $\hterm(p) \succ \min \{ \hterm(g) | g \in \ideal{}{}(F)\backslash \{ 0 \} \}$, by the
    same arguments used before we can commutatively reduce $p$ to a polynomial
    $q$ with $\hterm(q) \prec \hterm(p)$, and, thus, by our induction
    hypothesis we know that $q$ is commutatively top-reducible to zero.
   Therefore, as the reduction step $p \red{}{\myr}{c}{f_1} q$ takes
    place at the head term of $p$, the polynomial $p$ is also commutatively top-reducible
    to zero.
\\
\qed
\end{enumerate}\renewcommand{\baselinestretch}{1}\small\normalsize
Indeed, free commutative standard  bases and commutative
Gr\"obner bases are equivalent.
\begin{theorem}\label{theo.cgb.csr}~\\
{\sl
For a set $F$ of polynomials in $\myk[\m]$,
 the following statements are equivalent:
\begin{enumerate}
\item $F$ is a commutative Gr\"obner basis.
\item For all polynomials $g \in \ideal{}{}(F)$ we have $g \red{*}{\myr}{c}{F} 0$.
\item $F$ is a free commutative standard basis.
\lemend
\end{enumerate}
}
\end{theorem}
\Ba{}~\\
\mbox{$1 \R 2:$ }
  By (i) of definition \ref{def.cgb} we know that $g \in
  \ideal{}{}(F)$ implies $g \red{*}{\lr}{c}{F} 0$ and since
  $\red{}{\myr}{c}{F}$ is confluent and $0$ is irreducible,
 $g \red{*}{\myr}{c}{F} 0$ follows
  immediately.

\mbox{$2 \R 3:$ }
  This follows directly by adding up the polynomials used in the
   reduction steps  $g \red{*}{\myr}{c}{F} 0$.

\mbox{$3 \R 1:$ }
In order to show that $F$ is a commutative Gr\"obner basis we have
  to show two subgoals:
$\red{*}{\lr}{c}{F} = \;\; \equiv_{\ideal{}{}(F)}$ was already shown in lemma \ref{lem.csb.cong}.
It remains to show that $\red{}{\myr}{c}{F}$ is confluent.
             Since $\red{}{\myr}{c}{F}$ is Noetherian, we only have to prove
             local confluence.
             Suppose $g \red{}{\myr}{c}{F} g_1, g \red{}{\myr}{c}{F} g_2$
              and $g_1 \neq g_2$.
             Then $g_1 - g_2 \in \ideal{}{}(F)$ and, therefore, is commutatively
              top-reducible to zero as a result of lemma \ref{lem.csr.prop3}.
             Thus lemma \ref{lem.confluentc} provides the existence
              of a polynomial $h \in \myk[\m]$ such that
              $g_1 \red{*}{\myr}{c}{F} h$ and $g_2 \red{*}{\myr}{c}{F} h$, i.e.,
              $\red{}{\myr}{c}{F}$ is confluent.
\\
\qed
%
Since in general $\red{*}{\lr}{r}{} \neq \red{*}{\lr}{c}{}$ and 
  $\red{*}{\lr}{s}{} \neq \red{*}{\lr}{c}{}$ we again enrich our set of polynomials
 used for reduction to regain the desired expressiveness of strong reduction.
\begin{definition}\label{def.satc}~\\
{\rm
A set of polynomials $F\subseteq \{\alpha \skm p \mrm w \mid \alpha
\in \myk^*, w \in\m \}$
 is called a  
 \index{commutative!saturation}\index{saturation!commutative}\index{commutatively saturating set}\betonen{commutatively  saturating set}\/
 for a non-zero polynomial $p \in\myk[\m]$, if for all $ \alpha \in
 \myk^*$, $w \in \m$, with $\alpha \skm p \mrm w \neq 0$ we have  $\alpha \skm p \mrm w  \red{}{\myr}{c}{F} 0$\footnote{Since $\myk$ is a
   field it is sufficient to demand $p \mrm w \red{\leq 1}{\myr}{c}{F} 0$ for
   all $w \in \m$.}.
 $\SAT_c(p)$ denotes the family of all commutatively saturating sets for $p$.
We call a set $F \subseteq \myk[\m]$ 
 \index{commutatively saturated set}\index{saturated set!commutatively}\betonen{commutatively
  saturated}, if for all $\alpha \in \myk^*$, $f \in F$,  $w \in \m$, we
 have $\alpha \skm f \mrm w
 \red{}{\myr}{c}{F} 0$ in case $\alpha \skm f \mrm w \neq 0$.
\dend
}
\end{definition}
Note that commutatively saturating sets are also saturating sets.
Moreover, they give us additional information as they allow special
 representations of elements in the  ideals they generate comparable
 to lemma \ref{lem.prop1.1} for prefix saturated sets.
\begin{lemma}\label{lem.prop1.1c}~\\
{\sl
 Let $F$ be a commutatively saturated set of polynomials in $\myk[\m]$.
 Then every non-zero polynomial $g$ in $\ideal{}{}(F)$ has a representation $g =  \sum_{i=1}^{k} \alpha_i \skm f_i \mrm w_i$,
 where $\alpha_i \in \myk^*, f_i \in F, w_i \in \m$,
 and
 $\hterm(f_i \mrm w_i) = \hterm(f_i) \cm w_i$.
\lemend\ohnebeweis
}
\end{lemma}
Notice that these representations need not be free commutative standard
representations as we cannot conclude $\hterm(g) \succeq \hterm(f_i)
\cm w_i$.

As in the case of prefix saturation combined with prefix saturation,
commutative saturation enables us to do the same reduction steps we
can do using strong reduction or right reduction combined with
saturation (compare lemma \ref{lem.connection}).
\begin{lemma}\label{lem.connectionc}~\\
{\sl
Let $f,g,p$ be some polynomials in $\myk[\m], S \in\SAT(p)$,
 and $S_c \in\SAT_c(p)$. 
Then $$f  \red{}{\myr}{r}{S} g \mbox{ if and only if } f  \red{}{\myr}{c}{S_c} g.$$
\lemend\ohnebeweis
}
\end{lemma}
Further commutative reduction combined with commutative saturation is strong
enough to capture the ideal congruence.
\begin{lemma}\label{lem.congruencec}~\\
{\sl
Let $F$ be a commutatively saturated set of polynomials
and  $p,q$ some polynomials in $\myk[\m]$.
Then
  $$p \red{*}{\lr}{c}{F} q \mbox{ if and only if } p - q \in
  \ideal{}{}(F).$$
\lemend
}
\end{lemma}
\Ba{}~\\
This lemma follows directly from theorem \ref{theo.congruence}  and lemma
\ref{lem.connectionc}.
\\
\qed
The existence of finite commutatively saturating sets is guaranteed by
Dickson's lemma.
\begin{lemma}\label{lem.term}~\\
{\sl
Every commutatively saturating set for a polynomial contains a finite
 commutatively saturating set.
\lemend
}
\end{lemma}
\Ba{}~\\
For $p \in \myk[\m]$ let $S$ be an arbitrary saturating set.
Then comparable to the constructive approach for saturating sets in
definition \ref{def.sat} we can decompose $S$ into sets $S_t$ for $t
\in \terms(p)$ such that $S_t = \{ q \in S \mid \hterm(q) \mbox{ results
  from the term } \; t \}$.
Further let $Z_{t} = \{ \hterm(q) \mid q \in S_{t} \}\subseteq \m
 \subseteq \freecomm$.
Then $Z_{t}$ is a (possibly infinite) subset of $\freecomm$
 in the sense of Dickson's lemma and we can choose a finite subset $D_t$ of
 $Z_{t}$ such that for every $w \in Z_{t}$ there exist $u \in D_t$ and
 $v \in  \m \subseteq \freecomm$ such that $u \cm v = w$. 
Now for every $t \in \terms(p)$ we can associate every term $s \in  D_t$
with a polynomial $q_s \in S_t$ such that $\hterm(q_s) = s$.
The union of all such polynomials then is a finite saturating set for
$p$ which obviously is a subset of $S$.
\\
\qed
The property of being commutatively saturated is decidable.
\begin{lemma}~\\
{\sl 
It is decidable, whether a finite set $F \subseteq \myk[\m]$ is commutatively saturated.
}
\end{lemma}
\Ba{}~\\
We will show that the following procedure is correct:

\pagebreak

\renewcommand{\baselinestretch}{1}\small\normalsize
\procedure{Commutatively Saturated Check}%
{\vspace{-4mm}\begin{tabbing}
XXXXX\=XXXX \kill
\removelastskip
{\bf Given:} \>  A finite set $F \subseteq \myk[\m]$, and 
               $(\Sigma, T_c)$ a convergent  presentation of
               $\m$. \\
{\bf Find:} \> {\bf\em yes}, \=if $F$ is commutatively saturated \\
        \> {\bf\em no}, \>otherwise. 
\end{tabbing}
\vspace{-7mm}
\begin{tabbing}
XX\=XX\=XXX\=XXXXXXXX\=XXXXX\= \kill
answer := yes \\
{\bf for all}  $q \in F$ {\bf do} \\
\> $t := \hterm(q)$; \\
\> {\bf for all} \> \> $w \in C(t) = \{ $ \> $w \in \freecomm | \lcm(t,l) = t \cm
                               w = l \cm u \neq t \cm l \mbox{ for
                                 some rule }$ \\
\>   \>             \>             \> $(l,r) \in T_c, u \in \freecomm \}$ {\bf do}  \\
\>\> {\rm\kommentar \% $C(t)$ contains special 
  overlaps between  $t$ and  left hand sides of  rules in $T$} \\
\>\> {\rm\kommentar \% Multiplying $t$ with terms in $C(t)$
  leads to cancellation} \\
XX\=XX\=XX\=XXXX\=XXXXX\ \kill
\>  \> $q' := (q \mrm w_1)$; \\
\>  \> {\bf if}     \> $q'  \nred{}{\myr}{c}{F} 0$ and $q' \neq 0$ \\
\>  \> \>              {\bf then}   \> answer := no; \\
\>  \> {\bf endif}\\
\> {\bf endfor} \\
{\bf endfor} 
\end{tabbing}}
\renewcommand{\baselinestretch}{1.1}\small\normalsize\\
It remains to show that the answer of our procedure is ``no'' if and only if
 $F$ is not commutatively saturated.
Obviously, the answer ``no'' implies the existence of an element $w \in \m$ such that for some 
 $f \in F$ with $f \mrm w \neq 0$,  $f \mrm w \nred{}{\myr}{c}{F} 0$.
On the other hand let us assume that our procedure gives us ``yes'', but $F$ is not commutatively saturated.
Then there exists an element  $w \in \m$ and a polynomial $f \in F$ such that $\hterm(f) \cm w$ is minimal according to our
ordering on $\freecomm$, $f \mrm w \neq 0$ and $f \mrm w \nred{}{\myr}{c}{F} 0$.
In case $w \in C(\hterm(f))$ this would give us ``no'' contradicting our assumption.
Therefore, let us assume $w \not\in C(\hterm(f))$.
Furthermore, $\hterm(f) \cm w$ must be $T$-reducible as otherwise
 $f \mrm w \red{}{\myr}{c}{f \in F} 0$.
Let $ w=w_1 \cm w_2$ such that $\hterm(f) \cm w_1 = l \cm v =
\lcm(\hterm(f),l)$ for some rule $(l,r) \in T$, i.e., $w_1 \in C(\hterm(f))$.
Now $f \mrm w_1$ is considered by our procedure and since $w \succ w_1$ and the choice of
 $w,f$ was minimal, we either get $f \mrm w_1 = 0$ contradicting the
 fact that $f \mrm w \neq 0$ or $f \mrm w_1 \red{}{\myr}{c}{F} 0$.
Furthermore, this gives us the existence of $f' \in F$ such that
 $\hterm(f) \cm w_1 \succ \hterm(f \mrm w_1) = \hterm(f') \cm z$
 for some $z \in \m$.
This implies $\hterm(f) \cm w = \hterm(f) \cm (w_1 \cm w_2) = (\hterm(f)
\cm w_1) \cm w_2 \succ \hterm(f \mrm w_1) \cm w_2 = (\hterm(f') \cm z) \cm
w_2 =\hterm(f') \cm (z \cm w_2) \succeq \hterm(f') \cm (z \mm w_2)$ and gives us
 $f \mrm w = \alpha \skm f' \mrm (z \mm w_2) \red{}{\myr}{c}{F} 0$
 contradicting our assumption.
\\
Further this procedure terminates, as the sets $H$ and $T$ are
always finite.
\\
\qed
The next lemma states that minimal commutatively saturating sets
exist.
\begin{lemma}~\\
{\sl
Let $p$ be a polynomial in $\myk[\m]$ and $S \in \SAT_c(p)$.
Then if there is a polynomial $q \in S$ such that $q \red{}{\myr}{c}{S
 \backslash \{ q \}} 0$, the set $S \backslash \{ q \}$ is a commutatively saturating
 set for $p$.
\lemend
}
\end{lemma}
\Ba{}~\\
This is an immediate consequence as $p \red{}{\myr}{c}{q_1} 0$ and $q_1
\red{}{\myr}{c}{q_2} 0$ implies $p \red{}{\myr}{c}{q_2} 0$, i.e., item 3
of lemma \ref{lem.red} also holds for commutative reduction.
\\
\qed 
We proceed  to give a procedure, which actually computes a
simplified commutatively saturating set for a polynomial $p$ (compare
the approach for prefix saturation using simplification on page \pageref{prefix.sauration.using.simplification}).
The idea is to compute overlaps of terms with rules in $T_c$ using least
common multiples in $\freecomm$.
In giving a procedure for prefix saturating a polynomial, prefixes were used to
define a set of critical overlaps for a term $t$ of the form
$C(t) = \{ w \in \Sigma^* \mid tw \id t_1t_2w \id t_1l, t_2 \neq \lambda
\mbox{ for some rule } (l,r) \in T \}$.
Here we will look at the set
$C(t) = \{ w \in \freecomm \mid \lcm(t,l) = t \cm w = l \cm u \neq t \cm
l \mbox{ for some rule } (l,r) \in T_c, u \in \freecomm \}$
which corresponds to the overlaps between the term $t$ and the rules
in $T_c$ in terms of semi-Thue systems modulo commutativity.
Additionally, polynomials are removed from the computes set in case
 they are commutatively reducible to zero in one step
 by a polynomial computed later on.


\procedure{Commutative Saturation using Simplification\protect{\label{commutative.saturation}}}%
{\vspace{-4mm}\begin{tabbing}{ll}
XXXXX\=XXXX \kill
\removelastskip
{\bf Given:} \> A polynomial $p \in \myk[\m]$, and 
              $(\Sigma , T_c)$ a convergent  presentation of
               $\m$. \\
{\bf Find:} \> $S_c \in \SAT_c(p)$.
\end{tabbing}
\vspace{-7mm}
\begin{tabbing}
XX\=XX\=XXX\=XXXXXXXX\=XXXXX\=XXX\=XXX\= \kill
$S_0$:= $\{p \}$; \\
$H$ := $\{p \}$; \\
$i$ := $0$; \\
{\bf while} $H \neq \emptyset$ {\bf do} \\
\> $i$ := $i+1$; \\
\> $S_i$ := $S_{i-1}$; \\
\> $q := {\rm remove}(H)$; \\
\>{\rm\kommentar \% Remove an element from a set using a fair
  strategy}\\
\> $t := \hterm(q)$; \\
\> {\bf for all}  \>\> $w \in C(t) = \{ $ \> $w \in \freecomm | \lcm(t,l) = t \cm
                               w = l \cm u \neq t \cm l \mbox{ for
                                 some rule }$ \\
\>                \>\>             \> $(l,r) \in T_c, u \in \freecomm \}$ {\bf do}  \\
\>\> {\rm\kommentar \% $C(t)$ contains special
  overlaps between  $t$ and  left hand sides of  rules in $T$} \\
\> \> $q' := (q \mrm w_1)$; \\
XX\=XX\=XX\=XXXX\=XXXXX\=XXX\=XXX\= \kill
\>  \> {\bf if}  \>$q'  \nred{}{\myr}{c}{S_c} 0$ and $q' \neq 0$ \\
\>  \>  \>{\bf then} \> $S_i := {\rm simplify}(S_i, q') \cup \{q' \}$; \\
\>  \>  \>           \>{\rm\kommentar \% 
                     Simplify removes elements $s$ from $S_i$ in
                     case $s \red{}{\myr}{c}{q'} 0$}\\
\>  \> \>            \> $H := H \cup \{q' \}$; \\
\>  \> {\bf endif} \\
\> {\bf endfor} \\
{\bf endwhile}\\ 
$S_c$ := $S_i$
\end{tabbing}}

\begin{lemma}~\\
{\sl
Procedure {\sc Commutative Saturation using Simplification} is correct.
\theoend
}
\end{lemma}
\Ba{}~\\
We will prove  that for all $q \in S_c$,$w \in \m$ we have $ q \mrm w
\red{}{\myr}{c}{S_c} 0$ in case $q \mrm w \neq 0$.
Let us assume this is not true.
Then we can choose a  counter-example $q \mrm w$ such that
 $\hterm(q) \cm w$ is minimal (according to the ordering on $\freecomm$)  
 and $q \mrm w \nred{}{\myr}{c}{S_c} 0$.
Then $\hterm(q) \cm w$ must be $T$-reducible, as otherwise
 $q \mrm w \red{}{\myr}{c}{q} 0$
 and $q \in S_c$.
Let $ w=w_1 \cm w_2$ such that $w_1 \in \m$  causing
 $\hterm(q) \mm w_1 \neq \hterm(q) \cm w_1 = l \cm z = \lcm(\hterm(q),l)$ for
some $(l,r) \in T$, $z \in \freecomm$.
Since $q \in S_c$, the polynomial $q \mrm w_1$ is considered
 during the computation of $S_c$\footnote{We can assume that the
   elements are removed from the set $H$ using a fair strategy, e.g.,
   first in first out.}.
We have to consider two cases.
 If $q \mrm w_1 \in S_c$ then $q \mrm w = (q \mrm w_1) 
       \mrm w_2 \red{}{\myr}{c}{S_c} 0$ since $w_1 \in \m$ and
       $\hterm(q) \cm w \succ \hterm(q \mrm w_1) \cm w_2$ 
       contradicting our assumption.
On the other hand, if $q \mrm w_1 \not\in S_c$ 
       then for some iteration step $i_q$ we have
       $q \mrm w_1 \red{}{\myr}{c}{q' \in S_{i_q}} 0$ and first we show
       that even $q \mrm w_1 \red{}{\myr}{c}{S_c} 0$ holds.
      \auskommentieren{
      It remains to show that $q \mrm w_1 \red{}{\myr}{c}{S_c} 0$.
      \\
      This follows immediately, if we can show that for all
       $f \in \{ p \mrm w | w \in \m \}$, $f \red{}{\myr}{c}{S_i} 0$
       implies $f \red{}{\myr}{c}{S_{i+1}} 0$ for all iterations of the
       {\bf while} loop.
      \\
      $f \red{}{\myr}{c}{s \in S_i} 0$ implies $f = \alpha \skm s \mrm u$ for some
       $\alpha \in \myk^*$, $u \in \m$ and $\hterm(f) = \hterm(s) \cm u$.
      In case $s \in S_{i+1}$ we are done.
      \\
      Otherwise there exists a polynomial $s' \in S_{i+1}$ such that 
       $s \red{}{\myr}{c}{s'} 0$, i.e., $s = \alpha' \skm s' \mrm u'$
       for some $\alpha' \in \myk^*$, 
       $u' \in \m$, $\hterm(s) = \hterm(s') \mm_{\freecomm} u'$ and $u
       \mm u' = u \cm u'$.
      But then, 
       $$f = \alpha \skm s \mrm u = \alpha \skm (\alpha' \skm s' \mrm
       u') \mrm u = (\alpha \skm \alpha') \skm s' \mrm (u' \cm u)$$
      and
       $$\hterm(f) = \hterm(s) \cm u =
        (\hterm(s') \cm u') \cm u =
         \hterm(s') \cm u'\cm u$$
       imply that $f \red{}{\myr}{p}{s' \in S_{i+1}} 0$.
     \\}
     This follows since the fact that $q \red{}{\myr}{c}{q_1}
     0$ and $q_1 \red{}{\myr}{c}{q_2} 0$ implies $q \red{}{\myr}{c}{q_2}
     0$ yields that $q \mrm w_1 \red{}{\myr}{c}{S_i} 0$ implies $q \mrm
     w_1 \red{}{\myr}{c}{S_{i+n}} 0$ for all $n \in \n$.
     Furthermore, as removing a polynomial $s$ from a set $S_i$
     because of a polynomial $q'$ we know $\hterm(q') \tup
     \hterm(s)$ and this cannot occur infinitely often, we would get
     $q \mrm w_1 \red{}{\myr}{c}{S_c} 0$ contradicting our assumption.
     Thus we can assume $q \mrm w_1 \red{}{\myr}{c}{S_c} 0$ and moreover, we know $w_1 \in
       \m$ and 
       $\hterm(q)\mm_{\freecomm} w_1 \succ \hterm(q \mrm w_1) =
        \hterm(q')\mm_{\freecomm} z$ for some $z \in\m$.
      Further $q \mrm w_1 = \alpha' \skm q' \mrm z$ and 
       $\hterm(q \mrm w_1) = \hterm(q')\mm_{\freecomm} z$ give us 
       $q \mrm w = (q \mrm w_1) \mrm w_2 = (\alpha' \skm q' \mrm z)
       \mrm w_2 = \alpha' \skm q' \mrm (z \mm w_2)$,
       and $\hterm(q)\cm w = \hterm(q) \cm (w_1 \cm w_2) = (\hterm(q)
       \cm w_1) \cm w_2 \succ (\hterm(q') \cm z) \cm w_2 = \hterm(q') \cm
       (z \cm w_2) \succeq \hterm(q')\cm (z \mm w_2)$.
      Therefore,
       $q \mrm w =  (\alpha' \skm q' \mrm z) \mrm w_2 =
        \alpha' \skm q' \mrm (z \mm w_2)  \red{}{\myr}{c}{S_c} 0$,
       contradicting our assumption.
\\
\qed
%
\begin{lemma}~\\
{\sl
Procedure {\sc Commutative Saturation using Simplification} terminates.
}
\end{lemma}
\Ba{}~\\
This follows at once as the procedure is correct and by lemma
\ref{lem.term} the constructed set must contain a finite commutatively
saturated subset.
\\
\qed
Commutative saturation enriches a polynomial $p$ by adding a set of
polynomials $S \in \SAT_c(p)$ such that we can
 substitute a reduction step $q \red{}{\myr}{(s,r)}{p}q'$ by a
 commutative reduction step $q \red{}{\myr}{c}{p' \in S}q'$.
This additional information can be combined with special s-polynomials
to give a finite confluence test similar to the approach using prefix reduction. 
%
\begin{definition}\label{def.cpc}~\\
{\rm
Given two non-zero polynomials $p_{1}, p_{2} \in \myk[\m]$ such that
 $\hterm(p_{1}) \cm w_{1} = \hterm(p_{2}) \cm w_{2}= \lcm(\hterm(p_1),
 \hterm(p_2)) \in \m$ for some $w_{1}, w_{2} \in \m$, then  the 
 \betonen{commutative \linebreak s-polynomial}\index{commutative!s-polynomial}\index{s-polynomial!commutative}
 is defined as 

 $\mbox{}\hfill\spol{c}(p_{1}, p_{2}) = \hc(p_1)^{-1} \skm p_1 \mrm w_1 -
 \hc(p_2)^{-1} \skm p_2 \mrm w_2.\hfill\diamond$
}
\end{definition}
Obviously this definition implies $\hterm(\spol{c}(p_{1}, p_{2})) \pred  \lcm(\hterm(p_1),
 \hterm(p_2))$.
An s-polynomial is called non-trivial in case it is not zero.
Notice that a finite set $F \subseteq \myk[\m]$ only gives us finitely
 many commutative s-polynomials.
A commutative s-polynomial for two polynomials $p_1$, $p_2$
 corresponds to a term, namely $\lcm(\hterm(p_1),
 \hterm(p_2)) \in \m$, where both
polynomials $p_1$ and $p_2$ can be applied to perform a commutative
reduction step.
As in the previous section (commutative) Gr\"obner bases cannot be
 characterized by such restricted s-polynomials  alone.
However, this can be done in case we have a commutatively saturated set.
%
%
\auskommentieren{
\begin{lemma}\label{lem.redc}~\\
{\sl
Let $F \subseteq \myk[\m]$ and $p,q \in \myk[\m]$. 
Further, let $p \red{}{\myr}{c}{q} 0$ and $q \red{*}{\myr}{c}{F} 0$
 result in  the representations $p = \alpha \skm q \mrm w$ and 
 $q = \sum_{i=1}^{k} \alpha_i \skm g_i \mrm w_i$, where
 $\alpha, \alpha_i \in \myk^*, g_i \in F, w,w_i \in \m$.
\\
Then the following statements hold:
\begin{enumerate}
\item There exists $s \in \{ 1, \ldots , k \}$ such that
  $\hterm(p)=\hterm(g_s \mrm w_s \mrm w) = \hterm( g_s \mrm w_s) \cm w$.
\item For all remaining terms 
       $t \in (\bigcup_{i=1 \atop i \neq s}^k \terms(g_i \mrm w_i \mrm w)$
       we have $\hterm(p) \succ t$.
\lemend
\end{enumerate}
}
\end{lemma}
\Ba{}~\\
Since $p \red{}{\myr}{c}{q} 0$ and $p = \alpha \skm q \mrm w$ we know $\hterm(p)
= \hterm(q) \cm w$.
\\
Further, as $\sum_{i=1}^{k} \alpha_i \skm g_i \mrm w_i$ belongs to
 the reduction sequence $q \red{*}{\myr}{c}{F} 0$,
 there exists an index $s \in \{ 1, \ldots , k \}$ such
 that $\hterm(q)=\hterm(g_s) \cm w_s$,
 $\hc(q)=\alpha_s \skm \hc(g_s)$ and for all
 $v \in \bigcup_{i=1 \atop i \neq s}^{k} \terms(g_i \mrm w_i)$
 we have $\hterm(q) \succ v$ implying
 $\hterm(p)=\hterm(q) \cm w \succ v \cm w \succeq v \mm w$. \\
\qed}
\begin{theorem}\label{theo.ccp}~\\
{\sl
For a commutatively saturated set $F$ of polynomials in $\myk[\m]$, the
 following statements are equivalent:
\begin{enumerate}
\item For all polynomials $g \in \ideal{}{}(F)$ we have $g \red{*}{\myr}{c}{F} 0$. 
\item For all polynomials $f_{k}, f_{l} \in F$ we have 
  $\spol{c}(f_{k}, f_{l}) \red{*}{\myr}{c}{F} 0$.
\theoend
\end{enumerate}
}
\end{theorem}
\Ba{}~\\
\mbox{$1 \R 2:$ }
 In case $\hterm(f_{k}) \cm w_k = \hterm(f_{l}) \cm w_{l} =
 \lcm(\hterm(f_k),\hterm(f_l)) \in \m$ for some elements
  $w_{k}, w_{l} \in \m$,
 then by definition \ref{def.cpc} we get
 $$ \spol{c}(f_{k}, f_{l}) = \hc(f_k)^{-1} \skm f_k \mrm
 w_k - \hc(f_l)^{-1} \skm f_l \mrm w_l \:\in \ideal{}{}(F),$$ 
 and hence $\spol{c}(f_{k}, f_{l}) \red{*}{\myr}{c}{F} 0$.

\mbox{$2 \R 1:$ }
    We have to show that every non-zero element  $g \in
    \ideal{}{}(F)\backslash \{ 0 \}$ is $\red{}{\myr}{c}{F}$-reducible
    to zero.
    Remember that for
      $h \in \ideal{}{}(F)$, $ h \red{}{\myr}{c}{F} h'$ implies $h' \in \ideal{}{}(F)$.
     As  $\red{}{\myr}{c}{F}$ is Noetherian
      it suffices to show that every  $g \in \ideal{}{}(F) \backslash
      \{ 0 \}$ is $\red{}{\myr}{c}{F}$-reducible.
     Let $g = \sum_{j=1}^m \alpha_{j} \skm f_{j} \mrm w_{j}$ be a
      representation of a non-zero polynomial $g$ such that $\alpha_{j} \in \myk^*, f_j \in F,
      w_{j} \in \m$.
     By lemma \ref{lem.prop1.1c} we can assume  $\hterm(f_{i} \mrm w_{i})
     = \hterm(f_{i}) \cm w_{i}$.
This will be important to restrict ourselves to commutative
s-polynomials in order to modify this representation of $g$.
     Depending on the above  representation of $g$ and a well-founded total ordering $\succeq$ on $\m$ we define
      $t = \max \{ \hterm(f_{j}) \cm w_{j} \mid j \in \{ 1, \ldots m \}  \}$ and
      $K$ is the number of polynomials $f_j \mrm w_j$ containing
      $t$ as a term.
Then $t \succeq \hterm(g)$ and in case $\hterm(g) = t$ this immediately implies that $g$ is
$\red{}{\myr}{c}{F}$-reducible. 
So
by lemma \ref{lem.cprop} it is sufficient to  show that
$g$ has a free commutative standard representation, as this implies that $g$ is
top-reducible using $F$.
This will be done by induction
    on $(t,K)$, where
      $(t',K')<(t,K)$ if and only if $t' \prec t$ or $(t'=t$ and $K'<K)$\footnote{Note
        that this ordering is well-founded since $\succ$ is and $K \in\n$.}.
Now if $t \succ \hterm(g)$
      there are two polynomials $f_k,f_l$ in the corresponding
       representation\footnote{Not necessarily $f_l \neq f_k$.}
      such that  $\hterm(f_k) \cm w_k = \hterm(f_l) \cm w_l$.
     By definition \ref{def.cpc} we have a commutative s-polynomial
      $\spol{c}(f_k,f_l) = \hc(f_k)^{-1} \skm  f_k \mrm z_1-
      \hc(f_l)^{-1} \skm f_l \mrm z_2$ and  we further know 
      $\hterm(f_k) \cm w_k = \hterm(f_l)
      \cm w_l =  \hterm(f_k) \cm z_1 \cm w = \hterm(f_l) \cm z_2 \cm w$ for
      some $z_1,z_2,w \in \m$ such that $\hterm(f_k) \cm
      z_1=\hterm(f_l) \cm z_2= \lcm(\hterm(f_k), \hterm(f_l)) \in \m$.
Note that as in the case of prefix reduction in theorem \ref{theo.pcp}
we can localize ourselves to this part of the original overlap
$\hterm(f_k) \mm w_k = \hterm(f_l) \mm w_l$.
We will now change our representation of $g$ by using the additional
information on this s-polynomial in such a way that for the new
representation of $g$ we either have a smaller maximal term or the occurrences of the term $t$
are decreased by at least 1.
     Let us assume  $\spol{c}(f_k,f_l) \neq 0$\footnote{In case  $\spol{c}(f_k,f_l) = 0$,
               just substitute $0$ for the sum $\sum_{i=1}^n \delta_i \skm
               h_i \mrm v_i$ in the equations below.}. 
     Hence,  the reduction sequence  $\spol{c}(f_k,f_l)
     \red{*}{\myr}{c}{F} 0 $ results in a free commutative standard representation of the form 
     $\spol{c}(f_k,f_l) =\sum_{i=1}^n \delta_i \skm h_i \mrm v_i$ such that
     $\delta_i \in \myk^*,h_i \in F,v_i \in \m$
      and all terms occurring in the sum are bounded by $\hterm(\spol{c}(f_k,f_l))$.
     Then by lemma \ref{lem.redc} we can conclude that $t$ is a real bound for all terms occurring
      in the sum $\sum_{i=1}^n \delta_i \skm h_i \mrm v_i \mrm w$.
     Furthermore, we can assume that this representation is of the required form,
      as we can substitute all polynomials $h_i$, where
      $\hterm(h_i \mrm v_i \mrm w_k) \neq \hterm(h_i) \cm (v_i \mm w)$ 
      without increasing
      $t$ or $K$.
     This gives us: 
     \begin{eqnarray}
       &  & \alpha_{k} \skm f_{k} \mrm w_{k} + \alpha_{l} \skm f_{l} \mrm w_{l}  \nonumber\\ 
       &  &    \nonumber \\
       & = &  \alpha_{k} \skm f_{k} \mrm w_{k} + \underbrace{ \alpha'_{l} \skm \beta_k \skm f_{k} \mrm w_{k}
                   - \alpha'_{l} \skm \beta_k \skm f_{k} \mrm w_{k}}_{=\, 0} 
                   + \alpha'_{l}\skm \beta_l  \skm f_{l} \mrm w_{l} \nonumber\\
       &  &    \nonumber \\ 
       & = & (\alpha_{k} + \alpha'_{l} \skm \beta_k) \skm f_{k} \mrm w_{k} - \alpha'_{l} \skm
             \underbrace{(\beta_k \skm f_{k} \mrm w_{k}
             -  \beta_l \skm f_{l} \mrm w_{l})}_{=\, \spol{c}(f_k,f_l) \mrm w} \nonumber\\
       & = & (\alpha_{k} + \alpha'_{l} \skm \beta_k) \skm f_{k} \mrm w_{k} - \alpha'_{l} \skm
                   (\sum_{i=1}^n \delta_{i} \skm h_{i} \mrm v_{i} \mrm w) \label{s6}
     \end{eqnarray}
     where  $\beta_k = \hc(f_k)^{-1}$,$\beta_l = \hc(f_l)^{-1}$
      and  $\alpha'_l \skm \beta_l = \alpha_l$.
     By substituting (\ref{s6}) in our representation of $g$  
 either $t$ disappears  or in
 case $t$ remains maximal among the terms occurring in the new
 representation of $g$, $K$ is decreased.
\\
\qed
Now theorem \ref{theo.ccp} gives rise to the following procedure that computes a Gr\"obner
 basis with respect to $\red{}{\myr}{c}{}$ for a finitely generated
 ideal similar to procedure {\sc Prefix Gr\"obner Bases}
 on page \pageref{prefix.groebner.bases}.
The resulting set is additionally commutatively saturated.


\procedure{Commutative Gr\"obner Bases}%
{\vspace{-4mm}\begin{tabbing}
XXXXX\=XXXX \kill
\removelastskip
{\bf Given:} \> A finite set of polynomials $F \subseteq \myk[\m]$. \\
{\bf Find:} \> $\gb_c(F)$, a commutative Gr\"obner basis
                 of $F$.\\
{\bf Using:} \> $\s_c$ a commutative saturating procedure for polynomials. 
\end{tabbing}
\vspace{-7mm}
\begin{tabbing}
XX\=XX\=XXXX\=XX\=XXXX\=XXXX\kill
$G$ := $\bigcup_{f \in F} \s_c(f)$\\
$B$ := $\{ (q_{1}, q_{2}) \mid q_{1}, q_{2} \in G \}$; \\
{\rm\kommentar \% $B$ is used to check statement 2 of theorem \ref{theo.ccp}} \\
{\bf while} $B \neq \emptyset$ {\bf do} \\
\>      $(q_{1}, q_{2})$ := {\rm remove}$(B)$; \\
\>      {\rm\kommentar \% Remove an element from a set} \\
\>      {\bf if} \> $\spol{c}(q_{1}, q_{2})$ exists \\
\>               \> {\rm\kommentar \% Compute the s-polynomial  if
  it is not trivial} \\
\> \>     {\bf then} \> $h$ := {\rm normalform}$(\spol{c}(q_{1}, q_{2}),
                              \red{}{\myr}{c}{G})$;   \\
\> \> \>  {\rm\kommentar \% Compute a normal form using
        commutative reduction}\\
\>\>  \>    {\bf if} \> $h \neq 0$  \\
\>\>\>               \> {\rm\kommentar \% Statement 2 of theorem
  \ref{theo.ccp} does not hold} \\
\>\>  \> \> {\bf then}  \>$G$ := $G \cup \s_c(h)$; \\
\>\>  \> \>             \> {\rm\kommentar \% $G$ is extended to achieve
  statement 2} \\
\>\>  \> \>             \> {\rm\kommentar \% $G$ is commutatively
  saturated} \\
\>\>  \> \>             \>$B$ := $B \cup \{ (f, {\tilde h}) \mid f \in G, {\tilde h} \in  \s_c(h)  \}$; \\
\>\>  \>    {\bf endif}\\
\>      {\bf endif}\\
{\bf endwhile} \\
$\gb_c(F)$ := $G$
\end{tabbing}}

\begin{lemma}~\\
{\sl
Procedure {\sc Commutative Gr\"obner Bases} terminates.
\lemend
}
\end{lemma}
\Ba{}~\\
New polynomials are only added in case an s-polynomial does not reduce
 to zero.
Hence, if our procedure would not terminate there would be an infinite
 sequence of normal forms of s-polynomials added contradicting the fact
 that the set of head terms of
 these polynomials has a finite basis via $\cm$ as a subset of
 $\freecomm$ by Dickson's lemma and the ideal can be characterized by
 the head terms of the polynomials it contains as described
 in theorem \ref{theo.cequiv}. 
\\
\qed
\begin{lemma}~\\
{\sl
Procedure {\sc Commutative Gr\"obner Bases} is correct.
\lemend
}
\end{lemma}
\Ba{}~\\
This follows immediately from theorem \ref{theo.ccp}.
\\
\qed

The sets characterized in theorem \ref{theo.ccp} are commutative Gr\"obner
bases and hence right Gr\"obner bases, but they are required to be
commutatively saturated.
Reviewing example \ref{exa.fcsb} we see that there exist right
 Gr\"obner bases in commutative monoid rings which are not  commutatively saturated.
\begin{example}~\\
{\rm
Let $\Sigma = \{ a,b \}$ and $T_c = \{ a^2 \myr \lambda, b^2 \myr \lambda
\}$ be a presentation of a commutative monoid $\m$ (which is in fact a
 group)  with  a  length-lexicographical ordering induced by $a \succ b$.
\\
Then the set $\{ ab + \lambda \} \subseteq \q[\m]$ itself is
  a right Gr\"obner basis,
  but neither commutatively  saturated nor a commutative Gr\"obner basis. 
We have $a + b \in \ideal{}{}(ab + \lambda)$
  but $a+b \nred{}{\myr}{c}{ab + \lambda} 0$.
\exaend
}
\end{example}
Note that a commutative Gr\"obner basis need not be commutatively saturated.
\begin{example}\label{exa.ac=d.bc=e.c}~\\
{\rm
Let $\Sigma = \{ a,b,c,d, e \}$ and $T_c = \{ ac \myr d, bc \myr e \}$ be a
 presentation of a commutative monoid $\m$ with
  a length-lexicographical ordering induced by $a \succ b \succ c \succ d \succ
 e$.
\\
Then the set $F = \{ a+b,d+ \lambda , e - \lambda \}$ is
 a commutative Gr\"obner basis in ${\bf Q}[\m]$.
This can be seen by studying the right ideal generated by $F$,
 $\ideal{}{}(F) = \{ \alpha_1 \skm (a + b)\mrm w_1 + \alpha_2 \skm (d
 - \lambda)\mrm w_2
    + \alpha_3 \skm (e - \lambda)\mrm w_3 + \alpha_4 \skm (d + e)\mrm w_4| \alpha_i \in \q,
    w_i \in \m, w_1 \neq c \cm w' \}$.
But $F$ is not commutatively saturated,
 as $(a+b) \mrm c = d+ e \nred{}{\myr}{c}{F} 0$.
We only have $d+e \red{2}{\myr}{c}{F} 0$.
\exaend
}
\end{example}
Next we will give a characterization of commutative Gr\"obner bases
 as free commutative standard bases without demanding that the set of
 polynomials is commutatively saturated.
This is important as interreducing a set of commutatively saturated
 polynomials destroys this property.
But interreducing a free commutative standard basis again gives us a
free commutative standard basis.
Remember that this is also true for prefix standard bases but not for
stable standard bases in general.
\begin{theorem}\label{theo.altccp}~\\
{\sl
For a set  $F$ of polynomials in $\myk[\m]$,
 equivalent are:
\begin{enumerate}
\item Every polynomial $g \in \ideal{}{}(F)$ has a free commutative standard representation.
\item \begin{enumerate}
      \item For all polynomials $f \in F$ and all elements  $w \in\m$,
             the polynomial $f \mrm w$ has a  free commutative
             standard representation.
      \item For all polynomials $f_{k}, f_{l} \in F$ the non-trivial  
             commutative s-polynomials have  free commutative standard representations.
      \end{enumerate}
\end{enumerate}
}
\end{theorem}
\Ba{}~\\
\mbox{$1 \R 2:$ }
This follows immediately.

\mbox{$2 \R 1:$ }
    We have to show that every non-zero element  $g \in
    \ideal{}{}(F)\backslash \{ 0 \}$ has a  free commutative
             standard representation.
    Let $g = \sum_{j=1}^m \alpha_{j} \skm f_{j} \mrm w_{j}$ be an arbitrary
      representation of a non-zero polynomial $g$ such that $\alpha_{j} \in \myk^*, f_j \in F,
      w_{j} \in \m$.
By our assumption and lemma \ref{lem.csr.prop2}
 we can assume  $\hterm(f_{i} \mrm w_{i}) =
 \hterm(f_{i}) \cm w_{i}$ as $f_{i} \in F$ and every $f_{i} \mrm
 w_{i}$ has a free commutative standard representation\footnote{Note that these
 free commutative standard representations do not yield a free
 commutative  standard
 representation for the polynomial $g$,
 as $\hterm(g) \prec \hterm(f_i) \cm w_i$ is possible.}.
Now using statement (a) and (b) we can show as in theorem
\ref{theo.ccp} how the representation of $g$ can be transformed into a
free commutative standard representation.
\auskommentieren{
     Depending on this  representation of $g$ and a well-founded total ordering $\succeq$ on $\m$ we define
      $t = \max \{ \hterm(f_{j}) \cm w_{j} \mid j \in \{ 1, \ldots m \}  \}$ and
      $K$ is the number of polynomials $f_j \mrm w_j$ containing
      $t$ as a term.
     \\
Then $t \succeq \hterm(g)$ and
in case $\hterm(g) = t$ this immediately implies that the
representation of $g$ is a free commutative standard representation. 
\\
So it is sufficient to  show that
$g$ has such a representation with $\hterm(g) =t$.
This will be done by induction on $(t,K)$, where
      $(t',K')<(t,K)$ if and only if $t' \prec t$ or $(t'=t$ and $K'<K)$\footnote{Note
        that this ordering is well-founded since $\succ$ is and $K \in\n$.}.
\\
In case $t \succ \hterm(g)$
      there are two polynomials $f_k,f_l$ in the corresponding
       representation\footnote{Not necessarily $f_l \neq f_k$.}
      such that  $\hterm(f_k) \cm w_k = \hterm(f_l) \cm w_l$.
     By definition \ref{def.cpc} we have a commutative s-polynomial
      $\spol{c}(f_k,f_l) = \hc(f_k)^{-1} \skm  f_k \mrm z_1-
      \hc(f_l)^{-1} \skm f_l \mrm z_2$ and  we further know 
      $\hterm(f_k) \cm w_k = \hterm(f_l)
      \cm w_l =  \hterm(f_k) \cm z_1 \cm w = \hterm(f_l) \cm z_2 \cm w$ for
      some $z_1,z_2,w \in \m$ such that $\hterm(f_k) \cm
      z_1=\hterm(f_l) \cm z_2= \lcm(\hterm(f_k), \hterm(f_l))$.
\\
We will now change our representation of $g$ by using the additional
information on this s-polynomial in such a way that for the new
representation of $g$ we either have a smaller maximal term or the occurrences of the term $t$
are decreased by at least 1.
\\
     Let us assume  $\spol{c}(f_k,f_l) \neq 0$\footnote{In case  $\spol{c}(f_k,f_l) = 0$,
               just substitute $0$ for the sum $\sum_{i=1}^n \delta_i \skm
               h_i \mrm v_i$ in the equations below.}. 
\\
     By our assumption $\spol{c}(f_k,f_l)$ has a free commutative standard representation of the form
     $\spol{c}(f_k,f_l) =\sum_{i=1}^n \delta_i \skm h_i \mrm v_i$, where
     $\delta_i \in \myk^*$,$h_i \in F$, and $v_i \in \m$.
     By lemma \ref{lem.redc} we can conclude that $t$ is a real bound
      for all terms occurring 
      in the sum $\sum_{i=1}^n \delta_i \skm h_i \mrm v_i \mrm w$.
     Furthermore, we can assume that this representation is of the required form,
      as we can substitute all polynomials $h_i$, where
      $\hterm(h_i \mrm v_i \mrm w_k) \neq \hterm(h_i) \cm (v_i \mm w)$ 
      without increasing
      $t$ or $K$.
     \\
     This gives us: 
     \begin{eqnarray}
       &  & \alpha_{k} \skm f_{k} \mrm w_{k} + \alpha_{l} \skm f_{l} \mrm w_{l}  \nonumber\\ 
       &  &    \nonumber \\
       & = &  \alpha_{k} \skm f_{k} \mrm w_{k} + \underbrace{ \alpha'_{l} \skm \beta_k \skm f_{k} \mrm w_{k}
                   - \alpha'_{l} \skm \beta_k \skm f_{k} \mrm w_{k}}_{=\, 0} 
                   + \alpha'_{l}\skm \beta_l  \skm f_{l} \mrm w_{l} \nonumber\\
       &  &    \nonumber \\ 
       & = & (\alpha_{k} + \alpha'_{l} \skm \beta_k) \skm f_{k} \mrm w_{k} - \alpha'_{l} \skm
             \underbrace{(\beta_k \skm f_{k} \mrm w_{k}
             -  \beta_l \skm f_{l} \mrm w_{l})}_{=\, \spol{c}(f_k,f_l) \mrm w} \nonumber\\
       & = & (\alpha_{k} + \alpha'_{l} \skm \beta_k) \skm f_{k} \mrm w_{k} - \alpha'_{l} \skm
                   (\sum_{i=1}^n \delta_{i} \skm h_{i} \mrm v_{i} \mrm w) \label{s7}
     \end{eqnarray}
     where  $\beta_k = \hc(f_k)^{-1}$,$\beta_l = \hc(f_l)^{-1}$
      and  $\alpha'_l \skm \beta_l = \alpha_l$.
     By substituting (\ref{s7}) in our representation of $g$  
 either $t$ disappears  or in
 case $t$ remains maximal among the terms occurring in the new
 representation of $g$, $K$ is decreased.
}
\\
\qed
Note that this theorem gives us a stronger characterization of commutative
 Gr\"obner bases, as it does not require them to be commutatively saturated.
\begin{remark}~\\
{\rm
Let $F$ be a set of polynomials in $\myk[\m]$.
\begin{enumerate}
\item If $F$ is commutatively saturated, then for every polynomial
       $f \in F$ and every element $w \in \m$ the polynomial $f \mrm
       w$ has a free commutative standard representation.
      This follows immediately, since $f \mrm w \red{}{\myr}{c}{F} 0$
       implies that there exists a polynomial
       $f' \in F$ such that $f \mrm w \red{}{\myr}{c}{f'} 0$ and 
       $\hterm(f \mrm w) = \hterm(f') \cm u$ for some $u \in \m$.
      Note that $F$ need not be a free commutative standard basis.  
\item On the other hand, if for all polynomials $f \in F$ and all
       elements $w \in \m$ the polynomial $f \mrm w$ has a free commutative
       standard representation, this need not imply that $F$ is commutatively
       saturated. 
      To see this, let us review example \ref{exa.ac=d.bc=e.c}\\
      Let $\Sigma = \{ a,b,c,d, e \}$ and $T_c = \{ ac \myr d, bc \myr e \}$ be a
       semi-Thue system modulo commutativity presenting a commutative
       monoid $\m$ with  a length-lexicographical ordering induced by
       $a \succ b \succ c \succ d \succ
       e$.
      \\
      Then for every polynomials $f$ in 
       the set $F = \{ a+b,d+ \lambda , e - \lambda \}$ and every element
       $w \in \m$ the multiple $f \mrm w$ has a  free commutative
       standard representation.
      For the elements $w \in \m$ where $(a+b) \mrm w = a \cm w + b
       \cm w$, $(d+ \lambda) \mrm w = d \cm w+ w$ and
       $(e - \lambda) \mrm w = e \cm w -w$, these are  free commutative
       standard representations.
      It remains to check the case $(a+b) \mrm (c \cm w) = d \cm w + e
       \cm w$.
      Since $d+e \red{}{\myr}{c}{d + \lambda}
             e - \lambda \red{}{\myr}{c}{e - \lambda} 0$,
       we have a  free commutative standard representation 
       $d \cm w + e \cm w = (d+ \lambda) \mrm w + (e - \lambda ) \mrm w$,
       but $d \cm w + e \cm w$ does not commutatively reduce to zero using $F$
       in {\em one}\/ step.
      We only have $(d+e) \mrm w \red{2}{\myr}{c}{F} 0$.
      \remend      
\end{enumerate}
}
\end{remark}
\begin{definition}~\\
{\rm
We call a set of polynomials $F \subseteq \myk[\m]$ \index{weakly commutatively
  saturated}\index{commutative saturated!weakly}\betonen{weakly commutatively
  saturated}, if for all $f \in F$ and all $\alpha \in \myk^*$, $w \in
 \m$ we have $\alpha \skm f
 \mrm w \red{*}{\myr}{c}{F} 0$.
\dend
}
\end{definition}
If a set of polynomials $F$ is weakly commutatively saturated this implies
 that for all $f \in F$ and all $w \in \m$ the polynomial $f \mrm w$,
 in case it is non-zero, has a free commutative standard representation.
Thus we can give the following procedure to compute reduced 
 commutative Gr\"obner bases where weak saturation is ensured by the
 use of a
 saturating procedure for polynomials.


\procedure{Reduced Commutative Gr\"obner Bases\protect{\label{reduced.commutative.groebner.bases}}}
{\vspace{-4mm}\begin{tabbing}
XXXXX\=XXXX \kill
\removelastskip
{\bf Given:} \> A finite set $F \subseteq \myk[\m]$. \\
{\bf Find:} \> $\gb(F)$, a commutative Gr\"obner basis of $F$. \\
{\bf Using:} \> A procedure $\s_c$ for computing commutatively
saturating sets.
\end{tabbing}
\vspace{-7mm}
\begin{tabbing}
XX\=XX\= XXXX \= XXXX \=\kill
$G_0$ := $\emptyset$; \\
$S_0$ := $F$; \\
$i$ := $0$; \\
{\bf while} $S_i \neq \emptyset$ do\\
\> $i$ := $i+1$; \\
\> $q_i$ := {\rm remove}$(S_{i-1})$; \\
\>{\rm\kommentar \% Remove an element  using a fair strategy}\\
\> $q_i'$ := ${\rm normalform}(q_i, \red{}{\myr}{c}{G_{i-1}})$; \\
\>{\rm\kommentar \% Compute a normal form using
           commutative reduction} \\
\>{\bf if} \>$q_i' \neq 0$ \\
\>         \>{\rm\kommentar \% Statement 2 of theorem \ref{theo.altccp} does
  not hold} \\
\> \>{\bf then}  \>$H_i$ := $\{ g \in G_{i-1} \mid
                              \hterm(g) \mbox{ is commutatively reducible using }
                              q_i' \}$;\\
\> \>            \>{\rm\kommentar \% These polynomials would
  have new head terms after commutative} \\
\> \>            \>{\rm\kommentar \% reduction using $q_i'$} \\
\> \>\>$G_i$ := {\rm reduce}$((G_{i-1} \backslash H_i) \cup \{ q_i'
\})$; \\
\>\>\>{\rm\kommentar
  \% {\rm reduce}$(F) = \{ {\rm normalform}(f, \red{}{\myr}{c}{F \backslash \{
    f \}}) | f \in F \}$\footnotemark} \\ 
\> \>  \>{\rm\kommentar \% No head term of a polynomial in $G_i$ is
  commutatively reducible by the} \\
\> \>  \>{\rm\kommentar \% other polynomials in $G_i$} \\
\> \>\>$S_i$ := \> $S_{i-1}  \cup H_i \cup 
              \bigcup_{g \in (G_i \backslash G_{i-1})}(\s_c(g)
              \backslash \{ g \})$\\
\> \> \>        \> $\cup \{ \spol{c}(f,g) \mid f \in
              G_i, g \in G_i \backslash G_{i-1}\}$; \\
\> \>{\bf else}  \> $G_i$ := $G_{i-1}$; \\
\> \> \> $S_i$ := $S_{i-1}$; \\
\>{\bf endif}\\   
{\bf endwhile} \\
$\gb (F)$:= $G_i$
\end{tabbing}}
\footnotetext{Notice that only the reducts of the polynomials are touched in this procedure.}

The sets $S_i$ will contain saturating sets and s-polynomials
corresponding to polynomials added to some set $G_j$.
Hence, in general this set will contain polynomials that are no longer
necessary, e.g. in case a polynomial is removed from a set $G_j$
neither its s-polynomials nor all saturating polynomials need to be
considered.
This is due to the fact that the conditions named in theorem
\ref{theo.altccp} must hold for the final set $G_k$ only.
Hence it is possible to develop marking strategies in order to keep the sets
$S_i$ smaller.

The following lemmata ensure the correctness of our procedure.
\begin{lemma}\label{lem.idealpropertyc}~\\
{\sl
Let $F$ be a set of polynomials in $\myk[\m]$ and $G_i$, $S_i$, $i \in
\n$ the respective sets in procedure {\sc Reduced Commutative Gr\"obner
  Bases}.
Then  we have
$$\ideal{}{}(F) = \ideal{}{}(G_i \cup S_i).$$
\lemend}
\end{lemma}
\Ba{}~\\
For $i = 0$ we have $G_0 \cup S_0 = F$ and hence
 $\ideal{}{}(F) = \ideal{}{}(G_0 \cup S_0)$.
For $i > 0$ let $G_{i-1}$, $S_{i-1}$ be the respective sets before
 entering the {\bf while} loop for its $i$-th iteration.
Further let $q_i$ be the polynomial chosen from $S_{i-1}$ and 
 $q_i'$ a  normal from of $q_i$ with respect to commutative reduction using $G_{i-1}$.
Then in case $q_i' =0$ we know $q_i \in \ideal{}{}(G_{i-1})$ and thus
 as $G_i = G_{i-1}$ and $S_i = S_{i-1} \backslash \{ q_i \}$ we can
 conclude
 $$\ideal{}{}(G_i \cup S_i) = \ideal{}{}(G_{i-1} \cup S_{i-1})=\ideal{}{}(F).$$
In case $q_i' \neq 0$, we  know $\ideal{}{}({\rm reduce}((G_{i-1} \backslash H_i)
\cup \{ q_i' \})) = \ideal{}{}((G_{i-1} \backslash H_i)
\cup \{ q_i' \})$, $\ideal{}{}(G_{i-1} \cup \{ q_i' \} \cup (S_{i-1} \backslash \{ q_i \}) ) 
 = \ideal{}{}(G_{i-1} \cup S_{i-1})$ and
$\ideal{}{}(\bigcup_{g \in (G_i \backslash G_{i-1})}(\s_c(g) \backslash \{
g \})\cup \{ \spol{c}(f,g) \mid f \in
              G_i, g \in G_i \backslash G_{i-1}\}) \subseteq \ideal{}{}(G_i)$.\\
Abbreviating $\bigcup_{g \in (G_i \backslash G_{i-1})}(\s_c(g) \backslash \{
g \})$ by $S_c$ we can conclude
\begin{eqnarray}
  &\phantom{=}& \ideal{}{}(G_i \cup S_i) \nonumber\\  
  &=& \ideal{}{}(G_i \cup (S_{i-1} \backslash \{ q_i \}) \cup H_i \cup 
            S_c \cup \{ \spol{c}(f,g) \mid f \in
              G_i, g \in G_i \backslash G_{i-1}\}) \nonumber\\ 
  &=& \ideal{r}{}(G_i \cup (S_{i-1} \backslash \{ q_i \}) \cup H_i)  \nonumber\\
  &=& \ideal{r}{}({\rm reduce}((G_{i-1} \backslash H_i)
      \cup \{ q_i' \}) \cup (S_{i-1} \backslash \{ q_i \}) \cup H_i) \nonumber\\ 
  &=& \ideal{}{}((G_{i-1} \backslash H_i)
      \cup \{ q_i' \} \cup  (S_{i-1} \backslash \{ q_i \}) \cup H_i) \nonumber\\
  &=& \ideal{}{}(G_{i-1} 
      \cup \{ q_i' \} \cup (S_{i-1} \backslash \{ q_i \})) \nonumber\\ 
  &=& \ideal{}{}(G_{i-1} \cup S_{i-1}) \nonumber\\ 
  &=& \ideal{r}{}(F). \nonumber
\end{eqnarray}
\nopagebreak
\qed
\begin{remark}~\\
{\rm
As in the case of procedure {\sc Reduced Prefix Gr\"obner Bases} on
page \pageref{reduced.prefix.groebner.bases} the head terms of the
polynomials in the sets $G_i$ fulfill that every term that had a
divisor in a set $\hterm(G_k)$ will then also have a divisor in each
set $\hterm(G_{k+n})$, $n \in \n$.
Hence, as the polynomials added are always in normal form and the sets
$\hterm(G_i)$ can be regarded as subsets of ${\cal T}$, Dickson's
lemma yields termination.
\remend
}
\end{remark}
\begin{theorem}\label{theo.corr.rcgb}~\\
{\sl
Let $G$ be the set generated by procedure
 {\sc Reduced Commutative Gr\"obner Bases} on a finite input $F \subseteq \myk[\m]$.
Then the following statements hold:
\begin{enumerate}
\item $\ideal{}{}(G) = \ideal{}{}(F)$.
\item For all polynomials $f \in G$ and all elements  $w \in\m$,
             the polynomial $f \mrm w$ has a  free commutative
             standard representation in case it is non-zero.
\item For all polynomials $f_{k}, f_{l} \in F$ the non-trivial  
             commutative s-polynomials have  free commutative standard representations.
\item  $G$ is a reduced commutative Gr\"obner basis.
\end{enumerate}

\theoend
}
\end{theorem}
\Ba{}~\\
Since procedure {\sc Reduced Commutative Gr\"obner Bases} terminates we
 have $G = G_k$ for some $k \in \n$ and $S_k = \emptyset$.
Now as by lemma \ref{lem.idealpropertyc} we have $\ideal{}{}(F) =
\ideal{}{}(G_i \cup S_i)$ for all $i \in \n$, this implies
$\ideal{}{}(F) = \ideal{}{}(G)$.
Since according to our construction all s-polynomials corresponding to polynomials in
 $G$ are not commutatively reducible to zero, it remains to show that for all
$f \in G$, $w \in \m$ the right 
 multiple $f \mrm w$ has a free commutative standard representation with respect
 to $G$.
We will first show that for all $i \leq k$, $f \in S_i \cup G_i$ implies
 that $f$ has a free commutative standard representation with respect to $G$.
This can be done by induction on $j$ where $i = k-j$.
The case $j = 0$ is trivial, as this implies $S_k = 0$ and every $f \in G_k$ has
 a free commutative standard representation with respect to $G = G_k$.
Hence let $f \in S_{k-(j+1)} \cup G_{k-(j+1)}$ and suppose
 $f \not\in S_{k-j} \cup G_{k-j}$, as otherwise we are already done.
In case $f \in  G_{k-(j+1)}$, as $f \not\in  G_{k-j}$, $f$ must be
 commutatively reducible by $q_{k-j}'$.
Then either $f \in S_{k-j}$ if $\hterm(f)$ is commutatively reducible
 by $q_{k-j}'$ 
 and our induction hypothesis then yields the
 existence of a free commutative standard representation for $f$.
Or $f$ is the result of reduction during the computation of the
set $G_{k-j} = {\rm reduce}((G_{k-(j+1)} \backslash H_{k-j}) \cup
\{ q_{k-j}' \})$.
But then by lemma \ref{lem.csr.prop3}, $f$ has a free commutative
 standard representation with respect to $G_{k-j}$ yielding the
 existence of a free commutative
 standard representation with respect to $G$, as every polynomial in
 $G_{k-j}$ has one.
In case $f \in  S_{k-(j+1)}$, as $f \not\in  S_{k-j}$, $f$ is chosen
 to compute the polynomial
 $q_{k-j}' = {\rm normalform}(f, \red{}{\myr}{c}{G_{k-{j+1}}})$.
Now in case $q_{k-j}' = 0$ we have $G_{k-(j+1)} = G_{k-j}$ and
 hence $f$ has a free commutative standard representation with respect to
 $G_{k-j}$.
Furthermore, as every polynomial in $G_{k-j}$ by induction hypothesis
 has a free commutative standard representation with respect to $G$, by lemma
 \ref{lem.csr.prop3} we are done.
If $q_{k-j}' \neq 0$ we get that $f$ has a free commutative standard
 representation with respect to $G_{k-(j+1)} \cup \{ q_{k-j}' \}$ and
 even with respect to $G_{k-j} = {\rm reduce}((G_{k-(j+1)} \backslash H_{k-j}) \cup
\{ q_{k-j}' \})$.
Again by the induction hypothesis and by lemma
\ref{lem.csr.prop3}, $f$ has a free commutative standard representation
with respect to $G$.
\\
Now let us return to our initial goal, to prove that for all
 polynomials $f \in G$ and all elements $w \in \m$, the polynomial
 $f \mrm w$ has a free commutative standard representation in case it
 is non-zero.
Since $f \in G$, the polynomials in $\s_c(f) \backslash \{ f \}$ have been added to some
 set $S_i$, i.e., there exists a polynomial $f' \in S_i \cup G_i$ such
 that $f \mrm w \red{}{\myr}{c}{f'} 0$ and $f'$ has a free commutative standard
 representation.
Then by lemma \ref{lem.csr.prop3} we can conclude that $f \mrm w$ also
 has a free commutative standard representation.
\\
By construction, $G$ is a reduced set of polynomials and by
 theorem \ref{theo.altccp} we can conclude that $G$ is a reduced commutative Gr\"obner basis.
\\
\qed
We end  this section by sketching how Buchberger's algorithm can
 be lifted to quotient rings of polynomial rings and how these
 quotients are related to finitely generated commutative monoid
 rings. 
The ideas are mainly the same as stated by Dei{\ss} in \cite{De89}.
Let us start by recalling Buchberger's definition of reduction for a
 polynomial ring $\myk[X_1, \ldots , X_n]$.
Remember that $p \red{}{\myr}{b}{f} q$ at a monomial $\alpha \skm t$, if $\hterm(f) \cm u = t$ for some $u \in {\cal T}$, and $q = p - \alpha \skm \hc(f)^{-1} \skm f \mrm u$.

Henceforth we will assume that our finitely generated commutative
monoid $\m$ is presented by a finite convergent semi-Thue system modulo commutativity
 $T_c \subseteq \freecomm \times \freecomm$ where ${\cal T}$ is again the free commutative
 monoid generated now by an alphabet $\Sigma = \{ X_1, \ldots , X_n
 \}$ to show the connection to the ordinary polynomial ring.
Then we can consider the set $R$ as a set of polynomials 
 $P_{T_c} = \{ l - r | (l,r) \in T_c \}$.
Then $P_{T_c}$ is a Gr\"obner basis with respect to the
 well-founded admissible completion ordering $\succeq$ on $\freecomm$ related to
 the presentation $(\Sigma, T_c)$.

We can show that in this context the monoid ring $\myk[\m]$ is in fact
 isomorphic to a quotient\index{quotient!of a polynomial ring}
 of the ordinary polynomial ring $\myk[X_1, \ldots , X_n]$.
\begin{lemma}~\\
{\sl
Let $\m$ be a commutative monoid presented by $(\Sigma,T_c)$.
Then $\myk[\m]$ is isomorphic to
 $\myk[X_1, \ldots , X_n]/\ideal{}{}(P_{T_c})$.
\lemend
}
\end{lemma}
\Ba{}~\\
Let $\varphi : \myk[X_1, \ldots , X_n] \myr \myk[\m]$ be the natural ring homomorphisms
       defined by setting
      $\varphi(\sum_{i=1}^k \alpha_i \skm w_i) = 
       \sum_{i=1}^k \alpha_i \skm [w_i]_{\m}$ with  
      $\alpha_i \in\myk$, $w_i \in\freecomm$.
      Then as in lemma \ref{lem.prefix.quotient} one can show that
      the kernel of $\varphi$ is an ideal in
      $\myk[X_1, \ldots , X_n]$, namely $\ideal{}{\myk[X_1, \ldots , X_n]}(P_{T_c})$, and
      hence $\myk[\m]$ is isomorphic to the ring 
       $\myk[X_1, \ldots , X_n]/\ideal{}{\myk[X_1, \ldots , X_n]}(P_{T_c})$.
\\
\qed
Now let us state how computation in our quotient structure is done.
For a polynomial  $p =
 \sum_{i=1}^m \alpha_i \skm t_i \in \myk[X_1, \ldots, X_n]$ we define 
 $$[p]_{\ideal{}{}(P_{T_c})} = \sum_{i=1}^m \alpha_i \skm [t_i[_{\m} = \sum_{i=1}^m \alpha_i \skm t_i\nf{P_{T_c}}$$
 and  we will
 write $p$ instead of $[p]_{\ideal{}{}(P_{T_c})}$ to denote elements of 
 $\myk[X_1, \ldots , X_n]/\ideal{}{}(P_{T_c})$ or $p\nf{P_{T_c}}$ if we want
 to turn a polynomial in $\myk[X_1, \ldots, X_n]$ into an element of the quotient.
\begin{definition}~\\
{\rm 
Let $p,q$ be two polynomials in $\myk[X_1, \ldots ,
X_n]/\ideal{}{}(P_{T_c})$.
Then we can define \index{addition! in a commutative quotient}
 \betonen{addition}
 and \index{multiplication!in a commutative quotient}
 \betonen{multiplication} as follows:
\begin{enumerate}
\item $p \qadd q = [p +_{\myk[X_1, \ldots , X_n]}
  q]_{\ideal{}{}(P_{T_c})}= (p +_{\myk[X_1, \ldots , X_n]} q)\nf{P_{T_c}}$
\item $p \qmult  q = [p \mrm_{\myk[X_1, \ldots , X_n]} q]_{\ideal{}{}(P_{T_c})}=(p \mrm_{\myk[X_1, \ldots , X_n]} q)\nf{P_{T_c}}$
\end{enumerate}
such that  $+_{\myk[X_1, \ldots , X_n]}$ and $\mrm_{\myk[X_1, \ldots , X_n]}$
 are the corresponding ring operations in the polynomial ring
 $\myk[X_1, \ldots, X_n]$.
\mbox{\phantom{XX}}\dend
}
\end{definition}
Then $\myk[X_1, \ldots , X_n]/\ideal{}{}(P_{T_c})$ together with $\qadd$ and 
 $\qmult $ is a commutative ring with unit.
We will now introduce quotient reduction to this structure by lifting
 Buchberger's reduction.
\begin{definition}\label{def.redq}~\\
{\rm
Let $p,f$ be two non-zero polynomials in $\myk[X_1, \ldots , X_n]/\ideal{}{}(P_{T_c})$.
Then we set\index{reduction!in a commutative quotient}
$p \red{}{\myr}{}{f} q$  at a monomial $\alpha \skm t$ of $p$ if there exists a polynomial
 $q'$ in $\myk[X_1, \ldots, X_n]$ such that
 $p \red{}{\myr}{b}{f} q'$ at $\alpha \skm t$ and $q = q'\nf{P_{T_c}}$.
We can define $\red{*}{\myr}{}{}, \red{+}{\myr}{}{}$, and
 $\red{n}{\myr}{}{}$ as usual.
Reduction by a set $F \subseteq \myk[X_1, \ldots , X_n]/\ideal{}{}(P_{T_c})$
 is denoted by $p \red{}{\myr}{}{F} q$ and stands for 
 $p \red{}{\myr}{}{f} q$ for some $f \in F$,
 also written as  $p \red{}{\myr}{}{f \in F} q$.
\dend
}
\end{definition}
For this reduction we can now state:
\begin{lemma}~\\
{\sl
Let $F$ be a set of polynomials and
 $p,q$ some polynomials in $\myk[X_1, \ldots , X_n]/\ideal{}{}(P_{T_c})$.
\begin{enumerate}
\item $\red{}{\myr}{}{F} \subseteq \red{*}{\myr}{b}{F \cup P_{T_c}}$, i.e.,
  reduction in the quotient can be simulated by reduction in the polynomial ring
  using additional polynomials.
\item $p \red{}{\myr}{}{F} q$ implies $p > q$.
\item $\red{}{\myr}{}{F}$ is Noetherian.
\lemend
\end{enumerate}
}
\end{lemma}
\Ba{}~\\
\vspace{-8mm}
\begin{enumerate}
\item This follows by the definition of reduction, since
       $p \red{}{\myr}{}{f}q$ can be simulated by 
       $p \red{}{\myr}{b}{f}q' \red{*}{\myr}{b}{P_{T_c}} q$.
\item This follows immediately as Buchberger's reduction already has
       this property.
\item This follows from the fact that Buchberger's reduction is
  Noetherian.
\\
\qed
\end{enumerate}\renewcommand{\baselinestretch}{1}\small\normalsize
But although $\myk[X_1, \ldots , X_n]/\ideal{}{}(P_{T_c})$ still is a Noetherian ring,
 many other properties of reduction  in $\myk[X_1, \ldots , X_n]$ are lost as for example
 the quotient may contain zero-divisors.
\begin{lemma}~\\
{\sl
Let $p,q,h \in \myk[X_1, \ldots , X_n]/\ideal{}{}(P_{T_c})$ and $h \neq 0$.
\begin{enumerate}
\item $q < p$ no longer implies $q \qmult  h < p \qmult  h$.
\item $p \red{}{\myr}{}{p} 0$ no longer implies $p \qmult  h
  \red{}{\myr}{}{p} 0$.
\lemend
\end{enumerate}
}
\end{lemma}
\begin{example}\label{exa.qsdefect}~\\
{\rm
Let $\Sigma = \{ X_1 \}$ and $P_T = \{ X_1^2 -1 \}$ be a presentation
 of a commutative monoid.
\\
Then $p=X_1$, $q=1$ and $h = X_1$ gives us an appropriate
counter-example, as we have $q < p$, but $q \qmult h = X_1 > p \qmult
h = 1$, and $p \qmult h = 1$ cannot be quotient reduced to zero by $p$.
\exaend
}
\end{example}
In the polynomial ring a basis of an ideal is called a Gr\"obner basis
 if Buchberger's reduction using it is confluent.
We can easily extend this definition to our quotient structure.
\begin{definition}~\\
{\rm
A set of polynomials $G \subseteq \myk[X_1,\ldots,X_n]/\ideal{}{}(P_{T_c})$ is said to be a
 \index{Gr\"obner basis!in a commutative quotient}\betonen{Gr\"obner
   basis} with respect to  $\red{}{\myr}{}{}$, if
\begin{enumerate}
\item $\red{*}{\Longleftrightarrow}{}{G} = \;\;\equiv_{\ideal{}{}(G)}$, and
\item $\red{}{\myr}{}{G}$ is confluent.
\dend
\end{enumerate}
}
\end{definition}
Unfortunately, 
 a polynomial alone is no longer  a Gr\"obner basis, as in 
 example \ref{exa.qsdefect}  the set $\{ p \}$ is no
 Gr\"obner basis.
This is due to the fact that reduction in the quotient
 $\myk[X_1, \ldots ,X_n]/\ideal{}{}(P_{T_c})$ no longer captures the ideal
 congruence.
In our example we have $X_1 \equiv_{\ideal{}{}(p)} 1$ but
 $X_1 \nred{*}{\Longleftrightarrow}{}{p} 1$.

In order to describe Gr\"obner bases, let us continue  by giving
 a sufficient condition for confluence.
\begin{lemma}\label{lem.qs.confluent}~\\
{\sl
Let $P_{T_c}$ be a Gr\"obner basis in $\myk[X_1, \ldots ,X_n]$,
 $F \subseteq \myk[X_1, \ldots , X_n]/\ideal{}{}(P_{T_c})$.
Then if $\red{}{\myr}{b}{F \cup P_{T_c}}$ is confluent on $\myk[X_1, \ldots,
 X_n]$, $\red{}{\myr}{}{F}$ is confluent on $\myk[X_1, \ldots , X_n]/\ideal{}{}(P_{T_c})$.
\lemend
}
\end{lemma}
\Ba{}~\\
Suppose there exist $f,h_1,h_2 \in \myk[X_1, \ldots , X_n]/\ideal{}{}(P_{T_c})$
 such that we get $f \red{}{\myr}{}{F} h_1$ and $f \red{}{\myr}{}{F} h_2$.
Then we can view these polynomials as elements of 
 $\myk[X_1, \ldots , X_n]$ and substitute $\red{}{\myr}{}{F}$ by 
 $\red{*}{\myr}{b}{F\cup P_{T_c}}$ giving us
 $f \red{*}{\myr}{b}{F \cup P_{T_c}} h_1$ and 
 $f \red{*}{\myr}{b}{F \cup P_{T_c}} h_2$.
Hence, as $\red{}{\myr}{b}{F \cup P_{T_c}}$ is confluent, there exists
 a polynomial $g \in \myk[X_1, \ldots , X_n]$ such that
 $h_1 \red{*}{\myr}{b}{F \cup P_{T_c}} g$ and
 $h_2 \red{*}{\myr}{b}{F \cup P_{T_c}} g$.
Since $\red{}{\myr}{b}{F \cup P_{T_c}}$ is convergent we can use the following
 reduction strategy:
\begin{enumerate}
\item Do as many reduction steps as possible using $P_{T_c}$.
\item If possible apply one reduction step using $F$ and return to 1.
\end{enumerate}
We stop as soon as no more reduction steps are possible.
Note that in this fashion we can combine reduction steps as required in the
 definition of 
 reduction $\red{}{\myr}{}{F}$.
This gives us that $h_1 \red{*}{\myr}{}{F} \tilde{g}$ and
 $h_2 \red{*}{\myr}{}{F} \tilde{g}$, where $\tilde{g} = g\nf{P_{T_c}}$.
\\
\qed 
The converse is not true as in example \ref{exa.qsdefect}, 
 $\red{}{\myr}{}{X_1}$ is confluent on $\myk[X_1]/\ideal{}{}(X_1^2 -1)$,
 but $\red{}{\myr}{b}{\{ X_1 \} \cup \{ X_1^2 -1 \}}$
 is not confluent on $\myk[X_1]$.
\\
In order to use this lemma to sketch how a Gr\"obner basis
 with respect to quotient reduction can be computed, we use
  Buchberger's s-polynomials.
Remember that the s-polynomial 
  for two polynomials
 $p,q \in \myk[X_1, \ldots,X_n]$ is defined as
$ \spol{}(p,q) = \hc(p)^{-1} \skm p \mrm  u - \hc(q)^{-1} \skm q \mrm v,$
 where  $\lcm(\hterm(p),\hterm(q)) = \hterm(p) \cm u = \hterm(q) \cm v$.
We can thus compute an
 (even reduced) Gr\"obner basis by modifying
 Buchberger's algorithm as follows:

Compute the (reduced) Gr\"obner basis $G'$ of $F \cup P_{T_c}$ in
 $\myk[X_1, \ldots, X_n]$ with respect to Buchberger's
 reduction without changing the polynomials in $P_{T_c}$.
Then the set $G := G' \backslash P_T$ is the reduced Gr\"obner basis
 of $F$ in $\myk[X_1, \ldots, X_n]/\ideal{}{}(P_{T_c})$.

Computing the Gr\"obner basis of $F \cup P_{T_c}$ three kinds of
 s-polynomials can arise:
\begin{enumerate}
\item $f,g \in P_{T_c}$:
      Then the corresponding s-polynomial can be omitted, as
       $P_T$ is already a Gr\"obner basis.
\item $f,g \in F$:
      Then the s-polynomial corresponds to the commutative
       s-polynomial as defined in \ref{def.cpc}.
\item $f \in F$, $g \in P_{T_c}$:
      Then the s-polynomial corresponds to the process of
       saturating $f$, in particular to the step of overlapping
       the head term of $f$ with the rule $l \myr r$ where $g = l-r$
       (compare procedure {\sc Commutative Saturation} on page \pageref{commutative.saturation}).
\end{enumerate}
Thus we have a characterization of Gr\"obner bases in the quotient
structure as follows:
\begin{theorem}~\\
{\sl 
Let $P_{T_c} \subseteq \myk[X_1, \ldots, X_n]$ be a reduced Gr\"obner basis 
 in $\myk[X_1, \ldots, X_n]$ and
 $F$ a set of polynomials in $\myk[X_1, \ldots, X_n]/\ideal{}{}(P_{T_c})$.
Then the following statements are equivalent:
\begin{enumerate}
\item $F$ is a Gr\"obner basis in $\myk[X_1, \ldots, X_n]/\ideal{}{}(P_{T_c})$.
\item
\begin{enumerate}
\item For all $f,g \in F$ we have $(\spol{}(f,g))\nf{P_{T_c}} \red{*}{\myr}{}{F} 0$, and
\item for all $f \in F$, $g \in P_{T_c}$ we have $(\spol{}(f,g))\nf{P_{T_c}} \red{*}{\myr}{}{F} 0$.
\theoend
\end{enumerate}
\end{enumerate}
}
\end{theorem}
\Ba{}~\\
\mbox{$1 \R 2:$ }
This follows immediately as all s-polynomials
   lie in $\ideal{}{\myk[X_1, \ldots,
   X_n]/\ideal{}{}(P_{T_c})}(F)$ and therefore are congruent to zero.
Thus the confluence of $F$ implies that they can be reduced to zero
 using $F$.

\mbox{$2 \R 1:$ }
The statements (a) and (b) imply that $F \cup P_{T_c}$ is a Gr\"obner basis
 in the polynomial ring $\myk[X_1, \ldots, X_n]$, i.e., $\red{}{\myr}{b}{F \cup P_{T_c}}$ is
 confluent.
Hence, by lemma \ref{lem.qs.confluent} $\red{}{\myr}{}{F}$ is also
 confluent.
\\
It remains to show that 
 $\red{*}{\Longleftrightarrow}{}{F} =\;\; \equiv_{\ideal{}{}(F)}$.
Obviously, $\red{*}{\Longleftrightarrow}{}{F} \subseteq\;\;
\equiv_{\ideal{}{}(F)}$.
On the other hand let $p$ and $q$ be polynomials in 
 $\myk[X_1, \ldots, X_n]/\ideal{}{}(P_{T_c})$.
Then $p \equiv_{\ideal{}{}(F)} q$ in $\myk[X_1, \ldots, X_n]/\ideal{}{}(P_{T_c})$
 implies $p \equiv_{\ideal{}{}(F \cup P_{T_c})} q$ in
 $\myk[X_1, \ldots, X_n]$.
Further, as $F \cup P_{T_c}$ is a Gr\"obner basis in
 $\myk[X_1, \ldots, X_n]$, we know 
 $\equiv_{\ideal{}{}(F \cup P_{T_c})} = \red{*}{\lr}{b}{F \cup P_{T_c}}$,
 and as $\red{}{\myr}{b}{F \cup P_{T_c}}$ is confluent,
 $p \red{*}{\lr}{b}{F \cup P_{T_c}} q$ implies 
 $p \downarrow_{F \cup P_{T_c}} q$.
Thus, as in lemma \ref{lem.qs.confluent}, we can conclude
 $p \Downarrow_{F} q$ giving us $p \red{*}{\Longleftrightarrow}{}{F} q$.
\\
This completes the proof that $F$ is a Gr\"obner basis in 
 $\myk[X_1, \ldots, X_n]/\ideal{}{}(P_{T_c})$.
\\
\qed
\auskommentieren{
The following lemma proves that sets as characterized in this theorem
 indeed are Gr\"obner bases.
\begin{lemma}~\\
{\sl
Let $F \subseteq \myk[X_1, \ldots, X_n]/\ideal{}{}(P_T)$,
 $\;\;p,q,h \in \myk[X_1, \ldots, X_n]/\ideal{}{}(P_T)$.
\begin{enumerate}
\item
Let $p-q \red{}{\myr}{r}{F} h$.
Then there are  $p',q' \in \myk[X_1, \ldots, X_n]/\ideal{}{}(P_T)$ such that 
 $p  \red{*}{\myr}{}{F} p', q  \red{*}{\myr}{}{F} q'$ and $h=p'-q'$.
\item
Let $0$ be a normal form of $p-q$ with respect to $F$.
Then there exists a polynomial 
 $g \in \myk[X_1, \ldots,X_n]/\ideal{}{}(P_T)$
 such that
 $p  \red{*}{\myr}{}{F} g$ and $q  \red{*}{\myr}{}{F} g$.
\lemend
\end{enumerate}
}
\end{lemma}
\Ba{}~\\
\begin{enumerate}
\item  Let $p-q \red{}{\myr}{}{f \in F} h$,  where the reduction takes
        place at the monomial $c \skm t$, $c \in \myk$,
        $t \in \freecomm$.
       Then $\hterm(f) \cm w = t$, for some $w \in \freecomm$.
       Further let $a_1$ be the coefficient of $t$ in $p$ and
        $a_2$ the coefficient of $t$ in $q$.
       We have to distinguish three cases:
       \begin{enumerate}
         \item $t \in \terms(p)$ and $t \in \terms(q)$: \\
               Then we can set
                $p \red{}{\myr}{b}{f} p - a_1 \skm \hc(f)^{-1} \skm f \mrm w
                  =: p'$,
                $q \red{}{\myr}{}{f} q - a_2 \skm \hc(f)^{-1} \skm f \mrm w
                  =:  q'$,
                where $f \mrm w$ is the normal form of the
                ``commutative'' polynomial $f \cm w$ with respect
                to $P_T$ in the polynomial ring.
         \item $t \in \terms(p)$ and $t \not\in \terms(q)$: \\
               Then $p \red{}{\myr}{}{f} p - c \skm \hc(f)^{-1}\skm f \mrm w=: p'$ 
                and $q =: q'$.
         \item $t \in \terms(q)$ and $t \not\in \terms(p)$: \\
               Then $q \red{}{\myr}{}{f} q + c  \skm f \mrm w=: q'$
                and $p =: p'$.
       \end{enumerate}
      In all three cases we have $p' -q' =  p - q - a \skm f \mrm w = h$.
\item We show our claim by induction on $k$, where $p-q \red{k}{\myr}{}{F} 0$.
      \\
      In the base case $k=0$ there is nothing to show.
      \\
      Let $p-q \red{}{\myr}{}{F} h  \red{k}{\myr}{}{F} 0$.
      \\
      Then by (1) there are $p',q' \in \myk[X_1, \ldots,X_n]/\ideal{}{}(P_T)$ such that 
       $p \red{*}{\myr}{}{F} p', q  \red{*}{\myr}{}{F} q'$ and $h=p'-q'$.
      \\
      Now the induction hypothesis for $p'-q' \red{k}{\myr}{s}{F} 0$  yields 
       the existence of $g \in \myk[X_1, \ldots,X_n]/\ideal{}{}(P_T)$ such that
       $p  \red{*}{\myr}{}{F} p' \red{*}{\myr}{}{F} g$ and
       $q  \red{*}{\myr}{}{F} q' \red{*}{\myr}{}{F} g$.
\\
\qed
\end{enumerate}\renewcommand{\baselinestretch}{1}\small\normalsize
\begin{corollary}~\\
{\sl
If for all polynomials $g \in \ideal{}{}(F)$ we have
 $g \red{*}{\myr}{}{F} 0$,
 then $F$ is a Gr\"obner basis.
\ohnebeweis
}
\end{corollary}
\procedure{Completion in Commutative Quotient Rings}%
{\vspace{-4mm}\begin{tabbing}
XXXXX\=XXXX \kill
\removelastskip
{\bf Given:} \> A finite set of polynomials 
                $F \subseteq \myk[X_1, \ldots, X_n]/\ideal{}{}(P_R)$, \\
 \>               and $R$ a set of relations presenting $\m$. \\
{\bf Find:} \> $\gb(F)$, a  Gr\"obner basis of $F$ in
               $\myk[X_1, \ldots, X_n]/\ideal{}{}(P_R)$.
\end{tabbing}
\vspace{-7mm}
\begin{tabbing}
XX\=XX\=XXXX\= XX\= XXXX \=\kill
$P_T$ := Gr\"obner$(P_R)$; \hspace{1cm}{\rm\kommentar \%  compute
  a reduced Gr\"obner basis in {\bf K}$[X_1, \ldots, X_n]$} \\
$G$ := $F$; \\
$H$ := $\{ (q, l-r) \mid q \in G, l-r \in P_T \}$; \\
$B$ := $\{ (q_{1}, q_{2}) \mid q_{1}, q_{2} \in G, q_{1} \neq q_{2} \}$; \\
{\bf while} $B \cup H \neq \emptyset$ {\bf do} \\
\>{\bf if} \> $B \neq \emptyset$ \\
\>         \>{\bf then} \>$(q_{1}, q_{2})$ := remove$(B)$; \\
\>         \>           \>  $h$ := ${\rm normalform}(\spol{}(q_{1}, q_{2}),\red{}{\myr}{b}{G \cup P_T})$ \\
\>         \>           \>{\rm\kommentar \%  compute a normal form
              with respect to Buchberger's reduction}\\
\>         \>           \>  {\bf if} \>$h \neq 0$\\
\>         \>           \>           \>{\bf then} \>     $B$ := $B \cup \{ (f,h) \mid f \in G \}$; \\
\>         \>           \>           \>           \>$G$ := $G \cup \{ h \}$; \\
\>         \>           \>           \>           \>$H$ := $H \cup \{ (h, l-r) \mid l-r \in P_T \}$ \\
\>         \>{\bf else} \> $(q, l-r)$ := {\rm remove}$(H)$; \\
\>         \>           \>  $h$ := ${\rm normalform}(\spol{}(q, l-r),\red{}{\myr}{b}{G \cup P_T})$ \\
\>         \>           \>  {\bf if} \> $h \neq 0$ \\
\>         \>           \>           \>{\bf then} \>     $B$ := $B \cup \{ (f,h) \mid f \in G \}$; \\
\>         \>           \>           \>           \>$G$ := $G \cup \{ h \}$; \\
\>         \>           \>           \>           \>$H$ := $H \cup \{ (h, l-r) \mid l-r \in P_T \}$ \\
{\bf endwhile} \\
$\gb (F)$:= $G$
\end{tabbing}}
} 
We will close this section by comparing the reduction introduced
 here for a quotient structure to right and commutative reduction
 in the corresponding commutative monoid ring.
For a polynomial $p \in \myk[X_1, \ldots, X_n]/\ideal{}{}(P_{T_c})$ let
$\tilde{p}$ be the corresponding polynomial in the monoid ring $\myk[\m]$.
%
\begin{lemma}~\\
{\sl
Let $p,q,f$ be some polynomials in $\myk[X_1, \ldots, X_n]/\ideal{}{}(P_{T_c})$
 and let $\tilde{p},\tilde{q},\tilde{f}$ be the corresponding polynomials in $\myk[\m]$.
Then $p \red{}{\myr}{}{f} q$ if and only if $\tilde{p} \red{}{\myr}{c}{\tilde{f}} \tilde{q}$.
\lemend
}
\end{lemma}
\Ba{}~\\
Before entering the proof of our claim let us first take a closer look at
reduction in $\myk[X_1, \ldots, X_n]/\ideal{}{}(P_{T_c})$.
If $p \red{}{\myr}{}{f} q$ at a monomial $\alpha \skm t$ with
 $t = \hterm(f) \cm  u$, then we can express this reduction step by
 $p \red{}{\myr}{b}{f} q' \red{*}{\myr}{b}{P_{T_c}} q$ and we have
 $q' = p - \alpha \skm \hc(f)^{-1} \skm f \mrm_{\myk[X_1, \ldots, X_n]} u$ 
 and
 $q = p - \alpha \skm \hc(f)^{-1} \skm f \qmult u$.
In this context it is easy to see that $p \red{}{\myr}{}{f} q$
 implies $\tilde{p} \red{}{\myr}{c}{\tilde{f}} \tilde{q}$.
On the other hand, commutative reduction requires that the head term
 of the polynomial is a divisor with respect to $\mm_{\freecomm}$ of
 the term to be reduced, i.e., for $T_c$ the term corresponding
 to $t$ in $\tilde{p}$ we get $T_c = \hterm(\tilde{f}) \cm
 \tilde{u}$ and $\tilde{q} = \tilde{p} - \alpha \skm
 \hc(\tilde{f})^{-1} \skm \tilde{f} \mrm_{\myk[\m]} \tilde{u}$.
Hence, a commutative reduction step in our monoid ring can be split
 into first doing one step using Buchberger's reduction and then
 normalizing the new monomials using $P_{T_c}$.
\\
\qed
\begin{corollary}~\\
{\sl
Let $G \subseteq \myk[X_1, \ldots, X_n]/\ideal{}{}(P_{T_c})$ and let $\tilde{G}$ be
 the corresponding set in $\myk[\m]$.
Then $G$ is a Gr\"obner basis with respect to $\R$ if and only if
 $\tilde{G}$ is a Gr\"obner basis with respect to $\red{}{\myr}{c}{}$.
\ohnebeweis
}
\end{corollary}
Obviously, then $\R$ must be weaker than $\red{}{\myr}{r}{}$.
\begin{corollary}~\\
{\sl
Let $p,q,f$ be some polynomials in $\myk[X_1, \ldots, X_n]/\ideal{}{}(P_{T_c})$
 and let $\tilde{p},\tilde{q},\tilde{f}$ be the corresponding polynomials in $\myk[\m]$.
Then $p \red{}{\myr}{}{f} q$ implies 
 $\tilde{p} \red{}{\myr}{r}{\tilde{f}} \tilde{q}$ but not
 vice versa.
}
\end{corollary}
\begin{example}~\\
{\rm
Let $\Sigma = \{ X_1,X_2 \}$ and $T_c = \{ X_1^2 \myr \lambda,
 X_2^2 \myr \lambda \}$
 be a presentation of a commutative group $\g$ with  a
 length-lexicographical ordering induced by $X_1 \succ X_2$.
Further let $\tilde{p} = X_2$ and $\tilde{f} = X_1X_2$ be polynomials
in $\myk[\g]$.
\\
Then $\tilde{p} \red{}{\myr}{r}{\tilde{f}} 0$, but for the corresponding polynomials
 $p,f$ in  $\myk[X_1, \ldots, X_n]/\ideal{}{}(P_{T_c})$ we have
 $p \nred{}{\myr}{}{f}$, as $X_1 \nred{}{\myr}{b}{X_1X_2}$.
Note that the set $\{ \tilde{f} \}$ itself is a Gr\"obner basis with respect
 to $\red{}{\myr}{r}{}$ but not with respect to
 $\red{}{\myr}{c}{}$, nor is $\{ f \}$ a Gr\"obner basis with respect to $\R$.
\exaend
}
\end{example}
%

%% file: specialgroups.tex
\chapter{Group Rings}\label{chapter.grouprings}
\spruch{10}{7}{Longum iter est per praecepta,\\
Breve et efficax per exempla.}{Seneca}

In this chapter we want to apply the ideas of reduction developed in
the previous chapter to group
rings.
In groups many problems have easy solutions due to the existence of
inverses, e.g., the solvability of equations. 
Additional information provided by the presentations for different classes
 of groups is incorporated to give improved and terminating
 procedures to compute Gr\"obner bases for finitely generated right
 ideals.

{\bf Section 5.1:} Similar to the equivalence of certain restricted versions of
  the word problem for semi-Thue systems to restricted versions of the ideal
  congruence problem for free monoid respectively group rings in
  section \ref{section.undecidable}, we show that the subgroup problem
  is equivalent to a restricted version of the right ideal membership problem in a
  group ring.
  Thus only groups having solvable subgroup problem are candidates for
  allowing the computation of finite right Gr\"obner bases. 

{\bf Section 5.2:} The results on prefix reduction can be used to
  give a terminating procedure to compute finite reduced prefix
  Gr\"obner bases for finitely generated right ideals in free group
  rings.

{\bf Section 5.3:} Similar ideas as in the case of free groups can be
  carried over to the class of plain groups and a procedure is
  provided to compute finite reduced prefix
  Gr\"obner bases for finitely generated right ideals.

{\bf Section 5.4:} The results on free groups are combined with
  special presentations of context-free groups to give a procedure to 
  compute finite reduced prefix
  Gr\"obner bases for finitely generated right ideals.

{\bf Section 5.5:} The ideas of commutative reduction are generalized
  to the case of nilpotent groups resulting in the definition of
  quasi-commutative reduction.
  We give a terminating procedure to compute Gr\"obner bases for
  finitely generated right ideals in torsion-free nilpotent group
  rings and show how these ideas can be generalized for nilpotent
  group rings in a construction similar to context-free groups.

\section{The Subgroup Problem}
In section \ref{section.undecidable} we have shown that the word
problem for group presentations is equivalent to a  restricted version of the
ideal membership problem for a free group ring.
We will now show that a similar equivalence holds for the right ideal
membership problem in group rings.

\begin{definition}~\\
{\rm
Given a subset $S$ of a group $\g$ let $\left< S \right>$ denote
 the subgroup generated by $S$. 
The \index{generalized word problem}\index{word problem!generalized}\betonen{generalized word problem}\/ or
\index{subgroup problem}\betonen{subgroup problem}\/ is then to
 determine, given an element $w \in \g$, whether $w \in \left< S \right>$. 
\dend
} 
\end{definition}
The word problem for a group $\g$ is just the generalized word problem
for the trivial subgroup in $\g$.
Thus the existence of a group with undecidable word problem yields
undecidability for the subgroup problem.
On the other hand, decidable word problem for a subgroup does not
imply decidable generalized word problem.

The next theorem states that the subgroup problem for a group is
equivalent to a special instance of the right membership problem in
the corresponding group ring.
\begin{theorem}\label{theo.subgroup.problem}~\\
{\sl
Let $S$ be a finite subset of $\g$ and $\myk[\g]$ the  group ring
corresponding to $\g$. 
Further let $P_S
 = \{ s - 1 \mid s \in S \}$ be a set of polynomials associated to $S$\footnote{Note that we use
                                 $1 = 1 \skm \lambda = \lambda$.}. 
Then the following statements are equivalent:
\begin{enumerate}
\item $w \in  \left< S \right>$. 
\item $w-1 \in \ideal{r}{}(P_S)$.
\end{enumerate}
}
\end{theorem}
\Ba{}~\\
\mbox{$1 \R 2:$ }
      Let $w = u_1\mm \ldots\mm u_k \in \left< S \right>$, i.e., $u_1,
      \ldots, u_k \in S \cup \{ \inv{s} | s \in S \}$.
      We show $w-1 \in \ideal{r}{}(P_S)$ by induction on $k$.
      In the base case $k=0$ there is nothing to show, as
       $w=\lambda \in  \left< S \right>$ and $0 \in \ideal{r}{}(P_S)$.
      Hence, suppose  $w = u_1\mm \ldots\mm u_{k+1}$ and $
       u_1\mm \ldots\mm u_k -1 \in \ideal{r}{}(P_S)$.
      Then
      $(u_1\mm\ldots\mm u_k-1)\mrm u_{k+1} \in
      \ideal{r}{}(P_S)$ and, since 
      $ u_{k+1}-1 \in \ideal{r}{}(P_S)$\footnote{We either have
        $u_{k+1}-1 \in P_S$ or $\inv{u_{k+1}} \in S$, i.e., 
        $(\inv{u_{k+1}} -1) \mrm u_{k+1} =  u_{k+1}-1 \in \ideal{}{}(P_S)$.}, we get
      $( u_1\mm \ldots\mm u_k-1) \mrm u_{k+1} +
      (u_{k+1} - 1) = w-1 \in \ideal{r}{}(P_S)$.

\mbox{$2 \R 1:$ }
       We have to show that $w-1 \in \ideal{r}{}(P_T)$ implies
        $w \in  \left< S \right>$.
       We know $w-1 = \sum_{j=1}^{n} \alpha_j \skm (u_j -1) \mrm x_j$, where
        $\alpha_j \in \myk^*$, $u_j \in S \cup \{ \inv{s} | s \in S \}$,
        $x_j \in \g$. 
       Therefore, by  showing the following stronger result we are
       done:
       A representation
        $ w - 1 = \sum_{j=1}^m p_j$ where $p_j = \alpha_j \skm (w_j -w'_j)$,
        $\alpha_j \in \myk^*$,$w_j \neq w'_j$ and $w_j \mm \inv{w'_j}
        \in \left< S \right>$ implies 
        $w \in \left< S \right>$.
       Now, let $w - 1 = \sum_{j=1}^m p_j$ be such a representation and $\succeq$ be an arbitrary total well-founded ordering
       on $\g$. 
       Depending on this  representation 
        and $\succeq$ we define
        $t = \max \{w_j,w'_j   \mid j = 1, \ldots m  \}$  and
        $K$ is the number of polynomials $p_j$ containing $t$ as a term.
       We will show our claim by induction on $(m,K)$, where 
       $(m',K') < (m,K)$ if and only if $m'<m$ or $(m'=m$ and $K'<K)$.
       In case $m=0$, $w-1=0$ implies $w=1$ and hence $w \in  \left< S \right>$.
       Thus let us assume $m>0$.
\\
       In case $K=1$, let $p_k$ be the polynomial containing
        $t$.
       As we either have $p_k=\alpha_k \skm (t-w'_k)$ or
        $p_k=\alpha_k \skm (w_k - t)$, where $\alpha_k \in \{ 1, -1 \}$,
        without loss of generality we can assume $p_k=t-w'_k$.
       Using $p_k$ we can decrease $m$ by subtracting $p_k$
        from $w-1$ giving us
        $w'_k-1 = \sum_{j=1,j \neq k}^{m} p_j$.
       Since $t \mm   \inv{w'_k} \in \left< S \right>$
        and our induction hypothesis yields
        $w'_k \in \left< S \right>$,
        we can conclude $w=t=(t \mm  \inv{w'_k}) \mm  w'_k \in \left< S \right>$.
\\
       In case $K>1$ there are two polynomials $p_k,p_l$ in the
       corresponding representation 
        and without loss
        of generality we can assume $p_k = \alpha_k \skm( t - w'_k)$
        and $p_l = \alpha_l \skm (t -w'_l)$. 
       If then $w'_k = w'_l$  we can immediately decrease $m$ 
        by substituting the occurrence of $p_k+p_l$  by $(\alpha_k + \alpha_l) \skm p_l$.
       Otherwise we can proceed as follows:
       \begin{eqnarray*}
         p_k + p_l & = & p_k \underbrace{-  \alpha_k \skm \alpha_l^{-1} \skm p_l + \alpha_k  \skm \alpha_l^{-1} \skm p_l}_{=0} + p_l \\
                   & = & \underbrace{(-\alpha_k \skm w'_k + \alpha_k  \skm w'_l)}_{p'_k} + (\alpha_k \skm \alpha_l^{-1} +1) \skm p_l \\
       \end{eqnarray*}
       where $p'_k = \alpha_k \skm (w'_l - w'_k)$, $w'_k \neq w'_l$
       and $w'_k \mm  \inv{w'_l} \in\left< S \right>$, since
       $w'_k \mm  \inv{t},t \mm  \inv{w'_l} \in\left< S \right>$ and 
       $w'_k \mm  \inv{w'_l}=w'_k \mm  \inv{t} \mm  t \mm  \inv{w'_l}$.
      In case $\alpha_k \skm \alpha_l^{-1} + 1=0$, i.e., $\alpha_k=-\alpha_l$, $m$ is decreased.
      On the other hand $p'_k$ does not contain $t$, i.e., if $m$ is
      not decreased $K$ is.
\\
\qed
This theorem implies that we can only expect group rings over groups with solvable
generalized word problem to allow solvable membership problem for
right ideals.
On the other hand, solvable subgroup problem only implies the
solvability of a  restricted version of the right ideal membership problem.

The usage of right ideals corresponds to the fact that the set $S
\subseteq \g$ induces a left congruence, namely  $u \sim_{S} v$ if and
only if $\left< S \right> u =_{\g} \left< S \right> v$.
Different methods to express this left congruence by reduction
methods in order to solve the subgroup problem can be found in the literature.
For free groups there is Nielsen's approach  known as Nielsen
reduction (compare \cite{LySch77,AvMa84}).
Kuhn and Madlener have developed prefix reduction methods and applied
them successfully to the
class of plain groups (see \cite{KuMa89}).
Wi{\ss}mann solved the subgroup problem for  the class of polycyclic groups (compare
\cite{Wi84,Wi89}) and Cremanns and
Otto successfully treated the class of context-free groups (see
\cite{CrOt94}). 

We move on now to study the right ideal membership problem in special
classes of groups.
\section{Free Groups}
In group theory a particularly important role is played by groups that
 are free  in the class of all groups which themselves are rather
 simple groups.
In this section we state  how the ideas of prefix reduction can be applied
 to give a completion algorithm for finitely generated free group rings.

Let $\free$ be a free group generated by a finite set $X = \{ x_1,
 \ldots, x_n \}$.
Then  $\Sigma = X \cup X^{-1}$ and
 $T = \{ xx^{-1} \myr \lambda, x^{-1}x \myr \lambda \mid x \in X \}$ is a
 presentation of $\free$ with
 $x_1^{-1} \succ x_1 \succ \ldots \succ x_n^{-1} \succ x_n$ inducing a
 length-lexicographical ordering on $\free$.
Note that $(\Sigma, T)$ then  is   a convergent 2-monadic monoid
presentation of $\free$.
We call $x^{-1}$ the formal inverse of $x \in X$ and we will allow the
 following notations for such $x \in X$: $\inv{x} = x^ {-1}$ and
 $\inv{x^{-1}} = x$.
This can be extended to $\Sigma^*$ by setting $\inv{\lambda} = \lambda$,
 $\inv{wx}= x^{-1}\inv{w}$ and $\inv{wx^{-1}} = x \inv{w}$ for 
 $wx, wx^{-1} \in \Sigma^*$.

Let us start with some technical notions for
 polynomials in $\myk[\free]$ which will allow an immediate
 characterization of saturating sets for polynomials.
Note that the ideas used in the following definition for special
instances can be compared
to isolating prefixes as it is done in Nielsen's approach to solve the
generalized word problem in free groups.
\begin{definition}\label{def.freegroupsat}~\\
{\rm
For a polynomial $p \in \myk[\free]$ which has more than one monomial,
 we define
\begin{eqnarray}
\sigma_1(p) & = & \max \{ u \in \free \mid \inv{u} \mbox{ is a suffix of }
                       \hterm(p) \mbox{ and } \hterm(p \mrm u) = \hterm(p) \mm u
                       \}, \nonumber \\
\sigma_2(p) & = & \min \{ u \in \free \mid \inv{u} \mbox{ is a suffix of }
                       \hterm(p) \mbox{ and } \hterm(p \mrm u) \neq \hterm(p) \mm u
                       \}. \nonumber
\end{eqnarray}
Then we can set $\can(p)= p \mrm \sigma_1(p)$ and
$\satpoly(p) = p \mrm \sigma_2(p)$.
For a polynomial $\alpha \skm t \in \myk[\free]$ we set $\sigma_1(p) =\sigma_2(p)=
 \inv{t}$ and $\can(p) = \satpoly(p) = \lambda$.
\dend
}
\end{definition}
The polynomials $\can(p)$ and $\satpoly(p)$ will often be called
``mates'' of each other.
Note that $\sigma_1(p)$ is a prefix of $\sigma_2(p)$ and in case $p$
contains more than one monomial we have $\sigma_2(p) \id \sigma_1(p)a$
 where $a = \inv{\ell(\hterm(p \mrm \sigma_1(p)))}$, i.e., $a$ is the
 inverse of the last letter of the head term of the polynomial $\can(p)$.
Hence $|\sigma_2(p)| = |\sigma_1(p)| +1$ holds.
\begin{example}~\\
{\rm
Let $\Sigma = \{ x,x^{-1} \}$ and $\free$ the free group generated by
$x$.
\\
Then for the polynomial $p= x^4 + x^2 + \lambda \in \q[\free]$ we get
 $\sigma_1(p) = x^{-1}$, $\sigma_2(p) = x^{-2}$, $\can(p) = p \mrm \sigma_1(p) =
 \underline{x^3} + x + x^{-1}$, and $\satpoly(p) = p \mrm \sigma_2(p) = x^{2} +
 \lambda + \underline{x^{-2}}$.
\exaend
}
\end{example}
Notice that $\hterm(p \mrm \sigma_1(p))$ is a prefix of $\hterm(p)$ and
  hence    $p \mrm \sigma_1(p) \leq p$.
Furthermore, in case $p \neq \alpha \skm t$, we get $p = \can(p) \mrm \inv{\sigma_1(p)} = \satpoly(p) \mrm
\inv{\sigma_2(p)}$ yielding $\ideal{r}{}(p) = \ideal{r}{}(\can(p)) =
  \ideal{r}{}(\satpoly(p))$.

Next let us take a closer look at the special forms
$\can(p)$ and $\satpoly(p)$ of a
polynomial $p$ and their head terms respectively the last letters of
their head terms.\label{def.largerterm}
Let us associate a pair of terms $(t_1^p, t_2^p)$ to $p$ such that $t_1^p = \min \{ \hterm(\can(p)), \hterm(\satpoly(p)) \}$, $t_2^p = \max \{ \hterm(\can(p)),
\hterm(\satpoly(p)) \}$ and set 
$q_i \in \{ \can(p), \satpoly(p) \}$ such that $\hterm(q_i) = t_i^p$.
Then the following lemma holds.
\begin{lemma}\label{lem.lastletter}~\\
{\sl
Let $p$ be a non-zero  polynomial in $\myk[\free]$ with more than one
monomial, and $(t_1^p, t_2^p)$,
$q_1, q_2$ as described above.
Then the following statements hold:
\begin{enumerate}
\item $|t_2^p| - |t_1^p| \leq 1$.
\item For all terms $t' \in \terms(q_2)$ with $|t'| = 
       |t_2^p|$ we have
       $\ell(t') =
       \ell(t_2^p) = \inv{\ell(t_1^p)}$.
\item $q_1 \mrm \inv{\ell(t_1^p)} = q_2$ and
      $q_2 \mrm \inv{\ell(t_2^p)} = q_1$.
\end{enumerate}
\lemend
}
\end{lemma}
\Ba{}~\\
Let $\hterm(\can(p)) = \hterm(p \mrm \sigma_1(p)) = t$ and as
$\sigma_2(p) \id \sigma_1(p)a$ for some $a \in \Sigma$ we get
$\hterm(\satpoly(p)) = \hterm(p \mrm \sigma_2(p)) = \hterm(p \mrm
\sigma_1(p) \mrm a) = s \mm a \succ t \mm a$ for some $s \in
\terms(\can(p))$.
Then $s \prec t$ and $s \mm a \succ t \mm a$ implies $s \mm a \id sa$
and $|t \mm a| < |t|$.
Hence either $|sa| = |t|$ or $|sa| = |t|+1$, and hence, as $t_1^p,
t_2^p \in \{ sa, t \}$ we get $|t_2^p| - |t_1^p| \leq 1$.
\\
To see that for all $t' \in \terms(q_2)$ with $|t'| = |t_2^p|$ the
last letters coincide with $\ell(t_2^p)$, let us take a closer look at
the terms $t_1^p$ and $t_2^p$ respectively the terms in $q_1$ and $q_2$.
We know $t_1^p \preceq t_2^p \id ub$ for some $u \in \free$, $b \in
\Sigma$ and $t_2^p = v \mm c$ for some $v
\in \terms(q_1)$, $c \in \Sigma$, i.e., $v \preceq t_1^p$, but $ub \id
v \mm c \succ t_1^p \mm c$.
Hence for the last letter of $t_1^p$, $\ell(t_1^p) = \inv{c}$ must
hold.
On the other hand, $\ell(v) = \inv{c}$ is not possible, implying $ub
\id vc$ and in particular $b = c$ and $u \in \terms(q_1)$.
It remains to study those terms $v'\in \terms(q_1)$ with $|v' \mm b| =
|t_2^p|$.
In distinguishing the four possible cases we  find that for
the cases $|t_2^p| = |t_1^p| > |v'|$ and $|t_2^p| > |t_1^p| = |v'|$ we have $v' \mm b \id v'b$,
and for the  cases $|t_2^p| = |t_1^p| = |v'|$ and $|t_2^p| > |t_1^p| >
|v'|$ either $|v' \mm b| = |t_2^p| +1$ respectively
 $|v' \mm b| = |t_2^p| -1$ or $|v' \mm b| \leq |t_2^p| -1$
gives us a contradiction to $|v' \mm b| =
|t_2^p|$.
\\
Finally, since $\ell(t_1^p) = \inv{b}$ and $\ell(t_2^p) = b$,
 $q_1 \mm b = q_2$ and $q_2 \mrm \inv{b} = q_1$ follows immediately.
Furthermore, since $t_1^p, t_2^p \in \{ sa, t \}$, this implies $b \in
\{ a, \inv{a} \}$.
\\
\qed
\begin{corollary}\label{cor.ht.no.prefix}~\\
{\sl
Let $p$ be a polynomial in $\myk[\free]$ containing more than one
monomial.
Then neither $\hterm(\can(p))$ nor $\hterm(\satpoly(p))$ are prefixes
of one another.
In particular  we get
 $\can(p) = \satpoly(p)$ if and only if $p = \alpha \skm t$ for some
 $\alpha \in \myk^*$, $t \in \free$.
}
\end{corollary}
\Ba{}~\\
Let $\hterm(\can(p)) \id ta$ for some $a \in \Sigma$ and by
 lemma \ref{lem.lastletter} $\hterm(\satpoly(p)) \id s\inv{a}$
 for some $s \in \terms(\can(p))$.
Now suppose $ta$ is a prefix of $s\inv{a}$.
Then, as $||\hterm(\satpoly(p))| - |\hterm(\can(p))|| \leq 1$ and $a
\neq \inv{a}$, we get $s
\id ta$ contradicting the fact that $s \mm \inv{a} \id s\inv{a}$.
Likewise, if $s\inv{a}$ were a prefix of $ta$ we would get $t \id
s\inv{a}$ contradicting that $ta \in \free$.
\\
In particular $\can(p) = \satpoly(p)$ implies $\hterm(\can(p)) =
\hterm(\satpoly(p))$, and this is only possible in case  $p = \alpha \skm t$ for some
 $\alpha \in \myk^*$, $t \in \free$ and $\can(p) = \satpoly(p) = \lambda$.
\\
\qed
\begin{lemma}~\\
{\sl
Let $p$ be a polynomial in $\myk[\free]$ containing more than one
monomial.
Then the following statements hold:
\begin{enumerate}
\item $\can(\can(p)) = \can(p)$.
\item $\satpoly(\can(p)) = \satpoly(p)$.
\item $\can(\satpoly(p)) = \satpoly(p)$.
\item $\satpoly(\satpoly(p)) = \can(p)$.
\end{enumerate}
\lemend
}
\end{lemma}
\Ba{}~\\
Note that $\sigma_1(\can(p)) = \lambda$ and $\sigma_2(\can(p)) =
\inv{\ell(\hterm(\can(p)))} = \ell(\hterm(\satpoly(p)))$.
Therefore, we get $\can(\can(p)) = \can(p) \mrm \sigma_1(\can(p)) = \can(p)$ and 
 $\satpoly(\can(p) = \can(p) \mrm \sigma_2(\can(p)) = \satpoly(p)$.
On the other hand we find $\sigma_1(\satpoly(p)) = \lambda$ and
$\sigma_2(\satpoly(p))=\inv{\ell(\hterm(\satpoly(p)))} =
\ell(\hterm(\can(p)))$.
Thus, $\can(\satpoly(p)) = \satpoly(p) \mrm \sigma_1(\satpoly(p)) = \satpoly(p)$ and 
  $\satpoly(\satpoly(p)) = \satpoly(p) \mrm \sigma_2(\satpoly(p)) = \can(p)$.
\\
\qed
We can specify  prefix saturating sets for polynomials
 in $\myk[\free]$ in terms of  the polynomials defined in definition
 \ref{def.freegroupsat}.
\begin{lemma}\label{lem.freegroupsat}~\\
{\sl
If a  polynomial $p  \in \myk[\free]$ contains more than one monomial, the set
$\{ \can(p), \satpoly(p) \}$ is a prefix saturating set for $p$.
In particular, we find $\SAT_p(p) = \SAT_p(\can(p))= \SAT_p(\satpoly(p))$.
\lemend
}
\end{lemma}
\Ba{}~\\
\auskommentieren{We have to show that for all elements $w \in \free$, 
 $p \mrm w \red{}{\myr}{p}{\{ \can(p), \satpoly(p) \}} 0$.
\\
In case $\hterm(p)w$ is not $T$-reducible, we immediately get
 $p \mrm w \red{}{\myr}{p}{p \mrm \sigma_1(p)} 0$.
\\
Hence let $\hterm(p)w \id \hterm(p)w_1w_2$ be a decomposition in such a way that we have
 $w_1 = \max \{ u \in \free | \inv{u} \mbox{ is a suffix of }
                       \hterm(p) \mbox{ and } u \mbox{ is a prefix of }
                       w \}.$
\\
Then the case $w_1 \preceq \sigma_1(p)$ yields 
 $p \mrm w \red{}{\myr}{p}{p \mrm \sigma_1(p)} 0$.
\\
Otherwise let us assume $w_1 \succ \sigma_1(p)$, i.e.,
 $w_1 \id \sigma_1(p)az \id \sigma_2(p)z$ for some $a \in \Sigma$, $z \in \free$.
Note that $p \mrm \sigma_1(p)a = p \mrm \sigma_2(p)$ and that
 the head term of $p \mrm \sigma_2(p)$ ends with the letter $a$.
Further, $zw_2$ does not start with $\inv{a}$ and thus
 $\hterm(p \mrm \sigma_2(p))\mm zw_2 \id \hterm(p \mrm \sigma_2(p))zw_2$.
This implies $p \mrm w = p \mrm \sigma_2(p) \mrm zw_2 \red{}{\myr}{p}{p
  \mrm \sigma_2(p)} 0$.}
We can use procedure {\sc Prefix Saturated Check} on page
\ref{is.prefix.saturated} to show that the set $\{ \can(p), \satpoly(p) \}$ is
a prefix saturating set for a polynomial $p$.
For the polynomials $\can(p)$ and $\satpoly(p)$ we get the
corresponding sets $C(\hterm(\can(p))) = \{
\inv{\ell(\hterm(\can(p)))} \}$ respectively $C(\hterm(\satpoly(p))) = \{
\inv{\ell(\hterm(\satpoly(p)))} \}$.
Now by lemma \ref{lem.lastletter} we know $\can(p) \mrm
\inv{\ell(\hterm(\can(p)))} = \satpoly(p)$ and
$\satpoly(p) \mrm \inv{\ell(\hterm(\satpoly(p)))} = \can(p)$ and hence
the set $\{ \can(p), \satpoly(p) \}$ is prefix saturated.
Furthermore, as it is a subset of $\{ p \mrm w \mid w \in \free \}$ it
is also a prefix saturating set for $p$.
\\
\qed

\begin{remark}~\\
{\rm
For a non-zero polynomial $p$ in $\myk[\free]$ the set
 $\{ \can(p), \satpoly(p) \}$ even is a prefix Gr\"obner basis of the
 right ideal generated by $p$.
In case $p = \alpha \skm t \in \myk[\free]$ this is true as the set $\{
\lambda \}$ is a prefix Gr\"obner basis for $\ideal{r}{}(p) =
\myk[\free]$.
If $p$ contains more than one monomial, 
 the head terms of the polynomials $\can(p)$ and $\satpoly(p)$
 are no prefixes of each other and this set is prefix saturated.
As $\ideal{r}{}(p) = \ideal{r}{}(\can(p))=\ideal{r}{}(\satpoly(p))$, theorem \ref{theo.pcp}
 implies that $\{ \can(p), \satpoly(p) \}$ is a prefix Gr\"obner basis.
\remend
}
\end{remark}
Now we can give a completion procedure for $\myk[\free]$ by
 modifying procedure {\sc Reduced Prefix Gr\"obner Bases} 
(see page \pageref{reduced.prefix.groebner.bases}) in specifying the
saturating procedure for polynomials.
It remains to show that finite reduced prefix Gr\"obner basis exist,
as then this procedure will compute them.
This will be done for the more general case of plain groups in the next
section.
In section \ref{section.prefixreduction} we mentioned that the sets
$S_i$ in procedure {\sc Reduced Prefix Gr\"obner Bases}  in general
contain unnecessary polynomials.
Therefore, we next provide a procedure using additional information to
prevent this overhead.
The idea is to use  prefix
 reduction at head terms (this is comparable to the step of computing
 prefix s-polynomials) combined with saturating polynomials and
 to remove not only {\em one} polynomial from the set $G_i$ for reduction, but
 also its corresponding ``mate'', i.e., we remove not only a polynomial
 $q$ from $G_i$ but the set $\{ \can(q), \satpoly(q) \}$\footnote{Note that
  by construction we will have $q \in \{ \can(q), \satpoly(q) \}$.}.
Hence the sets $G_i$ contain at most $2 \skm |F|$ polynomials.
This algorithm can be compared to the results of Rosenmann's approach
to free group rings in \cite{Ro93}.
The procedure will use the definitions of the special terms $t_1^q$
and $t_2^q$ introduced for polynomials $q$ on page
\pageref{def.largerterm}.
Given two polynomials $q,q'$ we define the following tuple ordering on
the accompanying terms which is well-founded:
$(t_1^q, t_2^q)
 \succ (t_1^{q'}, t_2^{q'})$ if and only if $t_1^q \succ t_1^{q'}$ or
 $(t_1^q = t_1^{q'}$ and $t_2^q \succ t_2^{q'})$.

\procedure{Completion in Free Group Rings\protect{\label{free.group.rings}}}%
{\vspace{-4mm}\begin{tabbing}
XXXXX\=XXXX \kill
\removelastskip
{\bf Given:} \> A finite set $F \subseteq \myk[\free]$.  \\
{\bf Find:} \> $\gb(F)$, a  prefix Gr\"obner basis of $F$.
\end{tabbing}
\vspace{-7mm}
\begin{tabbing}
XX\=XX\=XXXX\= XX \=XXXX\=\kill
$i$ := $0$; \\
$G_0$ := $\{ \can(f), \satpoly(f) \mid f \in F\}$; \\
{\rm\kommentar \%
 $\ideal{r}{}(F) = \ideal{r}{}(G_0)$ and $G_0$ is prefix saturated}\\
{\bf while} there exists  $q \in G_i$ such that
$\hterm(q)$ is prefix reducible by  $g' \in  G_i\backslash
\{ q \}$ with \\
\>$(t_1^{g'}, t_2^{g'}) \preceq  (t_1^q, t_2^q)$ {\bf do}\\
\> $i$ := $i+1$; \\
\> $G_i$ := $G_{i-1} \backslash \{ \can(q), \satpoly(q) \}$; \\
\> $G'$ := $\{ g \in G_i \mid (t_1^{g}, t_2^{g}) \preceq  (t_1^q, t_2^q) \}$; \\
\> $q'$ := ${\rm headnormalform}(q,  \red{}{\myr}{p}{G'})$; \\
\> {\rm\kommentar \% Compute a normalform of a polynomial allowing only
  reduction steps at the} \\
\> {\rm\kommentar \% respective head terms} \\
\> {\bf if} \>$|\terms(q')| = 1$ {\rm\kommentar \% The right ideal
  generated by $F$ is trivial} \\
\>          \> {\bf then} \>       $G_i$ := $\{ \lambda \}$;\\
\>          \>{\bf else} \> {\bf if } \>$q' \neq 0$ \\
\>          \>           \>           \>{\bf then } \>$G_i$ := $G_{i} \cup \{ \can(q'), \satpoly(q')
\}$; \\ 
\>          \>           \>            \>           \>{\rm\kommentar \%
  $\ideal{r}{}(F) = \ideal{r}{}(G_i)$ and $G_i$ is prefix saturated}\\
\>          \> \>  {\bf endif } \\
\> {\bf endif}\\
{\bf endwhile} \\
$G$ := $G_i$
\end{tabbing}}

Notice that we always have $\ideal{r}{}(F) = \ideal{r}{}(G_i)$ and
  for the prefix reduced set
 $\tilde{G} =\{ {\rm normalform}(g, G \backslash \{ g \})
 \mid g \in G \}$ we again have $\ideal{r}{}(G) = \ideal{r}{}(\tilde{G})$.
Furthermore,
$\hterm(G) = \hterm(\tilde{G})$ holds since by construction no term in
$\hterm(G)$ is prefix of another term in $\hterm(G)$.
Thus $\tilde{G}$ is a reduced prefix
Gr\"obner basis of $\ideal{r}{}(G)$.

\begin{theorem}\label{theo.correct.free.group}~\\
{\sl
Procedure {\sc Completion in Free Group Rings} is totally correct.
}
\end{theorem}
\Ba{}~\\
In case the procedure terminates, correctness follows at once from the
 fact that the final set $G_k$ is a prefix Gr\"obner basis of
 $\ideal{r}{}(F)=\ideal{r}{}(G_k)$, as it is
 prefix saturated by construction and no prefix s-polynomials exist.
To see the latter, let us assume that although no polynomial $q \in
G_k$ exists such that
$\hterm(q)$ is prefix reducible by a polynomial $g' \in  G_k\backslash
\{ q \}$ with $(t_1^{g'}, t_2^{g'}) \preceq  (t_1^q, t_2^q)$, 
 there exist $g, g' \in G_k$ such that $\hterm(g) \id
 \hterm(g')u$ for some $u \in \free$.
Then $(t_1^{g'}, t_2^{g'}) \succ  (t_1^g, t_2^g)$, i.e., $t_1^{g'}
\succ t_1^g$ or $(t_1^{g'} = t_1^g$ and $t_2^{g'} \succ t_2^g)$, must
hold.
The case $u = \lambda$ is not possible as then  the {\bf while loop} would have to be executed using $q = g'$
contradicting our assumption.
Now let us distinguish the possible cases.
If $\hterm(g) = t_1^g$ we get $t_1^g = \hterm(g) \succ
\hterm(g') \succeq t_1^{g'}$ contradicting the fact that $t_1^{g'}
\succeq t_1^g$  holds.
It remains to look at the case $\hterm(g) = t_2^g$.
Now, if $|u| > 1$, this immediately implies $|t_1^g| \geq |t_2^g| - 1 > |\hterm(g')|$
giving us $t_1^g \succ \hterm(g') \succeq t_1^{g'}$ contradicting $t_1^{g'}
\succeq t_1^g$. 
Hence we can assume $\hterm(g) = t_2^g \id \hterm(g')a$ for some $a
\in \Sigma$.
By lemma \ref{lem.lastletter} we know $\hterm(g \mrm \inv{a}) = t_1^g \succ
t_2^g \mrm \inv{a} = \hterm(g')$ and again $t_1^g \succ \hterm(g')
\succeq t_1^{g'}$ contradicts the assumption that $t_1^{g'}
\succeq t_1^g$.
\\
It remains to show that the procedure does indeed terminate.
This is not trivial, as we  either remove a pair of polynomials or
replace a pair of polynomials by another pair of polynomials.
We will use the technique of multisets.
Let $C_i = \{\{ (t_1^g, t_2^g) | g \in G_i \}\}$ be a multiset of
 pairs of terms and set $C_i = \{\{ (\lambda, \lambda) \}\}$ in case
 $G_i = \{ \lambda \}$.
We will show that for all computed sets $G_i$ we have 
 $C_{i+1} \ll C_i$ according to the well-founded tuple-ordering $(t_1^g, t_2^g)
 \succ (t_1^{g'}, t_2^{g'})$ if and only if $t_1^g \succ t_1^{g'}$ or
 $(t_1^g = t_1^{g'}$ and $t_2^g \succ t_2^{g'})$\footnote{This induces
   a well-founded multiset ordering.}.
In case $G_{i+1} = \{ \lambda \}$ or the respective polynomial $q$ in the
algorithm prefix reduces to zero we are done.
Hence it remains to show that in case 
 $p_1 \red{}{\myr}{p}{g} p_2 = p_1 - \alpha \skm g \mrm w\neq 0$ at
 $\hterm(p_1)$ and for some polynomial $q$ we have $(t_1^g,t_2^g) \preceq
 (t_1^{q},t_1^{q})$, $(t_1^{p_1}, t_2^{p_1}) \preceq (t_1^q,t_2^q)$ then this implies 
 $(t_1^{p_2}, t_2^{p_2}) \prec (t_1^{q}, t_2^{q})$.
For the polynomial $q$ chosen by our procedure and $q' = {\rm
  headnormalform}(q, \red{}{\myr}{p}{G'})$ this then implies
$(t_1^q,t_2^q) \succ (t_1^{q'},t_2^{q'})$.
We show our claim by distinguishing the possible cases.
In case $\hterm(p_1) = t_1^{p_1}$ we find $\hterm(p_1) \succ \hterm(p_2) \succeq
t_1^{p_2}$ immediately implies $t_1^{p_1} \succ t_1^{p_2}$.
Hence let us assume $\hterm(p_1) = t_2^{p_1} \id ta$ for some $t \in \free$,
$a \in \Sigma$.
Then by lemma \ref{lem.lastletter} for all $s \in \terms(p_1)$ with $|s| = |t_2^{p_1}|$ we have $\ell(s) =
a$.
Let us take a closer look at $p_2 = p_1 - \alpha \skm g \mrm w$.
\\
In case $w = \lambda$ we get $\hterm(p_1) = \hterm(g)$.
Then if $t_1^g = \hterm(g)$, since $t_1^g \preceq t_1^q$ we get
 $t_1^q \succeq t_1^g = \hterm(g) = t_2^{p_1} = \hterm(p_1) \succ
 \hterm(p_2) \succeq t_1^{p_2}$.
Hence it remains to look at $\hterm(g) = t_2^g$ and since $w = \lambda$
we know $t_2^g \id ta \id t_1^{p_1}$ and this implies that as for all $s \in
\terms(p_1)$, for all $s' \in
\terms(g)$ with $|s'|= |t_2^g|$ we have $\ell (s') = a$.
Furthermore, for all $s \in \terms(p_1)$ and for all $s' \in \terms(g)$
we have $s \mm \inv{a} \preceq t_1^{p_1} \preceq t_1^q$ and $s' \mm \inv{a} \preceq
t_1^g \preceq t_1^q$.
Thus we have $t_1^{p_2} \preceq \hterm(p_2 \mrm \inv{a}) \preceq t_1^{q}$ and in case
$\hterm(p_2 \mrm \inv{a}) = t_1^{p_2} = t_1^q$ we find that either
$t_1^{p_1} = t_1^q$ or $t_1^g = t_1^q$ and hence $t_2^{p_2} \preceq \hterm(p_2) 
\prec \hterm(p_1) = t_2^{p_1} = t_2^g \preceq t_2^q$.
\\
In case $w \id ua$  for all $s \in \terms(\reductum(g
\mrm u))$ we have $s \prec \hterm(g)u  \id t \prec \hterm(p_1 \mrm
\inv{a}) = t_1^{p_1} \preceq t_1^q$.
Thus we get
$t_1^{p_2} \preceq \hterm(p_2 \mrm \inv{a}) = \hterm(p_1 \mrm \inv{a} -
\alpha \skm g \mrm u) = t_1^{p_1} \preceq t_1^q$
and in case $t_1^{p_2} = t_1^{p_1} = t_1^q$ we know, $t_1^{p_2} =
\hterm(p_2 \mrm \inv{a})$, i.e., $t_2^{p_2} \preceq \hterm(p_2) \prec
\hterm(p_1) = t_2^{p_1} \preceq t_2^q$.
\\
\qed 

We will end this section by showing how a special case of our approach
can be compared to the Nielsen method to solve the subgroup problem in
free groups.

Let us start by giving a short description of this method, which can
e.g. be found in \cite{LySch77}.
Let $\free$ be a free group with basis $X$.
We call a word $w \id w_1 \ldots w_k$, $w_i \in \free$,
\betonen{reduced}, in case $w = w_1 \mm \ldots \mm
w_k$, i.e., $|w| = \sum_{i=1}^{k} |w_i|$.
Subsets of $\free$ are written as $U = \{ u_i \mid i \in \n \}$ or $U =
\{ u_1, \ldots , u_n \}$ depending on whether they are finite or not.
Then we can define \index{elementary Nielsen transformations}\index{Nielsen transformation!elementary}\betonen{elementary Nielsen transformations} on a
set $U$ as follows:
\begin{itemize}
\item[(T1)] Replace some $u_i \in U$ by $\inv{u_i}$. 
\item[(T2)] Replace some $u_i \in U$ by $u_i \mm u_j$ where $j
\neq i$. 
\item[(T3)] Delete some $u_i \in U$ where $u_i = \lambda$.
\end{itemize}
In all three cases it is understood that the $u_l$ remain unchanged
for $l \neq i$.
A product of such elementary transformations is called a
\index{Nielsen transformation}\betonen{Nielsen transformation}.
\begin{lemma}~\\
{\sl
If a subset $U$ of $\free$ is carried into a set $U'$ by a Nielsen
transformation, then $U$ and $U'$ generate the same subgroup.
\ohnebeweis
}
\end{lemma}
We call a set $U$ \index{Nielsen reduced}\betonen{Nielsen reduced}, if
for all $v_1,v_2,v_3 \in U \cup \{ \inv{u_i} | u_i \in U \}$ the
following conditions hold:
\begin{itemize}
\item[(N0)] $v_1 \neq \lambda$;
\item[(N1)] $v_1 \mm v_2 \neq \lambda$ implies $|v_1
  \mm v_2| \geq \max \{ |v_1|, |v_2| \}$;
\item[(N2)] $v_1 \mm v_2 \neq \lambda$ and $v_2 \mm
  v_3 \neq \lambda$ imply $|v_1 \mm v_2 \mm v_3| >
  |v_1| - |v_2| + |v_3|$.
\end{itemize}
%
%
\auskommentieren{
\begin{theorem}~\\
{\sl
Given a finite Nielsen reduced set $U$ there is an algorithm which
decides the subgroup problem.
\theoend
}
\end{theorem}
\Ba{}~\\
The algorithm to decide the subgroup problem uses the following lemma
for Nielsen reduced sets $U$:
\\
If $w = w_1 \mm \ldots \mm w_n$ where $w_i \in U \cup \{ \inv{u} | u \in U
\}$ and $w_i \mm w_{i+1} \neq \lambda$ for all $1 \leq i < n$,
then $|w| \geq n$ (compare lemma 2.13 in \cite{LySch77}).
\\
Hence, to decide whether $w$ is in the subgroup generated by $U$ it is
sufficient to check whether $w = w_1\mm \ldots \mm w_m$ for some $w_i \in U
\cup \{ \inv{u} | u \in U \}$ where $m \leq |w|$.
This is effective, as $U$ is finite.
\\
\qed}
Nielsen reduced sets play an important role, as they are free
generating systems for the subgroup they generate.
The following theorem due to Ziechang states that freely reducing a
product of elements of a Nielsen reduced set cannot result in
arbitrary cancellations on the elements involved.
\begin{theorem}\label{theo.ziechang}~\\
{\sl
Let $U$ be a Nielsen reduced set.
Then for every $u \in U \cup \{ \inv{u} \mid u \in U \}$ there are words
$a(u)$ and $m(u)$ with $m(u) \neq \lambda$ such that 
$u \id a(u)m(u)\inv{a(\inv{u})}$ and if $w = u_1 \mm \ldots \mm u_n$
for some $u_i \in U \cup \{ \inv{u} \mid u \in U \}$, $u_{i} \mm u_{i+1}
\neq \lambda$, then the words $m(u_i)$ remain uncancelled in the
reduced form of $w$.
In particular we get $|w| \geq n$.
\theoend\ohnebeweis
}
\end{theorem}
This property can be used to solve the subgroup problem for Nielsen
reduced sets by computing appropriate right coset representations.
Therefore, it remains to find Nielsen reduced sets.
The following theorem gives an effective way to transform an arbitrary
finite set $U$ into a Nielsen reduced set.
To see how the necessary Nielsen transformation can be constructed we include a
proof.
\begin{theorem}~\\
{\sl
Let $U \subseteq \free$ be a finite set.
Then there is a Nielsen transformation from $U$ into some Nielsen reduced
set $V$.
\theoend
}
\end{theorem}
\Ba{}~\\
Let $U$ be a finite subset of $\free$.
We will show how $U$ can be carried over into a Nielsen reduced set by
using elementary Nielsen transformations.
\\
Condition (N0) can be achieved using finitely many transformation
steps (T3), and hence we can assume that $U$ satisfies (N0).
\\
Now suppose that $U$ does not fulfill (N1).
Then without loss of generality we can assume that there are $u_i,
u_j$ such that $|u_i \mm u_j| < |u_i|$.
Then $j \neq i$, as $|u^2| < |u|$ is not possible for $u \in \free$.
Using transformation (T2) we can replace $u_i$ by $u_i \mm
u_j$ and for the new set $U'$ we have $\sum_{u \in U'} |u| < \sum_{u
  \in U} |u|$.
Hence, we can assume that $U$ can be transformed using (T2) into a set
$U'$ with $\sum_{u \in U'} |u|$ minimal, i.e., no further applications
of (T2) are possible and hence condition (N1) must hold.
Since (T2) does not affect condition (N0) we can now assume that $U$
satisfies (N0) and (N1).
\\
Finally let us transform $U$ into a set additionally fulfilling (N2).
To see how this can be done let us consider a triple $x,y,z$ such that
$x \mm y \neq \lambda$ and $y \mm z \neq \lambda$.
As $U$ satisfies (N1) we know $|x \mm y| \geq |x|$ and $|y
\mm z| \geq |z|$, i.e., the part of $y$ which cancels in $x
\mm y$ is no more than half of $y$, and likewise the part that
cancels in $y \mm z$.
Now let $x \id w_1w_2$, $y \id \inv{w_2}w_3w_4$, $z \id \inv{w_4}w_5$
such that $x \mm y \id w_1w_3\inv{w_4}$ and $y \mm z
\id w_2w_3w_5$.
In case $w_3 \neq \lambda$ we find $x \mm y \mm z \id
w_1w_3w_5$ and hence $|x \mm y \mm z| = |x| - |y| +
|z| + |w_3| > |x| - |y| + |z|$ and thus (N2) holds for this triple.
Otherwise, we get $y \id w_2\inv{w_4}$ and (N2) is violated.
Note that, as $U$ satisfies (N1), we have $|w_2| = |w_4| = \frac{1}{2}
|y| \leq \min \{
\frac{1}{2} |x|, \frac{1}{2} |z| \}$ and $w_2 \neq w_4$.
We could now use transformation (T2) and either replace $\inv{x}$ by
$\inv{x \mm y}$ or $z$ by $y \mm z$ without changing
the sum of the lengths of the elements in $U$.
To decide which replacement should take place, we use the following technique: \\
Suppose there is a precedence on the letters $X \cup X^{-1}$ inducing a
length-lexicographical ordering on the reduced words presenting the
elements of $\free$.
We define the left half of a reduced word $w$ to be the initial
segment  $L(w)$ of length $[\frac{|w|+1}{2}]$.
This can be used to define a well-founded ordering on reduced words as
follows:
For two reduced words $w_1,w_2$ we set $w_1 \prec w_2$ if and only if
$\min \{ L(w_1),L(\inv{w_1}) \} < \min \{ L(w_2),L(\inv{w_2}) \}$ or 
$(\min \{ L(w_1),L(\inv{w_1}) \} = \min \{ L(w_2),L(\inv{w_2}) \}$ and
$\max \{ L(w_1),L(\inv{w_1}) \} < \max \{ L(w_2),L(\inv{w_2}) \})$.
Now suppose $x \id w_1\inv{w_2}$, $y \id w_2\inv{w_4}$ and $z \id
w_4w_5$ as above.
If $w_2 \prec w_4$ then $y \mm z \id w_2w_5 \prec z \id w_4w_5$
and if $w_4 \prec w_2$ then $x \mm y \id w_1\inv{w_4} \prec x \id
w_1\inv{w_2}$.
We can now suppose that the set $U$ is transformed using (T2)
according to the relation $\prec$ as far as possible.
Since this terminates and does not affect condition (N0) nor condition (N1) we
are done.
\\
\qed
There are well-known algorithms for performing this task and Avenhaus
and Madlener have provided one which works in polynomial time (see \cite{AvMa84}).
We will now proceed to show how Nielsen's method is related to
solving the generalized word problem in free groups using Gr\"obner
bases.
Applying theorem \ref{theo.subgroup.problem} we find that the  subgroup
problem related to a finite set $U \subseteq \free$ can be transformed into the
membership problem for the right ideal generated by the set of
polynomials $P_U = \{ u-1 \mid u \in U \}$.
The following lemma states that special prefix Gr\"obner bases of $P_U$ exist.
\begin{lemma}~\\
{\sl
For a finite subset $U$ of $\free$
let $G$ be  the reduced prefix Gr\"obner basis computed by
 procedure {\sc Completion in Free Group Rings} 
 on input $P_U = \{ u-1 \mid u \in U \}$ assuming  a length
 lexicographical ordering on $\free$ assuming that all
 polynomials are made monic.
Then the following conditions hold:
\begin{enumerate}
\item No head term of a polynomial in $G$ is a prefix of a head term
  of another polynomial in $G$.
\item For all $g \in G$ we have $g = u -v$ for some $u,v \in \free$.
\item For $u-v \in G$ we have $u \mm \inv{v} \id u\inv{v}$, in particular
      $\ell(u) \neq \inv{\ell(v)}$.
\item For $u-v \in G$ we have $|u| - |v| \leq 2$.
\item If $u-v \in G$ so is its `mate' $(-1) \skm (u-v) \mrm
  \inv{\ell(u)}$.
\end{enumerate}
\lemend
}
\end{lemma}
\Ba{}
\begin{enumerate}
\item This follows immediately from theorem \ref{theo.correct.free.group}.
\item The only possibilities of changing polynomials are the saturation 
       process and the normal form computation.
      By definition \ref{def.freegroupsat}, for $p = u-v$ the polynomials
       $\can(p)$ and $\satpoly(p)$ are also
       of this form in case they are made monic.
      Let us hence take a closer look at prefix reducing $p=u-v$ by
       a polynomial $u'-v'$ at $u$.
      Then $u \id u'z$ for some $z \in \free$ and the result is
       $u-v - (u'-v') \mrm z = -v + v' \mm z$ and again the monic
       version of this polynomial has the desired form.
\item This follows immediately from the fact that all polynomials in
       $G$ are either in $\can$ or $\satpoly$ form.
\item  The case $|u| - |v| > 2$ would contradict the previous statement.
\item  To see this let us assume $(-1) \skm (u-v) \mrm \inv{\ell(u)}
  \not\in G$ and set $a = \inv{\ell(u)}$.
  Remember that $u-v$ is either in its $\can$ or $\satpoly$ form and hence
   $v \mm a \id va$.
  Since we assume $va - u \mm a \not\in G$, then $va$ must be 
    prefix reducible by a polynomial $g \in G$.
  As $u-v \in G$, $v$ is not prefix reducible using $G$, hence
  $\hterm(g) \id va$, say $g = va - z$ for some $z \in \free$ and $z
  \neq u \mm a$.
  But then we have $(u-v) \mrm a, va-z \in \ideal{r}{}(G)$, implying
  $(u-v) \mrm a - (va-z) = - u \mm a + z \in \ideal{r}{}(G)$.
  Thus either $z$ or $u \mm a$, which is a prefix of $u$, must
  be prefix reducible using $G$ contradicting that $u-v$ and $va-z$
  are supposed to belong to $G$ and hence must be prefix reduced.
\\
\qed
\end{enumerate}\renewcommand{\baselinestretch}{1}\small\normalsize
\auskommentieren{
We will next state how the steps of the algorithm can be compared to
 elementary Nielsen transformations.
To do this we add a slightly modified transformation rule  as follows:
\begin{itemize}
\item[(T1')] If $u_i \in U$ but $\inv{u_i} \not\in U$ then add
    $\inv{u_i}$ to $U$.
\end{itemize}
Note that all properties of Nielsen transformations and Nielsen
reduced sets mentioned above remain true for this modification.
\\
The operations  involved in the procedure are saturating and prefix
reducing a polynomial.
Let $S_i,G_i$ be the respective sets of polynomials computed and
associate a set $X_i = \{ u \mm \inv{v} \mid u-v \in S_i \cup G_i
\}$ to them.
We will show that changes to the sets $S_i,G_i$ can be described by
a Nielsen transformation on the accompanying set $X_i$.
\\
Note that the set $S_0 \cup G_0$ contains the polynomials $\can(u-1)$
 where $u-1 \in P_U$.
Hence, $X_0 \subseteq U \cup \{ \inv{u} \in U \}$.
This follows, as if we let $L(u)$ be the initial segment of length
$[\frac{|u|+1}{2}]$ of a reduced word $u$ and $R(u)$ the remaining
segment, then the candidates for $\can(u-1)$ and $\satpoly(u-1)$
belong to the set $\{ L(u) - \inv{R(u)}, (L(u) - \inv{R(u)}) \mrm
\inv{\ell(L(u))}, (L(u) - \inv{R(u)}) \mrm \inv{\ell(R(u))} \}$, i.e.,
the accompanying elements in $X_0$ are of the form $L(u)R(u) \id u$ and
$\inv{R(u)} \inv{L(u)} \id \inv{u}$.
Further $X_0$ can be achieved from $U$ by a Nielsen transformation
involving the rules (T1) and (T1').
\\
It remains to show, that reducing a polynomial $f \in S_i$ to its
prefix normal form $f'$ and saturating $f$ by $\can(f')$ and
$\satpoly(f')$ can also be reflected by a Nielsen transformation.
\\
This will be done by investigating one reduction step, as we will find
that more reduction steps can be combined to multiple (T2) transformations.
\\
Let us assume $u-v \red{}{\myr}{p}{u'-v'} v' \mm z - v$ and $u \id u'z$ for some $z
\in \free$\footnote{The case $v\id u'z$ can be treated similarly.}.
\\
In case $v' \mm z > v$ we have to show that $v' \mm z
\mm \inv{v}$ can be achieved from elements in $X_i$ by a
Nielsen transformation. 
\\
As $u-v, u'-v' \in S_i \cup G_i$ we know $u \mm \inv{v}, u'
\mm \inv{v'} \in X_i$.
Hence, we can assume $\inv{u' \mm \inv{v'}} \in X_i$\footnote{
This can always be achieved by an application of either (T1) or (T1').}.
Since we have $\inv{u' \mm \inv{v'}} \mm (u
\mm \inv{v}) = (v' \mm \inv{u'}) \mm (u'
\mm z \mm \inv{v}) = v' \mm z \mm
\inv{v}$, we may replace $u \mm \inv{v}$ by $v' \mm z
\mm \inv{v}$ using transformation (T2).
\\
Hence, as for a prefix normal form $g$ of $u-v$ the polynomials
$\can(g)$ and $\satpoly(g)$ are added to $S_i$ this can be reflected by (T2)
transformations for the reduction process and a (T1') transformation
for the saturation process.
Remember that in case $w - z = g$ and $w'-z' = \can(g)$, $(w'-z') \mrm
a = \satpoly(g)$, then either $w \mm \inv{z} \id w'\inv{z'}$ or 
 $w \mm \inv{z} \id \inv{w'\inv{z'}}$ and
 $(z'a) \mm \inv{w' \mm a} = (z'a) \mm
(\inv{a} \mm \inv{w'}) = z' \mm \inv{w'}$.
\\
In case $v > v' \mm z$ we can show that $v \mm \inv{v'
  \mm z}\in X_i$ similarly, as $v \mm \inv{ v' \mm z}
= \inv{ v' \mm z \mm \inv{v}}$.
}
\begin{theorem}~\\
{\sl
Let $U$ be a finite subset of $\free$ and $G$ the monic  reduced
 prefix Gr\"obner of the right ideal generated by $\{ u - 1 \mid u \in U
 \}$ in $\myk[\free]$.
Then the set $X_G = \{ u\inv{v} \mid u-v \in G \}$ is Nielsen reduced for $U$.
\theoend
}
\end{theorem}
\Ba{}~\\
We have to show that the set $X_G$ satisfies the conditions (N0), (N1)
and (N2). 
\\
(N0) is valid, since $u\inv{v} \id \lambda$ would imply $u=v=\lambda$,
 but we assume that $0$ is not contained in $G$. 
\\
To show that (N1) is true, we prove that for two polynomials
 $u-v,u'-v' \in G$, $u\inv{v} \mm u'\inv{v'} \neq \lambda$ implies
  $|u\inv{v} \mm u'\inv{v'}| \geq \max \{
  |u\inv{v}|,|u'\inv{v'}| \}$.
Let us assume $u\inv{v} \id w_1w_2$ and $u'\inv{v'} \id \inv{w_2}w_3$ such
 that $u\inv{v} \mm u'\inv{v'} \id w_1w_3$ and $w_1w_3 \not \id
 \lambda$. 
Since $u-v \in G$ we know that its mate $(-1) \skm (u-v) \mrm a$ belongs
 to $G$, where $a = \inv{\ell(u)}$, and this polynomial has head term $va$.
Let us first assume that $\inv{va}$ is a suffix of $w_2$, i.e., 
 $va$ is a prefix of $\inv{w_2}$.
Then, in case  $u'$ is a prefix of
 $\inv{w_2}$, this would  imply that one of the terms $va$
 or $u'$ is prefix of the other which would contradict the
 fact that $G$ is a reduced prefix Gr\"obner basis unless
 we have $u' \id va$.
But then, as $G$ contains the mate of $u-v$,
 this mate must be $u'-v'$ contradicting the fact that we require $u\inv{v}
 \mm u'\inv{v'} \neq \lambda$. 
Assuming that $\inv{w_2}$ is a prefix of $u'$ would immediately give a
 contradiction as $va$ then would be a prefix of $u'$ implying that
 $u'-v'$ is prefix reducible by the mate of $u-v$.
Hence let us assume that $w_2$ is a suffix of $\inv{v}$, i.e., 
 $\inv{w_2}$ is a proper prefix of $va$.
Then, as $|\inv{w_2}| < |va|$ and $|va| \leq [\frac{|u\inv{v}|}{2} +
1]$ as $|u| - |v| \leq 1$.
Thus at most half of $u\inv{v}$ is cancelled by multiplication with $w_2$ and
$|w_1| \geq |u|$.
Now in case $\inv{w_2}$ is also a proper prefix of $u'$ this implies
$|w_2| < |u'|$ and $|w_3| > |v'|$.
Hence, as $|u| - |v| \leq 1$ and $|u'| - |v'| \leq 1$ this implies
 $|w_2| \leq  |w_1|$ and $|w_2| < |w_3|$, i.e.
 $|w_1w_3| = |w_1| + |w_3| > |w_1| + |w_2| = |w_1w_2| = |u \inv{v}|$ and $|w_1| + |w_3| \geq
 |w_2| + |w_3| = |\inv{w_2}w_3| = |u'\inv{v'}|$ and hence
 $|u\inv{v} \mm u'\inv{v'}| \geq \max \{
  |u\inv{v}|,|u'\inv{v'}| \}$. 
Note that $u'$ cannot be a  prefix of $\inv{w_2}$ as then $va$
 would be prefix reducible using $u'$ contradicting that the mate of
 $u-v$ belongs to $G$. 
\\
To show that (N2) holds, we prove that for three polynomials
$u-v$,$u'-v'$,$u''-v'' \in G$, $u\inv{v} \mm u'\inv{v'} \neq
\lambda$ and $u'\inv{v'} \mm u''\inv{v''} \neq \lambda$ imply that the equation
 $|u\inv{v} \mm u'\inv{v'} \mm u''\inv{v''}| > 
  |u\inv{v}| - |u'\inv{v'}| + |u''\inv{v''}|$ holds.
Let us assume $u\inv{v} \id w_1w_2$, $u'\inv{v'} \id \inv{w_2}w_3w_4$ 
 and  $u''\inv{v''} \id \inv{w_4}w_5$ such
 that $u\inv{v} \mm u'\inv{v'} \id w_1w_3w_4$ and 
 $u'\inv{v'} \mm u''\inv{v''} \id \inv{w_2}w_3w_5$. 
Then, as 
 $|u\inv{v} \mm u'\inv{v'} \mm u''\inv{v''}| 
  =|w_1w_3w_5| = |w_1w_2| - |\inv{w_2}w_3w_4| + |\inv{w_4}w_5| + |w_3| \geq 
  |u\inv{v}| - |u'\inv{v'}| + |u''\inv{v''}|$,
 in case $w_3 \neq \lambda$ we are done.
Hence let us assume $w_3 = \lambda$.
Then, as (N1) holds, we get 
 $|u\inv{v} \mm u'\inv{v'}| \geq \max \{ |u \inv{v}|,
 u'\inv{v'}| \}$ and
 $|u'\inv{v'} \mm u''\inv{v''}| \geq \max \{ |u' \inv{v'}|,
 u''\inv{v''}| \}$
 implying that $|w_2|=|w_4| = \frac{1}{2}|u'\inv{v'}|$ and $w_2 \neq
 w_4$\footnote{We have $u'\inv{v'} \neq \lambda$ since (N0) holds.}.
Thus, $u' \id \inv{w_2}$ and $v' \id \inv{w_4}$ since $|u'| - |v'|
\leq 1$.
This would imply
that either $v$ is prefix reducible, in case $\inv{w_2}$ is a prefix
of $v$ contradicting that $u-v \in G$, or $u'$ is prefix reducible by
the mate of $u-v$ contradicting that $u'-v' \in G$.
Therefore, $w_3 = \lambda$ is not possible and we are done.
\\
\qed
In particular the steps performed  in procedure {\sc
  Reduced Prefix Gr\"obner Bases} modified for free group rings 
  or in procedure {\sc Completion in
  Free Group Rings} can be compared to Nielsen transformations on a
set $X$ corresponding to the actual sets of polynomials, namely in the
first case $X = \{ u \mm \inv{v} \mid u-v \in S_i \cup G_i \}$ and in
the second case $X = \{ u \mm \inv{v} \mid u - v \in G_i \}$.
Both algorithms in changing the actual polynomials mainly involve 
 saturation and prefix reduction.
We close this section in sketching how these two operations are
related to Nielsen transformations.
Let us start with investigating the polynomials $\can(u-1)$ and
$\satpoly(u-1)$ as saturating an arbitrary polynomial of the form $u-v$
can be reduced to saturating the polynomial $u \mm \inv{v} - 1$.
We claim that replacing a polynomial $u-1$ by $\can(u-1)$ and
$\satpoly(u-1)$ corresponds to a Nielsen transformation involving the
rules (T1) and (T1') on the respective set $X$ involving the element
$u \in X$.
To see this let $L(u)$ be the initial segment of length
$[\frac{|u|+1}{2}]$ of a reduced word $u$ and $R(u)$ the remaining
segment, then the candidates for $\can(u-1)$ and $\satpoly(u-1)$
belong to the set $\{ L(u) - \inv{R(u)}, (L(u) - \inv{R(u)}) \mrm
\inv{\ell(L(u))}, (L(u) - \inv{R(u)}) \mrm \inv{\ell(R(u))} \}$, i.e.,
the accompanying elements are of the form $L(u)R(u) \id u$ and
$\inv{R(u)} \inv{L(u)} \id \inv{u}$ and hence we can use the rules
(T1) and (T1') to do the appropriate changes on the set $X$.
It remains to show how prefix reduction steps in this context are
related to Nielsen transformations.
We show that a single reduction step corresponds to an application of
a rule (T2).
Let us assume $u-v \red{}{\myr}{p}{u'-v'} v' \mm z - v$ and $u \id u'z$ for some $z
\in \free$\footnote{The case $v\id u'z$ can be treated similarly.}.
In case $v' \mm z > v$ we have to show that $v' \mm z
\mm \inv{v}$ can be achieved from elements in $X$ by a
Nielsen transformation. 
As $u-v, u'-v'$ must belong to the actual set of polynomials we know $u \mm \inv{v}, u'
\mm \inv{v'} \in X$.
Hence, we can assume $\inv{u' \mm \inv{v'}} \in X$\footnote{
This can always be achieved by an application of either (T1) or (T1').}.
Since we have $\inv{u' \mm \inv{v'}} \mm (u
\mm \inv{v}) = (v' \mm \inv{u'}) \mm (u'
\mm z \mm \inv{v}) = v' \mm z \mm
\inv{v}$, we may replace $u \mm \inv{v}$ by $v' \mm z
\mm \inv{v}$ using transformation (T2).
\section{Plain Groups}
A plain group is a free product of finite groups and a free group,
 and in \cite{AvMaOt86}]
 it has been shown that these groups allow finite 2-monadic, convergent
 (even reduced) group presentations.
Using these presentations and the syntactical information they
 provide we can show that a slight modification of 
  procedure  {\sc Prefix Gr\"obner Bases} (see page
  \pageref{prefix.groebner.bases}) terminates.
\begin{theorem}\label{theo.term.plain.groups}~\\
{\sl
Given  a 2-monadic confluent group presentation for a plain group 
 $\g$ and a finite set of polynomials $F \subseteq \myk[\g]$,
 the procedure {\sc Prefix Gr\"obner Bases}
 terminates.
\theoend
}
\end{theorem}
\Ba{}~\\
Note that if $( \Sigma, T)$ is a convergent interreduced  presentation of
 a cancellative monoid $\m$,
 then no rules of the form $wa \myr a$ or $aw \myr a$ appear in $T$ for $a \in\Sigma$.
This is of course always true if such  presentations are given for
groups.
\\
Let us assume that procedure {\bf normalform} computes a normalform
 of a polynomial allowing only prefix reduction steps at the respective
 head terms.
Then the  proof  is done in two steps:
first we show that all polynomials computed have a certain property that will be used in the
 second step to ensure termination.
 We say a polynomial $q$ has property $\p_F$ if and only if
      \begin{enumerate}
      \item[($\alpha$)] $|\hterm(q)| \leq K$, where $K = \max \{ |\hterm(f)| \mid f \in F \} +1$.
      \item[($\beta$)] If $|\hterm(q)|=K$ then there exists an element $a \in \Sigma$ such that
           \begin{enumerate}
             \item[(i)] all terms of length $K$ in $q$ have $a$ as a
               common suffix, and
             \item[(ii)] for all $s \in\terms(q)$ with  $|s|=K-1$ we
               either have $s \id s_1a$ or in case $s \id s_1d$,
                   $d \in \Sigma\backslash\{ a \}$ there is a rule $ea \myr d \in T, e \in \Sigma$.
           \end{enumerate}
      \end{enumerate}
      We will show that all polynomials $q$ computed by the
       procedure on input $F$ have property $\p_F$. \\
      By the choice of $K$ all input polynomials have $\p_F$.
      Hence, let  $G$ be the actual  set of polynomials having $\p_F$,
       and let $q$ be the next polynomial computed by our procedure.
      In case $q$ is due to computing the normal form of a
       polynomial $p$ having $\p_F$ using prefix reduction at head terms only the property is preserved.
      To see this we can  restrict ourselves to a single step reduction.
In case $|\hterm(p)|<K$ we are done.
            Therefore, suppose $|\hterm(p)|=K$ and $\hm(p)$
            is reduced in the reduction step $p \red{}{\myr}{p}{g \in G}
            q'$.
            We have to show that $q'$ satisfies $\p_F$.
            Let $\hterm(p)\id \hterm(g)w$ and
             $q'= p - \alpha  \skm g \mrm w$, $\alpha \in\myk^*$, $w \in\m$.
            Now $g \mrm w$ has $\p_F$ as $\hterm(g \mrm w) \id
            \hterm(g)w$ and for all $s \in \terms(\reductum(g))$
            we either have $|s \mm w| < |sw|$ or $sw$ and $\hterm(g)w$
            have the same last letter.
            Since 
             $\terms(q') \subseteq\terms(p) \cup\terms(g \mrm w)$, 
             $q'$ then likewise has $\p_F$.
      In case $q$ is due to saturating a polynomial as specified e.g.
      in procedure {\sc Prefix Saturation} on page
       \pageref{prefix.saturation} and results from a polynomial $q'$
       having $\p_F$ being overlapped with a rule $ab \myr c \in T$, $c
       \in \Sigma \cup \{ \lambda \}$\footnote{The polynomial $q'$
         here is said to overlap with the rule $ab \myr c \in T$ in case
         $\ell(\hterm(q')) = a$.},  we
       can also show that $\p_F$ is preserved.
      Note that only the case $|\hterm(q)|=K$ is critical.
               In case $|\hterm(q')| <K$ and $|\hterm(q)|=K$ we know
               $\hterm(q) \id tb$
               and for all $s \in\terms(q')$
               with $|s \mm b|=K-1$ either $s \mm b \id sb \in \irr(T)$ or
               $s  \id s_1e$ and $s \mm b = s_1e \mm b \id s_1d$,
               where $eb \myr d \in T$.
              Note that these are the only possibilities to gain a
               term of length $K-1$ from a term
               of length less or equal to $K-1$ by multiplication
               with a letter $b$.
              On the other hand, if  $|\hterm(q')| = K$ with $\hterm(q') \id ta$
              we can only violate $\p_F$ in case we have 
               $t_1,t_2 \in\terms(q')$ such that
               $|t_1|=K,|t_2|=K-1, t_1 \id t'_1a$ and 
               $t_1 \mm b  \id t'_1c, t_2 \mm b \id t_2b$ with $c \neq \lambda$.
              Therefore, we examine all $s \in\terms(q')$  with $|s|=K-1$.
              If there are none  $q$ must have $\p_F$, since then a
               term $s \in\terms(q')$ can only reach length $K-1$ 
               by multiplication with $b$ in case $|s|=K-2$ and
               $sb \in \irr(T)$.
              Since $ab \myr c\in T$ and $\g$ is
               a group including inverses of length 1 for the generators,
               $a$ has an inverse $\tilde{a}$ and  $b
               \red{*}{\lr}{}{T} \underline{\tilde{a}a}b \id
               \tilde{a}\underline{ab}\red{*}{\lr}{}{T} \tilde{a}c$ gives 
               us the existence of a rule $\tilde{a}c \myr b \in T$ as $T$ is
                confluent\footnote{This is no longer true in case $a$ has an
                         inverse $u_a$ of length $|u_a|>1$ or 
                         no inverse at all.}.
              Now let $s \in\terms(q')$ have length $K-1$.
              Then if $s\id s_1a$ there is nothing to show\footnote{Then
                       $s \mm b = s_1a \mm b = s_1 \mm c$ and
                        either $|s \mm b|<K-1$ or $s \mm b \id s_1c$.}. 
              On the other hand, in case $s \id s_1d, d \neq a$ we know that there is
               a rule $ea \myr d \in T$ as $q'$ has $\p_F$.
              Then we have $db \la \underline{ea}b \id e\underline{ab} \myr
               ec$ and, since $ea \myr d \in T$ gives us $e \neq d$, 
               there are rules $db \myr g, ec \myr g \in T, g \in \Sigma
               \cup \{ \lambda \}$.
       Finally let us assume that $q$ is due to s-polynomial
       computation.
       But computing s-polynomials can be compared to a single prefix reduction
       step on the head monomial of a polynomial and we have seen that
       prefix reduction preserves property $\p_F$.
      \auskommentieren{ 
       As we use polynomials having $\p_F$ and prefix overlaps
        of their head terms to get $q$, 
        in case $|\hterm(q)|=K$, $q$ inherits $\p_F$ from the involved
        polynomials $q_1, q_2$, where
        $\hterm(q_1) \id \hterm(q_2)z$ for some $z \in\m$.
\\
       In case $z = \lambda$ there is nothing to show.
\\
       Let $|\hterm(q_2)z|=K$ and $z \neq \lambda$, i.e., $z \id z'a$
        for some $a \in\Sigma$.
       Looking at $s \in\terms(q_2 \mrm z)$ we get that
       \begin{enumerate}
       \item $|s|<K-1$ is not critical. 
       \item $|s|=K$ gives us $s = s_1 \mm z$, where
             $s_1 \in\terms(q_2), |s_1|<K$, and as $s_1, z \in \irr(T)$
             and $T$ is monadic, we get $s \id s'_1bz''a$,
             where $s'_1$ is a prefix of $s_1$,
             $b \in\Sigma \cup \{ \lambda \}$, and $z''$ is a suffix of $z'$. 
       \item $|s|=K-1$ gives us $s = s_1 \mm z$, where 
             $s_1 \in\terms(q_2), |s_1|<K$, and as $s_1, z \in \irr(T)$
             and $T$ is monadic either $s \id s'_1bz''a$ as above or
             $s \id s'_1b$ and $|s_1|=K-1,s_1 \id s'_1e$ and
             $ez'a \red{*}{\myr}{}{T} fa \myr b$, as $T$ is 2-monadic,
             i.e., we have $fa \myr b \in T$.
       \end{enumerate}
       Since $\terms(q) \subseteq\terms(q_1) \cup\terms(q_2 \mrm z)$ we are done.}
\\
It remains to show that the procedure does terminate.
Thus let us assume the contrary.
       Then there are infinitely many polynomials $q_i, i\in\n$
       resulting from s-polynomial computations added to $G$.
       Note that every such polynomial is in prefix normal from with respect
       to all polynomials in $G$ so far.
       On the other hand, as $|\hterm(q_i)| \leq K$,
        this would mean that there is a term $t$,
        which occurs infinitely often as a head term among these
        polynomials $q_i$ contradicting the fact that the head terms
        of all added polynomials
        are in prefix normal form with respect to the polynomials
        added to the Gr\"obner set so far, and hence no head term can
        appear twice among the head terms of the polynomials ever
        added to the set $G$.
\\
\qed
An immediate consequence of this theorem is the existence of finite
 Gr\"obner bases for finitely generated right ideals in free and plain group
 rings.

Following the approach sketched for free group rings in the previous
 section, we can give a more efficient prefix saturating procedure
 for plain group rings.
We will use the following observation from Kuhn in \cite{Ku91}.
\begin{lemma}~\\
{\sl
Let $(\Sigma,T)$ be a 2-monadic, convergent, reduced group
 presentation of a plain group $\g$.
Then for $w \in \g$ and $b \in \Sigma$ we have
\[ w \mm b  \id \left\{ \begin{array}{r@{\quad\quad}l}
               wb &  wb \mbox{ is irreducible} \\
               w'c & w \id w'a, (ab,c) \in T,
              \end{array} \right.\]
in particular, $ |w|-1 \leq |w \mm b| \leq  |w|+1$.
\lemend
}
\end{lemma}
\Ba{}~\\
In case $w \mm b \id wb$ we immediately get $|w \mm b| = |w| +1$.
On the other hand, if $w \id w'a$ and $w \mm b \id w'c$ in case $c \in
 \Sigma$  we have
 $|w \mm b| = |w|$ or if $c = \lambda$, $|w \mm b| = |w| -1$.
It remains to show that multiplication with $b$ cannot result in a
 larger decrease of length.
Let us assume $w \id w'a_2a_1$ and $|w \mm b| < |w| -1$.
Then there must be at least two rules $a_1b \myr c$, $a_2c \myr d \in T$
 with $c \in \Sigma \backslash \{ a_1,b \}$, $d \in \Sigma \backslash
 \{ a_2, c \} \cup \{ \lambda \}$.
This implies there is  a rule $c \,\inv{b} \myr a_1 \in T$, as $(\Sigma,T)$
 is a group presentation and $|\inv{b}|=1$.
Thus $d \, \inv{b} \la \underline{a_2c\,} \inv{b} \id
      a_2\underline{c\,\inv{b}} \myr a_2a_1$ either implies $a_2a_1 \myr
 \inv{b} \in T$ in case $d = \lambda$ or there exists $e \in \Sigma
 \cup \{ \lambda \}$ such that $a_2a_1 \myr e$, $d \,\inv{b} \myr e \in T$ and
 in both cases $a_2a_1$ is reducible, contradicting our assumption
 that $w \id w'a_2a_1 \in \g$.
\\
\qed
\begin{lemma}\label{lem.sat.in.pgr}~\\
{\sl
Let $(\Sigma,T)$ be a 2-monadic, convergent, reduced group
 presentation of a plain group $\g$.
For a polynomial $p \in \myk[\g]$ containing more than one monomial 
 we define $\sigma_1(p)$ and $\sigma_2(p)$
 as in definition \ref{def.freegroupsat}.
Further for $q \in \myk[\g]$ let 
$$C_q = \{ b \in \Sigma \mid (ab,c) \in T, \mbox{ where }
        a=\ell(\hterm(q)), c \in \Sigma \}.$$
Then  
 $$\s_p(p) =  \{ \can(p), \satpoly(p) \} \cup \{ \can(p) \mrm d \mid d \in C_{\can(p)}\}
        \cup \{ \satpoly(p) \mrm d \mid  d \in C_{\satpoly(p)} \}$$
  is a prefix saturating set for $p$.
In case $p = \alpha \skm t \in \myk[\g]$ we can set $\s_p(p) = \{
\lambda \}$.
\lemend
}
\end{lemma}
\Ba{}~\\
We have to show that the polynomials in the set $\{ \alpha \skm p \mrm
w \mid \alpha \in \myk^*, w \in \g \}$ are prefix reducible to zero in one
step by $\s_p(p)$.
In case $p = \alpha \skm t$, $\alpha \in \myk^*$, $t \in \g$, we are done as
 $\s_p(p) = \{ \lambda \} \in \SAT(p)$.
In case the polynomial $p$ contains more than one monomial,
 we use  procedure {\sc Prefix
   Saturated Check} on page \pageref{is.prefix.saturated} to prove our
 claim by showing that for every polynomial  $q \in \s_p(p)$ and every
$w \in C(\hterm(q))$ the multiple $q \mrm w$ is prefix reducible to
zero in one step using $\s_p(p)$.
Let $\hterm(\can(p)) \id ta$ and $\hterm(\satpoly(p)) = t'\mm\inv{a}$
for some $t,t' \in \g$, $a \in \Sigma$.
In case $q \in \{ \can(p), \satpoly(p) \}$, the fact that
 $C(\hterm(q)) = C_q \cup \{ \inv{\ell(\hterm(q))} \}$ and 
 the definition of $\s_p(p)$  imply that for all $b \in C(\hterm(q)) $ we have $q \mrm
 b \red{}{\myr}{p}{S_p} 0$.
Now, let us assume that $q = \can(p) \mrm b$ for some $b \in
C_{\can(p)}$ and $(ab,c) \in T$, $c \in \Sigma$.
We have to distinguish the following two cases.
If $\hterm(q) \id tc$, then $C(tc) = \{ d \mid (cd, e) \in T, d \in \Sigma, e \in \Sigma
      \cup \{ \lambda \} \}$ and in case this set is
      not empty let us look at such a rule $(cd, e) \in T$.
      Since our presentation is a reduced  convergent group presentation, there
      exists a rule of the form $\inv{a}c \myr b \in T$ where $|\inv{a}| = 1$.
      Now this gives us 
      $$bd \la \underline{\inv{a}c}d \id
      \inv{a}\underline{cd} \myr \inv{a}e$$ and as $d \neq e$ and $b \neq
      \inv{a}$,
      there exists an element $f \in \Sigma \cup \{ \lambda \}$ 
      such that $bd \myr f$, $\inv{a}e \myr f \in T$.
      Again this results in the situation
      $$cd \la \underline{ab}d \id a\underline{bd} \myr af$$
      and we either have $b \mm d = \lambda$ in case $f = \lambda$
      or there exists a rule $af \myr e \in T$.
      In case $b \mm d = \lambda$ this implies $q \mrm d = (\can(p) \mrm b) \mrm d =
      \can(p) \mrm (b \mm d) =
      \can(p)$ and hence $q \mrm d \red{}{\myr}{p}{\sm_p(p)} 0$.
      Otherwise,
      $q \mrm d = (\can(p) \mrm b) \mrm d =\can(p) \mrm (b \mm d) = \can(p) \mrm f$
      implies $q \mrm d \red{}{\myr}{p}{\sm_p(p)} 0$ as $f \in
      C_{\can(p)}$ and hence $\can(p) \mrm f \in \s_p(p)$.
      On the other hand, if $\hterm(q) \not\id tc$ there exists a term $s \in \terms(\can(p))$ such that
      $\hterm(q) = s \mm b$ and $s \mm b \succ tc$.
      We have to distinguish two cases:
      In case $|s| < |ta|$ we know $s \mm b \id sb$, as $|s \mm b| = |tc|$.
      If $C(sb)$ is not empty let $be \myr f \in T$ be a corresponding
      rule.
      We get 
      $$ce \la \underline{ab}e \id a \underline{be} \myr af.$$
      As $c \neq a$ we either get $ b \mm e = \lambda$ in case $f =
      \lambda$ implying that $q \mrm e = (\can(p) \mrm b) \mrm e =
      \can(p) \mrm (b \mm e) = \can(p)$ and hence $q \mrm e
      \red{}{\myr}{p}{\sm_p(p)} 0$, or there exists an element $g \in \Sigma
      \cup \{ \lambda \}$ such that $ce \myr g$, $af \myr g \in T$, giving
      us $q \mrm e = (\can(p) \mrm b) \mrm e = \can(p) \mrm (b \mm e)
      = \can(p) \mrm f$ and thus $q \mrm e \red{}{\myr}{p}{\sm_p(p)} 0$
      as $f \in C_{\can(p)}$.
      On the other hand, if $|s| = |ta|$ with $s \id s'd$
      and $db \myr f \in T$, then $s
      \prec ta$ and $s \mm b \id s'f \succ tc$ implies $s'
      \id t$ and $f \succ c$.
      Now suppose $C(s'f) \neq \emptyset$ and let $fg \myr h \in T$ be a
      corresponding rule.
      Since $db \myr f \in T$ we also have $\inv{d}f \myr b \in T$,
      resulting in 
      $$bg \la \underline{\inv{d}f}g \id \inv{d}\underline{fg} \myr \inv{d}h.$$
      Since $g \neq h$  in case $h = \lambda$ we have $bg \myr \inv{d} \in T$
       giving us $cg \la \underline{ab}g=a\underline{bg} \myr a\,\inv{d}$.
      But then, as $a,c,g, \inv{d}$ all are not equal to $\lambda$, there
       exists $i \in \Sigma \cup \{ \lambda \}$ such that $cg \myr i$,
       $a\,\inv{d} \myr i \in T$, and thus $\inv{d} \in C_{\can(p)}$.
      This implies $q \mrm g = (\can(p) \mrm b) \mrm g = \can(p) \mrm \inv{d} \in \s_p(p)$.
      On the other hand, in case $h \neq \lambda$, there exists $i \in
       \Sigma \cup \{ \lambda \}$ such that $bg \myr i$, $\inv{d}h \myr i \in T$.
      Hence, $cg \la \underline{ab}g \id a\underline{bg} \myr ai$.
      In case $i = \lambda$, $bg \myr \lambda \in T$ immediately implies
       $q \mrm g = (\can(p) \mrm b) \mrm g = \can(p) \mrm (b \mm g) = \can(p)$.
      Otherwise there exists $j \in \Sigma \cup \{ \lambda \}$ such that 
       $cg \myr j$, $ai \myr j \in T$, and hence $i \in C_{\can(p)}$, giving us
       $q \mrm g = (\can(p) \mrm b) \mrm g = \can(p) \mrm (b \mm g) =
       \can(p) \mrm i \in \s_p(p)$.
      \\
      Hence in all these cases we have $q \mrm g \red{}{\myr}{p}{\sm_p(p)}
      0$.
\\
The case $q = \satpoly(p) \mrm b$ is similar in case
 $\hterm(\satpoly(p)) \id t' \inv{a}$.
Hence let us assume $\hterm(\satpoly(p)) = t' \mm\inv{a}\not\id t' \inv{a}$.
Then $t' \id t''k$, $t' \mm \inv{a} = t''l$ and $k\inv{a} \myr l \in T$.
The rule corresponding to $b \in C_{\satpoly(p)}$ then is $lb \myr c \in T$.
We have to distinguish the following two cases.
If $\hterm(q) \id tc$, then $C(tc) = \{ d \mid (cd, e) \in T, d \in \Sigma, e \in \Sigma
      \cup \{ \lambda \} \}$ and in case this set is
      not empty let us look at such a rule $(cd, e) \in T$.
      Since our presentation is a reduced  convergent group presentation, there
      exists a rule of the form $\inv{l}c \myr b \in T$ where $|\inv{l}| = 1$.
      Now this gives us 
      $$bd \la \underline{\inv{l}c}d \id
      \inv{l}\underline{cd} \myr \inv{l}e$$ and as $d \neq e$ and $b \neq
      \inv{l}$,
      there exists an element $f \in \Sigma \cup \{ \lambda \}$ 
      such that $bd \myr f$, $\inv{l}e \myr f \in T$.
      Again this results in the situation
      $$cd \la \underline{lb}d \id l\underline{bd} \myr lf$$
      and we either have $b \mm d = \lambda$ in case $f = \lambda$
      or there exists a rule $lf \myr e \in T$.
      In case $b \mm d = \lambda$ this implies $q \mrm d = (\satpoly(p) \mrm b) \mrm d =
      \satpoly(p) \mrm (b \mm d) =
      \satpoly(p)$ and hence $q \mrm d \red{}{\myr}{p}{\sm_p(p)} 0$.
      Otherwise,
      $q \mrm d = (\satpoly(p) \mrm b) \mrm d =\satpoly(p) \mrm (b \mm d) = \satpoly(p) \mrm f$
      implies $q \mrm d \red{}{\myr}{p}{\sm_p(p)} 0$ as $f \in
      C_{\satpoly(p)}$ and hence $\satpoly(p) \mrm f \in \s_p(p)$.
      On the other hand, if $\hterm(q) \neq tc$ there exists a term $s \in \terms(\satpoly(p))$ such that
      $\hterm(q) = s \mm b$ and $s \mm b \succ tc$.
      We have to distinguish two cases:
      In case $|s| < |tl|$ we know $s \mm b \id sb$, as $|s \mm b| = |tc|$.
      If $C(sb)$ is not empty let $be \myr f \in T$ be a corresponding
      rule.
      We get 
      $$ce \la \underline{lb}e \id l \underline{be} \myr lf.$$
      As $c \neq l$ we either get $ b \mm e = \lambda$ in case $f =
      \lambda$ implying that $q \mrm e = (\satpoly(p) \mrm b) \mrm e =
      \satpoly(p) \mrm (b \mm e) = \satpoly(p)$ and hence $q \mrm e
      \red{}{\myr}{p}{\sm_p(p)} 0$, or there exists an element $g \in \Sigma
      \cup \{ \lambda \}$ such that $ce \myr g$, $lf \myr g \in T$, giving
      us $q \mrm e = (\satpoly(p) \mrm b) \mrm e = \satpoly(p) \mrm (b \mm e)
      = \satpoly(p) \mrm f$ and thus $q \mrm e \red{}{\myr}{p}{\sm_p(p)} 0$
      as $f \in C_{\satpoly(p)}$.
      On the other hand, if $|s| = |tl|$ with $s \id s'd$ and $db \myr f \in T$ then $s
      \prec tl$ and $s \mm b \id s'f \succ tc$ implies $s'
      \id t$ and $f \succ c$.
      Now suppose $C(s'f) \neq \emptyset$ and let $fg \myr h \in T$ be a
      corresponding rule.
      Since $db \myr f \in T$ we also have $\inv{d}f \myr b \in T$,
      resulting in 
      $$bg \la \underline{\inv{d}f}g \id \inv{d}\underline{fg} \myr \inv{d}h.$$
      Since $g \neq h$  in case $h = \lambda$ we have $bg \myr \inv{d} \in T$
       giving us $cg \la \underline{lb}g=l\underline{bg} \myr l\,\inv{d}$.
      But then, as $l,c,g, \inv{d}$ all are not equal to $\lambda$, there
       exists $i \in \Sigma \cup \{ \lambda \}$ such that $cg \myr i$,
       $l\,\inv{d} \myr i \in T$, and thus $\inv{d} \in C_{\satpoly(p)}$.
      This implies $q \mrm g = (\satpoly(p) \mrm b) \mrm g = \satpoly(p) \mrm \inv{d} \in \s_p(p)$.
      On the other hand, in case $h \neq \lambda$, there exists $i \in
       \Sigma \cup \{ \lambda \}$ such that $bg \myr i$, $\inv{d}h \myr i \in T$.
      Hence, $cg \la \underline{lb}g \id l\underline{bg} \myr li$.
      In case $i = \lambda$, $bg \myr \lambda \in T$ immediately implies
       $q \mrm g = (\satpoly(p) \mrm b) \mrm g = \satpoly(p) \mrm (b \mm g) = \satpoly(p)$.
      Otherwise there exists $j \in \Sigma \cup \{ \lambda \}$ such that 
       $cg \myr j$, $li \myr j \in T$, and hence $i \in C_{\satpoly(p)}$, giving us
       $q \mrm g = (\satpoly(p) \mrm b) \mrm g = \satpoly(p) \mrm (b \mm g) =
       \satpoly(p) \mrm i \in \s_p(p)$.
      \\
      Hence in all these cases we have $q \mrm g \red{}{\myr}{p}{\sm_p(p)}
      0$.
\\
\qed
%
Notice that unlike in the case of free groups, the sets $\s_p$ as defined
in lemma \ref{lem.sat.in.pgr} need not be prefix Gr\"obner bases.
\begin{example}~\\
{\rm
Let $\Sigma = \{ a,b,c,d,b^{-1},c^{-1} \}$ and 
 $T = \{ a^2 \myr \lambda, d^2 \myr \lambda, bb^{-1} \myr \lambda, b^{-1}b
 \myr\lambda,\linebreak[4] cc^{-1} \myr \lambda, c^{-1}c \myr \lambda, ab \myr c, ac \myr b,
 c^{-1}b \myr d, c^{-1}a \myr b^{-1}, cb^{-1} \myr a,\linebreak[4] cd \myr b, db^{-1}\myr
 c^{-1}, dc^{-1} \myr b^{-1}, bc^{-1} \myr a, bd \myr c, b^{-1} a \myr c^{-1},
 b^{-1}c \myr d \}$ be a presentation of a plain group $\g$\footnote{This
 follows as the presentation is 2-monadic convergent and includes
 inverses of length 1 for all generators.} with  a length-lexicographical
 ordering induced by $a \succ b^{-1} \succ
 b \succ c^{-1} \succ c \succ d$.
\\
For the polynomial $p = ad + a + \lambda \in \q[\g]$ we get
$\sigma_1(p) = \lambda$, $\sigma_2(p) = d$,$\can(p) = p \mrm \lambda = ad + a + \lambda$
and $\satpoly(p) = p \mrm d = ad + a + d$.
In contrary to the case of free groups the set $\{ \can(p),
 \satpoly(p) \}$ alone is not prefix saturated and even not confluent
 since we have $\hterm(\can(p)) = \hterm(\satpoly(p))$.
By lemma
 \ref{lem.sat.in.pgr} we can give 
 a prefix saturating set for $p$, namely
 $\s_p(ad + a + \lambda) = \{ ad + a + \lambda, ad + a + d, ab^{-1} +
 ac^{-1} + b^{-1}, ab^{-1} + ac^{-1} + c^{-1} \}$ where the
 polynomials arise from the multiplications
 $(ad + d + \lambda) \mrm d = a + \underline{ad} + d$,
 $(ad + a + \lambda) \mrm b^{-1} = ac^{-1} + \underline{ab^{-1}} +
 b^{-1}$ and $(a + \underline{ad} + d) \mrm b^{-1} = ac^{-1} +
 \underline{ab^{-1}} + c^{-1}$.
Note that this set is no prefix Gr\"obner basis as for the prefix s-polynomial
 $$\spol{p}(ad + a + \lambda, ad + a + d) = d - \lambda$$
 we get that it is not prefix reducible by the polynomials in
 $\s_p(ad  + a + \lambda)$. 
\exaend
}
\end{example}
Specifying saturation of polynomials in procedure {\sc Reduced Prefix
  Gr\"obner Bases} on page \pageref{reduced.prefix.groebner.bases} we
  we can compute finite reduced prefix Gr\"obner bases in plain group rings.
\section{Context-free Groups}
As stated in the introductory chapter,
 a finitely generated  context-free group $\g$ is a group with a
 free normal subgroup of finite index.
Hence, let the group $\g$ be given by $X$  a finite set of generators
 for a free subgroup $\free$ and ${\cal E}$ a
 finite group such that $({\cal E}\backslash\{ \lambda \}) \cap (X \cup X^{-1}) =
 \emptyset$ and $\g/\free \cong {\cal E}$.
For all $e \in {\cal E}$ let $\phi_e : X \cup X^{-1} \myr \free$ be a function
 such that $\phi_{\lambda}$ is the inclusion and for all $x \in X \cup X^{-1}$,
 $\phi_e(x) = \inv{e} \mm_{\g} x \mm_{\g} e$.
For all $e_1,e_2 \in {\cal E}$ let $z_{e_1,e_2} \in \free$ such that
 $z_{e_1,\lambda} \id z_{\lambda,e_1} \id \lambda$ and for all
 $e_1,e_2,e_3 \in {\cal E}$
 with $e_1 \mm_{\cal E} e_2 =_{\cal E} e_3$, $e_1 \mm_{\g} e_2 \id
 e_3z_{e_1,e_2}$. 
Let $\Sigma = ({\cal E} \backslash \{ \lambda \}) \cup X \cup X^{-1}$
 and let $T$ contain the following rules: 
\begin{tabbing}
XX\=XXXX\=XXX\=XXXXXXX\= XXXXXXXXXXXXXXXXXXXXXXXX\= \kill
\>$xx^{-1}$ \> $\myr$ \> $\lambda$ \> and  \\
\>$x^{-1}x$ \> $\myr$ \> $\lambda$ \> for all $x \in X$,  \\
\>$e_1e_2$ \>  $\myr$ \>  $e_3z_{e_1,e_2}$     \> 
   for all $e_1,e_2 \in {\cal E} \backslash \{ \lambda \}, e_3
   \in {\cal E}$ such that $e_1 \mm_{{\cal E}} e_2 =_{\cal E} e_3$, \\
\>$xe$ \>  $\myr$ \>  $e \phi_e(x)$ \> and  \\
\>$x^{-1}e$ \>  $\myr$ \>  $e \phi_e(x^{-1})$ \> for all
        $e \in {\cal E} \backslash \{ \lambda \}, x \in X$.
\end{tabbing}
 $(\Sigma , T)$ then is a canonical and is called a virtually free
 presentation (compare \cite{CrOt94}).
\\
Presenting $\g$ in this way we find that the elements of the group are of
the form $eu$ where $e \in {\cal E}$ and $u \in \free$.
We can specify a total well-founded ordering on our group by combining a
 total well-founded ordering $\succeq_{\cal E}$ on ${\cal E}$ and a
 length-lexicographical 
 ordering $\geq_{\rm lex}$ on $\free$:
Let $w_1,w_2 \in \g$ such that $w_i \id e_iu_i$ where $e_i \in {\cal
E}$, $u_i \in \free$.
Then we define 
$w_1 \succ w_2$  if and only if   $|w_1| > |w_2|$ 
  or $(|w_1| = |w_2|$ and $e_1 \succ_{\cal E} e_2)$ 
  or $(|w_1| = |w_2|$ and $e_1 =_{\cal E} e_2$ and
             $u_1 >_{\rm lex} u_2)$.
This ordering is compatible with right concatenation using elements in
$\free$ in the following sense: Given
 $w_1, w_2 \in \g$ presented as described above, $w_1
 \succ w_2$ implies $w_1u \succ w_2u$ for all $u \in \free$ in case
 $w_1u, w_2u \in \g$.
\begin{example}\label{exa.cf.group}~\\
{\rm
Let ${\cal E}$ be the finite group presented by $\Sigma' = \{ a \}$
and $T' = \{ a^2 \myr \lambda \}$ and $\free$ the free group generated
by $X = \{ x \}$.
Further let $\phi_a(x) = x$ and $\phi_a(x^{-1}) = x^{-1}$ be a
conjugation homomorphism.
\\
Then 
$\Sigma = \{a, x, x^{-1} \}$ and
$T = \{ xx^{-1} \myr \lambda, x^{-1}x \myr \lambda  \} \cup \{ a^2 \myr
\lambda \} \cup \{ xa \myr ax, x^{-1}a \myr a x^{-1} \}$ is a virtually free presentation of $\g$.
\exaend
}
\end{example}
Let us take a closer look at prefix reduction in $\myk[\g]$.
\begin{example}~\\
{\rm
Let $\g$ be the group specified in example \ref{exa.cf.group}.
Further let $p = ax^2 + x + \lambda$, $q_1 = a + x$ and $q_2 =
x^2 + \lambda$ be polynomials in $\q[\g]$.
\\
Then the polynomial $p$ is prefix reducible at its head term $ax^2$ by
$q_1$ giving us $$p \red{}{\myr}{p}{q_1} p - q_1 \mrm x^2 = \underline{ax^2} + x +
\lambda - \underline{ax^2} - x^3 = x + \lambda + \underline{x^3}.$$
On the other hand, as $x^2$ is no prefix of $ax^2$, this is not true
for $q_2$.
\exaend
}
\end{example}
\begin{definition}~\\
{\rm
Let ${\cal H}$ be a subgroup of a group ${\cal N}$ and $p$ a
 non-zero polynomial in $\myk[{\cal N}]$.
A set $S \subseteq \{ \alpha \skm p \mrm w \mid \alpha \in \myk^*, w \in {\cal H} \}$ is called a 
 \index{saturating set!${\cal H}$-prefix}\betonen{${\cal H}$-prefix saturating set} for $p$,
 if for all $\alpha \in \myk^*$, $w \in {\cal H}$
 the polynomial $\alpha \skm p \mrm w$ is prefix reducible to zero using $S$ in
 one step.
A set of polynomials $F \subseteq \myk[{\cal N}]$ is called a
 \index{saturated set!${\cal H}$-prefix}\betonen{${\cal H}$-prefix saturated set}, if for all $f \in F$ and
 for all $\alpha \in \myk^*$, $w \in {\cal H}$ the polynomial $\alpha
 \skm f \mrm w$ is prefix
 reducible to zero using $F$ in one step.
\dend
}
\end{definition}
Reviewing the results on free groups, for a polynomial $p$ in $\myk[\g]$ we can specify $\can(p)$ and
$\satpoly(p)$ as follows.
\begin{definition}~\\
{\rm
For a non-zero polynomial $p \in \myk[\g]$
 we define
 $$\sigma_1(p) = \max \{ u \in \free \mid \inv{u} \mbox{ is a suffix of }
                       \hterm(p) \mbox{ and } \hterm(p \mrm u) = \hterm(p) \mm u
                       \},$$
In case $p$ contains more than one monomial or $\hterm(p) \neq e\inv{\sigma_1(p)}$ for $e \in {\cal E}$
we define
 $$\sigma_2(p) = \min \{ u \in \free \mid \inv{u} \mbox{ is a suffix of }
                       \hterm(p) \mbox{ and } \hterm(p \mrm u) \neq \hterm(p) \mm u
                       \}$$
and else $\sigma_2(p) = \sigma_1(p)$.
Then  we can set
$\can(p) =  p \mrm \sigma_1(p)$ and
$\satpoly(p)=  p \mrm
 \sigma_2(p)$.
\dend
}
\end{definition}
\begin{lemma}~\\
{\sl
Let $p \in \myk[\g]$ be a non-zero polynomial.
Then the set $\{ \can(p), \satpoly(p) \}$,
 is a $\free$-prefix saturating set for $p$.
\lemend
}
\end{lemma}
\Ba{}~\\
The proof is straightforward as in lemma \ref{lem.freegroupsat}.
We only have to consider the additional case $\hterm(p \mrm
\sigma_1(p)) \id e \in {\cal E}$.
Since for all $w \in
\free$ we get $\hterm(p \mrm w) = \hterm(p) \mrm w \id ew'$ for some
$w' \in \free$, then $\{ \can(p) = \satpoly(p) = p \mrm
\sigma_1(p) \}$ is a prefix saturating set.
\\
\qed
\begin{example}~\\
{\rm
Let $\g$ be the group specified in example \ref{exa.cf.group}.
Then for the polynomial $p =
ax^2 + x + \lambda$ in $\q[\g]$ we get $\sigma_1(p) = x^{-1}$ and $\sigma_2(p) = x^{-2}$ giving us
$\can(p) = p \mrm \sigma_1(p) = \underline{ax} +  \lambda + x^{-1}$
and $\satpoly(p) = p \mrm\sigma_2(p) = a + x^{-1} + \underline{x^{-2}}$.
\\
On the other hand $q = a+x$ gives us
 $\sigma_1(q) = \sigma_2(q) = \lambda$ and thus $\can(q) = \satpoly(q) = q$.
\mbox{\phantom{X}}\exaend
}
\end{example}
The following lemma will be used as an analogon to lemma
\ref{lem.redp} when we characterize prefix Gr\"obner bases by using
prefix reduction, prefix s-polynomials and now $\free$-prefix saturated
sets.
\begin{lemma}\label{lem.cf}~\\
{\sl
Let $p$ be a non-zero polynomial and $F$ a set of polynomials in
$\myk[\g]$.
Then $p \red{*}{\myr}{p}{F} 0$ gives us a prefix standard representation
 of $p = \sum_{i=1}^k \alpha_i \skm f_i \mrm w_i$, with $\alpha_i \in \myk^*, f_i
 \in F, w_i \in \g$ such that for all $w \in
 \free$ with $\hterm(p \mrm w) \id \hterm(p)w$, we get $ \hterm(p)w \succeq
 \hterm(f_i \mrm w_i \mrm w)$.
Note that additionally for all $t \in \m$ with $t \succeq \hterm(p)$,
 if $t \mm w \id tw$ for some $w \in \m$, then $tw \succeq
 \hterm(f_i \mrm w_i \mrm w)$.
\lemend\ohnebeweis
}
\end{lemma}
For every $e \in {\cal E}$ let the mapping 
 $\psi_e: \myk[\g] \myr \myk[\g]$ be defined
 by $\psi_e(f) = f \mrm e$ for $f \in \myk[\g]$.
We now can give a characterization of  prefix Gr\"obner bases by
 transforming a generating set for a right ideal using these finitely many
 mappings.
This will enable us to restrict ourselves to $\free$-prefix saturated
sets when characterizing prefix Gr\"obner bases.
\begin{theorem}\label{theo.compcf}~\\
{\sl
Let $
 F \subseteq \myk[\g]$ and $G \subseteq \myk[\g]$ such that
 \begin{itemize}
 \item $\ideal{r}{}(F) = \ideal{r}{}(G)$
 \item $\{ \psi_e(f) \mid f \in F, e \in {\cal E} \} \subseteq G$
 \item $G$ is $\free$-prefix saturated.
 \end{itemize}.
Then the following statements are equivalent:
\begin{enumerate}
\item For all $g \in \ideal{r}{}(F)$ we have $g \red{*}{\myr}{p}{G} 0$.
\item For all $f_{k}, f_{l} \in G$ we have 
  $\spol{p}(f_{k}, f_{l}) \red{*}{\myr}{p}{G} 0$.
\theoend
\end{enumerate}
}
\end{theorem}
\Ba{}~\\
\mbox{$1 \R 2:$ }
Let $\hterm(f_{k}) \id \hterm(f_{l})w$ for $f_k,f_l \in G$ and
 $w \in \g$.
Then by definition \ref{def.cpp} we get
 $$\spol{p}(f_{k}, f_{l}) = \hc(f_k)^{-1} \skm f_{k} 
    -\hc(f_l)^{-1} f_{l} \mrm w \:\in \ideal{r}{}(G)= \ideal{r}{}(F),$$
 and hence $\spol{p}(f_{k}, f_{l}) \red{*}{\myr}{p}{G} 0$.

\mbox{$2 \R 1:$ }
We have to show that every non-zero element $g \in \ideal{r}{}(F)$
 is $\red{}{\myr}{p}{G}$-reducible to zero.
Remember that for
 $h \in \ideal{r}{}(F)= \ideal{r}{}(G)$, $ h \red{}{\myr}{p}{G} h'$
 implies $h' \!\in \ideal{r}{}(G)= \ideal{r}{}(F)$.
Thus as  $\red{}{\myr}{p}{G}$ is Noetherian
 it suffices to show that every 
 $g \in \ideal{r}{}(F)\backslash\{ 0 \}$ is $\red{}{\myr}{p}{G}$-reducible.
Let $g = \sum_{j=1}^m \alpha_{j} \skm f_{j} \mrm w_{j}$ be a
  representation of a non-zero polynomial $g$ such that
  $\alpha_{j} \in \myk^*, f_j \in F, w_{j} \in \g$.
Further for all $1 \leq j \leq m$, let $w_j \id e_ju_j$, with
 $e_j \in {\cal E}$, $u_j \in \free$.
Then, we can modify our representation of $g$ to
 $g = \sum_{j=1}^m \alpha_j \skm \psi_{e_j}(f_j) \mrm u_j$.
Since $G$ is $\free$-prefix saturated and
 $\psi_{e_j}(f_j) \in G$ we can assume $g= \sum_{j=1}^m \alpha_j \skm
 g_j \mrm v_j$, where $\alpha_j \in \myk^*, g_j \in G, v_j \in \free$
 and  $\hterm(g_{j} \mrm v_{j}) \id \hterm(g_{j})v_{j}$.
Depending on this representation of $g$ and our well-founded
 total ordering $\succeq$ on $\g$ we define
 $t = \max \{ \hterm(g_{j})v_{j} \mid j \in \{ 1, \ldots m \}  \}$ and
 $K$ is the number of polynomials $g_j \mrm v_j$ containing $t$ as a term.
Then $t \succeq \hterm(g)$ and
in case $\hterm(g) = t$ this immediately implies that $g$ is
$\red{}{\myr}{p}{G}$-reducible. 
So
by lemma \ref{lem.srprop} it is sufficient to  show that
$g$ has a prefix standard representation, as this implies that $g$ is
top-reducible using $G$.
This will be done by induction on $(t,K)$, where
 $(t',K')<(t,K)$ if and only if $t' \prec t$ or
 $(t'=t$ and $K'<K)$\footnote{Note that this ordering is well-founded
                                    since $\succeq$ is and $K \in\n$.}.
If $t \succ \hterm(g)$
 there are two polynomials $g_k,g_l$ in the corresponding
 representation\footnote{Not necessarily $g_l \neq g_k$.}
 and  $t \id \hterm(g_k)v_k \id \hterm(g_l)v_l$.
Without loss of generality let us assume $\hterm(g_k)\id\hterm(g_l)z$ for
 some $z \in \free$ and  $v_l \id zv_k$.
Then by definition \ref{def.cpp} we have a prefix s-polynomial
 $\spol{p}(g_k,g_l) = \hc(g_k)^{-1} \skm  g_k -
 \hc(g_l)^{-1} \skm g_l \mrm z$.
We will now change our representation of $g$ by using the additional
information on this s-polynomial in such a way that for the new
representation of $g$ we either have a smaller maximal term or the occurrences of the term $t$
are decreased by at least 1.
Let us assume  $\spol{p}(g_k,g_l) \neq 0$\footnote{In case  $\spol{p}(g_k,g_l) = 0$,
 just substitute $0$ for $\sum_{i=1}^n \delta_i \skm h_i \mrm v'_i$ in 
 the equations below.}.
Hence,  the reduction sequence $\spol{p}(g_k,g_l) \red{*}{\myr}{p}{G} 0
$ yields a prefix standard representation of the form
 $\spol{p}(g_k,g_l) =\sum_{i=1}^n \delta_i \skm h_i \mrm v'_i$, $\delta_i \in
 \myk^*$,$h_i \in G$,$v'_i \in \free$
 and all terms occurring in the sum are bounded by $\hterm(\spol{p}(g_k,g_l))$.
By lemma \ref{lem.cf} we can conclude that $t$ is a proper
 bound for all terms occurring
 in the sum $\sum_{i=1}^n \delta_i \skm h_i \mrm v'_i \mrm v_k$ and
 again we can substitute all polynomials $h_i$, where
 $\hterm(h_i \mrm v'_i \mrm v_k) \neq \hterm(h_i)(v'_i \mm v_k)$  without increasing
 $t$ or $K$.
Similarly, in case $v_i' \in {\cal E}$, we can substitute $h_i$ by
 $\psi_{v_i'}(h_i) \in G$ by our assumption.
Therefore, without loss of generality we can assume that the
 representation has the required form.
This gives us: 
\begin{eqnarray}
 &  & \alpha_{k} \skm g_{k} \mrm v_{k} + \alpha_{l} \skm g_{l} \mrm v_{l}  \nonumber\\ 
 &  &  \nonumber \\
 & = &  \alpha_{k} \skm g_{k} \mrm v_{k} +
        \underbrace{ \alpha'_{l} \skm \beta_k \skm g_{k} \mrm v_{k}
                   - \alpha'_{l} \skm \beta_k \skm g_{k} \mrm
                   v_{k}}_{=\, 0\phantom{\spol{p}(g_k,g_l) \mrm v_k}} 
                   + \alpha'_{l}\skm \beta_l  \skm g_{l} \mrm v_{l} \nonumber\\ 
 & = & (\alpha_{k} + \alpha'_{l} \skm \beta_k) \skm g_{k} \mrm v_{k} - \alpha'_{l} \skm
        \underbrace{(\beta_k \skm g_{k} \mrm v_{k}
        -  \beta_l \skm g_{l} \mrm v_{l})}_{=\, \spol{p}(g_k,g_l) \mrm v_k}\nonumber\\
 & = & (\alpha_{k} + \alpha'_{l} \skm \beta_k) \skm g_{k} \mrm v_{k} - \alpha'_{l} \skm
       (\sum_{i=1}^n \delta_{i} \skm h_{i} \mrm v'_{i} \mrm v_k) \label{s5}
\end{eqnarray}
 where  $\beta_k = \hc(g_k)^{-1}$, $\beta_l = \hc(g_l)^{-1}$ and
 $\alpha'_l \skm \beta_l = \alpha_l$.
By substituting (\ref{s5}) in our representation of $g$ 
 either $t$ disappears  or in
 case $t$ remains maximal among the terms occurring in the new
 representation of $g$, $K$ is decreased.
\\
\qed
Next we give a procedure to compute reduced prefix Gr\"obner bases by
modifying procedure {\sc Reduced Prefix Gr\"obner Bases} on page \pageref{reduced.prefix.groebner.bases}.

\procedure{Reduced Gr\"obner Bases in Context-Free Group Rings}%
{\vspace{-4mm}\begin{tabbing}
XXXXX\=XXXX \kill
\removelastskip
{\bf Given:} \> $F \subseteq \myk[\g]$, 
              $(\Sigma , T)$  a virtually free presentation of
             $\g$. \\
{\bf Find:} \> $\gb(F)$, a  (prefix) Gr\"obner basis of $F$.
\end{tabbing}
\vspace{-7mm}
\begin{tabbing}
XX\=XX\= XXXX \= XX\=XXXX \=\kill
$G_0$ := $\emptyset$; \\
$S_0$ := $\{ \can(\psi_e(f)), \satpoly(\psi_e(f)) \mid e \in {\cal E}, f \in F \}$; \\
$i$ := $0$; \\
{\bf while} $S_i \neq \emptyset$ {\bf do} \hspace{1cm}{\rm\kommentar \%
  $\ideal{r}{}(F) = \ideal{r}{}(G_i \cup S_i)$}\\
\> $i$ := $i+1$; \\
\> $q_i$ := {\rm remove}$(S_{i-1})$; \\
\>{\rm\kommentar \%   Remove an element using a fair strategy}\\
\> $q_i'$ := ${\rm normalform}(q_i, \red{}{\myr}{p}{G_{i-1}})$; \\
\> {\rm\kommentar \% Compute a normal form using
           prefix reduction} \\
\>{\bf if} \>$q_i' \neq 0$ \\
\>         \>{\bf then} \>{\bf if} \>$|\terms(q_i')| = 1$ \\
\>         \>           \>         \> {\rm \kommentar \% The right
  ideal is trivial} \\
\>         \>           \>         \>{\bf then}\>$G_i$ := $\{ \lambda \}$; \\
\>         \>           \>         \>          \>$S_i$ := $\emptyset$; \\
\>         \>           \>         \>{\bf else} \\
\>         \>           \>         \>            \>$H_i$ := $\{ g \in G_i \mid \hterm(g) \mbox{ is prefix reducible using } q_i' \}$;\\
\>         \>            \>        \>            \>{\rm \kommentar \% These polynomials would
  have new head terms when prefix}\\
\>         \>            \>        \>            \>{\rm \kommentar \% reduced using $q_i'$} \\
\>         \>           \>         \>            \>$G_i$  := $(G_{i-1} \backslash H_i) \cup \{ q_i' \}$; \\
\>         \>           \>         \>            \>$S_i$ := $S_{i-1} \cup H_i \cup \{\can(q_i'),  \satpoly(q_i') \}$; \\
\>         \> \>{\bf endif}\\
\>{\bf endif}\\
{\bf endwhile} \\
$\gb (F)$:= ${\rm reduce}(G_i)$ \\
{\rm\kommentar
  \% reduce$(F) = \{ {\rm normalform}(f, \red{}{\myr}{p}{F \backslash \{
    f \}}) \mid f \in F \}$\footnotemark}
\end{tabbing}}
\footnotetext{Notice that only the reducts are touched in this procedure.}

\begin{remark}~\\
{\rm
For all sets computed we have $\ideal{r}{}(F) = \ideal{r}{}(G_i \cup
S_i) = \ideal{r}{}(G)$.
But note that the set $G$ need not fulfill the condition that 
 $\{ \psi_e(f) \mid e \in {\cal E}, f \in F \} \subseteq G$.
Hence theorem \ref{theo.compcf} cannot be applied to ensure that $G$ is a
 prefix Gr\"obner basis.
We will see later on that $\{ \can(\psi_a(f)) \mid a \in {\cal E}, f \in
F \} \subseteq G_0$ in fact is sufficient to ensure correctness. 
\remend
}
\end{remark}
Termination of the procedure follows by the same arguments used for  plain groups.
\begin{theorem}~\\
{\sl
Procedure {\sc Reduced Gr\"obner Bases in Context-Free Group Rings} terminates on
 finite input $F$.
\theoend
}
\end{theorem}
\Ba{}~\\
The proof is done in two steps:
 first we show that all polynomials computed have a certain property
 that will be used in the 
 second step to ensure termination.
 We say a polynomial $q$ has property $\p_F$ if and only if
      \begin{enumerate}
      \item[($\alpha$)] $|\hterm(q)| \leq K$, where \\
            $K = \max \{ |\hterm(f)| \mid f \in 
                 \{ \can(\psi_e(f)), \satpoly(\psi_e(f)) \mid e \in {\cal E}, f \in F \}\} +1$.
      \item[($\beta$)] If $|\hterm(q)|=K$ then there exists an element
            $a \in X \cup X^{-1}$ such that all terms of length $K$
            in $q$ have $a$ as a common suffix.
      \end{enumerate}
      We will show that all polynomials $q$ computed by the
       procedure on input $F$ have property $\p_F$. \\
      By the choice of $K$ all  polynomials in $G_0$ and $S_0$ 
       have $\p_F$.
      Let $G_i$ and $S_i$ be actual computed sets of polynomials
       having $\p_F$ and 
       let $q_i$ be the next polynomial chosen by our procedure.
      Then $q_i$ is first prefix reduced to normal form with respect to
       $G_i$ and only  polynomials having $\p_F$ are involved.
      By the definition of prefix reduction this operation preserves 
       property $\p_F$.
\\
      It remains to show that computing the polynomials $\can(p)$ and
      $\satpoly(p)$ of a polynomial having $\p_F$ does not destroy this property.
      In case $q = \alpha \skm t$ or $\hterm(q \mrm \sigma_1(q)) \in
      {\cal E}$ we are done.
      Hence let us assume $\hterm(q \mrm \sigma_1(q)) \id eu$ for some
      $e \in {\cal E}$, $u \in \free$ and let $|eu| = K$.
      Then $\sigma_1(p) = \lambda$ as $|\hterm(q \mrm \sigma_1(q))|
      \leq |\hterm(q)|$ and $q$ has $\p_F$.
      Furthermore, 
      all $t \in \terms(q)$ of length $K$ have a common last letter,
      say $a \in X \cup X^{-1}$.
      Since $\sigma_2(p) = \inv{a}$, the head term of $q \mrm
      \inv{a}$ must again have length less equal to $K$ and all terms
      of length $K$ must have $\inv{a}$ as last letter.
\\
It remains to show that the procedure terminates.
Let us assume the contrary.
       Then there are infinitely many polynomials $q_i', i\in\n$, with
        heads in normal form added. 
\\
       But since no term occurs more than once among the head terms of
       polynomials added to a set $G_i$, $i>0$, and for every such polynomial
       we have  $|\hterm(q_i')| \leq K$ this is not possible.
\\
\qed
\begin{theorem}~\\
{\sl
Procedure {\sc Reduced Gr\"obner Bases in Context-Free Group Rings}
 is correct.
\theoend
}
\end{theorem}
\Ba{}~\\
Let $G$ be the output of procedure {\sc Reduced Gr\"obner Bases in
  Context-Free Group Rings} on input $F$.
Without loss of generality let us assume that $\ideal{r}{}(F)$ is not
trivial, i.e., it is neither $\{ 0 \}$ nor $\myk[\g]$.
Then by theorem \ref{theo.pgb.psr} it is sufficient to show that every
non-zero polynomial $g$ in $\ideal{r}{}(F)$ has a prefix
standard representation with respect to $G$.
This will be done by transforming an arbitrary representation of $g$
with respect to $F$ into a prefix standard representation with respect
to $G$.
Let $g = \sum_{i=1}^{k} \alpha_i \skm f_i \skm w_i$ with $\alpha_i \in
\myk^*$, $f_i \in F$, $w_i \in \g$ and $w_i \id e_iu_i$ for $e_i \in
{\cal E}$, $u_i \in \free$.
Then since the set $\{ \can(\psi_e(f)) | e \in {\cal E}, f \in F \}
\subseteq G_0 \cup S_0$ we can represent $g$ as $g = \sum _{i=1}^k \alpha_i \skm
\psi_{e_i}(f_i) \mrm u_i$.
Furthermore, as the set $G_0 \cup S_0$ is $\free$-prefix saturated we even get
a representation $g = \sum_{i=1}^{m'} \beta_i \skm g_i \mrm v_i$ with
$\beta_i \in \myk^*$, $g_i \in G_0 \cup S_0$ and $v_i \in
\free$\footnote{This is due to the fact that $\ideal{r}{}(F)$ does not
  contain the identity.} such
that $\hterm(g_i \mrm v_i) \id \hterm(g_i)v_i$.
In showing that every $f \in G_0 \cup S_0$ has a prefix standard
representation with respect to $G$ by lemma \ref{lem.psr.prop3} then a
prefix standard representation for $g$ with respect to $G$ also
exists.
Remember that lemma \ref{lem.psr.prop3} states that in case a
polynomial has a prefix standard representation with respect to a set
and all polynomials in this set have prefix standard representations
with respect to another set, the polynomial itself again has a prefix
standard representation with respect to the latter set.
Let us proceed now in showing our last claim.
Since our procedure terminates there exists an index $k \in \n$ such
that $G = {\rm reduce}(G_k)$ and $S_k = \emptyset$.
We will now prove that our claim holds for every $f \in G_i \cup S_i$
by induction on $j$ where $i = k-j$.
In case $j = 0$ we are immediately done since by lemma
\ref{lem.psr.prop3} $f \in G_k$ has a prefix
standard representation with respect to $G$.
Hence let $f \in G_{k - (j+1)} \cup S_{k- (j+1)}$ and suppose $f \not
\in G_{k-j} \cup S_{k-j}$, as then our induction hypothesis can be
applied and therefore the claim holds.
Now if $f \in G_{k-(j+1)}$, $\hterm(f)$ must be prefix reducible using the
polynomial $q_{k-j}'$ computed in this iteration, as we assume $f \not
\in G_{k-j}$.
Therefore we get $f \in H_{k-j}$
and hence $f \in S_{k-j}$ and we are done.
It remains to study the case that $f \in S_{k-(j+1)}$.
Since $f \not \in S_{k-j}$, $f$ is chosen to compute the polynomial
$q_{k-j}'$ and again we can conclude that $f$ has a prefix standard
representation with respect to $G_{k-j} \cup S_{k-j}$ and hence with
respect to $G$.
\\
\qed

\auskommentieren{Note that this approach can be applied to every group that can be presented by a subgroup
 of finite index and in case the subgroup ring allows finite prefix Gr\"obner bases this carries over
 to the group ring.}

%% file: nilpotent.tex
\section{Nilpotent Groups}
Nilpotent groups and their presentations were briefly described
 in section \ref{section.presentations}.
The fact that  their elements can be presented by ordered group words
motivates an approach similar to the one for commutative monoids.
Since multiplication is no longer commutative, we will first restrict
 ourselves to right ideals and show
 the existence of finite right Gr\"obner bases for finitely generated
 right ideals.
Later on this approach will be extended to two-sided ideals.
Let us start by generalizing the concept of special divisors which
 can be interpreted as (commutative) ``prefixes'' in the set of ordered
group words as in the commutative case.
This will be done by extending the tuple ordering on ordered words (compare
definition \ref{def.tuple}) to a tuple ordering on the set of ordered
group words $\ord(\Sigma) = \{ a_1^{i_1} \ldots a_{n\phantom{1}}^{i_n}
| i_j \in \z \}$.
\begin{definition}~\\
{\rm
For $w \id a_1^{i_1} \ldots a_n^{i_n}, v \id a_1^{j_1}  \ldots
 a_n^{j_n} \in \ord(\Sigma)$, we define $w \tupeq v$ if for each $1 \leq l \leq n$ we
       have either $j_l = 0$ or $\sgn(i_l)=\sgn(j_l)$ and $|i_l| \geq
       |j_l|$.
      Further we define $w \tupgreater v$ if $w \tupeq v$ and $|i_l|>|j_l|$
       for some $1 \leq l \leq n$.
\dend
}
\end{definition}
According to this ordering we call $v$ a prefix of $w$ if $v \tupleq
w$.
Notice that an element then has finitely many prefixes.

Let us now start by investigating  the special case of torsion-free
 nilpotent groups.
As seen in section \ref{section.presentations}, such a group can be presented by
 a convergent CNI-system over an alphabet
 $\Sigma = \{a_1,a_1^{-1}, \ldots, a_n, a_n^{-1} \}$, and we show that for
 such presentations (which contain no P-rules)
 additional syntactical lemmata hold which enable a
 weakening of reduction comparable to commutative reduction.
\begin{example}\label{exa.fnp2CNI}~\\
{\rm
The free nilpotent group of class 2 with 2 generators
described in example \ref{exa.fnp2} can be presented by the
convergent CNI-system $\Sigma = \{ a_1, a_1^{-1}, a_2, a_2^{-1}, a_3, a_3^{-1} \}$ and $T
= \{ a_2a_1 \myr a_1a_2a_3, a_2^{-1}a_1^{-1} \myr a_1^{-1}a_2^{-1}a_3,
a_2^{-1}a_1 \myr a_1a_2^{-1}a_3^{-1}, a_2a_1^{-1} \myr
a_1^{-1}a_2a_3^{-1}, \linebreak a_3^{\delta}a_2^{\delta'} \myr
a_2^{\delta'}a_3^{\delta}, a_3^{\delta}a_1^{\delta'} \myr
a_1^{\delta'}a_3^{\delta} \mid \delta, \delta' \in \{ 1, -1 \}\}$.
\exaend
}
\end{example}
The ordering on our group $\g$ will be the syllable ordering  and for $w,v
\in \g$ we have that $w \tupgreater v$ implies $w \syll v$, but the converse is not true,
as $\tupeq$ is not total.
Reviewing example \ref{exa.fnp2CNI} we find $a_1a_2 \tupgreater a_1$ and
$a_1a_2 \syll a_1$, but $a_1^{-1}a_2 \syll a_1$ and $a_1^{-1}a_2 \not\tupgreater a_1$.
The following lemma reveals a connection between these two orderings
in  nilpotent groups having convergent CNI-presentations which
specifies under which conditions a term remains a proper bound for
$\syll$-smaller terms under restricted right multiplication.
\begin{lemma}\label{lem.order}~\\
{\sl
Let $\g$ be a finitely generated  nilpotent group
 presented by a convergent CNI-system on 
 $\Sigma =  \{ a_1,a_1^{-1}, \ldots, a_n, a_n^{-1} \}$.
Further let $w,v,\tilde{v} \in \g$ with $w \tupeq v$ and
 $v \syll \tilde{v}$.
Then for $u \in \g$ such that $w = v \mm u$, we get $w
\syll \tilde{v}\mm u$. 
Notice that since $\g$ is a group such an element $u$ always exists,
namely $u = \inv{v} \mm w$, and $u$ is unique.
\lemend
}
\end{lemma}
\Ba{}~\\
Let $w,v, \tilde{v},u \in \g$ be presented by ordered group words,
i.e., 
 $w \id a_1^{w_1} \ldots a_n^{w_n}$, 
 $v \id a_1^{v_1} \ldots a_n^{v_n}$,
 $\tilde{v} \id a_1^{\tilde{v}_1} \ldots a_n^{\tilde{v}_n}$,  and
 $u \id a_1^{u_1} \ldots a_n^{u_n}$
 with $w_i,v_i,\tilde{v}_i,u_i \in \z$.
Further let $a_d$ be the  distinguishing letter
 between $v$ and $\tilde{v}$, i.e.,
 $v_d >_{\z} \tilde{v}_d$.
Since the commutation system only includes rules of the form
 $a_j^\delta a_i^{\delta'} \myr  a_i^{\delta'} a_j^\delta z$, $j>i$,
 $z \in \ord(\Sigma_{j+1}), \delta, \delta' \in \{ 1, -1 \}$
 and we have no P-rules, we can conclude 
$$a_1^{v_1} \ldots a_{d-1}^{v_{d-1}} \mm a_1^{u_1} \ldots a_{d-1}^{u_{d-1}}
 = a_1^{\tilde{v}_1} \ldots a_{d-1}^{\tilde{v}_{d-1}} \mm a_1^{u_1}
 \ldots a_{d-1}^{u_{d-1}} \id a_1^{w_1} \ldots a_{d-1}^{w_{d-1}} a_d^{s_d} \ldots a_n^{s_n}$$
for some $s_i \in \z$.
Moreover, $a_d^{v_d} \mm a_d^{s_d} \mm a_d^{u_d} = a_d^{w_d}$, i.e.,
 $v_d + s_d + u_d = w_d$.
To prove $w_d >_{\z} \tilde{v}_d + s_d +u_d$ and hence $w \syll
\tilde{v} \mm u$, we have to
 take a closer look at $v_d$ and $\tilde{v}_d$.
\begin{enumerate}
\item In case $v_d >0$ this implies $w_d >0$ as $w \tupeq v$.
      Therefore, $v_d + s_d + u_d = w_d$ and $w_d \geq v_d >0$
       give us $s_d + u_d \geq 0$.
     Now $v \syll \tilde{v}$ and $v_d >0$ imply that
      $v_d > \tilde{v}_d \geq 0$, as otherwise $\tilde{v}_d \geq_{\z} v_d$
      would contradict our assumption.
     Hence we get $\tilde{v}_d + s_d + u_d < w_d$, implying
      $w \syll \tilde{v} \mm u$.
\item In case $v_d < 0$ this implies $w_d < 0$, $|w_d| \geq |v_d|$ and
      thus $v_d + s_d + u_d = w_d$ yields $s_d + u_d \leq 0$.
     Further we know $|v_d|+|s_d+u_d| = |w_d|$.
     We have to distinguish two cases:
  \begin{enumerate}
  \item In case $\tilde{v}_d \leq 0$, then $v \syll \tilde{v}$
         implies $|v_d| > |\tilde{v}_d|$.
       Therefore, we get $ |\tilde{v}_d|+|s_d+u_d| < |w_d|$ and
        $w \syll \tilde{v} \mm u$.
  \item In case $\tilde{v}_d > 0$, as $s_d + u_d \leq 0$ we
         have to take a closer look at $\tilde{v}_d +s_d+u_d$.
        In case  $\tilde{v}_d +s_d+u_d \geq 0$ we are done as
         this implies $w \syll  \tilde{v} \mm u$.
        In case  $\tilde{v}_d +s_d+u_d < 0$ we get that
         $\tilde{v}_d < |s_d + u_d|$
         implying  
         $|\tilde{v}_d+s_d+u_d|< |s_d+u_d| < |w_d|$
         and hence $w \syll  \tilde{v} \mm u$.
\\
\qed
  \end{enumerate}
\end{enumerate}
\begin{remark}\label{rem.multCNI}~\\
{\rm
In the previous proof it is very important that the modifications due to
 the changes of the
 occurrences of the letters $a_1, \ldots, a_{d-1}$ are the same in $v$
 and $\tilde{v}$.
Further all changes of the occurrence of the letter $a_d$ due to
 those modifications are actually the same and are expressed in the
 exponent $s_d$.
This is true since ``moving'' a letter with a
 smaller index past $a_d$ can only add a word in $\ord(\Sigma_{d+1})$
 and moving a newly introduced letter with an index $k > d$
 back past the letter $a_d$ can only add words in $\ord(\Sigma_{k+1})$
 (compare lemma \ref{lem.nilpotentmultiplication}).
Note that this argumentation need no longer hold for arbitrary
 CR-rules $a_j^\delta a_i^{\delta'} \myr  a_i^{\delta'}z$,
      $z \in \ord(\Sigma_{i+1})$, as a rule $a_da_i \myr a_iz$ can have
 additional influence on the occurrence of the distinguishing letter $a_d$.
\remend
}
\end{remark}
Notice that the observations of this remark on the behaviour  
of CNI-presentations are
essential and will be frequently used in the proofs of this section.

Henceforth, let $\g$ be a torsion-free nilpotent group with a
 convergent CNI-presentation $(\Sigma, T)$.
We are using the syllable ordering to induce the ordering on the group
 ring $\myk[\g]$ and the tuple ordering to restrict this ordering.
This is similar to the concepts of prefix and commutative reduction.
Note that these reductions correspond to the property described in
lemma \ref{lem.order} as follows:
If $\m$ is a monoid and $w,v, \tilde{v} \in \m$ such that $v$ is a prefix
of $w$ as a word, i.e., $w \id vu$ for some $u \in
\m$, then $v \succ \tilde{v}$ implies $vu \succ \tilde{v}u \succeq
\tilde{v} \mm u$.
On the other hand, if $\m$ is a commutative monoid, $w \tupeq v$ and $w = v \cm u$, then similarly
$v \succ \tilde{v}$ implies $w = v \cm u \succ \tilde{v} \cm u \succeq
\tilde{v} \mm u$.
In chapter \ref{chapter.reduction} these properties were used
 to characterize appropriate Gr\"obner bases by special
 s-polynomials.
We will now give a similar approach by introducing special standard
 representations and a reduction
corresponding to ``prefixes'' in the set of ordered group words
defined by the extended tuple-ordering.
\begin{definition}\label{def.qcsr}~\\
{\rm 
Let $F$ be a set of polynomials  and $p$ a non-zero
polynomial in $\myk[\g]$.
A representation 
$$ p = \sum_{i=1}^{n} \alpha_i \skm f_{i} \mrm w_i, \;\;
  \mbox{ with } \alpha_i \in \myk^*,
  f_{i} \in F, w_i \in \g $$
  is called a \betonen{quasi-commutative (qc-)standard representation}\/ 
 in case for all $1 \leq i \leq n$ we have $\hterm(p)   
 \succeq \hterm(f_{i}) \mm w_i \succeq \hterm(f_i \mrm w_i)$
 and $\hterm(f_i \mrm w_i) \tupeq \hterm(f_i)$.
\dend
}
\end{definition} 
\begin{definition}\label{def.redqc}~\\
{\rm
Let $p, f$ be two non-zero polynomials  in $\myk[\g]$. 
We say $f$ 
 \index{quasi-commutative!right reduction}\index{reduction!quasi-commutative right}\betonen{quasi-commutatively (right) reduces} $p$ to $q$ at
 a monomial $\alpha \skm t$ of $p$ in one step, denoted by $p \red{}{\myr}{qc}{f} q$, if
\begin{enumerate}
\item[(a)] $t \tupeq \hterm(f)$, and
\item[(b)] $q = p - \alpha \skm \hc(f)^{-1} \skm f \mrm
  (\inv{\hterm(f)} \mm t)$.
\end{enumerate}
We write $p \red{}{\myr}{qc}{f}$ if there is a polynomial $q$ as defined
above  and $p$ is then called  quasi-commutatively reducible by $f$. 
Further we can define $\red{*}{\myr}{qc}{}, \red{+}{\myr}{qc}{}$,
 $\red{n}{\myr}{qc}{}$ as usual.
Quasi-commutative reduction by a set $F \subseteq \myk[\g]$ is denoted by
 $p \red{}{\myr}{qc}{F} q$ and abbreviates $p \red{}{\myr}{qc}{f} q$
 for some $f \in F$,
 which is also written as  $p \red{}{\myr}{qc}{f \in F} q$.
\dend
}
\end{definition}
Notice that if $f$ quasi-commutatively reduces $p$ at $\alpha \skm t$
to $q$,
then  $t \not\in \terms(q)$.
This reduction is effective, as it is possible to decide, whether
$t \tupeq \hterm(f)$.
Further it is Noetherian, as by lemma \ref{lem.order}, $p \red{}{\myr}{qc}{f} q$ implies $p > q$. 

For a commutative group $\g$, quasi-commutative right reduction
 and commutative reduction  coincide as follows:  For $w,v \in
 \g$, $w \tupeq v$  implies $v \cm u = w$ where $u = \inv{v}
 \mm w$ and on the other hand $v \cm u = w$ implies $w \tupeq v$.

Furthermore, the translation lemma holds.
\begin{lemma} \label{lem.confluentqc}~\\
{\sl
Let $F$ be a set of polynomials in $\myk[\g]$ and $p,q,h \in\myk[\g]$ some
 polynomials.
\begin{enumerate}
\item
Let $p-q \red{}{\myr}{qc}{F} h$.
Then there are  $p',q' \in \myk[\g]$ such that 
 $p  \red{*}{\myr}{qc}{F} p', q  \red{*}{\myr}{qc}{F} q'$ and $h=p'-q'$.
\item
Let $0$ be a normal form of $p-q$ with respect to $\red{}{\myr}{qc}{F}$.
Then there exists a polynomial  $g \in \myk[\g]$ such that
 $p  \red{*}{\myr}{qc}{F} g$ and $q  \red{*}{\myr}{qc}{F} g$.
\lemend
\end{enumerate}
}
\end{lemma}
\Ba{}
\begin{enumerate}
\item  Let $p-q \red{}{\myr}{qc}{F} h = p-q-\alpha \skm f \mrm w$,
        where $\alpha \in \myk^*, f \in F, w \in \g$
        and $\hterm(f) \mm w = t \tupeq \hterm(f)$, i.e. $\alpha \skm \hc(f)$ is the
        coefficient of $t$ in $p-q$.
       We have to distinguish three cases:
       \begin{enumerate}
         \item $t \in \terms(p)$ and $t \in \terms(q)$: 
               Then we can eliminate the term $t$ in the polynomials 
                $p$ respectively $q$ by qc-reduction. We then
                get $p \red{}{\myr}{qc}{f} p - \alpha_1 \skm f \mrm w= p'$ and
                $q \red{}{\myr}{qc}{f} q - \alpha_2 \skm f \mrm w= q'$,
                with $\alpha_1  -  \alpha_2 =\alpha$,
                where $\alpha_1 \skm \hc(f)$ and 
                $\alpha_2 \skm \hc(f)$ are
                the coefficients of $t$ in $p$ respectively $q$.
         \item $t \in \terms(p)$ and $t \not\in \terms(q)$: 
               Then  we can eliminate the term $t$ in the polynomial 
                $p$  by  qc-reduction and  get
                $p \red{}{\myr}{qc}{f} p - \alpha \skm f \mrm w= p'$ 
                and $q = q'$.
         \item $t \in \terms(q)$ and $t \not\in \terms(p)$: 
               Then  we can eliminate the term $t$ in the polynomial 
                $q$ by  qc-reduction and  get
                $q \red{}{\myr}{qc}{f} q + \alpha \skm f \mrm w= q'$
                and $p = p'$.
       \end{enumerate}
      In all cases we have $p' -q' =  p - q - \alpha \skm f \mrm w = h$.
\item We show our claim by induction on $k$, where $p-q \red{k}{\myr}{qc}{F} 0$.
      In the base case $k=0$ there is nothing to show.
      Hence, let $p-q \red{}{\myr}{qc}{F} h  \red{k}{\myr}{qc}{F} 0$.
      Then by (1) there are polynomials $p',q' \in \myk[\g]$ such that 
       $p  \red{*}{\myr}{qc}{F} p', q  \red{*}{\myr}{qc}{F} q'$ and $h=p'-q'$.
      Now the induction hypothesis for $p'-q' \red{k}{\myr}{qc}{F} 0$  yields 
       the existence of a polynomial $g \in \myk[\g]$ such that
       $p  \red{*}{\myr}{qc}{F} p' \red{*}{\myr}{qc}{F} g$ and
       $q  \red{*}{\myr}{qc}{F} q' \red{*}{\myr}{qc}{F} g$.
\\
\qed
\end{enumerate}\renewcommand{\baselinestretch}{1}\small\normalsize
Let us continue by defining Gr\"obner bases with respect to
 quasi-commutative reduction.
\begin{definition}~\\
{\rm
A set $G \subseteq \myk[\g]$ is said to be a 
 \index{quasi-commutative!Gr\"obner basis}\betonen{(right) Gr\"obner basis} with respect to $\red{}{\myr}{qc}{}$ or a \betonen{quasi-commutative Gr\"obner basis}, if
\begin{enumerate}
\item[(i)]
 $\red{*}{\lr}{qc}{G} = \;\;\equiv_{\ideal{r}{}(G)}$, and
\item[(ii)] $\red{}{\myr}{qc}{G}$ is confluent.
\dend
\end{enumerate}
}
\end{definition}
As before, when weakening our reduction we have to do
 saturation to express the right ideal congruence.
Reviewing our previous approaches to saturation,
 we can restrict saturation for a polynomial $p$ in a torsion-free nilpotent
 group ring to computing appropriate representatives for the sets
 $Y_{t}$, $t \in \terms(t)$ as specified in definition \ref{def.sat}
 (compare procedure {\sc Saturation 2} on page \ref{saturation}).
We show that it is decidable whether $Y_t$ is empty and how in case $Y_t \neq \emptyset$ a representative
for $Y_t$ can be constructed, in particular that we can construct a
polynomial $q \in Y_t$ such
that for all $q' \in Y_t$ we have $q' \red{}{\myr}{qc}{q} 0$.
First let us look at an example to illustrate how a term can be
brought into head position and what candidates are likely to cause
such a situation.
\begin{example}~\\
{\rm 
Let $\g$ be the free nilpotent group of class 2 with 2 generators
with the CNI-presentation $(\Sigma, T)$ given in example
\ref{exa.fnp2CNI}.
Further let us take a look at the polynomial $p = a_1^{-1}a_2^2a_3^3 +
a_1^{-1}a_2a_3^{-2} + a_1^{-1}a_2a_3$. where the term $t \id
a_1^{-1}a_2a_3$ is the one to be brought in head position.
For $w \id a_1a_2^{-2}a_3$ we get
 $p \mrm w = a_3^5 + a_2^{-1} + \underline{a_2^{-1}a_3^3}$, i.e., 
 $Y_{t} \neq \emptyset$ as $\hterm(p \mrm w) = t \mm w$.
Now all terms in $p$ start with the prefix $a_1^{-1}$ and the
distinguishing letter between the head term $a_1^{-1}a_2^2a_3^3$ of
$p$ and the term $t \id a_1^{-1}a_2a_3$ is $a_2$.
Notice that in order to bring the latter term into head position one
has to change the occurrence of $a_2$.
One idea might be to split $w$ into a prefix in the letters smaller
than the distinguishing letter $a_2$ and a remaining part in
$\ord(\Sigma_2)$, i.e., $w \id w'w''$ with $w' \id a_1$, $w'' \id
a_2^{-2}a_3$.
Let us now look at $p \mrm w''$.
We get $p \mrm w'' = a_1^{-1}a_3^4 +
\underline{a_1^{-1}a_2^{-1}a_3^{-1}} + a_1^{-1}a_2^{-1}a_3^2$ and
unfortunately $\hterm(p \mrm w'') \neq t \mm w''$.
This is due to the fact that although multiplication with $a_1$ is the
same on the prefix  $a_1^{-1}$ for all three terms, it does
have influence on the remaining part of the terms as multiplication is
not commutative\footnote{For C-presentations $w''$ is the appropriate
  candidate to bring $t$ into head position.}.
On the other hand there is an element in $\ord(\Sigma_2)$ that can bring
$t$ into head position, namely $v \id a_2^{-2}a_3^2$.
We get $p \mrm v = a_1^{-1}a_3^5 + a_1^{-1}a_2^{-1} +
\underline{a_1^{-1}a_2^{-1}a_3^3}$.
\exaend
}
\end{example}
The next lemma shows that if a term  can be brought into head
position in a polynomial 
by right multiplication with a group element, then there is a group
element of a special form that will also do. 
Since this special element will only depend on the polynomial and is
given in a constructive way, we then
can decide whether there exists such an element and if this is not the
case the term cannot be brought into head position.
\begin{lemma}\label{lem.dletter}~\\
{\sl
Let $(\Sigma, T)$ be a convergent CNI-presentation of a group $\g$ and
$p$ a non-zero polynomial in $\myk[\g]$.
In case there exists an element $w \in \g$ such that $\hterm(p
\mrm w) = t \mm w$ for some $t \in \terms(p)$, let $a_d$ is the distinguishing
letter between $t$ and $\hterm(p)$.
Then one can construct an
element $v \in \ord(\Sigma_d)$ such that  $\hterm(p\mrm v) = t \mm v$.
\lemend
}
\end{lemma}
\Ba{}~\\
We show that for all polynomials $q \in \{ p \mrm u | u \in \g \}$ the
following holds:
In case $\hterm(q \mrm w) = t_i \mm w$ for some $w \in \g$, $t_i \in
\terms(q)$ then one can construct an element $v \in \ord(\Sigma_d)$ where $a_d$
is the distinguishing letter between $t_i$ and $\hterm(q)$, and
$\hterm(q \mrm v) = t_i \mm v$.
This will be done by induction on $k$ where $d = n-k$.
In the base case let $k=0$, i.e., $a_{n}$ is the
distinguishing letter between $\hterm(q) = t_1 \id a_1^{1_1} \ldots
a_{n\phantom{1}}^{1_n}$ and $t_i \id a_1^{i_1} \ldots
a_{n\phantom{1}}^{i_n}$.
Hence  $1_j = i_j$ for all $1 \leq j \leq n-1$ and $1_n >_{\z} i_n$.
By our assumption there exists $w \in \g$ such that $\hterm(q \mrm w)
= t_i \mm w$, with $w \id w'a_{n\phantom{1}}^{w_n}$, $w' \in
\ord(\Sigma \backslash \Sigma_n)$, and there exist $k_1, \ldots , k_{n-1}, x \in
\z$ such that 
\auskommentieren{ 
\begin{tabbing}
$t_1 \mm w$ \= = \= $a_1^{1_1} \ldots a_{n\phantom{1}}^{1_n} \mm w$ \\
            \> = \> $a_1^{1_1} \ldots a_{n-1\phantom{1}}^{1_{n-1}} \mm
            w \mm a_{n\phantom{1}}^{1_n}$ \\
            \> = \> $(a_1^{1_1} \ldots a_{n-1\phantom{1}}^{1_{n-1}}
            \mm w') \mm a_{n\phantom{1}}^{1_n + w_n}$\\
            \> $\id$ \> $a_1^{k_1} \ldots a_{n-1\phantom{1}}^{k_{n-1}}
            a_{n\phantom{1}}^{1_n  +  x}$ \\ 
and\\
$t_i \mm w$ \> = \> $a_1^{1_1} \ldots a_{n-1\phantom{1}}^{1_{n-1}}a_{n\phantom{1}}^{i_n} \mm w$ \\
            \> = \> $a_1^{1_1} \ldots a_{n-1\phantom{1}}^{1_{n-1}} \mm
            w \mm a_{n\phantom{1}}^{i_n}$ \\
            \> = \> $(a_1^{1_1} \ldots a_{n-1\phantom{1}}^{1_{n-1}}
            \mm w') \mm a_{n\phantom{1}}^{i_n + w_n}$ \\
            \> $\id$ \> $a_1^{k_1} \ldots a_{n-1\phantom{1}}^{k_{n-1}}
            a_{n\phantom{1}}^{i_n  +  x}$.
\end{tabbing}}
$t_1 \mm w = a_1^{1_1} \ldots a_{n\phantom{1}}^{1_n} \mm w = a_1^{1_1} \ldots a_{n-1\phantom{1}}^{1_{n-1}} \mm
            w \mm a_{n\phantom{1}}^{1_n} = (a_1^{1_1} \ldots a_{n-1\phantom{1}}^{1_{n-1}}
            \mm w') \mm a_{n\phantom{1}}^{1_n + w_n} \id a_1^{k_1} \ldots a_{n-1\phantom{1}}^{k_{n-1}}
            a_{n\phantom{1}}^{1_n  +  x}$
and 
$t_i \mm w =a_1^{1_1} \ldots
a_{n-1\phantom{1}}^{1_{n-1}}a_{n\phantom{1}}^{i_n} \mm w = a_1^{1_1} \ldots a_{n-1\phantom{1}}^{1_{n-1}} \mm
            w \mm a_{n\phantom{1}}^{i_n} =(a_1^{1_1} \ldots a_{n-1\phantom{1}}^{1_{n-1}}
            \mm w') \mm a_{n\phantom{1}}^{i_n + w_n} \id a_1^{k_1} \ldots a_{n-1\phantom{1}}^{k_{n-1}}
            a_{n\phantom{1}}^{i_n  +  x}$.
Thus  
 $1_n  + x <_{\z} i_n  + x$ must hold.
Let us  set  $v \id a_{n}^{-1_n}$.
We  show that for all $t_j \in T(q) \backslash \{t_i \}$ we
have $t_i \mm v \succ t_j \mm v$.
Note that for all $t_j$ with prefix $a_1^{j_1} \ldots
a_{n-1\phantom{1}}^{j_{n-1}} \pred a_1^{1_1} \ldots
a_{n-1\phantom{1}}^{1_{n-1}}$ we have $t_j \mm v \pred t_i \mm v$, as
right multiplication with $v$ only changes the exponent of $a_n$ in
the respective term.
It remains to look at those terms $t_j$ with $a_1^{j_1} \ldots
a_{n-1\phantom{1}}^{j_{n-1}} \id a_1^{1_1} \ldots
a_{n-1\phantom{1}}^{1_{n-1}}$.
Hence, let us assume that there exists  a term $t_j$ such that $t_j \mm v \succ t_i
\mm v$, i.e., $j_n - 1_n >_{\z} i_n - 1_n$.
Since  $\hterm(q \mrm w) = t_i \mm w$ we know 
$j_n + x <_{\z} i_n + x$ and $1_n +x <_{\z} i_n + x$.
Furthermore, as $t_1 = \hterm(q)$ we have $1_n >_{\z} i_n$ and 
$1_n >_{\z} j_n$.
We prove that $t_j \mm v \succ t_i \mm v$ yields $j_n +x >_{\z}
i_n +x$  contradicting our assumption by analysing the possible cases for
these exponents.
First suppose that $1_n < 0$ and thus $1_n +x <_{\z} i_n + x$ implies $x \geq |1_n| > 0$.
Then in case $i_n \leq 0$ this gives us $|1_n| > |i_n|$.
Now $j_n - 1_n >_{\z} i_n - 1_n > 0$ and $j_n - 1_n >0$ yields either
 $j_n > 0$ or  ($j_n \leq 0$ and  $|j_n| < |i_n|$), 
 both implying $j_n + y > i_n + y$ for all $y \geq |1_n|$,
 especially for $y=x$.
In case $i_n > 0$ as before $j_n - 1_n >_{\z} i_n - 1_n > 0$, $j_n -
1_n >0$ imply $j_n > i_n$ and
for all $y \geq |1_n|$ we get $j_n + y > i_n + y$, especially for $y = x$.
Hence let us assume that $1_n > 0$ and thus $1_n +x <_{\z} i_n + x$
implies $x < 0$ and $|x| > i_n$, since $1_n > i_n \geq 0$ and $1_n > j_n \geq 0$.
Now $j_n - 1_n >_{\z} i_n - 1_n$ and $i_n - 1_n < 0$  imply $j_n - 1_n
<0$ and $|i_n - 1_n| < |j_n - 1_n|$.
Hence we get $j_n
< i_n$ and for all  $y < 0$ with $|y| > j_n$ we  have
  $j_n + y >_{\z} i_n + y$, especially for $y = x$ as $|x| > i_n > j_n$.

In the induction step let us assume that for all polynomials $q \in \{
p \mrm u | u \in \g \}$ and $w \in \g$
with $\hterm(q \mrm w) = t_i \mm w$, if the distinguishing letter $a_d$
between $\hterm(q)$ and $t_i$ has index $d \geq n - (k-1)$ there exists
an element $v \in \ord(\Sigma_d)$ such that $\hterm(q \mrm v) = t_i \mm
v$.
Now for $q \in \{ p \mrm u | u \in \g \}$, $w \in \g$ with $\hterm(q
\mrm w) = t_i \mm w$ let us assume that the distinguishing letter
between $\hterm(q)$ and $t_i$ has index $d = n-k$.
Since $\hterm(q \mrm w) = t_i \mm w$, for $w \id
w'a_{d\phantom{1}}^{w_d}w''$ with $w' \in \ord(\Sigma \backslash \Sigma_d)$, $w''
\in \ord(\Sigma_{d+1})$, we know that there exist $k_1, \ldots, k_{d-1}, x
\in \z$ and $z_1, z_i, \tilde{z}_1 \in \ord(\Sigma_{d+1})$
such that 
\auskommentieren{
\begin{tabbing}
$t_1 \mm w$ \= = \= $a_1^{1_1} \ldots a_{n\phantom{1}}^{1_n} \mm w$ \\
            \> = \> $a_1^{1_1} \ldots a_{d-1\phantom{1}}^{1_{d-1}} \mm
            w' \mm a_{d\phantom{1}}^{1_d} \mm \tilde{z}_1$ \\
            \> $\id$ \> $a_1^{k_1} \ldots 
            a_{d-1\phantom{1}}^{k_{d-1}}a_{n\phantom{1}}^{1_d + x}z_1$\\
and similarly\\
$t_i \mm w$ \> = \>  $a_1^{k_1} \ldots a_{d-1\phantom{1}}^{k_{d-1}} 
                      a_{n\phantom{1}}^{i_d + x}z_i$.
\end{tabbing}}
$t_1 \mm w = a_1^{1_1} \ldots a_{n\phantom{1}}^{1_n} \mm w = a_1^{1_1} \ldots a_{d-1\phantom{1}}^{1_{d-1}} \mm
            w' \mm a_{d\phantom{1}}^{1_d} \mm \tilde{z}_1 \id a_1^{k_1} \ldots 
            a_{d-1\phantom{1}}^{k_{d-1}}a_{n\phantom{1}}^{1_d + x}z_1$
and similarly
$t_i \mm w = a_1^{k_1} \ldots a_{d-1\phantom{1}}^{k_{d-1}} 
                      a_{n\phantom{1}}^{i_d + x}z_i$.
As $1_d \neq i_d$ then $1_d + x <_{\z} i_d + x$ must hold  and we can
set $v \id a_{n\phantom{1}}^{-1_d}$. 
We have to show that for all $t_j \in T(q) \backslash \{t_i \}$ we
have $t_i \mm v \succ t_j \mm v$.
Note that for all $t_j$ with prefix $a_1^{j_1} \ldots
a_{d-1\phantom{1}}^{j_{d-1}} \pred a_1^{1_1} \ldots
a_{d-1\phantom{1}}^{1_{d-1}}$ we have $t_j \mm v \pred t_i \mm v$, as
right multiplication with $v$ has no influence on the prefix 
in $\ord(\Sigma \backslash \Sigma_d)$.
Therefore, it remains to look at those terms $t_j$ with $a_1^{j_1} \ldots
a_{d-1\phantom{1}}^{j_{d-1}} \id a_1^{1_1} \ldots
a_{d-1\phantom{1}}^{1_{d-1}}$.
Let us assume that there exists a term $t_j$ such that $t_j \mm v \succ t_i
\mm v$, i.e., $j_d - 1_d \geq_{\z} i_d - 1_d$.
We will show that then $j_d = i_d$ and hence our induction hypothesis
can be applied since for the polynomial $q \mrm v$ the distinguishing
letter between $\hterm(q \mrm v)$ and $t_i \mm v$ is of index $d' > d
=n-k$ and by our assumption there exists $\inv{v} \mm w \in \g$ such that
$\hterm((q \mrm v ) \mrm (\inv{v} \mm w)) = \hterm (q \mrm w) = t_i
\mm w = t_i \mm (\inv{v} \mm w)$. 
We know $j_d +x \leq_{\z} i_d + x$ and $1_d +x <_{\z} i_d + x$ since
$\hterm(q \mrm w) = t_i \mm w$.
Next we prove that $t_j \mm v \succ t_i \mm v$ implies $j_d =
i_d$ by analysing the possible cases.
First suppose that $1_d < 0$ and thus $1_d +x <_{\z} i_d + x$ implies $x \geq |1_d| > 0$.
Then in case $i_d \leq 0$ this gives us $|1_d| > |i_d|$.
Now $j_d - 1_d \geq_{\z} i_d - 1_d > 0$ and $j_d - 1_d >0$ yield either
 $j_d > 0$ or  ($j_d \leq 0$ and  $|j_d| \leq |i_d|$), 
 both implying $j_d + y \geq i_d + y$ for all $y \geq |1_d|$.
Thus as $x \geq |1_d|$ we get  $j_d + x \geq i_d + x$ yielding $j_d = i_d$.
In case $i_d > 0$ as before $j_d - 1_d \geq_{\z} i_d - 1_d > 0$, $j_d
- 1_d >0$ yield $j_d \geq i_d$ and
for all $y \geq |1_n|$ we get $j_d + y \geq i_d + y$.
Thus as $x \geq |1_d|$, $j_d + x \geq i_d + x$ again yields $j_d = i_d$.
Therefore, let us assume that $1_d > 0$ and thus $1_d +x <_{\z} i_d + x$
implies $x < 0$ and $|x| > i_d$, since $1_d > i_d \geq 0$ and $1_d \geq j_d \geq 0$.
Now $j_d - 1_d \geq_{\z} i_d - 1_d$ and $i_d - 1_d < 0$  imply $j_d - 1_d
<0$ and $|i_d - 1_d| \leq |j_d - 1_d|$.
Hence we get $j_d \leq i_d$ and for all $y < 0$ with $|y| > j_d$, we have
$j_d + y \geq_{\z} i_d + y$.
Thus as $|x| > i_d \geq j_d$, then $j_d + x \geq_{\z} i_d + x$ yields
$j_d = i_d$.
\\
\qed
\begin{corollary}~\\
{\sl
Let $(\Sigma, T)$ be a convergent CNI-presentation of a group $\g$ and
$p$ a non-zero polynomial in $\myk[\g]$.
Then it is decidable whether $Y_t \neq \emptyset$ for some $t \in \terms(p)$.
}
\end{corollary}
\Ba{}~\\
The proof of lemma \ref{lem.dletter} can be turned into a procedure
which, given a polynomial $p \in \myk[\g]$ and a term $t \in \terms(p)$,
tries to compute a representative of $Y_t$.
This can be done by starting with the distinguishing letter between
$t$ and $\hterm(p)$ and proceeding to compute an appropriate $v \in
\ord(\Sigma_{d})$ if possible. 
Remember that for each modification step the multiple is defined by
the exponent of the distinguishing letter of the respective
polynomial.
Hence let $p_0 = p$, $t_0 = t$ and $d_0$ be the index of the
distinguishing letter between $\hterm(p_0)$ and $t_0$ and further let
$i_0$ be the exponent of $a_{d_0}$ in the term $\hterm(p_0)$.
Now proceed for $p_j$, $t_j$, $d_j$ and $v_j \id a_{d_j}^{-i_j}$,
$j \geq 0$ as follows:
 If the exponent of $a_{d_j}$ in $\hterm(p_j \mrm v_j)$
  differs from its exponent in $t_j \mm v_j$ then $t_j$ cannot be
  brought to head position.
 If $\hterm(p_j \mrm v_j) = t_j \mm v_j$ then we are done.
Else we set $p_{j+1} = p_j \mrm v_j$, $t_{j+1} = t_j \mm v_j$ and
$i_{j+1}>i_j$ is the index of the distinguishing letter between 
$\hterm(p_{j+1})$ and $t_{j+1}$.
\\
Now either we find that $t$ cannot be brought to head position or for
some $0 \leq k < n$ we have $v \id v_0 \ldots v_k$.
\\
\qed
Hence we can  compute a saturating set for a polynomial
 in a torsion-free nilpotent group ring with respect to right reduction.
We will now move on to show how these polynomials can be modified in
order to compute a saturating set with respect to quasi-commutative reduction.
\begin{definition}\label{def.satqc}~\\
{\rm
A set of polynomials $F\subseteq \{\alpha \skm p \mrm w | \alpha \in
\myk,  w \in\g \}$
 is called  a 
 \index{quasi-commutative!saturation}\index{saturation!quasi-commutative}\index{quasi-commutatively!saturating set}\index{saturating set!quasi-commutatively}\betonen{quasi-commutatively saturating set}\/ for  
 a polynomial $p \in \myk[\g]$,
 if for all $\alpha \in \myk^*$, $w \in \g$, $\alpha \skm p \mrm w  \red{}{\myr}{qc}{F} 0$ holds.
$\SAT_{qc}(p)$ denotes the family of all quasi-commutatively
 saturating sets for $p$.
We call a set $F \subseteq \myk[\g]$ \index{quasi-commutatively!saturated set}\betonen{quasi-commutatively
 saturated}, if for all $f \in F$, $\alpha \in \myk^*$, $w \in \g$, we have 
 $\alpha \skm f \mrm w \red{}{\myr}{qc}{F} 0$.
\dend
}
\end{definition}
Important is that quasi-commutative saturated sets allow special
representations.
\begin{lemma}\label{lem.propqc}~\\
{\sl
Let $F$ be a quasi-commutatively saturated set of polynomials in $\myk[\g]$.
 Then every non-zero polynomial $g \in \ideal{r}{}(F)$ has a representation $g =  \sum_{i=1}^{k} \alpha_i \skm f_i \mrm w_i$,
 where $\alpha_i \in \myk^*, f_i \in F, w_i \in \g$,
 $\hterm(f_i \mrm w_i) = \hterm(f_i) \mm w_i$
 and $\hterm(f_{i} \mrm w_i) \tupeq \hterm(f_{i})$.
\lemend\ohnebeweis
}
\end{lemma}
For a polynomial $p$ and a term $t \in \terms(p)$ we call a term $s$
in a multiple $p \mrm w$ a \betonen{$t$-term}, if $s = t \mm w$.
The following lemma states that if in two multiples of a polynomial
the head terms result from the same term $t$, then there is also a
multiple of the polynomial with a $t$-term as head term which is in
some sense a greatest common divisor of the head terms of the
original two polynomials. 
This leads to the existence of one representative for each $t$-term
occurring as a head term comparable to the existence of a unique minimal
element in the set $Y_t$ with respect to the ordering $\tup$ which is
not obvious, as $\tup$ is no total ordering.
\begin{lemma}\label{lem.uniqueCNI}~\\
{\sl
Let $(\Sigma, T)$ be a convergent CNI-presentation of a group $\g$.
Further let $p \mrm u$ and $p \mrm v$ be two multiples of a non-zero
polynomial $p \in \myk[\g]$ such that for some term $t \in \terms(p)$
the head terms are $t$-terms, i.e.,
$\hterm(p \mrm u) = t \mm u \id a_1^{i_1} \ldots a_n^{i_n}$ and
$\hterm(p \mrm v) = t \mm v \id a_1^{j_1} \ldots a_n^{j_n}$.
Then there exists a term $\tilde{t} 
\tupleq a_1^{\rho_1} \ldots a_n^{\rho_n}$ where
\[ \rho_l = \left\{ \begin{array}{l@{\quad\quad}l}
         \sgn(i_l) \skm \min \{ |i_l|, |j_l| \} & \sgn(i_l) =
         \sgn(j_l) \\
         0 &  \mbox{otherwise}
 \end{array} \right. \]
and an element $\tilde{z} \in \g$ such that   
 $\hterm(p \mrm \tilde{z}) = t \mm \tilde{z} = \tilde{t}$. 
}
\end{lemma}
\Ba{}~\\
Let $p$, $p \mrm u$ and $p \mrm v$ be as described in the lemma and
 let $a_1, \ldots, a_n$ be the letters corresponding to our presentation.
We show the existence of $\tilde{z}$ by constructing a sequence
 $z_1, \ldots , z_n \in \g$, such that for $1 \leq l \leq n$ we have
 $\hterm(p \mrm z_l) = t \mm z_l \id a_1^{s_1} \ldots a_l^{s_l}r_l$
 with $r_l \in \ord(\Sigma_{l+1})$ and 
 $a_1^{s_1} \ldots a_l^{s_l} \tupleq a_1^{\rho_1} \ldots a_l^{\rho_l}$.
Then for $\tilde{z} = z_n$ our claim holds.
\\
Let us start by constructing an element $z_1 \in \g$ such that 
 $\hterm(p \mrm z_1) = t \mm z_1 \id a_1^{s_1}r_1$,
 $r_1 \in \ord(\Sigma_{2})$ and $a_1^{s_1} \tupleq a_1^{\rho_1}$.
In case $i_1 = j_1$ or $j_1 =0$ we can set $z_1 = v$ and
 $s_1 = j_1 = \rho_1$ since
 $\hterm(p \mrm v) = t \mm v \id a_1^{j_1} \ldots a_n^{j_n}$.
Similarly in case $i_1 = 0$ we can set  $z_1 = u$ and
 $s_1 = i_1 = 0 = \rho_1$
 since $\hterm(p \mrm u) = t \mm u \id a_2^{i_2}\ldots
 a_n^{i_n}\in\ord(\Sigma_2)$.
Hence let us assume $i_1 \neq j_1$ and both are non-zero.
First suppose that $\sgn(i_1) = \sgn(j_1)$.
Then if $|i_1| \geq |j_1|$ we again  set $z_1 = v$ 
 since for $s_1 = j_1 = \rho_1$ our claim holds.
In case $|j_1| > |i_1|$ we set $z_1 = u$ because for $s_1 = i_1 = \rho_1$ our claim holds.
Now let us proceed with the case $\sgn(i_1) \neq \sgn(j_1)$, i.e., we 
 construct $z_1 \in \g$ such that 
 $\hterm(p \mrm z_1) = t \mm z_1 \in \ord(\Sigma_2)$ as $\rho_1 = 0$.
\auskommentieren{Let us consider the polynomial $p \mrm v \mrm a_1^{-j_1}$.
In case $\hterm(p \mrm v \mrm a_1^{-j_1}) = t \mm v \mm a_1^{-j_1}$
 we can set $z_1 = v \mm a_1^{-j_1}$ and are done since
 $t \mm v \mm a_1^{-j_1} \in \ord(\Sigma_2)$.
\\
Hence suppose $\hterm(p \mrm v \mrm a_1^{-j_1}) = 
 s \mm a_1^{-j_1} \id a_1^{c_1} \ldots a_n^{c_n}$
 for some $s \id a_1^{b_s}x_s\in \terms(p \mrm v)$, $x_s \in \ord(\Sigma_2)$.
\\
We show that then we have $c_1 =0$.
Then for the polynomial $p \mrm v \mrm a_1^{-j_1}$ the distinguishing letter
 between $\hterm(p \mrm v \mrm a_1^{-j_1})$ and the term
 $t \mm v \mm a_1^{-j_1}$ has at least index 2 and since 
 $\hterm((p \mrm v \mrm a_1^{-j_1}) \mrm a_1^{-j_1}) = \hterm(p \mm v)
  = t \mm v$ by lemma \ref{lem.dletter} there exists an element 
 $r \in \ord(\Sigma_2)$ such that
 $\hterm((p \mrm v \mrm a_1^{-j_1}) \mrm r) = t \mm v \mm a_1^{-j_1} \mm r
 \in \ord(\Sigma_2)$ and we can set $z_1 = v \mm a_1^{-j_1} \mm r$.
\\
Hence it remains to show that $c_1 = 0$.
\\
Remember that $s \mm a_1^{-j_1} = a_1^{b_s}x_s \mm a_1^{-j_1} 
 \id a_1^{c_1} \ldots a_n^{c_n}$ yields $c_1 = b_s - j_1$ and furthermore
 as $s \in \terms(p \mrm v)$ we know $b_s \leq_{\z} j_1$.
\\
By multiplying $t \mm v \mm a_1^{-j_1}$ with
 $\inv{v \mm a_1^{-j_1}} \mm u = a_1^{j_1} \mm \inv{v} \mm u$ we can
 bring this term to head position as
 $(p \mrm v \mrm a_1^{-j_1}) \mrm (a_1^{j_1} \mm \inv{v} \mm u) = p \mrm u$.
Hence multiplying $\hterm(p \mrm v \mrm a_1^{-j_1})$ by
 $a_1^{j_1} \mm \inv{v} \mm u$ must result in a term smaller than
 $t \mm u$.
Now $\hterm(p \mrm v \mrm a_1^{-j_1}) \mm a_1^{j_1} \mm \inv{v} \mm u
 = a_1^{c_1} \ldots a_n^{c_n} \mm a_1^{j_1} \mm a_1^{f_1}u'
 = a_1^{c_1 + j_1+f_1}y$ where $\inv{v} \mm u \id a_1^{f_1}u'$ and
  $y, u' \in \ord(\Sigma_2)$, $f_1 \in \z$.
Remember that for $t \mm v \mm a_1^{-j_1} \id a_2^{d_2} \ldots a_n^{d_n}$
 we have
 $(t \mm v \mm a_1^{-j_1}) \mm a_1^{j_1} \mm \inv{v} \mm u
 = a_2^{d_2} \ldots a_n^{d_n} \mm a_1^{j_1} \mm a_1^{f_1}u'
 \id a_1^{i_1} \ldots a_n^{i_n}$ implying $j_1 + f_1 = i_1$.
\\
Thus we know $c_1 + j_1 + f_1 = b_s - j_1 + i_1 \leq_{\z} i_1$.
Let us now prove that $c_1 = 0$ by distinguishing the two possible cases.
\\
If $i_1 > 0$ and $j_1 < 0$ then $b_s - j_1 + i_1 \leq_{\z} i_1$
 implies $0 \leq b_s - j_1 + i_1 \leq i_1$ and as $-j_1 + i_1 > i_1$ we
 must have $b_s < 0$.
On the other hand, $s \id a_1^{b_s}z_s \in \terms(p \mrm v)$ implies
 $b_s \leq_{\z} j_1$ and thus $|b_s| \leq |j_1|$.
Then $b_s - j_1 \geq 0$ and hence $0 \leq b_s - j_1 + i_1 \leq i_1$
 yields $c_1 = b_s - j_1 = 0$.
\\
It remains to examine the case $i_1 < 0$ and $j_1 >0$.
Since $s \id a_1^{b_s}z_s \in \terms(p \mrm v)$ now we know $0 \leq
 b_s \leq j_1$ and this implies $b_s - j_1 \leq 0$.
Hence, $b_s - j_1 + i_1 < 0$ and thus $b_s - j_1 + i_1 \leq_{\z} i_1$
 yields $|b_s - j_1 + i_1| \leq |i_1|$ implying $c_1 =b_s -j_1 =
 0$.}
We claim that  the letter $a_1$ has the same exponent for all
 terms in $\terms(p)$, say $b$.
In case this holds, no term in the polynomial $p \mrm a_1^{-b}$ will
 contain the letter $a_1$ and the distinguishing letter between
 $\hterm(p \mrm a_1^{-b})$ and the term $t \mm a_1^{-b}$ is at least of
 index 2.
Furthermore we know $\hterm((p \mrm a_1^{-b}) \mrm (a_1^{b} \mm v)) =
 \hterm(p \mrm v) = t \mm v$.
Thus by lemma \ref{lem.dletter} there exists an element $r \in
 \ord(\Sigma_2)$ such that $\hterm((p \mrm a_1^{-b}) \mrm r) = t \mm
 a_1^{-b} \mm r \in \ord(\Sigma_2)$ and thus we can set $z_1 = a_1^{-b}r$ and $s_1 = 0 = \rho_1$.
\\
Hence it remains to prove our initial claim.
Suppose we have the representations $s' \id a_1^{b_{s'}}x_s$, $b_{s'}
\in \z$, $x_{s'} \in \ord(\Sigma_2)$ for the
terms $s' \in \terms(p)$ and $ \hterm(p)=s \id a_1^{b_s}x_s$.
Then we know $b_s \geq_{\z} b_t$ since $t \in \terms(p)$.
Hence in showing that the case $b_s >_{\z} b_t$ is not possible we
 find that the exponents of $a_1$ in $s$ and $t$ are equal.
To see this, let us study the possible cases.
If $b_s > 0$ we have $b_s > b_t \geq 0$ and hence there exists no $x
\in \z$ such that $b_t + x > b_s + x \geq 0$.
On the other hand $b_s < 0$ either implies $b_t > 0$ or
 $b_t \leq 0$ and $|b_s| > |b_t|$.
In both cases there exists no $x \in \z$ such that $b_t + x < 0$ and
 $|b_t + x| > |b_s + x|$.
Hence $b_t = b_s$ must hold as we know that $t$ can be brought to 
 head position by $u$ respectively $v$ 
 such that the exponents of $a_1$  in
 $\hterm(p \mrm u)$ respectively $\hterm(p \mrm v)$ have different
 sign.
It remains to show that there cannot exist a term $s' \in \terms(p)$
 with $b_{s'} <_{\z} b_s = b_t$.
Let us assume such an $s'$ exists.
Since $\hterm(p \mrm u) = t \mm u \id a_1^{i_1} \ldots a_n^{i_n}$ and
 $\hterm(p \mrm v) = t \mm v \id a_1^{j_1} \ldots a_n^{j_n}$ there
 then must exist $x_1, x_2 \in \z$ such that
 $b_{s'} + x_1 <_{\z} b_t + x_1 = i_1$ and $b_{s'} + x_2 <_{\z} b_t +
 x_2 = j_1$.
Without loss of generality let us assume $i_1 > 0$ and $j_1 < 0$ (the
 other case is symmetric).
In case $b_t < 0$ we get that $b_t + x_1 = i_1 > 0$ implies
 $x_1 > |b_t| > 0$.
Now, as $b_{s'} <_{\z} b_t$ either implies
 $b_{s'} > 0$ or $b_{s'} \leq 0$ and $|b_{s'}| < |b_t|$, we find
 $b_{s'} + x_1 > b_t + x_1$ contradicting $b_{s'} + x_1 <_{\z} b_t +
 x_1$.
On the other hand, in case $b_t > 0$ we know $b_t > b_{s'} \geq 0$.
Furthermore, $b_t + x_2 = j_1 < 0$ implies $x_2 <0$ and $|x_2| >
 b_t$.
Hence we get $b_{s'} + x_2 < 0$ and $|b_{s'} + x_2| > |b_t + x_2|$
 contradicting $b_{s'} + x_2 <_{\z} b_t + x_2$.
\\
Thus let us assume that for the letter $a_{k-1}$ we have
constructed  $z_{k-1} \in \g$ such that
 $\hterm(p \mrm z_{k-1}) = t \mm z_{k-1} \id a_1^{s_1} \ldots
 a_{k-1}^{s_{k-1}}r_{k-1} \id a_1^{s_1} \ldots
 a_{k-1}^{s_{k-1}}a_k^{l_k}r'$
 with $r_{k-1} \in \ord(\Sigma_{k})$, $r' \ord(\Sigma_{k+1})$ and 
 $a_1^{s_1} \ldots a_{k-1}^{s_{k-1}} \tupleq a_1^{\rho_1} \ldots a_{k-1}^{\rho_{k-1}}$.
We now show that we can find $z_k = z_{k-1} \mm \tilde{w} \in \g$ such that
 $\hterm(p \mrm z_{k}) = t \mm z_{k} \id a_1^{s_1} \ldots a_{k}^{s_{k}}r_{k}$
 with $r_{k} \in \ord(\Sigma_{k+1})$ and 
 $a_1^{s_1} \ldots a_{k}^{s_{k}} \tupleq a_1^{\rho_1} \ldots a_{k}^{\rho_{k}}$.
This will be done in two steps.
First we show that for the polynomials $p \mrm u$ and $p \mrm z_{k-1}$
 with head terms $a_1^{i_1}\ldots a_n^{i_n}$ respectively
 $a_1^{s_1} \ldots a_{k-1}^{s_{k-1}}a_k^{l_k}r'$
 we can find an element $w_1 \in \g$ such that
 $\hterm(p \mrm z_{k-1} \mrm w_1) = t \mm z_{k-1} \mm w_1
    \id a_1^{s_1} \ldots a_{k-1}^{s_{k-1}}a_{k}^{\tilde{s}_{k}}\tilde{r}$,
  $\tilde{r} \in \ord(\Sigma_{k+1})$ and 
 $a_k^{\tilde{s}_{k}} \tupleq a_k^{\tilde{\rho}_{k}}$ with
 \[\tilde{\rho}_k = \left\{ \begin{array}{l@{\quad\quad}l}
         \sgn(i_k) \skm \min \{ |i_k|, |l_k| \} & \sgn(i_k) =
         \sgn(l_k) \\
         0 &  \mbox{otherwise}.
 \end{array} \right. \]
Then in case $a_k^{\tilde{\rho}_{k}} \tupleq a_k^{\rho_{k}}$ we are done
 and set $z_k = z_{k-1} \mm w_1$ and $s_k = \tilde{s}_{k}$.
Else we can similarly proceed for the polynomials $p \mrm v$ and 
 $p \mrm z_{k-1} \mrm w_1$ with head terms $a_1^{j_1} \ldots a_n^{j_n}$
 respectively 
 $a_1^{s_1} \ldots a_{k-1}^{s_{k-1}}a_{k}^{\tilde{s}_{k}}\tilde{r}$
 and find an element $w_2 \in \g$ such that
 for $z_k = z_{k-1} \mm w_1 \mm w_2$ we have
 $\hterm(p \mrm z_k) = t \mm z_k \id a_1^{s_1} \ldots a_{k}^{s_{k}}r_{k}$,
 $r_k \in \ord(\Sigma_{k+1})$ and
 $a_k^{s_{k}} \tupleq a_k^{\tilde{\rho}'_{k}}$ with
 \[\tilde{\rho}'_k = \left\{ \begin{array}{l@{\quad\quad}l}
         \sgn(j_k) \skm \min \{ |j_k|, |\tilde{s}_k| \} & \sgn(j_k) =
         \sgn(\tilde{s}_k) \\
         0 &  \mbox{otherwise}.
 \end{array} \right. \]
Then we can conclude $a_k^{s_{k}} \tupleq a_k^{\rho_{k}}$
 as in case $s_k = 0$ we are immediately done and otherwise we get
 $\sgn(j_k) = \sgn(\tilde{s}_k) = \sgn(\tilde{\rho}_k) = \sgn(i_k)$ and
 $\min \{ |i_k|, |\tilde{s}_k|, |j_k| \} \leq
  \min \{ |i_k|, |j_k| \}$.
\\
Let us hence show how to construct $w_1$.
Remember that $\hterm(p \mrm u) = t \mm u \id a_1^{i_1} \ldots
a_n^{i_n}$ and $\hterm(p \mrm z_{k-1}) = t \mm z_{k-1} \id a_1^{s_1}
\ldots a_{k-1}^{s_{k-1}}a_k^{l_k}r'$ for some $r' \in \ord(\Sigma_{k+1})$.
\\
In case $i_k = l_k$ or $l_k = 0$ we can set $w_1 = \lambda$ and 
 $\tilde{s}_k = l_k = \tilde{\rho}_k$ as 
 $\hterm(p \mrm z_{k-1}) = t \mm z_{k-1} \id 
   a_1^{s_1} \ldots a_{k-1}^{s_{k-1}}a_{k}^{l_{k}}r'$.
Hence let $i_k \neq l_k$ and $l_k \neq 0$.
\\
First let us assume that $\sgn(i_k) = \sgn(l_k)$.
Then in case $|i_k| \geq |l_k|$ we are done by setting $w_1 = \lambda$
 as again 
 $\hterm(p \mrm z_{k-1}) = t \mm z_{k-1} \id 
   a_1^{s_1} \ldots a_{k-1}^{s_{k-1}}a_{k}^{l_{k}}r'$
 will do with $\tilde{s}_k = l_k = \tilde{\rho}_k$.
Therefore,  let us assume that  $|l_k| > |i_k|$.
Then we consider the multiple $p \mrm z_{k-1} \mrm a_k^{-l_k+i_k}$,
 i.e., the exponent of the letter $a_k$ in the term
 $t \mm z_{k-1} \mm a_k^{-l_k+i_k}$ will be $i_k$.
If $\hterm(p \mrm z_{k-1} \mrm a_k^{-l_k+i_k}) = t \mm z_{k-1} \mm
 a_k^{-l_k+i_k}$ we are done because then
 $t \mm z_{k-1} \mm a_k^{-l_k+i_k} \id a_1^{s_1} \ldots a_{k-1}^{s_{k-1}}a_{k}^{i_{k}}\tilde{r}_k$ for some $\tilde{r}_k \in \ord(\Sigma_{k+1})$
 and we can set $w_1 =  a_k^{-l_k+i_k}$ and 
 $\tilde{s}_k = i_k = \tilde{\rho}_k$.
Otherwise we show that the $t$-term $t \mm z_{k-1} \mm a_k^{-l_k+i_k}$
 in this multiple can be brought
 to head position using an element  $r \in \ord(\Sigma_{k+1})$ thus allowing
 to set $\tilde{s}_k = i_k = \tilde{\rho}_k$ and $w_1 = a_k^{-l_k+i_k}r$ as then we have
 $\hterm(p \mrm z_{k-1} \mrm w_1) = t \mm z_{k-1} \mm w_1 = a_1^{s_1}
  \ldots a_{k-1}^{s_{k-1}}a_k^{l_k}r' \mm a_k^{-l_k+i_k}r
  \id a_1^{s_1} \ldots a_{k-1}^{s_{k-1}}a_k^{i_k}\tilde{r}$
 where $a_k^{l_k}r' \mm a_k^{-l_k+i_k}r \id
 a_k^{i_k}\tilde{r}$\footnote{Note that the product of two elements in
   $\ord(\Sigma_i)$ is again an element in $\ord(\Sigma_i)$.}.
This follows immediately if we can prove that the exponent of $a_k$ in
the term $\hterm(p \mrm z_{k-1} \mrm a_k^{-l_k + i_k}$ is also $i_k$.
Then we can apply lemma \ref{lem.dletter} to the polynomial $p \mrm
z_{k-1} \mrm a_k^{-l_k+i_k}$ and the term $t \mm z_{k-1} \mm
a_k^{-l_k+i_k}$.
Note that $\hterm(p \mrm z_{k-1} \mrm a_k^{-l_k+i_k})$ and $t \mm z_{k-1} \mm
a_k^{-l_k+i_k}$ have then distinguishing letter  of at least index $k+1$  and
further $\hterm((p \mrm z_{k-1} \mrm a_k^{-l_k+i_k}) \mrm
a_k^{-l_k+i_k}) =  \hterm(p
 \mrm z_{k-1})= t \mm z_{k-1}$.
Therefore, we show that the exponent of $a_k$ in
the term $\hterm(p \mrm z_{k-1} \mrm a_k^{-l_k + i_k})$ is also $i_k$.
Let $a_1^{s_1} \ldots a_{k-1}^{s_{k-1}} a_k^{b_k}r''$ with $r'' \in
 \ord(\Sigma_{k+1})$ be the term in $p \mrm z_{k-1}$ that became head
 term\footnote{Note that a candidate in $\terms(p \mrm z_{k-1})$ 
 for the head term in $p \mrm z_{k-1}
   \mrm a_k^{-l_k + i_k}$ must have prefix $a_1^{s_1} \ldots
   a_{k-1}^{s_{k-1}}$ since $\hterm(p \mrm z_{k-1}) \id a_1^{s_1}
   \ldots a_{k-1}^{s_{k-1}}r_{k-1}$ and multiplication with $a_k^{-l_k +
     i_k}$ only involves $r_{k-1}$.},
 i.e., $a_1^{s_1} \ldots a_{k-1}^{s_{k-1}} a_k^{b_k}r'' \mm
 a_k^{-l_k+i_k} \id a_1^{s_1} \ldots a_{k-1}^{s_{k-1}} a_k^{c_k}x \succ 
 a_1^{s_1} \ldots a_{k-1}^{s_{k-1}} a_k^{i_k} y \id t \mm z_{k-1} \mm
 a_k^{-l_k+i_k}$ for some $x,y \in \ord(\Sigma_{k+1})$ and therefore $c_k
  \geq_{\z} i_k$.
\auskommentieren{
Then the $t$-term $t \mm z_{k-1} \mm a_k^{-l_k + i_k}$
 in $p \mrm z_{k-1} \mrm a_k^{-l_k + i_k}$ can be modified 
 by ``filling'' up its prefix
 $a_1^{s_1} \ldots a_{k-1}^{s_{k-1}}$ and some additional
 modifications to equal $\hterm(p \mrm u)$.
Notice that it is  brought to head position  by these multiplications
 and they can be motivated as follows:
\\
There exists $z \in \ord(\Sigma \backslash \Sigma_k)$ such that 
 $a_1^{s_1} \ldots a_{k-1}^{s_{k-1}} \mm z \id a_1^{i_1} \ldots 
 a_{k-1}^{i_{k-1}}w$ and $w \id a_k^{f_k}w'$, $ w' \in \ord(\Sigma_{k+1})$.
Then by lemma \ref{lem.nilpotentmultiplication} since $z \in \ord(\Sigma \backslash
\Sigma_k)$ and $y \in \ord(\Sigma_{k+1})$ we know 
 $a_k^{i_k}y \mm z \id za_k^{i_k}y'$ for some 
 $y' \in \ord(\Sigma_{k+1})$ and on the other hand 
 $w \mm a_k^{i_k}y' = a_k^{f_k}w' \mm a_k^{i_k}y' \id a_k^{f_k +
   i_k}y''$
 for some $y'' \in \ord(\Sigma_{k+1})$ as $w' \in \ord(\Sigma_{k+1})$.
Furthermore, there exists $z'' \in \ord(\Sigma_{k+1})$ such that for
 $z' \id a_k^{-f_k}z''$ we have 
 $y'' \mm y' \mm z' \id a_k^{-f_k}a_{k+1}^{i_{k+1}} \ldots a_n^{i_n}$.
Now these observations can be combined to
\begin{tabbing}
$a_1^{s_1} \ldots a_{k-1}^{s_{k-1}} a_k^{i_k} y \mm z \mm z'$ 
 \= = \= 
 $a_1^{s_1} \ldots a_{k-1}^{s_{k-1}} \mm za_k^{i_k} y' \mm z'$ \\
 \> = \> $a_1^{i_1} \ldots  a_{k-1}^{i_{k-1}}w \mm a_k^{i_k}y' \mm z'$ \\
 \> = \> $a_1^{i_1} \ldots  a_{k-1}^{i_{k-1}} a_k^{i_k +f_k}y'' \mm y' \mm z'$ \\
 \> $\id$ \> $a_1^{i_1} \ldots  a_{k-1}^{i_{k-1}} a_k^{i_k + f_k - f_k}a_{k+1}^{i_{k+1}} \ldots a_n^{i_n}.$ 
\end{tabbing}
On the other hand, multiplying $a_1^{s_1} \ldots a_{k-1}^{s_{k-1}}
 a_k^{c_k}x$ by $z \mm z'$, there exist elements
 $x', x'', \tilde{x} \in \ord(\Sigma_{k+1})$ such that 
 $a_k^{c_k}x \mm z \id za_k^{c_k}x'$,
 $w \mm a_k^{c_k}x' = a_k^{f_k}w' \mm a_k^{c_k}x' \id a_k^{c_k + f_k}x''$ and $x'' \mm z' = x'' \mm
 a_k^{-f_k}z'' \id
 a_k^{f_k}\tilde{x}$ and the following holds: 
\begin{tabbing}
$a_1^{s_1} \ldots a_{k-1}^{s_{k-1}} a_k^{c_k} x \mm z \mm z'$ 
 \= = \= 
 $a_1^{s_1} \ldots a_{k-1}^{s_{k-1}} \mm za_k^{c_k} x' \mm z'$ \\
 \> = \> $a_1^{i_1} \ldots  a_{k-1}^{i_{k-1}}w \mm a_k^{c_k}x' \mm z'$ \\
 \> = \> $a_1^{i_1} \ldots  a_{k-1}^{i_{k-1}} a_k^{c_k +f_k}x'' \mm a_k^{-f_k}z''$ \\
 \> $\id$ \> $a_1^{i_1} \ldots  a_{k-1}^{i_{k-1}} a_k^{c_k + f_k - f_k}\tilde{x}$ 
\end{tabbing}
implying $c_k \leq_{\z} i_k$ and thus $c_k = i_k$ because $c_k \geq_{\z} i_k$
 also holds.
}
Then by lemma \ref{lem.nilpotentmultiplication} 
 there exist $u_1 \in \ord(\Sigma \backslash \Sigma_{k-1})$ and
 $u_2 \in \ord(\Sigma_k)$ such that 
 $a_1^{s_1} \ldots a_{k-1}^{s_{k-1}}a_k^{i_k}y \mm u_1 \id a_1^{i_1} \ldots
 a_{k-1}^{i_{k-1}} a_k^{i_k + f_k}z$ for some $z \in \ord(\Sigma_{k+1})$ and
 $a_k^{i_k + f_k}z \mm u_2 \id a_k^{i_k}a_{k+1}^{i_{k+1}}$, i.e., $u_2 \id
 a_k^{-f_k}u_2'$ for some $u_2' \in \ord(\Sigma_{k+1})$.
Note that the $t$-term is brought to head position by this
multiplication.
Now multiplying $\hterm(p \mrm z_{k-1} \mrm a_k^{-l_k+i_k})$ by
$u_1u_2$ we find $a_1^{s_1} \ldots a_{k-1}^{s_{k-1}} a_k^{c_k} x \mm
u_1u_2 \id a_1^{i_1} \ldots  a_{k-1}^{i_{k-1}} a_k^{c_k + f_k -
  f_k}\tilde{x}$ for some $\tilde{x} \in \ord(\Sigma_{k+1})$.
This gives us $c_k \leq_{\z} i_k$ and thus $i_k \leq_{\z} c_k$ yields
$c_k = i_k$.
\\
Finally, we have to check the case that $\sgn(i_k) \neq \sgn(l_k)$ and
$l_k \neq 0$.
Let us take a look at the polynomial $p \mrm z_{k-1} \mrm a_k^{-l_k}$,
 i.e., the exponent of the letter $a_k$ in the term
 $t \mm z_{k-1} \mm a_k^{-l_k}$ will be $0$.
Suppose 
 $\hterm(p \mrm z_{k-1} \mrm a_k^{-l_k}) \id a_1^{s_1} \ldots
a_{k-1}^{s_{k-1}}a_{k}^{c_k}x$,  for some term 
 $s \id a_1^{s_1} \ldots a_{k-1}^{s_{k-1}}a_{k}^{b_s}x_s 
    \in \terms(p \mrm z_{k-1})$,
 $x, x_s \in \ord(\Sigma_{k+1})$, i.e., $c_k = b_s - l_k$.
In case this head term is already the corresponding $t$-term 
 $t \mm z_{k-1} \mm a_k^{-l_k}$, we are
 done and we set  $w_1 = a_k^{-l_k}$ and $\tilde{s}_k = 0 = \tilde{\rho}_k$.
Now if we can show $c_k = 0$,
 by lemma \ref{lem.dletter} the $t$-term 
 $t \mm z_{k-1} \mm a_k^{-l_k}$ can be
brought to head position using an element in $\ord(\Sigma_{k+1})$ since
 the distinguishing letter between
 $\hterm(p \mrm z_{k-1} \mrm a_k^{-l_k})$ and the term
 $t \mm z_{k-1} \mm a_k^{-l_k}$ then has at least index $k+1$ and  we know
 $\hterm((p \mrm z_{k-1} \mrm a_k^{-l_k}) \mrm  a_k^{l_k}) =  \hterm(p
 \mrm z_{k-1}) =t \mm z_{k-1}$.
Hence, in showing that $c_k=0$ we are done.
\auskommentieren{
Remember that the $t$-term 
 $t \mrm z_{k-1} \mrm a_k^{-l_k}\id a_1^{s_1} \ldots
  a_{k-1}^{s_{k-1}}y$, $y \in \ord(\Sigma_{k+1})$ can be
 modified to equal $\hterm(p \mrm u)$ as before, i.e., there exist
 $z \in \ord(\Sigma \backslash \Sigma_k)$,$w, z'\in \ord(\Sigma_{k})$, $w',z'', y', y'' \in \ord(\Sigma_{k+1})$
 such that
 $a_1^{s_1} \ldots a_{k-1}^{s_{k-1}} \mm z \id a_1^{i_1} \ldots
 a_{k-1}^{i_{k-1}}w$, $w \id a_k^{f_k}w'$, $y \mm z \id zy'$, $w \mm
 y' = a_k^{f_k}w' \mm y' \id a_k^{f_k}y''$ and $y'' \mm z' \id a_k^{-f_k +
   i_k}a_{k+1}^{i_{k+1}} \ldots a_n^{i_n}$ for $z' \id a_k^{-f_k +
   i_k}z''$ giving us
\begin{tabbing}
$a_1^{s_1} \ldots a_{k-1}^{s_{k-1}}y \mm z \mm z'$ 
 \= = \= 
 $a_1^{s_1} \ldots a_{k-1}^{s_{k-1}} \mm zy' \mm z'$ \\
 \> = \> $a_1^{i_1} \ldots  a_{k-1}^{i_{k-1}}w \mm y' \mm z'$ \\
 \> = \> $a_1^{i_1} \ldots  a_{k-1}^{i_{k-1}} a_k^{f_k}y''\mm z'$ \\
 \> $\id$ \> $a_1^{i_1} \ldots  a_{k-1}^{i_{k-1}}a_{k}^{i_{k}} \ldots a_n^{i_n}.$ 
\end{tabbing}
Similarly, multiplying $\hterm(p \mrm z_{k-1} \mrm a_k^{-l_k})$ by $z
 \mm z'$ there exist $x',x''  \in\ord(\Sigma_{k+1})$ such that
 $a_k^{c_k}x \mm z \id za_k^{c_k}x'$, $w \mm a_k^{c_k}x' = a_k^{f_k}w'
 \mm a_k^{c_k}x' \id a_k^{c_k
   + f_k}x''$ and $x'' \mm z' = x'' \mm a_k^{-f_k + i_k}z'' \id
 a_k^{-f_k + i_k}\tilde{x}$ and we get
\begin{tabbing}
$a_1^{s_1} \ldots a_{k-1}^{s_{k-1}}a_{k}^{c_{k}}x \mm z \mm z'$ 
 \= = \= 
 $a_1^{s_1} \ldots a_{k-1}^{s_{k-1}} \mm z a_{k}^{c_{k}}x' \mm z'$ \\
 \> = \> $a_1^{i_1} \ldots  a_{k-1}^{i_{k-1}}w \mm a_{k}^{c_{k}}x' \mm z'$ \\
 \> = \> $a_1^{i_1} \ldots  a_{k-1}^{i_{k-1}} a_k^{c_k + f_k}x''\mm z'$ \\
 \> $\id$ \> $a_1^{i_1} \ldots  a_{k-1}^{i_{k-1}}a_{k}^{c_k + i_{k}}\tilde{x}.$ 
\end{tabbing}
}
As before there exist $u_1 \in \ord(\Sigma \backslash \Sigma_{k-1})$
and $u_2 \in \ord(\Sigma_k)$ such that
 $t \mm z_{k-1} \mm a_k^{-l_k} \mm u_1u_2 \id a_1^{i_1} \ldots  a_n^{i_n}$.
Remember that this multiplication brings the $t$-term to head
position.
Since the exponent of $a_k$ in the term $t \mm z_{k-1} \mm a_k^{-l_k}$
is $0$,  multiplying $\hterm(p \mrm z_{k-1} \mrm a_k^{-l_k})$ by $u_1u_2$
we find
 $a_1^{s_1} \ldots a_{k-1}^{s_{k-1}}a_{k}^{c_{k}}x \mm u_1u_2 \id
 a_1^{i_1} \ldots  a_{k-1}^{i_{k-1}}a_{k}^{c_k + i_{k}}\tilde{x}$ for
 some
 $\tilde{x} \in \ord(\Sigma_{k+1})$. 
Thus we know $c_k + i_k \leq_{\z} i_k$.
To see that this implies $c_k = 0$ we have to distinguish three cases.
Remember that $c_k = b_s - l_k$ and since our head term is an $s$-term $s \mm a_k^{-l_k}$ for some
 $s \in \terms(p\mrm z_{k-1})$ we know $b_s \leq_{\z} l_k$.
In case $i_k = 0$, we have $c_k \leq_{\z} 0$
 implying $c_k = 0$.
In case $i_k > 0$ then $c_k + i_k = b_s - l_k + i_k \leq_{\z} i_k$ implies
 $0 \leq b_s - l_k + i_k \leq i_k$.
Furthermore, as $l_k < 0$ we have $-l_k + i_k > i_k$ implying $b_s < 0$ and hence 
 $|b_s| \leq |l_k|$.
But then $b_s - l_k \geq 0$ and $0 \leq b_s - l_k + i_k \leq i_k$ yields $c_k = b_s - l_k = 0$.
On the other hand, $i_k < 0$ and $l_k > 0$ imply $0 \leq b_s \leq l_k$
and hence $b_s - l_k + i_k <0$ yielding $|b_s - l_k + i_k| \leq |i_k|$.
Since $b_s - l_k \leq 0$ this inequation can only hold in case $c_k = b_s - l_k = 0$.
\\
\qed
%
\begin{corollary}~\\
{\sl
For $p \mrm v$, $p \mrm u$ and $q = p \mrm \tilde{z}$ in the previous lemma we have $p \mrm u \red{}{\myr}{qc}{q} 0$ and $p \mrm v \red{}{\myr}{qc}{q} 0$.
}
\end{corollary}
\Ba{}~\\
This follows immediately, as for the term $\tilde{t} = \hterm(q)$ we
have $\tilde{t} \tupleq \hterm(p \mrm u)$ and $\tilde{t} \tupleq
\hterm(p \mrm v)$. 
\\
\qed
\begin{corollary}~\\
{\sl
Let $p$ be a non-zero polynomial in $\myk[\g]$ and $t \in \terms(p)$ such that $Y_t \neq \emptyset$.
Then $Y_t$ contains a polynomial $q$ such that for all $q'
\in Y_t$ we have $q' \red{}{\myr}{qc}{q} 0$. 
}
\end{corollary}
\Ba{}~\\
The set of head terms that are $t$-terms can be
ordered with respect to $\tupeq$.
Now suppose there are two different polynomials  $p \mrm u, p \mrm v$ in $Y_t$, both having minimal
 head terms with respect to $\tupeq$.
Then lemma \ref{lem.uniqueCNI} yields $t \mm u = t \mm v$ and hence $u = v$ contradicting our assumption.
\\
\qed
\auskommentieren{
\begin{remark}~\\
{\rm
To compute a quasi-commutatively saturating set for a polynomial
 $p \in \myk[\g]$ we can proceed as follows:
\begin{enumerate}
\item For each $Y_t$ compute a representative $p_t$.
\item For each term $s$ with $\hterm(p_t) \tupgreater s $  compute $p_t
  \mrm (\inv{\hterm(p_t)} \mm s)$ and in case the head
  term of this multiple is the $t$-term and cannot be quasi-reduced to
  zero in one step using the already computed polynomials, add it.
  Reviewing the proof of lemma \ref{lem.uniqueCNI} gives additional information
  on reducing the test set of such terms.
\remend
\end{enumerate}
}
\end{remark}}
Now quasi-commutatively saturating sets can be constructed by
computing the appropriate minimal polynomials for the non-empty sets $Y_t$.

\procedure{Quasi-Commutative Saturation}%
{\vspace{-4mm}\begin{tabbing}{ll}
XXXXX\=XXXX \kill
\removelastskip
{\bf Given:} \> A polynomial $p \in
                \myk[\g]$ and \\
             \>  $(\Sigma, T)$ a convergent CNI-Presentation
                of $\g$. \\
{\bf Find:} \> $S_{qc} \in \SAT_{qc}(p)$.
\end{tabbing}
\vspace{-7mm}
\begin{tabbing}
XX\=XX\=XXXX\=XX\=XX\=XXXX\=XXX\=XXX\= \kill
{\bf for all} $t \in \terms(p)$ {\bf do} \\
\> $S_t$ := $\emptyset$; \\
\> {\bf if}   \> $Y_t \neq \emptyset$ \\
\>            \>{\bf then}  \> compute $q = p \mrm w \in Y_t$ as described in lemma
\ref{lem.dletter} \\
\>        \>    \>  $H_t$ := $\{ s \in \g | \hterm(q) \tupeq s
\}$;\\
\>       \>     \> {\rm\kommentar \% These are candidates for smaller
  polynomials in $Y_t$}\\
\>       \>     \> $q$ := $\min \{\{ p \mrm (\inv{t} \mm s) \mid s \in
H_t \} \cap Y_t\}$;\\
\>        \>    \>  $S_t$ := $\{ q \}$; \\
\> {\bf endif} \\
{\bf endfor} \\
$S_{qc}$ := $\bigcup_{t \in \terms(p)} S_t$
\end{tabbing}}

Quasi-commutative saturation enriches a polynomial $p$ by adding a set
 of polynomials  $S \in \SAT_{qc}(p)$ such that we can
 substitute $q \red{}{\myr}{(s,r)}{p}q'$ by $q \red{}{\myr}{qc}{p' \in S}q'$.
Therefore, we have more information on the reduction step than using
 (strong) right reduction and we define s-polynomials corresponding to
 this reduction.
\begin{definition}\label{def.cpqc}~\\
{\rm
Let $p_{1}, p_{2}$ be two polynomials  in $\myk[\g]$ such that 
 $\hterm(p_1) \id a_1^{i_1} \ldots a_n^{i_n}$ and $\hterm(p_2) \id
               a_1^{j_1} \ldots a_n^{j_n}$ with either
 $i_l = 0$ or $j_l = 0$ or $\sgn(i_l)=\sgn(j_l)$ for $1 \leq l \leq n$.
Then  setting
\[\rho_l = \left\{ \begin{array}{l@{\quad\quad}l}
         \sgn(j_l) & i_l = 0 \\
         \sgn(i_l) &  \mbox{otherwise}
 \end{array} \right. \]
 the  situation
$${\sf qclcm}(t_1,t_2) = a_1^{\rho_1 \cd \max \{ |i_1|,|j_1| \}}\ldots a_n^{\rho_n 
  \cd \max \{ |i_n|,|j_n| \}}= t_1 \mm w_1 = t_2 \mm w_2$$
 for some $w_1,w_2 \in \g$ defines a 
 \index{quasi-commutative!s-polynomial}\index{s-polynomial!quasi-commutative}
 \betonen{quasi-commutative s-polynomial}
 $$ \spol{qc}(p_{1}, p_{2}) = \hc(p_1)^{-1} \skm p_1 \mrm
 w_1 - \hc(p_2)^{-1} \skm p_2 \mrm w_2.$$
\dend
}
\end{definition}
A quasi-commutative s-polynomial is called non-trivial in case it is non-zero.
Furthermore, we get 
 $\hterm(\spol{qc}(p_{1}, p_{2})) \prec {\sf qclcm}(t_1,t_2)$ and
 $\hterm(p_i) \tupleq {\sf qclcm}(t_1,t_2)$.
Notice that a finite set $F \subseteq \myk[\g]$ only gives us finitely
many such s-polynomials.
As before  Gr\"obner bases cannot be
characterized by these s-polynomials alone unless they are quasi-commutatively
 saturated sets.
\begin{theorem}\label{theo.qcpc}~\\
{\sl
For a quasi-commutatively saturated set $F$ of polynomials in $\myk[\g]$, the
 following statements are equivalent:
\begin{enumerate}
\item For all polynomials $g \in \ideal{r}{}(F)$ we have $g \red{*}{\myr}{qc}{F} 0$.
\item For all polynomials $f_{k}, f_{l} \in F$ we have
  $\spol{qc}(f_{k}, f_{l}) \red{*}{\myr}{qc}{F} 0$.
\end{enumerate}
}
\end{theorem}
\Ba{}~\\
\mbox{$1 \R 2:$ }
 Let $\hterm(f_{k}) \mm w_k = \hterm(f_{l}) \mm w_{l}$ for
  $w_{k}, w_{l} \in \g$  such that we have an overlap as described in definition \ref{def.cpqc}.
 We get the following
 $$ \spol{qc}(f_{k}, f_{l}) = \hc(f_k)^{-1} \skm f_k \mrm
 w_k - \hc(f_l)^{-1} \skm f_l \mrm w_l \:\in \ideal{r}{}(F),$$
 and hence $\spol{qc}(f_{k}, f_{l}) \red{*}{\myr}{qc}{F} 0$.

\mbox{$2 \R 1:$ }
    We have to show that every non-zero element  $g \in
    \ideal{r}{}(F)$ is $\red{}{\myr}{qc}{F}$-reducible to zero.
    Remember that for
      $h \in \ideal{r}{}(F)$, $ h \red{}{\myr}{qc}{F} h'$ implies $h' \in \ideal{r}{}(F)$.
     As  $\red{}{\myr}{qc}{F}$ is Noetherian
      it suffices to show that every  $g \in \ideal{r}{}(F) \backslash \{ 0 \}$ is $\red{}{\myr}{qc}{F}$-reducible.
     Let $g = \sum_{j=1}^m \alpha_{j} \skm f_{j} \mrm w_{j}$  be a
      representation of a non-zero polynomial $g$ with $\alpha_{j} \in \myk^*, f_j \in F,
      w_{j} \in \g$.
     By lemma \ref{lem.propqc} we can assume  $\hterm(f_{i} \mrm w_{i})
     = \hterm(f_{i}) \mm w_{i} \tupeq \hterm(f_i)$.
     Depending on this  representation of $g$ and our well-founded total ordering $\syll$ on $\g$ we define
      $t = \max \{ \hterm(f_{j}) \mm w_{j} \mid j \in \{ 1, \ldots m \}  \}$ and
      $K$ is the number of polynomials $f_j \mrm w_j$ containing $t$ as a term.
Then $t \succeq \hterm(g)$ and
in case $\hterm(g) = t$ this immediately implies that $g$ is
$\red{}{\myr}{qc}{F}$-reducible. 
So we  show that
$g$ has a special  representation (a standard representation
corresponding to quasi-commutative reduction) where all terms are
bounded by $\hterm(g)$, as this implies that $g$ is
top-reducible using $F$ since then $\hterm(g) = t$.
This will be done by induction on $(t,K)$, where
      $(t',K')<(t,K)$ if and only if $t' \syllless t$ or $(t'=t$ and $K'<K)$\footnote{Note
        that this ordering is well-founded since $\syll$ is and $K \in\n$.}.
If $t \succ \hterm(g)$
      there are two polynomials $f_k,f_l$ in the corresponding
       representation\footnote{Not necessarily $f_l \neq f_k$.}
      with $\hterm(f_k) \mm w_k = \hterm(f_l) \mm w_l$ and $t \tupeq \hterm(f_k)$,$t \tupeq \hterm(f_l)$.
     By definition \ref{def.cpqc} we then have an s-polynomial
      \mbox{$\spol{qc}(f_k,f_l) = \hc(f_k)^{-1} \skm  f_k \mrm z_1-
      \hc(f_l)^{-1} \skm f_l \mrm z_2$} such that $t \tupeq \hterm(f_k) \mm
      z_1$, $t \tupeq \hterm(f_l) \mm z_2$
      and  $\hterm(f_k) \mm w_k = \hterm(f_l)
      \mm w_l =  \hterm(f_k) \mm z_1 \mm w = \hterm(f_l) \mm z_2 \mm w$ for
      some $z_1,z_2,w \in \g$.
     Let us assume  $\spol{qc}(f_k,f_l) \neq 0$\footnote{In case  $\spol{qc}(f_k,f_l) = 0$,
               just substitute $0$ for the sum $\sum_{i=1}^n \delta_i \skm h_i \mrm v_i$ in the equations below.}.
     Hence,  $\spol{qc}(f_k,f_l) \red{*}{\myr}{qc}{F} 0 $ implies
     $\spol{qc}(f_k,f_l) =\sum_{i=1}^n \delta_i \skm h_i \mrm v_i,\delta_i \in \myk^*,h_i \in F,v_i \in \g$,
      where the $h_i$ are due to the reduction of the s-polynomial
      and all terms occurring in the sum are bounded by $\hterm(\spol{qc}(f_k,f_l))$.
     Since $t \tupeq \hterm(f_k) \mm z_1$ and $t = \hterm(f_k) \mm z_1 \mm w$ by lemma \ref{lem.order}
      we can conclude that $t$ is a proper bound for all terms occurring
      in the sum $\sum_{i=1}^n \delta_i \skm h_i \mrm v_i \mrm w$.
We can assume that this representation is of the required form,
      as we can substitute all polynomials $h_i$ violating
      $\hterm(h_i \mrm v_i \mrm w_k) \not{\tupeq}  \hterm(h_i)$  without increasing
      $t$ or $K$.
     This gives us:
     \begin{eqnarray}
       &  & \alpha_{k} \skm f_{k} \mrm w_{k} + \alpha_{l} \skm f_{l} \mrm w_{l}  \nonumber\\
       &  &    \nonumber \\
       & = &  \alpha_{k} \skm f_{k} \mrm w_{k} + \underbrace{
         \alpha'_{l} \skm \beta_k \skm f_{k} \mrm w_{k}
                   - \alpha'_{l} \skm \beta_k \skm f_{k} \mrm w_{k}}_{=\, 0}
                   + \alpha'_{l}\skm \beta_l  \skm f_{l} \mrm w_{l} \nonumber\\
       &  &    \nonumber \\
       & = & (\alpha_{k} + \alpha'_{l} \skm \beta_k) \skm f_{k} \mrm w_{k} - \alpha'_{l} \skm
             \underbrace{(\beta_k \skm f_{k} \mrm w_{k}
             -  \beta_l \skm f_{l} \mrm w_{l})}_{=\, \spol{qc}(f_k,f_l) \mrm w} \nonumber\\
       & = & (\alpha_{k} + \alpha'_{l} \skm b_k) \skm f_{k} \mrm w_{k} - \alpha'_{l} \skm
                   (\sum_{i=1}^n \delta_{i} \skm h_{i} \mrm v_{i} \mrm w) \label{number.qc}
     \end{eqnarray}
     where  $\beta_k = \hc(f_k)^{-1}$, $\beta_l = \hc(f_l)^{-1}$ and
       $\alpha'_l \skm \beta_l = \alpha_l$.
     By substituting (\ref{number.qc}) in our representation of $g$ either $t$ disappears or in
 case $t$ remains maximal among the terms occurring in the new
 representation of $g$, $K$ is decreased.
\\
\qed
Due to Dickson's lemma finite quasi-commutative Gr\"obner bases exist.
\begin{lemma}\label{lem.term.np}~\\
{\sl
Let $\g$ be a nilpotent group with a convergent CNI-presentation.
Then every Gr\"obner basis with respect to $\red{}{\myr}{qc}{}$ of a 
 finitely generated right ideal contains a finite one.
\lemend
}
\end{lemma}
\Ba{}~\\
Let $F$ be a finite subset of $\myk[\g]$ and $G$ a infinite Gr\"obner
basis of $\ideal{r}{}(F)$ with respect to quasi-commutative reduction.
Further let
 $H = \{ \hterm(g) \mid g \in G \} \subseteq \g$.
Then for every polynomial $f \in \ideal{r}{}(F)$ there exists a term
$t \in H$ such that $\hterm(f) \tupeq t$.
Each element of $H$ can be viewed as an n-tuples over $\z$ as it is
presented by an ordered group word.
But we can also view it as a 2n-tuples over $\n$ by representing
 each element $u \in H$ by an extended ordered group word
 $u \id a_1^{-i_1} a_1^{j_1} \ldots  a_n^{-i_n} a_n^{j_n}$,
 where $i_l,j_l \in \n$ and the
 representing 2n-tuple is $(i_1,j_1, \ldots , i_n,j_n)$.
Notice that at most one of the two exponents $i_l$ and $j_l$ is non-zero.
Now $H$ can be considered as a (possibly infinite) subset of a free commutative
 monoid  ${\cal T}_{2n}$ with $2 \skm n$ generators.
Thus by  Dickson's lemma there exists 
 a finite subset $B$ of $H$
 such that for every $w \in H$ there is a $b \in B$ with
 $w = b \mm_{{\cal T}_{2n}} u$ for some $u \in {\cal T}_{2n}$,
 and hence $w \tupeq b$.
Now we can use the set $B$ to distinguish a finite Gr\"obner basis in
$G$ as follows.
To each term $t \in B$ we can assign a polynomial $g_t \in G$ such
that $\hterm(g_t) = t$.
Then the set $G_B = \{ g_t \mid t \in B \}$ is again a Gr\"obner basis
with respect to quasi-commutative reduction since for every polynomial
$f \in \ideal{r}{}(F)$ there still exists a polynomial $g_t$ now in $G_B$
such that $\hterm(f) \tupeq \hterm(g_t) = t$.
Hence all polynomials in $\ideal{r}{}(F)$ are quasi-commutatively
reducible to zero using $G_B$.
\\
\qed
Finite Gr\"obner bases with respect to $\red{}{\myr}{qc}{}$ can now be
computed as follows:

\procedure{Quasi-Commutative Gr\"obner Bases\protect{\label{quasi.commutative.groebner.bases}}}%
{\vspace{-4mm}\begin{tabbing}
XXXXX\=XXXX \kill
\removelastskip
{\bf Given:} \> A finite set of polynomials $F \subseteq \myk[\g]$. \\
{\bf Find:} \> $\gb(F)$, a  quasi-commutative Gr\"obner basis of $\ideal{r}{}(F)$. \\
{\bf Using:} \> $\s_{qc}$ a quasi-commutatively saturating procedure for polynomials.
\end{tabbing}
\vspace{-7mm}
\begin{tabbing}
XX\=XX\=XXXX\=XX\=XXXX\=XXXX\kill 
$G$ := $\bigcup_{f \in F} \s_{qc}(f)$; \\
$B$ := $\{ (q_{1}, q_{2}) \mid q_{1}, q_{2} \in G, q_{1} \neq q_{2} \}$; \\
{\bf while} $B \neq \emptyset$ {\bf do} \\
\>{\rm\kommentar \% Test if
  statement 2 of theorem \ref{theo.qcpc} is valid}\\
\>  $(q_{1}, q_{2})$ := {\rm remove}$(B)$; \\
\> {\rm\kommentar \% Remove an element using a fair strategy} \\
\> {\bf if}  \>  $\spol{qc}(q_{1}, q_{2})$ exists \\
\>                \> {\rm\kommentar \% The s-polynomial is not trivial}\\
\>           \> {\bf then}\>   $h$ := ${\rm normalform}(\spol{qc}(q_{1},
                              q_{2}),\red{}{\myr}{qc}{G})$; \\
\>           \>     \> {\rm\kommentar \% Compute a normal form using
        quasi-commutative reduction} \\
\>          \>       \>  {\bf if} \> $h \neq 0$ \\
\>          \>       \>\>{\rm\kommentar \% Statement 2 of theorem
  \ref{theo.qcpc} does not hold}\\
\>          \>       \>           \> {\bf then} \> $G$ := $G \cup \s_{qc}(h)$; \\
\> \>  \>\>                 \> {\rm\kommentar \% $G$ is quasi-commutatively saturated} \\
\>          \>       \>           \>            \> $B$ := $B \cup \{ (f, {\tilde h}), ( {\tilde h},f) \mid f \in G, {\tilde h} \in  \s_{qc}(h) \}$; \\
\>          \>       \> {\bf endif} \\
\> {\bf endif} \\
{\bf endwhile}  \\
$\gb(F):= G$
\end{tabbing}}

It is possible to use the concept of quasi-commutative representations
 to introduce
interreduction to this setting as well.

We will proceed to show how Gr\"obner bases can be presented for arbitrary
nilpotent group rings.
As stated in the introductory chapter,
 a finitely generated  arbitrary nilpotent group $\g$ is a group
 containing a torsion-free subgroup of finite index.
Therefore, we can apply the approach used for context free groups to this
 situation and see how completion in a nilpotent group can be
 reduced to completion in the torsion-free quotient.

Let our group be given by a
 torsion-free subgroup subgroup ${\cal N}$ with a convergent 
 CNI-presentation
 $(\Sigma , C \cup I)$ and ${\cal E}$ a
 finite group such that $({\cal E}\backslash\{ \lambda \}) \cap \Sigma =
 \emptyset$ and $\g/{\cal N} \cong {\cal E}$.
Then every element $g \in \g$ can be uniquely expressed
 in the form $g \id ew$ where $e \in \g/{\cal N}$ and $w$ is an
 ordered group word in ${\cal N}$.
For all $e \in {\cal E}$ let $\phi_e : \Sigma \myr {\cal N}$ be a function
such that $\phi_{\lambda}$ is the inclusion and for all $a \in \Sigma$,
 $\phi_e(a) = \inv{e}\mm_{\g} a \mm_{\g} e$.
For all $e_1,e_2 \in {\cal E}$ let $z_{e_1,e_2} \in {\cal N}$ such that
 $z_{e_1,\lambda} \id z_{\lambda,e_1} \id \lambda$ and for all
 $e_1,e_2, e_3 \in {\cal E}$
 with $e_1 \mm_{\cal E} e_2 =_{\cal E} e_3$, $e_1 \mm_{\g} e_2 \id e_3z_{e_1,e_2}$. 
Let $\Gamma = ({\cal E} \backslash \{ \lambda \}) \cup \Sigma$
 and let $T$ contain the sets of rules $C$ and $I$, and
 the additional rules: 
\begin{tabbing}
XX\=XXXX\=XXX\=XXXXXXX\= XXXXXXXXXXXXXXXXXXXXXXXX\= \kill
\>$e_1e_2$ \>  $\myr$ \>  $e_3z_{e_1,e_2}$     \> 
   for all $e_1,e_2 \in {\cal E} \backslash \{ \lambda \}, e_3
   \in {\cal E}$ such that $e_1 \mm_{\cal E} e_2 =_{\cal E} e_3$, \\
\>$ae$ \>  $\myr$ \>  $e\phi_e(a)$ \> for all
        $e \in {\cal E} \backslash \{ \lambda \}, a \in \Sigma$.
\end{tabbing}

Then $(\Gamma,T)$ is a canonical presentation of $\g$ as an extension of ${\cal
  N}$ by ${\cal E}$.
The elements of our group $\g$
 are words of the form $eu$ where
 $e \in {\cal E}$ and $u \in {\cal N}$.
We can specify a total well-founded ordering on our group by combining a
 total well-founded ordering $\succeq_{\cal E}$ on ${\cal E}$ and
 a syllable ordering $\sylleq$ on ${\cal N}$:
Let $e_1u_1,e_2u_2 \in \g$ such that
 $e_i \in {\cal E}$, $u_i \in {\cal N}$.
Then we define $e_1u_1 \succ e_2u_2$ if and only if $e_1 \succ e_2$ 
  or $(e_1 = e_2$ and 
             $u_1 \syll  u_2)$.

For every $e \in {\cal E}$ let the mapping
 $\psi_e: \myk[\g] \myr \g$ be defined
 by $\psi_e(f) := f \mrm e$ for $f \in \myk[\g]$.
We now can give a characterization of  Gr\"obner bases by
 transforming a finite generating set for a right ideal using these finitely many
 mappings and then applying our results for finitely generated
 torsion-free nilpotent groups to this modified generating set.
We start by modifying quasi-commutative reduction in order to enable a lifting of
 the characterization of Gr\"obner bases in terms of special
 s-polynomials and restricted saturation.
The tuple ordering can be extended by setting $e_1u_1 \tupeq e_2u_2$
if and only if $e_1 = e_2$ and $u_1 \tupeq u_2$.
Similarly the concept of quasi-commutative right reduction is extended.
\begin{definition}\label{def.redpqc}~\\
{\rm
Let $p, f$ be two non-zero polynomials  in $\myk[\g]$. 
We say $f$ 
 \index{quasi-commutatively!right reduction}\index{reduction!quasi-commutative}\betonen{quasi-commutatively  reduces} $p$ to $q$ at
 a monomial $\alpha \skm eu$ of $p$ in one step, denoted by $p \red{}{\myr}{qc}{f} q$, if
\begin{enumerate}
\item[(a)] $eu \tupeq eu'$, where $eu' \id \hterm(f)$, and
\item[(b)] $q = p - \alpha \skm \hc(f)^{-1} \skm f \mrm
  (\inv{\hterm(f)} \mm t)$.
\end{enumerate}
We write $p \red{}{\myr}{qc}{f}$ if there is a polynomial $q$ as defined
above  and $p$ is then called quasi-commutatively reducible by
$f$. 
Further we can define $\red{*}{\myr}{qc}{}, \red{+}{\myr}{qc}{}$,
 $\red{n}{\myr}{qc}{}$ as usual.
Quasi-commutative reduction by a set $F \subseteq \myk[\g]$ is denoted by
 $p \red{}{\myr}{qc}{F} q$ and abbreviates $p \red{}{\myr}{qc}{f} q$
 for some $f \in F$,
 which is also written as  $p \red{}{\myr}{qc}{f \in F} q$.
\dend
}
\end{definition}
Notice that if $f$ quasi-commutatively reduces $p$ at $\alpha \skm eu$
to $q$,
then  $eu$ is no longer in the set $\terms(q)$ and $p > q$.
\begin{definition}~\\
{\rm
Let ${\cal H}_1$ be a subgroup of a group ${\cal H}_2$ and $p$ a
 non-zero polynomial in $\myk[{\cal H}_2]\backslash \{ 0 \}$.
A set $S \subseteq \{ p \mrm w \mid w \in {\cal H}_1 \}$ is called a 
 \index{saturating set!${\cal H}_1$-quasi-commutatively}
 \betonen{${\cal H}_1$-quasi-commutatively saturating set} for $p$,
 if for all $w \in {\cal H}_1$, $p \mrm w \red{}{\myr}{qc}{S} 0$
A set of polynomials $F \subseteq \myk[{\cal H}_2]$ is called a
 \index{saturated set!${\cal H}_1$-quasi-commutatively}
 \betonen{${\cal H}_1$-quasi-commutatively saturated set}, if for all $f \in F$ and
 for all $w \in {\cal H}_1$, $f \mrm w \red{}{\myr}{qc}{F} 0$.
\phantom{XXX}\dend
}
\end{definition}
Applying this definition to the subgroup ${\cal N}$ of $\g$ we find
that the essential lemmata \ref{lem.dletter} and \ref{lem.propqc}
for describing and computing saturating sets in
torsion-free nilpotent groups can be applied to describe 
${\cal N}$-quasi-commutatively saturating sets. 
The definition of s-polynomials can be extended by demanding that the
${\cal E}$ part of the head terms must coincide.
\begin{definition}\label{def.cppqc}~\\
{\rm
Let $p_{1}, p_{2}$ be  two polynomials  in $\myk[\g]$ such that
 $\hterm(p_1) \id ea_1^{i_1} \ldots a_n^{i_n}$ and $\hterm(p_2) \id
               ea_1^{j_1} \ldots a_n^{j_n}$
 with either $i_l = 0$ or $j_l = 0$ or $\sgn(i_l)=\sgn(j_l)$
 for $1 \leq l \leq n$.
Setting
\[\rho_l = \left\{ \begin{array}{l@{\quad\quad}l}
         \sgn(j_l) & i_l = 0 \\
         \sgn(i_l) &  \mbox{otherwise}
 \end{array} \right. \]
 the  situation
$$ ea_1^{\rho_1 \cd \max \{ |i_1|,|j_1| \}}\ldots a_n^{\rho_n 
  \cd \max \{ |i_n|,|j_n| \}}= t_1 \mm w_1 = t_2 \mm w_2$$
 for some $w_1,w_2 \in {\cal N}$ defines a 
 \index{quasi-commutative!s-polynomial}\index{s-polynomial!quasi-commutative}
 \betonen{quasi-commutative s-polynomial}
 $$ \spol{qc}(p_{1}, p_{2}) = \hc(p_1)^{-1} \skm p_1 \mrm
 w_1 - \hc(p_2)^{-1} \skm p_2 \mrm w_2.$$
\dend
}
\end{definition}
We can now give a characterization of Gr\"obner bases in this setting.
Notice that  $\hterm(p_i)
\tupleq ea_1^{\rho_1 \cd \max \{ |i_1|,|j_1| \}}\ldots a_{n\phantom{1}}^{\rho_n 
  \cd \max \{ |i_n|,|j_n| \}}$ for $i \in \{ 1,2 \}$ holds in case such an s-polynomial exists.
Furthermore, if there exists a term $t$ such that $t \tupeq \hterm(p_1)\id ea_1^{i_1} \ldots a_n^{i_n}$ and
 $t \tupeq \hterm(p_2) \id ea_1^{j_1} \ldots a_n^{j_n}$ an s-polynomial always exists\footnote{Notice that the condition for the existence of an s-polynomial is fulfilled as the tuple-ordering requires that the exponent of a letter $a_i$ in the smaller term is either zero or has the same sign as the exponent of $a_i$ in the tuple-larger term.} and we even have $t \tupeq ea_1^{\rho_1 \cd \max \{ |i_1|,|j_1| \}}\ldots a_{n\phantom{1}}^{\rho_n 
  \cd \max \{ |i_n|,|j_n| \}}$.
For every $e \in{\cal E}$ let the mapping $\psi_e:\myk[\g] \myr \myk[\g]$ be defined
 by $\psi_e(f) = f \mrm e$ for $f \in \myk[\g]$.
We now can give a characterization of a right Gr\"obner basis in a familiar way after transforming a generating set for the right ideal using these mappings.
\begin{theorem}\label{theo.comp.arbitrary.nilpotent}~\\
{\sl
Let $
 F,G \subseteq \myk[\g]$ such that
 \begin{itemize}
 \item[(i)] $\ideal{r}{}(F) = \ideal{r}{}(G)$
 \item[(ii)] $\{ \psi_e(f) \mid f \in F, e \in {\cal E} \} \subseteq G$
 \item[(iii)] $G$ is ${\cal N}$- saturated.
 \end{itemize}
Then the following statements are equivalent:
\begin{enumerate}
\item For all polynomials $g \in \ideal{r}{}(F)$ we have $g \red{*}{\myr}{qc}{G} 0$.
\item For all polynomials $f_{k}, f_{l} \in G$ we have
  $\spol{}(f_{k}, f_{l}) \red{*}{\myr}{qc}{G} 0$.
\theoend
\end{enumerate}
}
\end{theorem}
\Ba{}~\\
\mbox{$1 \R 2:$ }
By definition \ref{def.cppqc} in case for $f_k,f_l \in G$ the s-polynomial exists we get
 $$\spol{}(f_{k}, f_{l}) = \hc(f_k)^{-1} \skm f_{k} \mrm w_1
    -\hc(f_l)^{-1} f_{l} \mrm w_2 \:\in \ideal{r}{}(G)= \ideal{r}{}(F),$$
   and then $\spol{}(f_{k}, f_{l}) \red{*}{\myr}{qc}{G} 0$.

\mbox{$2 \R 1:$ }
We have to show that every non-zero element $g \in \ideal{r}{}(F)$
 is $\red{}{\myr}{qc}{G}$-reducible to zero.
Without loss of generality we assume that $G$ contains no constant polynomials, as then we are done at once.
Remember that for
 $h \in \ideal{r}{}(F)= \ideal{r}{}(G)$, $ h \red{}{\myr}{qc}{G} h'$
 implies $h' \in \ideal{r}{}(G)= \ideal{r}{}(F)$.
Thus as  $\red{}{\myr}{qc}{G}$ is Noetherian
 it suffices to show that every 
 $g \in \ideal{r}{}(F)\backslash\{ 0 \}$ is $\red{}{\myr}{qc}{G}$-reducible.
Let $g = \sum_{j=1}^m \alpha_{j} \skm f_{j} \mrm w_{j}$ be a
  representation of a non-zero polynomial $g$ such that
  $\alpha_{j} \in \myk^*, f_j \in F, w_{j} \in \g$.
Further for all $1 \leq j \leq m$, let $w_j \id e_ju_j$, with
 $e_j \in {\cal E}$, $u_j \in {\cal N}$.
Then, we can modify our representation of $g$ to
 $g = \sum_{j=1}^m \alpha_j \skm \psi_{e_j}(f_j) \mrm u_j$.
Since $G$ is ${\cal N}$-saturated and 
 $\psi_{e_j}(f_j) \in G$ by definition \ref{def.satqc} there exists
 $g_j \in G$ such that $\psi_{e_j}(f_j) \mrm u_j \red{}{\myr}{qc}{g_j} 0$ and
 hence we can assume $g= \sum_{j=1}^m \alpha_j \skm
 g_j \mrm v_j$, where $\alpha_j \in \myk^*, g_j \in G, v_j \in {\cal N}$
 and  $\hterm(g_{j} \mrm v_{j}) = \hterm(g_{j}) \mm v_{j} \tupeq \hterm(g_j)$.
Depending on this representation of $g$ and our well-founded
 total ordering  on $\g$ we define
 $t = \max \{ \hterm(g_{j})\mm v_{j} \mid j \in \{ 1, \ldots m \}  \}$ and
 $K$ is the number of polynomials $g_j \mrm v_j$ containing $t$ as a term.
Then $t \succeq \hterm(g)$ and
in case $\hterm(g) = t$ this immediately implies that $g$ is
$\red{}{\myr}{qc}{G}$-reducible. 
Otherwise we  show that
$g$ has a special  representation (a standard representation
corresponding to qc-reduction) where all terms are
bounded by $\hterm(g)$, as this implies that $g$ is
top-reducible using $G$.
This will be done by induction on $(t,K)$, where
 $(t',K')<(t,K)$ if and only if $t' \prec t$ or
 $(t'=t$ and $K'<K)$\footnote{Note that this ordering is well-founded
                                    since $\sylleq$ is and $K \in\n$.}.
In case $t \succ \hterm(g)$
 there are two polynomials $g_k,g_l$ in the corresponding
 representation\footnote{Not necessarily $g_l \neq g_k$.}
 such that  $t = \hterm(g_k) \mm v_k = \hterm(g_l) \mm v_l$ and  we have
 $t \tupeq \hterm(g_k)$,$t \tupeq \hterm(g_l)$.
Hence by definition \ref{def.cppqc} there exists an s-polynomial
 \mbox{$\spol{}(g_k,g_l) = \hc(g_k)^{-1} \skm  g_k \mrm z_1-
 \hc(g_l)^{-1} \skm g_l \mrm z_2$} and  $\hterm(g_k) \mm v_k = \hterm(g_l)
 \mm v_l =  \hterm(g_k) \mm z_1 \mm w = \hterm(g_l) \mm z_2 \mm w \tupeq
 \hterm(g_k) \mm z_1= \hterm(g_l) \mm z_2$ for
 some $z_1,z_2,w \in {\cal N}$.
Let us assume  $\spol{}(g_k,g_l) \neq 0$\footnote{In case
 $\spol{}(g_k,g_l) = 0$,
 just substitute $0$ for $\sum_{i=1}^n \delta_i \skm h_i \mrm v'_i$ in
 the equations below.}.
Hence,  $\spol{}(g_k,g_l) \red{*}{\myr}{qc}{G} 0$ implies
 $\spol{}(g_k,g_l) =\sum_{i=1}^n \delta_i \skm h_i \mrm v'_i,\delta_i \in
 \myk^*,h_i \in G,v'_i \in {\cal N}$\footnote{Note that the
   case $v_i' \in {\cal E}$ cannot occur as it implies that $h_i$ is a constant
   polynomial and we assumed that $G$ does not contain constant polynomials.},
 where the $h_i$ are due to the qc-reduction of the s-polynomial
 and all terms occurring in the sum are bounded by $\hterm(\spol{}(g_k,g_l))$.
By lemma \ref{lem.order}, since $t = \hterm(g_k) \mm z_1 \mm w  \tupeq
 \hterm(g_k) \mm z_1$ and $\hterm(g_k) \mm z_1 \succ \hterm(\spol{}(g_k,g_l))$,  we can conclude that $t$ is a proper
 bound for all terms occurring
 in the sum $\sum_{i=1}^n \delta_i \skm h_i \mrm v'_i \mrm w$.
Since $w \in {\cal N}$ and $G$ is ${\cal N}$-saturated,
 without loss of generality we can assume that the
 representation has the the required form.
We now have:
\begin{eqnarray}
 &  & \alpha_{k} \skm g_{k} \mrm v_{k} + \alpha_{l} \skm g_{l} \mrm v_{l}  \nonumber\\
 &  &  \nonumber \\
 & = &  \alpha_{k} \skm g_{k} \mrm v_{k} +
        \underbrace{ \alpha'_{l} \skm \beta_k \skm g_{k} \mrm v_{k}
                   - \alpha'_{l} \skm \beta_k \skm g_{k} \mrm v_{k}}_{=\, 0}
                   + \alpha'_{l}\skm \beta_l  \skm g_{l} \mrm v_{l} \nonumber\\
 &  &  \nonumber \\
 & = & (\alpha_{k} + \alpha'_{l} \skm \beta_k) \skm g_{k} \mrm v_{k} - \alpha'_{l} \skm
        \underbrace{(\beta_k \skm g_{k} \mrm v_{k}
        -  \beta_l \skm g_{l} \mrm v_{l})}_{=\, \spol{}(g_k,g_l) \mrm w}\nonumber\\
 & = & (\alpha_{k} + \alpha'_{l} \skm \beta_k) \skm g_{k} \mrm v_{k} -
        \alpha'_{l} \skm
       (\sum_{i=1}^n \delta_{i} \skm h_{i} \mrm v'_{i} \mrm w) \label{arbitrary.nilpotent}
\end{eqnarray}
 where  $\beta_k = \hc(g_k)^{-1}$, $\beta_l = \hc(g_l)^{-1}$ and
 $\alpha'_l \skm \beta_l = \alpha_l$.
By substituting (\ref{arbitrary.nilpotent}) in our representation of $g$ either
 $t$ disappears or in
 case $t$ remains maximal among the terms occurring in the new
 representation of $g$, $K$ is decreased.
\\
\qed
On first sight this characterization might seem artificial. 
The crucial point is that in losing the property ``admissible'' for
our ordering, an essential lemma in Buchberger's context, namely that
$p \red{*}{\myr}{}{F} 0$ implies $p \mrm w \red{*}{\myr}{}{F} 0$ for any
term $w$ no longer holds.
Defining reduction by restricting ourselves to commutative prefixes we gain enough
structural information to weaken this lemma, but we have to do
additional work to still describe the right ideal congruence. 
One step is to close the set of polynomials generating the right ideal with
respect to the finite group ${\cal E}$:
For
a set of polynomials $F$  using the ${\cal E}$-closure $F_{\cal E} =
\{ \psi_e(f) \mid f \in F, e \in {\cal E} \}$ we can characterize the
right ideal generated by $F$ as a set of ${\cal N}$-right-multiples since
$\ideal{r}{}(F) = \{ \sum_{i = 1}^k \alpha_i \skm f_i \mrm u_i \mid
\alpha_i \in \myk, f_i \in F_{\cal E}, u_i \in {\cal N} \}$.
If we additionally incorporate the concept of ${\cal N}$-saturation,
qc-reduction can be used to express the right ideal congruence and then a right Gr\"obner basis
can be characterized as usual by s-polynomials.
Now, using the characterization given in theorem \ref{theo.comp.arbitrary.nilpotent} 
 we can state a procedure which
enumerates  right Gr\"obner bases  in nilpotent group rings:


\procedure{\sc Right Gr\"obner Bases in Nilpotent Group Rings}
{\begin{tabbing}
XXXXX\=XXXX \kill
\removelastskip
{\bf Given:} \> $F \subseteq \myk[\g]$ and a presentation of $\g$ by
                ${\cal E}$ and ${\cal N}$ as specified above\\
{\bf Find:} \> $\gb_{r}(F)$, a  right Gr\"obner basis of $\ideal{r}{}(F)$.\\
\\
XX\= XX\= XXXX \= XX \= XXXX  \=\kill
$G := \{ \psi_e(f) \mid f \in F, e \in {\cal E} \}$; 
{\rm\kommentar \phantom{X}\% $G$ contains $F_{\cal E}$} \\
$G$ := $\bigcup_{g \in G} \s_{}(g)$;
{\rm\kommentar \phantom{X}\% $G$ is ${\cal N}$-saturated and $\ideal{r}{}(F) = \ideal{r}{}(G)$}\\
$B$ := $\{ (q_{1}, q_{2}) \mid q_{1}, q_{2} \in G, q_{1} \neq q_{2} \}$; \\
{\bf while} $B \neq \emptyset$ {\bf do} 
{\rm\kommentar \phantom{X}\% Test if
  statement 2 of theorem \ref{theo.comp.arbitrary.nilpotent} is valid}\\
\>      $(q_{1}, q_{2})$ := ${\rm remove}(B)$;
{\rm\kommentar \phantom{X}\% Remove an element using a fair strategy}\\
\>      {\bf if} \> $h$ := $\spol{}(q_{1}, q_{2})$ exists  \\
\>               \>   {\bf then}  \> $h'$ := ${\rm normalform}(h, \red{}{\myr}{qc}{G})$; {\rm\kommentar
  \phantom{X}\% Compute a normal form}\\
\>               \>               \> {\bf if} \> $h' \neq 0$ 
{\rm\kommentar \phantom{X}\% The s-polynomial does not
  reduce to zero}\\
\> \>     \>               \>{\bf then}  \>$G$ := $G \cup \{ g \mid g \in \s_{}(h') \}$;\\
\>\>\>\>\>
{\rm\kommentar \% $G$ is ${\cal N}$-saturated and $\ideal{r}{}(F) = \ideal{r}{}(G)$}\\
\> \>       \>               \>          \>$B$ := $B \cup \{ (f, g ) \mid f \in G,  g \in \s_{}(h') \}$; \\
\>\>\> {\bf endif}\\
\> {\bf endif} \\
{\bf endwhile}\\
$\gb_{r}(F):= G$
\end{tabbing}
}

The set $G$ enumerated by this procedure fulfills the requirements of
 theorem \ref{theo.comp.arbitrary.nilpotent}, i.e., we have 
 $F_{\cal E} \subseteq G$ and the set $G$ at each stage generates 
 the right ideal $\ideal{r}{}(F)$ and is ${\cal N}$-saturated.
Using a fair strategy to remove elements from the test set $B$ ensures
 that for all polynomials entered into $G$ the s-polynomial is considered
 in case it exists.
Hence, in case the procedure terminates, it computes a right Gr\"obner basis.
Later on we will see that every right Gr\"obner basis contains a finite one and hence this procedure must terminate.
Let us first continue to show how similar to the case of solvable polynomial rings or skew
polynomial rings (\cite{Kr93,We92}), Gr\"obner bases of two-sided
ideals can be characterized by right Gr\"obner bases which have
additional properties.
We will call a set of polynomials a \betonen{Gr\"obner basis} with respect to
qc-reduction of the
two-sided ideal it generates, if it fulfills one of the equivalent
statements in the next theorem.
\begin{theorem}\label{theo.ideals}~\\
{\sl
For a set of polynomials $G \subseteq \myk[\g]$, assuming that $\g$ is
presented by $(\Gamma, T)$ as described above, the following properties
are equivalent:
\begin{enumerate}
\item $G$ is a right Gr\"obner basis and $\ideal{r}{}(G) =
  \ideal{}{}(G)$.
\item For all $g \in \ideal{}{}(G)$ we have $g \red{*}{\myr}{qc}{G} 0$.
\item $G$ is a right Gr\"obner basis and for all $w \in \g$, $g \in G$
  we have $w \mrm g \in \ideal{r}{}(G)$.
\item $G$ is a right Gr\"obner basis and for all $a \in \Gamma$, $g \in G$
  we have $a \mrm g \in \ideal{r}{}(G)$.
\theoend
\end{enumerate}
}
\end{theorem}
\Ba{}~\\
\mbox{$1 \R 2:$ }
Since $g \in \ideal{}{}(G) = \ideal{r}{}(G)$ and $G$ is a right
Gr\"obner basis, we are done.

\mbox{$2 \R 3:$ }
To show that $G$ is a right Gr\"obner basis we have to prove $\red{*}{\lr}{qc}{G} = \;\;\equiv_{\ideal{r}{}(G)}$ and for all $g \in \ideal{r}{}(G)$, $g \red{*}{\myr}{qc}{G} 0$.
The latter follows immediately since $\ideal{r}{}(G) \subseteq \ideal{}{}(G)$ and hence for all $g
\in \ideal{r}{}(G)$ we have $g \red{*}{\myr}{qc}{G} 0$.
The inclusion $\red{*}{\lr}{qc}{G} \subseteq\;\;\equiv_{\ideal{r}{}(G)}$
is obvious.
Hence let $f \equiv_{\ideal{r}{}(G)} g$, i.e., $f -g \in \ideal{r}{}(G)$.
But then we have $f-g \red{*}{\myr}{qc}{G} 0$ and hence by lemma
\ref{lem.confluentqc}  there exists a polynomial
$h \in\myk[\g]$ such that $f \red{*}{\myr}{qc}{G} h$ and $g \red{*}{\myr}{qc}{G} h$,
yielding $f \red{*}{\lr}{qc}{G} g$.
Finally, $w \mrm f \in \ideal{}{}(G)$ and $w \mrm f
\red{*}{\myr}{qc}{G} 0$ implies $w \mrm f \in \ideal{r}{}(G)$.

$3 \R 4:$ 
This follows immediately.

$4 \R 1:$ 
Since it is obvious that $\ideal{r}{}(G) \subseteq \ideal{}{}(G)$ it
remains to show that $\ideal{}{}(G) \subseteq \ideal{r}{}(G)$ holds.
Let $g \in \ideal{}{}(G)$, i.e., $g = \sum_{i=1}^n \alpha_i \skm u_i
\mrm g_i \mrm w_i$ for some $\alpha_i \in \myk$, $g_i \in G$ and $u_i,
w_i \in \g$.
We will show by induction on $|u_i|$ that for $u_i \in \g$, $g_i \in
G$, $u_i \mrm g_i \in \ideal{r}{}(G)$ holds. Then $g$ also has a
representation in terms of right multiples and hence lies in the right
ideal generated by $G$ as well.
In case $|u_i|=0$  we are immediately done.
Hence let us assume $u_i \id ua$ for some $a \in \Gamma$ and by our
assumption we know $a \mrm g_i \in \ideal{r}{}(G)$.
Let $a \mrm g_i = \sum_{j=1}^m \beta_j \skm g_j' \mrm v_j$ for some $\beta_j \in \myk$, $g_j' \in G$ and $v_j \in \g$.
Then we get $u_i \mrm g_i = ua \mrm g_i = u \mrm (a \mrm g_i) = u \mrm
(\sum_{j=1}^m \beta_j \skm g_j' \mrm v_j) = \sum_{j=1}^m \beta_j \skm
(u \mrm g_i') \mrm v_j$ and by our induction hypothesis $u \mrm g_j'
\in \ideal{r}{}(G)$ holds for every $1 \leq j \leq m$. 
Therefore, we can conclude $u_i \mrm g_i \in \ideal{r}{}(G)$.
\\
\qed
Statement 4 enables a constructive approach to use procedure {\sc
  Right Gr\"obner Bases in Nilpotent Group Rings} in order to compute
Gr\"obner bases of two-sided ideals and item 2 states that such bases
can be used to decide the membership problem for the two-sided ideal 
by using qc-reduction.
The following corollary of the previous two theorems will be the foundation
 of a procedure to compute two-sided Gr\"obner bases.

\begin{corollary}\label{cor.ideal}~\\
{\sl
Let $
 F,G \subseteq \myk[\g]$ such that
 \begin{itemize}
 \item[(i)] $\ideal{}{}(F) = \ideal{}{}(G)$
 \item[(ii)] $\{ \psi_e(f) \mid f \in F, e \in {\cal E} \} \subseteq G$
 \item[(iii)] $G$ is ${\cal N}$- saturated.
 \end{itemize}
Then the following statements are equivalent:
\begin{enumerate}
\item For all polynomials $g \in \ideal{}{}(F)$ we have $g \red{*}{\myr}{qc}{G} 0$.
\item \begin{enumerate}
       \item For all polynomials $f_{k}, f_{l} \in G$ we have
             $\spol{}(f_{k}, f_{l}) \red{*}{\myr}{qc}{G} 0$.
       \item For all $a \in \Gamma$, $g \in G$
             we have $a \mrm g \red{*}{\myr}{qc}{G} 0$.
       \end{enumerate}
\end{enumerate}
}
\end{corollary}
\Ba{}~\\
\mbox{$1 \R 2:$ }
By definition \ref{def.cppqc} we find that in case for $f_k,f_l \in G$ 
 an s-polynomial exists,
 $$\spol{}(f_{k}, f_{l}) = \hc(f_k)^{-1} \skm f_{k} \mrm w_1
    -\hc(f_l)^{-1} f_{l} \mrm w_2 \:\in \ideal{}{}(G)= \ideal{}{}(F),$$
  and then $\spol{}(f_{k}, f_{l}) \red{*}{\myr}{qc}{G} 0$.
Similarly, since $g \in G$ implies $a \mrm g \in \ideal{}{}(G) = \ideal{}{}(F)$ for all $a \in \Gamma$, we have $a \mrm g\red{*}{\myr}{qc}{G} 0$.

\mbox{$2 \R 1:$ }
We have to show that every non-zero element $g \in \ideal{}{}(F)$
 is $\red{}{\myr}{qc}{G}$-reducible to zero.
Without loss of generality we assume that $G$ contains no constant polynomials, as then we are done at once.
Let $g = \sum_{j=1}^m \alpha_{j} \skm u_j \mrm f_{j} \mrm w_{j}$ be a
  representation of such a non-zero polynomial $g$ such that
  $\alpha_{j} \in \myk^*, f_j \in F, u_j, w_{j} \in \g$ and suppose for $1 \leq j \leq m$ we have $w_j \id e_jv_j$ with $e_j \in {\cal E}$ and $v_j \in {\cal N}$.
Then we can modify this representation to $g = \sum_{j=1}^m \alpha_{j} \skm u_j \mrm \psi_{e_j}(f_{j}) \mrm v_{j}$ as $\psi_{e_j}(f_j) \in G$ by our assumption.
Next we will show that every multiple $u_j \mrm \psi_{e_j}(f_{j})$   has a representation $u_j \mrm \psi_{e_j}(f_{j}) =  \sum_{i=1}^{m_j} \beta_i \skm g_i \mrm v_i'$ with $\beta_i \in \myk^*$, $g_i \in G$ and  $v_i' \in {\cal N}$.
More general, we will show that this is true for every multiple $u \mrm g$, $u \in \g$, $g \in G$.
As in the previous theorem this will be done by  induction on $|u|$.
The case $|u| = 0$ is obvious.
Hence let $u  \id u'a$ for some $a \in \Gamma$.
By our assumption we know $a \mrm g \red{*}{\myr}{qc}{G} 0$ and as we assume that
 $G$ does not contain constant polynomials, this reduction sequence results in a representation $a \mrm g = \sum_{i=1}^k \gamma_i \skm g_i' \mrm v_i''$ with $\gamma_i \in \myk^*$, $g_i' \in G$ and $v_i'' \in {\cal N}$.
Hence, $u \mrm g = u' \mrm (a \mrm g) = u' \mrm (\sum_{i=1}^k \gamma_i \skm g_i' \mrm v_i'') = \sum_{i=1}^k \gamma_i \skm (u' \mrm g_i') \mrm v_i''$ and now our induction hypothesis can be applied to each multiple $u' \mrm g_i'$, and since products of elements in ${\cal N}$ are again in ${\cal N}$, we are done.
Therefore, we find that $g$ has a representation  $g = \sum_{j=1}^n \alpha_j' \skm f_j' \mrm w_j'$ where $\alpha_j' \in \myk^*, f_j' \in G, w_j' \in {\cal N}$ and now we can proceed as in theorem \ref{theo.comp.arbitrary.nilpotent} to prove our claim.
\\
\qed
\procedure{Gr\"obner Bases in Nilpotent Group Rings}
{
\begin{tabbing}
XXXXX\=XXXX \kill
\removelastskip
{\bf Given:} \> $F \subseteq \myk[\g]$ and a presentation $(\Gamma,T)$ of $\g$ by
                ${\cal E}$ and ${\cal N}$ as specified above.\\
{\bf Find:} \> $\gb_{}(F)$, a  Gr\"obner basis of $\ideal{}{}(F)$.\\
\\
XX\= XX\= XXXX \= XX \= XXXX \= XX \= XXXX \= XX \= XXXX\=\kill
$G := \{ \psi_e(f) \mid f \in F, e \in {\cal E} \}$; 
{\rm\kommentar \phantom{X}\% $G$ contains $F_{\cal E}$ and $\ideal{}{}(F) = \ideal{}{}(G)$}\\
$G$ := $\bigcup_{g \in G} \s_{}(g)$;
{\rm\kommentar \phantom{X}\% $G$ is ${\cal N}$-saturated}\\
$B$ := $\{ (q_{1}, q_{2}) \mid q_{1}, q_{2} \in G, q_{1} \neq q_{2}
\}$; \\
$M$ := $\{ a \mrm f \mid f \in G, a \in \Gamma \}$;\\
{\bf while} $M \neq \emptyset$ or $B \neq \emptyset$ {\bf do} \\
\> {\bf if} \> $M \neq \emptyset$\\
\>          \> {\bf then} \> $h$ := ${\rm remove}(M)$; {\rm\kommentar
  \phantom{X}\% Remove an element using a fair strategy}\\
\>          \>            \> $h'$ := ${\rm normalform}(h, \red{}{\myr}{qc}{G})$;\\
\>           \>            \>{\bf if} \> $h' \neq 0$ \\
\>           \>            \>         \>{\bf then}    \>$G$ := $G \cup \s_{}(h')$;\\
\>           \>            \>         \>    \>{\rm\kommentar \% $G$ is ${\cal N}$-saturated  and $\ideal{}{}(F) = \ideal{}{}(G)$}\\
\>           \>            \>         \>            \>$B$ :=
$B \cup \{ (f, g) \mid f \in G, g \in \s_{}(h')
\}$; \\
\>           \>            \>          \>            \>$M$ :=
$\{ a \mrm g \mid a \in \Gamma, g \in \s_{}(h')
\}$; \\
\>\>\>{\bf endif}\\
\>{\bf endif}\\
\>  {\bf if} \> $B \neq \emptyset$ \\
\>           \> {\bf then} \> $(q_{1}, q_{2})$ := ${\rm remove}(B)$;
{\rm\kommentar \phantom{X}\% Remove an element using a fair strategy}\\
\>           \>            \>{\bf if} \> $h$ := $\spol{}(q_{1}, q_{2})$ exists  \\
\>           \>            \>         \>{\bf then}\>$h'$ := ${\rm normalform}(h, \red{}{\myr}{qc}{G})$;\\
\>           \>            \>         \>          \>{\bf if} \> $h' \neq 0$ 
{\rm\kommentar \phantom{X}\% The s-polynomial does not
  reduce to zero}\\
\>           \>            \>         \>           \>         \>{\bf then}  \>$G$ := $G \cup  \s_{}(h')$;\\
\>           \>            \>         \>\>         \>    \>{\rm\kommentar\% $G$ is ${\cal N}$-saturated  and $\ideal{}{}(F) = \ideal{}{}(G)$}\\
\>           \>            \>         \>    \>         \>        \>$B$ :=
$B \cup \{ (f, g) \mid f \in G, g \in \s_{}(h')
\}$; \\
\>           \>            \>          \>    \>         \>        \>$M$ :=
$\{ a \mrm g \mid a \in \Gamma, g \in \s_{}(
h')
\}$; \\
\>\>\>\>\>{\bf endif}\\
\>           \>            \>{\bf endif}\\
\>{\bf endif}\\
{\bf endwhile}\\
$\gb_{}(F):= G$
\end{tabbing}  
}

Correctness of this procedure follows from  corollary \ref{cor.ideal}.
For the set $G$ enumerated by this procedure  we have 
 $F_{\cal E} \subseteq G$ and the set $G$ at each stage generates 
 the  ideal $\ideal{}{}(F)$ and is ${\cal N}$-saturated.
Using a fair strategy to remove elements from the test sets $B$ and $M$ ensures
 that for all polynomials entered into $G$ the existing s-polynomials
 and the critical left multiples are considered.
To show termination we need the following theorem which makes use of Dickson's lemma due to the special representatives of the group elements.
\begin{theorem}\label{theo.term.np}~\\
{\sl
Every (right) Gr\"obner basis  contains a finite one.
}
\end{theorem}
\Ba{}
Let $F$ be a  subset of $\myk[\g]$ and $G$ a Gr\"obner
basis\footnote{The proof for the existence of a finite right Gr\"obner
  basis for $\ideal{r}{}(F)$ is similar.} of $\ideal{}{}(F)$, i.e., $\ideal{}{}(F) = \ideal{}{}(G) =
\ideal{r}{}(G)$ and for all $g \in \ideal{}{}(F)$ we have $g
\red{*}{\myr}{qc}{G} 0$.
We can assume that $G$ is infinite as otherwise we are done.
Further let
 $H = \{ \hterm(g) \mid g \in G \} \subseteq \g$.
Then for every polynomial $f \in \ideal{}{}(F)$ there exists a term
$t \in H$ such that $\hterm(f) \tupeq t$.
$H$ can be decomposed into $H = \bigcup_{e \in {\cal E}} H_e$ where $H_e$ contains those terms in $H$ starting with $e$.
For each element of $eu \in H_e$ the element $u$ then can be viewed as an n-tuple over $\z$ as it is
presented by an ordered group word.
But we can also view it as a 2n-tuple over $\n$ by representing
 each element $u \in {\cal N}$ by an extended ordered group word
 $u \id a_1^{-i_1} a_1^{j_1} \ldots  a_n^{-i_n} a_n^{j_n}$,
 where $i_l,j_l \in \n$ and the
 representing 2n-tuple is $(i_1,j_1, \ldots , i_n,j_n)$.
Notice that at most one of the two exponents $i_l$ and $j_l$ is non-zero.
Now only considering the ordered group word parts of the terms, each set $H_e$ can be seen as a (possibly infinite) subset of a free commutative
 monoid  ${\cal T}_{2n}$ with $2 \skm n$ generators.
Thus by  Dickson's lemma there exists 
 a finite subset $B_e$ of $H_e$
 such that for every $w \in H_e$ there is a $b \in B_e$ with
  $w \tupeq b$.
Now we can use the sets $B_e$ to distinguish a finite Gr\"obner basis in
$G$ as follows.
To each term $t \in B_e$ we can assign a polynomial $g_t \in G$ such
that $\hterm(g_t) = t$.
Then the set $G_B = \{ g_t \mid t \in B_e, e \in {\cal E} \}$ is again a Gr\"obner basis since for every polynomial
$f \in \ideal{}{}(F)$ there still exists a polynomial $g_t$ now in $G_B$
such that $\hterm(f) \tupeq \hterm(g_t) = t$.
Hence all polynomials in $\ideal{}{}(F)$ are qc-reducible to zero using $G_B$.
\\
\qed

Since both procedures enumerate respective Gr\"obner bases and the
 sets enumerated contain finite Gr\"obner bases, the procedures terminate
 as soon as all polynomials of the contained bases are entered into $G$.
Therefore we now are able to solve problems related to right and two-sided
 ideals in nilpotent group rings using reduction similar to Buchberger's
 approach to commutative polynomial rings.

%% file: reductionrings.tex
\chapter{Monoid Rings over Reduction Rings}\label{chapter.reductionrings}
\spruch{10.8}{7}{More matter, with less art.}{Hamlet}

In this chapter we want to outline how the ideas of
 completion in a monoid ring over a field can be modified and carried
 over to monoid rings over\index{reduction ring}\index{ring!reduction}
 reduction rings as defined by Madlener in \cite{Ma86}.
The special case of monoid rings over the integers has been
 studied explicitly in \cite{MaRe93} and will be used to provide an example
 when  following  the more general approach given by Kapur and
 Narendran in \cite{KaNa85}.

Let $\rr$ be a commutative ring with  a reduction $\Longrightarrow_B$
associated with  subsets $B \subseteq\rr$ satisfying the following
axioms
\begin{enumerate}
\item[(A1)] $\Longrightarrow_B\; = \bigcup_{\beta \in B}
  \Longrightarrow_{\beta}$, $\Longrightarrow_B$ is terminating for all
  subsets $B \subseteq\rr$.
\item[(A2)] $\alpha \Longrightarrow_{\beta} \gamma$ implies $\alpha - \gamma
  \in \ideal{}{\rr}(\beta)$, i.e., $\gamma = \alpha - \beta \skm \rho$
  for some $\rho \in \rr$.
\item[(A3)] $\alpha \Longrightarrow_{\alpha} 0$ for all $\alpha \in \rr
  \backslash \{ 0 \}$.
\end{enumerate}
Notice that $0$ has to be irreducible for all $\R_{\alpha}$, $\alpha \in \rr$.
Therefore, $0$ will be chosen as the normal form of the ideal elements.
Let us recall the definition of G-bases (Gr\"obner bases)
 for ideals in such a ring $\rr$.
\begin{definition}~\\
{\rm
A finite subset $B$ of  $\rr$ is called a 
 \betonen{G-basis} of an ideal $\mswab{i}$, if
 $\red{*}{\Longleftrightarrow}{}{B} = \;\;\equiv_{\mswab{i}}$
 and $\red{}{\R}{}{B}$ is confluent.
\dend
}
\end{definition}
$\rr$ is called a \index{reduction ring}\betonen{reduction ring}
 if every finitely generated ideal has a G-basis.
It is often useful, if $\rr$ satisfies an additional axiom
strongly related to interreduction.
\begin{enumerate}
\item[(A4)] $\alpha \Longrightarrow_{\beta}$ and $\beta\Longrightarrow_{\gamma} \delta$ imply $\alpha \Longrightarrow_{\gamma}$ or $\alpha \Longrightarrow_{\delta}$.
\end{enumerate}
In the following we will always assume that the reduction ring fulfills
 the axioms (A1) to (A4).
\begin{lemma}\label{lem.BB'}~\\
{\sl 
Let $B \subseteq \rr$ be a  G-basis
 and $B' \subseteq B$ such that for all $\beta \in B$, $\beta \red{*}{\R}{}{B'} 0$ holds.
Then for all $\alpha \in \rr$, $\alpha \red{*}{\R}{}{B} 0$ implies
 $\alpha \red{*}{\R}{}{B'} 0$.
In particular, $B'$ is a G-basis of $\ideal{}{\rr}(B)$.
\lemend
}
\end{lemma}
\Ba{}~\\
Notice that by axiom (A4) and our assumptions on $B'$, all elements reducible
 using $B$ are also reducible using $B'$ and
 $\irr(\R_{B'}) \subseteq \irr(\R_B)$.
Assuming $\alpha \red{*}{\R}{}{B} 0$ but $\alpha \red{*}{\R}{}{B'} \alpha\rnf{B'} \neq 0$ we find $\alpha\rnf{B'} \in \irr(\R_{B'}) \subseteq \irr(\R_B)$ and
 $\alpha\rnf{B'} \in \ideal{}{\rr}(B)$, contradicting that $B$ is
 a  G-basis of $\ideal{}{\rr}(B)$.
\\
In particular, as $B$ is a G-basis we have
 $\alpha \red{*}{\R}{}{B'} \alpha\rnf{B}$ since
 $\red{*}{\Longleftrightarrow}{}{B'} \subseteq
 \red{*}{\Longleftrightarrow}{}{B} = \;\;
 \equiv_{\ideal{}{\rr}(B)}$ and $\R_B$ is confluent.
This implies that $\R_{B'}$ is also confluent, as $\alpha\rnf{B}$ is unique.
Now it remains to show that
 $\red{*}{\Longleftrightarrow}{}{B} \subseteq 
 \red{*}{\Longleftrightarrow}{}{B'}$ holds.
This follows immediately, as for 
 $\alpha \red{*}{\Longleftrightarrow}{}{B} \beta$ the confluence 
 of $\R_B$ yields
 $\alpha\rnf{B'} = \alpha\rnf{B} = \beta\rnf{B}  = \beta\rnf{B'}$ 
 which implies $\alpha \red{*}{\Longleftrightarrow}{}{B'} \beta$.
\\
\qed
Given a cancellative\footnote{In case we allow arbitrary monoids we
 have to be more careful in defining right reduction and critical situations corresponding
 to it.} monoid $\m$, 
 we call $\rr[\m]$ the monoid ring over $\rr$.
Using an appropriate ordering $\succ$ on the monoid, we can proceed as
 described in chapter \ref{chapter.reduction} and specify
 reduction in $\rr[\m]$ depending on reduction in $\rr$.
We additionally define a partial ordering on $\rr$ by setting for 
 $\alpha, \beta \in \rr$, $\alpha >_{\rr} \beta$ if and only if
 there exists a  set $B \subseteq \rr$ such that $\alpha
 \red{+}{\Longrightarrow}{}{B} \beta$.
Then we can define a Noetherian ordering on $\rr[\m]$ as follows:
 for $f,g \in \rr[\m]$, $f > g$ if and only if either
 $\hterm(f) \succ \hterm(g)$ or $(\hterm(f) = \hterm(g)$ and 
 $\hc(f) >_{\rr} \hc(g))$ or $(\hm(f) = \hm(g)$ and 
 $\reductum(f) > \reductum(g))$.
Notice that this ordering in general is not total 
 on $\rr[\m]$.
\begin{definition}\label{def.redrr}~\\
{\rm
Let $p, f$ be two non-zero polynomials in $\rr[\m]$. \\
We say $f$ 
 \index{reduction!right!in reduction rings}\index{right!reduction in reduction rings}\betonen{right reduces} $p$ to $q$ at a monomial 
 $\alpha \skm t$ in $p$ in one step, denoted by $p \red{}{\myr}{r}{f} q$, if
\begin{enumerate}
\item[(a)] $\hterm(f \mrm w) = \hterm(f) \mm w  = t$ for some $w \in \m$,
\item[(b)] $\alpha \Longrightarrow_{\hc(f)} \beta$  with  $\alpha = \gamma \skm \hc(f) +\beta$ for some $\beta, \gamma \in \rr$, and
\item[(c)] $q = p - \gamma \skm f \mrm w$.
\end{enumerate}
We write $p \red{}{\myr}{r}{f}$ if there is a polynomial $q$ as defined
above and $p$ is then called right reducible by $f$. 
Further we can define $\red{*}{\myr}{r}{}, \red{+}{\myr}{r}{}$,
 $\red{n}{\myr}{r}{}$ as usual.
Right reduction by a set $F \subseteq \rr[\m]$ is denoted by
 $p \red{}{\myr}{r}{F} q$ and abbreviates $p \red{}{\myr}{r}{f} q$
 for some $f \in F$,
 which is also written as  $p \red{}{\myr}{r}{f \in F} q$.
\dend
}
\end{definition}
Notice that in case $f$ right reduces $p$ to $q$ at a monomial $\alpha
\skm t$ this no longer implies $t \not\in\terms(q)$.
But when using a  set of polynomials  for reduction we know by (A1)
that reducing $\alpha$ in $\rr$ with respect to the head coefficients of
the applicable polynomials must terminate and then either the term $t$
disappears or is not further reducible.  
Hence the so-defined right reduction is Noetherian.

Analogous definitions can be introduced for strong, prefix and the other
 notions of reduction considered previously.
As before, for right reduction 
 $\alpha \skm p \mrm w \red{}{\myr}{r}{p} 0$ need not hold.
Hence, we have to introduce the concept of saturation as in section \ref{section.rightreduction}, extended
 to include possible problems caused by multiplication with coefficients.
\begin{definition}~\\
{\rm
A set of polynomials $F \subseteq \{\alpha \skm p \mrm w \mid \alpha
\in \rr^*, w \in\m \}$ is called
 a \betonen{saturating set}\/ for a polynomial $p\in \rr[\m]$,
 if for all $\alpha \in \rr^*$, $w \in \m$,  $\alpha \skm p \mrm w
 \red{}{\myr}{r}{F} 0$ holds in case $\alpha \skm p \mrm
 w \neq 0$.
Let $\SAT(p)$ denote the family of all saturating sets for $p$.
A set  $F$ of polynomials in $\rr[\m]$ is called \index{saturated set}\betonen{saturated},
 if  $\alpha \skm f \mrm w \red{}{\myr}{r}{F} 0$ holds for
 all $f \in F$ and all  $\alpha \in \rr^*$, $w \in \m$ in case $\alpha
 \skm f \mrm w \neq 0$.
\dend
}
\end{definition}
\begin{definition}~\\
{\rm
A  set $G \subseteq \rr[\m]$ is called a \betonen{Gr\"obner basis}\/
 with respect to
 the reduction $\red{}{\myr}{r}{}$ or a 
 \betonen{right (or stable) Gr\"obner basis}, if
\begin{enumerate}
\item[(i)] $\red{*}{\lr}{r}{G} = \;\; \equiv_{\ideal{r}{}(G)}$, and
\item[(ii)] $\red{}{\myr}{r}{G}$ is confluent.
\dend
\end{enumerate}
}
\end{definition}
We will assume that for the reduction ring $\rr$  there is an
 algorithm to compute   
 G-bases and a representation of 
 the elements of these bases in terms of the elements of the generating
 set.
Further we require that it is possible to compute a finite basis for
a module of solutions to linear homogeneous equations over $\rr$.
These sets are necessary to describe ``overlaps'' due to coefficients
 that we will need to characterize right Gr\"obner bases in the monoid
 ring later on.
Next we give
 a criterion for checking whether a given basis of a right ideal is a
 right Gr\"obner basis.
We start by defining special polynomials for finite subsets of
 polynomials, called G-polynomials and M-polynomials.
\begin{definition}\label{def.gpol}~\\
{\rm
Let $P = \{ p_1, \ldots, p_k \}$ be a set of polynomials in $\rr[\m]$ and
 $t$ an element in $\m$
 such that there are $w_1, \ldots, w_k \in \m$ with
 $\hterm(p_i \mrm w_i) = \hterm(p_i) \mm w_i = t$, for all $1 \leq i \leq k$.
Further let $\gamma_i = \hc(p_i)$ for  $1 \leq i \leq k$\footnote{Note that this definition would have to be modified for non-cancellative monoids, as then $\hterm(p \mrm w) = \hterm(p) \mm w$ does not imply $\hc(p \mrm w) = \hc(p)$.}.
\\
Let $\{ \alpha_1, \ldots, \alpha_n \}$ be a 
 G-basis of $\{ \gamma_1, \ldots, \gamma_k \}$ and
$$ \alpha_i = \beta_{i,1} \skm \gamma_1 + \ldots + \beta_{i,k} \skm \gamma_k $$
 for $\beta_{i,j} \in \rr$,  $1 \leq i \leq n$, and $1 \leq j \leq k$.
Notice that the $\alpha_i$ respectively the $\beta_{i,j}$ do not depend on $t$.
Then we define the 
 \index{G-polynomial!Gr\"obner polynomial}
 \index{Gr\"obner polynomial!G-polynomial}
 \betonen{G-polynomials (Gr\"obner polynomials)}
 corresponding to $P$ and $t$ by setting
$$ g_i = \sum_{j=1}^k \beta_{i,j} \skm p_j \mrm w_j
 \mbox{ for each } 1 \leq i \leq k.$$
Notice that $\hm(g_i)= \alpha_i \skm t$.
\\
For the module $M = \{ (\delta_1, \ldots, \delta_k) \mid  
 \sum_{i=1}^k \delta_i \skm \gamma_i = 0 \}$, let the set
 $\{A_i \mid 1 \leq i \leq r \}$ be a basis with
 $A_i = (\alpha_{i,1}, \ldots, \alpha_{i,k})$ for $\alpha_{i,j} \in \rr$,
 $1 \leq i \leq r$, and $1 \leq j \leq k$.
Notice that the $A_i$ do not depend on $t$.
Then we define the 
 \index{M-polynomial!module polynomial}\index{module polynomial!M-polynomial}\betonen{M-polynomials (module polynomials)}
 corresponding to $P$ and $t$ by setting
$$ m_i = \sum_{j=1}^k \alpha_{i,j} \skm p_j \mrm w_j
 \mbox{ for each } 1 \leq i \leq r.$$
Notice that $\hterm(m_i) \prec t$ for each $1 \leq i \leq r$.
\dend
}
\end{definition}
Given a set of polynomials $F$, the set of corresponding
 G- and M-polynomials contains those which are specified by
 definition \ref{def.gpol} for each finite subset $P \subseteq F$
 and each term $t \in \m$ fulfilling the respective conditions.
For a set consisting of one polynomial the corresponding
M-polynomials reflect the multiplication of the polynomial with
zero-divisors of the head coefficient, i.e., by a basis of the annihilator
 of the head coefficient.
This case is also treated by the idea of saturation.
Notice that given a finite set of polynomials the corresponding sets of
 G- and M-polynomials in general can be infinite (e.g., reviewing example
 \ref{exa.kippen}, for the polynomials $p = a+f$ and $q=b + \lambda$ infinitely
 many critical situations at the terms $ba^i$, $i \in \n^+$
 have to be considered).

We can use G- and M-polynomials to characterize right
Gr\"obner bases in monoid rings over a reduction ring in case they are
 additionally saturated.
\begin{theorem}\label{theo.rr.cp}~\\
{\sl
For a saturated subset $F$ of $\rr[\m]$ the following statements are equivalent:
\begin{enumerate}
\item For all polynomials $g \in \ideal{r}{}(F)$ we have $g
       \red{*}{\myr}{r}{F} 0$.
\item All G-polynomials and all M-polynomials corresponding to $F$
       right reduce to zero using $F$.
\end{enumerate}
\theoend
}
\end{theorem}
\Ba{}~\\
\mbox{$1 \R 2:$  }
This follows  from the fact that all G-polynomials
 and M-polynomials belong to the right ideal generated by $F$.

\mbox{$2 \R 1:$ }
We have to show that every element $g \in \ideal{r}{}(F) \backslash \{ 0 \}$
 is right reducible to zero using $F$.
Remember that for
 $h \in \ideal{r}{}(F)$, $ h \red{}{\myr}{r}{F} h'$ implies $h' \in \ideal{r}{}(F)$.
Thus as  $\red{}{\myr}{r}{F}$ is Noetherian
 it suffices to show that every 
 $g \in \ideal{r}{}(F) \backslash \{ 0 \}$ is right reducible
 using $F$.
This will be done by assuming the contrary.
Let $g = \sum_{j=1}^m \gamma_{j} \skm f_{j} \mrm w_{j}$ with $\gamma_{j} \in \rr^*, f_j \in F, w_{j} \in \m$ be a representation
 of a polynomial $g\in \ideal{r}{}(F) \backslash \{ 0 \}$.
As $F$ is saturated, we can always assume 
 $\hterm(\gamma_i \skm f_{i} \mrm w_{i}) = \gamma_i \skm \hterm(f_{i}) \mm w_{i}$.
Depending on this  representation of $g$ and
 the well-founded total ordering $\succeq$ on $\m$ we define
 the {\bf critical term} of $g$ to be
 $t = \max \{ \hterm(f_{j}) \mm w_{j} \mid j \in \{ 1, \ldots m \}  \}$.
We call another representation of $g$ ``smaller'' if for the
 corresponding critical term $\tilde{t}$ we have $\tilde{t} \prec t$.
Let us assume that our polynomial $g$ is not right reducible
 by $F$ and  that our representation of $g$ is a minimal one
 with respect to $t$.
We have to distinguish two cases:
In case $t \neq \hterm(g)$ without loss of generality let us assume
 that $t$ occurs in the first $k$ products of our representation.
Hence, we have $\hterm(f_i) \mm w_i = t$ for each $1 \leq i \leq k$ and
 $\sum_{i=1}^{k} \gamma_i \skm \hc(f_i) = 0$, i.e., the vector
 $(\gamma_1, \ldots, \gamma_k)$ is in the module
 $M = \{ (\alpha_1, \ldots, \alpha_k) \mid \sum_{i=1}^k \alpha_i \skm 
         \hc(f_i) = 0 \}$.
By our assumption this module has been considered when generating
 the M-polynomials for $\{ f_1, \ldots, f_k \}$ and $t$.
Let the set $\{ A_i = (\alpha_{i,1}, \ldots, \alpha_{i,k}) \mid 1 \leq i \leq n \}$
 be a basis of $M$.
Then for $1 \leq i \leq k$ we have $\gamma_i = \sum_{j=1}^n
\alpha_{j,i} \skm \delta_j$ for some $\delta_j \in R$.
Thus, we get
\begin{eqnarray}
 & & \sum_{i=1}^k \gamma_i \skm f_i \mrm w_i \nonumber\\
 &=& \sum_{i=1}^k (\sum_{j=1}^n \alpha_{j,i} \skm \delta_j) \skm f_i \mrm w_i \nonumber\\
 &=& \sum_{j=1}^n \delta_j \skm (\sum_{i=1}^k \alpha_{j,i} \skm f_i \mrm w_i)
     \label{reductionrings}
\end{eqnarray}
Taking a closer look at the last sum of these transformations
 in \ref{reductionrings},
 we see that we can express the sum of the first $k$ elements
 of our representation of $g$ by a sum of M-polynomials.
Since these M-polynomials belonging to  $\{ f_1, \ldots, f_k \}$ and
 $t$ all have head terms smaller than $t$ and are all right reducible
 to zero using $F$, we get a new representation of $g$ with a critical
 term smaller than $t$, contradicting our assumption that our
 chosen representation was minimal.
\\
In case $t = \hterm(g)$ we know that there exists a finite subset
 $P \subseteq F$ such that
 $\hc(g) \in \ideal{}{\rr}(\{ \hc(p) \mid p \in P \})$, and as this ideal
 is finitely generated it has a G-basis, say $G_P$.
Then $\hc(g)$ is reducible by an element $\alpha \in G_P$.
By our assumption now $\alpha \skm t$ is head monomial of a G-polynomial
 corresponding to $P$ and $t$ and since this G-polynomial is
 right reducible to zero using $F$, in particular
 there exist polynomials
 $f_1, \ldots, f_k \in F$ in volved in the reduction of $\alpha \skm t$
 such that $\alpha \R_{\hc(f_1)} \alpha_1 
 \R_{\hc(f_2)} \ldots \R_{\hc(f_k)} 0$. 
By lemma \ref{lem.BB'} this implies that  $\hm(g)$ is right reducible
 using $F$.
\auskommentieren{
We will now show that this implies that $\hm(g)$ is right reducible
 using $F$ by induction on $k$.
For $k=1$ we find $\hc(g) \R_{\alpha}$ and $\alpha \R_{\hc(f_1)} 0$,
 and hence axiom (A4) implies $\hc(g) \R_{\hc(f_1)}$, i.e., $g$ is
 right reducible at $\hm(g)$ using $f_1 \in F$.
Now let $k > 1$. 
Then by axiom (A4), $\hc(g)$ is either reducible by $\hc(f_1)$ or by
 $\alpha_1$.
This again gives us a contradiction, as either $g$ is right reducible
 at $\hm(g)$ using $f_1 \in F$ or the induction hypothesis can
 be applied to $\alpha_1$ and $f_2, \ldots, f_k$.}
\\
\qed

In case we additionally require that in $\rr$ every ideal is
 finitely generated we can even show a stronger result
 in the second part of this proof.
Given a set of polynomials $F$ and an element $s \in \m$, let
 $$C_{F}(s)= \{ \hc(f) \mid f \in F, s = \hterm(f) \mm z \mbox{ for some } z
 \in \m, \hterm(f \mrm z) = \hterm(f) \mm z \}.$$
Notice that if $F$ is finite then $C_F(s)$ is also finite.
Then the following lemma holds.

\begin{lemma}\label{lem.rr.congruence}~\\
{\sl
Let $\rr$ be a reduction ring such that every ideal is finitely generated.
Further let $F$ be a saturated set of polynomials fulfilling the conditions of theorem \ref{theo.rr.cp}.
Then for every $s \in \m$, the corresponding set $C_F(s)$ contains a
 G-basis of $\ideal{}{\rr}(C_{F}(s))$.
\lemend
}
\end{lemma}
\Ba{}~\\
Since $\ideal{}{\rr}(C_{F}(s))$ is finitely generated, there exists
 a finite set $A \subseteq C_{F}(s)$ such that 
 $\ideal{}{\rr}(A) = \ideal{}{\rr}(C_{F}(s))$.
Let $G_A$ be the finite  G-basis
 of this ideal.
Furthermore let $P \subseteq
 F_s = \{ f \in F \mid \mbox{ there exists an element } z 
 \in \m \mbox{ such that } \hterm(f \mrm z) = \hterm(f) \mm z = s \}$
 be a finite set of polynomials such that
 $A = \{ \hc(f) \mid f \in P, s = \hterm(f) \mm z=\hterm(f \mrm z)
 \mbox{ for some } z \in \m \}$.
Then, as for every $\alpha \in G_A$ there exists a G-polynomial corresponding
 to $P$ and $s$, by our assumption there exists a finite
 set $C \subset C_{F}(s)$ such that for all $\alpha \in G_A$ we have
  $\alpha \red{*}{\R}{}{C} 0$.
Moreover, since $C \subseteq C \cup G_A$ and the latter is a  G-basis
 of $\ideal{}{\rr}(A) = \ideal{}{\rr}(C_{F}(s))$,
 by lemma \ref{lem.BB'} $C$ is also a G-basis of the same ideal. 
\\
\qed

This lemma also holds in case $F$ is finite as then all ideals
 $\ideal{}{\rr}(C_{F}(s))$ in this proof are finitely generated, namely
 by the finite sets $C_F(s)$.

\begin{corollary}\label{lem.rridealcongruence}~\\
{\sl
Let $F$ be a saturated set fulfilling the conditions specified in
 theorem \ref{theo.rr.cp}.
In case $F$ is finite or $\rr$ is Noetherian, then $F$ is a right G-basis.
}
\end{corollary}
\Ba{}~\\
The inclusion $\red{*}{\lr}{r}{F} \subseteq \;\;\equiv_{\ideal{r}{}(F)}$
 is obvious.
Hence let us assume $f \equiv_{\ideal{r}{}(F)} g$, i.e., $f-g \in \ideal{r}{}(F)$ and, therefore,  $f-g \red{*}{\myr}{r}{F} 0$.
We show that this implies $f \red{*}{\lr}{r}{F} g$.
In case $f-g = 0$ we are immediately done.
Hence let us assume $f-g \neq 0$ and as any polynomial in $\ideal{r}{}(F)$
 is right reducible to zero using $F$, without loss of generality
 we can assume that the reduction sequence 
 $f-g \red{*}{\myr}{r}{F} 0$ uses top-reduction, i.e., all reductions take
 place at the respective head monomial.
Further let $t = \hterm(f-g)$ and let $\gamma_1$ respectively $\gamma_2$ be the
 coefficients of $t$ in $f$ respectively $g$.
We will now show that $f \red{*}{\lr}{r}{F} g$ holds  by induction on
 the term $t=\hterm(f-g)$.
In case $t = \lambda$ we find $f - g = \gamma_1 - \gamma_2 \red{*}{\myr}{r}{F} 0$ and as described in lemma \ref{lem.rr.congruence} there exists a subset
 $C \subseteq C_F(\lambda) = \{ \hc(f) \mid f \in F, \lambda = \hterm(f \mrm z) =
 \hterm(f) \mm z \mbox{ for some } z \in \m \}  \subseteq \rr$
 which is
 a  G-basis of $\ideal{}{\rr}(C_F(\lambda))$
 and then, as $\gamma_1 - \gamma_2 \in \ideal{}{\rr}(C_F(\lambda))$,
 $\gamma_1 - \gamma_2 \red{*}{\Longrightarrow}{}{C} 0$
  implies 
 $\gamma_1  \red{*}{\Longleftrightarrow}{}{C} \gamma_2$.
Using the respective
 polynomials belonging to the elements in $C$, we get $f \red{*}{\lr}{r}{F} g$.
Now let us assume $t \succ \lambda$ and $f-g \red{k}{\myr}{r}{F} h$
 where $k \in \n^+$ is minimal such that $\hterm(h) \neq t$.
Further let $ f_1, \ldots, f_k$ be the polynomials used in the respective
 reduction steps, i.e.,
 $\gamma_1 - \gamma_2 \red{*}{\Longrightarrow}{}{\{\hc(f_i) \mid 1 \leq i \leq k \}} 0$. 
Again there exists a set $C \subseteq C_F(t)$ which
 is a  G-basis of $\ideal{}{\rr}(C_F(t))$ and without loss of
 generaltity let $\{\hc(f_i) \mid 1 \leq i \leq k \} \subseteq C$. 
Then $\gamma_1 - \gamma_2 \red{*}{\Longrightarrow}{}{C} 0$ implies 
 $\gamma_1  \red{*}{\Longleftrightarrow}{}{C} \gamma_2$.
Now applying the polynomial multiples belonging to the elements
 of $C$ used in this last sequence 
 $\gamma_1  \red{*}{\Longleftrightarrow}{}{C} \gamma_2$
 to the monomial with term $t$ in
 $f$ we find an element
 $\tilde{f} \in \rr[\m]$ such that $f \red{*}{\lr}{r}{F} \tilde{f}$,
 $\hm(\tilde{f}) = \gamma_2 \skm t$, $\tilde{f} - g \in \ideal{r}{}(F)$,
 and $t \succ \hterm(\tilde{f} - g)$.
Hence our induction hypothesis yields $g \red{*}{\lr}{r}{F} \tilde{f}
  \red{*}{\lr}{r}{F} f$ and we are done. 
\\
It  remains to show that right reduction using $F$ is confluent.
Suppose there is a polynomial $g$ having two distinct normal
 forms with respect to $F$, say $p_1$ and $p_2$.
Let $t$ be the largest term on which $p_1$ and $p_2$ differ and
 let $\alpha_1$ respectively $\alpha_2$ be the coefficients of $t$ in $p_1$
 respectively $p_2$.
Since $p_1 - p_2 \in \ideal{r}{}(F)$ we know
 $p_1 - p_2 \red{*}{\myr}{r}{F} 0$ and 
 $\alpha_1 - \alpha_2 \in \ideal{}{\rr}(C_{F}(t))=\ideal{}{\rr}(C)$,
 where $C\subseteq C_{F}(t)$ is a  G-basis.
Hence, $\alpha_1 \red{*}{\Longleftrightarrow}{}{C} \alpha_2$, and
 either $\alpha_1$ or $\alpha_2$
 must be reducible using $C$, i.e., not both $p_1$ and
 $p_2$ can be in normal form with respect to $F$, contradicting our
 assumption.
\\
\qed
In \cite{De89} Dei{\ss} has shown that in case $\rr$ allows the
computation of Gr\"obner bases using pairs of polynomials, then
 Gr\"obner bases in commutative polynomial rings over $\rr$ can
also be characterized by G- and M-polynomials of pairs of
polynomials.
When working over the integers it is even possible to restrict
oneself to {\em one} special overlap called the s-polynomial.
We have done a similar characterization for monoid rings over the integers
 in \cite{MaRe93a} which is sketched in the following example.
Notice, that still a pair of polynomials can give rise to an infinite set of
 s-polynomials.


\begin{example}~\\
{\rm
We  give a definition of  $\z$ as a
 reduction ring fulfilling the axioms (A1) to (A4) by defining a
 reduction relation as follows:
\\
First we give a total ordering on $\z$ by $\alpha <_{\z} \beta$
 if and only if ($\alpha \geq 0$ and $\beta < 0$) or
 ($\alpha \geq 0$, $\beta > 0$ and $\alpha < \beta$) or
 ($\alpha < 0$, $\beta < 0$ and $\alpha > \beta$).
For $\alpha \in \z$ we call the elements $\rho$ with
 $0 \leq \rho < |\alpha|$ the remainders of $\alpha$. 
Reduction is now specified as $\alpha \R_{\beta} \gamma$ if and only if
 $\alpha \geq_{\z} \beta$, $\alpha = \delta \skm \beta + \gamma$
  and $\gamma$ is a remainder of $|\beta|$.
The axioms (A1) to (A3) are easily checked.
(A4) holds since $\alpha \R_{\beta}$ and $\beta \R_{\gamma} \delta$
 imply $\alpha  \geq_{\z} \beta  \geq_{\z} \gamma  \geq_{\z} \delta$
 and hence $\alpha$ is reducible by $\gamma$ as well as by $\delta$.

Notice that $\z$ contains pairs of non-zero associated
  elements\footnote{We call $\alpha$ and
  $\beta$ associated in case there exists a unit $\epsilon$ such that
  $\alpha = \epsilon \skm \beta$} which  only differ in sign.
We will call the positive element of such a pair the canonical element
 and define right reduction by restricting ourselves to the use of 
 polynomials with canonical head coefficients for reduction only (compare
 \cite{KaKa84}).

A polynomial $f$ right reduces a non-zero polynomial $p$ to $q$ at
a monomial $\alpha \cd t$ in $p$ in one step, denoted by $p \red{}{\myr}{r}{f} q$, if
\begin{enumerate}
\item[(a)] $\hterm (f \mrm w) = \hterm(f) \mm w  = t$ for some $w \in \m$.
\item[(b)] $\hc(f)>0$ and $\alpha = \gamma \skm \hc(f)+
  \beta$ for $\gamma, \beta \in\z$, $\gamma \neq 0$, $\beta$ a remainder of $\hc(f)$.
\item[(c)] $q = p - \gamma \skm f \mrm w$.
\end{enumerate}
Given two polynomials $p_{1}, p_{2} \in \z[\m]$ with $HT(p_{i}) = t_{i}$, $i = 1, 2$.
If there are $w_{1}, w_{2} \in \m$ with $t_{1} \mm w_{1} = t_{2} \mm 
 w_{2} =t$,
 and $\alpha_1,\alpha_2$ are the non-zero coefficients of $t$ in $p_1 \mrm w_1$
 respectively $p_2 \mrm w_2$, then
If $\alpha_2 \geq \alpha_1 >0$ and $\alpha_{2} = \beta\skm \alpha_{1} + \gamma$, where 
 $\beta, \gamma \in \z$, $ \gamma$ a remainder of $\alpha_1$, we get the following s-polynomial
 $$\spol{}{}(p_{1}, p_{2}, w_{1}, w_{2}) = \beta \skm p_1 \mrm w_1 - p_2 \mrm w_2.$$
Let $U_{\hm(p_1),\hm(p_2)}  \subseteq \m \times \m$ be the set containing all
pairs $w_1,w_2 \in\m$ as above.

Now we can characterize right Gr\"obner bases as follows:
\\
For a saturated set of polynomials $F$ in $\z[\m]$,
equivalent are:
\begin{enumerate}
\item $F$ is a  right Gr\"obner basis.
\item $\ideal{r}{}(F) \red{*}{\myr}{r}{F} 0$.
\item For all not necessarily different $f_{k}, f_{l} \in F, (w_{k}, w_{l}) \in U_{\hm(f_k),\hm(f_l)}$  we have 
      $\spol{}{}(f_{k}, f_{l}, w_{k}, w_{l}) \red{*}{\myr}{r}{F} 0$.
\exaend
\end{enumerate}
}
\end{example}
The existence of finite Gr\"obner bases for the special classes of monoids
and groups shown in the previous chapters can be transfered to
corresponding monoid rings over reduction rings.
This is on one hand due to the fact that with respect to prefix,
 commutative or quasi-commutative reduction, for a finite set of polynomials
 we can localize the corresponding G- and M-polynomials to 
 finitely many critical situations.
On the other hand the property of having a reduction ring as coefficient domain
 ensures the finiteness of the respective Gr\"obner bases.
Using the necessary procedures for calculations in $\rr$ one can
modify the given procedures as it has been done for the case of the
integers.

Let us close this section with a short remark on other possible
 definitions of ``Gr\"obner'' bases of ideals
 in reduction rings.
Recall that in order to decide the membership or the congruence problem
 of an ideal $\mswab{i}$ in a reduction ring it is sufficient to have
 a finite basis $G$ of $\mswab{i}$ such that for all $\alpha \in \mswab{i}$
 we have $\alpha \red{*}{\Longrightarrow}{}{G} 0$.
This definition of special ideal bases was used e.g. by Pan (\cite{Pa85})
 or Kapur and Narendran (\cite{KaNa85}).
Notice that for such a basis
 $\red{*}{\Longleftrightarrow}{}{G} = \;\;\equiv_{\mswab{i}}$ in 
 general need not hold.
In case one uses these bases to define G-polynomials, theorem \ref{theo.rr.cp}
 no longer characterizes right Gr\"obner bases since we cannot
 guarantee $\red{*}{\lr}{r}{F} = \;\;\equiv_{\ideal{}{r}(F)}$.

Another type of ideal bases studied in reduction rings are the so-called
 weak Gr\"obner bases.
\begin{definition}~\\
{\rm
A finite subset $B$ in $\rr$ is called a 
 \betonen{weak G-basis} of an  ideal $\mswab{i}$, if 
 $\red{*}{\Longleftrightarrow}{}{B} = \;\;\equiv_{\mswab{i}}$ and every element in the
  ideal  reduces to zero using $B$.
\dend
}
\end{definition}
\begin{example}~\\
{\rm
To give an example for a weak Gr\"obner basis, let us define reduction 
 using $4$ on
 $\z$ as follows:
 For $n > 7$ set $n \R_{4} n-4$, for $n<-4$ set $n \R_{4} n+4$ and further set
 $7 \R_{4} 3$, $7 \R_{4} -1$, $6 \R_{4} 2$, $6 \R_{4} -2$, $5 \R_{4} 1$,
 $5 \R_{4} -3$, $4 \R_{4} 0$, $-4 \R_{4} 0$.
Then we have that $\red{*}{\Longleftrightarrow}{}{\{ 4 \}} = \;\;\equiv_{\ideal{}{\rr}(4)}$ and  all elements in $\ideal{}{\rr}(4)$ reduce to zero.
But although $3 - (-1) = 4 \R_{4} 0$, $3$ and $-1$ are not joinable,
 i.e., the translation lemma does not hold and $\R_{\{ 4 \}}$ is not confluent.
\exaend
}
\end{example}
Although we now have that the weak G-bases in $\rr$ describe the ideal
 congruence, this does not carry over to the monoid ring: when using
 these bases to define G-polynomials, theorem \ref{theo.rr.cp} does not characterize
 right Gr\"obner bases, as lemma \ref{lem.BB'} no longer holds.
We still have
\begin{lemma}~\\
{\sl 
Let $B \subseteq \rr$ be a weak  G-basis
 and $B' \subseteq B$ such that for all $\beta \in B$, $\beta \red{*}{\R}{}{B'} 0$ holds.
Then for all $\alpha \in \rr$, $\alpha \red{*}{\R}{}{B} 0$ implies
 $\alpha \red{*}{\R}{}{B'} 0$.
\lemend
}
\end{lemma}
But we cannot show that $B'$ is a weak G-basis, as the axioms, especially (A2)
 and (A4), do not provide enough information on the reduction step to
 capture the ideal congruence.
Still, special bases can be characterized as in theorem \ref{theo.rr.cp}
 which allow to solve the membership and congruence problem in case they
 are finite.

%% file: conclusions.tex
\chapter{Concluding Remarks}
\spruch{8.5}{8}{So eine Arbeit wird eigentlich nie fertig.}{Goethe}

The aim of this thesis was to introduce reduction and the concept of
 Gr\"obner bases to monoid and group rings.

Since finitely generated ideals in general already do not have finite
 Gr\"obner bases in the free monoid or free group ring due to the fact
 that this would solve the word problem for monoids respectively
 groups, we have restricted our studies to right ideals except for the
 class of commutative monoids where  right ideals are of course ideals and 
 for the class of nilpotent groups.
Furthermore, as general monoid rings need not be right Noetherian, we
 have tried to localize monoids and groups where finitely generated
 right ideals have finite Gr\"obner bases.

In order to introduce reduction to a monoid ring, a well-founded
 ordering on the monoid elements is needed.
Such an ordering in general cannot be compatible with the monoid
 multiplication and hence 
 one of the main problems turned out to be that  monomial right
  multiples of a polynomial need no longer be reducible to zero
 by the polynomial itself.
Furthermore, Gr\"obner bases cannot be characterized by head terms as
 in the case of polynomial rings, solvable polynomial rings or
 free monoid rings.
This phenomenon of having to deal with non-stable orderings influenced
 the concepts of reduction we studied:
 strong reduction allowed to use all right multiples of a polynomial as
 rules,
 right reduction restricted the right multiples to the ``stable'' ones and
 in defining prefix respectively commutative\footnote{for commutative monoids
  only} or quasi-commutative\footnote{for nilpotent groups only}
 reduction additional
 syntactical restrictions for the multiplication of the head term of
 the polynomial were added.
While the weakenings of strong reduction on one side gave more
 information on the reduction steps, the expressiveness of the right
 ideal congruence by the reflexive transitive symmetric closure of such
 a weaker reduction was lost.
To recover this property of reduction the concept of saturation was
 used to enrich the sets of polynomials used for reduction.
But in general for none  of these reductions $p \red{*}{\myr}{}{F} 0$ implies
 $\alpha \skm p \mrm w \red{*}{\myr}{}{F} 0$,  which is a crucial lemma
 in Buchberger's original approach of characterizing
 Gr\"obner bases by s-polynomials and
 in most of the generalizations known in literature.
In fact, in defining prefix respectively commutative reduction, a 
 weaker condition to characterize Gr\"obner bases by corresponding
 s-polynomials could be proved.
These reductions are motivated  syntactically and correspond to
 polynomial multiples, where the
 multiplication of the head term can be interpreted by a multiplication not
 only in the monoid but also in the set of syntactical elements from which the representatives for
 the monoid elements are taken, e.g., a free monoid respectively a
 free commutative monoid.
These monoids now allow admissible well-founded orderings and in using
 such an ordering to induce the ordering on the respective monoid, we
 can characterize Gr\"obner bases by localized s-polynomials.
This approach resulted in a procedure to enumerate a prefix Gr\"obner
 basis for a finitely generated ideal in a monoid ring.
The procedure could be modified to give terminating algorithms for
 special classes of monoids and groups, e.g., finite monoids, free
 monoids, free groups, plain groups and
 context-free groups.
Similarly, for commutative monoid rings a terminating algorithm could
 be provided.

The key idea of introducing reduction and Gr\"obner bases to other
structures used in this work is as follows:
\begin{enumerate}
\item Define a weakening of strong reduction, say w-reduction,
  appropriate to the respective structure in the following sense:
  \\
  If for some polynomials $p, g \in \myk[\m]$ and a set of polynomials
  $F \subseteq \myk[\m]$ we have $p \red{}{\myr}{w}{g} 0$ and $g
  \red{*}{\myr}{w}{F} 0$, then there exists a w-representation of $p$
  such that one term in this representation equals the head term of
  $p$ and all other terms are smaller\footnote{Variations of this
    lemma are e.g. the lemmata \ref{lem.redp} and \ref{lem.redc}.}.
\item Define saturation with respect to w-reduction.
\item Define s-polynomials with respect to w-reduction.
\end{enumerate}
Then the following holds:

{\sl
For a w-saturated set $F \subset \myk[\m]$ the following statements are
equivalent:
\begin{enumerate}
\item For all polynomials $g \in \ideal{r}{}(F)$ we have $g \red{*}{\myr}{w}{F} 0$.
\item For all polynomials $f_{k}, f_{l} \in F$ we have 
  $\spol{w}(f_{k}, f_{l}) \red{*}{\myr}{w}{F} 0$.
\end{enumerate}
}

A further application of this approach is given in the section on
nilpotent groups, where we define quasi-commutative reduction and show
the existence of finite Gr\"obner bases for right ideals and ideals.

We close by giving possible further points of interest:
On the theoretical side, the class of polycyclic groups is a
 good candidate to extend our approach of generalizing Gr\"obner bases.
It is known from literature that polycyclic groups have solvable 
 subgroup problem \cite{Wi89} and that polycyclic group rings have
 solvable membership problem \cite{Ha59,BaCaMi81}.
On the practical side, the  algorithms provided in this thesis
should be improved, implemented and complexity bounds should be
investigated. A possible guide for the latter task can be found in the
literature on the subgroup problem for the respective classes of groups
\cite{AvMa84,CrOt94,KuMa89}.

%% file: references.tex
\addcontentsline{toc}{chapter}{Bibliography}
\bibliographystyle{alpha}

%% file: index.tex
\addcontentsline{toc}{chapter}{Index}
\printindex

%% file: lebenslauf.tex
\pagestyle{empty}
{\large\bf Curriculum vitae}
\\
\\
\\
\\
\begin{tabbing}
XXXXXXXXXXXXXXXXXX\= \kill\\
{\bf Pers{\"o}nliche Daten} \\
\\
Name \> Birgit Reinert, geborene Weber \\
Geburtsdatum \> 6.11.1964  \\
Geburtsort \>  Landstuhl \\
Familienstand \> verheiratet \\
Staatsangeh\"origkeit \> deutsch \\
\\
\\
{\bf Schulausbildung} \\
\\
 08/70 -- 12/71   \> Grundschule in El Paso/Texas, USA\\
 01/72 -- 03/73   \> Grundschule in Huntsville/Alabama, USA\\ 
 04/73 -- 07/75   \> Grundschule in Spesbach/Pfalz\\
 08/75 -- 06/84   \>Staatliches Gymnasium in Landstuhl \\
                  \>  Abschlu{\ss}: Abitur\\
\\
\\      
{\bf Hochschulstudium} \\ 
\\ 
 10/84 -- 09/89 \>   Mathematik mit Nebenfach Informatik an der\\
                \>   Universit{\"a}t Kaiserslautern \\
                \>   Abschlu{\ss}: Diplom-Mathematiker\\
\\
\\
{\bf Berufst\"atigkeit}  \\ 
\\
 10/89 -- 03/95 \>Wissenschaftlicher Mitarbeiter in der  Arbeitsgruppe\\
                \> ``Grundlagen der Informatik'' von 
                   Prof.~Dr.~K. Madlener \\
                \>Universit\"at Kaiserslautern 
\end{tabbing}